\documentclass[a4,10pt]{article}
\usepackage{amsmath,amscd,amssymb}
\usepackage{mathrsfs}

\newcommand{\nbiga}{\mathcal{A}}
\newcommand{\nbigb}{\mathcal{B}}
\newcommand{\nbigc}{\mathcal{C}}
\newcommand{\nbigd}{\mathcal{D}}
\newcommand{\nbige}{\mathcal{E}}
\newcommand{\nbigf}{\mathcal{F}}
\newcommand{\nbigg}{\mathcal{G}}
\newcommand{\nbigh}{\mathcal{H}}
\newcommand{\nbigi}{\mathcal{I}}
\newcommand{\nbigj}{\mathcal{J}}
\newcommand{\nbigk}{\mathcal{K}}
\newcommand{\nbigl}{\mathcal{L}}
\newcommand{\nbigm}{\mathcal{M}}
\newcommand{\nbign}{\mathcal{N}}
\newcommand{\nbigo}{\mathcal{O}}
\newcommand{\nbigp}{\mathcal{P}}
\newcommand{\nbigq}{\mathcal{Q}}
\newcommand{\nbigr}{\mathcal{R}}
\newcommand{\nbigs}{\mathcal{S}}
\newcommand{\nbigt}{\mathcal{T}}
\newcommand{\nbigu}{\mathcal{U}}
\newcommand{\nbigv}{\mathcal{V}}

\newcommand{\nbigx}{\mathcal{X}}
\newcommand{\nbigy}{\mathcal{Y}}
\newcommand{\nbigz}{\mathcal{Z}}

\newcommand{\proj}{\mathbb{P}}
\newcommand{\seisuu}{{\mathbb Z}}
\newcommand{\rnum}{{\mathbb Q}}

\newcommand{\cnum}{{\mathbb C}}
\newcommand{\real}{{\mathbb R}}

\newcommand{\hyperr}{\mathbb{R}}

\newcommand{\newTate}{\pmb{T}}

\newcommand{\DD}{\mathbb{D}}

\newcommand{\gbigd}{\mathfrak D}

\newcommand{\gbigp}{\mathfrak P}
\newcommand{\gbigq}{\mathfrak Q}

\newcommand{\gbigs}{\mathfrak S}
\newcommand{\gbigt}{\mathfrak T}

\newcommand{\gbigv}{\mathfrak V}

\newcommand{\gminib}{\mathfrak b}

\newcommand{\gminis}{\mathfrak s}

\newcommand{\gminiv}{\mathfrak v}

\newcommand{\vecxi}{{\boldsymbol \xi}}
\newcommand{\veceta}{{\boldsymbol \eta}}

\newcommand{\vecgamma}{{\boldsymbol \gamma}}

\newcommand{\veczero}{{\boldsymbol 0}}
\newcommand{\vecone}{{\boldsymbol 1}}
\newcommand{\vecalpha}{{\boldsymbol \alpha}}
\newcommand{\veca}{{\boldsymbol a}}
\newcommand{\vecb}{{\boldsymbol b}}
\newcommand{\vecbeta}{{\boldsymbol \beta}}
\newcommand{\vecdelta}{{\boldsymbol \delta}}

\newcommand{\vecc}{{\boldsymbol c}}

\newcommand{\veck}{{\boldsymbol k}}
\newcommand{\vecm}{{\boldsymbol m}}

\newcommand{\vecp}{{\boldsymbol p}}

\newcommand{\vecS}{{\boldsymbol S}}

\newcommand{\larr}{\leftarrow}
\newcommand{\llarr}{\longleftarrow}

\newcommand{\rarr}{\rightarrow}
\newcommand{\lrarr}{\longrightarrow}
\newcommand{\darr}{\downarrow}

\newcommand{\pf}{{\bf Proof}\hspace{.1in}}
\newcommand{\qed}{\mbox{\rule{1.2mm}{3mm}}}

\def\Hom{\mathop{\rm Hom}\nolimits}

\def\Cok{\mathop{\rm Cok}\nolimits}

\def\Image{\mathop{\rm Im}\nolimits}

\def\Re{\mathop{\rm Re}\nolimits}

\def\Gr{\mathop{\rm Gr}\nolimits}
\def\gr{\mathop{\rm Gr}\nolimits}

\def\Tot{\mathop{\rm Tot}\nolimits}

\def\rank{\mathop{\rm rank}\nolimits}
\def\Spec{\mathop{\rm Spec}\nolimits}

\def\Ker{\mathop{\rm Ker}\nolimits}
\def\ker{\mathop{\rm Ker}\nolimits}

\def\Gr{\mathop{\rm Gr}\nolimits}
\def\Sym{\mathop{\rm Sym}\nolimits}

\def\Res{\mathop{\rm Res}\nolimits}

\def\ord{\mathop{\rm ord}\nolimits}

\def\tr{\mathop{\rm tr}\nolimits}

\def\id{\mathop{\rm id}\nolimits}

\def\codim{\mathop{\rm codim}\nolimits}

\def\Supp{\mathop{\rm Supp}\nolimits}

\newcommand{\del}{\partial}
\newcommand{\delbar}{\overline{\del}}

\newcommand{\nhom}{{\mathcal Hom}}

\newcommand{\mbar}{\underline{m}}

\newcommand{\sankaku}{\triangle}

\newcommand{\barlambda}{\overline{\lambda}}
\newcommand{\lambdabar}{\barlambda}

\newcommand{\varphibar}{\overline{\varphi}}
\newcommand{\etabar}{\overline{\eta}}

\newcommand{\lefttop}[1]{{}^{#1}\!}

\def\reg{\mathop{\rm reg}\nolimits}

\newcommand{\tildepsi}{\widetilde{\psi}}
\newcommand{\psitilde}{\tildepsi}

\newcommand{\deldel}{\eth}

\newcommand{\distribution}{\gbigd\gminib}

\newcommand{\nbigvecs}{\boldsymbol{\mathcal S}}
\newcommand{\vecnbigs}{\boldsymbol{\mathcal S}}
\newcommand{\vecnbigp}{\boldsymbol{\mathcal P}}

\newcommand{\supp}{\gminis}

\newcommand{\mutilde}{\widetilde{\mu}}

\newcommand{\htilde}{\widetilde{h}}

\newcommand{\pitilde}{\widetilde{\pi}}

\newcommand{\ftilde}{\widetilde{f}}

\newcommand{\Ftilde}{\widetilde{F}}

\newcommand{\atilde}{\widetilde{a}}
\newcommand{\stilde}{\widetilde{s}}

\newcommand{\Psitilde}{\widetilde{\Psi}}

\newcommand{\Rtilde}{\widetilde{R}}

\newcommand{\Omegabar}{\overline{\Omega}}

\newcommand{\Dtilde}{\widetilde{D}}

\newcommand{\nbigrtilde}{\widetilde{\nbigr}}

\newcommand{\nbigmtilde}{\widetilde{\nbigm}}

\def\ord{\mathop{\rm ord}\nolimits}

\def\Hol{\mathop{\rm Hol}\nolimits}

\def\exchange{\mathop{\rm exchange}\nolimits}

\def\DR{\mathop{\rm DR}\nolimits}

\def\Ch{\mathop{\rm Ch}\nolimits}

\def\Conv{\mathop{\rm Conv}\nolimits}
\def\FL{\mathop{\rm FL}\nolimits}
\def\loc{\mathop{\rm loc}\nolimits}
\def\aff{\mathop{\rm aff}\nolimits}

\def\MTM{\mathop{\rm MTM}\nolimits}
\def\Eff{\mathop{\rm Eff}\nolimits}
\def\QDM{\mathop{\rm QDM}\nolimits}
\def\vir{\mathop{\rm vir}\nolimits}
\def\ev{\mathop{\rm ev}\nolimits}
\def\Proj{\mathop{\rm Proj}\nolimits}
\def\sp{\mathop{\rm sp}\nolimits}

\def\GKZ{\mathop{\rm GKZ}\nolimits}
\def\MHM{\mathop{\rm MHM}\nolimits}

\newcommand{\Wtilde}{\widetilde{W}}
\newcommand{\Ptilde}{\widetilde{P}}

\newcommand{\gtilde}{\widetilde{g}}

\newcommand{\gbigstilde}{\widetilde{\gbigs}}
\newcommand{\varphitilde}{\widetilde{\varphi}}
\newcommand{\Ctilde}{\widetilde{C}}

\newcommand{\Stilde}{\widetilde{S}}

\newcommand{\Ytilde}{\widetilde{Y}}

\newcommand{\gbigvtilde}{\widetilde{\gbigv}}

\newcommand{\nrhom}{R{\mathcal Hom}}
\newcommand{\DDD}{\boldsymbol D}

\newcommand{\Htilde}{\widetilde{H}}
\newcommand{\Omegatilde}{\widetilde{\Omega}}

\newcommand{\gammatilde}{\widetilde{\gamma}}

\newcommand{\vecq}{\boldsymbol q}
\newcommand{\nbigatilde}{\widetilde{\nbiga}}

\newcommand{\nbigctilde}{\widetilde{\nbigc}}

\newcommand{\nbigptilde}{\widetilde{\nbigp}}
\newcommand{\Zbar}{\overline{Z}}

\newcommand{\nbigubar}{\overline{\nbigu}}

\newcommand{\QDMbar}{\overline{\QDM}}

\newcommand{\Btilde}{\widetilde{B}}
\newcommand{\gbigptilde}{\widetilde{\gbigp}}
\newcommand{\chitilde}{\widetilde{\chi}}

\newcommand{\Bbar}{\overline{B}}

\newcommand{\distributionbar}{\overline{\distribution}}

\newcommand{\distributiontilde}{\widetilde{\distribution}}
\newcommand{\Thetatilde}{\widetilde{\Theta}}
\newcommand{\DDDtilde}{\widetilde{\DDD}}
\newcommand{\Thetabar}{\overline{\Theta}}
\newcommand{\nbigcbar}{\overline{\nbigc}}
\newcommand{\trtilde}{\widetilde{\tr}}

\newcommand{\MTMint}{\MTM^{\rm int}}
\newtheorem{thm}{Theorem}[section]
\newtheorem{cor}[thm]{Corollary}

\newtheorem{rem}[thm]{Remark}
\newtheorem{lem}[thm]{Lemma}
\newtheorem{prop}[thm]{Proposition}
\newtheorem{df}[thm]{Definition}

\newtheorem{example}[thm]{Example}

\newtheorem{assumption}[thm]{Assumption}

\begin{document}

\title{Twistor property of
GKZ-hypergeometric systems}
\author{Takuro Mochizuki}
\date{}
\maketitle

\begin{abstract}
We study the mixed twistor $\nbigd$-modules
associated to meromorphic functions.
In particular, we describe
their push-forward and specialization
under some situations.
We apply the results to study the twistor property of
a type of better behaved GKZ-hypergeometric systems,
and to study their specializations.
As a result, we obtain some isomorphisms 
of mixed TEP-structures in the local mirror symmetry.

{\footnotesize
\vspace{.1in}

}

\end{abstract}

\section{Introduction}
\subsection{Mixed twistor $\nbigd$-modules}

In \cite{s3},
C. Simpson introduced the concept of mixed twistor structure
as a generalization of mixed Hodge structure.
According to his principle, called Simpson's meta theorem,
most objects and theorems concerning 
with mixed Hodge structure
should have their counterparts in the context of
mixed twistor structure.

The concept of mixed Hodge module due to M. Saito \cite{saito1,saito2}
is one of the ultimate in the Hodge theory.
Roughly speaking,
mixed Hodge modules are regular holonomic $\nbigd$-modules
equipped with mixed Hodge structure.
One of the main results in the theory is 
the functoriality of mixed Hodge modules.
Namely, we have the standard operations 
on the derived category of algebraic mixed Hodge modules
such as push-forward, pull back, duality, inner homomorphism,
tensor product, nearby cycle functor, vanishing cycle functors,
which are compatible with the standard operations
for algebraic regular holonomic $\nbigd$-modules.

According to Simpson's meta theorem,
we should have a twistor version of mixed Hodge modules,
that is the concept of mixed twistor $\nbigd$-module.
The concept of pure twistor $\nbigd$-module
was introduced by C. Sabbah \cite{sabbah2,sabbah5},
and studied by himself and the author
\cite{Mochizuki-tame, Mochizuki-wild}.
The mixed case was studied 
in \cite{Mochizuki-MTM}.

As in the Hodge case,
we have the standard functors on the derived category 
of algebraic mixed twistor $\nbigd$-modules,
that is one of the most useful points in the theory of
mixed twistor $\nbigd$-modules.
Another interesting point is that 
we have the mixed twistor $\nbigd$-modules
associated to meromorphic functions.

Let $f$ be a meromorphic function on a complex manifold $X$
whose poles are contained in a hypersurface $H$.
Let $L_{\ast}(f,H)$ be the holonomic $\nbigd$-modules
given by $\nbigo_X(\ast H)$ with the flat connection $d+df$.
We also have the holonomic $\nbigd$-module
$L_!(f,H):=\DDD(L_{\ast}(-f,H))$,
where $\DDD$ denotes the duality of holonomic $\nbigd$-modules.
We have the natural mixed twistor $\nbigd$-modules
over $L_{\star}(f,H)$ $(\star=\ast,!)$.
Applying the standard functors to such mixed twistor $\nbigd$-modules,
we can obtain many mixed twistor $\nbigd$-modules.
In other words,
many important holonomic $\nbigd$-modules
are equipped with a natural mixed twistor structure.
It would be interesting to have applications of 
the mixed twistor structure.

In this paper,
we shall explain that
a type of better behaved GKZ-hypergeometric systems
are naturally equipped with the mixed twistor structure.
We shall also explain an application
in the study of toric local mirror symmetry.

\subsubsection{Ingredients for mixed twistor $\nbigd$-modules}

Originally, ingredients for twistor $\nbigd$-modules
are given as $\nbigr$-triples by Sabbah \cite{sabbah2}.
(See \S\ref{subsection;15.5.19.1} for a review.)
But, in this paper,
we consider only integrable mixed twistor $\nbigd$-modules
with real structure.
So, in this introduction,
we explain it in a slightly different way,
although we use the formalism of $\nbigrtilde$-triples
after \S\ref{section;15.11.22.10}.

Let $X$ be any complex manifold.
Let $\nbigx:=\cnum\times X$
and $\nbigx^0:=\{0\}\times X$.
Let $\lambda$ be the standard coordinate of $\cnum$.
We have the sheaf of differential operators
$\nbigd_{\nbigx}$ on $\nbigx$.
Let $\nbigrtilde_X\subset \nbigd_{\nbigx}$ 
denote the sheaf of subalgebras generated by
$\lambda\Theta_{\nbigx}(\log\nbigx^0)$.

An $\nbigrtilde_X$-module $\nbigm$ is called strict
if it is flat over $\nbigo_{\cnum}$.
Any coherent $\nbigrtilde_X$-module $\nbigm$ has 
the characteristic variety $\Ch(\nbigm)$
in $T^{\ast}\nbigx(\log\nbigx^0)$
as in the case of $\nbigd$-modules.
We say that $\nbigm$ is holonomic
if $\Ch(\nbigm)$ is contained in
$\cnum\times \Lambda$
for a Lagrangian subvariety
$\Lambda\subset T^{\ast}X$.
A real structure of strict holonomic $\nbigrtilde_X$-module
$\nbigm$ is a $\real$-perverse sheaf
$P_{\real}$ on $\cnum^{\ast}\times X$
with an isomorphism
$P_{\real}\otimes\cnum
\simeq
 \DR_{\cnum^{\ast}\times X}
 (\nbigm_{|\cnum^{\ast}\times X})$.

Then, an integrable mixed  twistor $\nbigd$-module
with real structure on $X$
is a strict holonomic $\nbigrtilde_X$-module $\nbigm$
with a real structure and a weight filtration
$(\nbigm,P_{\real},W)$
satisfying some conditions.
Here, $W$ is an increasing filtration
of strict holonomic $\nbigrtilde_X$-modules
with real structure.

In this introduction, 
``mixed twistor $\nbigd$-module''
means 
``integrable mixed twistor $\nbigd$-module
with real structure''.

\subsection{Better behaved GKZ-hypergeometric systems}
\label{subsection;15.5.16.20}

\subsubsection{$\nbigd$-modules}

L. Borisov and P. Horja \cite{Borisov-Horja}
introduced the concept of better behaved GKZ-hypergeometric system
as a generalization of GKZ-hypergeometric systems 
\cite{GKZ}.
(See also \cite{Iritani-quantum-cohomology-period}.)
Let us recall a special type of better behaved GKZ-hypergeometric
systems.

Let $\nbiga=\{\veca_1,\ldots,\veca_m\}\subset \seisuu^n$
be a finite subset generating $\seisuu^n$.
We have the cone
$K_{\real}(\nbiga):=\bigl\{
 \sum_{j=1}^m r_j\veca_j\,\big|\,
 r_j\geq 0
 \bigr\}
\subset\real^n$
generated by $\nbiga$.
We set $K(\nbiga):=K_{\real}(\nbiga)\cap\seisuu^n$.
The semigroup $K(\nbiga)$ is the saturation of
$\seisuu_{\geq 0}\nbiga
=\bigl\{
\sum_{j=1}^m n_{j}\veca_j\,\big|\, 
 n_{j}\in\seisuu_{\geq 0}
\bigr\}$.
We denote 
$\veca_j=(a_{j1},\ldots,a_{jn})$.

Let $\Gamma\subset K(\nbiga)$ be any subset
such that
$\Gamma+\veca\subset \Gamma$ for any $\veca\in\nbiga$.
Let $\vecbeta\in\cnum^n$.
The following system of differential equations
$\GKZ(\nbiga,\Gamma,\vecbeta)$
for tuples of holomorphic functions 
$(\Phi_{\vecc}\,|\,\vecc\in \Gamma)$
on any open subset in $\cnum^m$ 
is called the better behaved GKZ-hypergeometric system
associated to $(\nbiga,\Gamma,\vecbeta)$:
\[
 \del_{x_j}\Phi_{\vecc}=\Phi_{\vecc+\veca_j}
\quad(\forall \vecc\in \Gamma,\,\forall j=1,\ldots,m)
\]
\[
 \Bigl(
 \sum_{j}a_{ji}x_j\del_{x_j}+c_i-\beta_i
 \Bigr)
 \Phi_{\vecc}=0
\quad
 (\forall \vecc\in \Gamma,\,\forall i=1,\ldots,n)
\]
Here, $(x_1,\ldots,x_m)$ denotes 
the standard coordinate system
of $\cnum^m$.

\begin{rem}
In {\rm\cite{Borisov-Horja}},
any tuple in a finitely generated abelian group is 
considered 
instead of a finite subset in $\seisuu^n$.
In that sense,
the above system is a special case of
better behaved GKZ-hypergeometric systems.
Later, moreover, we shall mainly consider the case $\vecbeta=0$.

But, as in {\rm \cite{Adolphson}},
we omit the existence of an element
$\vecalpha\in(\seisuu^n)^{\lor}$
such that $\vecalpha(\veca_i)=1$ for any $i$.
\hfill\qed
\end{rem}

\paragraph{Expression as the twisted Gauss-Manin systems}

We can describe the corresponding $\nbigd$-modules
using the twisted Gauss-Manin system.
Let $T^n:=(\cnum^{\ast})^n$.
We consider the morphism
$\psi^{\aff}_{\nbiga}:T^n\lrarr \cnum^m$
given by 
$\psi^{\aff}_{\nbiga}(t_1,\ldots,t_n)
=(t^{\veca_1},\ldots,t^{\veca_m})$,
where $t^{\veca_j}=\prod_{i=1}^nt_i^{a_{ji}}$.
Let $X_{\nbiga}^{\aff}$ denote the closure
of the image of $\psi^{\aff}_{\nbiga}$.
Let $\check{X}^{\aff}_{\nbiga}\lrarr X^{\aff}_{\nbiga}$
be the normalization.
Let $\check{D}^{\aff}_{\nbiga}$
denote the complement of $T^n$ in
$\check{X}^{\aff}_{\nbiga}$.
Note that 
$\check{X}^{\aff}_{\nbiga}
=\Spec \cnum[K(\nbiga)]$
and 
$X^{\aff}_{\nbiga}
=\Spec\cnum[\seisuu_{\geq 0}\nbiga]$.

Let $\Omega^k_{\check{X}^{\aff}_{\nbiga}}(\log
\check{D}^{\aff}_{\nbiga})$
denote the sheaf of meromorphic differential $k$-forms
which are logarithmic along $\check{D}^{\aff}_{\nbiga}$,
studied in \cite{Batyrev-VMHS}.
Let $\Omega^k_{(\check{X}^{\aff}_{\nbiga},\check{D}^{\aff}_{\nbiga})}$
denote the sheaf of holomorphic $k$-forms on 
the normal variety $\check{X}^{\aff}_{\nbiga}$
which are $0$ along $\check{D}^{\aff}_{\nbiga}$,
studied in \cite{Danilov-de-Rham}.
Let $q:\check{X}^{\aff}_{\nbiga}\times\cnum^m\lrarr
 \check{X}^{\aff}_{\nbiga}$ be the projection.
We set
\[
 \Omega^k_{\check{X}^{\aff}_{\nbiga}\times\cnum^m/\cnum^m}
 (\log \check{D}^{\aff}_{\nbiga}\times\cnum^m):=
 q^{\ast}\Omega^k_{\check{X}^{\aff}_{\nbiga}}(\log
 \check{D}^{\aff}_{\nbiga}),
\quad\quad
  \Omega^k_{
 (\check{X}^{\aff}_{\nbiga},\check{D}^{\aff}_{\nbiga})
 \times\cnum^m
 /\cnum^m}:=
 q^{\ast}\Omega^k_{(\check{X}^{\aff}_{\nbiga},\check{D}^{\aff}_{\nbiga})}.
\]

The family of Laurent polynomials 
$\sum_{j=1}^m x_jt^{\veca_j}$
induces a meromorphic function
$F_{\nbiga}$ on 
$(\check{X}^{\aff}_{\nbiga},\check{D}^{\aff}_{\nbiga})\times\cnum^m$.
We also have the logarithmic closed one form
$\kappa(\vecbeta)=\sum_{i=1}^n \beta_i dt_i/t_i$.
We obtain the relative algebraic de Rham complexes:
\[
 \nbigc^{\bullet}(\nbiga,\vecbeta)_{\bullet}:=
 \Bigl(
 \Omega^{\bullet}_{(\check{X}^{\aff}_{\nbiga}\times \cnum^m)/\cnum^m}
 (\log\check{D}^{\aff}_{\nbiga}\times\cnum^m),
 d+dF_{\nbiga}-\kappa(\vecbeta)
 \Bigr)
\]
\[
 \nbigc^{\bullet}(\nbiga,\vecbeta)_{\circ}:=
 \Bigl(
 \Omega^{\bullet}_{(\check{X}^{\aff}_{\nbiga},\check{D}^{\aff}_{\nbiga})
 \times \cnum^m/\cnum^m},
 d+dF_{\nbiga}-\kappa(\vecbeta)
 \Bigr)
\]

Let $\check{\pi}_{\nbiga}:
 \check{X}^{\aff}_{\nbiga}\times\cnum^m
\lrarr \cnum^m$ denote the projection.
Each $\check{\pi}_{\nbiga\ast}\nbigc^{k}(\nbiga,\vecbeta)_{\star}$
is naturally a $\nbigd_{\cnum^m}$-module
by $\del_{x_j}\cdot(g)=\del_{x_j}g+(\del_{x_j}F_{\nbiga}) g$.
The differential $d+dF_{\nbiga}-\kappa(\vecbeta)$ of
the complexes are compatible
with the actions of $\nbigd_{\cnum^m}$.
So, we obtain the following $\nbigd_{\cnum^m}$-modules
\[
 M_{\nbiga,\vecbeta,\star}:=
 \hyperr^n\check{\pi}_{\nbiga\ast}
 \nbigc^{\bullet}(\nbiga,\vecbeta)_{\star}
\quad
 (\star=\bullet,\circ).
\]

Let $K(\nbiga)^{\circ}$ be the intersection of
$K(\nbiga)$ and the interior part of $K_{\real}(\nbiga)$.
The $\nbigd_{\cnum^m}$-modules
$M_{\nbiga,\vecbeta,\bullet}$
and $M_{\nbiga,\vecbeta,\circ}$
correspond to the systems
$\GKZ(\nbiga,K(\nbiga),\vecbeta)$
and $\GKZ(\nbiga,K(\nbiga)^{\circ},\vecbeta)$
in the sense of the following lemma.
(See Lemma \ref{lem;15.5.15.1}
and Proposition \ref{prop;15.4.25.2}
for more details.)
\begin{prop}
Let $U$ be any open subset in $\cnum^m$.
We have a natural bijective correspondence
between solutions of 
$\GKZ(\nbiga,K(\nbiga),\vecbeta)$
(resp.
 $\GKZ(\nbiga,K(\nbiga)^{\circ},\vecbeta)$)
and $\nbigd_{U}$-homomorphisms
$M_{\nbiga,\vecbeta,\bullet}\lrarr
 \nbigo_U$
(resp.
$M_{\nbiga,\vecbeta,\circ}\lrarr
 \nbigo_U$).
\end{prop}

In this paper,
we are mainly interested in the case $\vecbeta=0$.

\paragraph{Expressions as the push-forward of $\nbigd$-modules}

We take a toric resolution
$\varphi_{\Sigma_1}:X_{\Sigma_1}
 \lrarr \check{X}^{\aff}_{\nbiga}$,
i.e.,
$X_{\Sigma_1}$ be a smooth toric variety,
and $\varphi_{\Sigma_1}$ 
is a projective morphism which is $T^n$-equivariant.
Let $D_{\Sigma_1}:=X_{\Sigma_1}\setminus T^n$.
We have the meromorphic algebraic function
$F_{\nbiga,\Sigma_1}$ 
on $(X_{\Sigma_1},D_{\Sigma_1})\times\cnum^m$
given by $F_{\nbiga}$.
We have the algebraic holonomic 
$\nbigd_{X_{\Sigma_1}\times\cnum^m}$-modules
$L_{\star}(F_{\nbiga,\Sigma_1},D_{\Sigma_1}\times \cnum^m)$
($\star=\ast,!$).
Let $\pi_{\Sigma_1}:
 X_{\Sigma_1}\times\cnum^m\lrarr \cnum^m$
denote the projection.
By the results in \S\ref{subsection;15.11.22.11}
and \S\ref{subsection;14.11.23.221},
we have the following proposition,
which gives expressions of $\GKZ(\nbiga,K(\nbiga),0)$
and $\GKZ(\nbiga,K(\nbiga)^{\circ},0)$
in terms of $\nbigd$-modules.

\begin{prop}
\label{prop;15.5.15.10}
We have natural isomorphisms
\[
 \pi^0_{\Sigma_1\ast}
 L_{\ast}(F_{\nbiga,\Sigma_1},D_{\Sigma_1}\times\cnum^m)
\simeq
  \pi^0_{\Sigma_1!}
 L_{\ast}(F_{\nbiga,\Sigma_1},D_{\Sigma_1}\times\cnum^m)
\simeq
 M_{\nbiga,0,\bullet},
\]
\[
\pi^0_{\Sigma_1\ast}
 L_{!}(F_{\nbiga,\Sigma_1},D_{\Sigma_1}\times\cnum^m)
\simeq
 \pi^0_{\Sigma_1!}
 L_{!}(F_{\nbiga,\Sigma_1},D_{\Sigma_1}\times\cnum^m)
\simeq
 M_{\nbiga,0,\circ}.
\]
Here, 
$\pi^0_{\Sigma_1,\star}$ $(\star=\ast,!)$
denote the $0$-th cohomology
of the push-forward functors
of algebraic $\nbigd$-modules 
with respect to $\pi_{\Sigma_1}$.
\end{prop}

\paragraph{Special case}

For any $\vecp=(p_1,\ldots,p_m)\in \seisuu^m$,
we put $\supp_+(\vecp):=\{j\,|\,p_j\geq 0\}$
and $\supp_-(\vecp):=\{j\,|\,p_j\leq 0\}$.
We set
\[
 \square_{\vecp}
=\prod_{j\in\supp_+(\vecp)}\del_{x_j}^{p_j}
-\prod_{j\in\supp_-(\vecp)}\del_{x_j}^{-p_j}.
\]
We have the morphism
$\seisuu^m\lrarr \seisuu^n$
given by
$\vecp=(p_1,\ldots,p_m)\longmapsto
\sum_{j=1}^m p_j\veca_j$.
Let $L_{\nbiga}$ denote the kernel.

For $\vecc_0\in K(\nbiga)$,
we have the following ordinary GKZ-hypergeometric system
$\GKZ^{\ord}(\nbiga,\vecbeta-\vecc_0)$
for holomorphic functions $\Phi_{\vecc_0}$
on any open subset of $\cnum^m$:
\[
 \square_{\vecp}\Phi_{\vecc_0}=0
\quad
 (\forall\vecp\in L_{\nbiga})
\]
\[
 \Bigl(c_{0i}-\beta_i+\sum_{j=1}^ma_{ji}x_j\del_{x_j}\Bigr)
 \Phi_{\vecc_0}=0
\quad
 (i=1,\ldots,n)
\]
For any $\vecgamma\in\cnum^n$,
let $I(\nbiga,\vecgamma)$ denote
the left ideal of $\nbigd_{\cnum^m}$
generated by 
$\square_{\vecp}$ $(\vecp\in L_{\nbiga})$
and $-\gamma_i+\sum_{j=1}^ma_{ji}x_j\del_{x_j}$
$(i=1,\ldots,n)$.
Then, for $\Gamma=\vecc_0+\seisuu_{\geq 0}\nbiga$,
we have a natural isomorphism
$M^{\GKZ}(\nbiga,\Gamma,\vecbeta)
\simeq
 \nbigd_{\cnum^m}/I(\nbiga,\vecbeta-\vecc_0)$
which is the $\nbigd$-module corresponding
to the system $\GKZ^{\ord}(\nbiga,\vecbeta-\vecc_0)$.

Suppose that $K(\nbiga)=\seisuu_{\geq 0}\nbiga$.
Then, we have a natural isomorphism
\begin{equation}
 \label{eq;15.5.15.30}
 M_{\nbiga,0,\bullet}
\simeq
 \nbigd_{\cnum^m}/I(\nbiga,0).
\end{equation}
Suppose moreover 
$K(\nbiga)^{\circ}=\seisuu_{\geq 0}\nbiga+\vecc_0$
for an element $\vecc_0\in K(\nbiga)$.
Then, we have a natural isomorphism
\begin{equation}
 \label{eq;15.5.15.31}
 M_{\nbiga,0,\circ}
\simeq
 \nbigd_{\cnum^m}/I(\nbiga,-\vecc_0).
\end{equation}
For an expression 
$\vecc_0=\sum_{i=1}^m b_i\veca_i$,
we have the following commutative diagram:
\[
 \begin{CD}
 M_{\nbiga,0,\circ}
 @>{a}>>
 M_{\nbiga,0,\bullet}
 \\
 @V{\simeq}VV @V{\simeq}VV \\
 \nbigd_{\cnum^m}\big/I(\nbiga,-\vecc_0)
 @>{b}>>
 \nbigd_{\cnum^m}\big/I(\nbiga,0).
 \end{CD}
\]
Here, the vertical morphisms are
(\ref{eq;15.5.15.30}) and (\ref{eq;15.5.15.31}),
$a$ is the natural morphism,
and $b$ is the induced by the multiplication of 
$\prod \del_{x_i}^{b_i}$.

\subsubsection{Mixed twistor $\nbigd$-modules}

By the geometric expression in Proposition \ref{prop;15.5.15.10},
we naturally obtain mixed twistor $\nbigd$-modules
$\gbigt_{\nbiga,0,\star}$
over $M_{\nbiga,0,\star}$ $(\star=\bullet,\circ)$.
Namely, 
we have the algebraic mixed twistor $\nbigd$-modules
$\nbigt_{\star}(F_{\nbiga,\Sigma_1},D_{\Sigma_1}\times\cnum^m)$
$(\star=\ast,!)$ on $X_{\Sigma_1}\times \cnum^m$
whose underlying $\nbigd$-modules are
$L_{\star}(F_{\nbiga,\Sigma_1},D_{\Sigma_1}\times\cnum^m)$.
We obtain the algebraic mixed twistor $\nbigd$-modules
\[
 \gbigt_{\nbiga,0,\bullet}
=\pi^0_{\Sigma_1\ast}
 \nbigt_{\ast}(F_{\nbiga,\Sigma_1},D_{\Sigma_1}\times\cnum^m)
\simeq
 \pi^0_{\Sigma_1!}
 \nbigt_{\ast}(F_{\nbiga,\Sigma_1},D_{\Sigma_1}\times\cnum^m)
\]
\[
 \gbigt_{\nbiga,0,\circ}
=\pi^0_{\Sigma_1\ast}
 \nbigt_{!}(F_{\nbiga,\Sigma_1},D_{\Sigma_1}\times\cnum^m)
\simeq
 \pi^0_{\Sigma_1!}
 \nbigt_{!}(F_{\nbiga,\Sigma_1},D_{\Sigma_1}\times\cnum^m)
\]
Let $\nbigm_{\nbiga,0,\star}$
denote the underlying $\nbigrtilde_{\cnum^m}$-modules
of $\gbigt_{\nbiga,0,\star}$ $(\star=\ast,!)$.
Let us describe them as systems of differential equations.

Let $\Gamma\subset K(\nbiga)$ be any subset 
satisfying $\Gamma+\veca\subset \Gamma$
for any $\veca\in\nbiga$. 
We consider the following system of differential equations
$\GKZ_{\nbigrtilde}(\nbiga,\Gamma,\vecbeta)$
for a tuple $\Phi_{\Gamma}=(\Phi_{\vecc}\,|\,\vecc\in \Gamma)$
of holomorphic functions on any open subset of 
$\cnum\times\cnum^m
=\{(\lambda,x_1,\ldots,x_m)\}$:
\[
 \lambda\del_{x_j}\Phi_{\vecc}
=\Phi_{\vecc+\veca_j},
\quad
(\forall \vecc\in \Gamma,\,\,j=1,\ldots,m)
\]
\[
 \Bigl(
 \lambda^2\del_{\lambda}
+n\lambda
+\sum_{j=1}^m\lambda x_j\del_{x_j}
 \Bigr)
\Phi_{\vecc}
=0,
\quad
(\forall \vecc\in \Gamma)
\]
\[
 \Bigl(
 \lambda (c_i-\beta_i)+\sum_{j=1}^m a_{ji}\lambda x_j\del_{x_j}
 \Bigr)
 \Phi_{\vecc}=0,
\quad
(\forall\vecc\in \Gamma,\,\,\,i=1,\ldots,n)
\]
We have the $\nbigrtilde$-modules
$\nbigm^{\GKZ}(\nbiga,\Gamma,\vecbeta)$
corresponding to the system 
$\GKZ_{\nbigrtilde}(\nbiga,\Gamma,\vecbeta)$.
(See \S\ref{subsection;15.5.11.10},
 in particular Lemma \ref{lem;15.5.15.11}.)
According to Proposition \ref{prop;15.5.11.1}
and Remark \ref{rem;15.5.15.20},
we obtain the following proposition.

\begin{prop}
\label{prop;15.5.15.50}
$\nbigm^{\GKZ}(\nbiga,K(\nbiga),0)$
(resp. $\nbigm^{\GKZ}(\nbiga,K(\nbiga)^{\circ},0)$)
is naturally isomorphic to
$\nbigm_{\nbiga,0,\bullet}$
(resp. $\nbigm_{\nbiga,0,\circ}$).
\end{prop}

\paragraph{Special case}

For any $\vecp\in \seisuu^m$,
we set
\[
 \overline{\square}_{\vecp}:=
 \prod_{j\in \supp_+(\vecp)}
 (\lambda\del_j)^{p_j}
-\prod_{i\in\supp_-(\vecp)}
 (\lambda\del_j)^{-p_j}.
\]

For any element $\vecc_0\in K(\nbiga)$,
we consider the following system of differential equations
$\GKZ_{\nbigrtilde}^{\ord}(\nbiga,\vecbeta-\vecc_0)$
for $\Phi_{\vecc_0}$:
\[
 \overline{\square}_{\vecp}\Phi_{\vecc_0}=0
\quad
 (\forall\vecp\in L_{\nbiga})
\]
\[
 \Bigl(
 \lambda^2\del_{\lambda}
 +n\lambda+\sum_{j=1}^mx_j\lambda\del_{x_j}
 \Bigr)\Phi_{\vecc_0}=0
\]
\[
 \Bigl(
 \lambda(c_{0i}-\beta_i)
+\sum_{j=1}^ma_{ji}x_j\lambda\del_{x_j}\Bigr)
 \Phi_{\vecc_0}=0
\quad
 (i=1,\ldots,n)
\]
For any $\vecgamma\in\seisuu^m$,
let $\nbigi(\nbiga,\vecgamma)$ denote
the left ideal of $\nbigrtilde_{\cnum^m}$
generated by 
$\overline{\square}_{\vecp}$ $(\vecp\in L_{\nbiga})$,
$\lambda^2\del_{\lambda}+n\lambda
+\sum_{j=1}^m\lambda x_j\del_{x_j}$,
and $-\lambda\gamma_i+\sum_{j=1}^ma_{ji}\lambda x_j\del_{x_j}$
$(i=1,\ldots,n)$.
Then, we have a natural isomorphism
\[
 \nbigm^{\GKZ}(\nbiga,\Gamma,\vecbeta)
\simeq
 \nbigrtilde_{\cnum^m}/\nbigi(\nbiga,\vecbeta-\vecc_0)
\]
which is the $\nbigrtilde_{\cnum^m}$-module
corresponding to
$\GKZ^{\ord}(\nbiga,\vecbeta-\vecc_0)$.

Suppose that $K(\nbiga)=\seisuu_{\geq 0}\nbiga$.
Then, we have 
\begin{equation}
 \label{eq;15.5.15.40}
 \nbigm^{\GKZ}(\nbiga,K(\nbiga),0)
\simeq
 \nbigrtilde_{\cnum^m}\big/
 \nbigi(\nbiga,0).
\end{equation}
Suppose moreover that
$K(\nbiga)^{\circ}=K(\nbiga)+\vecc_0$
for an element $\vecc_0\in K(\nbiga)$.
Then, we have
\begin{equation}
 \label{eq;15.5.15.41}
 \nbigm^{\GKZ}(\nbiga,K(\nbiga)^{\circ},0)
\simeq
 \nbigrtilde_{\cnum^m}\big/
 \nbigi(\nbiga,-\vecc_0).
\end{equation}
For an expression $\vecc_0=\sum_{i=1}^m b_i\veca_i$,
we have the following commutative diagram:
\begin{equation}
 \label{eq;15.5.16.1}
 \begin{CD}
 \nbigm_{\nbiga,0,\circ}
 @>>>
 \nbigm_{\nbiga,0,\bullet}
 \\ 
 @V{\simeq}VV @V{\simeq}VV \\
 \nbigrtilde_{\cnum^m}\big/
  \nbigi(\nbiga,-\vecc_0)
 @>>>
  \nbigrtilde_{\cnum^m}\big/
  \nbigi(\nbiga,0)
 \end{CD}
\end{equation}
Here, the vertical arrows are given by
(\ref{eq;15.5.15.40}), (\ref{eq;15.5.15.41})
and Proposition \ref{prop;15.5.15.50},
and the upper horizontal arrow is
the morphism underlying
the natural morphism
of mixed twistor $\nbigd$-modules
$\gbigt_{\nbiga,0,\circ}
\lrarr
 \gbigt_{\nbiga,0,\bullet}$,
and the lower horizontal arrow is 
given by the multiplication of
$\prod_{i=1}^m(\lambda\del_{x_i})^{b_i}$.

\subsubsection{Relation with the reduced quantum
$\nbigd$-module of toric complete intersection}

Inspired by the work of E. Mann and T. Mignon
{\rm\cite{Mignon-Mann}},
T. Reichelt and C. Sevenheck {\rm\cite{Reichelt-Sevenheck2}}
introduced some systems of differential equations
to describe the reduced quantum $\nbigd$-modules
of complete intersections.

Let $X$ be an $n$-dimensional projective toric manifold
corresponding to a fan $\Sigma$.
Let $D_i$ $(i=1,\ldots,m)$ be the hypersurfaces of $X$
corresponding to the one dimensional cones in $\Sigma$.
Let $K_X$ denote the canonical bundle of $X$.
Let $\nbigl_i$ $(i=1,\ldots,r)$ be nef line bundles on $X$
such that 
$(K_X\otimes\bigotimes_{i=1}^{r}\nbigl_i)^{\lor}$ is nef.
We may assume that
$\nbigl_j=\nbigo\bigl(\sum_{i=1}^m \beta_{ji}D_i\bigr)$
for some $\beta_{ji}\geq 0$.
We can regard $\bigoplus \nbigl_i^{\lor}$ as a toric manifold.
Let $\nbiga\subset\seisuu^{n+r}$ be 
the set of primitive vectors in the one dimensional cones
of $\bigoplus \nbigl_i^{\lor}$.
We have the systems of differential equations
$\GKZ^{\ord}_{\nbigrtilde}(\nbiga,K(\nbiga),0)$
and 
$\GKZ^{\ord}_{\nbigrtilde}(\nbiga,K(\nbiga)^{\circ},-\vecc_0)$,
where 
$\vecc_0=
 (\overbrace{0,\ldots,0}^n,\overbrace{-1,\ldots,-1}^r)
 \in\seisuu^{n+r}$.
On $(\cnum^{\ast})^{m+r}$,
the systems
$\GKZ^{\ord}_{\nbigrtilde}(\nbiga,K(\nbiga),0)$
and 
$\GKZ^{\ord}_{\nbigrtilde}(\nbiga,K(\nbiga)^{\circ},-\vecc_0)$
are equivalent to 
the systems in Definition-Lemma 6.1 of \cite{Reichelt-Sevenheck2},
i.e.,
we have the following commutative diagram:
\[
 \begin{CD}
 \bigl(
 \nbigrtilde_{\cnum^{m+r}}
 \big/
 \nbigi(\nbiga,-\vecc_0)
 \bigr)_{|(\cnum^{\ast})^{m+r}}
 @>{\simeq}>>
  {}^{\ast}_0\widehat{\nbign}^{(0,0,0)}
 \\
 @VVV @V{b}VV \\
  \bigl(
 \nbigrtilde_{\cnum^{m+r}}
 \big/
 \nbigi(\nbiga,0)
 \bigr)_{|(\cnum^{\ast})^{m+r}}
 @>{\simeq}>>
 {}^{\ast}_0\widehat{\nbigm}^{(-r,0,0)}
 \end{CD}
\]
(See \cite[\S6]{Reichelt-Sevenheck2} for 
the $\nbigrtilde$-modules
${}^{\ast}_0\widehat{\nbign}^{(0,0,0)}$
and 
${}^{\ast}_0\widehat{\nbigm}^{(-r,0,0)}$.)
Hence, 
the commutative diagram (\ref{eq;15.5.16.1}) implies
that the $\nbigrtilde$-modules
${}^{\ast}_0\widehat{\nbign}^{(0,0,0)}$
and 
${}^{\ast}_0\widehat{\nbigm}^{(-r,0,0)}$
underlie mixed twistor $\nbigd$-modules.
We also obtain that the image of 
the morphism $b$ underlies a pure twistor $\nbigd$-module,
which is related with
the reduced quantum $\nbigd$-module
of the complete intersection.

\begin{rem}
In {\rm\cite{Reichelt-Sevenheck2}},
they constructed $\nbigrtilde$-modules
by using the partial Fourier-Laplace transform
of GKZ-hypergeometric systems
and the Brieskorn lattices associated to
the Hodge filtrations.
(See also {\rm\cite{Reichelt-Sevenheck1}}
and {\rm\cite{sabbah-Fourier-Laplace}}.)
They conjectured some relation of the $\nbigrtilde$-modules
with ${}^{\ast}_0\widehat{\nbign}^{(0,0,0)}$
and 
${}^{\ast}_0\widehat{\nbigm}^{(-r,0,0)}$
{\rm\cite[Conjecture 6.13]{Reichelt-Sevenheck2}}.

We review the construction in 
{\rm\S\ref{section;14.12.27.1}},
and we compare it with the $\nbigrtilde$-modules
underlying our mixed twistor $\nbigd$-modules
in Proposition {\rm\ref{prop;14.11.12.3}}.
(See also Proposition {\rm\ref{prop;14.12.15.1}}
and Corollary {\rm\ref{cor;14.12.3.102}}.)
In this way, we can verify their conjecture.
In {\rm\cite{Reichelt-Sevenheck3}},
they also verified their conjecture
in a different way.
\hfill\qed
\end{rem}

\subsection{Application to toric local mirror symmetry}

Recall that GKZ-hypergeometric systems
$\GKZ^{\ord}(\nbiga,\Gamma,0)$
and the variants
$\GKZ^{\ord}_{\nbigrtilde}(\nbiga,\Gamma,0)$
play important roles in the study of 
mirror symmetry of toric manifolds.
We shall apply the general theory of mixed twistor
$\nbigd$-modules to the study of toric local mirror symmetry.

One of the important goals in the study of mirror symmetry
is to obtain an isomorphism of Frobenius manifolds
associated to a mirror pair of
an $A$-model and a $B$-model.
Roughly, a Frobenius manifold
is a complex manifold 
with a holomorphic multiplication
and a holomorphic inner product
on the tangent bundle
satisfying some compatibility conditions.
The Frobenius manifold associated to the $A$-model
contains much information on 
the genus $0$ Gromov-Witten invariants.
The Frobenius manifold associated to the $B$-model
contains much information on 
the generalized Hodge structure 
of the Landau-Ginzburg model.
Hence, it is interesting to have such an isomorphism
of the Frobenius manifolds.
One of the most celebrated results is due to A. Givental
who established it in the case of toric weak Fano manifolds.
(See also \cite{Iritani,Reichelt-Sevenheck1}.)

Pursuing an analogue of such an isomorphism
in the study of local mirror symmetry,
Konishi and Minabe
\cite{Konishi-Minabe-cubic,
 Konishi-Minabe,Konishi-Minabe2,
 Konishi-Minabe3}
introduced the concept of 
mixed Frobenius manifold
as a generalization of Frobenius manifold.
Note that Frobenius manifolds do not appear
at least naively in the local case.
Roughly, a mixed Frobenius manifold 
is a complex manifold 
with a holomorphic multiplication,
a holomorphic filtration on the tangent bundle,
and inner products on the graded pieces
with respect to the filtration,
satisfying some compatibility conditions.
It looks natural to expect that
mixed Frobenius manifolds appear widely.

Konishi and Minabe particularly studied the case of 
any toric weak Fano surface $S$.
On the local $A$-side,
they constructed a mixed Frobenius manifold
from the genus $0$ local Gromov-Witten invariants of $S$.
On the local $B$-side,
they suggested that the expected mixed Frobenius manifolds
should be related with a variation of mixed Hodge structure
associated to the corresponding Landau-Ginzburg model.

Then, it is natural to ask whether there really exists
a mixed Frobenius manifold on the local B-side.
Recently, in his master thesis,
Y. Shamoto proved a reconstruction theorem of mixed Frobenius manifolds
as a generalization of a reconstruction theorem of
Frobenius manifolds due to C. Hertling and Y. Manin
\cite{Hertling-Manin}.
Together with the description of the variation of mixed Hodge structure
in \cite{Batyrev-VMHS, Konishi-Minabe, Stienstra},
he proved the existence of mixed Frobenius manifolds
associated to some toric local $B$-models,
up to the ambiguity of the choice of inner products
on the graded pieces.
It is still interesting to ask
how to choose pairings on the graded pieces,
and how to obtain an isomorphism of mixed Frobenius manifolds
associated to the local mirror pair.

\vspace{.1in}
In this paper,
we study the expectation in a rough level.
Instead of mixed Frobenius manifolds,
we shall study mixed TEP-structures
which will be explained in the next subsection.
We shall obtain an isomorphism of
mixed TEP-structures
associated to a mirror pair of
a local $A$-model and a local $B$-model.

\subsubsection{Mixed TEP-structure}

Recall that a TE-structure on a complex manifold $M$
in the sense of \cite{Hertling}
is a locally free 
$\nbigo_{\cnum_{\lambda}\times M}$-module $\nbigv$
with a meromorphic flat connection
$\nabla:\nbigv\lrarr\nbigv\otimes
 \Omega^1_{\cnum_{\lambda}\times M}
 \bigl(\log(\{0\}\times M)\bigr)
 \otimes
 \nbigo\bigl(\{0\}\times M\bigr)$.
If it is equipped with a perfect pairing
$P:\nbigv\otimes j^{\ast}\nbigv\lrarr
 \lambda^n\nbigo_{\cnum_{\lambda}\times M}$
such that 
$(j^{\ast}P)(j^{\ast}a\otimes b)
=(-1)^n P(b\otimes j^{\ast}a)$,
then $(\nbigv,\nabla,P)$ is called a TEP-structure 
or more precisely ${\rm TEP}(n)$-structure  \cite{Hertling} on $M$.
Here,
$j:\cnum_{\lambda}\times M
 \lrarr \cnum_{\lambda}\times M$
is given by
$j(\lambda,Q)=(-\lambda,Q)$.
We shall often omit to denote $\nabla$,
i.e., $(\nbigv,\nabla,P)$ is denoted by $(\nbigv,P)$.

A mixed TE-structure on a complex manifold $M$
consists of the following:
\begin{itemize}
\item
A TE-structure $\nbigv$ on $M$.
\item
An increasing filtration
$\Wtilde=\bigl(\Wtilde_m(\nbigv)\,|\,m\in\seisuu\bigr)$ 
on $\nbigv$
such that 
(i) $\Wtilde_m=0$ $(m<\!<0)$ and $\Wtilde_{m}=\nbigv$ $(m>\!>0)$,
(ii) $\Gr^{\Wtilde}_m(\nbigv)$ are
locally free $\nbigo_{\cnum_{\lambda}\times M}$-modules,
(iii) $\Wtilde_m$ are preserved by the connection of $\nbigv$.
\end{itemize}
If each $\Gr^{\Wtilde}_m(\nbigv)$ is equipped with 
a non-degenerate pairing
$P_m:\Gr^{\Wtilde}_m(\nbigv)
\otimes
 j^{\ast}\Gr^{\Wtilde}_m(\nbigv)
\lrarr
 \lambda^{-m}\nbigo_{\cnum_{\lambda}\times M}$
such that 
$(\Gr^{\Wtilde}_m\nbigv,P_m)$
is a ${\rm TEP}(-m)$-structure,
then $(\nbigv,W,\{P_m\})$ is called a mixed TEP-structure.

For example,
a graded polarized variation of mixed Hodge structure 
naturally induces a mixed TEP-structure
by the Rees construction.
In that case, we can reconstruct
the Hodge filtration, the weight filtration
and the polarization on the graded pieces
of the original graded polarized variation of mixed Hodge structure
from the mixed TEP-structure.
So, it still contains much interesting information.

\subsubsection{Mixed twistor $\nbigd$-modules
 with a graded polarization}

Let $(\nbigm,P_{\real},W)$ be a mixed twistor $\nbigd$-module
on a complex manifold $Y$.
If the $\nbigrtilde_Y$-module $\nbigm$ is 
a locally free $\nbigo_{\nbigy}$-module,
then $(\nbigm,W)$ is a mixed TE-structure by definition.

Suppose that 
$(\nbigm,P_{\real},W)$ is pure of weight $w$,
i.e., $\Gr^W_j=0$ unless $j=w$.
Then, a polarization of $(\nbigm,P_{\real},W)$
is equivalent to a morphism
$\nbigm\lrarr \lambda^{-w}j^{\ast}\DDD(\nbigm)$,
where $\DDD$ is the duality functor for
mixed twistor $\nbigd$-modules.
If $\nbigm$ is a locally free $\nbigo_{\nbigy}$-module,
the morphism gives a pairing
$\nbigp:\nbigm\otimes
 j^{\ast}\nbigm
\lrarr
 \lambda^{-w+\dim Y}\nbigo_{\nbigy}$,
and $(\nbigm,\nbigp)$ is 
a TEP$(-w+\dim Y)$-structure.

If $(\nbigm,P_{\real},W)$ is not necessarily pure,
a tuple of polarizations 
$\vecnbigp=(\nbigp_m\,|\,m\in\seisuu)$ on
$\Gr^W_m(\nbigm,P_{\real})$
is called a graded polarization of
$(\nbigm,P_{\real},W)$.
If $\nbigm$ is a locally free 
$\nbigo_{\nbigy}$-module,
each $\Gr^W_m(\nbigm)$ is a locally free
$\nbigo_{\nbigy}$-module.
We set 
$\Wtilde_{m-\dim Y}\nbigm:=
 W_{m}\nbigm$.
Then, 
$(\nbigm,\Wtilde,\vecnbigp)$
is a mixed TEP-structure.

\subsubsection{Mixed TEP-structures on 
$\GKZ(\nbiga,K(\nbiga),0)$
 and $\GKZ(\nbiga,K(\nbiga)^{\circ},0)$}

We use the notation in \S\ref{subsection;15.5.16.20}.
We have the open subset 
$(\cnum^m)^{\reg}\subset \cnum^m$
determined by the regularity condition of 
the Laurent polynomials $F_{\nbiga}$.
The restriction
$M_{\nbiga,0,\star|(\cnum^m)^{\reg}}$ $(\star=\ast,!)$
are locally free $\nbigo$-modules.
By the general theory of mixed twistor $\nbigd$-modules,
$\nbigm_{\nbiga,0,\star}$ are locally free $\nbigo$-modules
on $\cnum\times (\cnum^m)^{\reg}$.
As previously remarked,
we obtain the mixed TE-structure
$(\nbigm_{\nbiga,0,\star|(\cnum^m)^{\reg}},W)$.

If $0$ is an interior point of the convex hull of $\nbiga$,
then it turns out that 
we have a natural isomorphism
$(\nbigm_{\nbiga,0,!},P_{\real},W)
\simeq
 (\nbigm_{\nbiga,0,\ast},P_{\real},W)$,
that they are pure of weight $n+m$,
and that 
they are equipped with a natural polarization.
In particular,
the restriction to
$(\nbigm_{\nbiga,0,!},P_{\real},W)_{|(\cnum^m)^{\reg}}$
with the polarization
give a TEP$(-n)$-structure.
It is equivalent to the TEP-structure
previously studied in the mirror symmetry
\cite{Iritani,Reichelt-Sevenheck1}.

In the general case,
we need to choose an additional datum for the construction
of graded polarization on 
$(\nbigm_{\nbiga,0,!},P_{\real},W)$
and $(\nbigm_{\nbiga,0,\ast},P_{\real},W)$.
It turns out that
if we are given a point $\vecb\in\seisuu^n$
such that $0$ is an interior point of
the convex hull of $\nbiga\cup\{\vecb\}$,
then we obtain graded polarizations 
on
$(\nbigm_{\nbiga,0,!},P_{\real},W)$
and $(\nbigm_{\nbiga,0,\ast},P_{\real},W)$
depending on $\vecb$.
In particular,
we obtain 
mixed TEP-structures
$(\nbigm_{\nbiga,0,\star|(\cnum^m)^{\reg}},
 \Wtilde,\vecnbigp_{\vecb,\star})$
on $(\cnum^m)^{\reg}$.

We have the action of $T^n$ on $(\cnum^{\ast})^m$
given by 
$(s_1,\ldots,s_n)(z_1,\ldots,z_m)
=(s^{\veca_1}z_1,\ldots,s^{\veca_m}z_m)$.
Let $S_{\nbiga}:=(\cnum^{\ast})^m/T^n$.
Let $S_{\nbiga}^{\reg}$
denote the image of 
$(\cnum^{\ast})^m\cap(\cnum^m)^{\reg}$
by the projection
$(\cnum^{\ast})^m\lrarr S_{\nbiga}$.

It turns out that
$(\nbigm_{\nbiga,0,\star},\Wtilde,\vecnbigp_{\vecb,\star})$
on $(\cnum^{\ast})^m\cap(\cnum^m)^{\reg}$
is equivariant with respect to the action.
So, 
we obtain mixed TEP-structures
$(\nbigv_{\nbiga},\Wtilde,\vecnbigp_{\vecb,\star})$
$(\star=\ast,!)$
on $S_{\nbiga}^{\reg}$.

\subsubsection{An isomorphism}

Let $X$ be an $n$-dimensional smooth projective toric variety
corresponding to a fan $\Sigma$.
Let $D_i$ $(i=1,\ldots,m)$ be the hypersurfaces of $X$
corresponding to the one dimensional cones in $\Sigma$.
Let $K_X$ denote the canonical bundle of $X$.
Let $\nbigl_i$ $(i=1,\ldots,r)$ be nef line bundles on $X$
such that 
$(K_X\otimes\bigotimes_{i=1}^{r}\nbigl_i)^{\lor}$ is nef.
We may assume that
$\nbigl_j=\nbigo\bigl(\sum_{i=1}^m \beta_{ji}D_i\bigr)$
for some $\beta_{ji}\geq 0$.
Let $Y$ be the projective completion of 
$\nbige^{\lor}:=\bigoplus_{i=1}^r\nbigl^{\lor}_i$,
i.e.,
$Y=\proj\bigl(
 \bigoplus\nbigl_i\oplus\nbigo
 \bigr)$.

We have the degenerated quantum products 
$\bullet_{\sigma}$ $(\sigma\in\nbigu_X)$
on $H^{\ast}(X,\cnum)$
induced by the Gromov-Witten invariants of $Y$,
as explained in \S\ref{subsection;15.5.16.30}.
Here, $\nbigu_X$ is an appropriate open subset
in $H^{2}(X,\cnum)$.
As in the ordinary case,
we have the associated TE-structure 
$\QDM\bigl(X,\nbige^{\lor}\bigr)$
on $\nbigu_X$.
Let $\gamma\in H^2(Y)$ be the first Chern class of
the tautological line bundle of $Y$ over $X$.
As explained in \S\ref{subsection;15.5.16.41},
we introduce a filtration $\Wtilde$ on $H^{\ast}(X,\cnum)$,
and pairings on $\Gr^{\Wtilde}_jH^{\ast}(X,\cnum)$
by using the action of $\gamma$ on $H^{\ast}(Y)$
and inner product on $H^{\ast}(Y)$.
The construction was motivated by
both the general theory of mixed twistor $\nbigd$-modules
and the construction of Konishi and Minabe \cite{Konishi-Minabe2}.
Thus, we obtain a mixed TEP-structure
$(\QDM(X,\nbige^{\lor}),\Wtilde,\vecnbigp^{\nbige^{\lor}})$
on $\nbigu_X$.
It is equivariant with respect to 
the translation action of $2\pi\sqrt{-1}H^2(X,\seisuu)$ 
on $\nbigu_X$.
So, we obtain a mixed TEP-structure
$(\QDM(X,\nbige^{\lor}),\Wtilde,\vecnbigp^{\nbige^{\lor}})'$
on $\nbigu_X/2\pi\sqrt{-1}H^2(X,\seisuu)$.

\begin{rem}
As explained in 
{\rm\S\ref{subsection;15.5.16.50}},
the degenerated quantum products
are related with the local Gromov-Witten invariants
in some special cases.
It was essentially given in
{\rm \cite{Konishi-Minabe-cubic}},
and used in {\rm\cite{Konishi-Minabe2}}
for the construction of the weight filtration and 
the pairings.
\hfill\qed
\end{rem}

Let $\nbiga(\nbige^{\lor})\subset \seisuu^{n+r}$
be the set of the primitive vectors
in the one dimensional cones 
in the fan corresponding to $\nbige^{\lor}$.
Let $\nbiga(X)=\{[\rho_i]\,|\,i=1,\ldots,m\}\subset\seisuu^n$ 
be the set of the primitive vectors
in the one dimensional cones in the fan $\Sigma$
corresponding to $X$.
In this case,
$\nbiga(\nbige^{\lor})$
consists of elements
$\veca_i$ $(i=1,\ldots,m+r)$
given as follows:
\[
 \veca_i:=
 \left\{
 \begin{array}{ll}
 [\rho_i]+\sum_{j=1}^r\beta_{ji}n_j
 & (i=1,\ldots,m)\\
 n_{i-m}
 & (i=m+1,\ldots,m+r)
 \end{array}
 \right.
\]
Here, $n_i=(\overbrace{0,\ldots,0}^{m+i-1},1,0,\ldots,0)$.
We also set
$\veca_{m+r+1}:=-\sum_{i=1}^r n_i$.
Then,
we have the mixed TEP-structure
$(\nbigm_{\nbiga(\nbige^{\lor})},\Wtilde,\vecnbigp_{\veca_{m+r+1},\ast})$
on $(\cnum^{m+r})^{\reg}$.
We obtain a mixed TEP-structure
$(\nbigv_{\nbiga(\nbige^{\lor})},\Wtilde,\vecnbigp_{\veca_{m+r+1},\ast})$
on $S_{\nbiga(\nbige^{\lor})}^{\reg}$
as the reduction.

\begin{thm}
\label{thm;15.5.16.100}
We have the following.
\begin{itemize}
\item
An open subset $U_1\subset H^2(X,\cnum)/2\pi\sqrt{-1}H^2(X,\seisuu)$
which contains a neighbourhood of the large radius limit point.
\item
An open subset $U_2\subset S_{\nbiga(\nbige^{\lor})}^{\reg}$
which contains a neighbourhood of the large radius point.
\item
A holomorphic isomorphism $\varphi:U_1\simeq U_2$.
\item
Isomorphism of mixed TEP-structures 
\[
 \bigl(
 \QDM(X,\nbige^{\lor}),
 \Wtilde,
 \vecnbigp^{\nbige^{\lor}}
 \bigr)'_{|U_1}\,\,
\quad\mbox{\rm and}
\quad
 \varphi^{\ast}
 \bigl(
 \nbigv_{\nbiga(\nbige^{\lor}),\ast},
 \Wtilde,
 \vecnbigp_{\veca_{m+r+1},\ast}
 \bigr)_{|U_2}
\]
up to shift of weights.
\end{itemize}
\end{thm}

See Theorem \ref{thm;14.11.21.10}
for a more refined and precise claim.

\begin{rem}
In {\rm\cite{Reichelt-Sevenheck2}},
it is announced that
a result in {\rm \cite{Iritani-Mann-Mignon}}
implies the comparison of the TE-structures
in Theorem {\rm\ref{thm;15.5.16.100}}.
In {\rm\cite{Iritani-shift-operators}},
the mirror theorem of Givental
is generalized
for non-compact or non-semipositive toric manifolds.
At this moment, it is not clear to the author 
if we could also deduce 
the comparison of the weight filtrations
from 
{\rm \cite{Iritani-shift-operators,Iritani-Mann-Mignon}}.
\hfill\qed
\end{rem}

\subsubsection{Idea of the proof}

Roughly, our isomorphism is obtained as 
the specialization of the isomorphism of Givental
for the weak Fano toric manifold $Y$.
We naturally have 
$H^2(Y,\cnum)=H^2(X,\cnum)\times\cnum\gamma$.
We have the quantum $\nbigd$-module 
$(\QDM(Y),\nbigp_Y)$
on an appropriate open subset of 
$H^2(Y,\cnum)\big/2\pi\sqrt{-1}H^2(Y,\seisuu)$
associated to the Gromov-Witten invariants of $Y$.
The mixed TEP-structure
$\bigl(\QDM(X,\nbige^{\lor}),
 \Wtilde,\vecnbigp^{\nbige^{\lor}}\bigr)$
is described as the ``specialization'' of
$(\QDM(Y),\nbigp_Y)$.
(We explain the procedure ``specialization''
in \S\ref{subsection;14.12.9.2}.)
Let $\nbiga(Y):=
 \nbiga(\nbige^{\lor})\cup\{\veca_{m+r+1}\}$ 
in $\seisuu^{n+r}$,
which is the set of the primitive vectors
in the one dimensional cones in a fan corresponding to $Y$.
We naturally have
$S_{\nbiga(Y)}=S_{\nbiga(\nbige^{\lor})}\times\cnum^{\ast}$.
By using the results on the specialization of
mixed twistor $\nbigd$-modules,
we can describe the mixed TEP-structure
$\bigl(\nbigv_{\nbiga(\nbige^{\lor})\ast},
 \Wtilde,\vecnbigp_{\veca_{m+r+1}\ast}
 \bigr)$
as ``the specialization'' of 
the TEP-structure 
$(\nbigv_{\nbiga(Y)},\nbigp)$.
We have the isomorphism of Givental
between the TEP-structures $(\QDM(Y),\nbigp_Y)$
and $(\nbigv_{\nbiga(Y)},\nbigp)$.
(See \cite{Iritani, Reichelt-Sevenheck1}.)
Hence, we can obtain the isomorphism of
mixed TEP-structures
in Theorem \ref{thm;15.5.16.100}
as the specialization of the isomorphism of Givental.

\subsection{The mixed twistor $\nbigd$-modules
associated to meromorphic functions}

Motivated by the applications to the study of 
the twistor property of GKZ-hypergeometric systems,
mentioned in the previous subsections,
we shall consider the technical issues
on the mixed twistor $\nbigd$-modules.
Let $X$ be a complex manifold 
with a hypersurface $H$.
Let $f$ be a meromorphic function on $X$
whose poles are contained in $H$.
We have the naturally defined integrable mixed twistor 
$\nbigd$-modules with real structure
 $\nbigt_{\star}(f,H)$ $(\star=\ast,!)$ 
over $L_{\star}(f,H)$.
We also have the naturally induced polarizations on 
$\Gr^W_w\nbigt_{\star}(f,H)$
which depend on the choice of an effective divisor $D$
whose support is $H$.
The much part of this paper
is devoted to the study of such kind of
graded polarized mixed twistor $\nbigd$-modules
with real structure.
For instance, we consider the following issues.
\begin{itemize}
\item
Let $X^{(1)}:=\proj^1_{\tau}\times X$.
We set $H^{(1)}:=(\proj^1\times H)\cup(\{0,\infty\}\times X)$.
Let $\iota:X\lrarr X^{(1)}$ be the morphism
induced by $\{0\}\lrarr \proj^1_{\tau}$.
Let $f,g\in \nbigo_X(\ast H)$.
We have the graded polarized mixed twistor $\nbigd$-modules
with real structure $\nbigt_{\star}(g+\tau f,H^{(1)})$ on $X^{(1)}$
and $\nbigt_{\star}(g,H)$ on $X$,
where $\star=\ast,!$.
Under some assumptions, we shall relate
$\iota_{\dagger}\nbigt_{\star}(g,H)$
with 
the kernel and the cokernel of
$\nbigt_{!}(g+\tau f,H^{(1)})
\lrarr
 \nbigt_{\ast}(g+\tau f,H^{(1)})$.
It is not difficult to obtain a relation
in the level of mixed twistor $\nbigd$-modules
with real structures
(Proposition \ref{prop;14.10.18.11}).
We need more preliminaries
to obtain the relation between 
polarizations on the graded pieces
(Proposition \ref{prop;14.9.7.42}).
\item
Let $X^{(1)}$ and $H^{(1)}$ be as above.
Let $\pi:X^{(1)}\lrarr X$ be the projection.
We have the graded polarized mixed twistor $\nbigd$-modules
$\pi^0_{\dagger}\nbigt_{\star}(g+\tau f,H^{(1)})$
by taking the push-forward via $\pi$.
Let $Z_f$ be the zero-set of $f$.
We have the graded polarized mixed twistor $\nbigd$-modules
$\nbigt_{\ast}(g,H)[!(f)_0]$ and 
$\nbigt_{!}(g,H)[\ast(f)_0]$ on $X$.
Under some assumptions, we shall relate
$\pi^0_{\dagger}\nbigt_{\ast}(g+\tau f,H^{(1)})$
(resp.
 $\pi^0_{\dagger}\nbigt_{!}(g+\tau f,H)$)
and 
$\nbigt_{\ast}(g,D)[!Z_f]$
(resp.
$\nbigt_{!}(g,D)[\ast Z_f]$.)
Again, it is not difficult to obtain a relation
in the level of mixed twistor $\nbigd$-modules
with real structure
(Proposition \ref{prop;14.11.7.21}).
We need more preliminaries
to obtain the relation between 
polarizations on the graded pieces
(Proposition \ref{prop;14.9.7.210}).
\end{itemize}

We also need the compatibility 
of various standard functors for $\nbigd$-modules
and flat bundles,
which are explained in the appendix sections
\S\ref{section;14.12.26.10} and
\S\ref{section;14.12.26.11}.

\paragraph{Acknowledgement}

This study grows out of my attempt to 
understand the works of 
H. Iritani \cite{Iritani-convergence, Iritani},
Y. Konishi, S. Minabe
\cite{Konishi-Minabe, Konishi-Minabe2, Konishi-Minabe3}
and 
T. Reichelt, C. Sevenheck
\cite{Reichelt-Sevenheck1,Reichelt-Sevenheck2}.
I am grateful to them for helpful discussions
on various occasions.
In particular,
Reichelt and Sevenheck
explained their notation,
and remarked me the issue of the normality condition
in the study of GKZ-hypergeometric systems
and toric varieties.
I thank C. Sabbah for his kindness and for discussions
on many occasions.
I thank Y. Shamoto for some discussion.
I thank G. Wilkin for his interest to this study
and his kindness.
I am grateful to 
C. Hertling and U. Walther
for some discussions.
I thank A. Ishii,
A. Moriwaki,
M.-H. Saito,
Y. Tsuchimoto,
T. Xue and K. Vilonen for their kindness
and supports.

This work is partially supported by 
the Grant-in-Aid for Scientific
Research (C) (No. 22540078), 
the Grant-in-Aid for Scientific
Research (C) (No. 15K04843),
the Grant-in-Aid for Scientific Research
(A) (No. 22244003) and 
the Grant-in-Aid for Scientific Research (S)
(No. 24224001),
Japan Society for the Promotion of Science.
This research was partially completed 
while the author was visiting 
the Institute for Mathematical Sciences, 
National University of Singapore in 2014.
This research was partially completed 
while the author stayed in 
the Mathematisches Forschungsinstitut Oberwolfach,
partially supported by the Simons Foundation
and by the Mathematisches Forschungsinstitut Oberwolfach.

\section{$\nbigd$-modules associated to meromorphic functions}
\label{section;15.11.22.10}
\subsection{Purity condition}

\subsubsection{$\nbigd$-modules associated 
 to meromorphic functions}
\label{subsection;14.10.14.10}

Let $X$ be a complex manifold with a hypersurface $D$.
Let $\nbigo_X(\ast D)$ be the sheaf of meromorphic functions on $X$
whose poles are contained in $D$.
Let $\DDD$ denote the duality functor
on the category of holonomic $\nbigd$-modules on $X$.
For any coherent $\nbigd_X$-module $M$,
we set 
$M(\ast D):=\nbigo_X(\ast D)\otimes_{\nbigo_X}M$
and 
$M(!D):=\DDD\bigl(\DDD(M)(\ast D)\bigr)$.
The restriction of $M(\star D)$ $(\star =\ast,!)$
to $X\setminus D$ is naturally isomorphic to $M_{|X\setminus D}$.
We have the canonical morphism
$M(!D)\lrarr M(\ast D)$
whose restriction to $X\setminus D$ is the identity.

Let $f$ be a meromorphic function on $(X,D)$,
i.e.,
a section of $\nbigo_X(\ast D)$.
Let $(f)_0$ and $(f)_{\infty}$
denote the effective divisors
obtained as the zeroes and the poles of $f$.
The supports of the divisors are denoted by
$|(f)_0|$ and $|(f)_{\infty}|$.

We obtain the meromorphic flat bundle
$L_{\ast}(f,D):=\nbigo_X(\ast D)v$
with $\nabla v=v\,df$.
We naturally regard it as a $\nbigd_X$-module.
We set 
$L_!(f,D):=L_{\ast}(f,D)(!D)
=\DDD\bigl(L_{\ast}(-f,D)\bigr)$.
The image of the canonical morphism
$L_!(f,D)\lrarr L_{\ast}(f,D)$
is independent of the choice of $D$
such that $|(f)_{\infty}|\subset D$,
and denoted by $L(f)$.
We have natural isomorphisms
$L_{\star}(f,D)\simeq L(f)(\star D)$ for $\star=!,\ast$.
Indeed, 
$L(f)\lrarr L_{\ast}(f,D)$ 
naturally induces
$L(f)(\ast D)\lrarr L_{\ast}(f,D)$
which is clearly an isomorphism.
By using the duality,
we obtain
$L_!(f,D)\simeq L(f)(!D)$.

When $D=|(f)_{\infty}|$,
we set 
$L_{\star}(f):=L_{\star}(f,|(f)_{\infty}|)$.

\subsubsection{Purity}

We continue to use the notation in
\S\ref{subsection;14.10.14.10}.
We introduce a condition.

\begin{df}
\label{df;15.1.3.3}
We say that 
a meromorphic function $f$ on $X$
is pure at $P\in X$
if the canonical morphism
$L_{!}(f)\lrarr L_{\ast}(f)$ 
is an isomorphism
on a neighbourhood of $P$.
We say that $f$ is pure
if it is pure at any point of $X$.
\hfill\qed
\end{df}

Because 
$\DDD \bigl(L_{\ast}(f)\bigr)
\simeq
 L_!(-f)$
and 
$\DDD \bigl(L_{!}(f)\bigr)
\simeq
 L_{\ast}(-f)$,
we have the following easy lemmas.
\begin{lem}
\label{lem;14.11.23.10}
Let $f$ be a meromorphic function on $X$.
\begin{itemize}
\item
$f$ is pure if and only if $-f$ is pure.
\item
$f$ is pure if and only if
the morphisms
$L_!(f)\lrarr L_{\ast}(f)$
and $L_{!}(-f)\lrarr L_{\ast}(-f)$
are epimorphisms.
\hfill\qed
\end{itemize}
\end{lem}

Let $\varphi:Y\lrarr X$ be 
a proper morphism of complex manifolds
such that 
$\varphi$ induces an isomorphism
$Y\setminus \varphi^{-1}(D)\simeq
 X\setminus D$.
\begin{lem}
\label{lem;14.11.26.30}
Let $f$ be a meromorphic function on $(X,D)$.
Suppose that 
$\bigl|(\varphi^{\ast}f)_{\infty}\bigr|=\varphi^{-1}(D)$
and that $\varphi^{\ast}(f)$ is pure.
Then, $f$ is also pure,
and it satisfies 
$|(f)_{\infty}|=D$.
\end{lem}
\pf
We set $f_Y:=\varphi^{\ast}(f)$
and $D_Y:=\varphi^{-1}(D)$.
Because
$D=\varphi\bigl(|(f_Y)_{\infty}|\bigr)
\subset  |(f)_{\infty}|\subset D$,
we have $|(f)_{\infty}|=D$.
We have
$\varphi^i_{+}L_{\ast}(f_Y,D_Y)=0$ $(i\neq 0)$,
and 
$\varphi^0_{+}\bigl(
 L_{\ast}(f_Y,D_Y)\bigr)(\ast D)
\simeq 
\varphi^0_{+}\bigl(
 L_{\ast}(f_Y,D_Y)\bigr)$.
Hence, we have
$\varphi^0_{+}\bigl(
 L_{\ast}(f_Y,D_Y)\bigr)
\simeq
 L_{\ast}(f,D)$.
By the duality,
we have
$\varphi^0_{+}\bigl(
 L_!(f_Y,D_Y)\bigr)
\simeq
 L_!(f,D)$.
Then, from the purity of $f_Y$,
we obtain that 
$L_!(f,D)\lrarr L_{\ast}(f,D)$
is an isomorphism.
\hfill\qed

\subsubsection{Vanishing of cohomology 
for pure functions}

The purity condition sometimes implies
the vanishing of the cohomology.
We mention a typical case.
Let $X$ be a complex manifold
with a hypersurface $D$.
Let $f$ be a meromorphic function on $(X,D)$
such that
(i) $|(f)_{\infty}|=D$,
(ii) $f$ is pure.
Let $F:X\lrarr S$ be a proper morphism
of complex manifolds.

\begin{prop}
\label{prop;14.11.17.20}
Suppose that
$R^iF_{\ast}\bigl(
 \Omega_X^j\otimes\nbigo(\ast D)
 \bigr)=0$ 
for $(i,j)\in\seisuu_{\geq 0}^2$
such that 
$i+j>\dim X$.
Then, we have $F_{\dagger}^kL(f)=0$ $(k\neq 0)$.
\end{prop}
\pf
Because $\DDD L(f)\simeq L(-f)$
and because $-f$ is also pure,
it is enough to prove that
$F_{\dagger}^kL(f)=0$ for $k>0$.
We have only to prove the claim locally around any point of $S$.
We set 
$\omega_X:=\Omega_X^{\dim X}$
and 
$\omega_S:=\Omega_S^{\dim S}$.
Recall that $F_{\dagger}L(f)$ is 
obtained as
$RF_{\ast}\Bigl(
 \bigl[
 \omega_X
 \otimes_{F^{-1}\nbigo_S}
 F^{-1}\bigl(
 \nbigd_S\otimes \omega_S^{-1}
\bigr)
 \bigr]
 \otimes^L_{\nbigd_X}
 L(f)
 \Bigr)$.
We have the standard free resolution 
$\nbigd_X\otimes\Omega^{\bullet}_X[\dim X]$
of the right $\nbigd_X$-modules $\omega_X$.
By the assumption,
if $i+j>\dim X$, we have
\[
 R^jF_{\ast}\Bigl(
 \bigl[
 \Omega^i_X
 \otimes_{F^{-1}\nbigo_S}(F^{-1}\nbigd_S\otimes\omega_S)
 \bigr]
 \otimes_{\nbigo_X}
 L(f)
 \Bigr)
\simeq
 R^jF_{\ast}\bigl(
 \Omega^i_X\otimes L(f)
 \bigr)
\otimes_{\nbigo_S}
 (\nbigd_S\otimes\omega_S)
=0.
\]
Hence, we have 
the desired vanishing
$F_{\dagger}^kL(f)=0$ for $k>0$.
\hfill\qed

\begin{example}
The condition of Proposition {\rm\ref{prop;14.11.17.20}}
is satisfied
if $F$ is factorized into the composite of
morphisms of complex manifolds
$X\stackrel{\rho}{\lrarr} X'\stackrel{F'}{\lrarr}S$
such that 
(i) $D'=\rho(D)$ is also a hypersurface of $X'$,
(ii) $\rho$ induces $X\setminus D\simeq X'\setminus D'$,
(iii) $\nbigo_{X'}(D')$ is relatively ample with respect to $F'$.
\hfill\qed
\end{example}

\subsection{$\nbigd$-modules associated to 
non-degenerate meromorphic functions}

\subsubsection{A non-degeneracy condition}
\label{subsection;14.11.4.1}

Let $X$ be a complex manifold
with a simple normal crossing hypersurface $D$.
Let $f$ be a meromorphic function on $(X,D)$.
In this paper,
we shall often use the following non-degeneracy condition.

\begin{df}
\label{df;15.1.3.1}
$f$ is called non-degenerate along $D$
if the following holds
for a small  neighbourhood $N$ of $|(f)_{\infty}|$.
\begin{itemize}
\item
 $(f)_0\cap N$ is reduced and non-singular.
\item
$N\cap\bigl(|(f)_0|\cup D\bigr)$ is normal crossing.
\hfill\qed
\end{itemize}
\end{df}

Let $D=\bigcup_{i\in\Lambda}D_i$
and 
$|(f)_{\infty}|=\bigcup_{i\in\Lambda_f}D_i$
be the irreducible decompositions,
where $\Lambda_f\subset\Lambda$.
For any non-empty subset $I\subset \Lambda$,
we set 
$D_I:=\bigcap_{i\in\Lambda}D_I$
and $D_I^{\circ}:=D_I\setminus \bigcup_{j\not\in I}D_j$.
If $I=\emptyset$, we set $D_{\emptyset}=X$.
Then, the second condition can be reworded 
that 
$D_I^{\circ}$ is transversal with $|(f)_0|$
for any $I\subset \Lambda$
with $I\cap \Lambda_f\neq\emptyset$.

\begin{rem}
Suppose that a meromorphic function $f$
on $(X,D)$ is non-degenerate along $D$.
Let $D'$ be a hypersurface of $X$
such that $|(f)_{\infty}|\subset D'\subset D$.
Then, $f$ is non-degenerate along $D'$.
But, the converse does not hold in general.
Namely,
even if a meromorphic function $f$ on $(X,D')$
is non-degenerate along $D'$,
it is not necessarily non-degenerate along $D$.
For example,
set $X=\cnum^2$, $f=(z_1-z_2)/z_2$,
$D'=\{z_2=0\}$
and $D=\{z_1=0\}\cup\{z_2=0\}$.
\hfill\qed
\end{rem}

We reword the condition in terms of local coordinate systems.
We set $D^c_f:=\bigcup_{i\in\Lambda_f^c}D_i$,
where 
$\Lambda_f^c:=
 \Lambda\setminus \Lambda_f$.
We have $D=|(f)_{\infty}|\cup D_f^c$.
Let $Q\in |(f)_{\infty}|$.
We take a holomorphic coordinate neighbourhood
$(X_Q;z_1,\ldots,z_n)$ around $Q$
such that
\[
|(f)_{\infty}|\cap X_Q=
 \bigcup_{i=1}^{\ell_1}\{z_i=0\},
\quad
 D_f^c\cap X_Q=
 \bigcup_{i=\ell_1+1}^{\ell_1+\ell_2}\{z_i=0\}.
\]
Let $k_i$ denote the pole order of $f$
along $\{z_i=0\}$.
Then, we have an expression
$f=f_0\prod_{i=1}^{\ell_1}z_i^{-k_i}$
where $f_0$ is holomorphic.
Let $I\subset\Lambda_f$ be determined by
$Q\in D_I^{\circ}$.
Note $I\cap \Lambda_f\neq\emptyset$.
If $f$ is non-degenerate along $D$,
then the divisor $(f_{0|D_I^{\circ}})$ is reduced and non-singular.
Conversely, 
if the above holds for any $Q\in |(f)_{\infty}|$,
then $f$ is non-degenerate.

\begin{rem}
See Definition {\rm\ref{df;14.11.23.1}}
and Lemma {\rm\ref{lem;15.1.3.2}}
for the relation between the non-degeneracy condition in 
{\rm\cite{Kouchnirenko}}
and the condition in Definition {\rm\ref{df;15.1.3.1}}.
We postpone to discuss the exact relation 
between the (cohomologically) tameness condition for
algebraic functions 
(see {\rm\cite{Sabbah-tame-polynomial}})
and the conditions in 
Definitions {\rm\ref{df;15.1.3.3}} and {\rm\ref{df;15.1.3.1}}.
\hfill\qed
\end{rem}

\subsubsection{Convenient coordinate systems}
\label{subsection;14.6.27.2}

Suppose that a section $f$ of $\nbigo_X(\ast D)$ 
is non-degenerate along $D$.
We have a holomorphic coordinate system
$(X_Q;z_1,\ldots,z_n)$
around $Q\in|(f)_{\infty}|$ satisfying the following conditions.
\begin{itemize}
\item
$|(f)_{\infty}|\cap X_Q
=\bigcup_{i=1}^{\ell_1}\{z_i=0\}$
 and 
$D_f^c\cap X_Q
=\bigcup_{i=\ell_1+1}^{\ell_1+\ell_2}\{z_i=0\}$.
\item
 If $Q\in |(f)_0|$,
 we have 
 $f_{|X_Q}=z_n\prod_{i=1}^{\ell_1}z_i^{-k_i}$
 for some $\veck\in\seisuu_{> 0}^{\ell_1}$.
\item
 If $Q\not\in |(f)_0|$,
 we have
 $f_{|X_Q}=\prod_{i=1}^{\ell_1}z_i^{-k_i}$
 for some $\veck\in\seisuu_{> 0}^{\ell_1}$.
\end{itemize}
In this paper, such a coordinate system is called 
a convenient coordinate system.

\subsubsection{Families of non-degenerate functions}

Let $X\lrarr S$ be a smooth morphism
of complex manifolds.
Let $D$ be a simple normal crossing hypersurface
in $X$
with the irreducible decomposition
$D=\bigcup_{i\in\Lambda}D_i$.
Suppose that 
the induced morphisms 
$D^{\circ}_I\lrarr S$ are smooth
for any $I\subset\Lambda$.
For any $s\in S$,
let $(X_s,D_s)$
denote the fibers of $(X,D)$
over $s\in S$.
In such a situation,
we sometimes consider
a condition which is stronger than
the non-degeneracy along $D$.

\begin{df}
\label{df;14.12.9.10}
Let $f$ be a meromorphic function on $(X,D)$
which is non-degenerate along $D$.
If moreover $f_{|X_s}$ is non-degenerate along $D_s$
for any $s\in S$,
we say that $f$ is non-degenerate along $D$
over $S$.
\hfill\qed
\end{df}

\subsubsection{Purity in the non-degenerate case}

Let $X$ be a complex manifold 
with a normal crossing hypersurface $D$.
Let $f$ be a meromorphic function on $(X,D)$.

\begin{lem}
\label{lem;14.8.13.1}
Suppose that
$D=|(f)_{\infty}|$
and that $f$ is non-degenerate along $D$.
Then, $f$ is pure.
\end{lem}
\pf
We have only to check the claim
locally around any point of $|(f)_{\infty}|$.
We use a convenient coordinate system
as in \S\ref{subsection;14.6.27.2}.
We have a natural isomorphism
$L_!(f)(\ast D)
\simeq
 L_{\ast}(f)(\ast D)$.
Hence, for a large $N$,
$u:=\prod_{i=1}^{\ell_1}z_i^{N}v$
is a section of $L_!(f)$.

We have a natural morphism
$\varphi:\nbigd_X
\lrarr 
 L_!(f)$
given by
$P\longmapsto P u$.
Let us check that
the composite
of $\varphi$
and $L_!(f)\lrarr L_{\ast}(f)$
is an epimorphism.
It is enough to observe that
the image of
$u$ in $L_{\ast}(f)$ generates $L_{\ast}(f)$.
For $\vecm\in\seisuu^{m_1}$,
let $z^{\vecm}=\prod_{i=1}^{\ell_1}z_i^{m_i}$.
If $Q\in |(f)_{\infty}|$,
note 
$\del_{z_n}(z^{\vecm}v)
=z^{\vecm-\veck}v$.
If $Q\not\in|(f)_{\infty}|$,
note 
$z_1\del_{z_1}(z^{\vecm}v)
=-k_iz^{\vecm-\veck}v
+m_iz^{\vecm}v$.
Then, the claim easily follows.

Hence, the natural morphism
$\kappa:L_!(f)\lrarr L_{\ast}(f)$ is an epimorphism.
As remarked in Lemma \ref{lem;14.11.23.10},
because the dual of
$\kappa$
is the epimorphism
$L_{!}(-f)\lrarr L_{\ast}(-f)$,
we obtain that 
$\kappa$
is a monomorphism.
\hfill\qed

\subsubsection{Expressions 
of the $\nbigd$-modules
 associated to non-degenerate functions}

Let $V_D\nbigd_X$ denote the sheaf of 
subalgebras of $\nbigd_X$
generated by $\Theta_X(\log D)$
over $\nbigo_X$.
Let $f$ be a meromorphic function on $(X,D)$
which is non-degenerate along $D$.
Then, $L(f)$ and 
$L(f)(D)=L(f)\otimes_{\nbigo_X}\nbigo_X(D)$ are
naturally $V_D\nbigd_X$-modules.
\begin{lem}
\label{lem;14.11.23.101}
We naturally have
$L_!(f,D)\simeq
 \nbigd_X\otimes_{V_D\nbigd_X}
 L(f)$
and
$L_{\ast}(f,D)\simeq
 \nbigd_X\otimes_{V_D\nbigd_X}
 L(f)(D)$.
\end{lem}
\pf
Because 
$\bigl(
 \nbigd_X\otimes_{V_D\nbigd_X}
 L(f)\bigr)(\ast D)
\simeq
\bigl(
 \nbigd_X\otimes_{V_D\nbigd_X}
 L(f)(D)\bigr)(\ast D)
\simeq
 L(f)(\ast D)$,
we have natural morphisms:
\begin{equation}
\label{eq;14.6.30.1}
  \nbigd_X\otimes_{V_D\nbigd_X}L(f)(D)
\stackrel{\alpha_1}{\lrarr}
 L_{\ast}(f,D)
\end{equation}
\begin{equation}
\label{eq;14.6.30.2}
 L_!(f,D)
 \stackrel{\alpha_2}{\lrarr}
  \nbigd_X\otimes_{V_D\nbigd_X}L(f).
\end{equation}
It is enough to prove that the morphisms
are isomorphisms
locally around any point of $D$.
We use a convenient coordinate system
in \S\ref{subsection;14.6.27.2}.

Let $v$ be a frame of $L(f)$
over $\nbigo_X(\ast(f)_{\infty})$
such that $\nabla v=v\,df$.
Because $a_1:L_!(f,D)\lrarr L(f)$ is an epimorphism,
we can locally take a section $v'$ of $L_!(f,D)$
which is mapped to $v$ via $a_1$.
We consider the submodule 
$V_D\nbigd_X\,v'\subset L_!(f,D)$.
It is coherent over $V_D\nbigd_X$.
Note that 
$L(f)$ is coherent over $V_D\nbigd_X$.
Because
$\Ker a_1\cap V_D\nbigd_Xv'
=\Ker\bigl(
 a_{1|V_D\nbigd_Xv'}
 \bigr)$,
it is coherent over $V_D\nbigd_X$.
In particular, it is locally finitely generated 
over $V_D\nbigd_X$.
Take a generator $f_1,\ldots,f_m$
of $\Ker a_1\cap V_D\nbigd_Xv'$.
The supports of $f_j$ are contained in $D$.
We can take a large $N$ such that
$\prod_{i=1}^{\ell_1+\ell_2} z_i^{N}\,f_j=0$
in $L_!(f,D)$,
because $L_!(f,D)(\ast D)=L(f)(\ast D)$.
Then, it is easy to see that
$\prod_{i=1}^{\ell_1+\ell_2}z_i^N\,
 \bigl(
 V_D\nbigd_X\,v'\cap\Ker(a_1)
 \bigr)=0$.
We set 
$v'':=\prod_{i=1}^{\ell_1+\ell_2}z_i^Nv'$,
and then we have
$V_D\nbigd_Xv''\cap\Ker(a_1)=0$.
The morphism
$V_D\nbigd_Xv''\lrarr L(f)$
induced by $a_1$
gives an isomorphism
$V_D\nbigd_Xv''
\lrarr
 \nbigo_X(\ast (f)_{\infty})\prod_{i=1}^{\ell_1+\ell_2}z_i^Nv$
of $V_D\nbigd_X$-modules.
In particular,
$V_D\nbigd_Xv''$
is naturally 
an $\nbigo_X\bigl(\ast(f)_{\infty}\bigr)$-module.

We set 
$v^{(3)}:=\prod_{i=1}^{\ell_1}z_i^{-N}\,v''$.
We have
$z_i\del_iv^{(3)}
=-k_iv^{(3)}\prod_{i=1}^{\ell_1}z_i^{-k_i}$
for $i=1,\ldots,\ell_1$,
and
$z_i\del_iv^{(3)}
=Nv^{(3)}$
for $i=\ell_1+1,\ldots,\ell_1+\ell_2$.
Take $0\leq p_{\ell_1+1},\ldots,p_{\ell_1+\ell_2}\leq N$.
We set $\vecp:=(p_{\ell_1+1},\ldots,p_{\ell_1+\ell_2})$.
We set 
$v_{\vecp}:=
 \prod_{i=\ell_1+1}^{\ell_1+\ell_2}
 \del_i^{p_i}v^{(3)}$.
Note that
$a_1(v_{\vecp})
=C_{\vecp}\prod_{i=\ell_1+1}^{\ell_1+\ell_2}
 z_i^{N-p_i}v$
for a non-zero constant $C_{\vecp}$.
We consider the following morphism
induced by $a_1$:
\begin{equation}
\label{eq;14.7.16.1}
V_D\nbigd_X\cdot
 \prod_{i=\ell_1+1}^{\ell_1+\ell_2}
 \del_i^{p_i}v^{(3)}
\lrarr
 L(f)
\end{equation}
Let us observe that 
(\ref{eq;14.7.16.1}) is a monomorphism.
If $\vecp=(0,\ldots,0)$,
it has already been observed.
If $p_i>0$,
set $p_j':=p_j$ $(j\neq i)$
and $p_i':=p_i-1$.
We have 
$\del_iv_{\vecp'}=v_{\vecp}$.
Let $s\in\Ker(a_1)\cap V_D\nbigd_Xv_{\vecp}$.
We have
$s=\del_is'$ for some 
$s'\in V_D\nbigd_Xv_{\vecp'}$.
We have
$0=a_1(z_is)=a_1(z_i\del_is')=0$.
We have $z_i\del_is'=0$
in 
$V_D\nbigd_X v_{\vecp'}
\simeq
 \nbigo_X(\ast (f)_{\infty})
 \prod_{j=1}^{N-p_j'}v$.
But, because $N-p_i'>0$,
it is easy to see that
if a section $s'$
of $\nbigo_X(\ast (f)_{\infty})
 \prod_{j=1}^{N-p_j'}v$
satisfies $z_i\del_is'=0$
then $s'=0$.
Hence, we obtain $s=0$.
In particular,
the induced morphism
$V_D\nbigd_X\,v_{(N,\ldots,N)}
\lrarr L(f)$
is a monomorphism.
Because
$a_1(v_{(N,\ldots,N)})$ is
$v$ multiplied by a non-zero constant,
$V_D\nbigd_X\,v_{(N,\ldots,N)}
\lrarr L(f)$ is also an epimorphism,
i.e., an isomorphism.

Hence, we have
a $V_D\nbigd_X$-homomorphism
$a_2:L(f)\lrarr L_!(f,D)$
such that $a_1\circ a_2=\id$.
It induces a $\nbigd_X$-homomorphism
$\beta_2:\nbigd_X\otimes_{V_D\nbigd_X} L(f)
\lrarr
 L_!(f,D)$.
We set $g:=\beta_2\circ\alpha_2$
which is an endomorphism of $L_!(f,D)$.
For the dual,
$g$ induces the identity  on $\DDD L_!(f,D)$.
Hence, we obtain that $g$ is the identity.
In particular,
(\ref{eq;14.6.30.2}) is an epimorphism
and $\beta_2$ is a monomorphism.
Because $\nbigd_X\otimes_{V_D\nbigd_X}L(f)$
is generated by $v$,
 $\beta_2$ is an epimorphism.
Hence, $\alpha_2$ and $\beta_2$
are isomorphisms.

Let us study the morphism $\alpha_1$.
Because $L(f)=L(f)(\ast (f_{\infty}))$,
we naturally have
\[
 \nbigd_X\otimes_{V_D\nbigd_X}L(f)(D)
\simeq
\nbigd_X(\ast(f)_{\infty})
 \otimes_{V_D\nbigd_X(\ast (f)_{\infty})}L(f)(D).
\]
For any non-negative integer $m$,
we consider the $\nbigo_X(\ast(f)_{\infty})$-homomorphism
$\gamma_m:L(f)\bigl((m+1)D\bigr)
 \lrarr 
 \nbigd_X\otimes_{V_D\nbigd_X}L(f)(D)$
determined by
\[
 \gamma_m\Bigl(\prod_{i=1}^{\ell_1}z_i^{-1}
 \prod_{i=\ell_1+1}^{\ell_1+\ell_2}z_i^{-m-1}v\Bigr)
=(-1)^{\ell_2m}(m!)^{-\ell_2}
 \prod_{i=\ell_1+1}^{\ell_1+\ell_2}
 \del_i^{m}\otimes \prod_{i=1}^{\ell_1+\ell_2}z_i^{-1}v.
\]
Let $\iota_{m}:L(f)(mD)\lrarr L(f)((m+1)D)$
be the natural inclusion.
Then, we have
$\alpha_1\circ\gamma_{m}\circ\iota_m
=\alpha_1\circ\gamma_{m-1}$.
We obtain an $\nbigo_X(\ast (f)_{\infty})$-homomorphism
$\gamma:
 L(f)(\ast D)
\lrarr
 \nbigd_X\otimes_{V_D\nbigd_X}L(f)(D)$.
By the construction,
$\gamma$ is an epimorphism,
and $\alpha_1\circ\gamma$ is the identity.
Then, we obtain that $\alpha_1$ is an isomorphism.
\hfill\qed

\subsubsection{De Rham complexes}
\label{subsection;14.11.6.3}

We give some complexes
which are quasi-isomorphic to
the de Rham complexes of $L_{\star}(f,D)$
when $f$ is non-degenerate along $D$.
(See \cite{Esnault-Sabbah-Yu}, \cite{Sabbah-Yu},
\cite{Yu}
for the case $|(f)_0|\cap|(f)_{\infty}|=\emptyset$.)
We set $d_X:=\dim X$.

We have the natural complex of right $\nbigd_X$-modules
$\Omega_X^{\bullet}\otimes_{\nbigo_X} \nbigd_X[d_X]$
which is a right locally $\nbigd_X$-free resolution of 
$\Omega_X:=\Omega^{d_X}_X$.
We have the subcomplex
$\Omega_X^{\bullet}(\log D)(-D)\otimes_{\nbigo_X}
 V_D\nbigd_X[d_X]$,
which is a right locally $V_D\nbigd_X$-free resolution
of $\Omega_X$.
Indeed,
the natural morphism
$\Omega_X\otimes V_D\nbigd_X
\lrarr
 \Omega_X$
induces a quasi-isomorphism
$\Omega_X^{\bullet}(\log D)(-D)\otimes
 V_D\nbigd_X[d_X]
\simeq
 \Omega_X$
of complexes of right $V_D\nbigd_X$-modules.

Suppose that $f$ is non-degenerate along $D$.
Because 
$L_{\ast}(f,D)\simeq
 \nbigd_X\otimes_{V_D\nbigd_X} L(f)(D)$
according to Lemma \ref{lem;14.11.23.101},
we have the following natural isomorphisms:
\begin{multline}
\label{eq;14.12.9.2}
 \Omega_X\otimes^L_{\nbigd_X}
 L_{\ast}(f,D)
\simeq
 \Omega_X\otimes^L_{V_D\nbigd_X}
 L(f)(D)
\simeq
 \Omega_X^{\bullet}(\log D)(-D)\otimes L(f)(D)[d_X]
 \\
\simeq
 \Bigl(
 \Omega_X^{\bullet}(\log D)(\ast (f)_{\infty}),d+df
 \Bigr)[d_X]
\end{multline}
Because
$L_!(f,D)\simeq
 \nbigd_X\otimes_{V_D\nbigd_X}
L(f)$
according to Lemma \ref{lem;14.11.23.101},
we have the following isomorphisms:
\begin{multline}
\label{eq;14.12.9.3}
 \Omega_X\otimes^L_{\nbigd_X}
 L_!(f,D)
\simeq
 \Omega_X\otimes^L_{V_D\nbigd_X}
 L(f)
\simeq
\Omega_X^{\bullet}(\log D)(-D)
\otimes L(f)[d_X]
 \\
\simeq
 \Bigl(
 \Omega_X^{\bullet}(\log D)(-D)(\ast (f)_{\infty}),
 d+df
 \Bigr)[d_X]
\end{multline}

\paragraph{Coherent expression}
Let $H$ be any divisor in $X$.
Let us consider the complex
\[
 \Omega^{\ell-1}_X(\log D)\bigl(H\bigr)
\stackrel{a_1}{\lrarr}
 \Omega^{\ell}_X(\log D)\bigl(H+(f)_{\infty}\bigr)
\stackrel{a_2}{\lrarr}
 \Omega^{\ell+1}_X(\log D)\bigl(H+2(f)_{\infty}\bigr),
\]
where $a_i$ are induced by 
the multiplication of $df$.

\begin{lem}
If $f$ is non-degenerate along $D$,
we have $\Image a_1=\Ker a_2$
on a neighbourhood of 
$|(f)_{\infty}|$.
\end{lem}
\pf
We have only to check the claim
locally around any point $Q$ of $|(f)_{\infty}|$
by using a convenient local coordinate system.
Let us consider the case $Q\in |(f)_0|\cap|(f)_{\infty}|$.
Then, $f=z_n\prod_{i=1}^{\ell_1}z_i^{-k_i}$.
We have the following local section of
$\Omega_X^1(\log D)$:
\[
\tau:=
 df\cdot \prod_{i=1}^{\ell_1}z_i^{k_i}
=dz_n-\sum_{i=1}^{\ell_1} k_i z_n\frac{dz_i}{z_i}
\]
We have $\tau_{|Q}\neq 0$.
We consider the following 
on a neighbourhood of $Q$:
\[
 \Omega^{\ell-1}_X(\log D)
\stackrel{b_1}{\lrarr}
 \Omega^{\ell}_X(\log D)
\stackrel{b_2}{\lrarr}
 \Omega^{\ell+1}_X(\log D)
\]
Here, $b_i$ are induced by the multiplication of $\tau$.
Because $\tau_{|Q}\neq 0$,
we have $\Image(b_1)=\Ker(b_2)$.
Then, in the case $Q\in|(f)_0|$,
the claim of the lemma follows.
The case $Q\not\in|(f)_0|$
can be argued similarly,
and it is well known.
\hfill\qed

\vspace{.1in}

We set 
$\Omega_X^{\ell}(\log D,f):=
 \Omega_X^{\ell}(\log D)(\ell(f)_{\infty})$.
Together with the differential $d+df$,
we obtain a complex
$\bigl(
 \Omega_X^{\bullet}(\log D,f),d+df
 \bigr)$.
We also have
$\bigl(
 \Omega_X^{\bullet}(\log D,f)(-D),d+df
 \bigr)$.
We obtain the following from the previous lemma.
\begin{lem}
The following natural morphisms
are quasi-isomorphisms:
\[
  \bigl(
 \Omega_X^{\bullet}(\log D,f),d+df
 \bigr)
\lrarr
 \Bigl(
 \Omega_X^{\bullet}(\log D)(\ast (f)_{\infty}),d+df
 \Bigr)
\]
\[
   \bigl(
 \Omega_X^{\bullet}(\log D,f)(-D),d+df
 \bigr)
\lrarr
 \Bigl(
 \Omega_X^{\bullet}(\log D)(-D)(\ast (f)_{\infty}),d+df
 \Bigr)
\]
\hfill\qed
\end{lem}

Hence, we have the following natural quasi-isomorphisms:
\[
 \Omega_X\otimes^L_{\nbigd_X}L_{\ast}(f,D)
\simeq
 \bigl(
 \Omega_X^{\bullet}(\log D,f),d+df
 \bigr)[d_X]
\]
\[
 \Omega_X\otimes^L_{\nbigd_X}L_{!}(f,D)
\simeq
 \bigl(
 \Omega_X^{\bullet}(\log D,f)(-D),d+df
 \bigr)[d_X]
\]

\paragraph{Kontsevich complexes}

Let $\Omega_{X,f,D}^k$ denote the kernel of
the following morphism
induced by the multiplication of $df$:
\[
 \Omega^k_X(\log D)
\stackrel{df}{\lrarr}
 \frac{\Omega^{k+1}_X(\log D)\bigl(!(f)_{\infty}\bigr)}
 {\Omega^{k+1}_X(\log D)}.
\]
\begin{lem}
$\Omega_{X,f,D}^k$ are locally free.
\end{lem}
\pf
We have
$\Omega_{X,f,D}^k=\Omega_X^k(\log D)$
outside $|(f)_{\infty}|$.
Locally around any point of $(f)_{\infty}$,
we have local decompositions
$\Omega_X^k(\log D)
=A_k\oplus B_k$
such that the multiplication of 
$df$ induces $B_k\simeq A_{k+1}((f)_{\infty})$.
Hence, we have
$\Omega_{X,f,D}^k=A_k\oplus B_k(-(f)_{\infty})$,
and the claim of the lemma follows.
\hfill\qed

\vspace{.1in}
The multiplication of $df$
induces
$df:\Omega^{k}_{X,f,D}\lrarr \Omega^{k+1}_{X,f,D}$.
By the commutativity $[d,df]=0$,
the exterior derivative
induces
$d:\Omega^k_{X,f,D}\lrarr \Omega^{k+1}_{X,f,D}$.
Hence, we obtain the complex
$(\Omega^{\bullet}_{X,f,D},d+df)$.
We also obtain the complex
$(\Omega^{\bullet}_{X,f,D}(-D),d+df)$.

\begin{lem}
The following natural morphisms 
are quasi-isomorphisms:
\[
 \bigl(
 \Omega^{\bullet}_{X,f,D},d+df
 \bigr)
\lrarr
 \bigl(
 \Omega_X^{\bullet}(\log D,f),d+df
\bigr)
\]
\[
 \bigl(
 \Omega^{\bullet}_{X,f,D}(-D),d+df
 \bigr)
\lrarr
 \bigl(
 \Omega_X^{\bullet}(\log D,f)(-D),d+df
\bigr)
\]
\hfill\qed
\end{lem}

\subsubsection{The push-forward by a projection}
\label{subsection;14.11.6.2}

Let us consider the case 
$(X,D)=(X_0,D_0)\times S$.
Suppose that $X_0$ is compact.
Let $f$ be a meromorphic function
on $(X,D)$.
Let $\pi:X\lrarr S$ denote the projection.
Let $n:=\dim X_0$.

\begin{lem}
\label{lem;14.11.23.102}
Suppose that $f$ is non-degenerate along $D$.
We have the following natural isomorphisms:
\begin{equation}
 \label{eq;14.12.9.4}
 \pi_{+}^i\bigl(
 L_{\ast}(f,D)
 \bigr)
\simeq
 \hyperr^i\pi_{\ast}\bigl(
 \Omega^{\bullet+n}_{X/S}(\log D)(\ast (f)_{\infty}),
 d+df
 \bigr)
\end{equation}
\begin{equation}
\label{eq;14.12.9.5}
 \pi_{+}^i\bigl(
 L_{!}(f,D)
 \bigr)
\simeq
 \hyperr^i\pi_{\ast}\bigl(
 \Omega^{\bullet+n}_{X/S}(\log D)(-D)(\ast (f)_{\infty}),
 d+df
 \bigr)
\end{equation}
\end{lem}
\pf
Let us consider (\ref{eq;14.12.9.4}).
By Lemma \ref{lem;14.11.23.101},
we have 
$L_{\ast}(f,D)\simeq
 \nbigd_X\otimes_{V_D\nbigd_X} L(f)(D)$.
We have
$\nbigd_X=\nbigd_{X_0}\boxtimes \nbigd_{S}$
and 
$V_D\nbigd_X=V_{D_0}\nbigd_{X_0}\boxtimes \nbigd_{S}$.
Let $p:X\times X_0$ be the projection.
We obtain 
\begin{multline}
 \pi_{+}\bigl(
 L_{\ast}(f,D)
 \bigr)
\simeq
 R\pi_{\ast}\Bigl(
 p^{\ast}\Omega_{X_0}
 \otimes^L_{p^{\ast}\nbigd_{X_0}}
 L_{\ast}(f,D)
\Bigr)
\simeq
  R\pi_{\ast}\Bigl(
 p^{\ast}\Omega_{X_0}
 \otimes^L_{p^{\ast}V_D\nbigd_{X_0}}
 \bigl(
 L_{\ast}(f,D)\otimes\nbigo(D)
 \bigr)
\Bigr)
 \\
\simeq
 R\pi_{\ast}\Bigl(
 \Omega^{\bullet+n}_{X/S}(\log D)(-D)
 \otimes_{\nbigo_X}
 \bigl(
 L_{\ast}(f,D)\otimes\nbigo(D)
 \bigr)
\Bigr)
\simeq
  R\pi_{\ast}\Bigl(
 \Omega^{\bullet+n}_{X/S}(\log D)
 \otimes_{\nbigo_X}
 L_{\ast}(f,D)
\Bigr)
\end{multline}
It implies (\ref{eq;14.12.9.4}).
Similarly, we obtain (\ref{eq;14.12.9.5})
from the expression 
$L_!(f,D)\simeq \nbigd_X\otimes_{V_D\nbigo_X}
 \bigl(L(f)(-D)\bigr)$ in Lemma \ref{lem;14.11.23.101},
as in (\ref{eq;14.12.9.3}).
\hfill\qed

\vspace{.1in}

Suppose moreover that
$f$ is non-degenerate along $D$ over $S$.
(See Definition \ref{df;14.12.9.10}
for this stronger condition.)
As in \S\ref{subsection;14.11.6.3},
we can naturally define the complexes
$\bigl(\Omega^{\bullet}_{X/S}(\log D,f),d+df\bigr)$
and 
$\bigl(
 \Omega^{\bullet}_{X/S,f,D},d+df
 \bigr)$
in the relative setting.

\begin{lem}
If $f$ is non-degenerate along $D$ over $S$,
then we have the following natural isomorphisms:
\begin{equation}
\label{eq;14.6.27.10}
 \pi_{+}^i\bigl(
 L_{\ast}(f,D)
 \bigr)
\simeq
 \hyperr^i\pi_{\ast}
\bigl(\Omega^{\bullet+n}_{X/S}(\log D,f),d+df\bigr)
\simeq
  \hyperr^i\pi_{\ast}
 \bigl(
 \Omega^{\bullet+n}_{X/S,f,D},d+df
 \bigr)
\end{equation}
\begin{equation}
\label{eq;14.6.27.11}
 \pi_{+}^i\bigl(
 L_{!}(f,D)
 \bigr)
\simeq
 \hyperr^i\pi_{\ast}
\bigl(\Omega^{\bullet+n}_{X/S}(\log D,f)(-D),d+df\bigr)
\simeq
  \hyperr^i\pi_{\ast}
 \bigl(
 \Omega^{\bullet+n}_{X/S,f,D}(-D),d+df
 \bigr)
\end{equation}
\end{lem}
\pf
As in \S\ref{subsection;14.11.6.3},
we obtain the isomorphisms
from Lemma \ref{lem;14.11.23.102}.
\hfill\qed

\begin{cor}
\label{cor;14.7.10.10}
Suppose that $f$ is non-degenerate along $D$
over $S$.
Then, 
$\pi^i_{+}\bigl(
 L_{\star}(f,D)
 \bigr)$ are flat bundles on $S$,
i.e., locally free
 $\nbigo_S$-modules
 with an integrable connection.
\end{cor}
\pf
The right hand sides of
(\ref{eq;14.6.27.10})
and (\ref{eq;14.6.27.11})
are $\nbigo_S$-coherent.
Then, the claim of this corollary follows
from a well known result,
i.e.,
if a $\nbigd_X$-module is coherent over $\nbigo_X$,
then it is a flat bundle.
\hfill\qed

\subsection{Some functions
satisfying the purity condition}
\label{subsection;14.11.11.20}

We give some examples of meromorphic functions
which are not necessarily non-degenerate 
but satisfy the purity condition.
Let $X$ denote a complex manifold
with a normal crossing hypersurface $D$.

\subsubsection{Basic cases}

We set 
$X^{(1)}:=X\times\proj^1$
and 
$D^{(1)}:=(D\times\proj^1)
 \cup (\{\infty\}\times X)$.
Let $\tau$ be the standard coordinate on $\proj^1$.
Take $f,g\in\nbigo_X(\ast D)$,
and we consider a meromorphic function
$F:=\tau f+g$
on $(X^{(1)},D^{(1)})$.
We give some sufficient conditions
for $F$ to be pure on an open subset in $X^{(1)}$.

\begin{lem}
\label{lem;14.10.14.12}
Suppose the following.
\begin{itemize}
\item
 $g$ is non-degenerate along $D$.
\item
 $|(f)_0|\cap |(f)_{\infty}|=\emptyset$,
 and
 $|(f)_0|\subset |(g)_{\infty}|$.
 In particular,
 $|(f)_0|\subset D$.
\item
 $D=|(f)_{\infty}|\cup|(g)_{\infty}|$.
\end{itemize}
Then, $F$ is pure on $X^{(1)}$.
\end{lem}
\pf
Let us consider the 
$\nbigd_{X^{(1)}}(\ast D^{(1)})$-module
$\nbigv=\nbigo_{X^{(1)}}(\ast D^{(1)}) v$
with 
$\nabla v=v dF$.
It is enough to prove that
$\nbigv[!D^{(1)}]
\lrarr
 \nbigv$
is an epimorphism
around any point of 
$(P,\tau)\in X^{(1)}$.
We take a holomorphic coordinate system
$(x_1,\ldots,x_n)$ around $P$
such that
$D=\bigcup_{i=1}^{\ell}\{x_i=0\}$.

Let us consider the case $\tau\neq\infty$.
For a large $N$,
we have 
$v_1=\prod_{i=1}^{\ell} x_i^Nv
 \in\nbigv(\star D^{(1)})$
$(\star=\ast,!)$.
It is enough to prove that
$v_1$ generates $\nbigv$
as a $\nbigd_{X^{(1)}}$-module.
For any $\vecm\in\seisuu^{\ell}$,
we have 
\[
 \del_i
 (x^{\vecm}v)
=\bigl(x_i^{-1}m_i+\del_ig+\tau \del_if\bigr)
x^{\vecm}v.
\]
We also have
$\del_{\tau}v=fv$.
For $i=1,\ldots,\ell$,
we obtain
\begin{equation}
 \label{eq;14.11.4.2}
 x_i\del_i(x^{\vecm}v)
-(f^{-1}x_i\del_if)\cdot\tau\del_{\tau}
 (x^{\vecm}v)
=x^{\vecm}(x_i\del_ig+m_i)v.
\end{equation}
For $i=\ell+1,\ldots,n$,
we have
\begin{equation}
 \label{eq;14.11.4.3}
 \del_i(x^{\vecm}v)
-(f^{-1}\del_if)\tau\del_{\tau}
 (x^{\vecm}v)
=x^{\vecm}\del_i(g)v.
\end{equation}
By using (\ref{eq;14.11.4.2}) and (\ref{eq;14.11.4.3}),
we obtain that
$\nbigo_{X^{(1)}}(\ast (g)_{\infty})v_1$
is contained in 
$\nbigd_{X^{(1)}}v_1$.
Then, by using $\del_{\tau}v=fv$,
we obtain that
$\nbigv$ is contained in $\nbigd_{X^{(1)}}v_1$
around $(P,\tau)$.

Let us consider the case $\tau=\infty$.
Let $\kappa:=\tau^{-1}$.
It is enough to see that
$v_1=\kappa^N\prod_{i=1}^{\ell}x_i^Nv$
generates $\nbigv$ around $(P,\infty)$.
We have
$\kappa\del_{\kappa}v
=-\tau\del_{\tau}v
=-\kappa^{-1}fv$.
As in the case $\tau\neq\infty$,
we obtain 
$\nbigo_{X^{(1)}}(\ast (g)_{\infty})v_1$
is contained in 
$\nbigd_{X^{(1)}}v_1$
by using (\ref{eq;14.11.4.2}) and (\ref{eq;14.11.4.3}).
Then, by using 
$\kappa\del_{\kappa}v
=-\kappa^{-1}fv$,
we obtain that 
$\nbigv$ is contained in 
$\nbigd_{X^{(1)}}v_1$.
\hfill\qed

\vspace{.1in}
The following lemma is easy to see.
\begin{lem}
Suppose the following:
\begin{itemize}
\item $g=0$.
\item $D=|(f)_{\infty}|$.
\item $(f)_0$ is smooth and reduced,
 and $D\cup|(f)_0|$ is normal crossing.
\end{itemize}
Then, $F=\tau f$ is non-degenerate on
$\{\tau\neq 0\}\times X$.
In particular,
$F$ is pure on 
$\{\tau\neq 0\}\times X$.
\hfill\qed
\end{lem}

\subsubsection{A variant}

Let $X$ be a complex manifold
with a simple normal crossing hypersurface $D$.
We set $X^{(2)}:=\proj^1_{\tau}\times \proj^1_t\times X$.
Let $D^{(2)}$ denote the union of
$\{\infty\}\times\proj^1_t\times X$
and 
$\proj^1_{\tau}\times\{0,\infty\}\times X$
and
$\proj^1_{\tau}\times\proj^1_t\times D$.
Let $h$ be a meromorphic function on $(X,D)$
such that
(i) $D=|(h)_{\infty}|$,
(ii) $h$ is non-degenerate along $D$,
(iii) $(h)_0$ is reduced and non-singular on $X$.
We have the meromorphic function
$F=t^{-1}\tau+th$ on 
$(X^{(2)},D^{(2)})$.

We set $D_1:=\{t=0\}\times X$
and $D_2:=(\{t=\infty\}\times X)\cup(\proj^1_t\times D)$
in $\proj^1_t\times X$.
We take a projective birational morphism
$\varphi:Y\lrarr \proj^1_t\times X$
such that 
(i) $g:=\varphi^{\ast}(th)$ is non-degenerate,
(ii) $D_Y:=\varphi^{-1}(D_1\cup D_2)$ is normal crossing,
(iii) $Y\setminus \varphi^{-1}(D_1)\simeq 
 (\proj^1_t\times X)\setminus D_1$.
We set $f:=\varphi^{\ast}(t^{-1})$.
\begin{lem}
\label{lem;14.11.18.10}
$\Ftilde=\tau f+g$ is pure on $\Ytilde:=\proj^1_{\tau}\times Y$.
\end{lem}
\pf
Because $D_2\setminus D_1\subset |(th)_{\infty}|$,
we have
$\varphi^{-1}(D_2\setminus D_1)
\subset
 \bigl|(g)_{\infty}\bigr|$.
We clearly have
$\varphi^{-1}(D_1)=
 \bigl|
 (f)_{\infty}
 \bigr|$.
Hence, we have
$D_Y=|(f)_{\infty}|\cup|(g)_{\infty}|$.
We also have
$|(f)_0|\cap|(f)_{\infty}|=\emptyset$
and 
$|(f)_0|\subset|(g)_{\infty}|$.
Then, we obtain the claim of Lemma \ref{lem;14.11.18.10}
from Lemma \ref{lem;14.10.14.12}.
\hfill\qed

\begin{lem}
\label{lem;14.10.14.11}
$F$ is pure on 
$\proj^1_{\tau}\times\proj^1_t\times X$.
\end{lem}
\pf
Let $\varphitilde:\Ytilde\lrarr X^{(2)}$
be the induced morphism.
Let $\Dtilde_Y$ be the union of
$\{\infty\}\times Y$
and $\proj^1_{\tau}\times D_Y$
in $\Ytilde$.
By the previous lemma,
$\Ftilde:=\tau f+g$ is pure on $\Ytilde$.
We also have
$\Dtilde_Y=|(\Ftilde)_{\infty}|$.
Hence,
the natural morphism
$L_!(\Ftilde,\Dtilde_Y)
\lrarr
 L_{\ast}(\Ftilde,\Dtilde_Y)$
is an isomorphism.
Because 
$L_{\star}(F,D^{(2)})\simeq
 \varphitilde_{\dagger}
 L_{\star}(\Ftilde,\Dtilde_Y)$,
we obtain that
$L_{!}(F,D^{(2)})
\lrarr
 L_{\ast}(F,D^{(2)})$
is an isomorphism,
i.e., $F$ is pure on $X^{(2)}$.
\hfill\qed

\subsection{Push-forward}
\label{subsection;14.7.2.1}

\subsubsection{A statement}
\label{subsection;14.10.14.20}

Let $X$ be a complex manifold
with a simple normal crossing hypersurface $D$.
We set $Y:=X\times\proj^1$
and $D_Y^{(0)}:=
 (X\times\{\infty\})\cup
 (D\times \proj^1)$
and $D_Y^{(1)}:=D_Y^{(0)}\cup(X\times\{0\})$.
Let $f$ and $g$ be meromorphic functions on $(X,D)$.
We assume the following.
\begin{itemize}
\item
The divisor 
$(f)_0\cap (X\setminus D)$ is reduced
and non-singular.
\end{itemize}
We set $Z_f:=|(f)_0|$.
Let $[z_0:z_1]$ be a homogeneous coordinate system
of $\proj^1$,
and we set $t:=z_0/z_1$.
We obtain the meromorphic function $tf$ on $Y$.
The pull back of $g$ by the projection $Y\lrarr X$
is also denoted by $g$.
We set $F:=tf+g$.

We have the $\nbigd$-modules
$L_{\ast}(F,D_Y^{(0)})$ 
and $L_{\ast}(F,D_Y^{(1)})$
on $Y$.
We have the natural exact sequence:
\begin{equation}
 \label{eq;14.6.27.30}
 0\lrarr
 L_{\ast}(F,D_Y^{(0)})
\lrarr
 L_{\ast}(F,D_Y^{(1)})
\lrarr
 L_{\ast}(F,D_Y^{(1)})
 \big/
 L_{\ast}(F,D_Y^{(0)})
\lrarr 0
\end{equation}

Let $\pi:Y\lrarr X$ be the projection.
We shall prove the following proposition
in \S\ref{subsection;14.7.2.10}.
\begin{prop}
\label{prop;14.7.1.1}
We have
$\pi^{i}_{+}(M)=0$ $(i\neq 0)$
for 
\[
M=L_{\ast}(F,D_Y^{(0)}),\,\,
L_{\ast}(F,D_Y^{(1)}),\,\,
L_{\ast}(F,D_Y^{(1)})
 \big/L_{\ast}(F,D_Y^{(0)}).
\]
We have the following isomorphisms:
\begin{equation}
\label{eq;14.6.27.40}
 \pi_{+}^0L_{\ast}(F,D_Y^{(0)})
\simeq
 L_{\ast}(g,D)(\ast Z_f)/L_{\ast}(g,D)
\simeq
 \Ker\Bigl(
 L_{\ast}(g)(!Z_f)(\ast D)
\lrarr
 L_{\ast}(g,D)
 \Bigr)
\end{equation}
\begin{equation}
\label{eq;14.6.27.41}
 \pi_{+}^0L_{\ast}(F,D_Y^{(1)})
\simeq
 L_{\ast}(g)(!Z_f)(\ast D)
\end{equation}
\begin{equation}
\label{eq;14.6.27.42}
 \pi_{+}^0\bigl(
 L_{\ast}(F,D_Y^{(1)})\big/
 L_{\ast}(F,D_Y^{(0)})
 \bigr)
\simeq
 L_{\ast}(g,D)
\end{equation}
The push-forward of {\rm(\ref{eq;14.6.27.30})}
is isomorphic to the following standard exact sequence
\[
 0\lrarr
 \Ker\Bigl(
 L_{\ast}(g)(!Z_f)(\ast D)
\lrarr
 L_{\ast}(g,D)
 \Bigr)
\lrarr
  L_{\ast}(g)(!Z_f)(\ast D)
\lrarr 
 L_{\ast}(g,D)
\lrarr 0
\]
\end{prop}

By the duality,
we obtain the following
as a corollary of Proposition \ref{prop;14.7.1.1}.
\begin{cor}
\label{co;14.11.26.22}
We have 
$\pi_{+}^iM=0$
$(i\neq 0)$
for 
\[
 M=L_{!}(F,D_Y^{(0)}),\,\,
L_{!}(F,D_Y^{(1)}),\,\,
\Ker\bigl(
L_{!}(F,D_Y^{(1)})
\lrarr
L_{!}(F,D_Y^{(0)})
 \bigr).
\]
We have the following isomorphisms:
\[
 \pi_+^0\bigl(
 L_!(F,D_Y^{(0)})
 \bigr)
\simeq
 \Ker\Bigl(
 L_!(g,D)(!Z_f)
\lrarr
 L_!(g,D)
 \Bigr)
\simeq
 \Cok\Bigl(
 L_!(g,D)
\lrarr
 L_!(g)(\ast Z_f)(!D)
 \Bigr)
\]
\[
 \pi_+^0\bigl(
 L_!(F,D_Y^{(1)})
 \bigr)
\simeq
 L_!(g)(\ast Z_f)(!D)
\]
\[
 \pi_+^0\bigl(
 \Ker\bigl(
 L_!(F,D_Y^{(1)})
\lrarr
 L_!(F,D_Y^{(0)})
 \bigr)
 \bigr)
\simeq
 L_!(g,D)
\]
The push-forward of 
\begin{equation}
 \label{eq;14.11.26.10}
 0\lrarr
\Ker\bigl(
L_{!}(F,D_Y^{(1)})
\lrarr
L_{!}(F,D_Y^{(0)})
 \bigr)
\lrarr
L_{!}(F,D_Y^{(1)})
\lrarr
L_{!}(F,D_Y^{(0)})
\lrarr 0
\end{equation}
is isomorphic to the standard exact sequence:
\begin{equation}
\label{eq;14.11.26.11}
0\lrarr 
 L_!(g,D)
\lrarr
 L_!(g)(\ast Z_f)(!D)
\lrarr 
 \Cok\bigl(
  L_!(g,D)
\lrarr
 L_!(g)(\ast Z_f)(!D)
\bigr)
\lrarr 0
\end{equation}
\hfill\qed
\end{cor}

\begin{rem}
\label{rem;14.11.26.23}
Note that the restriction of
 the morphism
$L_{!}(F,D_Y^{(0)})
\lrarr
 L_{\ast}(F,D_Y^{(0)})$
to $\proj^1\times (X\setminus D)$
is an isomorphism.
Indeed, the restriction of $F$ to
$\proj^1\times (X\setminus D)$
is non-degenerate
along $\{\infty\}\times(X\setminus D)$.
In particular,
the restriction of the morphism
$\pi^0_{+}L_!(F,D_Y^{(0)})
\lrarr
 \pi^0_{+}L_{\ast}(F,D_Y^{(0)})$
to $X\setminus D$ is an isomorphism.
Indeed, 
$\pi^0_{+}L_{\star}(F,D_Y^{(0)})_{|X\setminus D}$
are isomorphic to
the push-forward of the $\nbigd_{Z_f\setminus D}$-module
$(\nbigo_{Z_f\setminus D},d+dg')$
via the inclusion $Z_f\setminus D\lrarr X\setminus D$,
where $g'$ is the restriction of $g$ to $Z_f\setminus D$.
\hfill\qed
\end{rem}

Let us consider the case that
$f$ is moreover non-degenerate along $D$.
In this case,
$Z_f$ is smooth, and
$Z_f\cup D$ is normal crossing.
Let $\iota:Z_f\lrarr X$ denote the inclusion.
We set $D_{Z_f}:=D\cap Z_f$.
We set $g_0:=g_{|Z_f}$.
We have the $\nbigd$-modules
$L_{\star}(g_0,D_{Z_f})$ on $Z_f$.

\begin{cor}
\label{cor;14.11.26.2}
If $f$ is moreover non-degenerate along $D$,
we have the following isomorphisms:
\begin{equation}
 \label{eq;14.11.26.20}
 \pi^0_{+}L_{\ast}(F,D_Y^{(0)})
\simeq
 L_{\ast}(g,D)(\ast Z_f)/L_{\ast}(g,D)
 \simeq
\iota_+
 L_{\ast}(g_0,D_{Z_f})
\end{equation}
\begin{equation}
 \label{eq;14.11.26.21}
 \pi^0_{+}L_{!}(F,D_Y^{(0)})
\simeq
\Ker\bigl(
 L_!(g,D)(!Z_f)
\lrarr
 L_!(g,D)
 \bigr)
\simeq
 \iota_+
 L_!(g_0,D_{Z_f})
\end{equation}
The image of 
$\pi^0_+L_!(F,D_Y^{(0)})
\lrarr 
 \pi^0_+L_{\ast}(F,D_Y^{(0)})$
is naturally isomorphic to
$\iota_+L(g_0)$.
Note that it is 
$\iota_{+}\nbigo_{Z_f}$
if $g=0$.
\end{cor}
\pf
The first isomorphism in (\ref{eq;14.11.26.20})
is given in Proposition \ref{prop;14.7.1.1}.
If $Z_f\cup D$ is normal crossing,
we clearly have the second isomorphism
in (\ref{eq;14.11.26.20}).
We obtain the isomorphisms in (\ref{eq;14.11.26.21}) 
by the duality.
According to Remark \ref{rem;14.11.26.23},
the restriction of the morphism
$\pi^0_+L_!(F,D_Y^{(0)})
\lrarr 
 \pi^0_+L_{\ast}(F,D_Y^{(0)})$
to $X\setminus D$ is an isomorphism.
Hence,
the image of 
$\pi^0_+L_!(F,D_Y^{(0)})
\lrarr 
 \pi^0_+L_{\ast}(F,D_Y^{(0)})$
is identified
with the image of
a non-zero morphism
$\iota_{+}L_{!}(g_0,D_{Z_f})
\lrarr
 \iota_+L_{\ast}(g_0,D_{Z_f})$,
which is isomorphic to 
$\iota_{+}L(g_0)$.
\hfill\qed

\subsubsection{Extensions}
\label{subsection;15.10.20.1}

We give a general preliminary
which is a variant of Beilinson's construction
\cite{beilinson2}.
Let $h$ be any meromorphic function
on a complex manifold $Z$.
We set $D_h:=|(h)_0|\cup|(h)_{\infty}|$.
Let $\nbigo_Z(\ast D_h)[s]:=
 \nbigo_Z(\ast D_h)\otimes_{\cnum}\cnum[s]$.
For any pair of integers $(a,b)$ with $a< b$,
we consider the meromorphic flat bundle
\[
 I_h^{a,b}:=
 s^a\nbigo_Z(\ast D_h)[s]
 \bigl/
 s^b\nbigo_Z(\ast D_h)[s]
\]
with the connection $\nabla$
determined by $\nabla s^j=s^{j+1}dh/h$.
For any holonomic $\nbigd_Z$-module $M$,
we set 
$\Pi^{a,b}_hM:=I_h^{a,b}\otimes_{\nbigo_Z} M$.

\vspace{.1in}

We return to the situation in \S\ref{subsection;14.10.14.20}.
We shall prove the following proposition in 
\S\ref{subsection;14.12.9.20}.
\begin{prop}
\label{prop;14.10.14.21}
Suppose that
 $(f)_0\cap X\setminus D$ is smooth and reduced.
Then, 
we have
$\pi_{+}^j\bigl(
  \Pi^{a,b}_{t}L_{\star}(F)(\star t)
 \bigr)=0$ $(\star=\ast,!)$
for $j\neq 0$,
and 
we have the following natural isomorphisms:
\[
 \Psi^{a,b}_{\ast}:
 \pi_{+}^0\bigl(
 \Pi^{a,b}_{t}L_{\ast}(F)(\ast t)
 \bigr)
\simeq
 \Bigl(
 \Pi^{a,b}_{f^{-1}}
 L_{\ast}(g)
 \Bigr)
 (!(f)_0)
 (\ast D)
\]
\[
 \Psi^{a,b}_{!}:
 \pi_{+}^0\bigl(
 \Pi^{a,b}_{t}L_{!}(F)(!t)
 \bigr)
\simeq
 \Bigl(
 \Pi^{a,b}_{f^{-1}}
 L_!(g)
 \Bigr)
 (\ast (f)_0)
 (!D)
\]
The following diagram is commutative:
\begin{equation}
\label{eq;15.11.6.100}
 \begin{CD}
\pi_{+}^0\bigl(
 \Pi^{a,b}_{t}L_{!}(F)(!t)
 \bigr)
@>>>
\pi_{+}^0\bigl(
 \Pi^{a,b}_{t}L_{\ast}(F)(\ast t)
 \bigr)
\\
 @V{\Psi^{a,b}_!}VV @V{\Psi^{a,b}_{\ast}}VV 
\\
 \Bigl(
 \Pi^{a,b}_{f^{-1}}
 L_!(g)
 \Bigr)
 (\ast (f)_0)
 (!D)
@>>>
 \Bigl(
 \Pi^{a,b}_{f^{-1}}
 L_{\ast}(g)
 \Bigr)
 (!(f)_0)
 (\ast D)
 \end{CD}
\end{equation}
Here, the upper horizontal arrow is the natural morphism,
and the lower horizontal arrow is given by 
the multiplication of $-s$.
\end{prop}

\subsubsection{Proof of Proposition \ref{prop;14.7.1.1}}
\label{subsection;14.7.2.10}

Let us obtain the first isomorphism in (\ref{eq;14.6.27.40}).
We set $\nbigo_X(\ast D)[t]:=
\nbigo_X(\ast D)\otimes_{\cnum}\cnum[t]$.
Let $R\pi_{\ast}$ denote the ordinary push-forward of sheaves
by $\pi$.
We have
$R^j\pi_{\ast}\bigl(L_{\ast}(F,D_Y^{(0)})\bigr)=0$ for $j>0$,
and 
$\pi_{\ast}\bigl(L_{\ast}(F,D_Y^{(0)})\bigr)
\simeq
 \nbigo_X(\ast D)[t]$
as an $\nbigo_X$-module.
Hence, 
$\pi_{+}L_{\ast}(F)$
is represented by the following complex:
\[
 \nbigo_X(\ast D)[t]
 \stackrel{\del_t+f}\lrarr
 \nbigo_X(\ast D)[t]
\]
Here, the second term sits in the degree $0$.
The action of vector fields $V$ on $X$
is given by
$V(t^i)=t^{i}(tV(f)+V(g))$.
It is easy to see that
the morphism
$\nbigo_X(\ast D)[t]
 \stackrel{\del_t+f}\lrarr
 \nbigo_X(\ast D)[t]$
is a monomorphism.
Let us look at the cokernel.
We consider the following morphism
of $\nbigo_X(\ast D)$-modules:
\[
 \nbigo_X(\ast D)[t]
\lrarr
 L_{\ast}(g,D)(\ast Z_f),
\quad\quad
 t^i\longmapsto
 (-1)^ii!f^{-i-1}
\]
Here, 
we use the natural identification
$L_{\ast}(g,D)(\ast Z_f)=\nbigo_X(\ast D)(\ast Z_f)$
as an $\nbigo_X$-module.
It is a morphism of $\nbigd_X$-modules.
Indeed, 
for a holomorphic vector field $V$ on $X$,
we have
\[
V(t^i)=t^{i+1}V(f)+t^iV(g),
\quad
V(f^{-i-1}(-1)^ii!)
=(-1)^{i+1}(i+1)!f^{-i-2}V(f)
+V(g)f^{-i-1}(-1)^i i!.
\]
For $i>0$, we have
$\del_tt^i+ft^i=it^{i-1}+ft^i$
which is mapped to
$i(-1)^{i-1}(i-1)!f^{-i}
+f(-1)^ii!f^{-i-1}
=0$.
We also have
$\del_tt^0+ft^0=f$
which is mapped to $1$.
Hence, we obtain the morphism of
the $\nbigd$-modules:
\begin{equation}
\label{eq;14.7.10.1}
 \Cok\Bigl(
 \nbigo_X(\ast D)[t]
\stackrel{\del_t+f}{\lrarr}
 \nbigo_X(\ast D)[t]
 \Bigr)
\lrarr
 L_{\ast}(g,D)(\ast Z_f)\big/
 L_{\ast}(g,D)
\end{equation}
We consider filtrations
$\nbigf_j(\nbigo_X(\ast D)[t])
=\bigoplus_{i\leq j}\nbigo_X(\ast D)t^iv$
and 
$\nbigf_jL_{\ast}(g,D)(\ast Z_f)
=L_{\ast}(g,D)\bigl((j+1)Z_f\bigr)$.
It is easy to see that
the induced morphism
on the graded modules is an isomorphism.
Hence, (\ref{eq;14.7.10.1}) is an isomorphism.
Thus, we obtain the first half of (\ref{eq;14.6.27.40}).

\vspace{.1in}

Let us prove (\ref{eq;14.6.27.42}).
Let $\nbigo_X(\ast D)[t,t^{-1}]:=
 \nbigo_X(\ast D)\otimes_{\cnum}
 \cnum[t,t^{-1}]$.
Let $\iota:X\times\{0\}\lrarr X\times\proj^1$
denote the inclusion.
We identify 
$L_{\ast}(F,D^{(1)}_Y)\big/L_{\ast}(F,D^{(0)}_Y)$
with the $\nbigo_{Y}$-module
\[
 \nbigm_1:=
\iota_{\ast}
\Bigl(
 \nbigo_X(\ast D)[t,t^{-1}]
 \Big/
\nbigo_X(\ast D)[t]
\Bigr)\,v
\]
and the connection $\nabla$
given by $\nabla v=v d(F)$.
Then,
$\pi_{+}
 \bigl(L_{\ast}(F,D^{(1)}_Y)\big/
 L_{\ast}(F,D^{(0)}_Y)
 \bigr)$
is represented by the complex
$\nbigm_1
 \stackrel{\del_t+f}{\lrarr}
 \nbigm_1$.
Here, the second term sits in the degree $0$.
The kernel and the cokernel
are denoted by $\Ker_1$ and $\Cok_1$,
respectively.
For the identification
$\nbigm_1=\bigoplus_{j=1}^{\infty}
 \nbigo_X(\ast D)t^{-j}v$,
we set
$\nbigf_n:=\bigoplus_{j\leq n}
 \nbigo_X(\ast D)t^{-j}v$.
We have
$\del_t+f:
 \nbigf_n\lrarr\nbigf_{n+1}$,
which induces an isomorphism
$\Gr^{\nbigf}_n\simeq\Gr^{\nbigf}_{n+1}$
for $n\geq 1$.
Hence, it is easy to check 
$\Ker_1=0$,
and that the natural 
$\nbigo_X$-morphism
$L_{\ast}(g,D)\lrarr
 \nbigo_X(\ast D)t^{-1}v=\nbigf_1$
given by $\varphi\longmapsto \varphi (t^{-1}v)$
induces an isomorphism
$L_{\ast}(g,D)\lrarr \Cok_1$
which is compatible with the flat connection.

\vspace{.1in}

Let us observe (\ref{eq;14.6.27.41}).
We consider
$\nbigm_2=\nbigo_X(\ast D)[t,t^{-1}]v$
with $\nabla v=vdF$.
As before,
$\pi_{+}\bigl(
 L_{\ast}(F,D_Y^{(1)})
 \bigr)$
is represented by
$\nbigm_2\stackrel{\del_t+f}\lrarr
 \nbigm_2$,
where the second term sits in the degree $0$.
The kernel and the cokernel
are denoted by $\Ker_2$ and $\Cok_2$,
respectively.
We set
$\nbigm_0:=\nbigo_X(\ast D)[t]v$.
We have the natural exact sequence
$0\lrarr
 \nbigm_0
\lrarr
 \nbigm_2
\lrarr
 \nbigm_1
\lrarr 0$.
Because 
$\nbigm_i\stackrel{\del_t+f}{\lrarr}\nbigm_i$
$(i=0,1)$ are monomorphisms,
we obtain that $\Ker_2=0$.
We also have the following exact sequence
of $\nbigd$-modules:
\[
 0\lrarr
L_{\ast}(g,D)(\ast Z_f)/L_{\ast}(g,D)
\lrarr
 \Cok_2
\lrarr
 L_{\ast}(g,D)
\lrarr 0
\]
It implies that
$\Cok_2(\ast Z_f)\simeq
 L_{\ast}(g,D)(\ast Z_f)$.
Hence, we have
a uniquely induced morphism
$ L_{\ast}(g,D)(!Z_f)
\lrarr
 \Cok_2$.
Because $\Cok_2(\ast D)\simeq \Cok_2$,
we have a uniquely induced morphism
\[
  \rho:
 L_{\ast}(g)(!Z_f)(\ast D)
\simeq
 L_{\ast}(g,D)(!Z_f)(\ast D)
\lrarr
 \Cok_2.
\]
Let us prove that the morphism
is an isomorphism.
Around any point of $Z_f\setminus D$,
we have
$\psi_f(\Cok_2)\simeq
 \psi_f\bigl(L_{\ast}(g)\bigr)$
and 
$\phi_f(\Cok_2)\simeq
 \phi_f\bigl(
 L_{\ast}(g)(\ast Z_f)/L_{\ast}(g)
 \bigr)$.
By a direct computation,
we can check that
the natural morphism
$\psi_f(\Cok_2)\lrarr
 \phi_f(\Cok_2)$ is non-zero.
It implies that 
$\psi_f(\Cok_2)\lrarr
 \phi_f(\Cok_2)$ is an isomorphism.
Hence, $\rho$ is an isomorphism
on $X\setminus D$,
and then it follows that
$\rho$ is an isomorphism on $X$.
We also obtain the second isomorphism
in (\ref{eq;14.6.27.40}),
and that the push-forward of
(\ref{eq;14.11.26.10})
is isomorphic to (\ref{eq;14.11.26.11}).
\hfill\qed

\subsubsection{Proof of Proposition \ref{prop;14.10.14.21}}
\label{subsection;14.12.9.20}

Let us prove the claim in the case $\star=\ast$.
We have the following exact sequence:
\[
 0\lrarr
 \Pi^{a-1,b}_t L_{\ast}(F)(\ast t)
\lrarr
  \Pi^{a,b}_t L_{\ast}(F)(\ast t)
\lrarr
 L_{\ast}(F)(\ast t)\,s^a
\lrarr 0
\]
Then, we obtain
$\pi_{+}^j
   \Pi^{a,b}_t L_{\ast}(F)(\ast t)=0$ $(j\neq 0)$
by applying Proposition \ref{prop;14.7.1.1}
and an easy induction.
We also have the exact sequence:
\[
 0\lrarr
 \pi_{+}^0
 \Pi^{a-1,b}_t L_{\ast}(F)(\ast t)
\lrarr
 \pi_{+}^0
  \Pi^{a,b}_t L_{\ast}(F)(\ast t)
\lrarr
 L_{\ast}(g)(!(f)_0)(\ast D)
\lrarr 0
\]
By an easy induction,
we obtain that 
the following natural morphisms
are isomorphisms:
\[
 \pi_{+}^0
 \bigl(
  \Pi^{a,b}_t L_{\ast}(F)(\ast t)
 \bigr)
\stackrel{\simeq}{\lrarr}
   \pi_{+}^0
 \bigl(
  \Pi^{a,b}_t L_{\ast}(F)(\ast t)
 \bigr)(\ast D)
\stackrel{\simeq}{\llarr}
    \pi_{+}^0
 \bigl(
  \Pi^{a,b}_t L_{\ast}(F)(\ast t)
 \bigr)(!(f)_0)(\ast D).
\]
Hence, for the construction of $\Psi^{a,b}_{\ast}$,
it is enough to obtain an isomorphism
$\pi_{+}^0\bigl(
 \Pi^{a,b}_tL_{\ast}(F)(\ast t)
\bigr)
\simeq
 \Pi_{f^{-1}}^{a,b}L_{\ast}(g)$
on $X\setminus (D\cup Z_f)$.

We have the following representative of
$\pi_{+}\bigl(
 \Pi_{t}^{a,b}L_{\ast}(F)(\ast t)
 \bigr)$:
\[
 \bigoplus_{j=a}^{b-1}
\nbigo_X(\ast D) [t]s^j
\stackrel{\kappa}{\lrarr}
 \bigoplus_{j=a}^{b-1}
\nbigo_X(\ast D)[t]t^{-1}s^j 
\]
Here, 
the morphism $\kappa$ is given by
$\del_{t}+f+st^{-1}$.
The action of any vector field $V\in\Theta_X$
is given by 
$V(s^j)=\bigl(V(g)+tV(f)\bigr)s^j$.
The kernel of $\kappa$ is clearly $0$.
The natural inclusion of
$\bigoplus_{j=a}^{b-1}
 L_{\ast}(g)t^{-1}s^j$
into 
$\bigoplus_{j=a}^{b-1}
 L_{\ast}(g)[t]t^{-1}s^j$
induces an isomorphism of
$\bigoplus_{j=a}^{b-1}
 L_{\ast}(g)t^{-1}s^j$
and the cokernel of $\kappa$.
We have
$V(s^jt^{-1})
=V(g)s^jt^{-1}+V(f)s^j$.
Because
$\kappa(s^j)
=fs^j+s^{j+1}t^{-1}$,
we have
$V(f)s^j\equiv
 fV(f^{-1})s^{j+1}t^{-1}$
on $X\setminus \bigl(Z_f\cup D\bigr)$.
Thus, 
we obtain the isomorphism
$\Psi^{a,b}_{\ast}:
\pi_{+}^0\bigl(
 \Pi^{a,b}_tL_{\ast}(F)(\ast t)
\bigr)
\simeq
 \Pi_{f^{-1}}^{a,b}L_{\ast}(g)$
on $X\setminus (D\cup Z_f)$
by setting
$\Psi^{a,b}_{\ast}(t^{-1}s^j):=s^j$.

Let us consider the case $\star=!$.
We can deduce it from the claim in the case
$\star=\ast$ by using the duality.
But, we give a more explicit construction 
which would be useful for our study later.
As in the case of $\star=\ast$,
we have 
$\pi_{+}^j
 \Pi_t^{a,b}L_!(F)(!t)=0$ $(j\neq 0)$
and the following natural isomorphisms:
\[
  \pi_{+}^0
 \bigl(
  \Pi^{a,b}_t L_{!}(F)(!t)
 \bigr)
\stackrel{\simeq}{\lrarr}
   \pi_{+}^0
 \bigl(
  \Pi^{a,b}_t L_{!}(F)(!t)
 \bigr)(!D)
\stackrel{\simeq}{\llarr}
    \pi_{+}^0
 \bigl(
  \Pi^{a,b}_t L_{!}(F)(!t)
 \bigr)(\ast Z_f)(!D).
\]
Hence, it is enough to obtain an isomorphism
$\pi_{+}^0\bigl(
 \Pi^{a,b}_tL_{!}(F)(!t)
\bigr)
\simeq
 \Pi_{f^{-1}}^{a,b}L_{!}(g)$
on $X\setminus (D\cup Z_f)$.

We have the following representative of
$\pi_{+}\bigl(
 \Pi_{t}^{a,b}L_{!}(F)(! t)
 \bigr)$
on $X\setminus\bigl(Z_f\cup D\bigr)$:
\[
 \bigoplus_{j=a}^{b-1}
 \nbigo_X[t]ts^j
\stackrel{\kappa}{\lrarr}
 \bigoplus_{j=a}^{b-1}
 \nbigo_X[t]s^j 
\]
The morphism $\kappa$
and the action of $V\in\Theta_X$
are given by the same formula.
The natural inclusion of
$\bigoplus_{j=a}^{b-1}
 L_{\ast}(g)s^j$
into 
$\bigoplus_{j=a}^{b-1}
 L_{\ast}(g)[t]s^j$
induces an isomorphism of
$\bigoplus_{j=a}^{b-1}
 L_{\ast}(g)s^j$
and the cokernel of $\kappa$.
On $X\setminus \bigl(Z_f\cup D\bigr)$,
we have 
$V(s^j)
=V(g)s^j+V(f)ts^j$
and 
$\kappa(s^{j}t)=fs^{j}t+s^{j+1}+s^j$.
Hence, we have
$V(fs^j)\equiv V(g)(fs^j)+fV(f^{-1})fs^{j+1}$.
Thus, we obtain the isomorphism
$\Psi^{a,b}_!$ by setting
$\Psi^{a,b}_!(fs^j):=s^j$.

Let us look at the restriction of the natural morphism
$\pi_{+}^0
 \Pi^{a,b}_tL_!(F)(!t)
\lrarr
\pi_{+}^0
 \Pi^{a,b}_tL_{\ast}(F)(\ast t)$
to $X\setminus (D\cup Z_f)$.
Because $\kappa(s^j)=fs^j+s^{j+1}t^{-1}$,
we have
$fs^j=-t^{-1}s^{j+1}$
in 
$\pi_{+}^0
 \Pi^{a,b}_tL_{\ast}(F)(\ast t)$.
Hence, under the above isomorphisms
$\Psi^{a,b}_!$ and $\Psi^{a,b}_{\ast}$,
it is identified with
$\Pi^{a,b}_tL(g)
\lrarr
 \Pi^{a,b}_tL(g)$
induced by the multiplication of $-s$,
i.e.,
the diagram (\ref{eq;15.11.6.100})
is commutative.
\hfill\qed

\begin{rem}
We also have the following representative
of $\pi^0_+\bigl(
 \Pi^{a,b}_{t}L_{\ast}(F)(!t)
\bigr)$:
\[
 \bigoplus_{j=a}^{b-1}
 \nbigo_X(\ast D)(!t)s^j
\stackrel{\kappa}\lrarr
 \bigoplus_{j=a}^{b-1}
 \nbigo_X(\ast D)(!t)s^j
\]
For any local section $c$ of
$\bigoplus_{j=a}^{b-1}
 \nbigo_X(\ast D)(!t)s^j$,
let $[c]_!$ denote the induced section of 
$\pi^0_+\bigl(
 \Pi^{a,b}_{t}L_{\ast}(F)(!t)
 \bigr)$.
Recall that we have a natural isomorphism
$L(F)(!t)\simeq L(F)\otimes \nbigo_X(!t)$.

\end{rem}

\subsubsection{A consequence}
\label{subsection;14.7.1.12}

We continue to use the notation in
\S\ref{subsection;14.10.14.20}.
For simplicity,
we consider the case $g=0$.
Motivated by the descriptions of some hypergeometric systems
in \cite{Batyrev-VMHS} and \cite{Stienstra},
we give a remark on another description of
$\pi_{+}^0\bigl(
 L_{\ast}(tf,D_Y^{(0)})
 \bigr)$.
We set 
$Y_0:=X\times\proj^1_t\times\proj^1_s$.
Let $p_i$ denote the projection of
$Y_0$ onto the $i$-th component.
We set
$D_0:=p_1^{-1}(D)
 \cup
 p_2^{-1}(\{0,\infty\})
 \cup
 p_3^{-1}(\{0,\infty\})$.
We regard $t$ and $s$ 
as meromorphic functions on $(Y_0,D_0)$.
Let us consider $L_{\ast}\bigl(t(f+s),D_0\bigr)$.

Let $\pi_{ts}:X\times\proj^1_t\times\proj^1_s\lrarr X$,
$\pi_t:X\times\proj^1_t\lrarr X$
and $\pi_s:X\times\proj^1_s\lrarr X$
denote the projections.
Set
$D_{X\times\proj^1_t}:=
 \bigl(
 X\times\{0,\infty\}
 \bigr)
 \cup
 \bigl(
 D\times\proj^1_t
 \bigr)$,
and 
$D_{X\times\proj^1_s}:=
 \bigl(
 X\times\{0,\infty\}
 \bigr)
 \cup
 \bigl(
 D\times\proj^1_s
 \bigr)$.

\begin{prop}
We have the following natural isomorphisms:
\begin{equation}
\label{eq;14.7.16.2}
 \pi_{s+}\nbigo_{X\times\proj^1_s}
 \bigl(!(f+s)_0\bigr)
 \bigl(\ast D_{X\times\proj^1_s}\bigr)
\simeq
 \pi_{ts+}L_{\ast}\bigl(
 t(f+s),D_0
 \bigr)
\simeq
 \pi_{t+}L_{\ast}(tf,D_{X\times\proj^1_t})
\simeq
  \nbigo_{X}(!(f)_0)(\ast D).
\end{equation}
\end{prop}
\pf
Let $p_{12}$ and $p_{13}$
denote the projection 
of $Y_0$ onto
$X\times\proj^1_t$
and 
$X\times\proj^1_s$,
respectively.
By Proposition \ref{prop;14.7.1.1},
we have
\[
 p^0_{13+}L_{\ast}\bigl(t(f+s),D_0\bigr)
\simeq
 \nbigo_{X\times\proj^1_s}
 \bigl(!(f+s)_0\bigr)
 \bigl(\ast D_{X\times\proj^1_s}\bigr),
\quad\quad
 p^i_{13+}L_{\ast}\bigl(t(f+s),D_0\bigr)=0\,\,\,
 (i\neq 0).
\]
We also have the following
by Proposition \ref{prop;14.7.1.1}:
\[
 p^0_{12+}L_{\ast}\bigl(t(f+s),D_0\bigr)
\simeq
 L_{\ast}\bigl(tf,D_{X\times\proj^1_t}\bigr),
\quad\quad
 p^i_{12+}L_{\ast}\bigl(t(f+s),D_0\bigr)=0\,\,\,
 (i\neq 0).
\]
Then, the claim of the proposition follows.
\hfill\qed

\subsubsection{Complement for a non-resonant case}

We give a remark on a non-resonant case
to compare the resonant case above.
Take $\alpha\in\cnum\setminus\seisuu$.
We consider the line bundle
$L_{\alpha,Y}(tf):=\nbigo_{Y}(\ast D_Y^{(1)})\,e$,
with the flat connection $\nabla$
determined by 
$\nabla e=e(d(tf)+\alpha dt/t)$,
where $e$ is a global frame.
Similarly, we consider the line bundle
$L_{-\alpha,X}:=
 \nbigo_X\bigl(\ast (D\cup (f)_0)\bigr)\,v$
with the connection $\nabla$
determined by $\nabla v=v\,(-\alpha df/f)$.
The following is 
essentially contained in Theorem {\rm 1.5}
of {\rm\cite{Adolphson-Sperber}}.
\begin{prop}
We have a natural isomorphism
$\pi^0_{+}L_{\alpha,Y}(tf)\simeq
 L_{-\alpha,X}$.
\end{prop}
\pf
We use the notation in {\rm\S\ref{subsection;14.7.2.10}}.
Indeed,
$\pi_{+}^0L_{\alpha,Y}(tf)$
is represented by the complex
\[
 \nbigo_X(\ast D)[t,t^{-1}]\,e
\stackrel{a}{\lrarr}
 \nbigo_X(\ast D)[t,t^{-1}]\,e.
\]
Here, the second term sits in the degree $0$,
and $a(g)=(t\del_t+\alpha+tf)g$.
The action of a holomorphic vector field $V$
on $t^j e$ is given as
$V(t^je)=t^{j+1}(Vf)e$.
We define an $\nbigo_X(\ast D)$-homomorphism
$\Phi:\nbigo_X(\ast D)[t,t^{-1}]\,e
\lrarr
 \nbigo_X(\ast D)\,v$
by
$\Phi(t^je)=
 \Gamma(-\alpha+1)\,
 \Gamma(-\alpha-j+1)^{-1}
 f^{-j}v$.
Then, we can check that
$\Phi$ is compatible with the connections,
and $\Phi\circ a=0$.
It induces a morphism of $\nbigd$-modules
$\pi^0_{+}L_{\alpha,Y}(tf)
\lrarr L_{-\alpha,X}$.
We can check that the induced morphism
is an isomorphism.
\hfill\qed

\subsection{Specialization}
\label{subsection;14.7.4.1}

\subsubsection{Statement}

Let $X$ be a complex manifold
with a simple normal crossing hypersurface $D$.
Let us consider meromorphic functions $f$ and $g$
on $(X,D)$.
We set $X^{(1)}:=X\times\cnum_{\tau}$
and $D^{(1)}:=D\times\cnum_{\tau}$.
We have the meromorphic function
$\tau f+g$ on $(X^{(1)},D^{(1)})$
and the associated $\nbigd$-modules
$\nbigm_{\star,f,g}:=L_{\star}(\tau f+g,D^{(1)})$ on $X^{(1)}$
for $\star=\ast,!$.
Let $K_{\star,f,g}$ and $C_{\star,f,g}$
denote the kernel and the cokernel
of 
$\nbigm_{\star,f,g}(!\tau)\lrarr
 \nbigm_{\star,f,g}(\ast \tau)$
for $\star=\ast,!$.
Let $\iota_0:X\lrarr X^{(1)}$ be given by
$\iota_0(Q)=(Q,0)$.
We shall prove the following proposition
in \S\ref{subsection;14.12.9.30}--\S\ref{subsection;14.12.9.31}.

\begin{prop}
\label{prop;14.6.28.10}
If $|(f)_0|\cap |(f)_{\infty}|=\emptyset$,
we have the following:
\[
C_{\ast,f,g}\simeq 
 \iota_{0+}L_{\ast}(g,D),
\quad
K_{\ast,f,g}\simeq
 \iota_{0+}L_{\ast}(g,D)(!(f)_{\infty}).
\]
\[
 K_{!,f,g}\simeq
 \iota_{0+}L_!(g,D),
\quad
 C_{!,f,g}\simeq
 \iota_{0+}L_!(g,D)(\ast (f)_{\infty}).
\]
\end{prop}

\subsubsection{The case $g=0$ and $D=|(f)_{\infty}|$}
\label{subsection;14.12.9.30}

First, we consider the case $g=0$
and $D=|(f)_{\infty}|$.
In this case, $\tau f$ is non-degenerate along  $D^{(1)}$.
Hence, we have
$\nbigm_{\ast,f,0}=\nbigm_{!,f,0}$
which we denote by $\nbigm_f$.
We also have
$K_{\ast,f,0}\simeq K_{!,f,0}$
and 
$C_{\ast,f,0}\simeq C_{!,f,0}$.
We have a global section $v$
of $\nbigm_f$ such that
$\nbigm_f=
\nbigo_{X^{(1)}}(\ast D^{(1)})v$
with $\nabla v=v\,d(\tau f)$.

We described the $V$-filtration 
$U_{\bullet}(\nbigm_f)$ of $\nbigm_f$
along $\tau=0$ in \cite{Mochizuki-K-complex},
which we recall here.
We use the convention that
$\tau\del_{\tau}+\alpha$
is locally nilpotent on
$U_{\alpha}\big/U_{<\alpha}$.
We describe the filtration
locally around any point of
$(Q,0)\in D\times\{0\}$.
Suppose $Q\not\in |(f)_{\infty}|$.
We have
$\Gr^U_j(\nbigm_f(\ast \tau))=0$
unless $j\in\seisuu$,
and 
$U_j\bigl( \nbigm_f(\ast\tau)\bigr)
=\tau^{-j}\nbigm_f$
for $j\in\seisuu$.
Suppose $Q\in |(f)_{\infty}|$.
We take a holomorphic coordinate system
$(z_1,\ldots,z_n)$ around $Q$
such that $D=\bigcup_{i=1}^{\ell}\{z_i=0\}$
and $f=z^{-\veck}$
for some $\veck\in\seisuu_{> 0}^{\ell}$.
For $0\leq\alpha<1$,
we set
$\vecp=[\alpha\veck]:=
 \bigl(
 [\alpha k_1],\ldots,[\alpha k_{\ell}]
 \bigr)$,
where $[a]:=\max\{n\in\seisuu\,|\,n\leq a\}$.
Let $\vecdelta=(1,\ldots,1)\in\seisuu_{\geq 0}^{\ell}$.
Let $\pi:X^{(1)}\lrarr X$ denote the projection.
We may naturally regard
$\pi^{\ast}\nbigd_X$
as a sheaf of subalgebras in $\nbigd_{X^{(1)}}$.
Locally around $(Q,0)$,
we have
\[
 U_{\alpha}\bigl(\nbigm_f(\ast \tau)\bigr)
=\pi^{\ast}\nbigd_X\bigl(
 \nbigo_{X^{(1)}} x^{-\vecdelta-\vecp}v
 \bigr)
=\pi^{\ast}\nbigd_X\Bigl(
 \sum_{j=0}^{\infty}
 \nbigo_{X^{(1)}}x^{-\vecdelta-\vecp}(\tau f)^j\,v
 \Bigr)
\]
in $\nbigo_{X^{(1)}}(\ast D^{(1)})v$.
We have
$U_1\bigl(
 \nbigm(\ast \tau)
 \bigr)
=\tau^{-1}U_0\nbigm$.
We obtain
\[
  U_{<0}\nbigm+\tau\del_{\tau}U_0\nbigm
=\pi^{\ast}\nbigd_X
\Bigl(
 \sum_{j=1}^{\infty}
 \nbigo_{X^{(1)}}x^{-\vecdelta}(\tau f)^j\,v
 \Bigr).
\]
Hence, 
$U_0\nbigm\big/
 (U_{<0}\nbigm+\tau\del_{\tau}U_0\nbigm)
\simeq
 \iota_{0\ast}\Bigl(
 \nbigd_X\bigl(
 \nbigo_X(D)
 \bigr)
\Bigr)
\simeq
 \iota_{0\ast}\nbigo_X(\ast D)$.
It implies  
$C_{\star,f,0}\simeq
 \iota_{0+}\nbigo_X(\ast D)$
if $D=|(f)_{\infty}|$.
By using the duality,
we obtain
$K_{\star,f,0}\simeq\iota_{0+}\nbigo_X(! D)$
if $D=|(f)_{\infty}|$.

\subsubsection{The case $|(g)_0|\cap|(g)_{\infty}|=\emptyset$}

Let us consider the case that
$|(g)_0|\cap|(g)_{\infty}|=\emptyset$.
We define $\nbigm_0:=L(\tau f)$.
Set $\nbigm:=\nbigm_{\ast,f,g}$.
We put $D_0:=|(f)_{\infty}|\cup|(g)_{\infty}|$.
We have the hypersurface $D_1\subset D$ such that
$D_0\cup D_1=D$
and $\codim_X (D_0\cap D_1)\geq 2$.
We set $D_2:=|(g)_{\infty}|\cup D_1$.
We have
$\nbigm\simeq
 \nbigm_0\otimes L_{\ast}(g,D_2^{(1)})$.
We naturally have
$\nbigm(\ast\tau)
\simeq
 \nbigm_0(\ast\tau)\otimes L_{\ast}(g,D_2^{(1)})$.
We shall observe that
$\nbigm(!\tau)
\simeq
 \nbigm_0(!\tau)\otimes L_{\ast}(g,D_2^{(1)})$.

We have the $V$-filtration 
$U_{\bullet}(\nbigm_0(\ast \tau))$
of $\nbigm_0(\ast\tau)$.
Set
$U_{\alpha}(\nbigm(\ast\tau))
:=U_{\alpha}(\nbigm_0(\ast\tau))\otimes L_{\ast}(g,D_2^{(1)})$
for any $\alpha\in\real$.

\begin{lem}
$U_{\bullet}(\nbigm(\ast \tau))$ is 
the $V$-filtration of $\nbigm(\ast\tau)$.
\end{lem}
\pf
By the construction,
$\tau\del_{\tau}+\alpha$
is locally nilpotent on
$U_{\alpha}\nbigm(\ast\tau)
 \big/
 U_{<\alpha}\nbigm(\ast\tau)$.
Let us prove that
$U_{\alpha}\nbigm$
is $V\nbigd_{X^{(1)}}$-coherent.
If $Q\not\in |(f)_{\infty}|$,
the claim is clear.
Let us consider the case
$Q \in |(f)_{\infty}|$.
We take a holomorphic coordinate system
$(x_1,\ldots,x_n)$ around $Q$
such that
$D=\bigcup_{i=1}^{\ell}\{x_i=0\}$
and 
$g=x^{-\veca}$
and
$f=f_1 x^{-\vecb}$,
where $f_1$ is nowhere vanishing.
Here $\veca,\vecb\in\seisuu^{\ell}_{\geq 0}$.
Put $\vecp:=[\alpha\vecb]$.
Let $\vecdelta=(1,\ldots,1)\in\seisuu^{\ell}$.

We use the identifications
$\nbigm_0=
 \nbigo_X^{(1)}\bigl(\ast (f)^{(1)}_{\infty}\bigr)\,v$
with $\nabla v=v\,d(\tau f)$,
and 
$L_{\ast}(g,D^{(1)}_2)=
 \nbigo_{X^{(1)}}\bigl(\ast D^{(1)}_2\bigr)\,e$
with $\nabla e=e\,dg$.
It is enough to prove that
$U_{\alpha}(\nbigm(\ast\tau))
\otimes L_{\ast}(g,D^{(1)}_2)$
is generated by $x^{-\vecdelta-\vecp}v\otimes e$
over $V\nbigd_{X^{(1)}}$.

Set $\supp(\veca):=\{i\,|\,a_i\neq 0\}$
and $\supp(\vecb):=\{i\,|\,b_i\neq 0\}$.
For $i\in \supp(\veca)\cup\supp(\vecb)$,
and for any $\vecq\in\seisuu_{\geq 0}^{\ell}$,
we have
\[
 \del_ix_i(x^{-\vecdelta-\vecp-\vecq}v\otimes e)
=x^{-\vecdelta-\vecp-\vecq}
 \Bigl(
 -p_i-q_i
-a_ix^{-\veca}
-b_i\tau f_1 x^{-\vecb}
+\tau x_ix^{-\vecb}\del_if_1
 \Bigr)\,v\otimes e.
\]
Because
$\tau\del_{\tau}(v\otimes e)
=f_1 \tau x^{-\vecb}(v\otimes e)$,
we have
\begin{equation}
\label{eq;14.12.9.40}
 \del_ix_i(x^{-\vecdelta-\vecp-\vecq}v\otimes e)
+(b_i-\tau f_1^{-1}x_i\del_if_1)
 \tau\del_{\tau}(x^{-\vecdelta-\vecp-\vecq}v\otimes e)
=x^{-\vecdelta-\vecp-\vecq}
(-p_i-q_i-a_ix^{-\veca})v\otimes e.
\end{equation}
By using (\ref{eq;14.12.9.40}),
we obtain that
$x^{-\vecdelta-\vecp}v\otimes 
 \nbigo_{X^{(1)}}(\ast D^{(1)}_2)e
\subset
 V\nbigd_{X^{(1)}}(v\otimes e)$.
Let $\vecdelta_1:=(1,\ldots,1)\in\seisuu^{\supp(\vecb)}$.
For any $Q\in \pi^{\ast}\nbigd_X$
and for any $h\in \nbigo_{X^{(1)}}(\ast D^{(1)}_2)$,
we have
$\del_i\bigl(
 Q(x^{-\vecdelta_1}v)
\otimes
 he
 \bigr)
=\del_iQ(x^{-\vecdelta_1}v)
 \otimes he
+Q(x^{-\vecdelta_1}v)
\otimes \del_i(he)$.
Hence,
we can easily deduce that
$U_{\alpha}(\nbigm)
\subset
 V\nbigd_{X^{(1)}}\,x^{-\vecdelta-\vecp}(v\otimes e)$.
Thus, we obtain the lemma.
\hfill\qed

\vspace{.1in}

We set 
$U_{\alpha}\bigl(
 \nbigm_0(!\tau)\otimes L_{\ast}(g,D^{(1)}_2)
 \bigr):=
 U_{\alpha}\bigl(\nbigm_0(!\tau)\bigr)
 \otimes L_{\ast}(g,D^{(1)}_2)$
for any $\alpha\in\real$.
\begin{lem}
$U_{\bullet}\bigl(\nbigm_0(!\tau)\otimes L_{\ast}(g,D^{(1)}_2)\bigr)$
is the $V$-filtration of
$\nbigm_0(!\tau)\otimes L_{\ast}(g,D^{(1)}_2)$.
\end{lem}
\pf
For $\alpha<1$,
we have
$U_{\alpha}\bigl(
 \nbigm_0(!\tau)\otimes L_{\ast}(g,D^{(1)}_2)
 \bigr)
=U_{\alpha}(\nbigm_0)\otimes L_{\ast}(g,D^{(1)}_2)$,
which is coherent over $V\nbigd_{X^{(1)}}$.
For $\alpha\geq 1$,
we have
$U_{\alpha}\bigl(\nbigm_0(!\tau)\otimes 
 L_{\ast}\bigl(g,D^{(1)}_2)\bigr)
=\sum_{\beta+n\leq \alpha,\,\,\beta<1}
 \del_{\tau}^n
 U_{\beta}\bigl(
 \nbigm_0(!\tau)\otimes L_{\ast}(g,D^{(1)}_2)\bigr)$.
We obtain that
$U_{\alpha}(\nbigm_0(!\tau)
 \otimes L_{\ast}\bigl(g,D^{(1)}_2)\bigr)$
are coherent over $V\nbigd_{X^{(1)}}$
for any $\alpha$.
We have
$\tau\del\tau+\alpha$
are nilpotent on $U_{\alpha}/U_{<\alpha}$.
Hence, 
$U_{\bullet}\bigl(\nbigm_0(!\tau)\otimes L_{\ast}(g,D^{(1)}_2)\bigr)$
is the $V$-filtration of 
$\nbigm_0(!\tau)\otimes L_{\ast}\bigl(g,D^{(1)}_2\bigr)$.
\hfill\qed

\vspace{.1in}
Because the induced morphism
\begin{multline}
 \del_{\tau}:
 U_{0}\bigl(\nbigm_0(!\tau)\otimes L_{\ast}(g,D^{(1)}_2)\bigr)
 \big/
U_{<0}\bigl(\nbigm_0(!\tau)\otimes L_{\ast}(g,D^{(1)}_2)\bigr)
 \\
\lrarr
  U_{1}\bigl(\nbigm_0(!\tau)\otimes L_{\ast}(g,D^{(1)}_2)\bigr)
 \big/
 U_{<1}\bigl(\nbigm_0(!\tau)\otimes L_{\ast}(g,D^{(1)}_2)\bigr)
\end{multline}
is an isomorphism,
we have
$\nbigm_0(!\tau)\otimes L_{\ast}(g,D^{(1)}_2)
\simeq
 \nbigm(!\tau)$.
We obtain
$C_{\ast,f,g}\simeq
\iota_{0+}\Bigl(
 \nbigo_X(\ast (f)_{\infty})
 \otimes L_{\ast}(g,D_2)
\Bigr)
\simeq 
\iota_{0+}L_{\ast}(g,D)$.
We also obtain
$K_{\ast,f,g}\simeq
 \iota_{0+}\Bigl(
 \nbigo_X(!(f)_{\infty})\otimes L_{\ast}(g,D_2)
 \Bigr)$.
Under the assumption
$|(g)_0|\cap|(g)_{\infty}|=\emptyset$,
it is naturally isomorphic to 
$\iota_{0+}L_{\ast}(g,D_2)(!(f)_{\infty})$.

\vspace{.1in}
Thus, we obtain the claims
for $K_{\ast,f,g}$
and $C_{\ast,f,g}$
in the case
$|(g)_{\infty}|\cap|(g)_0|=\emptyset$.
By using the duality,
we also obtain the claims for
$C_{!,f,g}$
and $K_{!,f,g}$ in this case. 

\subsubsection{The general case}
\label{subsection;14.12.9.31}

Let us consider the general case.
We take a projective birational morphism of complex manifolds
$G:X'\lrarr X$
such that 
(i) $D':=G^{-1}(D)$ is a simple normal crossing hypersurface,
(ii) $X'\setminus D'\simeq X\setminus D$,
(iii) $|(g')_0|\cap|(g')_{\infty}|=\emptyset$,
where $f':=G^{\ast}(f)$.
We set $g':=G^{\ast}(g)$.
We have
$(f')_{0}=G^{\ast}((f)_0)$
and 
$(f')_{\infty}=G^{\ast}((f)_{\infty})$.
In particular,
we have
$|(f')_0|\cap|(f')_{\infty}|=\emptyset$.
We have
$C_{\ast,f',g'}\simeq
 L_{\ast}(g',D')$
and 
$K_{\ast,f',g'}\simeq
 L_{\ast}(g',D')(!(f')_{\infty})$.

The induced morphism
$X^{\prime(1)}\lrarr X^{(1)}$
is also denoted by $G$.
We have
$G^0_{+}(\nbigm_{\ast,f',g'})
\simeq
 \nbigm_{\ast,f,g}$
and 
$G^i_{+}(\nbigm_{\ast,f',g'})=0$
for $i\neq 0$.
We obtain
$G^0_{+}(\nbigm_{\ast,f',g'}(\star\tau))
\simeq
 \nbigm_{\ast,f,g}(\star\tau)$
and 
$G^i_{+}(\nbigm_{\ast,f',g'}(\star\tau))=0$
for $i\neq 0$.
We have
$G^0_{+}\bigl(
 L_{\ast}(g',D')
 \bigr)
\simeq
 L_{\ast}(g,D)$
and 
$G^i_+\bigl(
 L_{\ast}(g',D')
 \bigr)=0$ for $i\neq 0$.
We also have
$G^0_+\bigl(
 L_{\ast}(g',D')\bigl(!(f')_{\infty}\bigr)
 \bigr)
\simeq
 L_{\ast}(g,D)\bigl((f)_{\infty}\bigr)
 $
and 
$G^i_+\bigl(
 L_{\ast}(g',D')\bigl(!(f')_{\infty}\bigr)
=0$ for $i\neq 0$.
Let $I$ be the image of
$\nbigm'_{\ast,f',g'}(!\tau)
\lrarr
 \nbigm'_{\ast,f',g'}(\ast\tau)$.
Because
$G^i_+\bigl(
 \nbigm_{\ast,f,g}(\ast\tau)
 \bigr)
=G^i_{+}\bigl(
 C_{\ast,f,g}
 \bigr)=0$ for $i\neq 0$,
we obtain that
$G^i_+(I)=0$ for $i\neq 0,-1$.
Because
$G^i_+\bigl(
 \nbigm_{\ast,f,g}(!\tau)
 \bigr)
=G^i_{+}\bigl(
 K_{\ast,f,g}
 \bigr)=0$ for $i\neq 0$,
we obtain that
$G^i_+(I)=0$ for $i\neq 0,1$.
Hence, we obtain 
$G^i_+(I)=0$ for $i\neq 0$.
It implies that
$G^0_+(C_{\ast,f',g'})
\simeq
 C_{\ast,f,g}$
and 
$G^0_+(K_{\ast,f',g'})
\simeq
 K_{\ast,f,g}$.
Thus, we obtain the claims for
$C_{\ast,f,g}$ and $K_{\ast,f,g}$.
By using the duality,
we obtain the claims for
$C_{!,f,g}$
and $K_{!,f,g}$.
Thus, the proof of Proposition \ref{prop;14.6.28.10}
is finished.
\hfill\qed

\subsection{Nearby cycle functor
and Push-forward}

\subsubsection{Beilinson's functors and variants}
\label{subsection;14.12.12.31}

Let $h$ be a meromorphic function on a complex manifold $Y$.
Let $M$ be a holonomic $\nbigd_Y$-module.
Suppose $|(h)_0|\cap|(h)_{\infty}|=\emptyset$.
We recall the functors of Beilinson \cite{beilinson2}:
\[
 \psi^{(a)}_h(M):=
 \varprojlim_{b}
 \Cok\bigl(
  \Pi^{a,b}_hM(!(h)_0)
\lrarr
 \Pi^{a,b}_hM(\ast(h)_0)
 \bigr)
\]
\[
 \Xi^{(a)}_h(M):=
 \varprojlim_{b}
 \Cok\bigl(
  \Pi^{a-1,b}_hM(!(h)_0)
\lrarr
 \Pi^{a,b}_hM(\ast(h)_0)
 \bigr).
\]
(\S\ref{subsection;15.10.20.1} for $\Pi^{a,b}_hM$.)
On any relatively compact subset $K$ in $Y$,
if $b$ is sufficiently large,
the cokernel of 
$\Pi^{a,b}_hM(!(h)_0)
\lrarr
 \Pi^{a,b}_hM(\ast(h)_0)$
is isomorphic to $\psi^{(a)}_h(M)$
and 
the  kernel is isomorphic to
$\psi^{(b)}_h(M)$.
On $K$, if $N$ is sufficiently large,
$\Xi^{(a)}_h(M)$ is isomorphic to
the cokernel of
$\Pi^{a+1,a+N}_hM(!(h)_0)
\lrarr
 \Pi^{a,a+N}_hM(\ast(h)_0)$
and the kernel of 
$\Pi^{a-N,a+1}_hM(!(h)_0)
\lrarr
 \Pi^{a-N,a}_hM(\ast(h)_0)$.
(See the argument in
\cite[\S4.1]{Mochizuki-MTM}, for example.)

\vspace{.1in}

In the case $|(h)_0|\cap|(h)_{\infty}|\neq\emptyset$
and $M=\nbigo_Y(\ast D)$
for a normal crossing hypersurface $D\subset Y$,
we also consider a variant.
For simplicity,
we assume that
(i) $(h)_{\infty}$ is smooth and reduced,
(ii) $D=|(h)_0|$,
(iii) $D\cup|(h)_{\infty}|$ is normal crossing.
We set $Z:=|(h)_{\infty}|$.

\begin{lem}
We have a canonical morphism
$\Pi^{a+1,b+1}_{h}\nbigo_Y
\lrarr
 \Pi^{a,b}_{h}\nbigo_Y(!Z)$
such that 
the composite with 
$\Pi^{a,b}_{h}\nbigo_Y(!Z)
\lrarr
 \Pi^{a,b}_{h}\nbigo_Y$
is the canonical one.
\end{lem}
\pf
Let $\nbigk$ and $\nbigc$ denote
the kernel and the cokernel
of the morphism
$\Pi^{a,b+1}_{h}\nbigo_Y(!Z)
\lrarr
 \Pi^{a,b+1}_{h}\nbigo_Y$.
We set $D_Z:=D\cap Z$.
Let $\iota:Z\lrarr Y$ denote the inclusion.
By a direct computation,
we can check
$\nbigc\simeq \iota_+\nbigo_Z(\ast D_Z)s^a$
and 
$\nbigk\simeq \iota_+\nbigo_Z(\ast D_Z)s^b$.
We have the induced monomorphism
$\Pi^{a,b+1}_{h}\nbigo_X(!Z)/\nbigk
\lrarr
 \Pi^{a,b+1}_{h}\nbigo_X$.

Let us consider the canonical morphism
$\Pi_{h}^{a+1,b+1}\nbigo_Y
\lrarr
 \Pi_{h}^{a,b+1}\nbigo_Y$.
The composite with 
$\Pi_{h}^{a,b+1}\nbigo_Y\lrarr \nbigc$
is $0$.
Hence, we have an induced morphism
$\Pi_{h}^{a+1,b+1}\nbigo_Y
\lrarr
 \Pi_{h}^{a,b+1}\nbigo_Y(!Z)/\nbigk$.
We have the following natural commutative diagram:
\[
 \begin{CD}
 \Pi^{a,b+1}_h\nbigo_Y(!Z)/\nbigk
 @>>>
 \Pi^{a,b+1}_h\nbigo_Y
 \\
 @VVV @VVV \\
 \Pi^{a,b}_h\nbigo_Y(!Z)
 @>>>
 \Pi^{a,b}_h\nbigo_Y
 \end{CD}
\]
Hence, we obtain
$\Pi_{h}^{a+1,b+1}\nbigo_Y
\lrarr
 \Pi_{h}^{a,b}\nbigo_Y(!Z)$
such that 
the composite with
$\Pi_{h}^{a,b}\nbigo_Y(!Z)
\lrarr
 \Pi_{h}^{a,b}\nbigo_Y$
is the canonical one.
\hfill\qed

\vspace{.1in}
Let us consider the following induced morphism
\begin{equation}
\label{eq;14.10.15.4}
 \rho_N:
  \Pi^{a+1,a+N+1}_h\nbigo_Y(\ast (h)_{\infty})(!(h)_0)
\lrarr
 \Pi^{a,a+N}_h\nbigo_Y(!(h)_{\infty})(\ast(h)_0).
\end{equation}
We set
$\Xi^{\prime(a)}_h\bigl(
 \nbigo_X[\ast D]
 \bigr):=
 \varprojlim_N
 \Cok(\rho_N)$.

\begin{lem}
\label{lem;15.10.20.2}
For any relatively compact subset
$K\subset X$,
there exists $N_0$ such that 
$\Cok(\rho_{N_2})_{|K}
\lrarr
 \Cok(\rho_{N_1})_{|K}$
is an isomorphism 
if $N_0\leq N_1\leq N_2$.
\end{lem}
\pf
As mentioned above,
the claim is well known outside of $Z$.
Hence,
the induced morphism
\[
 \Cok(\rho_{N_2})_{|K}(\ast (h)_{\infty})
\lrarr
 \Cok(\rho_{N_1})_{|K}(\ast (h)_{\infty})
\]
is an isomorphism.
Take $Q\in Z\cap K$,
and take a holomorphic local coordinate neighbourhood
$(Y_Q;x_1,\ldots,x_n)$ of $Y$ around $Q$
such that 
$h_{|Y_Q}=x_n^{-1}\prod_{i=1}^{\ell}x_i^{k_i}$,
where $k_i\in\seisuu_{>0}$.
It is enough to prove that 
$\phi_{x_n}^{(0)}
 \Cok(\rho_{N_2})_{|K}
\lrarr
 \phi_{x_n}^{(0)}
 \Cok(\rho_{N_1})_{|K}$
is an isomorphism.
Set $h_1:=\prod_{i=1}^{\ell}x_i^{k_i}$
on $Z_Q:=Y_Q\cap Z$.
Let $\iota_1:Z_Q\lrarr Y_Q$ be the inclusion.
We have the following commutative diagram:
\[
 \begin{CD}
 \phi_{x_n}^{(0)}
 \Pi_h^{a+1,a+N+1}
 \nbigo_Y(\ast (h)_{\infty})(!(h)_0)
 @>>>
 \phi_{x_n}^{(0)}
 \Pi_h^{a,a+N}
 \nbigo_Y(! (h)_{\infty})(\ast(h)_0)
 \\
 @V{\simeq}VV @V{\simeq}VV \\
  \psi_{x_n}^{(0)}
 \Pi_h^{a+1,a+N+1}
 \nbigo_Y(\ast (h)_{\infty})(!(h)_0)
 @.
 \psi_{x_n}^{(1)}
 \Pi_h^{a,a+N}
 \nbigo_Y(! (h)_{\infty})(\ast(h)_0)
\\
 @V{b_1}V{\simeq}V @V{b_2}V{\simeq}V \\
 \iota_{1+}\Pi^{a+1,a+N+1}_{h_1}\nbigo_{Z_Q}(!h_1)
 @>{c_1}>>
  \iota_{1+}\Pi^{a,a+N}_{h_1}\nbigo_{Z_Q}
 \end{CD}
\]
Here,
$b_i$ are given as in 
\cite[Proposition 4.3.1]{Mochizuki-MTM}.
We shall prove that
 $c_1$ is induced by the multiplication of $s^{-1}$ 
Lemma \ref{lem;15.11.11.1} below.
Hence, 
$\phi_{x_n}^{(0)}
 \Cok(\rho_{N_2})$
is identified with 
$\iota_{1+}\Cok\bigl(
 \Pi^{a,a+N}_{h_1}\nbigo_{Z_Q}(!h_1)
\lrarr
  \Pi^{a,a+N}_{h_1}\nbigo_{Z_Q}
 \bigr)$.
Then, the claim is reduced to the above standard case.
\hfill\qed

\begin{lem}
\label{lem;15.11.11.1}
$c_1$ is induced by the multiplication of $s^{-1}$.
\end{lem}
\pf
It is enough to consider the case $D=\emptyset$.
We may assume that $x_n=h_n^{-1}$.
Taking the limit of $\rho_N$ for $N\to\infty$,
we obtain
$\Pi^{a+1,\infty}_{x_n^{-1}}\nbigo_Y
\lrarr
 \Pi^{a,\infty}_{x_n^{-1}}\nbigo_Y[!x_n]$
such that the composite with 
the canonical morphism
$\kappa_1:
 \Pi^{a,\infty}_{x_n^{-1}}\nbigo_Y[!x_n]
\lrarr
 \Pi^{a,\infty}_{x_n^{-1}}\nbigo_Y$
is equal to the canonical morphism
$\kappa_2:
 \Pi^{a+1,\infty}_{x_n^{-1}}\nbigo_Y
\lrarr
 \Pi^{a,\infty}_{x_n^{-1}}\nbigo_Y$.

We have the following commutative diagram:
\[
\begin{CD}
 \phi_{x_n}^{(0)}
 \Pi_{x_n^{-1}}^{a+1,\infty}\nbigo_Y
 @>>>
 \phi_{x_n}^{(0)}
 \Pi_{x_n^{-1}}^{a,\infty}\nbigo_Y
 \\
 @V{\simeq}VV @V{\simeq}VV \\
 \iota_{1+}
 \Pi^{a+1,\infty}_{1}\nbigo_Z
@>{i_1}>>
 \iota_{1+}
 \Pi^{a,\infty}_{1}\nbigo_Z
\end{CD}
\]
Here, $i_1$ is the inclusion.

We have the following commutative diagram:
\[
 \begin{CD}
 \phi_{x_n}^{(0)}
 \Pi_{x_n^{-1}}^{a,\infty}\nbigo_Y[!x_n]
 @>>>
 \phi_{x_n}^{(0)}
 \Pi_{x_n^{-1}}^{a,\infty}\nbigo_Y
 \\
 @V{\simeq}VV @V{\simeq}VV \\
 \iota_{1+}
 \Pi^{a,\infty}_{1}\nbigo_Z
@>{i_2}>>
 \iota_{1+}
 \Pi^{a,\infty}_{1}\nbigo_Z
 \end{CD}
\]
Here, $i_2$ is induced by $-x_n\del_{x_n}$,
and it is equal to the multiplication of $s$.
Because $i_2\circ c_1=i_1$,
we obtain that $c_1$ is 
induced by the multiplication of $s^{-1}$.
\hfill\qed

\begin{rem}
\label{rem;15.11.10.400}
In the proof of Lemma {\rm\ref{lem;15.10.20.2}},
we also observed that
$\phi_{x_n}^{(0)}
 \Xi^{\prime(a)}_h(\nbigo_X[\ast D])_{|Y_Q}
\simeq
 \iota_{1+}\psi_{h_1}^{(a)}\nbigo_{Z_Q}$
around any point $Q\in Z\cap D$.
\hfill\qed
\end{rem}

\subsubsection{The push-forward of the nearby cycle sheaf}
\label{subsection;14.12.9.100}

Let $X$ be a complex manifold
with a normal crossing hypersurface $D$.
Let $f$ and $g$ be meromorphic functions on $(X,D)$.
We set $Y:=\proj^1_{\tau}\times X$
and $D^{(1)}:=(\proj^1_{\tau}\times D)\cup(\{\infty\}\times X)$.
Let $\pi:Y\lrarr X$ be the projection.
We consider $F:=\tau f+g$ on $(Y,D^{(1)})$.
We suppose the following.
\begin{itemize}
\item
 $F$ is pure on $\{\tau\neq 0\}\times X$.
\item
$g$ is pure on $X\setminus |(f)_{\infty}|$.
\item
$D=|(f)_{\infty}|\cup|(g)_{\infty}|$,
 $|(f)_0|\subset D$ and $|(f)_0|\cap |(f)_{\infty}|=\emptyset$.
\end{itemize}
Note that
we have
$\pi_+^0\psi_{\tau}^{(a)}L_{!}(F)
=\pi_+^0\psi_{\tau}^{(a)}L_{\ast}(F)$
and
$\Xi^{(a)}_{f^{-1}}L_{!}(g)
\simeq
 \Xi^{(a)}_{f^{-1}}L_{\ast}(g)$
by the assumptions.

\begin{prop}
\label{prop;14.10.15.1}
We have natural isomorphisms
$\Lambda^{(a)}:\pi_{+}^0\psi_{\tau}^{(a)}L_{\ast}(F)
\simeq
 \Xi^{(a)}_{f^{-1}}L_{\ast}(g)$.
\end{prop}
\pf
By the assumption of the purity of $F$ on 
$\{\tau\neq 0\}\times X$,
we have
$\Pi^{a,N}_{\tau}L_{\ast}(F)(\star \tau)
=\Pi^{a,N}_{\tau}L_!(F)(\star \tau)$.
Let $\nbigi$ be the image of
$\Pi^{a,b}_{\tau}L_{!}(F)(!\tau)
\lrarr
\Pi^{a,b}_{\tau}L_{\ast }(F)(\ast \tau)$.
The support of the kernel $K^{a,b}$
and the cokernel $C^{a,b}$
are contained in $\tau=0$.
Hence, we have
$\pi_{+}^iK^{a,b}=0$
and 
$\pi_{+}^iC^{a,b}=0$
for $i\neq 0$.
By Proposition \ref{prop;14.10.14.21},
we have 
$\pi_{+}^i
 \Pi^{a,N}_{\tau}L_{\star_1}(F)(\star_2 \tau)=0$
for $i\neq 0$
and for $\star_1,\star_2\in\{\ast,!\}$.
We can easily deduce that
$\pi_{+}^i\nbigi=0$ 
unless $i=0$.
Hence, we have the following exact sequence
\[
 0\lrarr
 \pi_{+}^0
 K^{a,b}
\lrarr
 \pi_{+}^0
 \Pi^{a,b}_{\tau} L_{!}(F)(!\tau)
\lrarr
 \pi_{+}^0
 \Pi^{a,b}_{\tau} L_{\ast}(F)(\ast\tau)
\lrarr
  \pi_{+}^0
 C^{a,b}
\lrarr 0
\]
If $|a-b|$ is sufficiently large,
we have
$C^{a,b}\simeq
 \psi^{(a)}_{\tau}L_{\ast}(F)$.
Hence, we have
\[
 \pi_{+}^0
 \psi^{(a)}_{\tau}L_{\ast}(F)
\simeq
 \Cok\Bigl(
  \pi_{+}^0
 \Pi^{a,b}_{\tau} L_{!}(F)(!\tau)
\lrarr
 \pi_{+}^0
 \Pi^{a,b}_{\tau} L_{\ast}(F)(\ast\tau)
 \Bigr)
\]

By Proposition \ref{prop;14.10.14.21},
we have the following natural isomorphisms:
\[
 \Psi^{a,N}_{\star}:
 \pi_{+}^0
 \Pi^{a,N}_{\tau}L_{\star}(F)(\star \tau)
\simeq
 \Pi^{a,N}_{f^{-1}}L_{\star}(g)(\star D)
\]
By $|(g)_{\infty}|\cup|(f)_{\infty}|=D$
and the purity of $g$ on 
$X\setminus|(f)_{\infty}|$,
we have 
\[
 \Pi^{a,N}_{f^{-1}}L_{\ast}(g)(\ast D)
\simeq
 \Pi^{a,N}_{f^{-1}}L_{\ast}(g)(\ast (f)_{\infty}),
\quad\quad
 \Pi^{a,N}_{f^{-1}}L_{!}(g)(!D)
\simeq
 \Pi^{a,N}_{f^{-1}}L_{\ast}(g)(!(f)_{\infty}).
\]
By the commutativity of (\ref{eq;15.11.6.100}),
the image of the morphism 
$\pi_{+}^0
 \Pi^{a,N}_{\tau}L_{!}(F)(!\tau)
\lrarr 
 \pi_{+}^0
 \Pi^{a,N}_{\tau}L_{\ast}(F)(\ast\tau)$
is identified with the image of
\[
 \Pi^{a+1,N+1}_{f^{-1}}L_{\ast}(g)(!(f)_{\infty})
\lrarr
 \Pi^{a,N}_{f^{-1}}L_{\ast}(g)(\ast (f)_{\infty}).
\]
Hence, we obtain the desired isomorphism.
\hfill\qed

\vspace{.1in}

Let $\iota_0:X\lrarr \proj^1\times X$
be the inclusion given by
$\iota_0(P)=(0,P)$.
\begin{cor}
Under the assumption,
we have a natural isomorphism
$\psi_{\tau}^{(a)}L_{\ast}(F)
\simeq
 \iota_{0+}
 \Xi^{(a)}_{f^{-1}}L_{\ast}(g)$.
\hfill\qed
\end{cor}

We have the canonical nilpotent map $N$
on $\psi_{\tau}^{(a)}L_{\ast}(F)$.
The following proposition is clear
by the above description.
\begin{prop}
\label{prop;14.10.18.20}
Under the same assumption,
we have
$\Ker N\simeq \iota_{0\dagger}L_!(g,D)$
and 
$\Cok N\simeq \iota_{0\dagger}L_{\ast}(g,D)$.
The canonical morphism
$\Ker N\lrarr \Cok N$
is identified with
the canonical morphism
$L_!(g,D)\lrarr L_{\ast}(g,D)$.
\hfill\qed
\end{prop}

\begin{rem}
Note that 
$\Ker N$ and $\Cok N$
are isomorphic to
the kernel and the cokernel of
$L_{\star}(F)(!\tau)
\lrarr
 L_{\star}(F)(\ast \tau)$.
The isomorphisms in
Proposition {\rm \ref{prop;14.6.28.10}}
and Proposition {\rm\ref{prop;14.10.18.20}}
are consistent.
\hfill\qed
\end{rem}

\subsubsection{A variant}

Let us give a similar statement
in a slightly different situation.
We continue to use the notation in
\S\ref{subsection;14.12.9.100}.
Suppose that
(i) $g=0$,
(ii) $D=|(f)_{\infty}|$,
(iii) $(f)_0$ is non-singular and reduced, and
 $|(f)_0|\cup D$ is normal crossing.
Note that $F=\tau f$ is non-degenerate on 
$\{\tau\neq 0\}\times X$.
We set $Z:=|(f)_0|$,
and let $\iota_Z:Z\lrarr X$ denote the inclusion.

\begin{prop}
\label{prop;14.10.15.4}
We have a natural isomorphism
$\pi_{+}^0\psi^{(a)}_{\tau}L_{\ast}(F)
\simeq
\Xi^{\prime(a)}_{f^{-1}}\nbigo_X$.
\end{prop}
\pf
As in the proof of Proposition \ref{prop;14.10.15.1},
we have
\begin{equation}
 \label{eq;15.11.11.2}
 \pi_{+}^0\psi^{(a)}_{\tau}L_{\ast}(F)
\simeq
 \Cok\Bigl(
 \pi_{+}^0\Pi^{a,a+N}_{\tau}L_{!}(F)(!\tau)
\stackrel{\kappa_N}{\lrarr}
  \pi_{+}^0\Pi^{a,a+N}_{\tau}L_{\ast}(F)(\ast \tau)
 \Bigr).
\end{equation}
According to Proposition \ref{prop;14.10.14.21},
we have the following isomorphisms:
\begin{equation}
\label{eq;14.10.15.2}
 \pi_{+}^0\Pi^{a,a+N}_{\tau}L_{!}(F)(!\tau)
\simeq
 \Pi^{a,a+N}_{f^{-1}}\nbigo_X(\ast Z)(!D)
\end{equation}
\begin{equation}
\label{eq;14.10.15.3}
 \pi_{+}^0\Pi^{a,a+N}_{\tau}L_{\ast}(F)(\ast\tau)
\simeq
 \Pi^{a,a+N}_{f^{-1}}\nbigo_X(!Z)(\ast D)
\end{equation}
Then, we obtain the claim of Proposition \ref{prop;14.10.15.4}
from the following lemma.

\begin{lem}
\label{lem;14.10.15.5}
Under the isomorphisms
{\rm(\ref{eq;14.10.15.2})}
and {\rm(\ref{eq;14.10.15.3})},
the image of $\kappa_N$ is equal to
the image $\rho_N$ in {\rm(\ref{eq;14.10.15.4})}
if $N$ is sufficiently large.
\end{lem}
\pf
We identify 
$\Pi^{a,a+N}_{f^{-1}}\nbigo_X(!Z)(\ast D)$
and 
$\Pi^{a+1,a+N+1}_{f^{-1}}\nbigo_X(!Z)(\ast D)$
by the multiplication of $-s$.
We can compare the restriction of 
$\rho_N$ and $\kappa_N$
to $X\setminus \bigl(Z\cup D\bigr)$
by the commutativity of (\ref{eq;15.11.6.100}).
Take any $N_1\geq N$.
We have the morphisms
$\kappa_{N_1}$ and $\rho_{N_1}$.
The morphisms $\kappa_{N}$
and $\rho_{N}$
are induced by
$\kappa_{N_1}$ and $\rho_{N_1}$.
Because 
$(\kappa_{N_1}-\rho_{N_1})_{|X\setminus (D\cup Z)}=0$,
the difference
$\kappa_{N_1}-\rho_{N_1}$
factors through
$\iota_{Z+}\nbigo_Z(\ast D_Z)s^{a+N_1}$,
where $D_Z:=Z\cap D$.
Hence, 
we obtain that
$\kappa_{N}-\rho_{N}=0$.
Thus, we finish the proof of
Lemma \ref{lem;14.10.15.5}
and Proposition \ref{prop;14.10.15.4}.
\hfill\qed

\begin{cor}
If $Z\cap D=\emptyset$,
we have
$\pi_{+}^0\psi^{(0)}_{\tau}L_{\ast}(\tau f)
\simeq
 \Xi^{(0)}_{f^{-1}}\nbigo_X$.
\hfill\qed
\end{cor}

Take a holomorphic local coordinate
$(x_1,\ldots,x_n)$
such that $f=x_n^{-1}\prod_{i=1}^{\ell} x_i^{k_i}$.
Set $x^{\veck}:=\prod_{i=1}^{\ell}x_i^{k_i}=f\cdot x_n$.
\begin{lem}
\label{lem;15.11.11.40}
We have the following commutative diagram:
\[
 \begin{CD}
 \phi^{(0)}_{x_n}
 \pi_+^0
 \Pi^{a,a+N}_{\tau}L_!(F)(!\tau)
 @>{\phi^{(0)}_{x_n}(\kappa_N)}>>
 \phi^{(0)}_{x_n}
 \pi_+^0
 \Pi^{a,a+N}_{\tau}L_{\ast}(F)(\ast\tau)
 \\
 @VVV @VVV \\
 \iota_{Z+}
 \Pi^{a,a+N}_{x^{\veck}}\nbigo_Z[!D]
 @>{c}>>
 \iota_{Z+}
 \Pi^{a,a+N}_{x^{\veck}}\nbigo_Z
 \end{CD}
\]
Here, $c$ is induced by the multiplication of $-1$.
\end{lem}
\pf
As observed in Lemma \ref{lem;14.10.15.5},
$\kappa_N$ is identified with $\rho$
under the isomorphism
$\Pi^{a,a+N}_{f^{-1}}\nbigo_X(!Z)(\ast D)$
and 
$\Pi^{a+1,a+N+1}_{f^{-1}}\nbigo_X(!Z)(\ast D)$
given by the multiplication of $-s$.
Then, the claim of the lemma follows from
Lemma \ref{lem;15.11.11.1}.
\hfill\qed

\section{Mixed twistor $\nbigd$-modules
 associated to meromorphic functions}
\subsection{Mixed twistor $\nbigd$-modules}
\label{subsection;15.5.19.1}

We recall some operations for $\nbigr$-modules
and twistor $\nbigd$-modules.
See \cite{Mochizuki-wild},
\cite{Mochizuki-MTM}
and \cite{sabbah2}
for more details.

\subsubsection{$\nbigr$-modules}

Let $X$ be a complex manifold.
We set $\nbigx:=\cnum_{\lambda}\times X$.
Let $p_{\lambda}:\nbigx\lrarr X$ denote the projection.
Let $\nbigr_X\subset\nbigd_{\nbigx}$
be the sheaf of subalgebras
generated by 
$\lambda p_{\lambda}^{\ast}\Theta_X$
over $\nbigo_{\nbigx}$.
We set $d_X:=\dim X$.
The pull back $p_{\lambda}^{\ast}\Omega_X$
is also denoted by $\Omega_X$.
For any left $\nbigr_X$-module $\nbigm$,
$\lambda^{-d_X}\Omega_X\otimes_{\nbigo_{\nbigx}}
 \nbigm$
is naturally a right $\nbigr_X$-module,
by which the category of left $\nbigr_X$-modules
and the category of right $\nbigr_X$-modules
are equivalent.
In this paper, $\nbigr_X$-modules mean
left $\nbigr_X$-modules.
Let $D^b_c(\nbigr_X)$ denote the derived category of
cohomologically bounded coherent complexes of $\nbigr_X$-modules.

An $\nbigr_X$-module is equivalent to 
an $\nbigo_{\nbigx}$-module $\nbigm$
with a meromorphic relative flat connection
\[
\DD^f:
 \nbigm\lrarr
\Omega_{\nbigx/\cnum_{\lambda}}^1(\nbigx^0)\otimes\nbigm
\]
where $\nbigx^0:=\{0\}\times X$.
The operator $\DD:=\lambda\DD^f$
is also often used, and called a family of flat $\lambda$-connections.

The easiest example of $\nbigr_X$-module
is the line bundle $\nbigo_{\nbigx}$
with the meromorphic relative flat connection
$\DD^f$ determined by $\DD^f(1)=0$.
It is just denoted by $\nbigo_{\nbigx}$.

Let $\nbigm_i$ $(i=1,2)$ be $\nbigr_X$-modules.
Then, $\nbigm_1\otimes_{\nbigo_{\nbigx}}\nbigm_2$
and $\nbigm_1\oplus\nbigm_2$
are naturally $\nbigr_X$-modules.
The tensor product 
$\nbigm_1\otimes_{\nbigo_{\nbigx}}\nbigm_2$
is denoted just by $\nbigm_1\otimes\nbigm_2$,
if there is no risk of confusion.

We define the duality functor
$\DDD_X:D^b_c(\nbigr_X)\lrarr D^b_c(\nbigr_X)$
by
\[
 \DDD_X(\nbigm):=
 \nrhom_{\nbigr_X}\bigl(\nbigm,
 \nbigr_X\otimes\Omega_X^{-1}\bigr)[\dim X].
\]
(We will review more details in \S\ref{subsection;15.1.4.20}.)
Note that 
$\DDD_X\nbigo_{\nbigx}$
is naturally isomorphic to
$\lambda^{d_X}\nbigo_{\nbigx}
=\nbigo_{\nbigx}(-d_X\nbigx^0)$.
If $\nbigm$ is an $\nbigr_X$-module
underlying a mixed twistor $\nbigd$-module,
then $\nbigh^j\bigl( \DDD_X\nbigm\bigr)=0$
unless $j=0$.
In that case,
we identify $\DDD_X\nbigm$
and $\nbigh^0\DDD_X\nbigm$.

Let $f:X\lrarr Y$ be any morphism of complex manifolds.
We set
\[
 \nbigr_{Y\larr X}:=
 \lambda^{-d_X}\Omega_X\otimes_{f^{-1}\nbigo_{\nbigy}}
 f^{-1}\bigl(
 \nbigr_Y\otimes \lambda^{d_Y}\Omega_Y^{-1}
 \bigr).
\]
It is naturally 
a left $f^{-1}\nbigr_Y$-module
and a right $\nbigr_X$-module.
For any object $\nbigm$ in $D^b_c(\nbigr_X)$,
we set
\[
 f_{\dagger}(\nbigm):=
 Rf_{!}
 \bigl(
 \nbigr_{Y\larr X}
\otimes^L_{\nbigr_X}
 \nbigm\bigr)
\]
in the derived category of $\nbigr_Y$-modules.
(We shall review more details on the push-forward
in \S\ref{subsection;15.1.4.21}.)
The $j$-th cohomology sheaves of $f_{\dagger}(\nbigm)$
are denoted by $f^j_{\dagger}(\nbigm)$.
If $\nbigm$ is good relative to $f$,
i.e., for any compact subset $K\subset\nbigy$,
then there exists good filtration of $\nbigm_{|f^{-1}(K)}$.
Moreover, suppose that the support of $\nbigm$
is proper over $f$.
Then, 
$f_{\dagger}\nbigm$
is an object in $D^b_c(\nbigr_Y)$,
and we have a natural isomorphism
$f_{\dagger}\DDD_X\nbigm
\simeq
 \DDD_Y f_{\dagger}\nbigm$.

Let $j:\cnum_{\lambda}\lrarr\cnum_{\lambda}$
be given by $j(\lambda)=-\lambda$.
The induced morphism
$j\times\id_X:\nbigx\lrarr\nbigx$
is also denoted by $j$.
For any $\nbigr_X$-module $M$,
we naturally regard $j^{\ast}M$
as an $\nbigr_X$-module
by the natural isomorphism
$j^{\ast}\nbigr_X\simeq\nbigr_X$.

\vspace{.1in}

Let $H$ be a hypersurface in $X$.
We set $\nbigh:=\cnum_{\lambda}\times H$.
For any $\nbigr_X$-module $\nbigm$,
we set $\nbigm(\ast H):=
\nbigm\otimes\nbigo_{\nbigx}(\ast\nbigh)$.
In particular,
$\nbigr_X(\ast H):=
\nbigr_X\otimes\nbigo_{\nbigx}(\ast\nbigh)$.
We may naturally consider
$\nbigr_X(\ast H)$-modules.
For any $\nbigm\in D^b_c(\nbigr_X(\ast H))$,
we set
\[
 \DDD_{X(\ast H)}(\nbigm)
:=\nrhom_{\nbigr_X(\ast H)}\bigl(
 \nbigm,\nbigr_X(\ast H)\otimes\Omega_X^{-1}
 \bigr)[\dim X].
\]
Let $f:X\lrarr Y$ be a proper morphism of complex manifolds.
Let $H_Y$ be a hypersurface in $Y$,
and we set $H_X:=f^{-1}(H_Y)$.
For any $\nbigm\in D^b_c(\nbigr_X(\ast H_X))$,
we have
\[
 f_{\dagger}(\nbigm):=
 Rf_{!}\bigl(
 \nbigr_{Y\rarr X}
 \otimes_{\nbigr_X}^L\nbigm
 \bigr)
\simeq
 Rf_{!}\bigl(
 \nbigr_{Y\rarr X}(\ast H_X)
 \otimes_{\nbigr_X(\ast H_X)}^L\nbigm
 \bigr)
\]
in the derived category of $\nbigr_Y(\ast H_Y)$-modules.
If $f$ induces an isomorphism
$X\setminus H_X\simeq Y\setminus H_Y$,
we have
$f^j_{\dagger}\nbigm=0$ unless $j=0$.
We shall identify $f_{\dagger}\nbigm$
and $f^0_{\dagger}\nbigm$
in that case.

\subsubsection{$\nbigr$-triples}

We set $\vecS:=\bigl\{\lambda\in\cnum\,\big|\,|\lambda|=1\bigr\}$.
Let $\distribution_{\vecS\times X/\vecS}$
denote the sheaf of distributions on $\vecS\times X$
which are continuous with respect to $\vecS$.
(See \cite{sabbah2}.)
For local sections
$P\in\nbigr_{X|\vecS\times X}$
and $Q\in\distribution_{\vecS\times X/X}$,
the local section
$P\bullet Q:=PQ$ of $\distribution_{\vecS\times X/\vecS}$
is naturally defined.
Let $\sigma:\vecS\times X\lrarr \vecS\times X$
be given by
$\sigma(\lambda,x)=(-\lambda,x)$.
For local sections
$\sigma^{\ast}P
 \in
 \sigma^{\ast}\nbigr_{X|\vecS\times X}$
and $Q\in\distribution_{\vecS\times X/X}$,
the local section
$\sigma^{\ast}(P)\bullet Q:=
 \sigma^{\ast}(\overline{P})Q$
is defined.
Thus, the sheaf $\distribution_{\vecS\times X/X}$
is naturally a
$(\nbigr_{X|\vecS\times X},
 \sigma^{\ast}\nbigr_{X|\vecS\times X})$-module.

Let $\nbigm_i$ $(i=1,2)$ be $\nbigr_X$-modules.
A sesqui-linear pairing of
$\nbigm_1$ and $\nbigm_2$ is a 
$(\nbigr_{X|\vecS\times X},
 \sigma^{\ast}\nbigr_{X|\vecS\times X})$-homomorphism
$C:
 \nbigm_{1|\vecS\times X}
\otimes
\sigma^{\ast} \nbigm_{2|\vecS\times X}
\lrarr
 \distribution_{\vecS\times X/\vecS}$.
Such a tuple
$\nbigt=(\nbigm_1,\nbigm_2,C)$
is called an $\nbigr_X$-triple.

Let $\nbigt_i=(\nbigm_{i1},\nbigm_{i2},C_i)$
be $\nbigr_X$-triples.
A morphism  of $\nbigr_X$-triples
$\varphi:\nbigt_1\lrarr \nbigt_2$
is a pair of $\nbigr_X$-homomorphisms
$\varphi_1:\nbigm_{21}\lrarr\nbigm_{11}$
and 
$\varphi_2:\nbigm_{12}\lrarr\nbigm_{22}$
such that
$C(\varphi_1(m_2),\sigma^{\ast}m_1)
=C(m_2,\sigma^{\ast}\varphi_2(m_1))$.
For any morphism,
we set
$\Ker(\varphi):=\bigl(
 \Cok(\varphi_1),\Ker(\varphi_2),C_{\Ker\varphi}
 \bigr)$,
where
$C_{\Ker\varphi}$ denotes the naturally induced
sesqui-linear pairing.
Similarly, we set
\[
\Cok(\varphi):=
 \bigl(
 \Ker(\varphi_1),\Cok(\varphi_2),
 C_{\Cok(\varphi)}
 \bigr),
\quad
\Image(\varphi):=
 \bigl(
 \Image(\varphi_1),
 \Image(\varphi_2),
 C_{\Image(\varphi)}
 \bigr).
\]

For an $\nbigr_X$-triple
$\nbigt=(\nbigm_1,\nbigm_2,C)$,
we set
$\nbigt^{\ast}:=
 (\nbigm_2,\nbigm_1,C^{\ast})$,
where
$C^{\ast}(m_2,\sigma^{\ast}(m_1)):=
 \overline{\sigma^{\ast}C(m_1,\sigma^{\ast}m_2)}$.
For a morphism $\varphi=(\varphi_1,\varphi_2):
\nbigt_1\lrarr\nbigt_2$ of $\nbigr_X$-triples,
we set 
$\varphi^{\ast}=(\varphi_2,\varphi_1)$.

We set
$\nbigu_X(a,b):=
 \bigl(
 \lambda^{a}\nbigo_{\nbigx},
 \lambda^{b}\nbigo_{\nbigx},
 C_0
 \bigr)$,
where
$C_0(f,\sigma^{\ast}g):=f\overline{\sigma^{\ast}g}$.
For any $\nbigr$-triple
$\nbigt:=(\nbigm_1,\nbigm_2,C)$,
we define
$\nbigt\otimes\nbigu_X(a,b):=
\bigl(
 \lambda^a\nbigm_1,\lambda^b\nbigm_2,C
\bigr)$,
where $C$ is induced by the natural identification
$(\lambda^{c}\nbigm_{i})_{|\vecS\times X}
=(\nbigm_i)_{|\vecS\times X}$.
Particularly,
$\nbigu_X(-n,n)$ is denoted by
$\newTate(n)$,
called the $n$-th Tate object,
and 
$\nbigt\otimes\newTate(n)$
is called the $n$-th Tate twist of $\nbigt$.
We use the identification
$\newTate(n)\simeq\newTate(-n)^{\ast}$
given by the morphism
$((-1)^n,(-1)^n)$.

A morphism
$\nbigs:\nbigt\lrarr \nbigt^{\ast}\otimes\newTate(-n)$
such that
$\nbigs^{\ast}=(-1)^n\nbigs$
is called a sesqui-linear duality of weight $n$.
For example,
we have
$\nbigu_X(p,q)
\lrarr
 \nbigu_X(p,q)^{\ast}
 \otimes
 \newTate(-(p-q))$
given by $\bigl((-1)^p,(-1)^p\bigr)$.

Let $f:X\lrarr Y$ be any morphism of complex manifolds.
Let $\nbigt$ be an $\nbigr_X$-triple.
We have a naturally defined push-forward
$f_{\dagger}$
in the derived category of $\nbigr_Y$-triples.
See \cite{sabbah2}.
It induces a functor 
from the derived category of $\nbigr_X$-triples
to the derived category of $\nbigr_Y$-triples.
The $j$-th cohomology $\nbigr_Y$-triple
of $f_{\dagger}\nbigt$
is denoted by $f_{\dagger}^j\nbigt$.
For $\nbigt=(\nbigm_1,\nbigm_2,C)$,
the $\nbigr$-triple
$f_{\dagger}^j\nbigt$ consists of
$(f_{\dagger}^{-j}\nbigm_1,f_{\dagger}^j\nbigm_2,f^j_{\dagger}C)$.

\vspace{.1in}
Let $H$ be a hypersurface.
The sheaf  $\distribution_{\vecS\times X/\vecS}(\ast H)$
is naturally
$(\nbigr_{X}(\ast H)_{|\vecS\times X},
 \sigma^{\ast}\nbigr_{X}(\ast H)_{|\vecS\times X})$-modules.
For $\nbigr_X(\ast H)$-modules
$\nbigm_i$ $(i=1,2)$,
a sesqui-linear pairing of $\nbigm_1$ and $\nbigm_2$
is a
$(\nbigr_{X}(\ast H)_{|\vecS\times X},
 \sigma^{\ast}\nbigr_{X}(\ast H)_{|\vecS\times X})$-homomorphism
$C:\nbigm_{1|\vecS\times X}
\times
 \sigma^{\ast}\nbigm_{2|\vecS\times X}
\lrarr
 \distribution_{\vecS\times X/\vecS}(\ast H)$.
For a given $\nbigr_X$-triple
$\nbigt=(\nbigm_1,\nbigm_2,C)$,
we obtain an $\nbigr_X(\ast H)$-triple
$\nbigt(\ast H):=
 \bigl(
 \nbigm_1(\ast H),
 \nbigm_2(\ast H),
 C(\ast H)
 \bigr)$,
where $C(\ast H)$ is a naturally induced
sesqui-linear pairing of 
$\nbigm_1(\ast H)$ and $\nbigm_2(\ast H)$.

\subsubsection{Mixed twistor $\nbigd$-modules}

A pure twistor $\nbigd_X$-module of weight $n$
is an $\nbigr_X$-triple $\nbigt=(\nbigm_1,\nbigm_2,C)$
satisfying some conditions.
For example,
we impose that $\nbigm_i$
are flat over $\nbigo_{\cnum_{\lambda}}$
and that the characteristic variety of $\nbigm_i$
are contained in the product of 
$\cnum_{\lambda}$
and Lagrangian varieties of $T^{\ast}X$.
A polarization of $\nbigt$
is a sesqui-linear duality
$\nbigs:\nbigt\lrarr\nbigt^{\ast}\otimes\newTate(-n)$
satisfying some conditions.
The precise conditions are given
based on the strategy due to M. Saito \cite{saito1}.
We refer the detail
to \cite{sabbah2,sabbah5}
and \cite{Mochizuki-tame,Mochizuki-wild}.
A pure twistor $\nbigd$-module
which admit a polarization is called a polarizable
pure twistor $\nbigd$-module.
In this paper, we consider only polarizable pure twistor
$\nbigd$-modules, and hence we often omit 
``polarizable''.

A mixed twistor $\nbigd$-module
is an $\nbigr$-triple $\nbigt$
with a finite increasing filtration $W$ indexed by integers
satisfying some conditions.
For example, 
we impose that $\Gr^W_n(\nbigt)$ are 
polarizable pure twistor $\nbigd$-modules of weight $n$.
The precise conditions are given
based on the strategy due to M. Saito \cite{saito2}.
We refer the detail to \cite{Mochizuki-MTM}.

We recall some important properties
of mixed twistor $\nbigd$-modules.
Let $\MTM(X)$ denote the category of
mixed twistor $\nbigd$-modules on $X$.

\begin{prop}
\label{prop;14.7.17.1}
The category $\MTM(X)$ is abelian. 
More concretely, we have the following.
\begin{itemize}
\item
Let $\varphi:(\nbigt_1,W_1)\lrarr (\nbigt_2,W_2)$
be a morphism in $\MTM(X)$.
Then, 
the $\nbigr$-triples
$\Ker(\varphi)$ and $\Image(\varphi)$
with the induced filtrations
are mixed twistor $\nbigd$-modules.
Moreover, 
$\varphi$ is strict with respect to 
the filtration $W$,
and $\Image(\varphi)$ with the induced filtration
is a mixed twistor $\nbigd$-module. 
\end{itemize}
Moreover, the category of polarizable pure twistor 
$\nbigd$-modules of weight $w$
is abelian and semisimple.
\hfill\qed
\end{prop}

\begin{prop}
\label{prop;14.7.17.2}
Let $f:X\lrarr Y$ be a projective morphism of 
complex manifolds.
\begin{itemize}
\item
Let $\nbigt$ be a polarizable pure twistor 
$\nbigd_X$-module of weight $w$.
Then,
$f^j_{\dagger}\nbigt$
are polarizable pure twistor $\nbigd_Y$-modules
of weight $w+j$.
\item
Let $(\nbigt,W)$ be a mixed twistor $\nbigd$-module on $X$.
Let $W_kf^j_{\dagger}\nbigt$ be the image of
$f^j_{\dagger}W_{k-j}\nbigt
\lrarr
 f^j_{\dagger}\nbigt$.
Then,
$(f^j_{\dagger}\nbigt,W)$
are mixed twistor $\nbigd$-modules.
\hfill\qed
\end{itemize}
\end{prop}

\begin{prop}
Let $(\nbigt,W)$ be a mixed twistor $\nbigd$-module on $X$.
Let $(\nbigm_1,\nbigm_2,C)$ be the underlying 
$\nbigr_X$-triple.
\begin{itemize}
\item
$\DDD_X\nbigm_i
\simeq
 \nbigh^0\bigl(\DDD_X\nbigm_i\bigr)$.
\item
We have the naturally induced sesqui-linear pairing
$\DDD C$ of
$\DDD_X\nbigm_1$
and 
$\DDD_X\nbigm_2$,
and the $\nbigr_X$-triple
$\DDD\nbigt=(\DDD\nbigm_1,\DDD\nbigm_2,\DDD C)$
with the induced filtration is a mixed twistor structure.
\item
For any projective morphism
$f:X\lrarr Y$,
we have a natural isomorphism
$\DDD f^j_{\dagger}\nbigt
\simeq
 f^{-j}_{\dagger}\DDD\nbigt$.
\hfill\qed
\end{itemize}
\end{prop}

\subsubsection{Underlying $\nbigd$-modules}

Let $(\nbigt,W)$ be a mixed twistor $\nbigd$-module
on $X$.
The $\nbigr$-triple $\nbigt$ is described as $(\nbigm_1,\nbigm_2,C)$.
Let $\iota_1:X\lrarr \nbigx$ be given by
$\iota_1(P)=(1,P)$.
We obtain the $\nbigd$-module
$\Xi_{DR}(\nbigt):=
 \iota_1^{-1}\bigl(\nbigm_2/(\lambda-1)\nbigm_2\bigr)$
which we call the $\nbigd$-module underlying $\nbigt$.
It is naturally equipped with the filtration $W$.
For any morphism of mixed twistor 
$\nbigt$-modules
$\varphi: 
 (\nbigt^{(1)},W)
\lrarr
 (\nbigt^{(2)},W)$,
we have the induced morphism of 
filtered $\nbigd$-modules
$\Xi_{DR}(\varphi):
 \Xi_{DR}(\nbigt^{(1)},W)
\lrarr
 \Xi_{DR}(\nbigt^{(2)},W)$.

The following lemma is easy to see.
\begin{lem}
Let $\varphi: 
 (\nbigt^{(1)},W)
\lrarr
 (\nbigt^{(2)},W)$
be any morphism of mixed twistor 
$\nbigt$-modules.
\begin{itemize}
\item
We have
$\Ker\Xi_{DR}(\varphi)
=\Xi_{DR}(\Ker(\varphi))$,
$\Image\Xi_{DR}(\varphi)
=\Xi_{DR}(\Image(\varphi))$,
and 
$\Cok\Xi_{DR}(\varphi)
=\Xi_{DR}(\Cok(\varphi))$.
\item
$\varphi$ is an epimorphism
(resp. monomorphism)
if and only if 
$\Xi_{DR}(\varphi)$ is 
an epimorphism (resp. monomorphism).
\hfill\qed
\end{itemize}
\end{lem}

The following lemma is easy to see by construction
of the functors.
\begin{lem}
\mbox{{}}
Let $(\nbigt,W)$ be a mixed twistor $\nbigd$-module on $X$.
\begin{itemize}
\item
Let $f:X\lrarr Y$ be a projective morphism of 
complex manifolds.
Then, we have
$f^i_{+}\Xi_{DR}(\nbigt)
\simeq
 \Xi_{DR}(f^i_{\dagger}\nbigt)$.
The induced filtrations $W$ are also equal.
\item
We have a natural isomorphism
$\DDD_X\Xi_{DR}(\nbigt)
\simeq
 \Xi_{DR}\DDD\nbigt$.
The induced filtrations $W$ are also equal.
\hfill\qed
\end{itemize}
\end{lem}

\subsubsection{Integrable case}

Let $\nbigd_{\nbigx}$ be the sheaf of holomorphic differential
operators on $\nbigx$.
Let $\nbigrtilde_X\subset\nbigd_{\nbigx}$
be the sheaf of subalgebras generated by
$\nbigr_X$ and $\lambda^2\del_{\lambda}$.
The correspondence
$\nbigm\lrarr
 \lambda^{-\dim X}\Omega_X
 \otimes\nbigm$
also induces an equivalence between 
the category of left $\nbigrtilde_X$-modules
and the category of right $\nbigrtilde_X$-modules.
An $\nbigrtilde_X$-module
is equivalent to an $\nbigo_{\nbigx}$-module
with a meromorphic flat connection
$\nabla:
\nbigm\lrarr
\Omega^1_{\nbigx}(\log \nbigx^0)\otimes\nbigo(\nbigx^0)
 \otimes\nbigm$.

Let $\nbigm_i$ $(i=1,2)$ be $\nbigrtilde_X$-modules.
Let $\theta$ be the polar coordinate of $\vecS$.
We have
$\del_{\theta}=
 \sqrt{-1}\bigl(
 \lambda\del_{\lambda}
-\lambdabar\del_{\lambdabar}
\bigr)$.
It naturally acts on $\nbigm_i$.
A sesqui-linear pairing of 
$\nbigm_1$ and $\nbigm_2$
is a 
$(\nbigr_{X|\vecS\times X},
 \sigma^{\ast}\nbigr_{X|\vecS\times X})$-homomorphism
$C:\nbigm_{1|\vecS\times X}
 \times
 \sigma^{\ast}\nbigm_{2|\vecS\times X}
\lrarr
 \distribution_{\vecS\times X/\vecS}$
which is compatible with the action of $\del_{\theta}$.
Such a tuple
$(\nbigm_1,\nbigm_2,C)$
is called an $\nbigrtilde_X$-triple
or integrable $\nbigr_X$-triple.
For any morphism of complex manifolds
$f:X\lrarr Y$,
the $\nbigr_X$-triples
$f^j_{\dagger}\nbigt$
are naturally integrable
if $\nbigt$ is integrable.

A mixed twistor $\nbigd$-module $(\nbigt,W)$ 
is called integrable
if the $\nbigr$-triple $\nbigt$ is integrable
and the filtration $W$ is compatible 
with the action of $\lambda^2\del_{\lambda}$.
Propositions \ref{prop;14.7.17.1}
and \ref{prop;14.7.17.2}
are naturally extended to the integrable case.

Let $(\nbigt,W)$ be an integrable mixed twistor $\nbigd_X$-module.
Take an injective resolution $\nbigi$
of $\nbigr_X\otimes\Omega_X^{-1}[\dim X]$
as a left $(\nbigrtilde_X,\nbigr_X)$-module.
Then,
$\DDD_X\nbigm_i\simeq
 \nrhom_{\nbigr_X}(\nbigm,\nbigi)$
are naturally $\nbigrtilde_X$-modules.
The action of $\lambda^2\del_{\lambda}$
is given by 
$(\lambda^2\del_{\lambda}f)(m)
=\lambda^2\del_{\lambda}(f(m))
-f(\lambda^2\del_{\lambda}m)$.
It is shown that $\DDD C$
is compatible with the action of $\del_{\theta}$.
(See \cite{Mochizuki-MTM}.)
Hence, $\DDD(\nbigt,W)$
is also naturally integrable.

We shall give some complements
on the duality and the push-forward 
of $\nbigrtilde$-modules
in \S\ref{section;14.12.26.12}.

\subsection{Integrable mixed twistor $\nbigd$-modules
associated to non-degenerate functions}

\subsubsection{Pure twistor $\nbigd$-modules
associated to meromorphic functions}
\label{subsection;14.12.11.1}

Let $X$ be a complex manifold with a hypersurface $D$.
We set $\nbigx:=\cnum_{\lambda}\times X$
and $\nbigd:=\cnum_{\lambda}\times D$.
Let $f$ be a meromorphic function on $(X,D)$.
We have a wild harmonic bundle $E(f)$ on $(X,D)$.
It consists of a Higgs bundle
$(\nbigo_{X\setminus D}\,e,df)$
with a pluri-harmonic metric $h$
determined by $h(e,e)=1$,
where $e$ denotes a global frame of the line bundle.
Recall that 
we have the associated 
polarizable pure twistor $\nbigd$-module
$(\nbigl(f),\nbigl(f),C(f))$
of weight $0$ on $X$,
where $\nbigl(f)$ is an $\nbigr_X$-module,
and $C(f)$ is a sesqui-linear pairing
$\nbigl(f)_{|\vecS\times X}
\times
 \sigma^{\ast}\nbigl(f)_{|\vecS\times X}
\lrarr
 \distribution_{\vecS\times X/\vecS}$.
Here $\vecS:=\{|\lambda|=1\}\subset\cnum_{\lambda}$.
A polarization is given by $(\id,\id)$.

\vspace{.1in}

Let $d_X:=\dim X$.
We have the polarizable pure twistor 
$\nbigd$-module
$\nbigt(f):=\bigl(
 \lambda^{d_X}\nbigl(f),
 \nbigl(f),
 C(f)
 \bigr)$
of weight $d_X$.
The natural polarization is given by
$\bigl((-1)^{d_X},(-1)^{d_X}\bigr)$.

\vspace{.1in}

The restriction $\nbigl(f)_{|\nbigx\setminus \nbigd}$
is equipped with a global frame $\upsilon$
such that 
$\nbigl(f)_{|\nbigx\setminus \nbigd}$
is isomorphic to
$\nbigo_{\nbigx\setminus\nbigd}\upsilon$
with the family of flat $\lambda$-connections $\DD$
determined by $\DD\upsilon=\upsilon df$.
We have
\begin{equation}
\label{eq;15.1.4.30}
 C(f)(\upsilon,\sigma^{\ast}\upsilon)
=\exp\bigl(
 -\lambda \overline{f}
+\lambda^{-1}f
 \bigr).
\end{equation}

In general,
it is not easy to describe $\nbigl(f)$ explicitly.
The following lemma is clear by the construction
of $\nbigl(f)$ in \cite{Mochizuki-wild},
and it will be used implicitly.
\begin{lem}
$\upsilon$ is naturally extended as a section of
$\nbigl(f)(\ast D)$,
and it gives a global frame of 
$\nbigl(f)(\ast D)$
as an $\nbigo_{\nbigx}(\ast D)$-module.
Moreover, the natural morphism
$\nbigl(f)\lrarr\nbigl(f)(\ast D)$
is a monomorphism.
\hfill\qed
\end{lem}

\begin{lem}
If $|(f)_0|\cap D=\emptyset$
and $|(f)_{\infty}|=D$,
then $\nbigl(f)$ is naturally isomorphic to 
$\nbigo_{\nbigx}(\ast D)\upsilon$.
\end{lem}
\pf
If $D$ is normal crossing,
the claim follows from the construction 
of $\nbigl(f)$ in \cite{Mochizuki-wild}.
Although we are mainly interested in the case
where $D$ is normal crossing,
we give a proof in the general case.
We take a projective morphism
$\varphi:X'\lrarr X$ of complex manifolds
such that
(i) $D':=\varphi^{-1}(D)$ is normal crossing,
(ii) $\varphi$ induces an isomorphism
$X'\setminus D'\simeq X\setminus D$.
We set $f':=\varphi^{\ast}(f)$.
We have
$\nbigl(f')=\nbigo_{\nbigx'}(\ast D')v'$
with a frame $v'$
such that $\DD v'=v'\,df'$.
Because 
$\nbigl(f')(\ast D')=\nbigl(f')$,
we obtain
$\varphi_{\dagger}\bigl(\nbigl(f')\bigr)
\simeq
 \varphi_{\dagger}\bigl(\nbigl(f')\bigr)(\ast D)$.
We also have that
$\nbigl(f)$ is a direct summand 
of $\varphi_{\dagger}\nbigl(f')$
and 
$\varphi_{\dagger}\bigl(\nbigl(f')\bigr)
 _{|X\setminus D}
=\nbigl(f)_{|X\setminus D}$.
Then, we can deduce the claim of the lemma.
\hfill\qed

\begin{prop}
Suppose that 
(i) $D$ is normal crossing,
(ii) $f$ is non-degenerate along $D$,
(iii) $D=|(f)_{\infty}|$.
Then, $\nbigl(f)$ is naturally isomorphic to 
$\nbigo_{\nbigx}(\ast D)\upsilon$.
\end{prop}
\pf
We have only to check the claim
locally around any point of 
$|(f)_{0}|\cap |(f)_{\infty}|$.
Let $(z_1,\ldots,z_n)$ be a convenient coordinate system
with $f=z_n\prod_{i=1}^{\ell}z_i^{-k_i}$
for some $k_i>0$.
Note that
$\prod_{i=1}^{\ell}z_i^N\upsilon$
is a section of $\nbigl(f)$.
We have
$\deldel_n\upsilon=\upsilon \prod_{i=1}^{\ell}z_i^{-k_i}$,
where $\deldel_n=\lambda\del_{x_n}$.
Then, it is easy to deduce that
$\nbigl(f)\supset
 \nbigr_{X}\upsilon=\nbigo_{\nbigx}(\ast \nbigd)\upsilon$.
\hfill\qed

\subsubsection{Mixed twistor $\nbigd$-modules
 associated to meromorphic functions}

For $\star=\ast,!$,
we have the $\nbigr_X$-modules
$\nbigl(f)[\star D]$
obtained as the localization
$\nbigl(f)$.
They are denoted by
$\nbigl_{\star}(f,D)$
in this paper.
If $D=|(f)_{\infty}|$,
they are also denoted by 
$\nbigl_{\star}(f)$.
We have the natural morphisms
of $\nbigr_X$-modules
$\nbigl_!(f,D)
\lrarr
 \nbigl(f)
\lrarr
 \nbigl_{\ast}(f,D)$.
It induces
$\nbigl_{\star}(f,D)(\ast D)
\simeq
 \nbigl(f)(\ast D)$.

For $\star=\ast,!$,
we have the mixed twistor $\nbigd$-modules
$\nbigt(f)[\star D]$
obtained as the localizations of
$\nbigt(f)$.
They are denoted by
$\nbigt_{\star}(f,D)$
in this paper.
If $D=|(f)_{\infty}|$,
they are also denoted by
$\nbigt_{\star}(f)$.
The underlying $\nbigr_X$-triples of
$\nbigt_{\ast}(f,D)$ 
and 
$\nbigt_!(f,D)$ are
$\bigl(
 \lambda^{d_X}\nbigl_!(f,D),
 \nbigl_{\ast}(f,D),
 C(f)[\ast D]
 \bigr)$
and 
$\bigl(
 \lambda^{d_X}\nbigl_{\ast}(f,D),
 \nbigl_{!}(f,D),
 C(f)[!D]
 \bigr)$,
respectively.
The weight filtrations of $\nbigt_{\star}(f,D)$
are denoted by $W$.
We have natural morphisms of mixed twistor 
$\nbigd$-modules
$\nbigt_!(f,D)\lrarr \nbigt(f)
\lrarr
 \nbigt_{\ast}(f,D)$.

The meromorphic family of flat connections $\DD^f$
on $\nbigl(f)(\ast D)$ is naturally extended to 
a meromorphic flat connection $\nabla$
given by
$\nabla \upsilon=\upsilon d(\lambda^{-1}f)$.
Let $C(f)(\ast D)$ be the sesqui-linear pairing on
$\nbigl(f)(\ast D)$ induced by $C(f)$.
Due to (\ref{eq;15.1.4.30}),
$C(f)(\ast D)$ is compatible with the actions of
$\lambda^{2}\del_{\lambda}$.

\begin{lem}
$\nbigt_{\star}(f,D)$ are integrable mixed twistor 
$\nbigd$-modules,
and $\nbigt(f,D)$ is an integrable 
pure twistor $\nbigd$-module.
The morphisms
$\nbigt_{!}(f,D)
\lrarr
 \nbigt(f,D)
\lrarr
 \nbigt_{\ast}(f,D)$
are integrable.
\end{lem}
\pf
By using \cite[Lemma 3.2.4]{Mochizuki-MTM},
we obtain that the underlying $\nbigr$-triples
of $\nbigt_{\star}(f,D)$ are integrable.
Then, we obtain that 
$\nbigt(f,D)$ is an integrable pure twistor $\nbigd$-module.
By using \cite[Lemma 7.1.35]{Mochizuki-MTM},
we obtain that the filtrations $W$
of $\nbigt_{\star}(f,D)$ are also integrable.
Hence, we obtain that
$\nbigt_{\star}(f,D)$ are integrable
mixed twistor $\nbigd$-modules.
By using \cite[Lemma 7.1.37]{Mochizuki-MTM},
we obtain the integrability of the morphisms.
\hfill\qed

\begin{lem}
\label{lem;15.1.5.1}
The $\nbigd$-modules
$\Xi_{DR}(\nbigt_{\star}(f,D))$
and $\Xi_{DR}(\nbigt(f))$
are naturally identified with
$L_{\star}(f,D)$ and $L(f)$.
The morphisms
$\Xi_{DR}(\nbigt_!(f,D))
\lrarr
 \Xi_{DR}(\nbigt(f))
\lrarr 
 \Xi_{DR}(\nbigt_{\ast}(f,D))$
are naturally identified with 
$L_!(f,D)\lrarr 
 L(f)\lrarr L_{\ast}(f,D)$.
\end{lem}
\pf
Let $\nbigm$ be an $\nbigr_X$-module
underlying $\nbigt_{\star}(f,D)$ $(\star=\ast,!)$
or $\nbigt(f,D)$.
Let $g$ be a holomorphic function
on an open subset $U\subset X$.
Because $\nbigt_{\star}(f,D)$ are integrable,
the KMS-spectrum of $\nbigm_{|U}$
along $g$ are contained in
$\real\times\{0\}$.
Hence, the specialization of 
the $V$-filtration along $g$ at $\lambda=1$
is equal to the $V$-filtration of 
$\iota_1^{-1}(\nbigm/(\lambda-1)\nbigm)$,
where $\iota_1:X\lrarr \nbigx$ given by
$\iota_1(P)=(1,P)$.
Then, the claim is clear.
\hfill\qed

\vspace{.1in}
In the proof of Lemma \ref{lem;15.1.5.1},
the KMS-spectrum of the $\nbigr_X$-modules
are contained in $\rnum\times\{0\}$, indeed.
We can directly check it by standard computations.

\begin{lem}
Suppose that (i) $D=|(f)_{\infty}|$,
(ii) $f$ is pure at $P\in X$.
Then, the canonical morphisms
$\nbigt_{!}(f,D)
\lrarr 
 \nbigt(f)
\lrarr
 \nbigt_{\ast}(f,D)$
are isomorphisms
on a neighbourhood of $P$.
\end{lem}
\pf
It follows from that 
the morphisms of the underlying $\nbigd$-modules
are isomorphisms.
\hfill\qed

\begin{cor}
\label{cor;14.10.18.2}
Suppose that 
(i) $D$ is normal crossing,
(ii) $D=|(f)_{\infty}|$,
(iii) $f$ is non-degenerate.
Then 
the canonical morphisms
$\nbigt_{!}(f)
\lrarr 
 \nbigt(f)
\lrarr
 \nbigt_{\ast}(f)$
are isomorphisms.
In particular,
we have isomorphisms
$\nbigl_{\star}(f,D)\simeq
 \nbigo_{\nbigx}(\ast D)\,\upsilon$.
\hfill\qed
\end{cor}

\subsubsection{Expression of $\nbigrtilde$-modules
in non-degenerate cases}

Suppose that $D$ is normal crossing.
Let us describe $\nbigr_X$-modules
$\nbigl_{\star}(f,D)$
in the case that 
$f$ is non-degenerate along $D$.
We do not assume $D=|(f)_{\infty}|$.
Let $V_D\nbigr_X$
denote the sheaf of subalgebras in $\nbigr_X$
generated by
$\lambda\cdot
 p_{\lambda}^{\ast}\Theta_{X}(\log D)$
over $\nbigo_{\nbigx}$.
Let $V_D\nbigrtilde_X$ denote the sheaf of subalgebras
in $\nbigrtilde_X$
generated by
$V_D\nbigr_X$ and $\lambda^2\del_{\lambda}$.
We have the $V_D\nbigrtilde_X$-modules
$\nbigl(f)$ and
$\nbigl(f)(\nbigd)
=\nbigl(f)\otimes_{\nbigo_{\nbigx}}\nbigo_{\nbigx}(\nbigd)$.

\begin{prop}
\label{prop;14.10.18.1}
We have natural isomorphisms of
$\nbigrtilde_X$-modules
\[
 \nbigl_{\ast}(f,D)\simeq
 \nbigr_X\otimes_{V_D\nbigr_X}
  \nbigl(f)\bigl(\nbigd\bigr),
\quad
 \nbigl_{!}(f,D)\simeq
 \nbigr_X\otimes_{V_D\nbigr_X}\nbigl(f).
\]
\end{prop}
\pf
We have natural isomorphisms
of $\nbigrtilde_X(\ast D)$-modules
\[
 \nbigl_{\ast}(f,D)(\ast D)
 \simeq
 \bigl(
 \nbigr_X\otimes_{V_D\nbigr_X}
  \nbigl(f)\bigl(\nbigd\bigr)
 \bigr)(\ast D),
\quad
 \nbigl_{!}(f,D)(\ast D)\simeq
 \bigl(
 \nbigr_X\otimes_{V_D\nbigr_X}\nbigl(f)
 \bigr)(\ast D).
\]
Note that if we are given an $\nbigr_X$-endomorphism
$g$ of $\nbigl_{\star}(f,D)$ such that 
$g_{|\cnum_{\lambda}\times (X\setminus D)}=0$,
then we have $g=0$.
Hence,
it is enough to prove
that there exist isomorphisms
as $\nbigr_X$-modules.

The claim is clear 
outside $|(f)_0|\cap|(f)_{\infty}|$
by the construction 
of $\nbigl_{\star}(f,D)$ in 
\cite[\S5.3]{Mochizuki-MTM}.
We set 
\[
 \nbigl'_{\ast}(f,D):=
  \nbigr_X\otimes_{V_D\nbigr_X}
  \nbigl(f)\bigl(\nbigd\bigr),
\quad\quad
 \nbigl'_{!}(f,D):=
 \nbigr_X\otimes_{V_D\nbigr_X}\nbigl(f).
\]
Let $Q\in|(f)_0|\cap|(f)_{\infty}|$.
We take a convenient coordinate system
$(z_1,\ldots,z_n)$
with 
$D=\bigcup_{i=1}^{\ell_1+\ell_2}\{z_i=0\}$
and $f=z_n\prod_{i=1}^{\ell_1}z_i^{-k_i}$.
We set
$X_1:=\{(z_1,\ldots,z_{\ell_1},z_n)\}$
and 
$X_2:=\{(z_{\ell_1+1},\ldots,z_{n-1})\}$.
We set $f_1:=z_n\prod_{i=1}^{\ell_1}z_i^{-k_i}$
and $D_1:=\bigcup_{i=1}^{\ell_1}\{z_i=0\}$
and $D_2:=\bigcup_{i=\ell_1+1}^{\ell_1+\ell_2}\{z_i=0\}$.
We have
$\nbigl'_{\star}(f_1,D_1)
=\nbigl(f_1)
=\nbigl_{\star}(f_1,D_1)$.
In particular,
we have
$\nbigl'_{\star}(f_1,D_1)[\star z_j]
=\nbigl'_{\star}(f_1,D_1)$
for $j=1,\ldots,\ell_1$.
We also have
$\nbigl'_{\star}(f_2,D_2)
\simeq
 \nbigl_{\star}(f_2,D_2)$.
Hence, 
$\nbigl'_{\star}(f_2,D_2)[\star z_i]
\simeq
 \nbigl'_{\star}(f_2,D_2)$
for $i=\ell_1+1,\ldots,\ell_1+\ell_2$.
We have
$\nbigl'_{\star}(f,D)
\simeq
 \nbigl'_{\star}(f_1,D_1)
\boxtimes
 \nbigl'_{\star}(f_2,D_2)$.
Then, it is easy to check that
$\nbigl'_{\star}(f,D)[\star z_i]
=\nbigl'_{\star}(f,D)$
for $i=\ell_1+1,\ldots,\ell_1+\ell_2$.
We obtain the desired isomorphism
by the characterization in 
\cite[Proposition 5.3.1]{Mochizuki-MTM}.
\hfill\qed

\subsubsection{De Rham complexes in  non-degenerate cases}
\label{subsection;14.11.6.11}

Let us consider the case that
$(X,D)=S\times (X_0,D_0)$,
where $X_0$ is a complex manifold
with a normal crossing hypersurface $D_0$.
Let $f$ be a meromorphic function on $(X,D)$
which is non-degenerate along $D$.
Let $\pi:X\lrarr S$ denote the projection.
We have the mixed twistor $\nbigd$-modules
$\pi_{\dagger}^i\nbigt_{\star}(f,D)$ on $S$.
Let us describe 
the underlying $\nbigr_{S}$-modules
$\pi_{\dagger}^i\nbigl_{\star}(f,D)$,
by assuming that $X_0$ is projective.

Set $n:=\dim X_0$.
Let $p_{\lambda}:\nbigx\lrarr X$
denote the projection.
For any $\nbigr_X$-module $\nbigm$,
we set 
\[
\DR_{X/S}^{m}\nbigm:=
\lambda^{-m-n}
 p_{\lambda}^{\ast}
 \bigl(
 \Omega^{m+n}_{X/S}
 \bigr)
\otimes_{\nbigo_{\nbigx}}
 \nbigm.
\]
Then, we obtain the complex
$\DR_{X/S}^{\bullet}\nbigm$.
We naturally have
$\pi_{\dagger}^i\nbigm
\simeq
 \hyperr^i\pi_{\ast}
 \DR^{\bullet}_{X/S}\nbigm$.

\vspace{.1in}
Suppose $|(f)_{\infty}|=D$.
We have 
$\nbigl_{\star}(f,D)=\nbigl(f)$
which is isomorphic to
$\nbigo_{\nbigx}
 (\ast D)\upsilon$
with $\DD\upsilon=\upsilon df$.
We set 
$\Omegabar^{\ell}_{X/S}:=
 \lambda^{-\ell}
 p_{\lambda}^{\ast}\Omega_{X/S}^{\ell}$.
Then, we immediately obtain the following natural isomorphism:
\[
 \DR^{\bullet}_{X/S}\nbigl(f)
\simeq
\Bigl(
 \Omegabar^{\bullet+n}_{X/S}(\ast D),
 d+\lambda^{-1}df
\Bigr)
\]

Let us consider the case that
$D$ is not necessarily $|(f)_{\infty}|$.
We set
$\Omegabar^{\ell}_{X/S}(\log D):=
 \lambda^{-\ell}
 p_{\lambda}^{\ast}\bigl(
 \Omega^{\ell}_{X/S}(\log D)\bigr)$.
We obtain the complexes
$\bigl(
 \Omegabar^{\bullet}_{X/S}(\log D),d+\lambda^{-1}df
 \bigr)$
and
$\bigl(
 \Omegabar^{\bullet}_{X/S}(\log D)(-\nbigd),
 d+\lambda^{-1}df
 \bigr)$.
\begin{lem}
We have the following natural quasi-isomorphisms:
\begin{equation}
\label{eq;14.11.23.200}
 \DR^{\bullet}_{X/S}\nbigl_{\ast}(f,D)
\simeq
 \Bigl(
 \Omegabar^{\bullet+n}_{X/S}(\log D)\bigl(\ast(f)_{\infty}\bigr)
 ,d+\lambda^{-1}df
\Bigr)
\end{equation}
\begin{equation}
\label{eq;14.11.23.201}
 \DR^{\bullet}_{X/S}\nbigl_{!}(f,D)
\simeq
 \Bigl(
 \Omegabar^{\bullet+n}_{X/S}(\log D)(-\nbigd)
 \bigl(\ast(f)_{\infty}\bigr)
 ,d+\lambda^{-1}df
\Bigr)
\end{equation}
Hence, 
$\pi^i_{\dagger}
 \nbigl_{\star}(f,D)$ $(\star=\ast,!)$
are expressed as the push-forward
of the right hand side of
{\rm (\ref{eq;14.11.23.200})} and 
{\rm(\ref{eq;14.11.23.201})}.
\end{lem}
\pf
It follows from the expression 
in Proposition \ref{prop;14.10.18.1}.
\hfill\qed

\vspace{.1in}
Suppose moreover that
$f$ is non-degenerate along $D$ over $S$.
We use the notation in \S\ref{subsection;14.11.6.3}.
We set
$\Omegabar^{\ell}_{X/S}(\log D,f):=
 \Omegabar^{\ell}_{X/S}(\log D)\bigl(\ell(f)_{\infty}\bigr)$,
and we obtain the following complexes:
\[
\bigl(
 \Omegabar^{\bullet}_{X/S}(\log D,f),d+\lambda^{-1}df
 \bigr),
\quad
 \bigl(
 \Omegabar^{\bullet}_{X/S}(\log D,f)(-\nbigd)
 ,d+\lambda^{-1}df
 \bigr).
\]
Similarly, we set 
$\Omegabar^{\ell}_{X/S,f,D}:=
 \lambda^{-\ell}
 p_{\lambda}^{\ast}
 \Omega^{\ell}_{X/S,f,D}$,
and we obtain the following complexes:
\[
  \bigl(
 \Omegabar^{\bullet}_{X/S,f,D},
 d+\lambda^{-1}df
 \bigr),
\quad\quad
 \bigl(
 \Omegabar^{\bullet}_{X/S,f,D}(-\nbigd),
 d+\lambda^{-1}df
 \bigr)
\]
As in the case of $\nbigd$-modules
in \S\ref{subsection;14.11.6.3} and 
\S\ref{subsection;14.11.6.2},
we have the following lemma.
\begin{lem}
Suppose that $f$ is non-degenerate along $D$ over $S$.
Then, we have the following natural quasi-isomorphisms:
\begin{equation}
 \label{eq;14.11.23.202}
 \DR^{\bullet}_{X/S}\nbigl_{\ast}(f,D)
\simeq
 \Bigl(
 \Omegabar^{\bullet+n}_{X/S}(\log D,f)
 ,d+\lambda^{-1}df
\Bigr)
\simeq
 \bigl(
 \Omegabar^{\bullet}_{X/S,f,D},
 d+\lambda^{-1}df
 \bigr)
\end{equation}
\begin{equation}
 \label{eq;14.11.23.203}
 \DR^{\bullet}_{X/S}\nbigl_{!}(f,D)
\simeq
 \Bigl(
 \Omegabar^{\bullet+n}_{X/S}(\log D,f)(-\nbigd)
 ,d+\lambda^{-1}df
\Bigr)
\simeq
 \bigl(
 \Omegabar^{\bullet}_{X/S,f,D}(-\nbigd),
 d+\lambda^{-1}df
 \bigr)
\end{equation}
Hence,
$\pi_{\dagger}^i\nbigl_{\star}(f,D)$
are expressed as the push-forward of the complexes
in {\rm(\ref{eq;14.11.23.202})}
and {\rm(\ref{eq;14.11.23.203})}.
\hfill\qed
\end{lem}

We also remark the following lemma.
\begin{lem}
Suppose that $f$ is non-degenerate along $D$
over $S$.
Then, 
$\pi^i_{\dagger}\nbigl_{\star}(f,D)$
are locally free 
$\nbigo_{\cnum_{\lambda}\times S}$-modules.
\end{lem}
\pf
The specialization at $\lambda=1$
are locally free $\nbigo_S$-modules
by Corollary \ref{cor;14.7.10.10}.
So, by the general property of 
mixed twistor $\nbigd$-modules,
we obtain that 
$\pi^i_{\dagger}\nbigl_{\star}(f,D)$
are locally free 
$\nbigo_{\cnum_{\lambda}\times S}$-modules.
\hfill\qed

\begin{cor}
If $f$ is non-degenerate along $D$ over $S$,
the $\nbigo_{\cnum_{\lambda}\times S}$-modules
\[
 \hyperr^i\pi_{\ast}
 \bigl(
 \Omegabar^{\bullet}_{X/S}(\log D)(\ast (f)_{\infty}),d+\lambda^{-1}df
 \bigr)
\simeq
  \hyperr^i\pi_{\ast}
 \bigl(
 \Omegabar^{\bullet}_{X/S}(\log D,f),d+\lambda^{-1}df
 \bigr)
\simeq
 \hyperr^i\pi_{\ast}
 \bigl(
 \Omegabar^{\bullet}_{X/S,f,D},d+\lambda^{-1}df
 \bigr)
\]
are locally free.
The $\nbigo_{\cnum_{\lambda}\times S}$-modules
\begin{multline*}
 \hyperr^i\pi_{\ast}
 \bigl(
 \Omegabar^{\bullet}_{X/S}(\log D) (\ast (f)_{\infty})(-\nbigd),
 d+\lambda^{-1}df
 \bigr)
\simeq
  \hyperr^i\pi_{\ast}
 \bigl(
 \Omegabar^{\bullet}_{X/S}(\log D,f)(-\nbigd),
 d+\lambda^{-1}df
 \bigr)
 \\
\simeq
 \hyperr^i\pi_{\ast}
 \bigl(
 \Omegabar^{\bullet}_{X/S,f,D}(-\nbigd),
 d+\lambda^{-1}df
 \bigr)
\end{multline*}
are also locally free.
\hfill\qed
\end{cor}

\subsection{The real structure of 
the localization of some mixed twistor $\nbigd$-modules}
\label{subsection;14.12.10.1}

Let $X$ be a complex manifold.
For a mixed twistor $\nbigd$-module
$\nbigt=(\nbigm_1,\nbigm_2,C)$ on $X$,
we have
\[
j^{\ast}\nbigt=(j^{\ast}\nbigm_1,j^{\ast}\nbigm_2,j^{\ast}C),
\quad
\DDD\nbigt=
 (\DDD\nbigm_1,\DDD\nbigm_2,\DDD C).
\]
We set
$\gammatilde^{\ast}\nbigt:=
 (j^{\ast}\DDD\nbigm_2,j^{\ast}\DDD\nbigm_1,
 j^{\ast}\DDD C^{\ast})$.
We will naturally identify
$j^{\ast}\DDD\nbigm_i$
and $\DDD j^{\ast}\nbigm_i$.
Recall that a real structure of $\nbigt$ is 
an isomorphism
$\kappa:
 \gammatilde^{\ast}
\nbigt
\simeq
\nbigt$
satisfying
$\gammatilde^{\ast}\kappa\circ\kappa=\id$.

\subsubsection{Basic case}
\label{subsection;14.11.8.50}

Let $d_X:=\dim X$.
We have the isomorphism
$\nu_X:
 \DDD\nbigo_{\nbigx}\simeq
 \nbigo_{\nbigx}\lambda^{d_X}$,
whose specialization at $\{\lambda_0\}\times X$
$(\lambda_0\neq 0)$
is equal to the morphism in \S\ref{subsection;14.11.14.1}.
We have the natural identification
$j^{\ast}\nbigo_{\nbigx}=
 \nbigo_{\nbigx}$
given by the pull back of functions.
We have
$\DDD\nu_X=(-1)^{d_X}\nu_X$.
As shown in \cite{Mochizuki-MTM},
the isomorphism
\[
 \bigl(
 \nu_X^{-1},
 (-1)^{d_X}\nu_X
 \bigr):
 \gammatilde^{\ast}\nbigu_X(d_X,0)
\simeq
 \nbigu_X(d_X,0)
\]
gives a real structure.

Let $Y$ be a smooth hypersurface of $X$.
We set $d_Y:=\dim Y=d_X-1$.
Let $\iota:Y\lrarr X$ be the inclusion.
We have the following natural morphisms
of integrable mixed twistor $\nbigd$-modules:
\begin{equation}
\label{eq;14.10.4.30}
\nbigu_X(d_X,0)[\ast Y]
\lrarr
 \iota_{\dagger}\nbigu_{Y}(d_Y,0)
 \otimes\newTate(-1)
\end{equation}
\begin{equation}
\label{eq;14.10.4.31}
\iota_{\dagger}\nbigu_{Y}(d_Y,0)
\lrarr
\nbigu_X(d_X,0)[!Y]
\end{equation}
The morphisms are induced 
by $\nbigr_X$-homomorphisms
$\nbigo_{\nbigx}[\ast Y]
\lrarr
 \iota_{\dagger}\nbigo_{\nbigy}\lambda^{-1}$
and 
$\iota_{\dagger}\nbigo_{\nbigy}
\lrarr 
\nbigo_{\nbigx}[!Y]$.
Locally, 
if $X$ is equipped with a holomorphic coordinate system
$(x_1,\ldots,x_n)$ such that $Y=\{x_1=0\}$,
the morphisms are given by
$x_1^{-1}
\longmapsto
 \lambda^{-1}
\iota_{\ast}(dx_1/\lambda)^{-1}$
and 
$\iota_{\ast}(dx_1/\lambda)^{-1}
\longmapsto
-\deldel_{1}(1)$,
where $\deldel_1=\lambda\del_{x_1}$.

\begin{prop}
\label{prop;14.9.12.11}
The natural morphisms
{\rm(\ref{eq;14.10.4.30})}
and 
{\rm(\ref{eq;14.10.4.31})}
are compatible with the real structures.
\end{prop}
\pf
Let us look at the compatibility of
{\rm(\ref{eq;14.10.4.30})}
which is the commutativity of the following diagrams:
\[
 \begin{CD}
 \nbigo_{\nbigx}[!Y]
 \lambda^{d_X}
 @<{\nu_X}<<
 \DDD\bigl(\nbigo_{\nbigx}[\ast Y]\bigr)
 \\
 @AAA @AAA \\
 \iota_{\dagger}\nbigo_{\nbigy}\lambda^{d_Y+1}
 @<{-\nu_Y}<<
 \DDD\bigl(
 \iota_{\dagger}\nbigo_{\nbigy}\lambda^{-1}
 \bigr)
 \end{CD}
\quad\quad\quad
 \begin{CD}
 \nbigo_{\nbigx}[\ast Y]
 @<{(-1)^{d_X}\nu_X}<<
 \DDD\bigl(\nbigo_{\nbigx}[!Y]\lambda^{d_X}\bigr)
 \\
 @VVV @VVV \\
 \iota_{\dagger}\nbigo_{\nbigy}\lambda^{-1}
 @<{(-1)^{d_Y+1}\nu_Y}<<
 \DDD\bigl(
 \iota_{\dagger}\nbigo_{\nbigy}\lambda^{d_Y+1}
 \bigr)
 \end{CD}
\]
To check the commutativity,
it is enough to compare the specialization
along $\{\lambda_0\}\times X$
for any $\lambda_0\neq 0$.
Then, it is reduced to Proposition \ref{prop;14.9.12.10}.
We can check the compatibility of 
(\ref{eq;14.10.4.31}) similarly.
\hfill\qed

\subsubsection{Mixed twistor $\nbigd$-modules
associated to holomorphic functions}
\label{subsection;15.11.10.10}

Let $F$ be any holomorphic function on $X$.
We have the $\nbigr_X$-module $\nbigl(F)$
given by $\nbigo_{\nbigx}\,e$ with $\DD e=e\,dF$,
where $e$ is a global frame.
Let $\nbigl(F)^{\lor}=
 \nhom_{\nbigo_{\nbigx}}(\nbigl(F),\nbigo_{\nbigx})$ 
is naturally isomorphic to
$\nbigo_{\nbigx}e^{\lor}$
with $\DD e^{\lor}=e^{\lor}(-dF)$,
where $e^{\lor}$ is the dual frame of $e$.
Hence, we have the isomorphism
$j^{\ast}\nbigl(F)^{\lor}
\simeq 
 \nbigl(F)$
given by $j^{\ast}e^{\lor}\longleftrightarrow e$.
We have the smooth
$\nbigr_X$-triple
$\nbigt_{sm}(F):=(\nbigl(F),\nbigl(F),C_F)$.
The isomorphism
$\gammatilde_{sm}^{\ast}\nbigt_{sm}(F):=
 (j^{\ast}\nbigl(F)^{\lor},j^{\ast}\nbigl(F),j^{\ast}C_F^{\lor})
\simeq
 \nbigt(F)$
is a real structure of $\nbigt_{sm}(F)$
as a smooth $\nbigr_X$-triple.
(See \cite[\S2.1.7.2]{Mochizuki-MTM}.)
Because 
$\nbigt(F)=
 \nbigt_{sm}(F)\otimes
 \nbigu_X(d_X,0)$,
we have the induced real structure 
on $\nbigt(F)$
as a mixed twistor $\nbigd$-module
\cite[Proposition 13.4.6]{Mochizuki-MTM}.

More explicitly,
as in the case of $\nbigo_{\nbigx}$,
we have a natural isomorphism
$\nu:j^{\ast}\DDD\nbigl(F)\simeq
 \lambda^{d_X}j^{\ast}\nbigl(F)^{\lor}
\simeq \lambda^{d_X}\nbigl(F)$.
The real structure of 
$\nbigt(F)=(\lambda^{d_X}\nbigl(F),\nbigl(F),C_F)$
is given by
$(\nu^{-1},(-1)^{d_X}\nu)$.

\begin{lem}
We have natural isomorphisms of
integrable mixed twistor $\nbigd$-modules
with real structure
\begin{equation}
 \label{eq;15.11.10.20}
 \nbigt(F)[\star Y]
\simeq
 \nbigt_{sm}(F)
 \otimes\nbigu_X(d_X,0)[\star Y].
\end{equation}
\end{lem}
\pf
We have the natural isomorphisms
of the integrable $\nbigr_X$-modules
$\nbigl(F)[\star Y]
\simeq
 \nbigl(F)\otimes\nbigo_{\nbigx}[\star Y]$ $(\star=\ast,!)$
such that the restriction to $X\setminus Y$ is the identity.
It is enough to check the claim locally around any point of $Y$.
We may assume that $Y$ is given as $\{t=0\}$
for a holomorphic function $t$ such that $dt$ is nowhere vanishing.
Then,
it is easy to see that
the $V$-filtration of $\nbigo_{\nbigx}[\star Y]$ along $t$
induces a $V$-filtration of
$\nbigl(F)\otimes\nbigo_{\nbigx}[\star Y]$
along $t$,
and that it satisfies the characterization condition
for $\nbigl(F)[\star t]$ in \cite[\S3.1]{Mochizuki-MTM}.
Because the sesqui-linear pairings
of $\nbigt(F)[\star Y]$
and
$\nbigt_{sm}(F)
 \otimes\nbigu_X(d_X,0)[\star Y]$
are the same on $X\setminus Y$,
we obtain that they are the same on $X$.
(See \cite[Proposition 3.2.1]{Mochizuki-MTM}.)
Because the real structures of
$\nbigt(F)[\star F]$ are uniquely determined
by their restriction to $X\setminus Y$,
the isomorphisms (\ref{eq;15.11.10.20})
are compatible with the real structures.
\hfill\qed

\vspace{.1in}

We also obtain the integrable mixed twistor $\nbigd$-modules
$\nbigt_{sm}(F)\otimes\iota_{\dagger}\nbigu_Y(d_Y,0)$
with real structure.
We have the isomorphisms 
of integrable mixed twistor $\nbigd$-modules
with real structure:
\begin{equation}
 \label{eq;15.11.10.40}
 \nbigt(F)[\ast Y]/\nbigt(F)
\simeq
 \nbigt_{sm}(F)\otimes\iota_{\dagger}\nbigu_Y(d_Y,0)
 \otimes\newTate(-1)
\end{equation}
\begin{equation}
\label{eq;15.11.10.41}
 \Ker\bigl(
 \nbigt(F)[!Y]\lrarr \nbigt(F)
 \bigr)
\simeq
 \nbigt_{sm}(F)\otimes
 \iota_{\dagger}\nbigu_Y(d_Y,0)
\end{equation}

Let $F_Y$ denote the restriction of $F$ to $Y$.
We have the integrable mixed twistor $\nbigd$-module 
$\nbigt(F_Y)$ on $Y$
with real structure on $Y$.
It is given as
$\nbigt(F_Y)=
\bigl(\lambda^{d_Y}\nbigl(F_Y),\nbigl(F_Y),C_{F_Y}
\bigr)$.

Let $\omega_{Y/X}$ denote the conormal bundle of $Y$ in $X$.
We set $\omega_{\nbigy/\nbigx}:=
 \lambda \cdot p_{\lambda}^{\ast}\omega_{Y/X}$,
where $p_{\lambda}:\nbigy\lrarr X$ denotes the projection.
We have the $\nbigo_{\nbigx}$-submodules
\[
 \nbige_0:=\nbigl(F)\otimes 
 \iota_{\ast}(\omega_{\nbigy/\nbigx})
\subset \nbigl(F)\otimes\iota_{\dagger}\nbigo_{\nbigx},
\quad 
\nbige_1:=\iota_{\ast}(\nbigl(F_Y)\otimes\omega_{\nbigy/\nbigx})
\subset
 \iota_{\dagger}\nbigl(F_Y).
\]
We have the natural isomorphism
of $\nbigo_{\nbigx}$-modules
$\varphi_0:\nbige_0\simeq\nbige_1$.

\begin{lem}
\label{lem;15.11.10.30}
The $\nbigo$-isomorphism
$\varphi_0$
is uniquely extended to
an isomorphism of integrable $\nbigr_X$-modules
$\varphi:
  \iota_{\dagger}\nbigl(F_Y)
\simeq
  \nbigl(F)\otimes\iota_{\dagger}\nbigo_{\nbigy}$.
They give the following isomorphism
of integrable pure twistor $\nbigd$-modules
with real structure:
\[
 \nbigt_{sm}(F)
 \otimes
 \iota_{\dagger}\nbigu_Y(d_Y,0)
\simeq
 \iota_{\dagger}\nbigt(F_Y)
\]
\end{lem}
\pf
Let us observe that $\varphi_0$ is uniquely extended.
Because 
$\nbige_0$ and $\nbige_1$ generate
$\nbigl(F)\otimes\iota_{\dagger}\nbigo_{\nbigx}$
and 
$\iota_{\dagger}\nbigl(F_Y)$
over $\nbigr_X$,
the uniqueness is clear.
Hence, it is enough to obtain $\varphi$
locally around any point of $Y$.
We may assume that $X$ is the product of
$Y$ and a neighbourhood of $0$ in $\cnum_t$.
We have $Y=\{t=0\}$.
Let $V$ denote the $V$-filtration of
$\nbigl(F)\otimes\iota_{\dagger}\nbigo_{\nbigy}$
and $\iota_{\dagger}\nbigl(F_Y)$.
We may naturally regard 
$V_0\bigl(\nbigl(F)\otimes\iota_{\dagger}\nbigo_{\nbigy}\bigr)$
and 
$V_0\bigl(
 \iota_{\dagger}\nbigl(F_Y)
 \bigr)$
as integrable $\nbigr_Y$-modules,
and the natural isomorphism $\varphi_0$
is an isomorphism of the integrable $\nbigr_Y$-modules
$V_0\bigl(\nbigl(F)\otimes\iota_{\dagger}\nbigo_{\nbigy}\bigr)
\simeq
 V_0\bigl(
 \iota_{\dagger}\nbigl(F_Y)
 \bigr)$.
For both $\nbigrtilde_X$-modules,
we have
$V_m=\bigoplus_{j\leq m}\deldel_t^jV_0$.
We can easily construct the isomorphism $\varphi$
inductively on the $V$-filtration.

Because $\varphi_0$ is compatible
with the sesqui-linear pairings on $V_0$,
we obtain that $\varphi$ gives
an isomorphism of the integrable $\nbigr_X$-triples.
The comparison of the real structures
is easy by construction.
\hfill\qed

\begin{prop}
\label{prop;14.11.8.100}
We have the following isomorphisms
\[
 \nbigt(F)[\ast Y]/\nbigt(F)
\simeq
 \iota_{\dagger}\nbigt(F_Y)\otimes\newTate(-1)
\]
\[
 \Ker\bigl(
 \nbigt(F)[! Y]\lrarr\nbigt(F)
 \bigr)
\simeq
 \iota_{\dagger}\nbigt(F_Y) 
\]
\end{prop}
\pf
It follows from (\ref{eq;15.11.10.40}),
(\ref{eq;15.11.10.41})
and Lemma \ref{lem;15.11.10.30}.
\hfill\qed

\vspace{.1in}

Let us describe the isomorphisms
of the underlying integrable $\nbigr_X$-modules
more explicitly.
For simplicity, we assume that 
$X$ is the product of $Y$ and
an open subset in $\cnum_t$.
Let $\upsilon$ be the frame of $\nbigl(F)$
such that $\nabla\upsilon=\upsilon d(F/\lambda)$.
Let $\upsilon_Y$ be the frame of $\nbigl(F_Y)$
obtained as the restriction of $\upsilon$ to $Y$.
The morphism
$\nbigl(F)[\ast t]\lrarr
 \iota_{\dagger}\nbigl(F_Y)$ is given by
\begin{equation}
 \label{eq;15.11.10.100}
 t^{-1}\upsilon \longmapsto
 \iota_{\ast}(\upsilon_Y \lambda^{-1})(dt/\lambda)^{-1}
=\iota_{\ast}(\upsilon_Y(dt)^{-1}).
\end{equation}
Note that we have a natural inclusion
$\eta:\nbigl(F)\lrarr 
 \nbigl(F)[!t]$
as $\nbigo_{\nbigx}$-module.
The morphism
$\iota_{\dagger}\nbigl(F_Y)
\lrarr
 \nbigl(F)[!t]$
is given by
\begin{equation}
 \label{eq;15.11.10.101}
 \iota_{\ast}(\upsilon_Y(dt/\lambda)^{-1})
\longmapsto
-\deldel_t\eta(\upsilon)
+\eta(\del_t(F)\cdot\upsilon)
\end{equation}

\subsection{Push-forward}

We use the notation in \S\ref{subsection;14.10.14.20}.
We have the twistor $\nbigd$-modules
$\nbigt_{\star}\bigl(F,D_Y^{(0)}\bigr)$ 
and 
$\nbigt_{\star}\bigl(F,D_Y^{(1)}\bigr)$ 
$(\star=\ast,!)$.
\begin{prop}
\label{prop;14.11.7.21}
For $\nbigt=
 \nbigt_{\ast}(F,D_Y^{(1)}),
 \nbigt_{\ast}(F,D_Y^{(0)}),
  \nbigt_{\ast}(F,D_Y^{(1)})\big/
 \nbigt_{\ast}(F,D_Y^{(0)})$,
we have
$\pi_{\dagger}^i(\nbigt)=0$ $(i\neq 0)$.
We also have the following natural isomorphisms
of the integrable mixed twistor $\nbigd$-modules with real structure:
\begin{equation}
\label{eq;14.11.29.300}
 \pi_{\dagger}^0\nbigt_{\ast}(F,D_Y^{(1)})
\simeq
 \Bigl(
\nbigt_{\ast}(g,D)\otimes \newTate(-1)
 \Bigr)
 [!(f)_0][\ast D]
\end{equation}
\begin{equation}
\label{eq;14.11.29.301}
 \pi_{\dagger}^0\bigl(
 \nbigt_{\ast}(F,D_Y^{(1)})
 \big/
 \nbigt_{\ast}(F,D_Y^{(0)})
 \bigr)
\simeq
 \nbigt_{\ast}(g,D)\otimes \newTate(-1)
\simeq
 \nbigt_{\ast}(g,D)\otimes \newTate(-1)
 [\ast D]
\end{equation}
\begin{equation}
\label{eq;14.11.29.302}
 \pi_{\dagger}^0\bigl(
 \nbigt_{\ast}(F,D_Y^{(0)})
 \bigr)
 \simeq
 \Ker\left(
  \Bigl(
\nbigt_{\ast}(g,D)\otimes \newTate(-1)
 \Bigr)
 [!(f)_0][\ast D]
\lrarr
 \nbigt_{\ast}(g,D)\otimes \newTate(-1)
 [\ast D]
\right)
\end{equation}
\end{prop}
\pf
The vanishing
$\pi_{\dagger}^i\nbigt=0$ $(i\neq 0)$
follows from Proposition \ref{prop;14.7.1.1}
and the twistor property.
We set $D^{(2)}:D\times\proj^1$.
Let $\iota_0:X\lrarr Y$ be induced by
$\{0\}\lrarr \proj^1$.
By the argument for Proposition \ref{prop;14.11.8.100},
we have a natural isomorphism of 
integrable $\nbigr_Y(\ast D^{(2)})$-triples:
\[
 \bigl(
 \nbigt_{\ast}(F,D_Y^{(1)})
 \big/\nbigt_{\ast}(F,D_Y^{(0)})
 \bigr)
(\ast D^{(2)})
\simeq
 \iota_{0\dagger}
 \nbigt_{\ast}(g,D)(\ast D)
\]
Hence, we have
\[
 \pi_{\dagger}^0
  \bigl(
 \nbigt_{\ast}(F,D_Y^{(1)})
 \big/\nbigt_{\ast}(F,D_Y^{(0)})
 \bigr)
(\ast D)
\simeq
 \nbigt_{\ast}(g,D)(\ast D)
\]
Because we naturally have
$\bigl(
 \nbigt_{\ast}(F,D_Y^{(1)})
 \big/\nbigt_{\ast}(F,D_Y^{(0)})
 \bigr)
 [\ast D^{(2)}]
\simeq
 \nbigt_{\ast}(F,D_Y^{(1)})
 \big/\nbigt_{\ast}(F,D_Y^{(0)})$
and 
$\nbigt_{\ast}(g,D)[\ast D]
\simeq
 \nbigt_{\ast}(g,D)$,
we obtain the isomorphism
(\ref{eq;14.11.29.301})
of integrable mixed twistor $\nbigd$-modules.
The comparison of the real structures is reduced to
Proposition \ref{prop;14.11.8.100}.

We have the natural morphisms of
integrable mixed twistor $\nbigd$-modules
with real structure:
\begin{equation}
\label{eq;14.7.17.10}
\pi_{\dagger}^0\nbigt_{\ast}\bigl(F,D_Y^{(1)}\bigr)
\lrarr
\pi_{\dagger}^0\nbigt_{\ast}\bigl(F,D_Y^{(1)}\bigr)[\ast D]
\llarr
 \Bigl(
 \pi_{\dagger}^0\nbigt_{\ast}\bigl(F,D_Y^{(1)}\bigr)
[!(f)_0]
\Bigr)
[\ast D].
\end{equation}
According to Proposition \ref{prop;14.7.1.1},
the induced morphisms of the underlying
$\nbigd$-modules are isomorphisms.
Hence,
we obtain that the morphisms in (\ref{eq;14.7.17.10})
are isomorphisms.
We have
\[
 \pi_{\dagger}^0\nbigt_{\ast}(F,D_Y^{(1)})
 \bigl(\ast (f)_0\bigr)(\ast D)
\simeq
 \pi_{\dagger}^0\Bigl(
 \nbigt_{\ast}(F,D_Y^{(1)})
 \big/
 \nbigt_{\ast}(F,D_Y^{(0)})
 \Bigr)
 \bigl(\ast (f)_0\bigr)(\ast D).
\]
Hence, we have
\[
  \pi_{\dagger}^0\nbigt_{\ast}(F,D_Y^{(1)})
\simeq
 \Bigl(
 \pi_{\dagger}^0 \bigl(
 \nbigt_{\ast}(F,D_Y^{(1)})
 \big/\nbigt_{\ast}(F,D_Y^{(0)})
 \bigr)
(\ast D)
 \Bigr)[!(f)_0][\ast D]
\simeq
 \Bigl(
 \nbigt_{\ast}(g,D)\otimes\newTate(-1)
 \Bigr)[!(f)_0][\ast D]
\]
Thus, we obtain (\ref{eq;14.11.29.300}).
We obtain (\ref{eq;14.11.29.302})
from the others.
\hfill\qed

\begin{cor}
\label{cor;14.11.7.22}
For $\nbigt=
 \nbigt_{!}(F,D_Y^{(1)}),
 \nbigt_{!}(F,D_Y^{(0)}),
\Ker\bigl(
  \nbigt_{!}(F,D_Y^{(1)})\lrarr
 \nbigt_{!}(F,D_Y^{(0)})
 \bigr)$,
we have
$\pi_{\dagger}^i(\nbigt)=0$ $(i\neq 0)$.
We also have the following natural isomorphisms
of the integrable mixed twistor $\nbigd$-modules with real structure:
\[
 \pi_{\dagger}^0\nbigt_{!}(F,D_Y^{(1)})
\simeq
\nbigt_{!}(g,D)
 [\ast(f)_0][!D]
\]
\[
 \pi_{\dagger}^0\Bigl(
\Ker\Bigl(
 \nbigt_{!}(F,D_Y^{(1)})
\lrarr
 \nbigt_{!}(F,D_Y^{(0)})
 \Bigr)
\Bigr)
\simeq
 \nbigt_{!}(g,D)
\simeq
 \nbigt_{!}(g,D)[!D]
\]
\[
 \pi_{\dagger}^0\bigl(
 \nbigt_{!}(F,D_Y^{(0)})
 \bigr)
 \simeq
 \Cok\Bigl(
 \nbigt_{!}(g,D)[!D]
\lrarr
\nbigt_{!}(g,D)
 [\ast (f)_0][!D]
\Bigr)
\]
\end{cor}
\pf
Because
$\nbigt_{\ast}(F,D_Y^{(a)})^{\ast}
\simeq
 \nbigt_!(F,D_Y^{(a)})
 \otimes\newTate(\dim Y)$,
the claim follows from Proposition 
\ref{prop;14.11.7.21}.
\hfill\qed

\vspace{.1in}
Let us consider the case that $f$ is moreover
non-degenerate along $D$.
In this case $Z_f$ is smooth,
and $Z_f\cup D$ is normal crossing.
As in \S\ref{subsection;14.10.14.20},
let $\iota:Z_f\lrarr X$ denote the inclusion,
and we set $D_{Z_f}:=D\cap Z_f$
and $g_0:=g_{|Z_f}$.
We have the integrable mixed twistor $\nbigd$-modules
$\nbigt_{\star}(g_0,D_{Z_f})$ with real structure on $Z_f$.
By using Proposition \ref{prop;14.11.8.100},
we obtain the following corollary.
\begin{cor}
If $f$ is moreover non-degenerate along $D$,
we have the following isomorphisms
of the integrable mixed twistor $\nbigd$-modules
with real structure:
\begin{equation}
\label{eq;14.11.29.200}
\pi^0_{\dagger}\nbigt_{\ast}(F,D_Y^{(0)})
\simeq
\iota_{\dagger} \nbigt_{\ast}(g_0,D_{Z_f})\otimes\newTate(-1)
\end{equation}
\begin{equation}
\label{eq;14.11.29.201}
\pi^0_{\dagger}\nbigt_!(F,D_Y^{(0)})
\simeq
\iota_{\dagger} \nbigt_{!}(g_0,D_{Z_f})\otimes\newTate(-1)
\end{equation}
The image of the morphism 
$\pi^0_{\dagger}\nbigt_!(F,D_Y^{(0)})
\lrarr
 \pi^0_{\dagger}\nbigt_{\ast}(F,D_Y^{(0)})$
is 
$\iota_{\dagger}\nbigt(g_0)\otimes\newTate(-1)$
under the isomorphisms.
\hfill\qed
\end{cor}

\subsection{Specialization}

Let $X$ be a complex manifold
with a simple normal crossing hypersurface $D$.
We set $X^{(1)}:=X\times\cnum_{\tau}$
and $D^{(1)}:=D\times\cnum_{\tau}$.
Let $f$ and $g$ be meromorphic functions on $(X,D)$.
We have the meromorphic function
$F=\tau f+g$ on $(X^{(1)},D^{(1)})$.
We have the associated integrable mixed twistor $\nbigd$-modules
$\nbigt_{\star}(F,D^{(1)})$ $(\star=\ast,!)$
with real structure on $X^{(1)}$.
Let $\nbigk_{\star,f,g}$
and $\nbigc_{\star,f,g}$
denote the kernel and the cokernel of
$\nbigt_{\star}(F,D^{(1)})[!\tau]
\lrarr
 \nbigt_{\star}(F,D^{(1)})[\ast\tau]$.
Let $\iota_0:X\lrarr X^{(1)}$ be given by
$\iota_0(Q)=(Q,0)$.

\begin{prop}
\label{prop;14.10.18.11}
If $|(f)_0|\cap|(f)_{\infty}|=\emptyset$,
we have the following isomorphism of
the integrable mixed twistor $\nbigd$-modules
with real structure:
\begin{equation}
 \label{eq;15.10.21.10}
 \nbigc_{\ast,f,g}\simeq
 \iota_{0\dagger}
 \nbigt_{\ast}(g,D)
 \otimes\newTate(-1),
\quad\quad
 \nbigk_{\ast,f,g}\simeq
 \iota_{0\dagger}
 \nbigt_{\ast}(g,D)[!(f)_{\infty}],
\end{equation}
\begin{equation}
 \label{eq;15.10.21.11}
 \nbigk_{!,f,g}\simeq
  \iota_{0\dagger}
 \nbigt_!(g,D),
\quad\quad
 \nbigc_{!,f,g}\simeq
  \iota_{0\dagger}
 \nbigt_!(g,D)[\ast(f)_{\infty}]
\otimes\newTate(-1).
\end{equation}
\end{prop}
\pf
By Proposition \ref{prop;14.11.8.100},
if $D=\emptyset$,
we already have the isomorphisms 
of integrable mixed twistor $D$-modules with real structure
in (\ref{eq;15.10.21.10}, \ref{eq;15.10.21.11}).
In other words,
we already have the isomorphisms
of integrable mixed twistor $D$-modules with real structure
on $X^{(1)}\setminus D^{(1)}$.
It is enough to prove that they are extended
isomorphisms 
of integrable mixed twistor $D$-modules with real structure
on $X^{(1)}$.
By construction of the morphisms 
in Proposition \ref{prop;14.11.8.100},
we have the following isomorphisms
of integrable filtered $\nbigr$-triples:
\[
\Bigl(
 \nbigc_{\ast,f,g}
 \Bigr)(\ast D^{(1)})
\stackrel{a_1}{\simeq}
 \Bigl(
 \iota_{0\dagger}
 \nbigt_{\ast}(g,D)
 \otimes\newTate(-1)
\Bigr)(\ast D^{(1)}),
\quad\quad
 \Bigl(
 \nbigk_{\ast,f,g}
 \Bigr)(\ast D^{(1)})
\stackrel{a_2}
\simeq
 \Bigl(
 \iota_{0\dagger}
 \nbigt_{\ast}(g,D)[!(f)_{\infty}]
 \Bigr)(\ast D^{(1)}),
\]
We have the following morphism
of integrable mixed twistor $D$-modules
induced by $a_1$:
\begin{equation}
\label{eq;15.10.21.20}
 \nbigc_{\ast,f,g}
\stackrel{b_1}\lrarr
 \iota_{0\dagger}
 \nbigt_{\ast}(g,D)
 \otimes\newTate(-1)
\end{equation}
By Proposition \ref{prop;14.6.28.10},
the morphism of the underlying $D$-modules
is an isomorphism.
Hence, (\ref{eq;15.10.21.20})
is an isomorphism.
We have the following morphisms:
\begin{equation}
\label{eq;15.10.21.21}
 \nbigk_{\ast,f,g}
 \stackrel{c_1}{\llarr}
 \nbigk_{\ast,f,g}[!(f)_{\infty}]
 \stackrel{c_2}{\lrarr}
 \nbigk_{\ast,f,g}[\ast D^{(1)}][!(f_{\infty})]
 \stackrel{c_3}{\lrarr}
 \iota_{0\dagger}\nbigt_{\ast}(g,D)[!D^{(1)}][!(f)_{\infty}]
\end{equation}
Here, $c_i$ $(i=1,2)$ are canonical morphisms,
and $c_3$ is induced by $a_2$.
By Proposition \ref{prop;14.6.28.10},
the morphisms of the underlying $D$-modules
are isomorphisms.
Hence, we obtain the isomorphism
$\nbigk_{\ast,f,g}
\simeq
 \iota_{0\dagger}\nbigt_{\ast}(g,D)[!(f)_{\infty}]$
from (\ref{eq;15.10.21.21}).
We obtain the isomorphisms in (\ref{eq;15.10.21.11}) similarly.
\hfill\qed

\subsection{$\cnum^{\ast}$-homogeneous $\nbigrtilde$-modules}

\subsubsection{Homogeneity}
\label{subsection;14.11.25.10}

Suppose that a complex manifold $X$ is equipped with
a $\cnum^{\ast}$-action,
i.e., a morphism
$\mu:\cnum^{\ast}\times X\lrarr X$
satisfying
$\mu(a_1a_2,x)=\mu(a_1,\mu(a_2,x))$
and $\mu(1,x)=x$.
We consider the action of 
$\cnum^{\ast}$ on $\cnum_{\lambda}$
given by the multiplication.
Set $\nbigx:=\cnum_{\lambda}\times X$.
We have the diagonal $\cnum^{\ast}$-action
$\mutilde$ on $\nbigx$.

Let $\nbigm$ be an $\nbigrtilde_X$-module.
It is equivalent to 
an $\nbigo_{\cnum_{\lambda}\times X}$-module $\nbigm$
with a meromorphic flat connection
$\nabla:\nbigm\lrarr
 \nbigm\otimes\Omega^1_{\nbigx}(\log \nbigx^0)
 \otimes\nbigo_{\nbigx}(\nbigx^0)$,
where $\nbigx^0:=\{0\}\times X$.
We have the $\nbigo_{\cnum^{\ast}\times\nbigx}$-module
$\mutilde^{\ast}\nbigm$ on $\cnum^{\ast}\times\nbigx$
equipped with the meromorphic flat connection
$\mutilde^{\ast}\nabla$.
We can easily check that 
\[
 (\mutilde^{\ast}\nabla)\bigl(
 \mutilde^{\ast}\nbigm
 \bigr)
\subset
 \mutilde^{\ast}\nbigm\otimes
 \Omega^1_{\cnum^{\ast}\times\nbigx}
 \bigl(
 \log(\cnum^{\ast}\times\nbigx^0)
 \bigr)
 \otimes\nbigo(\cnum^{\ast}\times\nbigx^0).
\]
Hence, $\mutilde^{\ast}\nbigm$ is naturally
an $\nbigrtilde_{\cnum^{\ast}\times X}$-module.

Let $p:\cnum^{\ast}\times \nbigx
 \lrarr \nbigx$ denote the projection.
Let $p_i:
 \cnum^{\ast}\times \cnum^{\ast}\times \nbigx
\lrarr \nbigx$ $(i=1,2,3)$
be given by
\[
 p_1(a_1,a_2,x)=x,
\quad
 p_2(a_1,a_2,x)=\mutilde(a_2,x),
\quad
 p_3(a_1,a_2,x)=\mutilde(a_1a_2,x).
\]
We have the morphisms 
$p_{23},\id\times\mutilde,\mu_{\cnum^{\ast}}\times\id:
  \cnum^{\ast}\times \cnum^{\ast}\times \nbigx
\lrarr 
 \cnum^{\ast}\times \nbigx$
given by
\[
 p_{23}(a_1,a_2,x)=(a_2,x),
\quad
 (\id\times\mutilde)(a_1,a_2,x)=(a_1,\mutilde(a_2,x)),
\quad
 (\mu_{\cnum^{\ast}}\times\id)(a_1,a_2,x)
=(a_1a_2,x).
\]

\begin{df}
\label{df;14.11.25.1}
An $\nbigrtilde_X$-module $\nbigm$ is called 
$\cnum^{\ast}$-homogeneous
if we have an isomorphism
 of $\nbigrtilde_{\cnum^{\ast}\times X}$-modules
 $\kappa:
 p^{\ast}\nbigm
\simeq
 \mutilde^{\ast}\nbigm$
 satisfying the cocycle condition
 $(\id\times\mutilde)^{\ast}\kappa\circ
 p_{23}^{\ast}\kappa
=(\mu_{\cnum^{\ast}}\times\id)^{\ast}\kappa$.
\hfill\qed
\end{df}

The restriction of $\mutilde$ to 
$\cnum_{\lambda}^{\ast}\times X\subset\nbigx$
is a free action, and the quotient space is $X$.
The condition in Definition \ref{df;14.11.25.1}
implies that there exists a $\nbigd$-module $M$
such that the restriction $\nbigm_{|\cnum^{\ast}\times X}$
is isomorphic to the pull back of $M$
by the quotient map
$\cnum^{\ast}_{\lambda}\times X\lrarr 
 (\cnum^{\ast}_{\lambda}\times X)/\cnum^{\ast}\simeq
 X$.
Indeed, $M$ is given as the specialization of $\nbigm$
to $\{1\}\times X$,
i.e., $M=\Xi_{DR}(\nbigm):=
 \iota_1^{-1}\bigl(\nbigm/(\lambda-1)\nbigm\bigr)$,
where $\iota_1:X\lrarr\nbigx$ is given by
$\iota_1(P)=(1,P)$.

\begin{lem}
Suppose
that $\nbigm$ is a locally free $\nbigo_{\nbigx}$-module,
for simplicity.
Then, the torus action is uniquely determined 
by the connection $\nabla$.
\end{lem}
\pf
Let $\underline{\gminiv}$ be the holomorphic fundamental vector field 
of the action $\mu$ on $X$,
i.e., $\underline{\gminiv}_{|Q}=T_{(1,Q)}\mu_{\ast}(\del/\del a)$
for any $Q\in X$.
The holomorphic fundamental vector field of $\mutilde$ on $\nbigx$
is given by $\lambda\del_{\lambda}+\underline{\gminiv}$.
The $\cnum^{\ast}$-action on $\nbigm$
induces a differential operator
$L:\nbigm\lrarr \nbigm$
such that
$L(fs)=fL(s)
+(\lambda\del_{\lambda}+\underline{\gminiv})f\cdot s$.
On $\cnum_{\lambda}^{\ast}\times X$,
we can easily check that 
$L(s)=\nabla_{\lambda\del_{\lambda}+\underline{\gminiv}}(s)$.
The equality holds on $\cnum_{\lambda}\times X$.
Because the $\cnum^{\ast}$-action is determined by $L$,
it is uniquely determined by the connection $\nabla$.
\hfill\qed

\subsubsection{$\nbigr$-modules associated to 
homogeneous meromorphic functions}

Let $X$ be a complex manifold
with a $\cnum^{\ast}$-action $\mu$ as above.
For any $a\in\cnum^{\ast}$,
the morphism $\mu(a,\bullet):X\lrarr X$
is denoted by $\mu_a$.
Let $D$ be a hypersurface of $X$
such that $\mu_a(D)=D$ for any $a\in\cnum^{\ast}$.
Let $f$ be a meromorphic function on $(X,D)$
such that $\mu_a^{\ast}(f)=a f$ for any $a\in\cnum^{\ast}$.

\begin{prop}
\label{prop;14.11.25.20}
The $\nbigr$-modules
$\nbigl_{\star}(f,D)$ are $\cnum^{\ast}$-homogeneous.
\end{prop}
\pf
First, let us consider 
the $\nbigrtilde_{X}(\ast D)$-module
$\nbigq(f)$
induced by 
$\nbigo_{\nbigx}(\ast D)$
with the meromorphic flat connection
$d+d(\lambda^{-1}f)$.
We set $X_1:=\cnum^{\ast}\times X$
and $D_1:=\cnum^{\ast}\times D$.
We set $f_1:=p^{\ast}(f)$.
We naturally have
$p^{\ast}\nbigo_{\nbigx}(\ast D)
\simeq
 \nbigo_{\nbigx_1}(\ast D_1)
\simeq
 \mutilde^{\ast}\nbigo_{\nbigx}(\ast D)$.
Because
$\mutilde^{\ast}(\lambda^{-1}f)
=p^{\ast}(\lambda^{-1}f)
 =\lambda^{-1}f_1$,
the pull back of the connections are also equal.
Hence, we have
$p^{\ast}\nbigq(f)
\simeq
 \mutilde^{\ast}\nbigq(f)$.
Similarly, we can check the cocycle condition for
$\nbigq(f)$
as in Definition \ref{df;14.11.25.1}.
It is enough to check that we have 
natural isomorphisms
$p^{\ast}\nbigl_{\star}(f)
\simeq
 \mutilde^{\ast}\nbigl_{\star}(f)$
which induce the above isomorphism
$p^{\ast}\nbigq(f)
\simeq \mutilde^{\ast}\nbigq(f)$.

Take $(a_0,\lambda_0,Q_0)\in 
 \nbigx_1=\cnum^{\ast}\times \cnum_{\lambda}\times X$.
Let us check that there exists a unique isomorphism
$p^{\ast}\nbigl_{\star}(f)
\simeq
 \mutilde^{\ast}\nbigl_{\star}(f)$
around $(a_0,\lambda_0,Q_0)$ which induces
$p^{\ast}\nbigq(f)
\simeq \mutilde^{\ast}\nbigq(f)$.
Let $g$ be a holomorphic function defined on a neighbourhood
of $\mu(a_0,Q_0)$ such that $g^{-1}(0)=D$.
We set $g_1:=g\circ\mu$ defined on a neighbourhood of $(a_0,Q_0)$.
We have 
$p^{\ast}\nbigl_{\star}(f)[\star g_1]
=p^{\ast}\nbigl_{\star}(f)$ 
on a neighbourhood of $(\lambda_0,a_0,Q_0)$.
It is enough to prove that
$\mutilde^{\ast}\nbigl_{\star}(f)[\star g_1]
=\mutilde^{\ast}\nbigl_{\star}(f)$ 
on a neighbourhood of $(\lambda_0,a_0,Q_0)$.

We set $Q_1:=\mu(a_0,Q_0)$.
We take a small neighbourhood 
$U=U_{\lambda_0}\times U_{Q_1}\subset \nbigx$
of $(\lambda_0,Q_1)$
such that $g$ is defined on $U_{Q_1}$.
Let $\iota_{g}:U_{Q_1}\lrarr U_{Q_1}\times\cnum_t$
be the graph.
We may assume that 
we have the $V$-filtration
of $\iota_{g\dagger}\nbigl_{\star}(f)$
on $U\times\cnum_t$ as $\nbigr$-modules.
Note that each 
$V_a\bigl(
 \iota_{g\dagger}\nbigl_{\star}(f)
 \bigr)$
are naturally $V\nbigrtilde_{X\times\cnum_t}$-modules.
In the case $\star=\ast$,
the induced morphism
$t:\Gr^V_{0}\lrarr \Gr^{V}_{-1}$ is an isomorphism.
In the case $\star=!$,
the induced morphism
$\deldel_t:\Gr^V_{-1}\lrarr \Gr^{V}_{0}$
is an isomorphism.

Let $U'\subset \nbigx_1$ be a small neighbourhood of
$(a_0,\lambda_0,Q_0)$ such that $\mutilde(U')\subset U$.
Let $\mutilde_1:U'\times\cnum_t\lrarr U\times\cnum_t$
be the induced morphism.
We naturally have
$\iota_{g\dagger}\nbigl_{\star}(f)
=\bigoplus_{j=0}^{\infty}
 \deldel_t^j\iota_{g\ast}\nbigl_{\star}(f)$
and 
$\iota_{g_1\dagger}\mutilde^{\ast}\nbigl_{\star}(f)
=\bigoplus_{j=0}^{\infty}
 \deldel_t^j\iota_{g_1\ast}\mutilde^{\ast}\nbigl_{\star}(f)$.
The natural isomorphisms
$\mutilde_1^{\ast}
 \iota_{g\ast}\nbigl_{\star}(f)
\simeq
 \iota_{g_1\ast}\mutilde^{\ast}\nbigl_{\star}(f)$
induce isomorphisms
$\mutilde_1^{\ast}(\iota_{g\dagger}\nbigl_{\star}(f))
\simeq
 \iota_{g_1\dagger}\mutilde^{\ast}\nbigl_{\star}(f)$.

We set 
$V_a\bigl(
 \iota_{g_1\dagger}\mutilde^{\ast}\nbigl_{\star}(f)
 \bigr)
:= 
 \mutilde_1^{\ast}
 V_a\bigl(
 \iota_{g\dagger}\nbigl_{\star}(f)
 \bigr)$
as $\nbigo_{U'\times\cnum_t}$-modules.
By the construction,
$V_a\bigl(
 \iota_{g_1\dagger}\mutilde^{\ast}\nbigl_{\star}(f)
 \bigr)$
are $V\nbigrtilde_{X_1\times\cnum_t}$-modules
on $U'\times\cnum_t$.
Because 
$V_a\bigl(
 \iota_{g\dagger}\nbigl_{\star}(f)
 \bigr)$ are coherent over
$V\nbigr_{X\times\cnum_t}$,
we obtain that 
$V_a\bigl(
 \iota_{g_1\dagger}\mutilde^{\ast}\nbigl_{\star}(f)
 \bigr)$
are pseudo-coherent over $\nbigo_{U'\times\cnum_t}$,
and we can easily check that
they are finitely generated over 
$V\nbigr_{X_1\times\cnum_t}$.
Hence they are coherent over $V\nbigr_{X_1\times\cnum_t}$.
(See \cite{kashiwara_text}.)
Then, we obtain that 
$V_{\bullet}\bigl(
 \iota_{g_1\dagger}\mutilde^{\ast}\nbigl_{\star}(f)
 \bigr)$
is a $V$-filtration of 
$\iota_{g_1\dagger}\mutilde^{\ast}\nbigl_{\star}(f)$
along $t$.
In the case $\star=\ast$,
the morphism
$t:\Gr^V_0\lrarr \Gr^V_{-1}$ is an isomorphism.
In the case $\star=!$,
the morphism
$\deldel_t:\Gr^V_{-1}\lrarr \Gr^V_{0}$
is an isomorphism.
Hence, we have
$\mutilde^{\ast}\nbigl_{\star}(f)
\simeq
 \nbigq(f_1)[\star g_1]
\simeq
 p^{\ast}\nbigl_{\star}(f)[\star g_1]$.
Thus, the proof of Proposition \ref{prop;14.11.25.20}
is finished.
\hfill\qed

\subsubsection{Good filtrations on holonomic $\nbigd$-modules
and $\cnum^{\ast}$-homogeneity}
\label{subsection;15.11.17.10}

Let $X$ be a complex manifold.
We regard that $X$ is equipped with 
the trivial $\cnum^{\ast}$-action $\mu_0$,
i.e., $\mu_0(a,x)=x$ for any $(a,x)\in \cnum^{\ast}\times X$.
The induced  action $\mutilde_0$ on $\cnum_{\lambda}\times X$
is just the multiplication on the $\cnum_{\lambda}$-component.
In this case,
we have a relation between
$\cnum^{\ast}$-homogeneity and good filtrations.

\vspace{.1in}

Let $\nbigc_1(X)$ denote the category of holonomic
$\nbigd_X$-modules with a good filtration $(M,F)$.
Morphisms $(M_1,F_1)\lrarr (M_2,F_2)$
in $\nbigc_1(X)$ are morphisms $\varphi:M_1\lrarr M_2$
of holonomic $\nbigd_X$-modules
such that $\varphi(F_jM_1)\subset F_jM_2$.

Let $\nbigc_2(X)$ denote the category of
coherent $\nbigrtilde$-modules $\nbigm$
such that 
(i) it is $\cnum^{\ast}$-homogeneous
with respect to $\mu_0$,
(ii) $\nbigm$ is coherent as an $\nbigr_X$-module,
(iii) $\nbigm$ is strict as an $\nbigr_X$-module
i.e.,
the multiplication of $\lambda-\lambda_0$
is a monomorphism for any $\lambda_0\in\cnum$,
(iv) $\Xi_{\DR}(\nbigm)$ is a holonomic $\nbigd_{X}$-module.
Morphisms $\nbigm_1\lrarr\nbigm_2$ in $\nbigc_2(X)$
are morphisms of $\nbigrtilde_X$-modules.

Let $(M,F)$ be an object in $\nbigc_1(X)$.
We have the Rees module
$R(M,F)=\sum F_jM\lambda^{j}$
of $(M,F)$,
which is naturally an $\nbigo_X[\lambda]$-module.
By taking the analytification 
we obtain an $\nbigo_{\cnum_{\lambda}\times X}$-module
$\Rtilde(M,F)$
which is naturally an $\nbigrtilde_X$-module.
This construction gives a functor
$\Rtilde:\nbigc_1(X)\lrarr\nbigc_2(X)$.

\begin{prop}
\label{prop;15.11.17.101}
The functor $\Rtilde$ is an equivalence.
\end{prop}
\pf
Set $\nbigx:=\cnum_{\lambda}\times X$.
Let $p_{\lambda}:\nbigx\lrarr X$ denote the projection.
Let $\Hol(X)$ denote the category of 
holonomic $\nbigd$-modules on $X$.
Let $\nbigc_3(X)$ denote the category of
$\nbigd_{\nbigx}(\ast\lambda)$-modules 
$\nbigl$ satisfying the conditions
(i) $\nbigl$ is $\cnum^{\ast}$-homogeneous
 with respect to $\mu_0$,
(ii) $\nbigl$ is coherent over 
 $p_{\lambda}^{\ast}(\nbigd_X)(\ast\lambda)$,
(iii) $\nbigl$ is strict as an $\nbigr$-module,
(iv) $\Xi_{\DR}(\nbigl)$ is a holonomic $\nbigd_X$-module.
Morphisms $\nbigl_1\lrarr\nbigl_2$ in $\nbigc_3(X)$
are $\cnum^{\ast}$-equivariant morphisms
of $\nbigd_{\nbigx}$-modules.

For any $M\in \Hol(X)$,
we set $G_1(M):=p_{\lambda}^{\ast}(M)(\ast\lambda)$.
Then, we obtain a functor
$G_1:\Hol(X)\lrarr \nbigc_3(X)$.
Clearly $G_1$ is an exact functor.

\begin{lem}
\label{lem;15.11.17.100}
$G_1$ is an equivalence.
\end{lem}
\pf
Clearly $G_1$ is faithful.
Let us prove that $G_1$ is full.
Let $M_i\in\Hol(X)$.
Let $f:G_1(M_1)\lrarr G_1(M_2)$
be a morphism in $\nbigc_3(X)$.
By applying $\Xi_{\DR}$,
we obtain a morphism
$f_1:=\Xi_{\DR}(f):M_1\lrarr M_2$.
Because  the restriction of
$G_1(f_1)-f:G_1(M_1)\lrarr G_1(M_2)$
to $\cnum^{\ast}\times X$ is $0$
because of the $\cnum^{\ast}$-equivariance of $f$.
Hence, we obtain $G_1(f_1)-f=0$
on $\cnum\times X$,
i.e., $G_1$ is full.

Let us prove that $G_1$ is essentially surjective.
Let $\nbigl$ be an object in $\nbigc_3(X)$.
We set $L:=\Xi_{\DR}(\nbigl)$.
It is enough to prove that
$G_1(L)\simeq \nbigl$.
Note that we have a natural isomorphism
$\kappa:G_1(L)_{|\cnum^{\ast}\times X}\simeq
 \nbigl_{|\cnum^{\ast}\times X}$
induced by the $\cnum^{\ast}$-homogeneity.
It is enough to prove that $\kappa$
is extended to an isomorphism on $\cnum\times X$.
Take any point $P\in X$.
It is enough to prove that 
for a small neighbourhood $X_P$ of $P$ in $X$,
$\kappa_{|\cnum^{\ast}\times X_P}$ is 
extended on $\cnum\times X_P$.
We have the subvariety $Z(\nbigl)\subset X_P$
such that
the support of the sheaf $\nbigl_{|\nbigx_P}$ on $\nbigx_P$
is $\cnum\times Z(\nbigl)$.
We use an induction on the dimension of 
$Z(\nbigl)$.
To simplify the description,
we set $Z:=Z(\nbigl)$.

We have a stratification
$Z=\coprod Z_i$ into locally closed smooth subvarieties
such that the characteristic variety of $L$
is the union of the conormal bundle 
$T_{Z_i}^{\ast}X$ of $Z_i$ in $X$.
By shrinking $X_P$,
we may assume to have holomorphic function $g$
such that 
$Z_i\subset g^{-1}(0)$
if and only if $\dim Z_i<\dim Z$.
We may assume to have a complex manifold $Y$
with a projective morphism
$\varphi:Y\lrarr X$
such that 
(i) $\varphi(Y)\subset Z$ and
 $\dim(Z\setminus \varphi(Y))<\dim Z$,
(ii) we set $g_Y:=g\circ\varphi$,
 and then $D_Y:=g_Y^{-1}(0)$ is a normal crossing,
(iii) $\varphi$ induces
 $Y\setminus g_Y^{-1}(0)\simeq Z\setminus g^{-1}(0)$.
We set $\gtilde:=g\circ p_{\lambda}$
and $\gtilde_Y:=g_Y\circ p_{\lambda}$.
We have a meromorphic flat bundle $L_1$ on $(Y,D_Y)$
with an isomorphism
$\varphi_{+}L_1\simeq L(\ast g)$
of $\nbigd_{X_P}$-modules.
We have a $\nbigd_{\nbigy}(\ast \gtilde_Y)$-module
$\nbigl_1$ with an isomorphism
$(\id\times\varphi)_{+}\nbigl_1\simeq
 \nbigl(\ast \gtilde)$.
Note that $\nbigl_1$ is a meromorphic flat bundle
on $\nbigy$ with the pole
$(\{0\}\times Y)\cup \gtilde_Y^{-1}(0)$.

Because $\nbigl(\ast \gtilde)$ is $\cnum^{\ast}$-homogeneous
with respect to the trivial action on $X$,
$\nbigl_1$ is $\cnum^{\ast}$-homogeneous
with respect to the trivial action on $Y$.
It is easy to observe that
$\nbigl_{1|\cnum\times (Y\setminus D_Y)}$
is regular along $\{0\}\times(Y\setminus D_Y)$.
Hence, the isomorphism
$p_{\lambda}^{\ast}(L_1)_{|\cnum^{\ast}\times(Y\setminus D_Y)}
\simeq
 \nbigl_{1|\cnum^{\ast}\times (Y\setminus D_Y)}$
is extended to an isomorphism
$G_1(L_1)_{|\cnum\times(Y\setminus D_Y)}
\simeq
 \nbigl_{1|\cnum\times (Y\setminus D_Y)}$.
By using the Hartogs property,
the isomorphism is extended to
$G_1(L_1)
\simeq
 \nbigl_1$.
We also obtain
$\nbigl_1[!\gtilde_Y]
\simeq
 G_1(L_1)[!\gtilde_Y]
\simeq
 G_1(L_1[!g_Y])$.
Hence, we have a natural isomorphism
$G_1\bigl(L(\star g)\bigr)
\simeq \nbigl(\star \gtilde)$.

Considering the Beilinson functors,
we have
$G_1\bigl(
 \Pi^{a,b}_g(L)(\star g)
 \bigr)
\simeq
 \Pi^{a,b}_{\gtilde}(\nbigl)[\star \gtilde]$.
We obtain
$G_1\bigl(\psi_g^{(a)}(L) \bigr)
\simeq
 \psi_{\gtilde}^{(a)}(\nbigl)$
and 
$G_1\bigl(
 \Xi_g^{(a)}(L) \bigr)
\simeq
 \Xi_{\gtilde}^{(a)}(\nbigl)$.
By the assumption of the induction,
we have 
$G_1\bigl(
 \phi_g^{(a)}(L)\bigr)
\simeq
 \phi_{\gtilde}^{(a)}(\nbigl)$.
The following is commutative:
\[
 \begin{CD}
 G_1\bigl(
 \psi^{(1)}_g(L)
 \bigr)
 @>>>
 G_1\bigl(
 \phi^{(0)}_g(L)
 \bigr)
\oplus
 G_1\bigl(
 \Xi^{(0)}_g(L)
 \bigr)
 @>>>
 G_1\bigl(
 \psi^{(0)}_g(L)
 \bigr)
 \\
 @V{\simeq}VV @V{\simeq}VV @V{\simeq}VV \\
 \psi^{(1)}_{\gtilde}(\nbigl)
 @>>>
 \phi^{(0)}_{\gtilde}(\nbigl)
\oplus
  \Xi^{(0)}_{\gtilde}(\nbigl)
 @>>>
  \psi^{(0)}_{\gtilde}(\nbigl)
 \end{CD}
\]
We obtain
$G_1(L)\simeq \nbigl$.
Thus,
Lemma \ref{lem;15.11.17.100}
is finished.
\hfill\qed

\vspace{.1in}

Let $\nbigc_4(X)$ denote the category
whose objects are the same as those of $\nbigc_3(X)$,
but morphisms $f:\nbigm_1\lrarr\nbigm_2$ in $\nbigc_4(X)$
are morphisms of $\nbigd_{\nbigx}$-modules,
i.e., we do not assume the $\cnum^{\ast}$-equivariance.
\begin{lem}
\label{lem;15.11.17.110}
The natural functor
$\nbigc_3(X)\lrarr \nbigc_4(X)$
is an equivalence.
Namely,
any morphism $f:\nbigm_1\lrarr\nbigm_2$
in $\nbigc_4(X)$ is $\cnum^{\ast}$-equivariant.
\end{lem}
\pf
Let $M_1\lrarr M_2$ be a morphism in $\Hol(X)$.
Let $f:G_1(M_1)\lrarr G_1(M_2)$ be a morphism
in $\nbigc_4(X)$.
Applying $\Xi_{\DR}$,
we obtain a morphism $f_0:M_1\lrarr M_2$ 
of $\nbigd_X$-modules.
It is easy to see that
$\DR_{\cnum^{\ast}\times X}
 G_1(f_0)_{|\cnum^{\ast}\times X}=
 \DR_{\cnum^{\ast}\times X}f_{|\cnum^{\ast}\times X}$.
Hence, we obtain
$G_1(f_0)_{|\cnum^{\ast}\times X}=
  f_{|\cnum^{\ast}\times X}$,
and $G_1(f_0)=f$ on $\nbigx$.
\hfill\qed

\vspace{.1in}

Let us return to the proof of Proposition \ref{prop;15.11.17.101}.
Let $\nbigm$ be an object in $\nbigc_2(X)$.
Set $M:=\Xi_{\DR}(\nbigm)$.
For any local section $s$ of $M$,
we have the section $p_{\lambda}^{\ast}(s)$
of $\nbigm(\ast\lambda)\simeq G_1(M)$
by Lemma \ref{lem;15.11.17.110}.
The number 
$i(s):=\min\big\{i\,\big|\,
 \lambda^{i}p_{\lambda}^{\ast}(s)\in\nbigm
 \}$
exists because $\nbigm$ is a coherent $\nbigr_X$-module.
It determines the filtration $F$ on $M$,
i.e.,
$F_j(M)=\{s\in M\,|\,i(s)\leq j\}$.
We have 
$F_j\nbigd_X \cdot
 F_k(M)\subset F_{j+k}(M)$,
where $F_j\nbigd_X$ denote the sheaf of 
differential operators whose orders are less than $j$.

\begin{lem}
\label{lem;14.11.25.30}
We have a natural isomorphism
$\nbigrtilde$-modules
$\Rtilde(M,F)\simeq
 \nbigm$.
\end{lem}
\pf
By Lemma \ref{lem;15.11.17.100},
we naturally have
$\Rtilde(M,F)(\ast\lambda)
=G_1(M)
\simeq
 \nbigm(\ast\lambda)$.
By the construction of $F$,
we have
$\Rtilde(M,F)
\subset
 \nbigm$.
Let $\overline{p}_{\lambda}:
 \proj^1_{\lambda}\times X\lrarr X$
be the projection.
Let $\nbigm'$ be the sheaf on
$\proj^1_{\lambda}\times X$
such that
(i) the restriction to 
$\bigl(\proj^1_{\lambda}
 \setminus\{0\}\bigr)\times X$
is equal to the restriction of
$\overline{p}_{\lambda}^{\ast}(M)(\ast\infty)$,
(ii) the restriction to
 $\cnum_{\lambda} \times X$
is equal to $\nbigm$.
Because $\nbigm$ is $\cnum^{\ast}$-equivariant,
$\nbigm'$ is also $\cnum^{\ast}$-equivariant.
It is easy to observe that
$R^i\overline{p}_{\lambda\ast}
 \nbigm'=0$ for $i>0$,
and 
that the natural morphism
$\overline{p}_{\lambda\ast}\nbigm'
\lrarr
 \nbigm_{|\{0\}\times X}$
is an epimorphism.
For any local section 
$s$ of $\nbigm/\lambda\nbigm$,
we take a section $\stilde$ of 
$\overline{p}_{\lambda\ast}\nbigm'$
which is mapped to $s$.
Because 
$\nbigm'\subset
 \overline{p}_{\lambda}^{\ast}(M)
 \bigl(\ast(\{0,\infty\}\times X)\bigr)$,
we have
$\overline{p}_{\lambda\ast}\nbigm'
\subset
 M[\lambda,\lambda^{-1}]$.
Hence, we have
$\stilde=\sum_{j=-N}^N\lambda^j\stilde_j$
for a large $N$,
where $\stilde_j$ are sections of $M$.
By the $\cnum^{\ast}$-equivariance of $\nbigm$,
we obtain that each $\lambda^j\stilde_j$
is a section of $\nbigm$.
By construction,
$\lambda^j\stilde_j$
are sections of $\Rtilde(M,F)$.
It implies that
$\Rtilde(M,F)\lrarr \nbigm/\lambda\nbigm$
is an epimorphism,
and hence we obtain
$\Rtilde(M,F)= \nbigm$.
\hfill\qed

\vspace{.1in}
Let $\iota_0:X\lrarr \nbigx$ be
given by $\iota_0(P)=(0,P)$.
Because 
$\Gr^F(M)\simeq
 \iota_0^{-1}\bigl(
 \nbigm/\lambda\nbigm\bigr)$,
the $\Sym^{\bullet}\Theta_X$-module
$\Gr^F(M)$ is coherent.
Hence, $F$ is a good filtration,
i.e.,
$(M,F)$ is an object in $\nbigc_1(X)$.
We obtain that 
$\Rtilde$ is essentially surjective.

\vspace{.1in}
Clearly $\Rtilde$ is faithful.
Let us prove that 
$\Rtilde$ is full.
Let $(M^{(i)},F^{(i)})$ $(i=1,2)$
be objects in $\nbigc_1(X)$.
Let $f:\Rtilde(M^{(1)},F^{(1)})\lrarr \Rtilde(M^{(2)},F^{(2)})$
be a morphism in $\nbigc_2(X)$.
By applying $\Xi_{\DR}$,
we have the morphism $f_1:M^{(1)}\lrarr M^{(2)}$
of $\nbigd_{X}$-modules.
By Lemma \ref{lem;15.11.17.110},
$f$ is the restriction of 
$G_1(f_1)$ to $\Rtilde(M^{(1)},F^{(1)})$.
Then, it is easy to deduce that
$f_1$ preserves the filtrations,
i.e., $f_1$ gives a morphism
$(M^{(1)},F^{(1)})\lrarr (M^{(2)},F^{(2)})$
in $\nbigc_1(X)$,
and $f=\Rtilde(f_1)$.
Thus, Proposition \ref{prop;15.11.17.101}
is proved.
\hfill\qed

\vspace{.1in}
Let us give some complements.
Let $\nbigc_2'(X)$ denote the category 
whose objects are the same as those of
$\nbigc_2(X)$,
and whose morphisms
$\nbigm_1\lrarr\nbigm_2$
are $\cnum^{\ast}$-equivariant $\nbigrtilde_X$-homomorphisms.

\begin{cor}
\label{cor;15.11.18.2}
The natural functor
$\nbigc_2'(X)\lrarr \nbigc_2(X)$
is an equivalence,
i.e.,
any morphism in $\nbigc_2(X)$
is $\cnum^{\ast}$-equivariant.
\end{cor}
\pf
It follows from Lemma \ref{lem;15.11.17.110}.
It also follows from the equivalence of $\Rtilde$
and the $\cnum^{\ast}$-equivariance
of any morphism
$\Rtilde(f):\Rtilde(M^{(1)},F^{(1)})\lrarr\Rtilde(M^{(2)},F^{(2)})$
induced by $f:(M^{(1)},F^{(1)})\lrarr(M^{(2)},F^{(2)})$
in $\nbigc_1(X)$.
\hfill\qed

\subsubsection{$\cnum^{\ast}$-homogeneity and specialization
(Appendix)}

Let $X$ be a complex manifold 
with a $\cnum^{\ast}$-action $\mu$.
Let $\nbigm$ be a coherent and strict $\nbigrtilde_X$-module
which is $\cnum^{\ast}$-homogeneous with respect to $\mu$.
Let $f$ be a holomorphic function on $X$
such that there exists an integer $m$
such that $f(\mu(a,x))=a^mx$  for any $x\in X$.
Suppose that $\nbigm$ is strictly specializable along $f$.
Recall that the $V$-filtration of 
any $\nbigrtilde_X$-modules
are indexed by $\real\times\{0\}$,
which is naturally identified with $\real$.
So, we have the $\nbigrtilde_X$-modules
$\psi_{t,u}(\nbigm)$.
They are not necessarily $\cnum^{\ast}$-homogeneous.
Let us observe that 
we can obtain homogeneous $\nbigrtilde$-modules
by modifying $\psi_{f,b}(\nbigm)$ $(b\in\real)$.
By considering the graph construction of $f$,
it is enough to consider the case where
$X=X_0\times\cnum_t$, $f=t$
and the action $\mu$ is given as
the diagonal action of 
the $\cnum^{\ast}$-action $\mu_0$ on $X_0$
and the $\cnum^{\ast}$-action on $\cnum_t$
given by $a\bullet t=a^mt$,
and to study the $\nbigrtilde_{X_0}$-modules
$\psi_{t,b}(\nbigm)$.

Let $\nabla$ denote the meromorphic flat connection
of $\nbigm$ induced by the $\nbigrtilde_X$-module structure,
i.e.,
$\nabla:
 \nbigm\lrarr
 \nbigm\otimes
 \lambda^{-1}\Omega_{\nbigx}(\log \nbigx^0)$.
Let $\mutilde:\cnum^{\ast}\times \nbigx\lrarr\nbigx$
be given by
$\mutilde(a,\lambda,x,t)=(a\lambda,\mu_0(a,x),a^mt)$.
Let $p_2:\cnum^{\ast}\times\nbigx\lrarr\nbigx$
be the projection.
We have the isomorphism
$p_2^{\ast}\nbigm\simeq \mutilde^{\ast}\nbigm$
satisfying the cocycle condition
under which
$p_2^{\ast}(\nabla)=\mutilde^{\ast}\nabla$.

Let $V_{\bullet}$ denote the $V$-filtration of $\nbigm$ along $t$.
Clearly $p^{\ast}(V_{\bullet}\nbigm)$
is the $V$-filtration of $p^{\ast}\nbigm$ along $t$.
It is easy to check that 
$\mutilde^{\ast}(V_{\bullet}\nbigm)$
is also the $V$-filtration of $\mutilde^{\ast}\nbigm=p^{\ast}\nbigm$
along $t$.
We have
$p^{\ast}V_{\bullet}=
 \mutilde^{\ast}(V_{\bullet})$.
Hence,
the $\nbigo_{\nbigx_0}$-modules
$\psi_{t,b}(\nbigm):=V_b\nbigm/V_{<b}\nbigm$
are $\mutilde_0^{\ast}$-equivariant,
i.e.,
we have the isomorphism
$\mutilde_0^{\ast}\psi_{t,b}(\nbigm)
\simeq
 p_{0,2}^{\ast}\psi_{t,b}(\nbigm)$
satisfying the cocycle condition,
where $p_{0,2}:\cnum^{\ast}\times\nbigx_0\lrarr\nbigx_0$
denotes the projection.
Note that 
it is not necessarily $\cnum^{\ast}$-homogeneous
as an $\nbigrtilde_{X_0}$-module.
We have
$\nabla:
 V_b\nbigm\lrarr
 V_b\nbigm\otimes
 \lambda^{-1}
 \Omega^1_{\nbigx}\bigl(\log(\nbigx^0\cup\nbigx_0)\bigr)$,
which induces
\[
 G^{(b)}:
 \psi_{t,b}(\nbigm)
\lrarr
 \psi_{t,b}(\nbigm)
 \otimes
 \Bigl(
 \lambda^{-1}
 \Omega^1_{\nbigx_0}(\log \nbigx_0^0)
\oplus
 \lambda^{-1}\nbigo_{\nbigx_0}dt/t
 \Bigr)
\]
We have the decomposition
$G^{(b)}=\nabla^{(b)}+\nbigb^{(b)}dt/t$,
where $\nabla^{(b)}$ is a meromorphic flat connection 
of the $\nbigrtilde_{X_0}$-module $\psi_{t,b}(\nbigm)$,
and $\nbigb^{(b)}$ is an $\nbigrtilde_{X_0}$-endomorphism
of $\psi_{t,b}(\nbigm)$.
By the relation
$\mutilde^{\ast}\nabla=p_2^{\ast}\nabla$,
we obtain
\[
 \mutilde_0^{\ast}(\nabla^{(b)})
+\mutilde^{\ast}\nbigb^{(b)}a^{-1}da
=p_{0,2}^{\ast}(\nabla^{(b)}),
\quad\quad
\mutilde_0^{\ast}\nbigb^{(b)}
=p_{0,2}^{\ast}(\nbigb^{(b)})
\]
Then, the meromorphic flat connection
$\nabla^{(b)}+m\nbigb^{(b)}d\lambda/\lambda$
gives an $\nbigrtilde_{X_0}$-module structure on $\psi_{t,b}(\nbigm)$
which is $\cnum^{\ast}$-homogeneous.

\begin{rem}
Note that $\nbigb^{(-1)}$ on $\psi_{t,-1}(\nbigm)$
is nilpotent,
and the induced action on
$\Cok(\nbigb^{(-1)})$
and $\Ker(\nbigb^{(-1)})$ are $0$.
Hence, 
$\Cok(\nbigb^{(-1)})$
and $\Ker(\nbigb^{(-1)})$ are 
naturally $\cnum^{\ast}$-homogeneous.
It also follows from the fact that
$\nbigm[\star t]$ $(\star=\ast,!)$ 
are $\cnum^{\ast}$-homogeneous,
and the kernel and the cokernel
of $\nbigm[!t]\lrarr\nbigm[\ast t]$
are also $\cnum^{\ast}$-homogeneous.
\hfill\qed
\end{rem}

\subsection{Hodge modules}
\label{subsection;15.11.18.1}

\subsubsection{Hodge modules and $\cnum^{\ast}$-homogeneity}

Let $\MHM(X,\real)$ denote the category of
graded polarizable mixed $\real$-Hodge modules on $X$.
Let $\MTMint_{\rnum\times\{0\}}(X,\real)$ denote the category of
integrable mixed twistor $\nbigd$-modules with real structure
on $X$
whose KMS-spectrum are contained in $\rnum\times\{0\}$ .
As explained in \cite{Mochizuki-MTM},
we have the fully faithful functor 
$\Psi_X:\MHM(X,\real)\lrarr
 \MTMint_{\rnum\times\{0\}}(X,\real)$.
It is exact.
It is compatible with the other operations
such as the duality and the push-forward.
Let us identify the essential image
by the $\cnum^{\ast}$-homogeneity condition.
Because we consider only graded polarizable mixed Hodge modules,
we omit to distinguish ``graded polarizable''.
We shall prove the following theorem
in \S\ref{subsection;15.11.17.120}--\ref{subsection;15.11.17.121}.

\begin{thm}
\label{thm;15.11.17.20}
Let $(\nbigt,W)$ be an object in
$\MTMint_{\rnum\times\{0\}}(X,\real)$.
Then, 
$(\nbigt,W)$ is contained in
the essential image of $\Psi_X$
if and only if 
the underlying $\nbigrtilde_X$-modules
$W_k\nbigt$ $(k\in\seisuu)$
are $\cnum^{\ast}$-homogeneous 
with respect to the trivial $\cnum^{\ast}$-action
on $X$.
\end{thm}

\subsubsection{Smooth case}
\label{subsection;15.11.17.120}

Let $\nbigt\in\MTMint_{\rnum\times\{0\}}(X,\real)$.
Suppose that the underlying $\nbigrtilde_X$-modules
are $\cnum^{\ast}$-homogeneous
with respect to the trivial action.
Let us consider the case where
the underlying $\nbigrtilde_X$-modules  of $\nbigt$
are smooth,
i.e., locally free $\nbigo_{\cnum\times X}$-modules.
We set
$\nbigt_0=(\nbigv_1,\nbigv_2,C_{\nbigv})
:=\nbigt\otimes\nbigu_X(-d_X,0)$.
It is naturally a variation of integrable twistor structure.
It is also equipped with the real structure
as a variation of integrable twistor structure,
i.e.,
an isomorphism 
$\kappa:
\gammatilde^{\ast}_{sm}(\nbigt_0)
\simeq
\nbigt_0$
such that 
$\gammatilde_{sm}^{\ast}\kappa\circ\kappa=\id$,
where 
$\gammatilde_{sm}^{\ast}(\nbigt_0):=
 (j^{\ast}\nbigv_2^{\lor},j^{\ast}\nbigv_1^{\lor},C_{\nbigv}^{\lor})$.
(Note that the notion of real structure 
for variations of integrable twistor structure
is not equal to the notion of real structure
for integrable mixed twistor $\nbigd$-modules.)

Recall that we have 
the $C^{\infty}$-vector bundle
on $\proj^1\times X$ with some differential operator
$\DD^{\sankaku}$
associated to the smooth $\nbigr_X$-triple $\nbigt_0$.
We have the $\nbigr_X$-module
$\nbige:=\nbigv_1^{\lor}$
on $\cnum_{\lambda}\times X$
which is a locally free 
$\nbigo_{\cnum_{\lambda}\times X}$-module.
Let $X^{\dagger}$ denote the conjugate complex manifold.
Let $\sigma:\cnum_{\mu}\lrarr \cnum_{\lambda}$
be given by $\sigma(\mu)=-\overline{\mu}$.
We have the anti-holomorphic isomorphism
$\sigma:
 \cnum_{\mu}\times X^{\dagger}
 \lrarr
 \cnum_{\lambda}\times X$.
We have the sheaf
$\nbige^{\dagger}:=\sigma^{\ast}\nbigv_2$
which is naturally an $\nbigr_{X^{\dagger}}$-module,
and a locally $\nbigo_{\cnum_{\mu}\times X^{\dagger}}$-module.
Note that  $\nbige$
is naturally a $C^{\infty}$-bundle
with a family of flat connections
\[
 \DD^f:\nbige
\lrarr
 \nbige
 \otimes
\Bigl(
 \lambda^{-1}p_{\lambda}^{\ast}\Omega^{1,0}_X
 \oplus p_{\lambda}^{\ast}\Omega^{0,1}_X
 \Bigr)
\]
induced by the structure of $\nbigr_X$-module,
where $p_{\lambda}:\cnum_{\lambda}\times X\lrarr X$
denote the projection.
The family of flat connections 
is holomorphic with respect to $\lambda$.
Note also that
$\nbige^{\dagger}$
is naturally a $C^{\infty}$-bundle
with a family of flat connections
\[
 \DD^{\dagger\,f}:
 \nbige^{\dagger}_{|\cnum_{\mu}^{\ast}\times X^{\dagger}}\lrarr
  \nbige^{\dagger}_{|\cnum_{\mu}^{\ast}\times X^{\dagger}}
 \otimes 
 \Bigl(
 \mu^{-1}p_{\mu}^{\ast}\Omega^{1,0}_{X^{\dagger}}
 \oplus
 \Omega^{0,1}_{X^{\dagger}}
 \Bigr)
\]
induced by the structure of $\nbigr_{X^{\dagger}}$-module.
The family of flat connections is holomorphic with respect to $\mu$.
We identify $\cnum^{\ast}_{\lambda}$ and $\cnum^{\ast}_{\mu}$
by $\lambda=\mu^{-1}$.
Then, we have the unique isomorphism of  family of flat connections
$\Phi:
 (\nbige,\DD^{f})_{|\cnum_{\lambda}^{\ast}\times X}
\simeq
 (\nbige^{\dagger},\DD^{\dagger\,f})
 _{|\cnum_{\mu}^{\ast}\times X^{\dagger}}$,
such that
(i) it is holomorphic with respect to $\lambda$,
(ii) $\Phi_{|\vecS\times X}$ is 
 compatible with the sesqui-linear pairing $C_{\nbigv}$.
Let $\nbige^{\sankaku}$ be the $C^{\infty}$-bundle
on $\proj^1\times X$
obtained as the gluing of
$\nbige$ and $\nbige^{\dagger}$
by $\Phi$.
It is equipped with the holomorphic structure
in the $\proj^1$-direction $d''_{\lambda}$,
and the differential operator
$\DD^{\sankaku}:
 \nbige^{\sankaku}
\lrarr
 \nbige^{\sankaku}
 \otimes
 \Omega^1_{\proj^1\times X/\proj^1}
 \otimes\nbigo_{\proj^1}(1)$
induced by the $\nbigr$-modules structure on $\nbige$
and $\nbige^{\dagger}$.
We have the commutativity of
$d''_{\lambda}$ and $\DD^{\sankaku}$.

We set 
$W_{k}(\nbigt_0):=
 W_{k+d_X}(\nbigt)\otimes\nbigu_X(-d_X,0)$.
By applying the above construction to
$W_k(\nbigt_0)$,
we obtain $W_k(\nbige^{\sankaku})$.
We have natural inclusions
$W_k(\nbige^{\sankaku})\subset \nbige^{\sankaku}$.
In this way, we obtain an increasing filtration $W$
on $\nbige^{\sankaku}$
compatible with $d''_{\lambda}$ and $\DD^{\sankaku}$.
By the construction,
$(\nbige^{\sankaku},W)$
gives a variation of mixed twistor structure in the sense of 
\cite{s3,Mochizuki-tame,mochi8}.

Because the sesqui-linear pairings of the underlying
$\nbigr$-modules are compatible with
the actions of $\lambda\del_{\lambda}$,
it is preserved by the action of
$S^1=\bigl\{t\in\cnum^{\ast}\,\big|\,|t|=1\bigr\}$
obtained as the restriction of 
the $\cnum^{\ast}$-action.
We obtain that the gluing $\Phi$ is $\cnum^{\ast}$-equivariant.
Hence, $(\nbige^{\sankaku},\DD^{\sankaku})$
is a $\cnum^{\ast}$-equivariant variation of mixed twistor structure.

Let $\gamma:\proj^1\lrarr\proj^1$ be given by 
$\gamma(\lambda)=(\lambdabar)^{-1}$.
The induced map
$\proj^1\times X\lrarr \proj^1\times X$
is also denoted by $\gamma$.
We have the naturally induced
variation of mixed twistor structure
$\gamma^{\ast}(\nbige^{\sankaku},\DD^{\sankaku},W)$
as in \cite{mochi8}.
The isomorphism 
$\kappa:\gammatilde_{sm}^{\ast}\nbigt_0\simeq\nbigt_0$
induces an isomorphism
$\kappa_1:\gamma^{\ast}(\nbige^{\sankaku},\DD^{\sankaku},W)
\simeq
 (\nbige^{\sankaku},\DD^{\sankaku},W)$
such that
(i) $\gamma^{\ast}\kappa_1\circ\kappa_1=\id$,
(ii) $\kappa_1$ is equivariant with respect to $\cnum^{\ast}$-action.
Let $\sigma:\proj^1\lrarr\proj^1$ be given by
$\sigma(\lambda)=-(\lambdabar)^{-1}$.
The induced morphism $\proj^1\times X\lrarr \proj^1\times X$
is also denoted by $\sigma$.
We have the naturally defined variation of mixed twistor structure
$\sigma^{\ast}(\nbige^{\sankaku},\DD^{\sankaku},W)$
as in \cite{s3,Mochizuki-tame,mochi8}.
We have $\sigma=j\circ\gamma$.
By the restriction of the $\cnum^{\ast}$-action,
we naturally have
$(\nbige^{\sankaku},\DD^{\sankaku},W)
\simeq
 j^{\ast}(\nbige^{\sankaku},\DD^{\sankaku},W)$.
Hence, we have
the isomorphism
$\kappa_2:\sigma^{\ast}(\nbige^{\sankaku},\DD^{\sankaku},W)
\simeq
 (\nbige^{\sankaku},\DD^{\sankaku},W)$
such that
(i) $\sigma^{\ast}(\kappa_2)\circ\kappa_2=\id$,
(ii) $\kappa_2$ is equivariant with respect to 
the $\cnum^{\ast}$-action.
Then, by using \cite[Corollary 3.72]{Mochizuki-tame},
we obtain a variation of mixed $\real$-Hodge structure
$(L_{\real},F,W)$
which induces $(\nbige^{\sankaku},\DD^{\sankaku},W)$
by the Rees construction,
where $L_{\real}$ denotes the local system over $\real$,
$F$ denotes the Hodge filtration of 
$H\otimes_{\real}\nbigo_X$,
and $W$ denotes the weight filtration.
Because each integrable variation of pure twistor structure
$\Gr^W_k(\nbigt_0)$ is assumed to have an integrable polarization
compatible with the real structure,
we obtain that
each $\Gr^W_k(\nbige^{\sankaku},\DD^{\sankaku})$
has a polarization which is equivariant
with respect to the $\cnum^{\ast}$-action.
Hence, by \cite[Corollary 3.72]{Mochizuki-tame},
the variation of mixed $\real$-Hodge structure 
$(L_{\real},F,W)$
is polarizable.

Let $(P_{\real},M,F,W)$ be the mixed Hodge module
associated to the graded polarizable variation of 
mixed Hodge structure $(L_{\real},F,W)$,
i.e.,
$P_{\real}=L_{\real}[d_X]$,
$M=L_{\real}\otimes_{\real}\nbigo_X$,
$F$ is the Hodge filtration,
$W_j(P)=W_{j+d_X}(L)[d_X]$.
Then, by construction,
we can observe that
$\Psi_X(P_{\real},M,F,W)$
is naturally isomorphic to $(\nbigt,W)$.

\subsubsection{Admissible variations of mixed twistor structure}
\label{subsection;15.11.17.123}

Let $H$ be a normal crossing hypersurface of $X$.
Let $(\nbigt_0,W)$ be an admissible variation of
integrable mixed twistor structure
with an integrable real structure
and integrable graded polarization.
Suppose that
the underlying $\nbigr_X$-modules of 
 $(\nbigt_0,W)$ are $\cnum^{\ast}$-homogeneous
with respect to the trivial $\cnum^{\ast}$-action.
By the result in \S\ref{subsection;15.11.17.120},
we already have 
a graded polarizable variation of mixed $\real$-Hodge structure
$(L_{\real},F,W)$ on $X\setminus H$
corresponding to
$(\nbigt_0,W)_{|X\setminus H}$.
By the assumption on the KMS-spectrum,
the monodromy automorphisms of $L_{\real}$ 
along the loop around the irreducible components of $H$
are quasi-unipotent.

In particular,
the variation of pure twistor structure 
corresponding to $\Gr^W_k(\nbigt_0)$
comes from a polarized variation of pure Hodge structure
of weight $k$.
Hence, the underlying harmonic bundle of
$\Gr^W_k(\nbigt_0)$ is tame,
and the $\nbigr$-modules of $\Gr^W_k\nbigt_0$
are regular singular along $H$.
We obtain that the $\nbigrtilde$-modules
of $\nbigt_0$ are regular singular.

Let $(\nbigm_1,\nbigm_2,C)$ be the underlying
$\nbigrtilde_X(\ast H)$-triple of $\nbigt_0$.
It is easy to observe that
the $\nbigrtilde_X$-modules
$\nbigm_2[\ast H]$ are also 
$\cnum^{\ast}$-homogeneous.
It gives an object in $\nbigc_2(X)$
in Proposition \ref{prop;15.11.17.101},
and hence 
we have a good filtration $F$ on
$M:=\Xi_{\DR}(\nbigm_2[\ast H])$
with an isomorphism
$\Rtilde(M,F)\simeq
 \nbigm_2[\ast H]$.
Let $F_j^{\circ}(M):=
 F_j(M)\otimes\nbigo_X(\ast H)$.
Then, we have
$\Rtilde(M,F^{\circ})
\simeq
 \nbigm_2$.
In particular,
$F^{\circ}_jM$ are locally free 
$\nbigo_X(\ast H)$-submodules of $M$
such that
$F^{\circ}_j(M)/F^{\circ}_k(M)$
are also locally free $\nbigo_X(\ast H)$-modules.

We set $W_k(M):=\Xi_{\DR}(\nbigm_2)$.
We obtain the filtration $W$ on $M$,
for which 
we naturally have
$\Gr^W_k(M)\simeq
 \Xi_{\DR}(\Gr^W_k\nbigm_2)$.
By applying the above construction to
$W_k(\nbigt_0)$
and $\Gr^W_k(\nbigt_0)$,
we obtain good filtrations
$F^{\circ}$
on $W_k(M)$ and $\Gr^W_k(M)$.

Because $\nbigm_2$ is
an integrable $\nbigrtilde_X(\ast H)$-module
which is regular singular along $H$,
$W_kM$ are meromorphic flat bundles
obtained as the extension of 
$W_kL_{\real}\otimes_{\real}\nbigo_{X\setminus H}$.
Similarly,
$\Gr^W_k(M)$
are regular singular meromorphic flat bundles
obtained as the extension of 
$\Gr^W_k(L_{\real})\otimes_{\real}\nbigo_{X\setminus H}$.
The restriction of $F^{\circ}$ on $W_k(M)$ and $\Gr^W_k(M)$
to $X\setminus H$
are the Hodge filtrations $F$
on $W_k(L_{\real})\otimes_{\real}\nbigo_{X\setminus H}$
and $\Gr^W_k(L_{\real})\otimes_{\real}\nbigo_{X\setminus H}$.

Let us prove that 
$(L_{\real},F,W)$ is admissible along $H$.
It is enough to study it on a small neighbourhood $X_P$
of any smooth point $P\in H$,
according to \cite{kashiwara-mixed-Hodge}.
We may also assume that the monodromy of $L_{\real}$
is unipotent, after the pull back by an appropriate covering
ramified along $H$.
Let $z$ denote the defining equation of $H$ on $X_P$.
Let $V_{\bullet}$ denote the $V$-filtration along $z$ for 
$W_k(\nbigm_2)$, $\Gr^W_k(\nbigm_2)$,
$W_k(M)$ and $\Gr^W_k(M)$.
Because the monodromy is unipotent,
$V_{\bullet}$ is indexed by $\seisuu$,
and $V_{-1}M$ is equal to the Deligne extension
of $M_{|X\setminus H}$.
We set 
$F_j^{\circ}V_{-1}:=F_j^{\circ}\cap V_{-1}$
on $W_k(M)$ and $\Gr^W_k(M)$.
Note that
$V_{-1}W_k(\nbigm_2)
=\Rtilde(V_{-1}W_k(M),F^{\circ})$.
Because
$V_{-1}W_k(\nbigm_2)$
and $V_{-1}\Gr^W_k(\nbigm_2)$
are locally free $\nbigo_{\nbigx_P}$-modules
by the admissibility of $\nbigt_0$,
we obtain that
$\Gr^{F^{\circ}}V_{-1}W_k(M)$
and
$\Gr^{F^{\circ}}V_{-1}\Gr^W_k(M)$
are locally free $\nbigo_{X_P}$-modules.
Because the morphisms
$V_{-1}W_k\nbigm_2
 \lrarr
 V_{-1}\Gr^W_{k}\nbigm_2$
are surjective for any $k\in\seisuu$
by the admissibility of $\nbigt_0$,
the morphisms
$F^{\circ}_jV_{-1}W_kM
 \lrarr
 F^{\circ}_jV_{-1}\Gr^W_{k}M$
are surjective for any $k,j\in\seisuu$,
i.e.,
$F^{\circ}$ on $V_{-1}\Gr^W_k(M)$
is equal to the filtration induced by $V_{-1}W_k(M)$.
By the admissibility of $\nbigt_0$,
we have the relative weight filtration of the action of 
the nilpotent part of $z\del_z$
on $(V_{-1}(M),W)_{|P}$.
Hence, 
$(L_{\real},F,W)$ is admissible.
(See \cite{kashiwara-mixed-Hodge}.)

\vspace{.1in}

We have objects
$\nbigt_{\star}:=\Bigl(
 (\nbigt_0,W)\otimes\nbigu_X(d_X,0)
 \Bigr)[\star H]$
in $\MTMint_{\rnum\times\{0\}}(X,\real)$
for $\star=\ast,!$.
We also have the mixed $\real$-Hodge modules
$(P_{\real\star},M(\star H),F,W)$ $(\star=\ast,!)$
on $X$ such that 
(i) $W_jP_{\real\star|X\setminus H}
 =W_{j+d_X}L_{\real}[d_X]$,
(ii) $F_{|X\setminus H}$ is the Hodge filtration of 
 $(L_{\real},F,W)$,
(iii)
 $W_k(P_{\real\star})_{|X\setminus H}
=W_{k-d_X}L_{\real}[d_X]$.
We have objects
$\Psi_X(P_{\real\star},M(\star H),F,W)$
in $\MTMint_{\rnum\times\{0\}}(X,\real)$.
By the construction,
we have
$\Psi_X(P_{\real\star},M(\star H),F,W)(\ast H)
 \simeq
 \nbigt_0\otimes\nbigu_X(d_X,0)
 \simeq
 \nbigt_{\star}(\ast H)$.
We also have
$\Psi_X(P_{\real\star},M(\star H),F,W)[\star H]
\simeq
\Psi_X(P_{\real\star},M(\star H),F,W)$.
Hence, we have
$\Psi_X(P_{\real\star},M(\star H),F,W)
\simeq
 (\nbigt_{\star},W)$
in $\MTMint_{\rnum\times\{0\}}(X,\real)$,
i.e.,
$(\nbigt_{\star},W)$
are contained in the essential image of
$\Psi_X$.

\subsubsection{End of Proof of Theorem \ref{thm;15.11.17.20}}
\label{subsection;15.11.17.121}

Let us consider the general case.
Let $(\nbigt,W)$ be an object in $\MTMint_{\rnum\times\{0\}}(X,\real)$
such that the underlying $\nbigrtilde$-modules
are $\cnum^{\ast}$-homogeneous
with respect to the trivial action.
We shall construct a mixed $\real$-Hodge module
$(P,M,F,W)$ with an isomorphism
$\Psi_X(P,M,F,W)\simeq (\nbigt,W)$.
Because the functor $\Psi_X$ is fully faithful,
it is enough to study the issue locally around
any point $P$ of $X$.
Note that an integrable polarization compatible with real structure
on $\Gr^W(\nbigt)$ induces
a polarization of $\Gr^W(P,M,F)$.
We take a small neighbourhood $X_P$ of $X$.
We use an induction on the dimension 
of the support $\Supp(\nbigt)$ of $\nbigt$.

We set $Z:=\Supp(\nbigt)$.
By shrinking $X_P$,
we may assume to have 
a holomorphic function $g$ on $X$,
a complex manifold $Y$
with a projective morphism
$\varphi:Y\lrarr X$
and a normal crossing hypersurface $H_Y\subset Y$,
and a variation of mixed twistor structure
$\nbigt_1$ on $(Y,H_Y)$
such that the following holds:
\begin{itemize}
\item
 $\varphi(Y)\subset Z$,
 $\dim(Z\setminus\varphi(Y))<\dim Z$,
and
 $Z\setminus \varphi(Y)\subset g^{-1}(0)$.
\item
 $H_Y=(g\circ \varphi)^{-1}(0)$.
\item
 $\varphi_{\dagger}(\nbigt_1,W)
\simeq
 (\nbigt,W)(\ast g)$. 
\end{itemize}

We have the $\cnum^{\ast}$-homogeneity of $(\nbigt_1,W)$
with respect to the trivial action.
We set $g_Y:=g\circ\varphi$.
Applying the results in \S\ref{subsection;15.11.17.123},
$ \Pi^{a,b}_{g_Y}(\nbigt_1,W)[\star H_Y]$
are contained in the essential image of $\Psi_Y$.
Hence, we obtain that
$\Pi^{a,b}_g(\nbigt,W)[\star g]$ 
are contained in the essential image of $\Psi_X$.
It follows that
$\psi^{(a)}_g(\nbigt,W)$
and 
$\Xi^{(a)}_g(\nbigt,W)$
are contained in the essential image of $\Psi_X$.
By construction,
the underlying $\nbigrtilde_X$-modules of $\phi_g^{(0)}(\nbigt,W)$
are also $\cnum^{\ast}$-homogeneous.
By the assumption of the induction,
 $\phi_g^{(0)}(\nbigt,W)$
is also contained in the essential image of
$\Psi_X$.
Because 
$(\nbigt,W)$
is reconstructed as the cohomology of
$\psi^{(1)}_g(\nbigt,W)
\lrarr
 \phi^{(0)}_g(\nbigt,W)\oplus
 \Xi^{(0)}_g(\nbigt,W)
\lrarr
 \psi^{(0)}_g(\nbigt,W)$,
we obtain that
$(\nbigt,W)$ 
is also contained in the essential image of
$\Psi_X$.
Thus, Theorem \ref{thm;15.11.17.20}
is proved.
\hfill\qed

\section{Graded sesqui-linear dualities}
\subsection{Sesqui-linear dualities and graded sesqui-linear dualities}

Let $X$ be a complex manifold.
Let $(\nbigt,W)$ be a mixed twistor $\nbigd$-module on $X$.
\begin{itemize}
\item
A sesqui-linear duality 
of weight $w$ on $(\nbigt,W)$ is a morphism 
$\nbigs:(\nbigt,W)\lrarr (\nbigt,W)^{\ast}\otimes\newTate(-w)$
such that 
$\nbigs^{\ast}=(-1)^w\nbigs$.
\item
A graded sesqui-linear duality on $(\nbigt,W)$
is a tuple of sesqui-linear dualities 
$\nbigs_w$ $(w\in\seisuu)$
of weight $w$ on $\Gr^W_w\nbigt$.
\item
A graded sesqui-linear duality $(\nbigs_w\,|w\in\seisuu)$
on $(\nbigt,W)$ is called a graded polarization
if each $\nbigs_w$ is a polarization 
of $\Gr^W_w(\nbigt)$.
\end{itemize}
\begin{rem}
The notions of sesqui-linear duality and graded sesqui-linear duality
are the same if $(\nbigt,W)$ is pure.
But, in general, they are not directly related.
A sesqui-linear duality $\nbigs$ of weight $w$
on a mixed twistor $\nbigd$-module $(\nbigt,W)$
induces just morphisms
$\Gr^W_{\ell}(\nbigt)\lrarr 
 \Bigl(
 \Gr^W_{-\ell+2w}(\nbigt)
 \Bigr)^{\ast}
 \otimes\newTate(-w)$.
\hfill\qed
\end{rem}

If $(\nbigt,W)$ is pure of weight $w$,
a sesqui-linear duality of weight $w$
is called just a sesqui-linear duality.

\subsection{Induced graded sesqui-linear dualities}

\subsubsection{Pure case}
\label{subsection;14.11.24.1}

Let $X$ be a complex manifold.
Let $\nbigt$ be a pure twistor $\nbigd$-module of weight $w$
on $X$.
Let $D$ be an effective divisor of $X$.
Recall that we have 
the mixed twistor $\nbigd$-modules
$\nbigt[\star D]$ $(\star=\ast,!)$
obtained as the localizations.
Note that a polarization $\nbigs$ of $\nbigt$
induces a graded polarization 
$\nbigs[\star D]$ of $\nbigt[\star D]$
as explained in \cite{Mochizuki-MTM}.
By the same construction,
from a sesqui-linear duality $\nbigs$ of $\nbigt$,
we obtain a graded sesqui-linear duality
$\nbigs[\star D]=(\nbigs[\star D]_j\,|\,j\in\seisuu)$ 
on $\nbigt[\star D]$.
Let us recall the local construction,
for which  we can take a holomorphic function
$f$ such that $D=(f)_0$.
(The graded sesqui-linear duality is eventually independent
from the choice of the function $f$.
See \cite{Mochizuki-MTM}.)

Let $\nbign:\psi_f^{(a)}(\nbigt)\lrarr \psi_f^{(a-1)}(\nbigt)$
be the canonical morphism.
Recall that the weight filtration of
$\psi_f^{(a)}(\nbigt)$
is the shift of 
the monodromy weight filtration of $\nbign$,
i.e.,
$W(\nbign)_j\psi_f^{(a)}(\nbigt)
=W_{w+1-2a+j}\psi_f^{(a)}(\nbigt)$.
In particular, the induced morphisms
\[
 \Gr^W_{w+1-2a+j}\psi_f^{(a)}(\nbigt)
\stackrel{\nbign^j}{\lrarr}
 \Gr^W_{w+1-2a+j}\psi_f^{(a-j)}(\nbigt)
=
 \Gr^W_{w+1-2a-j}\bigl(\psi_f^{(a)}(\nbigt)\bigr)
\otimes\newTate(-j)
\]
are isomorphisms for $j\geq 0$.
The primitive part
$P\Gr^W_{w+1-2a+j}\psi_f^{(a)}(\nbigt)$
$(j\geq 0)$
is defined to be the kernel of 
$\nbign^{j+1}:
 \Gr^W_{w+1-2a+j}\psi_f^{(a)}(\nbigt)
\lrarr
 \Gr^W_{w-1-2a+j}\psi_f^{(a-j-1)}(\nbigt)$.
We formally set
$P\Gr^W_{w+1-2a+j}\psi_f^{(a)}(\nbigt)=0$
for $j<0$,
in this paper.

We have natural isomorphisms:
\[
 \Gr^W_{w+j}\nbigt[\ast f]
\simeq
\left\{
 \begin{array}{ll}
 0 & (j<0)\\
 \nbigt & (j=0)\\
 P\Gr^W_{w+j}\psi^{(0)}_f\nbigt & (j>0)
 \end{array}
\right.
\]
A sesqui-linear duality $\nbigs[\ast (f)_0]_{w+1+\ell}$
on 
$P\Gr^W_{w+1+\ell}\psi_f^{(0)}\nbigt$ $(\ell\geq 0)$
is induced by the composite of the following morphisms:
\begin{multline}
 \Gr^W_{w+1+\ell}\psi_f^{(0)}\nbigt
\stackrel{a_1}{\lrarr}
 \Gr^W_{w+1+\ell}\psi_f^{(-\ell)}\nbigt
=\Gr^W_{w+1-\ell}\bigl(
 \psi_f^{(0)}\nbigt\bigr)
 \otimes\newTate(-\ell)
\stackrel{a_2}{\lrarr}
 \Gr^W_{-w+1-\ell}\bigl(
 \psi_f^{(0)}(\nbigt^{\ast})\bigr)
 \otimes\newTate(-w-\ell)
 \\
\stackrel{a_3}{\simeq}
 \Bigl(
 \Gr^W_{w-1+\ell}\psi_f^{(1)}(\nbigt)
 \Bigr)^{\ast}
 \otimes\newTate(-w-\ell)
\stackrel{a_4}{\simeq}
  \Bigl(
 \Gr^W_{w+1+\ell}
 \bigl(
 \psi_f^{(0)}(\nbigt)
 \bigr)
 \otimes
 \newTate(1)
 \Bigr)^{\ast}
 \otimes\newTate(-w-\ell) 
 \\
\stackrel{a_5}{\simeq}
 \Bigl(
 \Gr^W_{w+1+\ell}\psi_f^{(0)}(\nbigt)
 \Bigr)^{\ast}
 \otimes\newTate(-w-1-\ell)
\end{multline}
Here, $a_1$ is induced by
$(-\nbign)^{\ell}$,
$a_2$ is induced by $\nbigs$,
and $a_i$ $(i=3,4,5)$
are the isomorphisms in \cite{Mochizuki-MTM}.

For $j<0$,
let 
$P'\Gr^W_{w+1-2a+j}\psi^{(a)}_f\nbigt$
denote the image of
$P\Gr^W_{w+1-2a+j}\psi^{(a-j)}_f\nbigt
\lrarr
 \Gr^W_{w+1-2a+j}\psi^{(a)}_f\nbigt$.
We have the following natural isomorphisms:
\[
 \Gr^W_{w+j}\nbigt[!f]
\simeq
 \left\{
 \begin{array}{ll}
 0 & (j>0)\\
 \nbigt & (j=0)\\
 P'\Gr^W_{w+j}\psi^{(1)}_f(\nbigt)
 & (j<0)
 \end{array}
 \right.
\]
Because
$P'\Gr^W_{w+1-2a+j}\psi^{(a)}_f\nbigt
\simeq
P\Gr^W_{w+1-2a+j}\psi^{(a-j)}_f\nbigt$,
the pure twistor $\nbigd$-modules
$P'\Gr^W_{w-1-\ell}\psi^{(1)}_f(\nbigt)$ $(\ell\geq 0)$
are equipped with the induced sesqui-linear dualities
$\nbigs[!(f)_0]_{w-1-\ell}$.
They are induced by the composite of the following:
\begin{multline}
\Gr^W_{w-1-\ell}\psi_f^{(1)}\nbigt
\stackrel{b_1}{\lrarr}
\Gr^W_{-w-1-\ell}\psi_f^{(1)}(\nbigt^{\ast})\otimes\newTate(-w)
\\
\stackrel{b_2}{\simeq}
\Gr^{W}_{-w+1-\ell}\psi_f^{(0)}(\nbigt^{\ast})\otimes\newTate(-w+1)
=\Bigl(
\Gr^{W}_{w-1+\ell}\psi_f^{(1)}(\nbigt)
\Bigr)^{\ast}\otimes\newTate(-w+1)
 \\
\stackrel{b_3}{\lrarr}
\Bigl(
\Gr^{W}_{w-1-\ell}\psi_f^{(1)}(\nbigt)\otimes\newTate(-\ell)
\Bigr)^{\ast}\otimes\newTate(-w+1)
\stackrel{b_4}{\simeq}
\Bigl(
\Gr^{W}_{w-1-\ell}\psi_f^{(1)}(\nbigt)
\Bigr)^{\ast}
\otimes\newTate(-w+1+\ell)
\end{multline}
Here, $b_1$ is induced by $\psi_f^{(1)}(\nbigs)$,
$b_3$ is the inverse of 
the induced morphism 
of $(-1)^{\ell}\bigl(\nbign^{\ell}\bigr)^{\ast}$,
and $b_i$ $(i=2,4)$ are the natural isomorphisms.
It is also induced by the composite of the following:
\begin{multline}
 \Gr^W_{w-1-\ell}\psi_f^{(1)}(\nbigt)
 \stackrel{c_1}{\lrarr}
 \Gr^{W}_{w-1+\ell}\psi_f^{(1)}(\nbigt)
 \otimes\newTate(\ell)
 \stackrel{c_2}{\lrarr}
  \Gr^{W}_{-w-1+\ell}\psi_f^{(1)}(\nbigt^{\ast})\otimes\newTate(-w+\ell)
\\
= \Gr^{W}_{-w+1+\ell}\psi_f^{(0)}(\nbigt^{\ast})\otimes\newTate(-w+1+\ell)
\stackrel{c_3}{\simeq}
 \Bigl(
   \Gr^{W}_{w-1-\ell}\psi_f^{(1)}(\nbigt)
 \Bigr)^{\ast}
 \otimes\newTate(-w+1+\ell)
\end{multline}
Here, $c_1$ is the inverse of the induced morphism of $\nbign^{\ell}$,
$c_2$ is induced by $\nbigs$,
and $c_3$ is the natural morphism.

The graded sesqui-linear dualities
$\nbigs[!D]$
and $\nbigs[\ast D]$
induce graded sesqui-linear dualities
on the kernel 
and the cokernel of
the morphism
$\varphi:
 \nbigt[!D]
\lrarr
 \nbigt[\ast D]$,
denoted as
$\vecnbigs_{\Ker(\varphi)}=
 \bigl(\nbigs_{\Ker(\varphi),j}\,\big|\,j\in\seisuu\bigr)$
and 
$\vecnbigs_{\Cok(\varphi)}=
 \bigl(\nbigs_{\Cok(\varphi),j}\,\big|\,j\in\seisuu\bigr)$.

Suppose that 
the sesqui-linear duality
$\nbigs:\nbigt\lrarr\nbigt^{\ast}\otimes\newTate(-w)$
is an isomorphism.
Note that we have the induced isomorphism
$\Cok(\varphi)\simeq
 \Ker(\varphi)^{\ast}\otimes\newTate(-w)$,
and hence we have the following isomorphisms
for $j>0$:
\begin{equation}
\label{eq;15.11.9.1}
 \Gr^W_{w+1+j}\Cok(\varphi)
\simeq
 \Bigl(
 \Gr^W_{w-1-j}\Ker(\varphi)
 \Bigr)^{\ast}
 \otimes\newTate(-w)
\end{equation}

The following can be checked by a direct computation.
\begin{lem}
\label{lem;15.11.9.2}
Under the isomorphism {\rm(\ref{eq;15.11.9.1})},
we have
$\nbigs_{\Cok(\varphi),w+1+j}
=\bigl(
 \nbigs_{\Ker(\varphi),w-1-j}
 \bigr)^{-1}$ for $j\geq 0$.
\hfill\qed
\end{lem}

\begin{rem}
If $\nbigs$ is a polarization,
then $\nbigs[\star D]$ are graded polarizations.
\hfill\qed
\end{rem}

\subsubsection{Mixed case}
The construction was also generalized in the mixed case
\cite{Mochizuki-MTM}.
Let $(\nbigt,W)$ be a mixed twistor $\nbigd$-module
on $X$.
Let $D$ be an effective divisor of $X$.
A graded sesqui-linear duality
$\nbigvecs=(\nbigs_w\,|\,w\in\seisuu)$
of $(\nbigt,W)$
induces graded sesqui-linear dualities 
$\vecnbigs[\star D]$
of the mixed twistor $\nbigd$-modules $(\nbigt[\star D],W)$.
We recall the local construction.
If we are given a holomorphic function $f$
such that $D=(f)_0$,
then $\psi_f^{(a)}(\nbigt)$ is equipped with
the two filtrations.
One is the filtration $L$ induced by 
the weight filtration of $\nbigt$.
The other is the relative monodromy weight filtration
$W$ of $\nbign$ with respect to $L$,
which is equal to the weight filtration of 
the mixed twistor $\nbigd$-module $\psi_f^{(a)}(\nbigt)$.
The induced filtration $L$
on $\Gr^W\psi_f^{(a)}(\nbigt)$
has a canonical splitting due to Kashiwara:
\[
 \Gr^W\psi_f^{(a)}(\nbigt)
=\bigoplus_w\Gr^L_w\Gr^W\psi_f^{(a)}(\nbigt)
=\bigoplus_w\Gr^W\Gr^L_w\psi_f^{(a)}(\nbigt)
\]
Hence, the sesqui-linear dualities $\nbigs_w$
on $\Gr^L_w(\nbigt)$ $(w\in\seisuu)$
induce sesqui-linear dualities of
$\Gr^W_j\psi_f^{(a)}(\nbigt)$ $(j\in\seisuu)$.
We have the decomposition:
\[
 \Gr^L_k(\nbigt[\ast D])=A_{1,k}\oplus A_{2,k}
\]
Here, $A_{1,k}$ is the sum of 
the direct summands of $\Gr^L_k(\nbigt)$
whose strict supports are not contained in $D$,
and the support of $A_{2,k}$ is contained in $D$.
As shown in \cite{Mochizuki-MTM},
$A_{2,k}$ is naturally isomorphic to 
a subobject in $\Gr^W_k\psi^{(0)}_f(\nbigt)$.
Hence, it is equipped with the induced sesqui-linear duality,
which is the $k$-th entry of
$\vecnbigs[\star D]$.
If $\vecnbigs$ is a graded polarization,
then $\vecnbigs[\star D]$ are graded polarizations.

\subsection{Push-forward}

\subsubsection{A condition}
\label{subsection;14.11.16.10}

We introduce a condition on the push-forward of 
mixed twistor $\nbigd$-modules
equipped with a graded sesqui-linear duality.
Let $F:X\lrarr Y$ be a projective morphism
of complex manifolds.
Let $(\gbigt,W)$ be a mixed twistor $\nbigd$-module
on $X$
with a graded sesqui-linear duality 
$\vecnbigs=(\nbigs_j\,|\,j\in\seisuu)$.
Recall that 
we have the induced complex
\[
\begin{CD}
 F^{i-1}_{\dagger}\Gr^W_{j+1}\gbigt
@>{a^{i-1}_{j+1}}>>
 F^i_{\dagger}\Gr^W_{j}\gbigt
@>{a^i_j}>>
 F^{i+1}_{\dagger}\Gr^W_{j-1}\gbigt,
\end{CD}
\]
and $\Ker a^i_j/\Image a^{i-1}_{j+1}$
is naturally isomorphic to
$\Gr^W_{j}F^i_{\dagger}\gbigt$.
Here, $a^i_j$ are induced by 
the extensions
$0\lrarr
 \Gr^W_{j-1}\lrarr W_{j}/W_{j-2}
 \lrarr \Gr^W_{j}\lrarr 0$.
We set $\alpha_j:=a_{j+1}^{-1}$
and $\beta_j:=a_{j}^0$.
We have
\[
 \begin{CD}
 F^{-1}_{\dagger}\Gr^W_{j+1}\gbigt
@>{\alpha_j}>>
 F^0_{\dagger}\Gr^W_{j}\gbigt
@>{\beta_j}>>
 F^{1}_{\dagger}\Gr^W_{j-1}\gbigt.
 \end{CD}
\]
As the Hermitian adjoint, we have the following:
\[
\begin{CD}
\Bigl(
  F^{-1}_{\dagger}\Gr^W_{j+1}\gbigt
\Bigr)^{\ast}
@<{\alpha_j^{\ast}}<<
\Bigl(
 F^0_{\dagger}\Gr^W_{j}\gbigt
\Bigr)^{\ast}
@<{\beta_j^{\ast}}<<
\Bigl(
 F^{1}_{\dagger}\Gr^W_{j-1}\gbigt
\Bigr)^{\ast}
\end{CD}
\]
We also have the induced isomorphism
$F_{\dagger}^0\nbigs_{j}:
 F_{\dagger}^0\Gr^W_j\gbigt
\stackrel{\simeq}{\lrarr}
\bigl(
 F_{\dagger}^0\Gr^W_j\gbigt
\bigr)^{\ast}\otimes\newTate(-j)$.

\begin{lem}
\label{lem;14.11.11.1}
Let $\nbigi_j$ denote the image of 
$\Ker\beta_j\cap
 \Ker\bigl(
 \alpha_j^{\ast}\circ F_{\dagger}^0\nbigs_{j}
 \bigr)
\lrarr
\Ker\beta_j/\Image\alpha_j$.
Then, the morphism 
\begin{equation}
 \label{eq;14.11.7.20}
 \Ker\beta_j\cap
 \Ker\bigl(
 \alpha_j^{\ast}\circ F^0_{\dagger}\nbigs_{j}
 \bigr)
\lrarr
\Bigl(
  \Ker\beta_j\cap
 \Ker\bigl(
 \alpha_j^{\ast}\circ F^0_{\dagger}\nbigs_{j}
 \bigr)
\Bigr)^{\ast}
 \otimes\newTate(-j)
\end{equation}
induced by 
$F_{\dagger}^0\nbigs_{j}$
is factorized as follows:
\begin{equation}
 \label{eq;14.11.11.2}
\Ker\beta_j\cap
 \Ker\bigl(
 \alpha_j^{\ast}\circ F_{\dagger}^0\nbigs_{j}
 \bigr)
\lrarr
 \nbigi_j
\stackrel{\nu_j}{\lrarr}
 \nbigi_j^{\ast}\otimes\newTate(-j)
\lrarr
\Bigl(
  \Ker\beta_j\cap
 \Ker\bigl(
 \alpha_j^{\ast}\circ F_{\dagger}^0\nbigs_{j}
 \bigr)
\Bigr)^{\ast}
 \otimes\newTate(-j)
\end{equation}
Namely, we have an induced 
sesqui-linear duality $\nu_j$
of weight $j$ on $\nbigi_j$.
\end{lem}
\pf
Because
$\Bigl(
  \Ker\beta_j\cap
 \Ker\bigl(
 \alpha_j^{\ast}\circ F_{\dagger}^0\nbigs_{j}
 \bigr)
\Bigr)^{\ast}$
is the quotient of
$(F_{\dagger}^0\Gr^W_j\gbigt)^{\ast}$
by 
$\Image\beta_j^{\ast}
+
\Image\bigl(
 F_{\dagger}^0\nbigs_{j}\circ\alpha_j
 \bigr)$.
Hence, the morphism (\ref{eq;14.11.7.20})
factors through $\nbigi_j$.
Because $\nbigi_j^{\ast}$
is the image of 
$\Ker(\alpha_j^{\ast})$
to the quotient,
the morphism (\ref{eq;14.11.7.20})
factors through $\nbigi^{\ast}\otimes\newTate(-j)$.
The condition
$\nu_j^{\ast}=(-1)^{j}\nu_{j}$
is clear by the construction.
\hfill\qed

\begin{df}
We say that Condition {\bf (A)} is satisfied for
the morphism $F:X\lrarr Y$
and the mixed twistor $\nbigd$-module
$(\gbigt,W)$ with the graded sesqui-linear duality
$\vecnbigs$
if the following holds:
\begin{description}
\item[(A)]
 The morphisms
 $\Ker\beta_j\cap
 \Ker\bigl(
 \alpha_j^{\ast}\circ F^0_{\dagger}\nbigs_{j}
 \bigr)
\lrarr
\Ker\beta_j/\Image\alpha_j$
 are epimorphisms
for any $j$.
\hfill\qed
\end{description}
\end{df}

\begin{rem}
If $F_{\dagger}^0\nbigs_j$ is a polarization
of $F_{\dagger}^0\Gr^W_j\gbigt$,
then 
 $\Ker\beta_j\cap
 \Ker\bigl(
 \alpha_j^{\ast}\circ F_{\dagger}^0\nbigs_{j}
 \bigr)
\lrarr
\Ker\beta_j/\Image\alpha_j$
is an isomorphism.
\hfill\qed
\end{rem}

According to Lemma \ref{lem;14.11.11.1},
if the condition {\bf(A)} is satisfied,
we have the induced sesqui-linear duality
$[F^0_{\dagger}\nbigs_j]$ of weight $j$
on $\Gr^W_jF^0_{\dagger}\gbigt$
such that the following diagram is commutative:
\[
 \begin{CD}
 F_{\dagger}^0\Gr^W_j\gbigt
 @>{F_{\dagger}^0\nbigs_j}>>
 \bigl(F_{\dagger}^0
 \Gr^W_j\gbigt\bigr)^{\ast}\otimes\newTate(-j)\\
 @AAA @VVV\\
  \Ker\beta_j\cap
 \Ker\bigl(
 \alpha_j^{\ast}\circ F_{\dagger}^0\nbigs_{j}
 \bigr)
 @>>>
 \Bigl(
   \Ker\beta_j\cap
 \Ker\bigl(
 \alpha_j^{\ast}\circ F_{\dagger}^0\nbigs_{j}
 \bigr)
\Bigr)^{\ast}\otimes\newTate(-j)
 \\
 @VVV @AAA\\
  \Gr^W_jF^0_{\dagger}\gbigt
 @>{[F^0_{\dagger}\nbigs_j]}>>
 \Bigl(
 \Gr^W_jF^0_{\dagger}\gbigt
 \Bigr)^{\ast}\otimes\newTate(-j)
 \end{CD}
\]
The tuple
$\bigl(
 [F_{\dagger}^0\nbigs_j]\,\big|\,
 j\in\seisuu
 \bigr)$
is denoted by
$[F_{\dagger}^0\vecnbigs]$.

\subsubsection{Statements}
\label{subsection;14.11.16.120}

Let $F:X\lrarr Y$ be a projective morphism
of complex manifolds.
Let $\nbigt$ be a pure twistor $\nbigd$-module 
of weight $w$ on $X$ with a sesqui-linear duality $\nbigs_X$.
Let $D_Y$ be an effective divisor of $Y$.
We set $D_X:=F^{\ast}D_Y$.
As explained in \S\ref{subsection;14.11.24.1},
we have the mixed twistor $\nbigd$-modules
$\nbigt[\star D_X]$
with the induced graded sesqui-linear duality
$\nbigs_X[\star D_X]=
 \bigl(\nbigs_X[\star D_X]_{m}\,\big|\,m\in\seisuu\bigr)$.
We shall prove the following theorem
in \S\ref{subsection;14.11.16.101}.
\begin{thm}
\label{thm;14.11.16.110}
Condition {\bf(A)} is satisfied
for the projective morphism $F$
and the mixed twistor $\nbigd$-module
$\nbigt[\star D_X]$ with 
the graded sesqui-linear duality $\nbigs_X[\star D_X]$.
\end{thm}

The pure twistor $\nbigd$-module
$F_{\dagger}^0\nbigt$ of weight $w$
is equipped with the induced sesqui-linear duality
$\nbigs_Y:=F_{\dagger}^0\nbigs_X$.
Moreover, it induces a graded sesqui-linear duality
$\nbigs_Y[\star D_Y]=(\nbigs_Y[\star D_Y]_m\,|\,m\in\seisuu)$
of $F_{\dagger}^0(\nbigt)[\star D_Y]
 \simeq
 F_{\dagger}^0(\nbigt[\star D_X])$.
We shall prove the following theorem
in \S\ref{subsection;14.11.16.101}.
\begin{thm}
\label{thm;14.11.16.111}
We have
$\bigl[F_{\dagger}^0\nbigs_{X}[\star D_X]\bigr]
=\nbigs_{Y}[\star D_Y]$.
\end{thm}

\begin{cor}
Suppose that $\nbigs$ is a polarization,
and that $F_{\dagger}^i\nbigt[\star D_X]=0$ $(i\neq 0)$.
Then, 
$[F_{\dagger}^0\nbigs_X[\star D_X]]$ are graded polarizations.
\end{cor}
\pf
The assumptions imply that $\nbigs_Y$ is a polarization.
Then, the claim follows  from
$[F_{\dagger}^0\nbigs_X[\star D_X]]
=\nbigs_Y[\star D_Y]$.
\hfill\qed

\vspace{.1in}
Let $\nbigk$ and $\nbigc$
denote the kernel and the cokernel of
the natural morphism
$\nbigt[!D_X]\lrarr\nbigt[\ast D_X]$.
They are naturally equipped with
the induced graded sesqui-linear dualities
$\vecnbigs_{\nbigk}=(\nbigs_{\nbigk,j}\,|\,j\in\seisuu)$
and 
$\vecnbigs_{\nbigc}=(\nbigs_{\nbigc,j}\,|\,j\in\seisuu)$.

\begin{cor}
\label{cor;14.11.17.10}
Suppose that 
$F_{\dagger}^i\nbigt[\star D_X]
=F_{\dagger}^i\nbigk
=F_{\dagger}^i\nbigc=0$ for $i\neq 0$.
Then, 
Condition {\bf(A)} is satisfied for
the morphism $F$
and the mixed twistor $\nbigd$-modules 
$\nbigk$ (resp. $\nbigc$)
with $\vecnbigs_{\nbigk}$
(resp. $\vecnbigs_{\nbigc}$).
Moreover, we have
$[F_{\dagger}^0\nbigs_{\nbigk,j}]=\nbigs_{Y}[!D_Y]_j$
for $j<w$
and 
$[F_{\dagger}^0\nbigs_{\nbigc,j}]=\nbigs_{Y}[\ast D_Y]_j$
for $j>w$.
Under the assumptions,
if $\nbigs$ is a polarization,
$[F_{\dagger}^0\vecnbigs_{\nbigk}]$
and 
$[F_{\dagger}^0\vecnbigs_{\nbigc}]$
are graded polarizations.
\end{cor}
\pf
Under the assumption of the corollary,
we also have
$F_{\dagger}^i\nbigt=0$ $(i\neq 0)$.
Then, the claims immediately follow from
Theorem \ref{thm;14.11.16.110}.
\hfill\qed

\vspace{.1in}

For the proof of the theorems
we give an argument in the case $\star=\ast$,
and the other case is similar.
So, to simplify the description,
we denote 
$\nbigs_{X}[\ast D_X]_j$
and 
$\nbigs_{Y}[\ast D_Y]_j$
by 
$\nbigs_{X,j}$ and $\nbigs_{Y,j}$
in the following proof.
Because it is enough to consider the issue locally
on $Y$,
we may assume to have a holomorphic function $g_Y$
such that $D_Y=(g_Y)_{0}$.
The pull back $g_Y\circ F$ is denoted by $g$.
We shall use the notation in \S\ref{subsection;14.11.16.10}
with $\gbigt=\nbigt[\ast D_X]$.

\subsubsection{Preliminary}

We have a natural isomorphism
$F_{\dagger}^i\psi_{g}^{(0)}\nbigt
\simeq
 \psi_{g_Y}^{(0)}F_{\dagger}^i\nbigt$
of mixed twistor $\nbigd$-modules.
By the spectral sequence 
for $\psi_{g}^{(0)}(\nbigt)$
with the weight filtration, 
we have the complex
\begin{equation}
 \label{eq;14.9.6.2}
 F_{\dagger}^{-1}\Gr^W_{j+1}\psi_{g}^{(0)}(\nbigt)
\stackrel{\xi_j}{\lrarr}
 F_{\dagger}^{0}\Gr^W_j\psi_{g}^{(0)}(\nbigt)
\stackrel{\eta_j}{\lrarr}
 F_{\dagger}^{1}\Gr^W_{j-1}\psi_{g}^{(0)}(\nbigt),
\end{equation}
and  $\ker \eta_j/\Image\xi_j$ is naturally isomorphic to
$\Gr^W_jF_{\dagger}^0\psi_{g}^{(0)}(\nbigt)
\simeq
 \Gr^W_j\psi_{g_Y}^{(0)}F_{\dagger}^0(\nbigt)$.

For $\ell\geq 0$,
let $P_{\ell}\Gr^W_{w+1+j}\psi^{(0)}_{g}(\nbigt)$
denote the image of 
$P\Gr^W_{w+1+j}\psi^{(\ell)}_{g}(\nbigt)
\lrarr
 \Gr^W_{w+1+j}\psi^{(0)}_{g}(\nbigt)$.
Note
$P_{\ell}\Gr^W_{w+1+j}\psi_{g}^{(\ell)}(\nbigt)=0$
if $j+2\ell<0$.
We have the primitive decomposition
\[
 \Gr^W_{w+1+j}\psi_g^{(0)}(\nbigt)
=\bigoplus_{\ell\geq 0}
 P_{\ell}\Gr^W_{w+1+j}\psi_g^{(0)}(\nbigt).
\]

\begin{lem}
For $j\geq 0$,
the morphism
$F_{\dagger}^0P\Gr^W_{w+1+j}\psi^{(0)}_g(\nbigt)
\lrarr
 F_{\dagger}^1
 \Gr^W_{w+j}\psi^{(0)}_g(\nbigt)$
factors through
\[
  F_{\dagger}^1
 P_0\Gr^W_{w+j}\psi^{(0)}_g(\nbigt)
\oplus
 F_{\dagger}^1
 P_1\Gr^W_{w+j}\psi^{(0)}_g(\nbigt).
\]
\end{lem}
\pf
Note that $\xi_j$ and $\eta_j$ are compatible 
with the canonical morphisms
\[
 \nbign^{\ell}:
 F_{\dagger}^i\Gr^W_m\psi_g^{(0)}\nbigt
\lrarr
 F_{\dagger}^i\Gr^W_m\psi_g^{(-\ell)}\nbigt
=F_{\dagger}^i\Gr^W_{m-2\ell}\psi_g^{(0)}\nbigt
 \otimes\newTate(-\ell).
\]
For $j\geq 0$,
the morphisms
$\Gr^W_{w+1+j}\psi_g^{(0)}(\nbigt)
\lrarr
 \Gr^W_{w+1+j}\psi_g^{(-j)}(\nbigt)
=\Gr^W_{w+1-j}\psi_g^{(0)}(\nbigt)
 \otimes\newTate(-j)$
are isomorphisms.
Then, the claim easily follows.
\hfill\qed

\vspace{.1in}
The restriction of $\eta_{w+1+j}$ to 
$F_{\dagger}^0P\Gr^W_{w+1+j}\psi^{(0)}_g(\nbigt)$
induces the following morphisms:
\[
 \eta_{k,w+1+j}:
F_{\dagger}^0P\Gr^W_{w+1+j}\psi^{(0)}_g(\nbigt)
\lrarr
  F_{\dagger}^1
 P_k\Gr^W_{w+j}\psi^{(0)}_g(\nbigt),
\quad (k=0,1).
\]

Note that
$\Gr^W_j\Cok\bigl(\psi^{(1)}_g(\nbigt)
 \lrarr
 \psi^{(0)}_g(\nbigt)\bigr)
\simeq
 P\Gr^W_j\psi^{(0)}_g\nbigt$.
By using the spectral sequence
for the cokernel,
we obtain the following complex
\begin{equation}
\label{eq;14.11.11.3}
\begin{CD}
  F_{\dagger}^{-1}P\Gr^W_{j+1}\psi_{g}^{(0)}(\nbigt)
@>{\kappa_{1j}}>>
 F_{\dagger}^{0}P\Gr^W_j\psi_{g}^{(0)}(\nbigt)
@>{\kappa_{2j}}>>
 F_{\dagger}^{1}P\Gr^W_{j-1}\psi_{g}^{(0)}(\nbigt).
\end{CD}
\end{equation}
We have 
$\kappa_{2,w+1+j}
=\eta_{0,w+1+j}$
by construction.
\begin{lem}
\label{lem;14.9.7.4}
The following diagram is commutative
up to signature:
\begin{equation}
 \label{eq;14.9.7.3}
 \begin{CD}
 F_{\dagger}^0P\Gr^W_{w+1+j}\psi_g^{(0)}(\nbigt)
 @>{\eta_{1,w+1+j}}>>
  F_{\dagger}^1
 P_1\Gr^W_{w+j}\psi^{(0)}_g(\nbigt).
 \\
 @V{\simeq}VV @V{\simeq}VV \\
 \Bigl(
 F_{\dagger}^0P\Gr^W_{w+1+j}\psi_g^{(0)}(\nbigt)
 \Bigr)^{\ast}
 \otimes
 \newTate(-w-j-1)
@>{\kappa_{1,w+1+j}^{\ast}}>>
  \Bigl(
  F_{\dagger}^{-1}
  P\Gr^W_{w+j+2}\psi^{(0)}_g(\nbigt)
 \Bigr)^{\ast}
 \otimes\newTate(-w-j-1)
 \end{CD}
\end{equation}
The vertical arrows are induced by 
the induced sesqui-linear dualities
of $P\Gr^W_{w+k}\psi_g^{(0)}(\nbigt)$ $(k> 0)$.
\end{lem}
\pf
Note that we have the following diagram
which is commutative up to signatures.
\begin{equation}
\label{eq;14.9.6.1}
 \begin{CD}
 F_{\dagger}^0
 \Gr^W_{w+1+j}\psi_g^{(0)}(\nbigt)
 @>{\eta_{w+1+j}}>>
 F_{\dagger}^1
 \Gr^W_{w+j}\psi_g^{(0)}(\nbigt)
 \\
 @V{\nbign^j}V{\simeq}V @V{\nbign^j}VV \\
 F_{\dagger}^0
  \Gr^W_{w+1-j}\psi_g^{(0)}(\nbigt)
 \otimes\newTate(-j)
 @>{\eta_{w+1-j}}>>
 F_{\dagger}^1
 \Gr^W_{w-j}\psi_g^{(0)}(\nbigt)
 \otimes\newTate(-j)
 \\
 @V{F_{\dagger}^0\psi_g^{(0)}\nbigs_X}V{\simeq}V 
 @V{F_{\dagger}^1\psi_g^{(0)}\nbigs_X}VV \\
  F_{\dagger}^0
  \Gr^W_{-w+1-j}\psi^{(0)}_g(\nbigt^{\ast})
 \otimes\newTate(-w-j)
@>{\eta_{-w+1-j}'}>>
  F_{\dagger}^1
  \Gr^W_{-w-j}\psi^{(0)}_g(\nbigt^{\ast})
 \otimes\newTate(-w-j)
 \\
 @V{a_1}V{\simeq}V @V{a_2}VV \\
 \Bigl(
 F_{\dagger}^0
 \Gr^W_{w+1+j}\psi^{(0)}_g(\nbigt)
 \Bigr)^{\ast}
 \otimes\newTate(-w-j-1)
@>{\xi_{w+1+j}^{\ast}}>>
 \Bigl(
  F_{\dagger}^{-1}
  \Gr^W_{w+j+2}\psi^{(0)}_g(\nbigt)
 \Bigr)^{\ast}
 \otimes\newTate(-w-j-1)
 \end{CD}
\end{equation}
Here, $a_i$ are morphisms
induced by the natural isomorphisms
$\psi^{(0)}_g(\nbigt^{\ast})
\simeq
 \psi^{(1)}_g(\nbigt)^{\ast}
 \simeq
\psi^{(0)}_g(\nbigt)^{\ast}
\otimes\newTate(-1)$,
and $\eta'_{-w+1-j}$ denotes 
$\eta_{-w+1-j}$ for 
$\nbigt^{\ast}$.
Let 
$\mu:F_{\dagger}^1
 \Gr^W_{w+j}\psi_g^{(0)}(\nbigt)
\lrarr 
 \Bigl(
  F_{\dagger}^{-1}
  \Gr^W_{w+j+2}\psi^{(0)}_g(\nbigt)
 \Bigr)^{\ast}
 \otimes\newTate(-w-j-1)$
denote 
the composite of the right vertical arrows.
The restriction of 
$\mu$ to 
$F_{\dagger}^1P_0\Gr^W_{w+j}\psi^{(0)}_g(\nbigt)$
is $0$,
and the restriction of $\mu$ to 
$F_{\dagger}^1
 P_1\Gr^W_{w+j}\psi^{(0)}_g(\nbigt)$
induces 
an isomorphism
\[
\begin{CD}
 F_{\dagger}^1
 P_1\Gr^W_{w+j}\psi^{(0)}_g(\nbigt)
@>{\simeq}>>
  \Bigl(
  F_{\dagger}^{-1}
  P\Gr^W_{w+j+2}\psi^{(0)}_g(\nbigt)
 \Bigr)^{\ast}
 \otimes\newTate(-w-j-1)
\end{CD}
\]
Thus, we obtain (\ref{eq;14.9.7.3}).
\hfill\qed

\vspace{.1in}

For $j\geq 0$,
we have the natural morphism
\begin{equation}
\label{eq;14.10.10.1}
 \Ker\Bigl(
\begin{CD}
 F_{\dagger}^0P\Gr^W_{w+1+j}\psi^{(0)}_g\nbigt
@>{\eta_{w+1+j}}>>
 F_{\dagger}^1\Gr^W_{w+j}\psi^{(0)}_g\nbigt
\end{CD}
 \Bigr)
\lrarr
 P\Gr^W_{w+1+j}
 F_{\dagger}^0\psi^{(0)}_g\nbigt.
\end{equation}

\begin{lem}
\label{lem;14.10.10.11}
The morphism 
{\rm(\ref{eq;14.10.10.1})}
is an epimorphism.
\end{lem}
\pf
Let $\nbigt=(\nbigm_1,\nbigm_2,C)$.
Let $W$ denote the filtration of
$\psi_g^{(0)}\nbigm_2$
associated to the weight filtration of $\psi_g^{(0)}\nbigt$.
It is enough to prove that 
\begin{equation}
\label{eq;14.10.10.2}
 \Ker\Bigl(
 F_{\dagger}^0P\Gr^W_{w+1+j}\psi^{(0)}_g\nbigm_2
\lrarr
 F_{\dagger}^1\Gr^W_{w+j}\psi^{(0)}_g\nbigm_2
 \Bigr)
\lrarr
 P\Gr^W_{w+1+j}
 F_{\dagger}^0\psi^{(0)}_g\nbigm_2
\end{equation}
is an epimorphism.
Let $f_1$ be a section of
$\Ker\Bigl(
 F_{\dagger}^0\Gr^W_{w+1+j}\psi^{(0)}_g\nbigm_2
\lrarr
 F_{\dagger}^1\Gr^W_{w+j}\psi^{(0)}_g\nbigm_2
\Bigr)$
such that
$\nbign^{j+1}f_1=\del f_2$
for $f_2\in
 F_{\dagger}^{-1}\Gr^W_{w-j}\psi^{(0)}_g\nbigm_2
 \lambda^{j+1}$.
We have
$f_2'\in 
 F_{\dagger}^{-1}\Gr^W_{w+1+j+1}\psi^{(0)}_g\nbigm_2$
such that 
$\nbign^{j+1}f_2'=f_2$.
Then, $f_1-\del f_2'$ is a section of
$\Ker\Bigl(
 F_{\dagger}^0P\Gr^W_{w+1+j}\psi^{(0)}_g\nbigm_2
\lrarr
 F_{\dagger}^1\Gr^W_{w+j}\psi^{(0)}_g\nbigm_2
\Bigr)$.
Hence, (\ref{eq;14.10.10.2}) is an epimorphism.
\hfill\qed

\begin{lem} 
\label{lem;14.11.11.11} 
Under the identification
$\Gr^W_{w+1}(\nbigt[\ast g])
\simeq
 P\Gr^W_{w+1}\psi_g^{(0)}(\nbigt)$,
the kernel of 
\begin{equation}
\label{eq;14.11.11.5}
\begin{CD}
 F_{\dagger}^0P\Gr^W_{w+1}\psi^{(0)}_g\nbigt
@>{\eta_{w+1}}>>
 F_{\dagger}^1\Gr^W_w\psi^{(0)}_g\nbigt
\end{CD}
\end{equation}
is contained in $\Ker\beta_{w+1}$.
(Recall that we use the notation 
in {\rm\S\ref{subsection;14.11.16.10}}
with $\gbigt=\nbigt[\ast g]$.)
\end{lem}
\pf 
Recall that
$\phi^{(0)}_g\nbigt$
is equal to the image of
$\psi^{(1)}_g\nbigt
\lrarr
 \psi^{(0)}_g\nbigt$.
Let $\nbigc$ denote the cokernel.
Let $\nbigc_1\subset\psi^{(0)}_g\nbigt$
denote the inverse image of
$P\Gr^W_{w+1}\psi^{(0)}_g\nbigt\subset\nbigc$
by the projection $\psi_g^{(0)}\nbigt\lrarr\nbigc$.
The extension $0\lrarr 
 \phi^{(0)}_g\nbigt
\lrarr
 \nbigc_1\lrarr
 P\Gr^W_{w+1}\psi_g^{(0)}\nbigt
\lrarr 0$
induces 
the following morphism:
\begin{equation}
 \label{eq;14.11.11.10}
 F_{\dagger}^0P\Gr^W_{w+1}\psi^{(1)}_g\nbigt
\lrarr
 F_{\dagger}^1\phi^{(0)}_g(\nbigt).
\end{equation}
Because $\phi^{(0)}_g$ is an exact functor,
the kernel of $\beta_{w+1}$
is equal to the kernel of (\ref{eq;14.11.11.10}).
Moreover, the kernel of (\ref{eq;14.11.11.10})
is equal to the kernel of the following induced morphism
\begin{equation}
\label{eq;14.12.11.20}
F_{\dagger}^0P\Gr^W_{w+1}\psi^{(1)}_g\nbigt
\stackrel{c_1}{\lrarr}
 \Gr^W_{w+1}F_{\dagger}^1\phi^{(0)}_g(\nbigt).
\end{equation}
We have the following complex
associated to $\phi^{(0)}_g(\nbigt)$
with the weight filtration,
and 
$\Gr^W_{w+1}F_{\dagger}^1\phi^{(0)}_g(\nbigt)$
is $\Ker c_3\big/\Image c_2$:
\[
 F_{\dagger}^0 \Gr^W_{w+1}\phi^{(0)}_g(\nbigt)
\stackrel{c_2}{\lrarr}
 F_{\dagger}^1\Gr^{W}_w\phi^{(0)}_g(\nbigt)
\stackrel{c_3}{\lrarr}
 F_{\dagger}^2\Gr^{W}_{w-1}\phi^{(0)}_g(\nbigt).
\]
We have 
$\Gr^W_{w}\psi^{(0)}_g(\nbigt)
=\Gr^W_w\phi^{(0)}_g(\nbigt)$.
The image of $\eta_{w+1}$ is contained in 
$\Ker c_3$,
and $c_1$ is the composite of
\[
\begin{CD}
 F^0_{\dagger}P\Gr^W_{w+1}\psi_g^{(1)}\nbigt
@>{\eta_{w+1}}>>
 \Ker c_3
@>>>
 \Ker c_3\big/\Image c_2.
\end{CD}
\]
Then, the claim of the lemma follows.
\hfill\qed

\subsubsection{Proof of Theorem \ref{thm;14.11.16.110}
 and Theorem \ref{thm;14.11.16.111}}

\label{subsection;14.11.16.101}

For $j>w$,
we have natural isomorphisms
$\Gr^W_{j}F^i_{\dagger}\nbigt[\ast D_X]
\simeq
 \Gr^W_{j}F^i_{\dagger}\psi_g\nbigt$.
For $j>w$,
the morphism $\kappa_{1j}$
is identified with $\alpha_j$.
For $j>w+1$,
the morphism $\kappa_{2j}$
in (\ref{eq;14.11.11.3})
is identified with $\beta_j$.
We also have Lemma \ref{lem;14.11.11.11}.
By Lemma \ref{lem;14.9.7.4},
we have
\[
 \Ker\Bigl(
 F_{\dagger}^0P\Gr^W_{j}\psi^{(0)}_g\nbigt
\lrarr
 F_{\dagger}^1\Gr^W_{j-1}\psi^{(0)}_g\nbigt
 \Bigr)
\simeq
 \Ker \beta_j
 \cap
 \Ker\bigl(
 \alpha_j^{\ast}\circ
 F_{\dagger}^0\nbigs_{X,j}
 \bigr).
\]
Hence, Lemma \ref{lem;14.10.10.11}
implies Theorem \ref{thm;14.11.16.110}.

\vspace{.1in}

Let us consider the morphism
\[
\begin{CD}
 F_{\dagger}^0\Gr^W_{w+1+j}\psi_g^{(0)}(\nbigt)
@>{F_{\dagger}^0\psi_g^{(0)}\nbigs_X\circ(-\nbign)^j}>>
 \Bigl(
 F_{\dagger}^0\Gr^W_{w+1+j}\psi_g^{(0)}(\nbigt)
 \Bigr)^{\ast}
 \otimes\newTate(-w-j-1)
\end{CD}
\]
It induces a sesqui-linear duality for
$\Ker \eta_{w+1+j}
\subset
  F_{\dagger}^0\Gr^W_{w+1+j}\psi_g^{(0)}(\nbigt)$:
\[
 \Ker \eta_{w+1+j}
\lrarr
 \Bigl(
  \Ker \eta_{w+1+j}
 \Bigr)^{\ast}\otimes\newTate(-w-j-1)
\]
It is factorized as follows:
\begin{multline}
  \Ker\eta_{w+1+j}
\lrarr
 \Gr^W_{w+1+j}
 F_{\dagger}^0\psi_g^{(0)}(\nbigt)
\stackrel{b}{\lrarr}
 \\
 \Bigl(
 \Gr^W_{w+1+j}
 F_{\dagger}^0\psi_g^{(0)}(\nbigt)^{\ast}
 \Bigr)^{\ast}\otimes\newTate(-w-j-1)
\lrarr
 \Bigl(
  \Ker\eta_{w+1+j}
 \Bigr)^{\ast}\otimes\newTate(-w-j-1)
\end{multline}
By construction,
the restriction of $b$ to 
$P\Gr^W_{w+1+j}
 F_{\dagger}^0\psi_g^{(0)}(\nbigt)
\simeq
 \Gr^W_{w+1+j}
 F_{\dagger}^0\nbigt[\ast D_X]
=\Gr^W_{w+1+j}F_{\dagger}^0\nbigt$
is 
the induced sesqui-linear duality
$\nbigs_{Y,w+1+j}$.
Then, the claim of Theorem \ref{thm;14.11.16.111}
follows.
\hfill\qed

\subsection{Basic examples of induced sesqui-linear dualities}

Let $X$ be a complex manifold.
Set $d:=\dim X$.
We have the pure twistor $\nbigd$-module
$\nbigu_X(d,0)=\bigl(
 \nbigo_{\nbigx}\lambda^d,\nbigo_{\nbigx},C_0
 \bigr)$.
Here, $C_0$ is given by
$C_0(s_1,\sigma^{\ast}s_2)=
 s_1\cdot \overline{\sigma^{\ast}(s_2)}$.
The canonical polarization 
$\nbigu_X(d,0)\lrarr
 \nbigu_X(d,0)^{\ast}\otimes\newTate(-d)$
is given by
$\bigl(
 \lambda^d,1
 \bigr)
\longleftrightarrow
 \bigl(
 (-1)^d\cdot 1\cdot\lambda^{d},\,
 (-1)^d\cdot\lambda^d\cdot \lambda^{-d}
 \bigr)$.

Let $D=\sum k_iD_i$ be an effective divisor of $X$
such that $\bigcup D_i$ is normal crossing.
We describe the induced graded polarization
on $\nbigu_X(d,0)[\ast D]$.

\subsubsection{The simplest case}
\label{subsection;14.9.7.11}

Let us consider the case that  $D=(t)_0$
for a coordinate function $t$.
We will not distinguish $D$ and $|(t)_0|$.
We describe the induced graded polarizations
of $\nbigu_X(d,0)[\star t]$ $(\star=\ast,!)$.
It is enough to describe the induced polarizations
on $\psi^{(0)}_t\nbigu_X(d,0)\simeq \nbigu_X(d,0)[\ast t]\big/\nbigu_X(d,0)$
and 
$\psi^{(1)}_t\nbigu_X(d,0)
\simeq
 \Ker\bigl(
 \nbigu_X(d,0)[!t]
\lrarr
  \nbigu_X(d,0)
 \bigr)$.

\begin{lem}
\label{lem;15.11.10.210}
Let $\iota:D\lrarr X$ be the inclusion as above.
The natural isomorphisms
\[
\psi_t^{(0)}\nbigu_X(d,0)
\simeq
 \iota_{\dagger}
 \nbigu_D(d-1,0)\otimes\newTate(-1),
\quad
 \psi_t^{(1)}\nbigu_X(d,0)
\simeq
 \iota_{\dagger}
 \nbigu_D(d-1,0)
\]
are compatible with the polarizations.
\end{lem}
\pf
The natural isomorphism
$\psi_{t,-\vecdelta}\nbigu_X(d,0)\otimes\nbigu_D(-1,0)
\simeq
 \nbigu_D(d-1,0)$
is compatible with the polarizations.
Then, the claim follows from
\cite[Proposition 4.3.2]{Mochizuki-MTM}.
\hfill\qed

\vspace{.1in}
We give a more explicit description of the polarizations.

\begin{lem}
The induced polarization
$\psi_t^{(1)}\nbigu_X(d,0)
\lrarr
 \bigl(
\psi_t^{(1)}
 \nbigu_X(d,0)
 \bigr)^{\ast}
\otimes\newTate(-d+1)$
is given by
\begin{equation}
\label{eq;14.12.13.10}
\bigl(
 [t^{-1}\lambda^d],
 [\deldel_t]
\bigr)
\longleftrightarrow
\bigl(
 (-1)^d[\deldel_t]\lambda^{d-1},
 (-1)^d[t^{-1}\lambda]
\bigr).
\end{equation}
The induced polarization
$\psi_t^{(0)}\nbigu_X(d,0)
\lrarr
 \bigl(
\psi_t^{(0)}
 \nbigu_X(d,0)
 \bigr)^{\ast}
\otimes\newTate(-d-1)$
is given by
\begin{equation}
\label{eq;14.12.13.11}
 \bigl(
 [\deldel_t] \lambda^{d},[t^{-1}]
\bigr)
\longleftrightarrow
 \bigl(
 (-1)^{d+1}[t^{-1}]\lambda^{d+1},
 (-1)^{d+1}[\deldel_t]\lambda^{-1}
 \bigr).
\end{equation}
\end{lem}
\pf
Recall that we have a natural isomorphism
of 
$\psi^{(1)}_t\nbigu_X(d,0)
\simeq
 \iota_{\dagger}
 \psi_{-\vecdelta}\nbigu_X(d,0)
 \otimes\nbigu_X(-1,0)$
studied in \cite[Proposition 4.3.1]{Mochizuki-MTM}.
In this case,
the isomorphisms of the $\nbigr$-modules
are given as follows:
\[
 \psi^{(0)}_t(\nbigo_{\nbigx})
\simeq
 \iota_{\dagger}\psi_{-\vecdelta}(\nbigo_{\nbigx})
 \lambda^{-1},
\quad
 [t^{-1}]
\longleftrightarrow
 (dt/\lambda)^{-1}\lambda^{-1}
\]
\[
 \psi^{(1)}_t(\nbigo_{\nbigx})
\simeq
 \iota_{\dagger}\psi_{-\vecdelta}(\nbigo_{\nbigx}),
 \quad
 [\deldel_t]
 \longleftrightarrow
 -(dt/\lambda)^{-1}
\]
Then, 
according to \cite[Proposition 4.3.2]{Mochizuki-MTM},
the induced polarization
of $\psi^{(1)}_t\nbigu_X(d,0)$
is given by (\ref{eq;14.12.13.10}).

Similarly,
the isomorphism
$\psi^{(0)}_t\nbigu_X(d,0)
=\nbigu_X(d,0)[\ast t]/\nbigu_X(d,0)
\simeq
 \iota_{\dagger}
 \psi_{-\vecdelta}\nbigu_X(d,0)\otimes \nbigu_X(0,-1)$
is given by
\[
 \bigl(
 -[\deldel_t]\lambda^{d},
 [t^{-1}]
 \bigr)
\longleftrightarrow
 \bigl(
 (dt/\lambda)^{-1}\lambda^d,\,\,
 (dt/\lambda)^{-1}\lambda^{-1}
 \bigr).
\]
Then, we obtain the claim for 
$\psi^{(0)}_t\nbigu_X(d,0)$.
\hfill\qed

\vspace{.1in}

We have an isomorphism of 
mixed twistor $\nbigd$-modules
\[
 \nbigu_X(d,0)[\ast t]/\nbigu_X(d,0)
\simeq
 \Ker\bigl(
 \nbigu_X(d,0)[!t]
\lrarr
 \nbigu_X(d,0)
 \bigr)\otimes\newTate(-1)
\]
given by
$\bigl(
 [\deldel_t]\lambda^d,
  [t^{-1}]
\bigr)
\longleftrightarrow
\bigl(
 [t^{-1}]\lambda^{d+1},
 [\deldel_t]\lambda^{-1}
\bigr)$.
It is compatible with the polarizations
(\ref{eq;14.12.13.10}) and (\ref{eq;14.12.13.11}).

\subsubsection{Normal crossing case}

Let us consider the case
that $X$ is equipped with a holomorphic coordinate system
$(x_1,\ldots,x_d)$
such that 
$D=\sum_{i=1}^m k_iD_i$
for some $(k_1,\ldots,k_m)\in\seisuu_{>0}^m$,
where $D_i=\{x_i=0\}$.
We have the decomposition
\[
\Gr^W\bigl(\nbigu_X(d,0)[\ast D]\bigr)
=\bigoplus_{j\geq d}
 \Gr^W_j\bigl(
 \nbigu_X(d,0)[\ast D]
 \bigr).
\]
For $J\subset\{1,\ldots,m\}$,
we set $D_J:=\bigcap_{i\in J}D_i$.
Let $\iota_J:D_J\lrarr X$
denote the inclusion.
We have 
\begin{equation}
 \label{eq;14.10.3.1}
 \Gr^W_j\bigl(
 \nbigu_X(d,0)[\ast D]
 \bigr)
\simeq
 \bigoplus_{|J|=j}
 \iota_{J\dagger}
 \nbigu_{D_J}(d-j,0)
 \otimes\newTate(-j).
\end{equation}
Note that 
$ \iota_{J\dagger}
 \nbigu_{D_J}(d-j,0)
 \otimes\newTate(-j)$
is equipped with the natural polarization
$\bigl((-1)^d,(-1)^d\bigr)$.

\begin{prop}
\label{prop;14.9.28.1}
Under the natural isomorphism
{\rm(\ref{eq;14.10.3.1})},
the induced polarization on
$ \iota_{J\dagger}
 \nbigu_{D_J}(d-j,0)
 \otimes\newTate(-j)$
is equal to
$ \Bigl(
 (-1)^{d}\prod_{i\in J}k_i^{-1},
  (-1)^{d}\prod_{i\in J}k_i^{-1}
 \Bigr)$.
\end{prop}
\pf
It is enough to consider 
the induced polarization
on $\Gr^W_m\nbigu_X(d,0)$.
For $\veca=(a_i)\in\seisuu^{m}$,
we set 
$\phi^{(\veca)}:=
 \phi^{(a_1)}_{x_1}\circ\cdots
 \circ\phi^{(a_m)}_{x_m}$
and 
$\psi^{(\veca)}:=
 \psi^{(a_1)}_{x_1}\circ\cdots
 \circ\psi^{(a_m)}_{x_m}$.
Let $\veczero=(0,\ldots,0)\in\seisuu^m$
and $\vecone=(1,\ldots,1)\in\seisuu^m$.
We set $f:=\prod_{i=1}^m x_i^{k_i}$.
We have
$\psi_f^{(0)}\nbigu(d,0)
=\bigl(
 \psi_f^{(1)}\nbigo_{\nbigx}\lambda^d,
 \psi_f^{(0)}\nbigo_{\nbigx},
 \psi_f^{(0)}C
 \bigr)$.
We consider
$\phi^{(\veczero)}
 \psi^{(0)}_f\nbigu_X(d,0)$.

Recall that
$\psi_f^{(a)}(\nbigo_{\nbigx})$
is isomorphic to 
$\Cok\Bigl(
 \Pi^{a,\infty}_f\nbigo_{\nbigx}[!f]
\lrarr
 \Pi^{a,\infty}_f\nbigo_{\nbigx}[\ast f]
 \Bigr)$.
Then,
$\phi^{(\veczero)}
\psi_f^{(a)}(\nbigo_{\nbigx})$
is isomorphic to
\[
\Cok\Bigl(
 \Pi^{a,\infty}_0\psi^{(\vecone)}
 \nbigo_{\nbigx}
\stackrel{\rho}{\lrarr}
 \Pi^{a,\infty}_0\psi^{(\veczero)}
 \nbigo_{\nbigx}
 \Bigr).
\]
The morphism $\rho$ is induced by
$(\lambda s)^j
\prod_{i=1}^m\deldel_{x_i}
\longmapsto
 (\lambda s)^{j+m}
 \prod_{i=1}^m k_ix_i^{-1}$.
Hence, we have the following isomorphism:
\begin{equation}
 \label{eq;14.9.6.10}
 \phi^{(\veczero)}
 \psi_f^{(0)}
 (\nbigo_{\nbigx})
\simeq
 \bigoplus_{0\leq j\leq m-1}
 \psi^{(\veczero)}(\nbigo_{\nbigx})
 (\lambda s)^j
\end{equation}
Recall that
$\psi_f^{(a)}(\nbigo_{\nbigx})$
is also isomorphic to 
$ \Ker\Bigl(
 \Pi^{-\infty,a}_f\nbigo_{\nbigx}[!f]
\lrarr
 \Pi^{-\infty,a}_f\nbigo_{\nbigx}[\ast f]
 \Bigr)$,
and that
$\phi^{(\veczero)}
\psi_f^{(a)}(\nbigo_{\nbigx})$
is isomorphic to
\[
 \Ker\Bigl(
 \Pi^{-\infty,a}_0\psi^{(\vecone)}
 \nbigo_{\nbigx}
\lrarr
 \Pi^{-\infty,a}_0\psi^{(\veczero)}
 \nbigo_{\nbigx}
 \Bigr).
\]
We obtain the following isomorphism:
\[
\phi^{(\veczero)}
 \psi_f^{(1)}(\nbigo_{\nbigx})
\simeq
 \bigoplus_{-m+1\leq j\leq 0}
 \psi^{(\vecone)}(\nbigo_{\nbigx})
 (\lambda s)^j
\]
Hence,
$\phi^{(\veczero)}\psi_f^{(0)}\nbigu_X(d,0)$
is naturally isomorphic to
$\bigoplus_{j=0}^{m-1}
 \psi^{(\veczero)}\nbigu_X(d,0)
 \otimes
 \newTate(j)$,
where 
$\psi^{(\veczero)}\nbigu_X(d,0)$
is 
\[
 \bigl(
 \psi^{(\vecone)}\nbigo_{\nbigx}\lambda^d,
 \psi^{(\veczero)}\nbigo_{\nbigx},
 \psi^{(\veczero)}C
 \bigr)
\simeq
 \iota_{\mbar\dagger}
 \nbigu_{D_{\mbar}}(d-m,0)
 \otimes\newTate(-m).
\]
The canonical morphism
$\nbign^j:\phi^{(\veczero)}\psi_f^{(a)}\nbigu_X(d,0)
\lrarr
 \phi^{(\veczero)}\psi_f^{(a-j)}\nbigu_X(d,0)$
is given as follows;
the second component is given as 
$b(\lambda s)^i\longmapsto
 b(\lambda s)^i$
if $i\leq a-j+m-1$,
and 
$b(\lambda s)^i
 \longmapsto 0$
otherwise.
The first component is given similarly.

Let us describe the morphism
$\phi^{(\veczero)}\psi_f^{(0)}\nbigu_X(d,0)
\lrarr
 \phi^{(\veczero)}\psi_f^{(0)}\nbigu_X(d,0)^{\ast}
 \otimes\newTate(-d-1)$
induced by the polarization of $\nbigu_X(d,0)$.
We have the isomorphism
$\psi_f^{(0)}(\nbigo_{\nbigx}\lambda^d)
\simeq
 \psi_f^{(1)}(\nbigo_{\nbigx}\lambda^{d})\lambda^{-1}$
induced by
$(\lambda s)^j\longmapsto
 -(\lambda s)^{j+1}\lambda^{-1}$.
It induces an identification
$\phi^{(\veczero)}
 \psi_f^{(0)}(\nbigo_{\nbigx}\lambda^d)
\simeq
 \phi^{(\veczero)}
 \psi_f^{(1)}(\nbigo_{\nbigx}\lambda^d)
 \lambda^{-1}$.
Note that the isomorphism
\[
 \Cok\bigl(
  \Pi^{1,\infty}_f\nbigo_{\nbigx}
\lrarr
 \Pi^{1,\infty}_f\nbigo_{\nbigx}
 \bigr)
 \\
\simeq
\Ker\bigl(
  \Pi^{-\infty,1}_f\nbigo_{\nbigx}
\lrarr
 \Pi^{-\infty,1}_f\nbigo_{\nbigx}
 \bigr)
\]
induces the isomorphism
\begin{multline}
\bigoplus_{1\leq j\leq m}
 \psi^{(\veczero)}(\nbigo_{\nbigx})
 (\lambda s)^j
\simeq
\Cok\Bigl(
 \Pi^{1,\infty}_0\psi^{(\vecone)}\nbigo_{\nbigx}
\lrarr
 \Pi^{1,\infty}_0\psi^{(\veczero)}\nbigo_{\nbigx}
 \Bigr)
 \\
\simeq
 \Ker\Bigl(
 \Pi^{-\infty,1}_0\psi^{(\vecone)}\nbigo_{\nbigx}
\lrarr
 \Pi^{-\infty,1}_0\psi^{(\veczero)}\nbigo_{\nbigx}
 \Bigr)
\simeq
\bigoplus_{-m+1\leq j\leq 0}
 \psi^{(\vecone)}(\nbigo_{\nbigx})
 (\lambda s)^j
\end{multline}
which is given by
$(\lambda s)^{p+m}\prod_{i=1}^m(k_ix_i^{-1})
\longmapsto
 (\lambda s)^p\prod_{i=1}^m\deldel_i$.
In all,
the induced morphism
\[
 \phi^{(\veczero)}\psi^{(0)}(\nbigs):
 \bigoplus_{0\leq j\leq m-1}
 \psi^{(\veczero)}(\nbigo_{\nbigx})
 (\lambda s)^j
\lrarr
 \bigoplus_{-m+1\leq j\leq 0}
  \psi^{(\vecone)}(\nbigo_{\nbigx})
 (\lambda s)^j
 \lambda^{-1}
\]
is given by
$(\lambda s)^{p+m}
 \prod_{i=1}^m(k_ix_i^{-1})
\longmapsto
 (-1)^{d+1}\lambda^{-1}(\lambda s)^{p+1}
 \prod_{i=1}^m\deldel_i$.
Hence, the second component of the induced polarization
of $\Gr^W_{d+m}\phi^{(\veczero)}\psi^{(0)}_f\nbigu(d,0)$
is given by the isomorphisms
$\prod_{i=1}^m (k_ix_i^{-1})
 \longmapsto
 (-1)^{d+m}\lambda^{-m}\prod_{i=1}^m\deldel_i$.
The first component is obtained in the same way.

The isomorphism
$\psi^{(\veczero)}\nbigo_{\nbigx}
\simeq
 \iota_{\dagger}\nbigo_{\nbigd_{\mbar}}
 \lambda^{-m}$
is given by
$\prod_{i=1}^m (x_i^{-1}\lambda)
\longleftrightarrow
 \prod_{i=1}^m(dx_i/\lambda)^{-1}$,
and the isomorphism 
$\psi^{(\vecone)}\nbigo_{\nbigx}
\simeq
 \iota_{\dagger}\nbigo_{\nbigd_{\mbar}}
 \lambda^{-m}$
is given by 
$ \prod_{i=1}^m(-\deldel_i)
 \longleftrightarrow
 \prod_{i=1}^m(dx_i/\lambda)^{-1}$.
Then, the claim of Proposition \ref{prop;14.9.28.1}
follows.

\hfill\qed

\subsection{Nearby cycle functors and maximal functors}
\label{subsection;14.9.6.20}

We give a relation between the induced graded sesqui-linear
dualities on the nearby cycle sheaves
and the maximal sheaves.
Although we do not use the results directly in this paper,
the argument will be useful in \S\ref{section;14.12.12.1}.

Let $\nbigt$ be a pure twistor $\nbigd$-module
of weight $w$
on $X$.
Let $f$ be a holomorphic function on $X$.
We have the following exact sequences
of mixed twistor $\nbigd$-modules:
\[
 0\lrarr
 \psi_f^{(a+1)}\nbigt
 \lrarr
 \Xi^{(a)}_{f}\nbigt
 \lrarr 
 \nbigt^{(a)}[\ast f]
 \lrarr 0
\]
\[
 0\lrarr
 \nbigt^{(a)}[!f]
 \lrarr
 \Xi^{(a)}_f\nbigt
 \lrarr 
 \psi_f^{(a)}\nbigt
 \lrarr 0
\]
Hence, we have the following isomorphisms:
\[
 \Gr^W_{w-2a+j}\psi_f^{(a+1)}\nbigt
 \stackrel{\simeq}\lrarr
 \Gr^W_{w-2a+j}\Xi^{(a)}_f\nbigt
\quad
 (j<0)
\]
\[
 \Gr^W_{w-2a+j}\Xi^{(a)}_f\nbigt
\stackrel{\simeq}{\lrarr}
 \Gr^W_{w-2a+j}\psi^{(a)}_f\nbigt
\quad
 (j>0)
\]
The following exact sequences
have the unique splittings:
\[
 0\lrarr
 \Gr^W_{w-2a}
 \psi^{(a+1)}_f(\nbigt)
\lrarr
 \Gr^W_{w-2a}\Xi^{(a)}_f\nbigt
\lrarr
 \nbigt^{(a)}
 \lrarr 0
\]
\[
 0\lrarr \nbigt^{(a)}
\lrarr
 \Gr^W_{w-2a}\Xi^{(a)}_f\nbigt
\lrarr
 \Gr^W_{w-2a}\psi_f^{(a)}\nbigt
 \lrarr 0
\]
In other words,
we have
$\Gr^W_{w-2a}\Xi^{(a)}_f\nbigt
=\nbigt^{(a)}
\oplus
  \Gr^W_{w-2a}\psi_f^{(a)}\nbigt$,
and 
$\Gr^W_{w-2a}\psi^{(a)}_f\nbigt
=\Gr^W_{w-2a}\psi^{(a+1)}_f\nbigt$
in $\Gr^W_{w-2a}\Xi^{(a)}_f\nbigt$.
The isomorphism
$\Gr^W_{w-2a}\psi^{(a+1)}_f\nbigt
\simeq
 \Gr^W_{w-2a}\psi^{(a)}_f\nbigt$
is induced by
the canonical morphism
$\nbign:\psi^{(a+1)}_f\nbigt
\lrarr
 \psi^{(a)}_f\nbigt$.

\vspace{.1in}

Let  $\nbign:\Xi^{(0)}_f\nbigt
 \lrarr \Xi^{(0)}_f\nbigt\otimes\newTate(-1)$
be the canonical morphism.
Let $W(\nbign)$ denote the monodromy weight filtration
of $\nbign$.
For $j\geq 0$,
let $P\Gr^{W(\nbign)}_{j}\Xi_f^{(0)}\nbigt$
denote the primitive part,
i.e.,
the kernel of
$\nbign^{j+1}:
 \Gr^{W(\nbign)}_{j}(\Xi_f^{(0)}\nbigt)
\lrarr
 \Gr^{W(\nbign)}_{-j-2}(\Xi_f^{(0)}\nbigt)
 \otimes
 \newTate(-j-1)$.

\begin{lem}
We have $W(\nbign)_j=W_{w+j}$ $(j\in\seisuu)$,
and
\begin{equation}
 \label{eq;14.9.7.1}
 P\Gr^{W(\nbign)}_{j}\Xi_f^{(0)}\nbigt
\simeq
 \left\{
 \begin{array}{ll}
 \nbigt & (j=0)\\
 P\Gr^W_{w+j}\psi_f^{(0)}(\nbigt)
 & (j>0)
 \end{array}
 \right.
\end{equation}
\end{lem}
\pf
It is enough to consider the case $w=0$.
Because
the morphism 
$\nbign:\Xi_f^{(0)}\nbigt\lrarr 
 \Xi_f^{(0)}\nbigt\otimes \newTate(-1)$
induces
$W_j\Xi_f^{(0)}\nbigt\lrarr
 \bigl(
 W_{j-2}\Xi_f^{(0)}\nbigt
 \bigr)\otimes\newTate(-1)$.
Let us observe that the induced morphism
$\Gr^W_{j}\Xi_f^{(0)}\nbigt
\lrarr 
 \Gr^W_{-j}\Xi_f^{(0)}\nbigt$
is an isomorphism
for $j>0$.
Note that 
$\nbign:\Xi^{(0)}_f\nbigt\lrarr
 \Xi^{(-1)}_f\nbigt$
factors through
$\Xi^{(0)}_f\nbigt
\lrarr
 \psi^{(0)}_f\nbigt
\lrarr 
 \Xi^{(-1)}_f\nbigt$.
Hence, for $j>0$,
we have the following commutative diagram:
\[
 \begin{CD}
 \Gr^{W}_j\Xi^{(0)}_f\nbigt
 @>{\nbign^j}>>
 \Gr^{W}_{j}\Xi^{(-j)}_f\nbigt
 @>{\nbign}>>
 \Gr^W_{j}\Xi^{(-j-1)}_f\nbigt
 \\ @V{\simeq}VV @A{\simeq}AA @A{\simeq}AA\\
 \Gr^W_{j}\psi^{(0)}_f\nbigt
 @>{\nbign^{j-1}}>{\simeq}>
 \Gr^{W}_{j}\psi^{(-j+1)}_f\nbigt
 @>{\nbign}>>\Gr^W_{j-1}\psi^{(-j)}\nbigt
 \end{CD}
\]
Hence, we have $W=W(\nbign)$ on $\Xi^{(0)}\nbigt$.
We also obtain (\ref{eq;14.9.7.1})
for $j>0$.
We have the following commutative diagram:
\[
 \begin{CD}
 \Gr^W_2\Xi^{(0)}_f(\nbigt)
 @>>>
  \Gr^W_0\Xi^{(0)}_f(\nbigt)\otimes\newTate(-1)
=\nbigt\oplus \Gr^W_0\psi^{(0)}_f(\nbigt)
 @>>>
  \Gr^W_{-2}\Xi^{(0)}_f(\nbigt)\otimes\newTate(-2)
 \\
 @V{\simeq}VV @VVV @A{\simeq}AA \\
 \Gr^W_2\psi^{(0)}_f(\nbigt)
 @>{\simeq}>>
 \Gr^W_0\psi^{(0)}_f(\nbigt)
=\Gr^W_0\psi^{(1)}_f(\nbigt)
 @>{\simeq}>>
 \Gr^W_{-2}\psi^{(1)}_f(\nbigt)
 \end{CD}
\]
Then, we obtain (\ref{eq;14.9.7.1})
for $j=0$.
\hfill\qed

\vspace{.1in}

Let $\nbigs$ be a Hermitian sesqui-linear duality of $\nbigt$.
For $j\geq 0$,
let $\nbigs^{\Xi^{(0)}}_j$ denote
the composite of the following morphisms:
\[
 \Xi^{(0)}_f(\nbigt)
\stackrel{(-\nbign)^{j}}{\lrarr}
 \Xi^{(-j)}_f(\nbigt)
=\Xi^{(0)}_f(\nbigt)\otimes\newTate(-j)
\lrarr
 \Xi^{(0)}_f(\nbigt^{\ast})\otimes\newTate(-j-w)
 \\
\simeq
 \Xi^{(0)}_f(\nbigt)^{\ast}\otimes\newTate(-j-w)
\]
The middle morphism is induced by $\nbigs$.
Let $\nbigs^{\psi^{(0)}}_j$
denote the composite of the following morphisms:
\begin{multline}
\label{eq;14.9.7.2}
 \psi_f^{(0)}(\nbigt)
 \stackrel{(-\nbign)^j}{\lrarr}
 \psi_f^{(-j)}(\nbigt)
=\psi_f^{(0)}(\nbigt)\otimes\newTate(-j)
\lrarr
 \psi_f^{(0)}(\nbigt^{\ast})\otimes\newTate(-w-j)
\simeq
 \psi_f^{(1)}(\nbigt)^{\ast}\otimes\newTate(-w-j)
 \\
\simeq
 \Bigl(
 \psi_f^{(0)}(\nbigt)\otimes\newTate(1)
 \Bigr)^{\ast}\otimes\newTate(-w-j)
\simeq
 \psi_f^{(0)}(\nbigt)^{\ast}\otimes\newTate(-w-j-1)
\end{multline}

\begin{prop}
For $j\geq 0$,
the following diagram is commutative:
\[
\begin{CD}
 \Xi_f^{(0)}(\nbigt)
 @>{\nbigs^{\Xi_f^{(0)}}_{j+1}}>>
 \Xi_f^{(0)}(\nbigt)^{\ast}\otimes\newTate(-w-j-1)
 \\
 @VVV @AAA \\
 \psi_f^{(0)}(\nbigt)
 @>{\nbigs^{\psi_f^{(0)}}_{j}}>>
 \psi_f^{(0)}(\nbigt)^{\ast}
 \otimes\newTate(-w-j-1)
\end{CD}
\]
Here, the vertical morphisms
are natural ones.
\end{prop}
\pf
We can check the claim by a direct computation.
We remark that,
in (\ref{eq;14.9.7.2}),
the isomorphism
$\newTate(1)^{\ast}
\simeq
 \newTate(-1)$
is given by
$(-1,-1)$.
\hfill\qed

\begin{cor}
Suppose that $\nbigs$ is a polarization.
For $j\geq 1$,
the morphism
$\nbigs^{\Xi_f^{(0)}}_{j}$
induces a polarization of
$P\Gr^W_{w+j}\Xi^{(0)}_f(\nbigt)$.
It is equal to the polarization of
$P\Gr^W_{w+j}\psi^{(0)}_f(\nbigt)$
under the natural isomorphism.
In particular,
$\Xi^{(0)}_f\nbigt$
is equipped with an induced 
graded polarization.
\hfill\qed
\end{cor}

\section{Comparisons of polarizations}
\label{section;14.12.12.1}
\subsection{A specialization}
\label{subsection;14.9.6.30}

\subsubsection{Statements}
\label{subsection;14.11.24.10}

Let $X$ be a complex manifold
with a normal crossing hypersurface $D$.
Set $d:=\dim X$.
Let $f,g\in\nbigo_X(\ast D)$.
We consider a meromorphic function
$F:=\tau f+g$
on $X\times\proj^1_{\tau}$.

\begin{assumption}
\label{assumption;14.11.24.11}
We assume the following.
\begin{itemize}
\item
 $|(f)_0|\cap |(f)_{\infty}|=\emptyset$,
 and
 $|(f)_0|\subset |(g)_{\infty}|$.
 We also have $D=|(f)_{\infty}|\cup|(g)_{\infty}|$.
\item
 $F$ is pure on $X\times\{\tau\neq 0\}$,
and $g$ is pure on $X\setminus|(f)_{\infty}|$. 
\hfill\qed
\end{itemize}
\end{assumption}
For example,
the assumptions are satisfied
in the cases 
of Lemma \ref{lem;14.10.14.12} and 
Lemma \ref{lem;14.10.14.11}.

We set 
$D^{(1)}:=(\proj^1\times D)
 \cup (X\times\{\infty\})$.
We have the associated mixed twistor $\nbigd$-modules
$\nbigt_{\star}(F,D^{(1)})$ $(\star=\ast,!)$
which are equipped with natural real structure.
Recall that $\nbigt(F)$ denotes the image of
$\nbigt_{!}(F,D^{(1)})
\lrarr
 \nbigt_{\ast}(F,D^{(1)})$,
which is a pure twistor $\nbigd$-module
of weight $d+1$.
We have 
$\nbigt_{\star}(F,D^{(1)})[\star(\tau)_0]
=\nbigt(F)[\star(\tau)_0]
=\nbigt(F)[\star\tau]$ for $\star=\ast,!$
by the assumption of purity.

Let $\iota:X\lrarr 
 X\times\proj^1_{\tau}$
be given by $\iota(x)=(x,0)$.
According to Proposition \ref{prop;14.10.18.11},
for the morphism
$\varphi:
 \nbigt(F)[!\tau]\lrarr
\nbigt(F)[\ast\tau]$,
we have the following isomorphisms:
\[
 \Cok(\varphi)\simeq
 \iota_{\dagger}\nbigt_{\ast}(g,D)
 \otimes\newTate(-1),
\quad
 \Ker(\varphi)\simeq
 \iota_{\dagger}\nbigt_!(g,D)
\]

The polarization
$\nbigs_F=\bigl((-1)^{d+1},(-1)^{d+1}\bigr)$
of $\nbigt(F)$
induces graded polarizations of
$\nbigt(F)[\star \tau]$,
and they induce graded polarizations of
$\Cok(\varphi)$ and $\Ker(\varphi)$.
The induced graded polarizations are denoted by
$\vecnbigs_{\Cok(\varphi)}
 =(\nbigs_{\Cok(\varphi),w}\,|\,w\in\seisuu)$
and 
$\vecnbigs_{\Ker(\varphi)}
 =(\nbigs_{\Ker(\varphi),w}\,|\,w\in\seisuu)$.
Note that by applying $\phi_{\tau}^{(0)}$
we have the following exact sequence:
\[
 0\lrarr
 \Ker(\varphi)\lrarr
 \psi^{(1)}_{\tau}\nbigt(F)
\lrarr
 \psi^{(0)}_{\tau}\nbigt(F)
 \lrarr\Cok(\varphi)
 \lrarr 0
\]
For the weight filtrations of
$\Cok(\varphi)$ and $\Ker(\varphi)$,
we have
\[
 \Gr^W_{d+1+j}
 \Cok(\varphi)\simeq
 P\Gr^W_{d+1+j}
 \psi_{\tau}^{(0)}(\nbigt(F))\,\,\,
 (j>0),
\quad\quad
  \Gr^W_{d+1+j}\Cok(\varphi)=0\,\,\,
 (j\leq 0),
\]
\[
 \Gr^W_{d+1-j}\Ker(\varphi)
\simeq
 P'\Gr^W_{d+1-j}\psi_{\tau}^{(1)}(\nbigt(F))\,\,\,
 (j>0),
\quad\quad
  \Gr^W_{d+1-j}\Ker(\varphi)=0\,\,\,
 (j\leq 0).
\]
Here,
$P'\Gr^W_{d+1-j}\psi_{\tau}^{(1)}(\nbigt(F))$ is the image of
$P\Gr^W_{d+1-j}\psi_{\tau}^{(1-j)}(\nbigt(F))\lrarr
 \Gr^W_{d+1-j}\psi_{\tau}^{(1)}(\nbigt(F))$.
The isomorphisms are compatible
with the induced polarizations.

We also have other polarizations.
Let $\nbigt(g)$ denote the image of
$\nbigt_!(g,D)\lrarr\nbigt_{\ast}(g,D)$,
which is a pure twistor $\nbigd$-module of weight $d$.
It is equipped with the natural polarization
$\nbigs_g=((-1)^{d},(-1)^{d})$.
Because
$\nbigt(g)[\star (f)_{\infty}]
=\nbigt_{\star}(g,D)$,
we have the induced graded polarizations
$\iota_{\dagger}\nbigs_{g}[\ast (f)_{\infty}]$
of $\iota_{\dagger}\nbigt_{\ast}(g,D)\otimes\newTate(-1)$,
and 
$\iota_{\dagger}\nbigs_{g}[!(f)_{\infty}]$
of $\iota_{\dagger}\nbigt_{!}(g,D)$.
We shall prove the following proposition
in \S\ref{subsection;14.10.3.2}--\ref{subsection;15.11.6.1}.
\begin{prop}
\label{prop;14.9.7.42}
Under the isomorphisms
$\Cok(\varphi)
\simeq
 \iota_{\dagger}\nbigt_{\ast}(g,D)\otimes\newTate(-1)$
and 
$\Ker(\varphi)\simeq
 \iota_{\dagger}\nbigt_!(g,D)$,
we have
$\vecnbigs_{\Cok(\varphi)}=
 \iota_{\dagger}\nbigs_{g}[\ast (f)_{\infty}]$
and 
$\vecnbigs_{\Ker(\varphi)}=
 \iota_{\dagger}\nbigs_{g}[!(f)_{\infty}]$.
\end{prop}

\subsubsection{Some consequences}

Before going to the proof of Proposition \ref{prop;14.9.7.42},
we give a consequence.
Let $X$, $D$, $f$ and $g$ be 
as in \S\ref{subsection;14.11.24.10}.
Let $\rho:X\lrarr Y$ be 
a projective morphism of complex manifolds.
We assume that 
$R^i\rho_{\ast}\bigl(\Omega^j_X(\ast D)\bigr)=0$
for any $i>0$ and $j=0,\ldots,\dim X$.
Set
$\rho_1:=\rho\times\id_{\proj^1_{\tau}}:
 X\times\proj^1_{\tau}
 \lrarr
 Y\times\proj^1_{\tau}$.
Note $\rho_{1\dagger}^i\nbigt(F)[\star(\tau)_0]=0$
$(i\neq 0)$ 
by Proposition {\rm\ref{prop;14.11.17.20}}.
Let $\iota_Y:Y\lrarr Y\times\proj^1_{\tau}$ 
be the inclusion
induced by $\{0\}\lrarr \proj^1_{\tau}$.

\begin{cor}
\label{cor;14.11.24.20}
Suppose that
$\rho_{\dagger}^i\nbigt_{\star}(g,D)=0$ 
for $\star=\ast,!$ and for $i\neq 0$.
\begin{itemize}
\item
Condition {\bf (A)} is satisfied for
the morphism $\rho$ and 
the mixed twistor $\nbigd$-module
$\nbigt_{\star}(g,D)$
with the graded polarization
$\nbigs_g[\star(f)_{\infty}]$.
The induced graded sesqui-linear duality
$\bigl[\rho_{\dagger}^0\nbigs_g[\star (f)_{\infty}]\bigr]$
on $\rho^0_{\dagger}\nbigt_{\star}(g,D)$
is a graded polarization.
\item
Let $\nbigk_Y$ and $\nbigc_Y$ denote
the kernel and the cokernel of
$\rho_{1\dagger}^0\nbigt(F)[!\tau]
\lrarr
 \rho_{1\dagger}^0\nbigt(F)[\ast\tau]$.
Let $\vecnbigs_{\nbigc_Y}$
and $\vecnbigs_{\nbigk_Y}$
denote the graded polarizations
induced by 
the graded polarization
$\rho_{1\dagger}^0\nbigs_F$ 
of $\rho_{1\dagger}^0\nbigt(F)$.
Then, under the natural isomorphisms
$\iota_{Y\dagger}
 \rho^0_{\dagger}\nbigt_{!}(g,D)
\simeq
 \nbigk_Y$
and 
$\iota_{Y\dagger}
 \rho^0_{\dagger}\nbigt_{\ast}(g,D)
 \otimes\newTate(-1)
\simeq
 \nbigc_Y$,
we have
$\iota_{Y\dagger}
 \bigl[\rho_{\dagger}^0\nbigs_g[!(f)_{\infty}]\bigr]
=\vecnbigs_{\nbigc_Y}$
and 
$\iota_{Y\dagger}
 \bigl[\rho_{\dagger}^0\nbigs_g[\ast(f)_{\infty}]\bigr]
=\vecnbigs_{\nbigk_Y}$.
\end{itemize}
\end{cor}
\pf
It follows from Corollary \ref{cor;14.11.17.10}
and Proposition \ref{prop;14.9.7.42}.
\hfill\qed

\vspace{.1in}

Let $\kappa:X'\lrarr X$ be a projective morphism
of complex manifolds
such that 
(i) $D'=\kappa^{-1}(D)$ is normal crossing,
(ii) $\kappa$ induces $X'\setminus D'\simeq X\setminus D$,
(iii) Assumption \ref{assumption;14.11.24.11} is also satisfied 
for $f'=\kappa^{\ast}(f)$ and $g'=\kappa^{\ast}(g)$.
We have
$\kappa_{\dagger}^i\nbigt_{\star}(g',D')=0$ $(i\neq 0)$,
and we have natural isomorphisms
$\kappa_{\dagger}^0\nbigt_{\star}(g',D')
\simeq
 \nbigt_{\star}(g,D)$.
We set $\rho':=\rho\circ\kappa$.
We have natural isomorphisms
$\rho_{\dagger}^{\prime\,0}
 \nbigt_{\star}(g',D')
\simeq
 \rho_{\dagger}^0\nbigt_{\star}(g,D)$.

\begin{cor}\mbox{{}}
\label{cor;14.11.24.31}
\begin{itemize}
\item
Condition {\bf(A)} holds
for $\kappa$ 
and $\nbigt_{\star}(g',D')$
with $\nbigs_{g'}[\star (f')_{\infty}]$.
We have 
$\bigl[
 \kappa^0_{\dagger}\nbigs_{g'}[\star (f')_{\infty}]
 \bigr]
=\nbigs_g[\star(f)_{\infty}]$
under the natural isomorphisms.
\item
Condition {\bf(A)} holds
for  $\rho'$ 
and $\nbigt_{\star}(g',D')$
with $\nbigs_{g'}[\star (f')_{\infty}]$.
We have 
$\bigl[
 \rho^{\prime 0}_{\dagger}\nbigs_{g'}[\star (f')_{\infty}]
 \bigr]
=\bigl[
 \rho^0_{\dagger}\nbigs_g[\star(f)_{\infty}]
\bigr]$
under the natural isomorphisms.
\hfill\qed
\end{itemize}
\end{cor}

\subsubsection{Preliminary}
\label{subsection;14.10.3.2}

Let us return to the proof of
Proposition \ref{prop;14.9.7.42}.
By Lemma \ref{lem;15.11.9.2},
it is enough to study 
the polarizations on
$\Cok(\varphi)
\simeq
 \iota_{\dagger}
 \nbigt_{\ast}(g,D)\otimes\newTate(-1)$.
The comparison of the polarizations on 
$\Gr^W_{d+2}\Cok(\varphi)
\simeq
 \Gr^W_{d+2}\bigl(
 \iota_{\dagger}\nbigt_{\ast}(g,D)
 \otimes\newTate(-1)\bigr)$
can be given as in \S\ref{subsection;14.9.7.11}.
Hence, we shall argue the other parts.

We have the identification
$\Gr^W_{d+2+j}\Cok(\varphi)
=P\Gr^W_{d+2+j}\psi^{(0)}_{\tau}\nbigt(F)$
for $j\geq 0$.
For $j>0$,
we have the isomorphism
$\Gr^W_{d+2+j}\bigl(\nbigt_{\ast}(g,D)\otimes\newTate(-1)\bigr)
\simeq
 P\Gr^W_{d+2+j}\psi_{f^{-1}}^{(0)}(\nbigt(g)\otimes\newTate(-1))$.
Let $\pi:X\times\proj^1_{\tau}\lrarr X$ be the projection.
To compare $\nbigs_{\Cok(\varphi),d+j+2}$
and $\iota_{\dagger}\nbigs_g[\ast(f)_{\infty}]_{d+j+2}$
for $j> 0$,
it is enough to compare 
$\pi_{\dagger}^0\nbigs_{\Cok(\varphi),d+j+2}$
and $\nbigs_g[\ast(f)_{\infty}]_{d+j+2}$
under the isomorphisms:
\begin{equation}
 \label{eq;15.11.6.2}
 \pi_{\dagger}^0
 P\Gr^W_{d+2+j}\psi_{\tau}^{(0)}\nbigt(F)
\simeq
 P\Gr^W_{d+2+j}
 \psi_{f^{-1}}^{(0)}\bigl(\nbigt(g)\otimes\newTate(-1)\bigr)
\end{equation}
\begin{equation}
 \label{eq;15.11.6.200}
 \Bigl(
 P\Gr^W_{d+2+j}
 \psi_{f^{-1}}^{(0)}\bigl(\nbigt(g)\otimes\newTate(-1)\bigr)
 \Bigr)^{\ast}
\simeq
 \Bigl(
 \pi_{\dagger}^0
 P\Gr^W_{d+2+j}\psi_{\tau}^{(0)}\nbigt(F)
 \Bigr)^{\ast}
\end{equation}
Recall that any morphisms of pure twistor $\nbigd$-modules
are uniquely determined by their underlying morphisms of
$\nbigd$-modules.
Hence, it is enough to compare 
the morphisms of the underlying $\nbigd$-modules
for $\pi_{\dagger}^0\nbigs_{\Cok(\varphi),d+j+2}$ and
$\nbigs_g[\ast(f)_{\infty}]_{d+j+2}$.

The polarization 
$\nbigs_{\Cok(\varphi),d+2+j}$
of
$P\Gr^W_{d+2+j}\psi_{\tau}^{(0)}\nbigt(F)$
is induced by the composite of the following morphisms:
\begin{multline}
\label{eq;14.9.7.30}
\psi_{\tau}^{(0)}\nbigt(F)
\stackrel{(-\nbign)^j}{\lrarr}
\psi_{\tau}^{(-j)}(\nbigt(F))
=\psi_{\tau}^{(0)}\nbigt(F)\otimes\newTate(-j)
 \\
\stackrel{\nbigs_F}{\lrarr}
 \psi_{\tau}^{(0)}(\nbigt(F)^{\ast})
 \otimes\newTate(-j-d-1)
\lrarr
 \psi_{\tau}^{(1)}\bigl(
 \nbigt(F)
 \bigr)^{\ast}\otimes\newTate(-d-j-1)
 \\
\simeq
 \psi_{\tau}^{(0)}(\nbigt(F))^{\ast}
 \otimes\newTate(-d-j-2)
\end{multline}
The polarization $\nbigs_g[\ast(f)_{\infty}]_{d+2+j}$ on 
$\Gr^W_{d+2+j}\bigl(\nbigt_{\ast}(g,D)\otimes\newTate(-1)\bigr)$
is induced by the composite of 
the following morphisms:
\begin{multline}
\label{eq;14.9.7.41}
 \psi_{f^{-1}}^{(0)}\bigl(
 \nbigt(g)\otimes\newTate(-1)
\bigr)
\stackrel{(-\nbign)^{j-1}}
\lrarr
  \psi_{f^{-1}}^{(0)}\bigl(
 \nbigt(g)\otimes\newTate(-1)
\bigr)\otimes\newTate(-j+1)
 \\
\stackrel{\nbigs_{g}}
{\lrarr}
 \psi_{f^{-1}}^{(0)}\bigl(
 \nbigt(g)^{\ast}
 \otimes\newTate(-1)
\bigr)
 \otimes\newTate(-j-d-1)
\simeq
  \psi_{f^{-1}}^{(1)}\bigl(
 \nbigt(g)
 \otimes\newTate(-1)
\bigr)^{\ast}
 \otimes\newTate(-j-d-1)
 \\
 \simeq
  \psi_{f^{-1}}^{(0)}\bigl(
 \nbigt(g)
 \otimes\newTate(-1)
\bigr)^{\ast}
 \otimes\newTate(-j-d-2)
\end{multline}

Let $L(F)$ denote the $\nbigd_{X\times\proj^1}$-module
underlying $\nbigt(F)$.
Let $L(g)$ denote the $\nbigd_X$-module
underlying $\nbigt(g)$.
By the assumption of purity,
we have 
$L(F)[\star\tau]=L_{!}(F,D^{(1)})[\star\tau]
=L_{\ast}(F,D^{(1)})[\star\tau]$
and 
$L_{\star}(g,D)=L(g)[\star (f)_{\infty}]$.

By the above consideration,
Proposition \ref{prop;14.9.7.42}
is reduced to the following lemma.

\begin{lem}
\label{lem;15.11.6.300}
We have the following morphisms of $\nbigd$-modules
\[
 B_1:
 \pi_{\dagger}^0
 \psi_{\tau}^{(0)}L(F)
\lrarr
 \psi^{(0)}_{f^{-1}}L(g),
\quad
 B_2:
 \psi^{(1)}_{f^{-1}}L(g)
\lrarr
 \pi_{\dagger}^0
 \psi_{\tau}^{(1)}L(F)
\]
with the following properties:
\begin{description}
\item[(P1)]
 $B_1$ and $B_2$ induce
 the isomorphisms of $\nbigd$-modules
 underlying {\rm(\ref{eq;15.11.6.2})}
 and {\rm(\ref{eq;15.11.6.200})}.
\item[(P2)]
 The following diagrams are commutative:
\begin{equation}
 \label{eq;15.11.6.11}
 \begin{CD}
 \pi_{\dagger}^0\psi_{\tau}^{(0)}(L(F))
 @>{\epsilon}>>
  \pi_{\dagger}^0\psi_{\tau}^{(1)}L(F)
 \\
 @V{B_1}VV @A{B_2}AA \\
\psi_{f^{-1}}^{(0)}L(g)
 @>{\delta}>>
\psi_{f^{-1}}^{(1)}L(g)
 \end{CD}
\end{equation}
Here,
$\epsilon$ gives the morphism
underlying the composite of
the morphisms in {\rm(\ref{eq;14.9.7.30})}
applied $\pi_{\dagger}^{0}$,
and $\delta$ gives the morphism
underlying the composite of the morphisms
in {\rm(\ref{eq;14.9.7.41})}. 
\end{description}
\end{lem}

\subsubsection{The construction of the morphisms $B_1$ and $B_2$}

As in Proposition \ref{prop;14.10.15.1},
we have a natural isomorphism
$\Lambda^{(a)}:\pi_{+}^0\psi_{\tau}^{(a)}L_{\ast}(F)
\simeq
 \Xi^{(a)}_{f^{-1}}L_{\ast}(g)$.
Let $B_1$ denote the composite of
$\Lambda^{(0)}$
and the natural morphism
$\Xi^{(0)}_{f^{-1}}L_{\ast}(g)
\lrarr
 \psi^{(0)}_{f^{-1}}L_{\ast}(g)$.
Let $B_2$ be the composite of the following morphisms:
\[
\begin{CD}
 \psi^{(1)}_{f^{-1}}L_{\ast}(g)
 @>{c_1}>>
 \Xi^{(0)}_{f^{-1}}L_{\ast}(g)
 @>{c_2}>{\simeq}>
 \Xi^{(1)}_{f^{-1}}L_{\ast}(g)
 @>{c_3}>{\simeq}>
 \pi_{+}^0\psi_{\tau}^{(1)}L_{\ast}(F)
\end{CD}
\]
Here, $c_1$ is the canonical morphism,
$c_2$ is induced by the multiplication of $-s$,
and $c_3$ is $(\Lambda^{(1)})^{-1}$.

\subsubsection{Proof of {\bf (P1)}}

Let us check that 
$B_1$ induces the morphism of $\nbigd$-modules
underlying (\ref{eq;15.11.6.2}).
We have the following natural morphisms:
\begin{multline}
 \pi_+
 \psi^{(0)}_{\tau}L(F)
\stackrel{a_1}{\lrarr}
 \pi_+
 \Cok\bigl(
 \psi^{(1)}_{\tau}L(F)\lrarr
 \psi^{(0)}_{\tau}L(F)
 \bigr)
\stackrel{a_2}{\simeq}
 \pi_+
 \Cok\bigl(
 L(F)[!\tau]\lrarr L(F)[\ast\tau]
 \bigr)
\stackrel{a_3}{\simeq}
L_{\ast}(g,D)
 \\
 \lrarr
 \Cok\Bigl(
 L_!(g,D)\lrarr
 L_{\ast}(g,D)
 \Bigr)
\end{multline}
(See (\ref{eq;15.11.10.100}) for $a_3$.)
The composite induces the morphism of
$\nbigd$-modules underlying (\ref{eq;15.11.6.2}).
Under the natural isomorphism
\[
 \Cok\Bigl(
 L_!(g,D)\lrarr
 L_{\ast}(g,D)
 \Bigr)
\simeq
 \Cok\Bigl(
 \psi_{f^{-1}}^{(1)}L(g)
 \lrarr
  \psi_{f^{-1}}^{(0)}L(g)
 \Bigr)=:\nbiga,
\]
the following diagram is commutative:
\[
 \begin{CD}
 \Xi^{(0)}_{f^{-1}}L(g)
 @>>>
 L_{\ast}(g,D)\\
 @VVV @VVV  \\
 \psi^{(0)}_{f^{-1}}L(g)
 @>>> \nbiga
 \end{CD}
\]
Hence, it is enough to prove that
$b_1:=a_3\circ a_2\circ a_1$
is equal to
the composite $b_2$ of
$\Lambda^{(0)}$ and 
$\Xi^{(0)}_{f^{-1}}L(g)\lrarr L_{\ast}(g,D)$.
The kernel of both of $b_i$ are the image of
$\pi_+\psi^{(1)}_{\tau}L(F)
\lrarr
 \pi_+\psi^{(0)}_{\tau}L(F)$.
The automorphisms of $L_{\ast}(g,D)$
are only the multiplication of locally constant functions.
Hence, it is enough to check the claim
in the case where $D$ is empty,
in which both of the morphisms are given by
the correspondence 
$\tau^{-1}d\tau\longmapsto 1$.
Therefore, we obtain that 
$B_1$ induces the morphism of $\nbigd$-modules
underlying (\ref{eq;15.11.6.2}).

Let us check that $B_2$ induces the morphism
of $\nbigd$-modules
underlying (\ref{eq;15.11.6.200}).
We have the following natural morphisms:
\begin{multline}
 \Ker\Bigl(
 L_!(g,D)\lrarr L_{\ast}(g,D)
 \Bigr)
 \lrarr
 L_!(g,D)
\stackrel{c_1}{\simeq}
 \pi_+
 \Ker\bigl(
 L(F)[!\tau]\lrarr L(F)[\ast\tau]
 \bigr)
 \\
\stackrel{c_2}{\simeq}
 \pi_+
 \Ker\bigl(
 \psi^{(1)}_{\tau}L(F)\lrarr
 \psi^{(0)}_{\tau}L(F)
 \bigr)
\stackrel{c_3}{\lrarr}
 \pi_+
 \psi^{(1)}_{\tau}L(F)
\end{multline}
(See (\ref{eq;15.11.10.101})
for $c_1$.)
As before,
it is enough to check that
$d_1:=c_3\circ c_2\circ c_1$
is equal to
the composite $d_2$ of
the following:
\[
 L_!(g,D)\lrarr
 \Xi^{(0)}_{f^{-1}}L(g)
\stackrel{-s}{\simeq}
 \Xi^{(1)}_{f^{-1}}L(g)
\stackrel{(\Lambda^{(1)})^{-1}}{\lrarr}
 \pi_+^0\psi_{\tau}^{(1)}L_{\ast}(F)
\]
The images of $d_i$ are equal,
and the automorphisms of $L_!(g,D)$
are multiplications of locally constant functions.
Hence, it is enough to consider the case $D$ is empty.
We have the following:
\begin{equation}
 \label{eq;15.11.6.210}
\psi^{(1)}_{\tau}L(F)
\simeq
 \Cok\Bigl(
 \Pi^{1,N}_{\tau}L(F)[!\tau]
\lrarr
 \Pi^{1,N}_{\tau}L(F)[\ast\tau]
 \Bigr)
\simeq
 \Ker\Bigl(
 \Pi^{-N,1}_{\tau}L(F)[!\tau]
\lrarr
 \Pi^{-N,1}_{\tau}L(F)[\ast\tau]
 \Bigr)
\end{equation}
Under the isomorphism (\ref{eq;15.11.6.210}),
$-\del_{\tau}(1)+f$
in $\Pi^{-N,1}_{\tau}L(F)[!\tau]$
corresponds to
$-s/\tau$
in $\Pi^{1,N}_{\tau}L(F)[\ast\tau]$.
By using it, we can check $d_1=d_2$.
Hence, we obtain that
$B_2$ induces the isomorphism of $\nbigd$-modules
underlying (\ref{eq;15.11.6.200}).

\subsubsection{Proof of {\bf(P2)}}
\label{subsection;15.11.6.1}

Let us look at the morphisms of the $\nbigd$-module
underlying (\ref{eq;14.9.7.30}):
\begin{equation}
 \psi_{\tau}^{(0)}L(F)
\stackrel{(-s)^j}\lrarr
 \psi_{\tau}^{(0)}L(F)
\stackrel{(-1)^{d+1}}{\lrarr}
 \psi_{\tau}^{(0)}L(F)
\stackrel{-s}{\lrarr}
 \psi_{\tau}^{(1)}L(F) 
\end{equation}
By applying $\pi_+^0$,
we obtain the following morphisms:
\[
 \Xi^{(0)}_{f^{-1}}\bigl(
 L(g)
 \bigr)
\stackrel{(-s)^j}{\lrarr}
 \Xi^{(0)}_{f^{-1}}\bigl(
 L(g)
 \bigr)
\stackrel{(-1)^{d+1}}{\lrarr}
 \Xi^{(0)}_{f^{-1}}\bigl(
 L(g)
 \bigr) 
 \\
\stackrel{-s}{\simeq}
 \Xi^{(1)}_{f^{-1}}\bigl(
 L(g)
 \bigr)
\]

Let us look at the morphism of $\nbigd$-modules
underlying (\ref{eq;14.9.7.41}):
\[
 \psi^{(0)}_{f^{-1}}(L(g))
\stackrel{(-s)^{j-1}}{\lrarr}
 \psi^{(0)}_{f^{-1}}(L(g))
\stackrel{(-1)^{d}}{\simeq}
  \psi^{(0)}_{f^{-1}}(L(g))
\stackrel{-1}{\simeq}
  \psi^{(0)}_{f^{-1}}(L(g))
\stackrel{-s}{\simeq}
  \psi^{(1)}_{f^{-1}}(L(g))
\]
We have the following commutative diagram:
\[
 \begin{CD}
 \Xi^{(0)}_{f^{-1}}L(g)
 @>{(-1)^{d+1}(-s)^j}>>
 \Xi^{(0)}_{f^{-1}}L(g)
 @>{-s}>> 
 \Xi^{(1)}_{f^{-1}}L(g)
 \\
 @VVV @AAA \\
 \psi^{(0)}_{f^{-1}}L(g)
 @>{(-1)^{d+1}(-s)^{j}}>>
 \psi^{(1)}_{f^{-1}}L(g)
 \end{CD}
\]
It implies the commutativity of (\ref{eq;15.11.6.11}).
Thus, 
we finish the proof of 
Lemma \ref{lem;15.11.6.300}
and Proposition \ref{prop;14.9.7.42}.
\hfill\qed

\subsection{A push-forward}
\label{subsection;14.10.6.10}

\subsubsection{Statement}

Let $X$ be a complex manifold
with a normal crossing hypersurface $D$.
Set $d:=\dim X$.
\begin{assumption}
Let $f$ be a meromorphic function on $(X,D)$
such that
(i) $|(f)_{\infty}|=D$,
(ii) $f$ is non-degenerate along $D$,
(iii) $(f)_0$ is smooth.
\hfill\qed
\end{assumption}
We set $Z:=|(f)_0|$.
We have the meromorphic function $tf$
on $X\times\proj^1_t$.
We have the associated mixed twistor $\nbigd$-modules
$\nbigt(tf)[\star t]$
and 
$\nbigt_{\star}(tf)$.
Let $\pi:X\times\proj^1_t\lrarr X$
be the projection.
According to Proposition \ref{prop;14.11.7.21}
and Corollary \ref{cor;14.11.7.22},
we have
$\pi_{\dagger}^i\nbigt(tf)[\star t]=0$
for $i\neq 0$,
and we have the isomorphisms
\begin{equation}
\label{eq;14.11.7.23}
\Psi_{\ast}:
\pi_{\dagger}^0\nbigt(tf)[\ast t]
\simeq
 \nbigu_X(d,0)[!Z][\ast D]\otimes\newTate(-1),
\quad
\Psi_!:
 \pi_{\dagger}^0\nbigt(tf)[! t]
\simeq
 \nbigu_X(d,0)
 [\ast Z][! D].
\end{equation}

Recall that $\nbigt(tf)$ is the image of
$\nbigt_{!}(tf)\lrarr  \nbigt_{\ast}(tf)$.
We have the description
$\nbigt(tf)
=\bigl(
 \nbigl(tf)\lambda^{d+1},
 \nbigl(tf),C_{tf}
 \bigr)$
as an $\nbigr$-triple.
The polarization is given by
$\nbigs_{tf}=\bigl((-1)^{d+1},(-1)^{d+1}\bigr)$.
It induces graded polarizations
$\vecnbigs_{tf\star}$ 
(resp. $\nbigs_{tf}[\star t]$)
on $\nbigt_{\star}(tf)$
(resp. $\nbigt(tf)[\star t]$).
The polarization  $\nbigs_0=((-1)^d,(-1)^d)$
of $\nbigu_X(d,0)$ induces
graded polarizations
$\nbigs_0[\ast (f)_0][!(f)_{\infty}]$
on $\nbigu_X(d,0)[\ast Z][!D]$,
and $\nbigs_0[!(f)_0][\ast (f)_{\infty}]$ on
$\nbigu_X(d,0)[! Z][\ast D]\otimes\newTate(-1)$.

\begin{prop}
\label{prop;14.9.7.210}
\mbox{{}}
\begin{itemize}
\item
We have
$\pi_{\dagger}^i\Gr^W\nbigt(tf)[\star t]=0$
and 
$\pi_{\dagger}^i\gr^w\nbigt_{\star}(tf)=0$
unless $i\neq 0$.
In particular,
we have
\[
 \Gr^W
 \pi^0_{\dagger}\nbigt_{\star}(tf)
\simeq
 \pi^0_{\dagger}
 \Gr^W
\nbigt_{\star}(tf),
\quad\quad
  \Gr^W
 \pi^0_{\dagger}\nbigt(tf)[\star t]
\simeq
 \pi^0_{\dagger}
 \Gr^W
\nbigt(tf)[\star t].
\]
Note that it implies that
Condition {\bf (A)} holds for
the projection $\pi$ and 
the mixed twistor $\nbigd$-modules
$\nbigt(tf)[\star t]$ with $\nbigs_{tf}[\star t]$,
and the mixed twistor $\nbigd$-modules
$\nbigt_{\star}(tf)$
with $\vecnbigs_{tf\star}$.
\item
The isomorphisms
$\Gr^W\Psi_{\star}$ $(\star=\ast,!)$
are compatible with the graded polarizations.
\end{itemize}
\end{prop}

Before going to the proof of Proposition \ref{prop;14.9.7.210},
we give a consequence.
Let $\rho:X\lrarr Y$ be
any projective morphism of complex manifolds
such that 
$R^j\rho_{\ast}\Omega_X^i(\ast D)=0$
for any $j>0$ and 
$i=0,1,\ldots,\dim X$.
We also suppose
$\rho_{\dagger}^i\nbigu_X(d,0)[!Z_f][\ast D]=0$
for $i\neq 0$.
Note that it implies
$\rho_{\dagger}^i\nbigu_X(d,0)[\ast Z_f][! D]=0$
for $i\neq 0$.

\begin{cor}\mbox{{}}
Condition {\bf (A)} holds for $\rho$
and $\nbigu_X(d,0)[!Z_f][\ast D]$
with $\nbigs_0[!(f)_0][\ast (f)_{\infty}]$.
We have
\[
 \bigl[
 \rho_{\dagger}^0
 \nbigs_{0}[!(f)_0][\ast (f)_{\infty}]
 \bigr]
=\bigl[
 (\rho\circ\pi)_{\dagger}^0
 \nbigs_{tf}[\ast (t)_0]
 \bigr].
\]
Condition {\bf (A)} holds for $\rho$
and $\nbigu_X(d,0)[\ast Z_f][!D]$
with $\nbigs_0[\ast(f)_0][!(f)_{\infty}]$.
We have
\[
 \bigl[
 \rho_{\dagger}^0
 \nbigs_{0}[\ast(f)_0][!(f)_{\infty}]
 \bigr]
=\bigl[
 (\rho\circ\pi)_{\dagger}^0
 \nbigs_{tf}[!(t)_0]
 \bigr].
\]
\end{cor}
\pf
Set $X^{(1)}:=X\times\proj^1_t$
and $D^{(1)}:=(X\times\{0,\infty\})\cup(D\times\proj^1)$.
We have
$R^j(\rho\circ\pi)_{\ast}
 \Omega_{X^{(1)}}^i(\ast D^{(1)})=0$
for any $j>0$ and $i=0,\ldots,\dim X^{(1)}$.
By the assumption and Proposition \ref{prop;14.7.1.1},
we have
$(\rho\circ\pi)_{\dagger}^j\nbigt_{\star}(tf,D^{(1)})=0$
for $j\neq 0$.
According to Corollary \ref{cor;14.11.24.20} and 
Lemma \ref{lem;14.10.14.11},
Condition {\bf(A)} holds for
$\rho\circ\pi$
and $\nbigt(tf)[\star t]=\nbigt_{\star}(tf,D^{(1)})$
with $\nbigs_{tf}[\ast t]$.
Then, the claim follows from Proposition 
\ref{prop;14.9.7.210}.
\hfill\qed

\subsubsection{First claim of Proposition \ref{prop;14.9.7.210}}

\label{subsection;15.11.9.10}

Let us prove the first claim of Proposition \ref{prop;14.9.7.210}.
We argue the case $\star=\ast$.
We obtain the other case by using the Hermitian duality.

We have 
$\Gr^W_j\nbigt(tf)[\ast t]=
\Gr^W_{j}\nbigt_{\ast}(tf)=0$
for $j<d+1$,
and 
$\Gr^W_{d+1}\nbigt(tf)[\ast t]=
\Gr^W_{d+1}\nbigt_{\ast}(tf)=
\nbigt(tf)$.
The support of
$\Gr^W_{j}\nbigt(tf)[\ast t]$
and 
$\Gr^W_{j}\nbigt_{\ast}(tf)$ $(j>d+1)$
are contained in $X\times\{t=0\}$.
Hence, it is enough to prove
$\pi_{\dagger}^i\nbigt(tf)=0$
for $i\neq 0$.

Let $A_1$ and $A_2$
denote the kernel of
$\nbigt(tf)[!t]\lrarr
 \nbigt(tf)$
and the cokernel of 
$\nbigt(tf)\lrarr
 \nbigt(tf)[\ast t]$.
The support of $A_i$
are contained in $\{t=0\}$.
Hence, we obtain
$\pi_{\dagger}^i\nbigt(tf)=0$
$(i\neq 0)$,
and the following exact sequences:
\[
 0\lrarr
 \pi_{\dagger}^0A_1
\lrarr
 \pi_{\dagger}^0\nbigt(tf)[!t]
\lrarr
 \pi_{\dagger}^0\nbigt(tf)
\lrarr 0
\]
\[
 0\lrarr
 \pi_{\dagger}^0\nbigt(tf)
\lrarr
 \pi_{\dagger}^0\nbigt(tf)[\ast t]
\lrarr
 \pi_{\dagger}^0A_2
\lrarr 0
\]
In particular,
we obtain the first claim of Proposition 
\ref{prop;14.9.7.210}.

\subsubsection{Preliminary for the second claim}

Let us study the second claim.
Let $A_3$ denote the cokernel of 
$\nbigt_{\ast}(tf)\lrarr \nbigt(tf)[\ast t]$.
Let $\iota_Z:Z\lrarr X$ denote the inclusion,
and let $D_Z:=D\cap Z$.
According to Proposition \ref{prop;14.11.7.21},
we have the following isomorphisms
of integrable mixed twistor $D$-modules
with real structure:
\[
 \Psi_1:
 \pi_{\dagger}^0\nbigt_{\ast}(tf)
\simeq
 \iota_{Z\dagger}\nbigu_Z(d-1,0)[\ast D_Z]
 \otimes\newTate(-1),
\quad\quad
\Psi_2:
 \pi_{\dagger}^0A_3\simeq
 \nbigu_X(d,0)\otimes\newTate(-1)[\ast D]
\]
Note that 
the strict support of
$\Gr^W
\Bigl(
 \iota_{Z\dagger}\nbigu_Z(d-1,0)[\ast D_Z]\otimes\newTate(-1)
\Bigr)$
are contained in $Z$,
and 
that the strict support of any direct summand of
$\Gr^W\nbigu_X(d,0)\otimes\newTate(-1)[\ast D]$
are not contained in $Z$.
Hence, we have the unique splittings
\begin{equation}
 \label{eq;15.11.10.1}
 \Gr^W\pi^0_{\dagger}\nbigt(tf)[\ast t]
=\Gr^W\pi^0_{\dagger}\nbigt_{\ast}(tf)
\oplus
 \Gr^W\pi^0_{\dagger}A_3,
\end{equation}
\begin{equation}
 \label{eq;15.11.10.2}
 \Gr^W\Bigl(
 \nbigu_X(d,0)\otimes\newTate(-1)[!Z][\ast D]
 \Bigr)
=
 \Gr^W\Bigl(
 \iota_{Z\dagger}\nbigu_Z(d-1,0)[\ast D_Z]\otimes\newTate(-1)
 \Bigr)
\oplus
 \Gr^W\nbigu_X(d,0)[\ast D]\otimes\newTate(-1),
\end{equation}
and we have
$\Gr^W\Psi_{\ast}=\Gr^W\Psi_1\oplus \Gr^W\Psi_2$.
The graded polarizations are compatible
with the decompositions
(\ref{eq;15.11.10.1}, \ref{eq;15.11.10.2}).

By Proposition \ref{prop;14.9.7.42},
$\Gr^W\Psi_{2|X\setminus Z}$ is compatible
with the graded polarizations.
Hence, we obtain that
$\Gr^W\Psi_{2}$ is compatible
with the graded polarizations.
Therefore, it is enough to study that
$\Gr^W\Psi_1$ is compatible with the graded polarizations.

We remark that
the graded polarization
on 
$\Gr^W\Bigl(
 \iota_{Z\dagger}\nbigu_Z(d-1,0)[\ast D_Z]\otimes\newTate(-1)
 \Bigr)$
is equal to the graded polarization
induced by the canonical polarization 
of $\nbigu_Z(d-1,0)$,
which follows from Lemma \ref{lem;15.11.10.210}.

\subsubsection{The degree $(d+1)$-part}
\label{subsection;14.11.7.24}

Let us study the isomorphism
$\Gr^W_{d+1}\Psi_1:
 \pi_{\dagger}^0\nbigt(tf)
\simeq
 \iota_{Z\dagger}\nbigu_Z(d-1,0)\otimes\newTate(-1)$.
We may assume that $D$ is empty.
To describe $\Gr^W_{d+1}\Psi_1$,
let us look at
the isomorphisms of the $\nbigr$-modules
underlying $\Psi_{\ast}$ more closely,
by assuming $D=\emptyset$.

Set $Y:=X\times\proj^1$.
We have the following:
\[
 \pi_{\dagger}^0\nbigl(tf)[\ast t]
\simeq
  \Cok\Bigl(
\pi_{\ast}
 \nbigl(tf)[\ast t]
\stackrel{\varphi_1}{\lrarr}
\pi_{\ast}
 \nbigl(tf)[\ast t]
 \cdot (dt/\lambda)
 \Bigr)
\simeq
 \Cok\Bigl(
 \nbigo_{\nbigx}[t]
\stackrel{\varphi_2}{\lrarr}
 \nbigo_{\nbigx}[t]\cdot (dt/t\lambda)
 \Bigr)
\]
Here, 
$\varphi_1(c)=\deldel_t(c)\cdot (dt/\lambda)$
and 
$\varphi_2(h\,t^j)
=\bigl(
 jht^{j-1}+hft^j/\lambda
 \bigr)dt$.
For any local section $c$ of 
$\nbigo_{\nbigx}[t]\cdot dt/t\lambda$,
let $[c]_{\ast}$ denote the induced section of
$\pi_{\dagger}^0\nbigl(tf)[\ast t]$.
The isomorphism
$\pi_{\dagger}^0\bigl(
 \nbigl(tf)[\ast t]\bigr)(\ast Z)
\simeq
 \lambda^{-1}\nbigo_{\nbigx}(\ast Z)$
is given by
$[dt/\lambda t]_{\ast}
 \longleftrightarrow \lambda^{-1}$,
which we can check by using (\ref{eq;15.11.10.100}).

Similarly,
we have the following:
\[
 \pi_{\dagger}^0\bigl(
 \nbigl(tf)[!t]
 \bigr)
\simeq
 \Cok\Bigl(
 \pi_{\ast}
 \nbigl(tf)[!t]
 \stackrel{\varphi_1}{\lrarr}
 \pi_{\ast}
 \nbigl(tf)[!t]\cdot(dt/\lambda)
 \Bigr)
\simeq
 \Cok\Bigl(
 \nbigo_{\nbigx}[t]\cdot t
 \stackrel{\varphi_2}{\lrarr}
 \nbigo_{\nbigx}[t]\cdot(dt/\lambda)
 \Bigr)
\]
Here, $\varphi_i$ are given in the same ways.
For any local section $c$ of
$\nbigo_{\nbigx}[t]\cdot(dt/\lambda)$,
let $[c]_!$ denote the induced section of
$\nbigo_{\nbigx}[t]$.
The isomorphism
$\pi_{\dagger}^0(\nbigl(tf)[!t])(\ast Z)
\simeq
 \nbigo_{\nbigx}(\ast Z)$
is given by 
$[fdt/\lambda]_!\longleftrightarrow 1$,
which we can check by using (\ref{eq;15.11.10.101}).

We have the following:
\[
 \pi_{\dagger}^0\nbigl(tf)
\simeq
 \Cok\Bigl(
 \pi_{\ast}\nbigl(tf)
\stackrel{\varphi_1}{\lrarr}
 \pi_{\ast}\nbigl(tf)
 \cdot(dt/\lambda)
 \Bigr)
\simeq
  \Cok\Bigl(
 \nbigo_{\nbigx}[t]
\stackrel{\varphi_2}{\lrarr}
 \nbigo_{\nbigx}[t]
 \cdot(dt/\lambda)
 \Bigr)
\]
For any local section of
$\nbigo_{\nbigx}[t](dt/\lambda)$,
let $[c]$ denote the induced local section of
$\pi_{\dagger}^0\nbigl(tf)$.
We remark
$\lambda^2\del_{\lambda}[dt/\lambda]=0$.

\begin{lem}
\label{lem;15.11.10.200}
The isomorphism
$\Gr^W_{d+1}\Psi_1:
 \pi_{\dagger}\nbigt(tf)\simeq
 \iota_{Z\dagger}\nbigu_Z(d-1,0)\otimes\newTate(-1)$
is given as a pair of the following morphisms:
\[
 \pi_{\dagger}\nbigl(tf)\lambda^{d+1}
\stackrel{a_1}{\llarr}
 \iota_{Z\dagger}\nbigo_{\nbigz}\lambda^{d},
\quad
 a_1\bigl(
 \iota_{Z\ast}(\lambda^d(dt/\lambda)^{-1})
 \bigr)
=[dt/\lambda]\cdot\lambda^{d+1}
\]
\[
 \pi_{\dagger}\nbigl(tf)
\stackrel{a_2}{\lrarr}
 \iota_{Z\dagger}\nbigo_{\nbigz}\lambda^{-1},
\quad
a_2\bigl([dt/\lambda]\bigr)
=-\iota_{Z\ast}(\lambda^{-1}(dt/\lambda)^{-1})
\]
\end{lem}
\pf
We have the epimorphism
$b_1:
 \pi_{\dagger}^0\nbigl(tf)[!t]
\lrarr
 \pi_{\dagger}^0\nbigl(tf)$
and the monomorphism
$b_2:
 \pi_{\dagger}^0\nbigl(tf)
\lrarr
 \pi_{\dagger}^0\nbigl(tf)[\ast t]$.
We have
$b_1([c]_!)=[c]$
and $b_2([c])=[c]_{\ast}$.
We have
$[dt/\lambda]_{!}=f^{-1}[fdt/\lambda]_{!}$.
For a holomorphic local coordinate
$(x_1,\ldots,x_n)$ for $x_1=f$,
we have 
$[dt/\lambda]_{\ast}
=\deldel_{x_1}[dt/t\lambda]_{\ast}$.
Then, we obtain the claim of the lemma
by construction.
\hfill\qed

\vspace{.1in}
We can easily check that
the isomorphism
$\Gr^W_{d+1}\Psi_1$
is compatible with the polarizations
by the explicit description in 
Lemma \ref{lem;15.11.10.200}.

\subsubsection{Preliminary for the other parts}

Let us study 
$\Gr^W_{d+2+\ell}\Psi_1$ for $\ell\geq 0$.
We have only to consider the issue
locally around any point of $D\cap Z$.
We may assume to have a holomorphic coordinate
$(x_0,\ldots,x_{d-1})$ on $X$ such that
$Z=\{x_0=0\}$
and $D=\bigcup_{i=1}^{m}\{x_i=0\}$,
such that
$f=x_0\prod_{i=1}^mx_i^{k_i}$.
We set $x^{\veck}:=\prod_{i=1}^mx_i^{k_i}$.

We have the isomorphism
\begin{equation}
 \label{eq;15.11.10.300}
 \phi^{(0)}_{x_0}
 \pi_{\dagger}^0\Gr^W_{d+\ell+2}(\nbigt(tf)[\ast t])
\simeq
 P\Gr^W_{d+\ell+2}
 \phi^{(0)}_{x_0}\pi_{\dagger}\psi_t^{(0)}(\nbigt(tf))
\end{equation}
The polarization on 
$\phi^{(0)}_{x_0}
 \pi_{\dagger}^0\Gr^W_{d+\ell+2}(\nbigt(tf)[\ast t])$
is induced by the composite of the following:
\begin{multline}
\label{eq;15.11.10.301}
 \phi_{x_0}^{(0)}
 \pi_{\dagger}
 \psi_t^{(0)}(\nbigt(tf))
\stackrel{(-\nbign)^{\ell}}{\lrarr}
 \phi_{x_0}^{(0)}
 \pi_{\dagger}
 \psi_t^{(0)}(\nbigt(tf))
 \otimes\newTate(-\ell)
 \\
\stackrel{a_1}{\lrarr}
  \phi_{x_0}^{(0)}
 \pi_{\dagger}
 \psi_t^{(0)}(\nbigt(tf)^{\ast})
 \otimes\newTate(-\ell-d-1)
=\phi_{x_0}^{(0)}
 \pi_{\dagger}
 \psi_t^{(1)}(\nbigt(tf))^{\ast}
 \otimes\newTate(-\ell-d-1)
\\
\stackrel{a_2}{\lrarr}
 \phi_{x_0}^{(0)}
 \pi_{\dagger}
 \psi_t^{(0)}(\nbigt(tf))^{\ast}
 \otimes\newTate(-\ell-d-2) 
\end{multline}
Here, $a_1$ is induced by the polarization
on $\nbigt(tf)$,
and $a_2$ is the canonical isomorphism.

We have the isomorphism
\begin{equation}
\label{eq;15.11.10.302}
 \Gr^W_{d+\ell+2}
 \iota_{Z\dagger}\nbigu_Z(d,-1)[\ast D_Z]
 \simeq
 P\Gr^W_{d+\ell+2}\psi^{(0)}_{x^k}
 \iota_{Z\dagger}\nbigu_Z(d,-1)
\end{equation}
The polarization on 
$\Gr^W_{d+\ell+2}
 \iota_{Z\dagger}\nbigu_Z(d,-1)[\ast D_Z]$
is induced by the composite of the following morphisms:
\begin{multline}
 \label{eq;15.11.10.303}
 \psi^{(0)}_{x^k}
 \iota_{Z\dagger}\nbigu_Z(d,-1)
\stackrel{(-\nbign)^{\ell}}{\lrarr}
  \psi^{(0)}_{x^k}
 \iota_{Z\dagger}\nbigu_Z(d,-1)
 \otimes\newTate(-\ell)
 \\
\stackrel{b_1}{\lrarr}
  \psi^{(0)}_{x^k}\bigl(
 \iota_{Z\dagger}\nbigu_Z(d,-1)^{\ast}
 \bigr)
 \otimes\newTate(-\ell-d-1)
=\psi^{(1)}_{x^k}\bigl(
 \iota_{Z\dagger}\nbigu_Z(d,-1)
 \bigr)^{\ast}
 \otimes\newTate(-\ell-d-1)
 \\
\stackrel{b_2}{\lrarr}
 \psi^{(0)}_{x^k}\bigl(
 \iota_{Z\dagger}\nbigu_Z(d,-1)
 \bigr)^{\ast}
 \otimes\newTate(-\ell-d-2)
\end{multline}
Here, $b_1$ is induced by the polarization of
$\nbigu_Z(d,-1)$,
and $b_2$ is the canonical isomorphism.

We have the following isomorphism
induced by (\ref{eq;15.11.10.300}),
(\ref{eq;15.11.10.302})
and $\Gr^W_{d+\ell+2}\Psi_1$:
\begin{equation}
 \Lambda:
  P\Gr^W_{d+\ell+2}
 \phi^{(0)}_{x_0}\pi_{\dagger}\psi_t^{(0)}(\nbigt(tf))
\simeq
 P\Gr^W_{d+\ell+2}
  \psi^{(0)}_{x^k}
 \iota_{Z\dagger}\nbigu_Z(d,-1)
\end{equation}

We obtain the compatibility of $\Gr^W_{d+\ell+2}\Psi_1$
and the polarizations from the following lemma.

\begin{lem}
\label{lem;15.11.11.30}
We have isomorphisms of $\nbigd$-modules
\begin{equation}
 \Phi_0:
 \phi^{(0)}_{x_0}\pi_{+}\psi_t^{(0)}L(tf)
\simeq
 \psi_{x^k}^{(0)}\iota_{Z+}\nbigo_Z
\end{equation}
\begin{equation}
 \Phi_1:
 \phi^{(0)}_{x_0}\pi_+\psi_t^{(1)}L(tf)
\simeq
 \psi_{x^k}^{(1)}\iota_{Z+}\nbigo_Z
\end{equation}
with the following property.
\begin{description}
\item[(P1)]
 $\Phi_0$ induces the isomorphism of the $\nbigd$-modules
 underlying $\Lambda$,
 and $\Phi_1^{-1}$ induces the isomorphism of the $\nbigd$-modules
 underlying $\Lambda^{\ast}$. 
\item[(P2)]
 The following diagram is commutative:
\begin{equation}
 \label{eq;15.11.10.310}
 \begin{CD}
 \phi^{(0)}_{x_0}\pi_{+}\psi_t^{(0)}L(tf)
 @>{d_1}>>
 \phi^{(0)}_{x_0}\pi_{+}\psi_t^{(1)}L(tf)\\
 @V{\Phi_0}VV @V{\Phi_1}VV \\
 \psi_{x^k}^{(0)}\iota_{Z+}\nbigo_Z 
@>{d_2}>>
 \psi_{x^k}^{(1)}\iota_{Z+}\nbigo_Z
 \end{CD}
\end{equation}
Here, $d_1$ and $d_2$ are the morphisms
of $\nbigd$-modules
underlying 
{\rm(\ref{eq;15.11.10.301})}
and 
{\rm(\ref{eq;15.11.10.303})},
respectively.
\end{description}
\end{lem}

\subsubsection{Construction of the morphisms $\Phi_i$}

According to Proposition \ref{prop;14.10.15.4},
we have the isomorphism:
\begin{multline}
 \pi_+\psi^{(a)}_tL(tf)
\simeq
 \Xi^{\prime(a)}_{f^{-1}}\nbigo_X
\simeq
 \Cok\Bigl(
 \Pi^{a,N}_{f^{-1}}\nbigo_X[\ast Z][!D]
\stackrel{-s}{\lrarr}
 \Pi^{a,N}_{f^{-1}}\nbigo_X[!Z][\ast D]
 \Bigr)
 \\
\simeq
  \Ker\Bigl(
 \Pi^{-N,a}_{f^{-1}}\nbigo_X[\ast Z][!D]
\stackrel{-s}{\lrarr}
 \Pi^{-N,a}_{f^{-1}}\nbigo_X[!Z][\ast D]
 \Bigr)
\end{multline}
We have the isomorphism 
of  $\nbigd$-modules
(see Lemma \ref{lem;15.11.11.40}
for the signature):
\begin{equation}
\label{eq;15.11.10.410}
 \phi_{x_0}^{(0)}
 \pi_+\psi^{(0)}_tL(tf)
\simeq
\Cok\Bigl(
 \iota_{Z+}
 \Pi^{0,N}_{x^{\veck}}
 \nbigo_Z[!D_Z]
\stackrel{-1}{\lrarr}
 \iota_{Z+}
 \Pi^{0,N}_{x^{\veck}}
 \nbigo_Z[\ast D_Z]
 \Bigr)
\simeq
 \iota_{Z+}
 \psi_{x^{\veck}}^{(0)}\nbigo_Z
\end{equation}
Let $\Phi_0$ be the composite of the isomorphisms
in (\ref{eq;15.11.10.410}).

We have the isomorphism of $\nbigd$-modules:
\begin{equation}
 \label{eq;15.11.11.3}
 \phi_{x_0}^{(0)}\pi_+\psi_t^{(1)}L(tf)
\simeq
  \Ker\Bigl(
 \iota_{Z+}
 \Pi^{-N,1}_{x^{\veck}}\nbigo_Z[!D_Z]
\stackrel{-1}{\lrarr}
 \Pi^{-N,1}_{x^{\veck}}\nbigo_Z[\ast D_Z]
 \Bigr)
\simeq
 \iota_{Z+}
 \psi^{(1)}_{x^{\veck}}\nbigo_Z
\end{equation}
We define $\Phi_1$ as the composite of
(\ref{eq;15.11.11.3}).

\subsubsection{Property {\bf (P1)}}

Let $\Cok$ denote the cokernel of
$L(tf)[!t]\lrarr L(tf)[\ast t]$.
We have the following commutative diagram:
\begin{equation}
 \label{eq;15.11.11.11}
 \begin{CD}
 \Pi^{0,N}_tL(tf)[\ast t]
 @>>>
 L(tf)[\ast t] \\
 @VVV @VVV \\
 \psi_t^{(0)}L(tf)
 @>>>
 \Cok
 \end{CD}
\end{equation}
We apply $\phi_{x_0}^{(0)}\pi_+$
to (\ref{eq;15.11.11.11}).
By the construction of the isomorphism $\Phi_0$,
 (\ref{eq;15.11.11.11}) is transformed to the following:
\[
 \begin{CD}
 \iota_{Z+}\Pi^{0,N}_{x^{\veck}}\nbigo_Z[\ast D_Z]
 @>>>
 \iota_{Z+}\nbigo_Z[\ast D_Z]
 \\
 @VVV @VVV \\
 \iota_{Z+}\psi^{(0)}_{x^{\veck}}\nbigo_Z
 @>>>
 \iota_{Z+}\bigl(
 \nbigo_Z[\ast D_Z]/\nbigo_Z
 \bigr)
 \end{CD}
\]
It means that $\Phi_0$
is the morphism of $\nbigd$-modules
underlying $\Lambda$.

Similarly, 
let $\Ker_1$ denote the kernel of 
$L(tf)[!t]\lrarr L(tf)[\ast t]$.
Then, we have the following commutative diagram:
\begin{equation}
\label{eq;15.11.11.12}
 \begin{CD}
 \Ker @>>> \psi_t^{(1)} L(tf)
 \\
 @AAA @AAA\\
 L(tf)[!t] @>>> \Pi^{-N,1}L(tf)[!t]
 \end{CD}
\end{equation}
We apply $\phi_{x_0}^{(0)}\pi_+$
to (\ref{eq;15.11.11.12}).
Let $\Ker_2$ denote the kernel of
 $\iota_+\nbigo_Z[! D_Z]\lrarr \nbigo_Z$.
Then, $(\ref{eq;15.11.11.12})$ is transformed to
the following:
\[
 \begin{CD}
 \Ker_2
 @>>>
 \iota_{Z+}\psi_{x^{\veck}}^{(1)}\nbigo_Z
 \\
 @AAA @AAA \\
 \iota_{Z+}\nbigo_Z[!D_Z]
 @>>>
 \iota_{Z+}
 \Pi^{-N,1}_{x^{\veck}}\nbigo_Z[!D_Z]
 \end{CD}
\]
It implies that $\Phi_1$ underlies $\Lambda^{\ast}$.
Thus, we obtain that
{\bf(P1)} is satisfied.

\subsubsection{Property {\bf (P2)}}

The morphism $d_1$ is induced by the following
morphisms:
\begin{equation}
 \psi_{t}^{(0)}L(tf)
 \stackrel{(-s)^{\ell}}{\lrarr}
 \psi_t^{(0)}L(tf)
 \stackrel{(-1)^{d+1}}{\lrarr}
 \psi_t^{(0)}L(tf)
\stackrel{-s}\simeq
 \psi_t^{(1)}L(tf)
\end{equation}
The morphism $d_2$ is induced by the following:
\begin{equation}
 \psi^{(0)}_{x_k}\iota_{Z+}\nbigo_Z
 \stackrel{(-s)^{\ell}}{\lrarr}
 \psi^{(0)}_{x_k}\iota_{Z+}\nbigo_Z
 \stackrel{(-1)^{d}}{\lrarr}
 \psi^{(0)}_{x_k}\iota_{Z+}\nbigo_Z
 \stackrel{-s}{\lrarr}
 \psi^{(1)}_{x_k}\iota_{Z+}\nbigo_Z
\end{equation}
We remark that 
the following diagram is commutative:
\[
\begin{CD}
 \phi_{x_0}^{(0)}
 \pi_+\psi^{(0)}_{t}L(tf)
 @>{s}>>
  \phi_{x_0}^{(0)}
 \pi_+\psi^{(1)}_{t}L(tf)
 \\
 @V{\Phi_0}VV @V{\Phi_1}VV \\
 \psi^{(0)}_{x^{\veck}}
 \iota_{Z+}\nbigo_Z
 @>{-s}>>
 \psi^{(1)}_{x^{\veck}}
 \iota_{Z+}\nbigo_Z
\end{CD}
\]
The signature comes from $-1$ in
(\ref{eq;15.11.10.410}) and (\ref{eq;15.11.11.3}).
Hence, we obtain the desired commutativity of
(\ref{eq;15.11.10.310}).
Thus, we finish the proof of Lemma \ref{lem;15.11.11.30}
and Proposition \ref{prop;14.9.7.210}.
\hfill\qed

\section{Mixed twistor $\nbigd$-modules
and GKZ-hypergeometric systems}
\label{section;14.12.8.11}
\subsection{Better behaved GKZ-hypergeometric systems}
\label{subsection;14.7.4.10}

\subsubsection{The associated toric varieties}
\label{subsection;14.7.8.1}

We consider a finite subset
$\nbiga:=\{\veca_1,\ldots,\veca_m\}\subset 
 \seisuu^n\setminus\{(0,\ldots,0)\}$,
where $\veca_i=(a_{ij}\,|\,j=1,\ldots,n)$.
We assume that $\nbiga$ generates $\seisuu^n$.
Formally, we set $\veca_0:=(0,\ldots,0)$.

We set $T^n:=\Spec\cnum[t_1^{\pm 1},\ldots,t_n^{\pm 1}]$.
The associated complex manifold
is also denoted by $T^n$
if there is no risk of confusion.
We consider a morphism
$\psi_{\nbiga}:T^n\lrarr \proj^{m}$
given by
$\psi_{\nbiga}(t_1,\ldots,t_n)
=[1:t^{\veca_1}:\ldots:t^{\veca_m}]$,
where
 $t^{\veca_i}:=\prod_{j=1}^nt_j^{a_{ij}}$.
Let $X_{\nbiga}$ denote the closure of
$\psi_{\nbiga}(T^n)$ in $\proj^{m}$.

For any subset $\nbigb\subset\real^n$,
let $\Conv(\nbigb)$ denote the convex hull of
$\nbigb\subset\real^n$.
Note that
$0$ is an interior point of
$\Conv(\nbiga\cup\{0\})$
if and only if
$0$ is an interior point of
$\Conv(\nbiga)$.

For any face $\sigma$ of $\Conv(\nbiga\cup\{0\})$,
we put $I_{\sigma}:=
 \sigma\cap\bigl(\nbiga\cup\{0\}\bigr)$
and $J_{\sigma}:=(\nbiga\cup\{0\})\setminus I_{\sigma}$.
We have the subspace
\[
 \proj_{\sigma}:=\bigl\{
 [z_0:z_1:\cdots:z_m]\,\big|\,
 z_{j}=0\,\,(a_j\in J_{\sigma})
 \bigr\}.
\]
We set
$\proj_{\sigma}^{\ast}:=
 \bigl\{
 [z_0:\ldots:z_m]\in
 \proj_{\sigma}\,\big|\,
 z_i\neq 0\,\, (i\in I_{\sigma}),
 z_i=0\,\,(i\in J_{\sigma})
 \bigr\}$.
Recall the following.
\begin{prop}[Proposition 1.9 
 \cite{GKZ-book}]
For any face $\sigma$ of $\Conv(\nbiga\cup\{0\})$,
we have the non-empty intersection
$X_{\nbiga}^{\sigma}=X_{\nbiga}\cap \proj_{\sigma}^{\ast}$,
and it is a $T^n$-orbit
of $[z_0^{\sigma}:\cdots:z_m^{\sigma}]$,
where $z_i^{\sigma}=1$ $(i\in I_{\sigma})$
or $z_i^{\sigma}=0$ $(i\in J_{\sigma})$.
We have the decomposition into the orbits
$X_{\nbiga}
=\coprod_{\sigma} X_{\nbiga}^{\sigma}$
where $\sigma$ runs through
the set of faces of $\Conv(\nbiga\cup\{0\})$.
\hfill\qed
\end{prop}

Recall that we have a desingularization
of $X_{\nbiga}$ in the category of toric varieties.
(See \cite{Cox-et-al-toric} for example.)
Namely, there exist a smooth fan $\Sigma$
and a toric birational projective morphism $\varphi_{\Sigma}:
 X_{\Sigma}\lrarr X_{\nbiga}$,
where $X_{\Sigma}$
denotes the smooth toric variety associated to $\Sigma$.
In that situation,
we set $D_{\Sigma}:=\bigcup_{\tau\in\Sigma(1)} D_{\tau}$,
where $\Sigma(1)$ denotes the set of
the $1$-dimensional cones of $\Sigma$,
and $D_{\tau}\subset X_{\Sigma}$ 
denotes the hypersurface corresponding to $\tau$.
We have $T^n=X_{\Sigma}\setminus D_{\Sigma}$.

\subsubsection{$\nbigd$-modules associated to
families of Laurent polynomials}
\label{subsection;14.11.23.210}

Suppose that
we are given an algebraic map
$\gamma:
 S\lrarr H^0(\proj^m,\nbigo(1))$.
It determines a family of Laurent polynomials
$(\psi_{\nbiga}\times\id_S)^{\ast}(s/z_0)$ on $T^n\times S$
denoted by $F_{\gamma}$.
We obtain the algebraic $\nbigd$-module 
$L_{T^n\times S}(F_{\gamma})$
on $T^n\times S$
given as the line bundle
$\nbigo_{T^n\times S}\cdot e$
with $\nabla e=e\,dF_{\gamma}$,
where $e$ denotes a global frame.
We obtain the following algebraic $\nbigd$-modules
on $\proj^m\times S$:
\[
 L_{\star}(\nbiga,S,\gamma):=
 (\psi_{\nbiga}\times\id_S)_{\star}L_{T^n\times S}(F_{\gamma})
\quad
(\star=!,\ast)
\]
The image of 
$L_{!}(\nbiga,S,\gamma)
\lrarr
 L_{\star}(\nbiga,S,\gamma)$
is denoted by
$L_{\min}(\nbiga,S,\gamma)$,
which is the minimal extension of
$L_{T^n\times S}(F_{\gamma})$
on $\proj^m\times S$
via $\psi_{\nbiga}\times\id_S$.
Let $\pi:\proj^m\times S\lrarr S$ denote the projection.
We obtain the algebraic $\nbigd$-modules
$\pi^i_{+}L_{\star}(\nbiga,S,\gamma)$ $(\star=\ast,!)$
and 
$\pi^i_{+}L_{\min}(\nbiga,S,\gamma)$.

\vspace{.1in}

It is convenient for us to take a toric desingularization
$\varphi_{\Sigma}:X_{\Sigma}\lrarr X_{\nbiga}$
when we work in the complex analytic setting
or when we are interested in the twistor property
of the $\nbigd$-modules.
We set 
$(X_{\Sigma,S},D_{\Sigma,S})
 :=(X_{\Sigma},D_{\Sigma})\times S$.
Let $\varphibar_{\Sigma,S}$
denote $X_{\Sigma,S}\lrarr \proj^m\times S$
which is the composite of
$\varphi_{\Sigma}\times\id_S$
and the inclusion 
$X_{\nbiga}\times S\lrarr \proj^m\times S$.
We have the meromorphic function
$F_{\gamma,\Sigma}:=
 \varphibar_{\Sigma,S}^{\ast}(s)/\varphibar_{\Sigma,S}^{\ast}(z_0)$
on $(X_{\Sigma},D_{\Sigma})\times S$.
We obtain the $\nbigd_{X_{\Sigma,S}}$-modules
$L_{\star}(F_{\gamma,\Sigma},D_{\Sigma,S})$ 
$(\star=\ast,!)$ on $X_{\Sigma,S}$.
Then, we have
$\varphibar_{\Sigma,S+}^i
 L_{\star}(F_{\gamma,\Sigma},D_{\Sigma,S})=0$
for $i\neq 0$.
In the algebraic case,
we have
$\varphibar_{\Sigma,S+}^0
 L_{\star}(F_{\gamma,\Sigma},D_{\Sigma,S})
\simeq
 L_{\star}(\nbiga,S,\gamma)$.
In the complex analytic setting,
we adopt it as the definition,
i.e.,
$L_{\star}(\nbiga,S,\gamma):=
 \varphibar_{\Sigma,S+}^0
 L_{\star}(F_{\gamma,\Sigma},D_{\Sigma,S})$.
Let $\pi_{\Sigma}:X_{\Sigma}\times S\lrarr S$
denote the projection.
Then, we have the following natural isomorphisms
for $\star=\ast,!$:
\begin{equation}
\label{eq;14.11.6.1}
 \pi_{\Sigma+}
 L_{\star}(F_{\gamma,\Sigma},D_{\Sigma,S})
\simeq
 \pi_{+}
  L_{\star}(\nbiga,S,\gamma)
\end{equation}

If we are given another desingularization
$\varphi_{\Sigma'}:X_{\Sigma'}\lrarr X_{\nbiga}$,
we can find a smooth fan $\Sigma''$
and toric morphisms
$\psi_1:X_{\Sigma''}\lrarr X_{\Sigma}$
and $\psi_2:X_{\Sigma''}\lrarr X_{\Sigma'}$
such that
$\varphi_{\Sigma}\circ\psi_1
=\varphi_{\Sigma'}\circ\psi_2=:
\varphi_{\Sigma''}$.
We have natural isomorphisms
$(\psi_1\times\id_S)_{+}
 L_{\star}(F_{\gamma,\Sigma''},D_{\Sigma'',S})
\simeq
(\psi_1\times\id_S)^0_{+}
 L_{\star}(F_{\gamma,\Sigma''},D_{\Sigma'',S})
\simeq
 L_{\star}(F_{\gamma,\Sigma},D_{\Sigma,S})$.
We have similar isomorphisms
for $L_{\star}(F_{\gamma,\Sigma'},D_{\Sigma',S})$.
Hence, we have the following natural isomorphisms
on $\proj^m\times S$:
\[
 \varphibar^0_{\Sigma+}
 L_{\star}(F_{\gamma,\Sigma},D_{\Sigma,S})
\simeq
 \varphibar^0_{\Sigma''+}
 L_{\star}(F_{\gamma,\Sigma''},D_{\Sigma'',S})
\simeq
 \varphibar^0_{\Sigma'+}
 L_{\star}(F_{\gamma,\Sigma'},D_{\Sigma',S})
\]
In this sense,
the $\nbigd$-modules
$L_{\star}(\nbiga,S,\gamma)$
are independent of 
the choice of a resolution.
We also have the following natural isomorphisms:
\begin{equation}
\label{eq;14.7.10.11}
 \pi^i_{\Sigma+}L_{\star}(F_{\gamma,\Sigma},D_{\Sigma,S})
\simeq
\pi^i_{\Sigma''+}L_{\star}(F_{\gamma,\Sigma''},D_{\Sigma'',S})
\simeq
 \pi^i_{\Sigma'+}L_{\star}(F_{\gamma,\Sigma'},D_{\Sigma',S})
\simeq
  \pi^i_{+}L_{\star}(\nbiga,S,\gamma)
\end{equation}

\begin{rem}
Although it is convenient for us to 
take a toric desingularization $X_{\Sigma}$ of $X_{\nbiga}$,
we can also use any toric completion of $T^n$
for the study of
$\pi^i_{+}L_{\star}(\nbiga,S,\gamma)$.
Let $X_{\Sigma_1}$ be any $n$-dimensional smooth toric variety.
We fix an inclusion
$T^n\subset X_{\Sigma_1}$.
The family of Laurent polynomials $F_{\gamma}$
gives a meromorphic function
$F_{\gamma,\Sigma_1}$
on $(X_{\Sigma_1,S},D_{\Sigma_1,S})$.
Then, we naturally have
$\pi_{\Sigma_1+}L_{\star}(F_{\gamma,\Sigma_1},D_{\Sigma_1})
\simeq
 \pi_{\Sigma+}L_{\star}(F_{\gamma,\Sigma},D_{\Sigma})
\simeq
 \pi_+L_{\star}(\nbiga,S,\gamma)$.
\hfill\qed
\end{rem}

Let $\gamma_i:S\lrarr H^0(\proj^m,\nbigo(1))$
$(i=1,2)$ be holomorphic maps.
We give some conditions 
under which we have
$\pi^i_{+}L_{\star}(\nbiga,S,\gamma_1)
\simeq
 \pi^i_{+}L_{\star}(\nbiga,S,\gamma_2)$.
We have the $T^n$-action on $\proj^m$
given by
$(t_1,\ldots,t_n)[z_0:\cdots:z_m]
=[t^{\veca_0}z_0:\cdots:t^{\veca_m}z_m]$.
For $b\in T^n$,
let $m_b$ denote the induced action on $\proj^m$,
and let $m_b^{\ast}$ denote the induced action on
$H^0\bigl(\proj^m,\nbigo(1)\bigr)$.

\begin{lem}
\label{lem;14.12.30.10}
Suppose that there exists a holomorphic map
$\kappa:S\lrarr T^n$ 
such that 
$\gamma_1=m_{\kappa}^{\ast}\gamma_2$,
i.e.,
$\gamma_1(x)=m_{\kappa(x)}^{\ast}\gamma_2(x)$
for any $x\in S$.
Then, we have the natural isomorphisms
$\pi_{+}L_{\star}(\nbiga,S,\gamma_1)
\simeq
 L_{\star}(\nbiga,S,\gamma_2)$
for $\star=\ast,!$.
If $\gamma_i$ and $\kappa$ are algebraic,
then we have the isomorphisms 
of algebraic $\nbigd$-modules.
\end{lem}
\pf
We have the isomorphism
$\widetilde{m}_{\kappa}:
 \proj^m\times S\lrarr \proj^m\times S$
induced by $\kappa$ and the identity of $S$.
We have an isomorphism
$ \widetilde{m}_{\kappa+}L_{\star}(\nbiga,S,\gamma_1)
\simeq
 L_{\star}(\nbiga,S,\gamma_2)$.
Applying $\pi_{+}$,
we obtain the desired isomorphisms.
\hfill\qed

\vspace{.1in}
The following lemma is clear by construction.

\begin{lem}
Let $\rho:S\lrarr \cnum$ be a holomorphic map
such that 
$\gamma_1=\gamma_2+\rho\cdot z_0$.
Then, we have the isomorphisms of 
the analytic $\nbigd$-modules
$L_{\star}(\nbiga,S,\gamma_1)
\simeq
 L_{\star}(\nbiga,S,\gamma_2)$.
In particular,
we have the isomorphisms of 
the analytic $\nbigd$-modules
$\pi_{+}^iL_{\star}(\nbiga,S,\gamma_1)
\simeq
 \pi_{+}^iL_{\star}(\nbiga,S,\gamma_2)$.
\hfill\qed
\end{lem}

\subsubsection{Non-degenerate families of sections}
\label{subsection;14.11.23.240}

Let $p_1:\proj^m\times S\lrarr S$ be the projection.
Let $\gamma:S\lrarr H^0\bigl(\proj^m,\nbigo_{\proj^m}(1)\bigr)$
be any holomorphic map.
We have the expression
$\gamma(x)=\sum_{i=0}^m \gamma_i(x) z_i$.
It determines a family of Laurent polynomials
$F_{\gamma}(x,t)=\sum_{i=0}^m\gamma_i(x)t^{\veca_i}$.
For any face $\sigma$ of $\Conv(\nbiga\cup\{0\})$,
we consider the associated family of Laurent polynomials
$F_{\gamma,\sigma}(x,t)=
 \sum_{\veca_i\in\sigma}\gamma_i(x)t^{\veca_i}$.
We can regard it as a function on $T^n\times S$.
Let $dF_{\gamma,\sigma}$ denote the exterior derivative 
of $F_{\gamma,\sigma}$,
which is a section of $\Omega^1_{T^n\times S}$.
We use the following non-degeneracy condition
which is a minor generalization
of the classical non-degeneracy condition \cite{Kouchnirenko}
(see also \cite{Adolphson}, \cite{Reichelt-Sevenheck2})
in the sense that we also consider the derivatives
in the $S$-direction.

\begin{df}
\label{df;14.11.23.1}
$\gamma$ is called non-degenerate at $\infty$ for $X_{\nbiga}$
if we have 
$(F_{\gamma,\sigma})^{-1}(0)\cap
 (dF_{\gamma,\sigma})^{-1}(0)=\emptyset$
for any face $\sigma$ of $\Conv(\nbiga\cup\{0\})$
such that $0\not\in\sigma$.
\hfill\qed
\end{df}
The condition is equivalent to that
the zero-divisors of $F_{\gamma,\sigma}$
are smooth and reduced 
for any face $\sigma$ of $\Conv(\nbiga\cup\{0\})$
such that $0\not\in\sigma$.
We remark the following.
\begin{lem}
If $\sigma\subsetneq\Conv(\nbiga\cup\{0\})$
and $0\not\in\sigma$,
we have 
$(F_{\gamma,\sigma})^{-1}(0)\cap
 (dF_{\gamma,\sigma})^{-1}(0)=\emptyset$
if and only if $(dF_{\gamma,\sigma})^{-1}(0)=\emptyset$.
\end{lem}
\pf
Under the assumption,
there exists
$\beta=(\beta_1,\ldots,\beta_n)\in(\rnum^n)^{\lor}$
such that 
$\langle
 \beta,\veca_i
 \rangle=1$
for any $\veca_i\in I_{\sigma}$.
Because
$\Bigl(
 \sum_{j=1}^n\beta_jt_j\del_{t_j}
 \Bigr)t^{\veca_i}
=\langle
 \beta,\veca_i
 \rangle\cdot
 t^{\veca}
=t^{\veca}$,
we have
$F_{\gamma,\sigma}^{-1}(0)
\supset
 \bigl(
 dF_{\gamma,\sigma}
 \bigr)^{-1}(0)$.
Hence, the condition 
$F_{\gamma,\sigma}^{-1}(0)\cap
 (dF_{\gamma,\sigma})^{-1}(0)=\emptyset$
is equivalent to
that $\bigl(dF_{\gamma,\sigma}\bigr)^{-1}(0)=\emptyset$.
\hfill\qed

\vspace{.1in}

Let us reword the condition.
We have the section $s_{\gamma}$
of the line bundle $p_1^{\ast}\nbigo_{\proj^m}(1)$
such that
$s_{\gamma|\proj^m\times \{x\}}
 =\gamma(x)\in H^0(\proj^m,\nbigo_{\proj^m}(1))$
for any $x\in S$.
We set 
$H_{\gamma}:=s_{\gamma}^{-1}(0)
\subset \proj^m\times S$.
We also put $H_{\infty}:=\{z_0=0\}\subset\proj^m$.

\begin{lem}
\label{lem;14.7.1.10}
Let $\gamma:S\lrarr H^0(\proj^m,\nbigo(1))$ be a morphism.
For any face $\sigma$ of $\Conv(\nbiga\cup\{0\})$,
the following conditions are equivalent.
\begin{itemize}
\item
$(F_{\gamma,\sigma})^{-1}(0)
\cap
 (dF_{\gamma,\sigma})^{-1}(0)=\emptyset$.
\item
There exists an open neighbourhood 
$\nbigu_1$ of $X_{\nbiga}^{\sigma}\times S$ in $\proj^m\times S$
such that 
(i)
$H_{\gamma}\cap \nbigu_1$
is smooth
(ii) $H_{\gamma}$ is transversal with 
$X_{\nbiga}^{\sigma}\times S$.
\end{itemize}
\end{lem}
\pf
Let $Q\in X_{\nbiga}^{\sigma}$ be the point 
$[z^{\sigma}_0:\cdots:z^{\sigma}_m]$
as in \S\ref{subsection;14.7.8.1}.
By using the $T^n$-action on $X_{\nbiga}^{\sigma}$,
we consider the map 
$\psi:T^n\times S\lrarr X_{\nbiga}^{\sigma}\times S$
given by $\psi(t,x)=(t\cdot Q,x)$.
Take any $i_0\in I_{\sigma}$.
We have
\[
 \psi^{\ast}(s_{\gamma})/\psi^{\ast}(z_{i_0})
=\sum_{\veca_i\in I_{\sigma}} \gamma_i(x) t^{\veca_i-\veca_{i_0}}
=:G_{\gamma,\sigma}(t,x).
\]
Then, 
the second condition holds
if and only if
(i) $G_{\gamma,\sigma}$ is not constantly $0$,
(ii) the divisor 
$(G_{\gamma,\sigma})=(F_{\gamma,\sigma})$ is smooth.
It is equivalent to that
$(dF_{\gamma,\sigma})^{-1}(0)\cap
 F_{\gamma,\sigma}^{-1}(0)=\emptyset$.
\hfill\qed

\begin{cor}
$\gamma$ is
non-degenerate at $\infty$ for $X_{\nbiga}$
if and only if the following holds:
\begin{itemize}
\item 
There exists a neighbourhood $\nbigu$ of 
$H_{\infty}\times S$ in $\proj^m\times S$
such that $H_{\gamma}\cap \nbigu$ is a smooth hypersurface 
in $\nbigu$.
\item
If $\sigma$ is a face of $\Conv(\nbiga\cup\{0\})$
such that $0\not\in\sigma$,
then $H_{\gamma}$ is transversal with
$X_{\nbiga}^{\sigma}\times S$.
\hfill\qed
\end{itemize}
\end{cor}

The first condition in this corollary is trivial
if the image of $\gamma$ does not contain
$0\in H^0(\proj^m,\nbigo(1))$.
The conditions particularly imply
that $\gamma(U)\not\subset \cnum\cdot z_0$
for any open subset $U\subset S$.

\vspace{.1in}

Let $\varphi_{\Sigma}:X_{\Sigma}\lrarr X_{\nbiga}$
be any toric desingularization.
Let $\varphibar_{\Sigma}:X_{\Sigma}\lrarr \proj^m$
be the composite of $\varphi_{\Sigma}$
and the inclusion $X_{\nbiga}\lrarr \proj^m$.

\begin{lem}
\label{lem;15.1.3.2}
Let $\gamma:S\lrarr H^0(\proj^m,\nbigo(1))$
be a holomorphic map
which is non-degenerate at $\infty$ for $X_{\nbiga}$.
Then, $F_{\gamma,\Sigma}$
is non-degenerate along $D_{\Sigma}\times S$.
We also have
$|(F_{\gamma,\Sigma})_{\infty}|=
 \varphibar_{\Sigma}^{-1}(H_{\infty})\times S$.
\end{lem}
\pf
Let $T_{\rho}$ be a $T^n$-orbit in $X_{\Sigma}$
contained in 
$\varphi_{\Sigma}^{-1}(H_{\infty})$.
Note that 
$\varphi_{\Sigma}(T_{\rho})$
is contained in $X_{\Sigma}\cap H_{\infty}$,
and the morphism
$T_{\rho}\lrarr \varphi_{\Sigma}(T_{\rho})$
is smooth.
Let $\varphibar_{\Sigma,S}:=\varphibar_{\Sigma}\times\id_S$.
The restriction of
$\varphibar_{\Sigma,S}^{\ast}(s_{\gamma})$
to $T_{\rho}\times S$ is not $0$,
and the $0$-divisor of
$\varphibar_{\Sigma,S}^{\ast}(s_{\gamma})_{|T_{\rho}\times S}$
is smooth.
Then, the claims of the lemma are easy to see.
\hfill\qed

\vspace{.1in}

Any element $s\in H^0(\proj^m,\nbigo(1))$
determines a morphism $\gamma_s$
from the one-point set to $H^0(\proj^m,\nbigo(1))$.
\begin{df}
We say that $s$ is non-degenerate at $\infty$ for $X_{\nbiga}$
if $\gamma_s$ is non-degenerate at $\infty$ for $X_{\nbiga}$
in the sense of Definition {\rm\ref{df;14.11.23.1}}.
\hfill\qed
\end{df}

It is standard that
there exist non-empty Zariski open subsets
$U\subset H^0(\proj^m,\nbigo_{\proj^m}(1))$
such that any $s\in U$ is non-degenerate
at $\infty$ for $X_{\nbiga}$.
Let $U_{\nbiga}^{\reg}$ be the union
of such open subsets.

\subsubsection{Basic examples of non-degenerate
family of sections}
\label{subsection;15.11.22.11}

The following lemma is clear.
\begin{lem}
A holomorphic map $\gamma:S\lrarr U_{\nbiga}^{\reg}$
is non-degenerate at $\infty$ for $X_{\nbiga}$.
\hfill\qed
\end{lem}

Note that 
even if $\gamma:S\lrarr H^0(\proj^m,\nbigo(1))$
is non-degenerate at $\infty$ for $X_{\nbiga}$,
the image of $\gamma$ is not necessarily contained in $U_{\nbiga}^{\reg}$.
The following lemma is clear.
\begin{lem}
\label{lem;14.12.30.20}
Let $\gamma:S\lrarr H^0(\proj^m,\nbigo(1))$
be a holomorphic map.
Let $\Phi_{\gamma}:H_{\gamma}\lrarr \proj^m$
be the induced morphism,
i.e., the composite of the inclusion
$H_{\gamma}\lrarr \proj^m\times S$
and the projection $\proj^m\times S\lrarr \proj^m$.
Suppose that $H_{\gamma}$ is a smooth hypersurface,
and that any critical value of $\Phi_{\gamma}$ 
is not contained in 
$X_{\nbiga}\cap H_{\infty}$.
Then, $\gamma$ is 
non-degenerate at $\infty$
for $X_{\nbiga}$.
\hfill\qed
\end{lem}

\begin{example}
\label{example;14.11.24.100}
We set 
$W:=\bigl\{
 \sum_{i=1}^m\alpha_iz_i
 \bigr\}\subset
 H^0(\proj^m,\nbigo(1))$.
Then,
the inclusion
$\iota:W\lrarr H^0(\proj^m,\nbigo(1))$
is non-degenerate at $\infty$ for $X_{\nbiga}$
for any $\nbiga$.

Indeed, we set $\Htilde_{\iota}\subset
 \bigl(
 \cnum^{m+1}\setminus\{(0,\ldots,0)\}
 \bigr)\times W:=
 \Bigl\{(z_0,\ldots,z_m;\alpha_1,\ldots,\alpha_m)\,\Big|\,
 \sum_{i=1}^m\alpha_iz_i=0
 \Bigr\}$.
Then, $\Htilde_{\iota}$ is clearly smooth.
Because $H_{\iota}$ is the quotient of $\Htilde_{\iota}$
by the natural free $\cnum^{\ast}$-action,
$H_{\iota}$ is also smooth.
Let $\Htilde_{\iota}\lrarr
 \cnum^{m+1}\setminus\{(0,\ldots,0)\}$
be the projection.
Then, the critical values are 
$\{(z_0,0,\ldots,0)\,|\,z_0\in\cnum^{\ast}\}$.
Hence, the set of the critical values of
$\Phi_{\iota}:H_{\iota}\lrarr \proj^m$
is $\{[1:0:\cdots:0]\}$.
Then, by Lemma {\rm\ref{lem;14.12.30.20}},
we obtain that $\iota$ is non-degenerate.
\hfill\qed
\end{example}

Set
$W^{\ast}:=\bigl\{
 \sum_{i=1}^{m}\alpha_iz_i\,\big|\,
 \alpha_i\in\cnum^{\ast}
 \bigr\}$.
We have the action of $T^n$ on $W^{\ast}$ by 
$t\cdot (\alpha_1,\ldots,\alpha_m)
=(t^{\veca_1}\alpha_1,\ldots,t^{\veca_m}\alpha_m)$.
The action is free because $\nbiga$ generates $\seisuu^n$.
Let $q:W^{\ast}\lrarr W^{\ast}/T^n$
be the projection.
\begin{lem}
Let $\gamma:S\lrarr W^{\ast}$ be any holomorphic map
such that 
$q\circ\gamma$ is submersive.
Then, 
$\gamma$ is non-degenerate at $\infty$
for $X_{\nbiga}$.
\end{lem}
\pf
Let $\sigma$ be any face of 
$\Conv(\nbiga\cup\{0\})$
such that $0\not\in\sigma$.
Take any $(P,x)\in H_{\gamma}\subset \proj^m\times S$
such that $P\in X_{\nbiga}^{\sigma}$.
We have $P\in H_{\gamma(x)}$.
Let $\Delta_{\epsilon}:=\bigl\{z\in\cnum\,\big|\,|z|<\epsilon\bigr\}$
for $\epsilon>0$.
For any holomorphic map
$g:\Delta_{\epsilon}\lrarr W^{\ast}$,
we have the family of hyperplanes
$H_g\subset\proj^m\times \Delta_{\epsilon}$,
and we obtain the induced map
$\Phi_g:H_{g}\lrarr \proj^m$.
Because $P\neq [1:0:\cdots:0]$,
we can find a holomorphic map
$g_0:\Delta\lrarr W^{\ast}$
such that (i) $g_0(0)=\gamma(x)$,
(ii) the tangent map of $\Phi_{g_0}$ 
at $(P,0)$ is an isomorphism.
For some $0<\epsilon'<\epsilon$,
we can find holomorphic maps
$g_1:\Delta_{\epsilon'}\lrarr S$
and 
$h_1:\Delta_{\epsilon'}\lrarr T^n$
such that
(i) $g_1(0)=x$,
(ii) $h_1(0)=1$,
(iii) $\gamma\circ g_1(x)=h_1(x)\cdot g_0(x)$
 for $x\in \Delta_{\epsilon'}$.
Then, it is easy to see that
the tangent map of
$\Phi_{\gamma}$
at $(P,x)$ is surjective.
\hfill\qed

\begin{example}
\label{example;14.11.24.101}
Let $r:W^{\ast}/T^n\lrarr W^{\ast}$ be any section of $q$.
Then, 
the induced morphism
$r:W^{\ast}/T^n\lrarr 
 H^0\bigl(\proj^m,\nbigo(1)\bigr)$
is non-degenerate.
\hfill\qed
\end{example}

\subsubsection{$\nbigd$-modules
associated to non-degenerate families of sections}
\label{subsection;14.10.15.10}

Let $\gamma:S\lrarr H^0\bigl(\proj^m,\nbigo(1)\bigr)$
be any holomorphic map.
Suppose that $\gamma$ is non-degenerate at $\infty$
for $X_{\nbiga}$.
We have the $\nbigd$-modules
$L_{\star}(\nbiga,S,\gamma)$
on $\proj^m\times S$.
We set $H_{\infty,S}:=H_{\infty}\times S$.

\begin{lem}
\label{lem;14.6.28.1}
We naturally have 
$L_{\star}(\nbiga,S,\gamma)(\ast H_{\infty,S})
\simeq
  L_{\star}(\nbiga,S,\gamma)$
$(\star=\ast,!)$.
\end{lem}
\pf
For a desingularization
$\varphi_{\Sigma}:X_{\Sigma}\lrarr X_{\nbiga}$,
we obtain the $\nbigd$-modules
$L_{\star}(F_{\gamma,\Sigma},D_{\Sigma,S})$.
Because $F_{\gamma,\Sigma}$ is non-degenerate
along $D_{\Sigma,S}$,
we have
\[
 L_{\star}(\nbiga,S,\gamma)(\ast H_{\infty,S})
\simeq
 \varphibar^0_{\Sigma,S+}\Bigl(
 L_{\star}(F_{\gamma,\Sigma},D_{\Sigma,S})
  \bigl(\ast (F_{\gamma,\Sigma})_{\infty}\bigr)
 \Bigr)
\simeq
 \varphibar^0_{\Sigma,S+}\bigl(
 L_{\star}(F_{\gamma,\Sigma},D_{\Sigma,S})
 \bigr)
\simeq
  L_{\star}(\nbiga,S,\gamma).
\]
Thus, we are done.
\hfill\qed

\vspace{.1in}

Let $K_{\nbiga,S,\gamma}$ and $C_{\nbiga,S,\gamma}$
denote the kernel and the cokernel of
$L_{!}(\nbiga,S,\gamma)
\lrarr
 L_{\ast}(\nbiga,S,\gamma)$.
\begin{cor}
\label{cor;14.6.28.10}
We have
$M(\ast H_{\infty,S})
\simeq M$
for the $\nbigd$-modules
$M=K_{\nbiga,S,\gamma},C_{\nbiga,S,\gamma},
 L_{\min}(\nbiga,S,\gamma)$.
\hfill\qed
\end{cor}

\begin{prop}
\label{prop;14.7.4.20}
We have
$\pi_+^iM=0$ $(i\neq 0)$
for the $\nbigd$-modules
\[
 M=L_{\star}(\nbiga,S,\gamma)\,\,(\star=\ast,!),
 L_{\min}(\nbiga,S,\gamma),
 K_{\nbiga,S,\gamma},
 C_{\nbiga,S,\gamma}.
\]
\end{prop}
\pf
Let us prove
$\pi_+^iL_{\star}(\nbiga,S,\gamma)=0$
$(i\neq 0)$.
By using the duality,
it is enough to prove that
$\pi_{+}^iL_{\star}(\nbiga,S,\gamma)=0$
for $i>0$ and for $\ast,!$.
By Lemma \ref{lem;14.6.28.1},
we have
$R^j\pi_{\ast}\bigl(
 L_{\star}(\nbiga,S,\gamma)
 \otimes\Omega^i_{\proj^m}
 \bigr)=0$
for any $j>0$ and any $i$.
By using the expression
\[
 \pi^j_{+}L_{\star}(\nbiga,S,\gamma)
\simeq
 \hyperr^j\pi_{\Sigma\ast}\bigl(
 \Omega^{\bullet+m}_{\proj^m\times S/S}
 \otimes L_{\star}(\nbiga,S,\gamma)
 \bigr),
\]
we obtain the claim.

Let us prove the vanishings for 
$K_{\nbiga,S,\gamma}$,
$C_{\nbiga,S,\gamma}$
and 
$L_{\min}(\nbiga,S,\gamma)$.
By Corollary \ref{cor;14.6.28.10},
we obtain 
$\pi^i_{+}K_{\nbiga,S,\gamma}=0$,
$\pi^i_{+}C_{\nbiga,S,\gamma}=0$
and
$\pi_{+}^iL_{\min}(\nbiga,S,\gamma)=0$
for $i>0$.
Let $-\gamma:S\lrarr H^0\bigl(\proj^m,\nbigo(1)\bigr)$
denote the composite of
$\gamma$ and the multiplication of $-1$
on $H^0\bigl(\proj^m,\nbigo(1)\bigr)$.
We naturally have
$\DDD K_{\nbiga,S,\gamma}\simeq
 C_{\nbiga,S,-\gamma}$,
$\DDD C_{\nbiga,S,\gamma}\simeq
 K_{\nbiga,S,-\gamma}$,
and 
$\DDD
 L_{\min}(\nbiga,S,\gamma)
\simeq
 L_{\min}(\nbiga,S,-\gamma)$.
Then, 
we obtain 
$\pi^i_{+}K_{\nbiga,S,\gamma}=0$,
$\pi^i_{+}C_{\nbiga,S,\gamma}=0$
and
$\pi_{+}^iL_{\min}(\nbiga,S,\gamma)=0$
for $i<0$.
\hfill\qed

\begin{cor}
\label{cor;14.7.4.21}
The following are exact sequences:
\[
 0\lrarr
 \pi_{+}^0K_{\nbiga,S,\gamma}
\lrarr
 \pi_+^0L_!(\nbiga,S,\gamma)
\lrarr
  \pi_+^0L_{\min}(\nbiga,S,\gamma)
\lrarr 0
\]
\[
0\lrarr  \pi_+^0L_{\min}(\nbiga,S,\gamma)
\lrarr
  \pi_+^0L_{\ast}(\nbiga,S,\gamma)
\lrarr
 \pi_{+}^0C_{\nbiga,S,\gamma}
\lrarr 0
\]
In particular,
$\pi_+^0L_{\min}(\nbiga,S,\gamma)$
is isomorphic to the image of
$\pi_+^0L_!(\nbiga,S,\gamma)
\lrarr
  \pi_+^0L_{\ast}(\nbiga,S,\gamma)$.
\hfill\qed
\end{cor}

\begin{cor}
For any desingularization
$\varphi_{\Sigma}:X_{\Sigma}\lrarr X_{\nbiga}$,
we have
$\pi_{\Sigma+}^i
 L_{\star}(F_{\gamma,\Sigma},D_{\Sigma,S})
 =0$ for $i\neq 0$.
\hfill\qed
\end{cor}

\begin{cor}
\label{cor;14.11.23.220}
Suppose that $\gamma$ satisfies the stronger condition
that $\gamma(S)\subset U_{\nbiga}^{\reg}$.
Then, the $\nbigd$-modules
$\pi_+^0L_{\star}(\nbiga,S,\gamma)$,
$\pi_+^0L_{\min}(\nbiga,S,\gamma)$,
$\pi_+^0K_{\nbiga,S,\gamma}$
and 
$\pi_+^0C_{\nbiga,S,\gamma}$
are locally free $\nbigo_S$-modules
with a flat connection.
\end{cor}
\pf
We obtain the claims for 
$\pi_+^0L_{\star}(\nbiga,S,\gamma)$
from Corollary \ref{cor;14.7.10.10}
and the isomorphisms (\ref{eq;14.11.6.1}).
Then, the claims for the others follows.
\hfill\qed

\vspace{.1in}

Suppose moreover that $0$ is an interior point of $\Conv(\nbiga)$.
In this case,
we have
$L_{!}(\nbiga,S,\gamma)
\simeq
 L_{\min}(\nbiga,S,\gamma)
\simeq 
 L_{\ast}(\nbiga,S,\gamma)$,
which is denoted by $L(\nbiga,S,\gamma)$.
We have
$\pi_+^iL(\nbiga,S,\gamma)=0$ for $i\neq 0$.
For a desingularization
$\varphi_{\Sigma}:X_{\Sigma}\lrarr X_{\nbiga}$,
we have the $\nbigd$-module
$L(F_{\Sigma,\gamma})
=L_{!}(F_{\gamma,\Sigma},D_{\Sigma,S})
=L_{\ast}(F_{\gamma,\Sigma},D_{\Sigma,S})$.
We have the following corollary
as a special case.

\begin{cor}
\label{cor;14.7.9.100}
If $0$ is an interior point of $\Conv(\nbiga)$,
we have
$\pi^i_{\Sigma+}L(F_{\gamma,\Sigma})
=0$ for $i\neq 0$.
We naturally have
$\pi^0_{\Sigma+}L(F_{\gamma,\Sigma})
\simeq
 \pi^0_{+}L(\nbiga,S,\gamma)$.
\hfill\qed
\end{cor}

\begin{example}
The $\nbigd$-modules associated to 
the non-degenerate section in
Example {\rm\ref{example;14.11.24.100}}
are called the GKZ-hypergeometric systems.
The $\nbigd$-modules associated to 
the non-degenerate sections in
Example {\rm\ref{example;14.11.24.101}}
are called the reduced GKZ-hypergeometric systems.
Note that it is independent of the choice of a section,
according to Lemma {\rm \ref{lem;14.12.30.10}}.
\hfill\qed
\end{example}

\subsubsection{De Rham complexes}
\label{subsection;14.11.23.221}

Let $\gamma:S\lrarr H^0\bigl(\proj^m,\nbigo(1)\bigr)$
be a holomorphic map
which is non-degenerate at $\infty$ for $X_{\nbiga}$.
Take a toric desingularization
$\varphi_{\Sigma}:X_{\Sigma}\lrarr X_{\nbiga}$.

\begin{prop}
If $0$ is contained in the interior part of $\Conv(\nbiga)$,
we have the following natural isomorphism
for 
$L(F_{\gamma,\Sigma})=
L_{!}(F_{\gamma,\Sigma},D_{\Sigma,S})
=L_{\ast}(F_{\gamma,\Sigma},D_{\Sigma,S})$:
\begin{equation}
\label{eq;14.11.5.10}
 \pi^0_{\Sigma+}\bigl(
 L(F_{\gamma,\Sigma})
 \bigr)
\simeq
 \hyperr^n\pi_{\Sigma\ast}\bigl(
 \Omega^{\bullet}_{X_{\Sigma,S}/S}(\ast D_{\Sigma,S}),
 d+dF_{\gamma,\Sigma}
 \bigr)
\end{equation}
If $\gamma$ is algebraic,
it means
$\pi_{+}^0\bigl(
 L(\nbiga,S,\gamma)
 \bigr)
\simeq
 \hyperr^n\pi_{\ast}\bigl(
 \Omega^{\bullet}_{T^n\times S/S},
 d+dF_{\gamma}
 \bigr)$.

If $0$ is a boundary point of 
$\Conv(\nbiga\cup\{0\})$,
we have the following isomorphisms:
\begin{equation}
 \label{eq;14.11.6.12}
 \pi^0_{\Sigma+}
 L_{\ast}(F_{\gamma,\Sigma},D_{\Sigma,S})
\simeq
 \hyperr^n\pi_{\Sigma\ast}\Bigl(
 \Omega_{X_{\Sigma,S}/S}^{\bullet}
 (\log D_{\Sigma,S})(\ast (F_{\gamma,\Sigma})_{\infty}),
 d+dF_{\gamma,\Sigma}
 \Bigr)
\end{equation}
\begin{equation}
 \label{eq;14.11.6.13}
 \pi^0_{\Sigma+}
 L_!(F_{\gamma,\Sigma},D_{\Sigma,S})
\simeq
 \hyperr^n\pi_{\Sigma\ast}\Bigl(
 \Omega_{X_{\Sigma,S}/S}^{\bullet}
 (\log D_{\Sigma,S})
 (-D_{\Sigma,S})
 (\ast (F_{\gamma,\Sigma})_{\infty}),
 d+dF_{\gamma,\Sigma}
 \Bigr)
\end{equation}
\end{prop}
\pf
We immediately (\ref{eq;14.11.5.10}) from 
$L(F_{\gamma,\Sigma},D_{\Sigma,S})
\simeq
 \bigl(
 \nbigo_{X_{\Sigma,S}}(\ast D_{\Sigma,S}),
 d+dF_{\gamma,\Sigma}
 \bigr)$.
By applying the arguments in
\S\ref{subsection;14.11.6.3}
and 
\S\ref{subsection;14.11.6.2},
we obtain (\ref{eq;14.11.6.12}) and 
 (\ref{eq;14.11.6.13}).
\hfill\qed

\vspace{.1in}

Let $\check{\varphi}_{\nbiga}:
 \check{X}_{\nbiga}\lrarr X_{\nbiga}$ denote the normalization of
$X_{\nbiga}$,
which is the normal toric variety.
Let us rewrite (\ref{eq;14.11.6.12})
and (\ref{eq;14.11.6.13})
in terms of $\check{X}_{\nbiga}$.
We put
$\check{D}_{\nbiga}:=
 \check{X}_{\nbiga}\setminus T^n$.
We also set
$\check{D}_{\nbiga,\infty}:=
 \check{\varphi}_{\nbiga}^{-1}(X_{\nbiga}\cap H_{\infty})$.

Let $\Omega^i_{\check{X}_{\nbiga}}(\log \check{D}_{\nbiga})$
denote the sheaf of logarithmic $i$-forms 
$(\check{X}_{\nbiga},\check{D}_{\nbiga})$
as in \cite{Batyrev-VMHS}.
Let $\Omega^i_{\check{X}_{\nbiga},\check{D}_{\nbiga}}$
denote the sheaf of $i$-forms on $X_{\nbiga}$
whose restrictions to each stratum
of $\check{D}_{\nbiga}$ is $0$
as in \cite{Danilov-de-Rham}.
Let 
$\Omega_{\check{X}_{\nbiga}\times S/S}^{\bullet}
 (\log \check{D}_{\nbiga,S})$
and 
$\Omega^{\bullet}_{(\check{X}_{\nbiga},\check{D}_{\nbiga})\times S/S}$
denote the pull back of
$\Omega_{\check{X}_{\nbiga}}^{\bullet}
 (\log \check{D}_{\nbiga})$
and 
$\Omega^{\bullet}_{(\check{X}_{\nbiga},\check{D}_{\nbiga})}$
by the projection $\check{X}_{\nbiga}\times S\lrarr \check{X}_{\nbiga}$,
respectively.
\begin{prop}
We have the following isomorphisms:
\begin{equation}
\label{eq;14.11.6.10}
 \pi^0_{+}L_{\ast}(\nbiga,S,\gamma)
\simeq
 \hyperr^n\pi_{\ast}\Bigl(
 \Omega_{\check{X}_{\nbiga}\times S/S}^{\bullet}
 (\log \check{D}_{\nbiga,S})(\ast \check{D}_{\nbiga,\infty,S}),
 d+dF_{\gamma}
 \Bigr)
\end{equation}
\begin{equation}
\label{eq;14.11.6.11}
 \pi^0_{+}L_{!}(\nbiga,S,\gamma)
\simeq
 \hyperr^n\pi_{\ast}\Bigl(
 \Omega_{(\check{X}_{\nbiga},\check{D}_{\nbiga})\times S/S}^{\bullet}
 (\ast \check{D}_{\nbiga,\infty,S}),
 d+dF_{\gamma}
 \Bigr)
\end{equation}
\end{prop}
\pf
According to  \cite[Lemma 6.1]{Batyrev-VMHS},
we have
\begin{equation}
 \label{eq;14.7.4.30}
 R\varphi_{\Sigma\ast}\Omega^i_{X_{\Sigma}}(\log D_{\Sigma})
\simeq
 \Omega^i_{\check{X}_{\nbiga}}(\log \check{D}_{\nbiga}).
\end{equation}
For a smooth toric variety $X_{\Sigma}$,
we have
$\Omega^i_{X_{\Sigma},D_{\Sigma}}\simeq
 \Omega^i_{X_{\Sigma}}(\log D_{\Sigma})(-D_{\Sigma})$.
According to \cite[Proposition 1.8]{Danilov-de-Rham},
we have 
\begin{equation}
 \label{eq;14.7.4.31}
 R\varphi_{\Sigma\ast}\Omega^i_{X_{\Sigma}}(\log D_{\Sigma})(-D_{\Sigma})
\simeq
 \Omega^i_{\check{X}_{\nbiga},\check{D}_{\nbiga}}.
\end{equation}
Then, we obtain
(\ref{eq;14.11.6.10})
and 
(\ref{eq;14.11.6.11})
from 
(\ref{eq;14.11.6.12})
and 
(\ref{eq;14.11.6.13}).
\hfill\qed

\vspace{.1in}
In the algebraic setting,
we set 
$X^{\aff}_{\nbiga}:=X_{\nbiga}\setminus H_{\infty}$.
Let $\check{\varphi}^{\aff}:\check{X}^{\aff}\lrarr X^{\aff}$
denote the normalization.
Let $\check{D}^{\aff}_{\nbiga}:=
 \check{X}^{\aff}_{\nbiga}\setminus
 T^n$.
We obtain the following similarly.
\begin{prop}
We have the following isomorphisms:
\begin{equation}
 \label{eq;15.5.11.3}
 \pi^0_{+}L_{\ast}(\nbiga,S,\gamma)
\simeq
 \hyperr^n\pi_{\ast}\Bigl(
 \Omega_{\check{X}^{\aff}_{\nbiga,S}/S}^{\bullet}
 (\log \check{D}^{\aff}_{\nbiga,S}),
 d+dF_{\gamma}
 \Bigr)
\end{equation}
\begin{equation}
\label{eq;15.5.11.4}
 \pi^0_{+}L_{!}(\nbiga,S,\gamma)
\simeq
 \hyperr^n\pi_{\ast}\Bigl(
 \Omega_{(\check{X}^{\aff}_{\nbiga,S},
 \check{D}^{\aff}_{\nbiga,S})/S}^{\bullet},
 d+dF_{\gamma}
 \Bigr)
\end{equation}
\hfill\qed
\end{prop}

\begin{rem}
In the case of Example {\rm\ref{example;14.11.24.100}},
by using Proposition {\rm\ref{prop;15.4.25.2}},
we can describe 
{\rm(\ref{eq;15.5.11.3})}
and {\rm(\ref{eq;15.5.11.4})}
as the $\nbigd$-modules
associated to better behaved GKZ-systems
in {\rm\cite{Borisov-Horja}}.
\end{rem}

\subsection{Mixed twistor $\nbigd$-modules}
\label{subsection;14.11.21.100}

\subsubsection{Mixed twistor $\nbigd$-modules
associated to families of Laurent polynomials}

We use the notation in \S\ref{subsection;14.7.4.10}.
Let $\gamma:S\lrarr H^0(\proj^m,\nbigo(1))$
be a holomorphic map.
We take a toric desingularization
$\varphi_{\Sigma}:
 X_{\Sigma}
\lrarr
 X_{\nbiga}$.
Let $d_S:=\dim S$.
We have the integrable mixed twistor $\nbigd$-modules
$\nbigt_{\star}(F_{\gamma,\Sigma},D_{\Sigma,S})$
$(\star=\ast,!)$
with real structure 
on $X_{\Sigma}\times S$:
\[
 \nbigt_{\ast}(F_{\gamma,\Sigma},D_{\Sigma,S})
=\Bigl(
\lambda^{n+d_S}\nbigl_{!}(F_{\gamma,\Sigma},D_{\Sigma,S}),
\nbigl_{\ast}(F_{\gamma,\Sigma},D_{\Sigma,S}),
\nbigc_{\ast}(F_{\gamma,\Sigma},D_{\Sigma,S})
 \Bigr).
\]
\[
  \nbigt_{!}(F_{\gamma,\Sigma},D_{\Sigma,S})
=\Bigl(
\lambda^{n+d_S}
 \nbigl_{\ast}(F_{\gamma,\Sigma},D_{\Sigma,S}),
\nbigl_{!}(F_{\gamma,\Sigma},D_{\Sigma,S}),
\nbigc_{!}(F_{\gamma,\Sigma},D_{\Sigma,S})
 \Bigr).
\]
The weight filtration is denoted by $W$.
The image of
$\nbigt_!(\nbiga,S,\gamma)
\lrarr
 \nbigt_{\ast}(\nbiga,S,\gamma)$
is denoted by
$\nbigt_{\min}(\nbiga,S,\gamma)$.

On $\proj^m\times S$,
we obtain the $\nbigrtilde_{\proj^m\times S}$-modules
$\nbigl_{\star}(\nbiga,S,\gamma):=
 \varphi^0_{\Sigma,S+}
 \nbigl_{\star}(F_{\gamma,\Sigma},D_{\Sigma,S})$
and the integrable mixed twistor $\nbigd$-modules
with induced real structure 
$\nbigt_{\star}(\nbiga,S,\gamma):=
\varphibar^0_{\Sigma,S+}
 \nbigt_{\star}(F_{\gamma,\Sigma},D_{\Sigma,S})$:
\[
 \nbigt_{\ast}(\nbiga,S,\gamma)
=\Bigl(
 \lambda^{n+d_S}\nbigl_{!}(\nbiga,S,\gamma),
 \nbigl_{\ast}(\nbiga,S,\gamma),
 \nbigc_{\ast}(\nbiga,S,\gamma)
 \Bigr)
\]
\[
  \nbigt_{!}(\nbiga,S,\gamma)
=\Bigl(
 \lambda^{n+d_S}\nbigl_{\ast}(F_{\gamma},D_{\nbiga,S}),
 \nbigl_{!}(F_{\gamma},D_{\nbiga,S}),
 \nbigc_!(\nbiga,S,\gamma)
 \Bigr)
\]
Here, $\nbigc_{\star}(\nbiga,S,\gamma)$ $(\star=\ast,!)$
are obtained as the push-forward of
$\nbigc_{\star}(F_{\gamma,\Sigma},D_{\Sigma,S})$.
They are twistor enhancement
of $L_{\star}(\nbiga,S,\gamma)$.
As in the case of $\nbigd$-modules
(\S\ref{subsection;14.11.23.210}),
they are independent of the choice of
a toric desingularization.
We also obtain integrable mixed twistor $\nbigd$-modules
with real structure 
$\pi^i_{\dagger}\bigl(
 \nbigt_{\star}(\nbiga,S,\gamma)
 \bigr)$
$(\star=\ast,!)$
on $S$
which are naturally isomorphic to
$\pi^i_{\Sigma\dagger}\bigl(
 \nbigt_{\star}(F_{\gamma,\Sigma},D_{\Sigma,S})
 \bigr)$.

\begin{rem}
Let $X_{\Sigma_1}$ be any $n$-dimensional
smooth toric variety
with a fixed inclusion $T^n\subset X_{\Sigma_1}$.
We have the meromorphic function
$F_{\gamma,\Sigma_1}$
on $(X_{\Sigma_1,S},D_{\Sigma_1,S})$
associated to the family of Laurent polynomials
$F_{\gamma}$.
We have the associated integrable mixed twistor $\nbigd$-modules
$\nbigt_{\star}(F_{\gamma,\Sigma_1},D_{\Sigma_1,S})$.
As in the case of $\nbigd$-modules,
we naturally have
$\pi^i_{\Sigma_1\dagger}
 \nbigt_{\star}(F_{\gamma,\Sigma_1},D_{\Sigma_1,S})
\simeq
 \pi^i_{\dagger}
 \nbigt_{\star}(\nbiga,S,\gamma)$.
\hfill\qed
\end{rem}

\subsubsection{Non-degenerate family of sections}
\label{subsection;14.11.6.10}

Let $\gamma:S\lrarr 
 H^0\bigl(\proj^m,\nbigo_{\proj^m}(1)\bigr)$
be a holomorphic map
which is non-degenerate at $\infty$ for $X_{\nbiga}$.

\begin{prop}
We naturally have
$\nbigl_{\star}(\nbiga,S,\gamma)
 (\ast H_{\infty,S})
\simeq
 \nbigl_{\star}(\nbiga,S,\gamma)$.
\end{prop}
\pf
We obtain the claim by using the argument 
in Lemma \ref{lem;14.6.28.1}
together with
Corollary \ref{cor;14.10.18.2}
and Proposition \ref{prop;14.10.18.1}.
\hfill\qed

\vspace{.1in}

Let $\nbigk_{\nbiga,S,\gamma}$
and $\nbigc_{\nbiga,S,\gamma}$
denote the kernel and the cokernel of
$\nbigt_!(\nbiga,S,\gamma)
\lrarr
 \nbigt_{\ast}(\nbiga,S,\gamma)$.

\begin{cor}
Let $\nbigm$ be the $\nbigrtilde_{\proj^m\times S}$-modules
underlying 
$\nbigk_{\nbiga,S,\gamma}$, $\nbigc_{\nbiga,S,\gamma}$
or $\nbigt_{\min}(\nbiga,S,\gamma)$.
Then, we have $\nbigm(\ast H_{\infty,S})\simeq\nbigm$.
\hfill\qed
\end{cor}

We obtain the following 
from Proposition \ref{prop;14.7.4.20}.
\begin{prop}
\label{prop;14.11.8.110}
We have
$\pi^i_+\nbigt=0$ $(i\neq 0)$
for the integrable mixed twistor $\nbigd$-modules
for 
$\nbigt=\nbigt_{\star}(\nbiga,S,\gamma)$
$(\star=\ast,!)$
$\nbigt_{\min}(\nbiga,S,\gamma)$,
$\nbigk_{\nbiga,S,\gamma}$,
and 
$\nbigc_{\nbiga,S,\gamma}$.
\hfill\qed
\end{prop}

We obtain the following corollary.
\begin{cor}
\label{cor;14.11.29.400}
We have the following exact sequences
of integrable mixed twistor $\nbigd$-modules on $S$:
\[
 0\lrarr
 \pi_+^0\nbigk_{\nbiga,S,\gamma}
\lrarr
 \pi_+^0\nbigt_!(\nbiga,S,\gamma)
\lrarr
 \pi_+^0\nbigt_{\min}(\nbiga,S,\gamma)
\lrarr 0
\]
\[
 0\lrarr
 \pi_+^0\nbigt_{\min}(\nbiga,S,\gamma)
\lrarr
 \pi_+^0\nbigt_{\ast}(\nbiga,S,\gamma)
\lrarr
 \pi_+^0\nbigc_{\nbiga,S,\gamma}
\lrarr 0
\]
As a result,
the image of the morphism
$\pi_+^0\nbigt_!(\nbiga,S,\gamma)
\lrarr
 \pi_+^0\nbigt_{\ast}(\nbiga,S,\gamma)$
is naturally isomorphic to 
the integrable pure twistor $\nbigd$-module
$\pi_+^0\nbigt_{\min}(\nbiga,S,\gamma)$.
Moreover, we have
\[
 \Gr^W_{n+d_S}
 \pi_{\dagger}^0\nbigt_{!}(\nbiga,S,\gamma)
\simeq
  \pi_{\dagger}^0\nbigt_{\min}(\nbiga,S,\gamma)
\simeq
\Gr^W_{n+d_S}
 \pi_{\dagger}^0\nbigt_{\ast}(\nbiga,S,\gamma).
\]
\hfill\qed
\end{cor}

\begin{cor}
For any desingularization
$\varphi_{\Sigma}:X_{\Sigma}\lrarr X_{\nbiga}$,
we have
$\pi^i_{\Sigma\dagger}\bigl(
 \nbigt_{\star}(F_{\gamma,\Sigma},D_{\Sigma,S})
 \bigr)=0$ $(i\neq 0)$.
\hfill\qed
\end{cor}

\subsubsection{Descriptions in terms of the de Rham complexes}
\label{subsection;14.12.15.100}

Let 
$\gamma:S\lrarr H^0\bigl(\proj^m,\nbigo_{\proj^m}(1)\bigr)$
be a morphism which is non-degenerate at $\infty$
for $X_{\nbiga}$.
We give descriptions of
the $\nbigrtilde$-modules
$\pi_{\dagger}^0\nbigl_{\star}(\nbiga,S,\gamma)$.
We use the notation in 
\S\ref{subsection;14.11.6.11}
and \S\ref{subsection;14.11.23.221}.

\begin{prop}
\label{prop;14.11.6.20}
Take a toric desingularization
$\gamma:X_{\Sigma}\lrarr X_{\nbiga}$.
If $0$ is an interior point of $\Conv(\nbiga)$,
we have
\[
 \pi^0_{\Sigma\dagger}\bigl(
 \nbigl(F_{\gamma,\Sigma})
 \bigr)
\simeq
 \hyperr^n\pi_{\Sigma\ast}\bigl(
 \Omegabar^{\bullet}_{X_{\Sigma,S}/S}(\ast D_{\Sigma,S}),
 d+\lambda^{-1}dF_{\gamma,\Sigma}
 \bigr).
\]
If moreover $\gamma$ is algebraic,
it means
$\pi_{\dagger}^0\bigl(
 \nbigl(F_{\gamma})
 \bigr)
\simeq
 \hyperr^n\pi_{\ast}\bigl(
 \Omegabar^{\bullet}_{T^n\times S/S},
 d+\lambda^{-1}dF_{\gamma}
 \bigr)$.

If $0$ is a boundary point of
$\Conv(\nbiga\cup\{0\})$,
we have the following natural isomorphisms:
\begin{equation}
 \label{eq;14.11.24.40}
 \pi^0_{\Sigma\dagger}
 \nbigl_{\ast}(F_{\gamma,\Sigma},D_{\Sigma,S})
\simeq
 \hyperr^n\pi_{\Sigma\ast}\Bigl(
 \Omegabar_{X_{\Sigma,S}/S}^{\bullet}
 (\log D_{\Sigma,S})(\ast (F_{\gamma,\Sigma})_{\infty}),
 d+\lambda^{-1}dF_{\gamma,\Sigma}
 \Bigr)
\end{equation}
\begin{equation}
 \label{eq;14.11.24.41}
 \pi^0_{\Sigma\dagger}
 \nbigl_!(F_{\gamma,\Sigma},D_{\Sigma,S})
\simeq
 \hyperr^n\pi_{\Sigma\ast}\Bigl(
 \Omegabar_{X_{\Sigma,S}/S}^{\bullet}
 (\log D_{\Sigma,S})
 (-D_{\Sigma,S})
 (\ast (F_{\gamma,\Sigma})_{\infty}),
 d+\lambda^{-1}dF_{\gamma,\Sigma}
 \Bigr)
\end{equation}
\end{prop}
\pf
They are the special cases
of the isomorphisms given in \S\ref{subsection;14.11.6.11}.
\hfill\qed

\vspace{.1in}

Let $q:\cnum_{\lambda}\times \check{X}_{\nbiga}\times S
  \lrarr \check{X}_{\nbiga}$ be the projection.
Let $\Omegabar^i_{\check{X}_{\nbiga,S}/S}(\log \check{D}_{\nbiga,S})$
denote 
$\lambda^{-i}q^{\ast}
\Omega^i_{\check{X}_{\nbiga}}(\log \check{D}_{\nbiga})$.
Let $\Omegabar^i_{\check{X}_{\nbiga},\check{D}_{\nbiga}}$
denote
$\lambda^{-i}q^{\ast}
 \Omega^i_{\check{X}_{\nbiga},\check{D}_{\nbiga}}$.
\begin{prop}
\label{prop;14.12.3.20}
We have the following natural isomorphisms:
\begin{equation}
\label{eq;14.11.24.42}
 \pi^0_{\dagger}\nbigl_{\ast}(\nbiga,S,\gamma)
\simeq
 \hyperr^n\pi_{\ast}\Bigl(
 \Omegabar_{\check{X}_{\nbiga}\times S/S}^{\bullet}
 (\log \check{D}_{\nbiga,S})(\ast \check{D}_{\nbiga,\infty,S}),
 d+\lambda^{-1}dF_{\gamma}
 \Bigr)
\end{equation}
\begin{equation}
\label{eq;14.11.24.43}
 \pi^0_{\dagger}\nbigl_{!}(\nbiga,S,\gamma)
\simeq
 \hyperr^n\pi_{\ast}\Bigl(
 \Omegabar_{(\check{X}_{\nbiga},\check{D}_{\nbiga})\times S/S}^{\bullet}
 (\ast \check{D}_{\nbiga,\infty,S}),
 d+\lambda^{-1}dF_{\gamma}
 \Bigr)
\end{equation}
In the algebraic setting,
we have the following isomorphisms:
\begin{equation}
 \label{eq;14.11.24.44}
 \pi^0_{\dagger}\nbigl_{\ast}(\nbiga,S,\gamma)
\simeq
 \hyperr^n\pi_{\ast}\Bigl(
 \Omegabar_{\check{X}^{\aff}_{\nbiga}\times S/S}^{\bullet}
 (\log \check{D}^{\aff}_{\nbiga,S}),
 d+\lambda^{-1}dF_{\gamma}
 \Bigr)
\end{equation}
\begin{equation}
 \label{eq;14.11.24.45}
 \pi^0_{\dagger}\nbigl_{!}(\nbiga,S,\gamma)
\simeq
 \hyperr^n\pi_{\ast}\Bigl(
 \Omegabar_{(\check{X}^{\aff}_{\nbiga},
 \check{D}^{\aff}_{\nbiga})\times S/S}^{\bullet},
 d+\lambda^{-1}dF_{\gamma}
 \Bigr)
\end{equation}
\end{prop}
\pf
We obtain the isomorphisms
from Proposition \ref{prop;14.11.6.20}
with the isomorphisms (\ref{eq;14.7.4.30}) 
and (\ref{eq;14.7.4.31}).
\hfill\qed

\begin{cor}
The image of 
the natural morphism
\begin{equation}
\label{eq;14.11.24.46}
 \hyperr^n\pi_{\ast}\Bigl(
 \Omegabar_{(\check{X}^{\aff}_{\nbiga},
 \check{D}^{\aff}_{\nbiga})\times S/S}^{\bullet},
 d+\lambda^{-1}dF_{\gamma}
 \Bigr)
\lrarr
 \hyperr^n\pi_{\ast}\Bigl(
 \Omegabar_{\check{X}^{\aff}_{\nbiga}\times S/S}^{\bullet}
 (\log \check{D}^{\aff}_{\nbiga,S}),
 d+\lambda^{-1}dF_{\gamma}
 \Bigr)
\end{equation}
underlies a pure twistor $\nbigd$-module.
\hfill\qed
\end{cor}

Let $U_{\nbiga}^{\reg}\subset H^0(\proj^m,\nbigo(1))$
be the open subset as in \S\ref{subsection;14.11.23.240}.
Suppose that the image of $\gamma$ is contained in
$U_{\nbiga}^{\reg}$.
In this case,
$\pi^0_{\dagger}\bigl(
 \nbigl_{\star}(\nbiga,S,\gamma)
 \bigr)$
are locally free $\nbigo_{\cnum_{\lambda}\times S}$-modules.
Let 
$\pi^0_{\dagger}\bigl(
 \nbigl_{\star}(\nbiga,S,\gamma)
 \bigr)^{\lor}$
be the dual as $\nbigo_{\cnum_{\lambda}\times S}$-modules,
which are naturally equipped with the meromorphic connection.
They are naturally isomorphic to
$\lambda^{-d_S}\DDD\pi^0_{\dagger}\bigl(
 \nbigl_{\star}(\nbiga,S,\gamma)
 \bigr)$
up to signatures.
The real structure of 
$\pi_{\dagger}\nbigt_{!}(\nbiga,S,\gamma)$
gives
\[
 \pi_{\dagger}^0
 \nbigl_{!}(\nbiga,S,\gamma)
\simeq
 j^{\ast}\DDD\pi^0_{\dagger}
 \bigl(
  \lambda^{n+d_S}
 \nbigl_{\ast}(\nbiga,S,\gamma)
 \bigr)
\simeq
\lambda^{-n-d_S}
  j^{\ast}\DDD\pi^0_{\dagger}
 \nbigl_{\ast}(\nbiga,S,\gamma)
\simeq
 \lambda^{-n}
 \Bigl(
 j^{\ast}\pi_{\dagger}^0
 \nbigl_{\ast}(\nbiga,S,\gamma)
 \Bigr)^{\lor}
\]

\begin{cor}
$\pi_{\dagger}^0
 \nbigl_{!}(\nbiga,S,\gamma)$
is naturally identified with
$\lambda^{-n}j^{\ast}\hyperr^n\pi_{\ast}\Bigl(
 \Omegabar_{\check{X}^{\aff}_{\nbiga}\times S/S}^{\bullet}
 (\log \check{D}^{\aff}_{\nbiga,S}),
 d+\lambda^{-1}dF_{\gamma}
 \Bigr)^{\lor}$.
The image of the natural morphism
\[
 \lambda^{-n}
 j^{\ast}\hyperr^n\pi_{\ast}\Bigl(
 \Omegabar_{\check{X}^{\aff}_{\nbiga}\times S/S}^{\bullet}
 (\log \check{D}^{\aff}_{\nbiga,S}),
 d+\lambda^{-1}dF_{\gamma}
 \Bigr)^{\lor}
\lrarr
  \hyperr^n\pi_{\ast}\Bigl(
 \Omegabar_{\check{X}^{\aff}_{\nbiga}\times S/S}^{\bullet}
 (\log \check{D}^{\aff}_{\nbiga,S}),
 d+\lambda^{-1}dF_{\gamma}
 \Bigr)
\]
underlies an integrable pure twistor $\nbigd$-module.
\hfill\qed
\end{cor}

\begin{rem}
If $S$ and $\gamma$ are as in Example 
{\rm\ref{example;14.11.24.100}},
we can describe {\rm(\ref{eq;14.11.24.44})}
and {\rm(\ref{eq;14.11.24.45})}
in terms of systems of differential equations
as in {\rm\S\ref{subsection;15.5.2.1}}.
\hfill\qed
\end{rem}

\subsubsection{Graded polarizations}
\label{subsection;14.11.24.121}

Let $\nbiga=\{\veca_1,\ldots,\veca_m\}$
and $\veca_0:=(0,\ldots,0)$.
Let $\gamma:S\lrarr H^0\bigl(\proj^m,\nbigo(1)\bigr)$
be a morphism which is non-degenerate at $\infty$
for $X_{\nbiga}$.
We obtain the family of Laurent polynomials
\[
F_{\gamma}(x,t)
=\sum_{i=0}^m\gamma_i(x)t^{\veca_i}
\]
Take $\vecb\in\seisuu^{n}$ such that $0$ is contained in 
the interior part of $\Conv(\nbiga\cup\{\vecb\})$.
It gives a monomial $t^{\vecb}$.
We shall observe that
the mixed twistor $\nbigd$-modules
$\pi_{\dagger}^0\nbigt(\nbiga,S,\gamma)$
are equipped with graded sesqui-linear dualities
depending on the choice of $\vecb$.

Let $\varphi_{\Sigma}:X_{\Sigma}\lrarr X_{\nbiga}$
be a toric desingularization.
We have the meromorphic function
$t^{\vecb}_{\Sigma}$ on $(X_{\Sigma},D_{\Sigma})$
induced by $t^{\vecb}$.
We also have the meromorphic functions
$F_{\gamma,\Sigma}$ 
and $t_{\Sigma,S}^{\vecb}$ on $(X_{\Sigma,S},D_{\Sigma,S})$
induced by $F_{\gamma}$
and $t^{\vecb}$.

\begin{lem}
\label{lem;14.11.24.30}
Suppose that 
$|(t_{\Sigma}^{\vecb})_0|
 \cap
 |(t_{\Sigma}^{\vecb})_{\infty}|
=\emptyset$.
Then, the other conditions in Lemma {\rm\ref{lem;14.10.14.12}}
are satisfied for $g=F_{\gamma,\Sigma}$ and 
$f=t^{\vecb}_{\Sigma,S}$.
\end{lem}
\pf
By the non-degeneracy assumption on $\gamma$,
$F_{\gamma,\Sigma}$ is non-degenerate along $D_{\Sigma,S}$.
We clearly have
$|(t_{\Sigma,S}^{\vecb})_0|\cup
 |(t_{\Sigma,S}^{\vecb})_{\infty}|
 \subset D_{\Sigma,S}$.
It is enough to prove that 
$D_{\Sigma,S}=\bigl|(t_{\Sigma,S}^{\vecb})_{\infty}\bigr|
 \cup\bigl|(F_{\gamma,\Sigma})_{\infty}\bigr|$.

We set $\veca_{m+1}:=\vecb$
and $\nbigatilde:=\nbiga\cup\{\veca_{m+1}\}$.
We set $\Stilde:=S\times\cnum_{\tau}$.
We can naturally regard 
$H^0\bigl(\proj^{m+1},\nbigo(1)\bigr)
=H^0\bigl(\proj^m,\nbigo(1)\bigr)
 \times
 \bigl(
 \cnum\cdot z_{m+1}
 \bigr)$.
We consider the morphism
$\gammatilde:=\gamma\times\id:
\Stilde\lrarr 
 H^0\bigl(\proj^{m+1},\nbigo(1)\bigr)$.
It is enough to prove that
$F_{\gammatilde}(x,\tau,t)=
  F_{\gamma}(x,t)+\tau t_{\Sigma,S}^{\vecb}$
is non-degenerate at $\infty$
for $X_{\nbigatilde}$.
Let $\sigma$ be any face of
$\Conv(\nbigatilde)$.
If $\vecb\not\in\sigma$,
then $\sigma$ is a face of 
$\Conv(\nbiga\cup\{0\})$
such that $0\not\in\sigma$.
Hence, we have
$(dF_{\gammatilde,\sigma})^{-1}(0)
=(dF_{\gamma,\sigma})^{-1}(0)
=\emptyset$
by applying Lemma \ref{lem;14.7.1.10} to $\nbiga$.
If $\vecb\in\sigma$,
then 
$\del_{\tau}F_{\gammatilde,\sigma}
=t_{\Sigma,S}^{\vecb}$
is nowhere vanishing.
Thus, we obtain the claim of the lemma.
\hfill\qed

\vspace{.1in}

We have the canonical polarization
$\bigl(
 (-1)^{n+d_S},(-1)^{n+d_S}
 \bigr)$
on the pure twistor $\nbigd$-module
$\nbigt(F_{\gamma,\Sigma},D_{\Sigma,S})$
of weight $n+d_S$ on $X_{\Sigma,S}$,
which we denote by
$\nbigs_{\gamma,\Sigma}$.
By Lemma \ref{lem;14.11.24.30},
we have
\[
\nbigt_{\star}(F_{\gamma,\Sigma},D_{\Sigma,S})
=\nbigt(F_{\gamma,\Sigma},D_{\Sigma,S})
 [\star (t_{\Sigma,S}^{\vecb})_{\infty}].
\]
As explained in \S\ref{subsection;14.11.24.1},
we have the induced graded polarization
$\nbigs_{\gamma,\Sigma}[\star(t^{\vecb}_{\Sigma,S})_{\infty}]$
on 
$\nbigt_{\star}(F_{\gamma,\Sigma},D_{\Sigma,S})$,
which we denote by
$\nbigs_{\gamma,\Sigma,\star,\vecb}$.
By Corollary \ref{cor;14.11.24.20},
Condition {\bf(A)} holds
for the projection $\pi_{\Sigma}$
and the mixed twistor $\nbigd$-modules
$\nbigt_{\star}(F_{\gamma,\Sigma},D_{\Sigma,S})$
with 
$\nbigs_{\gamma,\Sigma,\star,\vecb}$,
and we obtain the graded polarization
$\bigl[
 \pi_{\Sigma\dagger}^0
 \nbigs_{\gamma,\Sigma,\star,\vecb}
 \bigr]$
on 
$\pi_{\Sigma\dagger}^0
 \nbigt_{\star}(F_{\gamma,\Sigma},D_{\Sigma,S})
\simeq 
 \pi_{\dagger}^0\nbigt_{\star}(\nbiga,S,\gamma)$.
We obtain the following lemma from
Corollary \ref{cor;14.11.24.31}.
\begin{lem}
The graded polarizations
$\bigl[
 \pi_{\Sigma\dagger}^0
 \nbigs_{\gamma,\Sigma,\star,\vecb}
 \bigr]$
on 
$\pi_{\dagger}^0\nbigt_{\star}(\nbiga,S,\gamma)$
are independent of the choice of
a toric desingularization $\varphi_{\Sigma}$.
We denote them by
$\vecnbigs_{\star}(\nbiga,S,\gamma,\vecb)$.
\hfill\qed
\end{lem}

If $0$ is contained in the interior part of $\Conv(\nbiga)$,
then $\pi_{\dagger}^0\nbigt_{\star}(\nbiga,S,\gamma)$
is pure,
and $\vecnbigs_{\star}(\nbiga,S,\gamma,\vecb)$
is equal to
$\pi_{\dagger}^0\nbigs_{\gamma,\Sigma}$
for any toric desingularization 
$\varphi_{\Sigma}:X_{\Sigma}\lrarr X_{\nbiga}$.

\subsubsection{The smooth part and the induced graded pairings}
\label{subsection;14.11.26.131}

Recall that we have the open subset $U_{\nbiga}^{\reg}$
in $H^0\bigl(\proj^m,\nbigo_{\proj^m}(1)\bigr)$.
(See {\rm\S\ref{subsection;14.11.23.240}}.)
For any $\gamma:S\lrarr H^0\bigl(\proj^m,\nbigo_{\proj^m}(1)\bigr)$,
we set $S^{\reg}:=\gamma^{-1}(U_{\nbiga}^{\reg})$.

\begin{prop}
\label{prop;14.11.25.31}
Let $\nbigm$ be the $\nbigr_{S}$-modules
underlying $\pi_{\dagger}^0\nbigt$
for $\nbigt=
\nbigt_{\star}(\nbiga,S,\gamma)$ $(\star=\ast,!,\min)$,
$\nbigk_{\nbiga,S,\gamma}$
and $\nbigc_{\nbiga,S,\gamma}$.
Then, $\nbigm_{|S^{\reg}}$ is a locally free
$\nbigo_{\cnum_{\lambda}\times S^{\reg}}$-module.
In particular,
the mixed twistor $\nbigd$-module
$\pi_{\dagger}^0\nbigt_{|S^{\reg}}$ comes from 
a graded polarizable variation of mixed twistor structure
on $S^{\reg}$.
It is admissible along $S\setminus S^{\reg}$.
\end{prop}
\pf
It follows from Corollary \ref{cor;14.11.23.220}
and the general property of mixed twistor $\nbigd$-modules.
\hfill\qed

\begin{cor}
The restriction of 
{\rm(\ref{eq;14.11.24.40})}
and {\rm(\ref{eq;14.11.24.41})}
to $S^{\reg}$ are locally free
$\nbigo_{\cnum_{\lambda}\times S^{\reg}}$-modules.
Equivalently,
the restriction of
{\rm(\ref{eq;14.11.24.42})},
{\rm(\ref{eq;14.11.24.43})},
{\rm(\ref{eq;14.11.24.44})}
and {\rm(\ref{eq;14.11.24.45})}
to $S^{\reg}$
are locally free
$\nbigo_{\cnum_{\lambda}\times S^{\reg}}$-modules.
The restriction of 
the image of {\rm(\ref{eq;14.11.24.46})}
to $S^{\reg}$
is also an $\nbigo_{\cnum_{\lambda}\times S^{\reg}}$-module.
\hfill\qed
\end{cor}

Let $\gamma:S\lrarr H^0\bigl(\proj^m,\nbigo_{\proj^m}(1)\bigr)$
be a morphism which is non-degenerate at $\infty$
for $X_{\nbiga}$.
Suppose that 
(i) $S^{\reg}\neq\emptyset$,
(ii) we are given a hypersurface $Y$
which contains $S\setminus S^{\reg}$.
Take $\vecb\in\seisuu^n$ such that
$0$ is an interior point of $\Conv(\nbiga\cup\{\vecb\})$.
We have the mixed twistor $\nbigd$-modules
$\pi_{\dagger}^0\nbigt_{\star}(\nbiga,S,\gamma)$
with the induced graded polarization
$\vecnbigs_{\star}(\nbiga,S,\gamma,\vecb)$
on $S$.
They are also equipped with the natural real structure.
We set 
\[
 \nbigv_{\star}(\nbiga,S,\gamma):=
 \pi_{\dagger}^0\nbigl_{\star}(\nbiga,S,\gamma)(\ast Y).
\]
We obtain the filtration $\Wtilde$
and the graded pairing 
$P_{\star}(\nbiga,S,\gamma,\vecb)=\bigl(
 P_{\star}(\nbiga,S,\gamma,\vecb)_k\,\big|\,
 k\in\seisuu
 \bigr)$
by the procedure in \S\ref{subsection;14.11.28.10}.
We explain it in this situation.
The $\nbigr$-modules
$\pi_{\dagger}^0\nbigl_{\star}(\nbiga,S,\gamma)$
are equipped with the filtration $W$
underlying the weight filtration of 
$\pi_{\dagger}^0\nbigt_{\star}(\nbiga,S,\gamma)$.
We set 
\[
 \Wtilde_{k}\nbigv_{\star}(\nbiga,S,\gamma):=
\Bigl(
 W_{k+d_S}
 \pi_{\dagger}^0\nbigl_{\star}(\nbiga,S,\gamma)
\Bigr)(\ast Y)
\]
The polarization and the real structure
of $\Gr^W_{k+d_S}\pi_{\dagger}^0\nbigt_{\star}(\nbiga,S,\gamma)$
induce a pairing
$P_{\star}(\nbiga,S,\gamma)_k$
of weight $k$
on $\Gr^{\Wtilde}_k\nbigv_{\star}(\nbiga,S,\gamma)$.
Thus, we obtain a graded pairing
$P_{\star}(\nbiga,S,\gamma)=
 \bigl(P_{\star}(\nbiga,S,\gamma)_k\,\big|\,k\in\seisuu\bigr)$
on 
$\bigl(
 \nbigv_{\star}(\nbiga,S,\gamma),\Wtilde
 \bigr)$.
In this way,
we obtain mixed TEP-structures
$\bigl(
 \nbigv_{\star}(\nbiga,S,\gamma),
 \Wtilde,
 P_{\star}(\nbiga,S,\gamma)
 \bigr)$.

\subsubsection{Algebraicity}
\label{subsection;14.11.24.120}

Let $Z$ be any smooth complex quasi-projective variety.
An integrable mixed twistor $\nbigd$-module 
$(\nbigt,W)$ on $Z$
is called algebraic if the following holds:
\begin{itemize}
\item
 Let $\Zbar$ be a smooth projective manifold
 with an open immersion $Z\lrarr \Zbar$.
 Then, there exists an integrable mixed twistor $\nbigd$-module
 $\nbigt'$ on $\Zbar$
 such that 
 $\nbigt'_{|Z}=\nbigt$.
\end{itemize}

We also have the notion of the algebraicity
of the underlying $\nbigrtilde$-modules
of integrable mixed twistor $\nbigd$-modules.
Let $\nbigd^a_{\cnum_{\lambda}\times Z}$
denote the sheaf of algebraic linear differential operators
on $\cnum_{\lambda}\times Z$.
Let $\Theta_Z^a$ denote the algebraic tangent sheaf of $Z$.
Let $p:\cnum_{\lambda}\times Z\lrarr Z$ denote the projection.
Let $\nbigrtilde^{a}_{Z}\subset 
 \nbigd^a_{\cnum_{\lambda}\times Z}$
denote the sheaf of subalgebras
generated by  $\lambda\cdot p^{\ast}\Theta_Z$
and $\lambda^2\del_{\lambda}$
over $\nbigo_{\cnum_{\lambda}\times Z}$.
We say that an $\nbigrtilde_Z$-module $\nbigm$
is algebraic
if there exists an $\nbigrtilde^a_{Z}$-module
$\nbigm^a$ such that
$\nbigm$ is isomorphic to the analytification of $\nbigm^a$.

\vspace{.1in}

We naturally regard $H^0\bigl(\proj^m,\nbigo(1)\bigr)$
as a quasi-projective variety.
Let $\gamma:S\lrarr H^0\bigl(\proj^m,\nbigo(1)\bigr)$
be an algebraic morphism of 
smooth quasi-projective varieties.

\begin{prop}
\mbox{{}}\label{prop;14.11.25.40}
\begin{itemize}
\item
The integrable mixed twistor $\nbigd$-modules
$\nbigt_{\star}(\nbiga,S,\gamma)$
and 
$\pi_{\dagger}^i\nbigt_{\star}(\nbiga,S,\gamma)$
are algebraic.
The underlying $\nbigrtilde$-modules of 
$\nbigt_{\star}(\nbiga,S,\gamma)$
and 
$\pi_{\dagger}^i\nbigt_{\star}(\nbiga,S,\gamma)$
are algebraic.
\item
For any toric desingularization
$\varphi_{\Sigma}:X_{\Sigma}\lrarr X_{\nbiga}$,
the integrable mixed twistor $\nbigd$-modules
$\nbigt_{\star}(F_{\gamma,S},D_{\Sigma,S})$
are algebraic.
The underlying $\nbigrtilde$-modules of 
$\nbigt_{\star}(F_{\gamma,S},D_{\Sigma,S})$
are algebraic.
\end{itemize}
\end{prop}
\pf 
The claims for 
$\nbigt_{\star}(F_{\gamma,S},D_{\Sigma,S})$
is clear.
Then, we obtain the claims for the others
by construction.
\hfill\qed

\subsection{Description as a specialization}
\label{subsection;14.11.21.21}

\subsubsection{$\nbigd$-modules}
\label{subsection;14.7.4.40}

Let $\nbiga=\{\veca_1,\ldots,\veca_m\}\subset\seisuu^{n}$
be a finite subset which generates $\seisuu^n$.
We take $\veca_{m+1}\in\seisuu^{n}$ 
such that $0$ is contained in $\Conv(\nbigatilde)$,
where $\nbigatilde:=\nbiga\cup\{\veca_{m+1}\}$.
Let us compare the $\nbigd$-modules
associated to $\nbiga$ and $\nbigatilde$.

We identify $\proj^{m}$
with the subspace of $\proj^{m+1}$
determined by $z_{m+1}=0$.
We naturally regard
$H^0\bigl(
 \proj^{m+1},\nbigo(1)
 \bigr)
=H^0\bigl(
 \proj^{m},\nbigo(1)
 \bigr)
\times
 \cnum\cdot z_{m+1}$.

Let $\gamma:S\lrarr H^0\bigl(\proj^m,\nbigo(1)\bigr)$
be a morphism which is non-degenerate at $\infty$
for $X_{\nbiga}$.
We set $\Stilde:=S\times\cnum_{\tau}$.
We have the induced map
$\gammatilde=\gamma\times\id:
 \Stilde\lrarr
 H^0\bigl(
 \proj^{m+1},\nbigo(1)
 \bigr)$
given by
$\gammatilde(x,\tau)
=\gamma(x)+\tau z_{m+1}$.
Let $\pi:\proj^m\times S\lrarr S$ 
and $\pitilde:\proj^{m+1}\times \Stilde\lrarr \Stilde$
denote the projections.
We have the following proposition
on the associated $\nbigd$-modules.

\begin{prop}
\mbox{{}}
\label{prop;14.10.18.10}
\begin{itemize}
\item
We have
$L_{\ast}(\nbigatilde,\Stilde,\gammatilde)
\simeq
 L_{!}(\nbigatilde,\Stilde,\gammatilde)
\simeq
 L_{\min}(\nbigatilde,\Stilde,\gammatilde)$.
We also have $\pitilde^i_+L_{\star}(\nbigatilde,\Stilde,\gammatilde)=0$
for $i\neq 0$.
\item
Let $\iota:S\lrarr \Stilde$ be the inclusion
induced by $\{0\}\lrarr \cnum$.
The kernel and the cokernel of the morphism
\[
 \pitilde^0_+L_{\min}(\nbigatilde,\Stilde,\gammatilde)(!\tau)
\lrarr
 \pitilde^0_+L_{\min}(\nbigatilde,\Stilde,\gammatilde)(\ast\tau)
\]
are naturally isomorphic to
$\iota_{+}
 \pi_{+}^0L_{!}\bigl(
\nbiga,S,\gamma
 \bigr)$
and
$\iota_{+}
 \pi_{+}^0L_{\ast}\bigl(
\nbiga,S,\gamma
 \bigr)$,
respectively.
\end{itemize}
\end{prop}
\pf
Let $\varphi_{\Sigma}:X_{\Sigma}\lrarr X_{\nbigatilde}$
be a toric desingularization.
We may assume to have a toric morphism
$\varphi_{\Sigma,\nbiga}:X_{\Sigma}\lrarr X_{\nbiga}$
which is a desingularization.
For the meromorphic function $f=t^{\veca_{m+1}}$ on $X_{\Sigma}$,
we may assume $|(f)_0|\cap|(f)_{\infty}|=\emptyset$.

We have 
$F_{\Stilde,\gammatilde}
=F_{S,\gamma}+\tau t^{\veca_{m+1}}$.
We can check that the assumption in Lemma \ref{lem;14.10.14.12}
is satisfied for $F_{S,\gamma}$ and $t^{\veca_{m+1}}$.
(See \S\ref{subsection;14.11.24.121}.)
Hence, we have the purity of 
$F_{\Stilde,\gammatilde}$
on $(X_{\Sigma,\Stilde},D_{\Sigma,\Stilde})$.
Let $\pitilde_{\Sigma}:X_{\Sigma,\Stilde}\lrarr \Stilde$
be the projection.
By using Proposition \ref{prop;14.11.17.20},
we obtain 
$\pitilde^i_{\Sigma+}
 L(F_{\Stilde,\gammatilde},D_{\Sigma,\Stilde})=0$
for $i\neq 0$.
Thus, we obtain the first claim
of Proposition \ref{prop;14.10.18.10}.

Let $\iota_{1}:X_{\Sigma,S}\lrarr X_{\Sigma,\Stilde}$
denote the inclusion
induced by $\{0\}\lrarr \cnum$.
According to Proposition \ref{prop;14.6.28.10},
the kernel and the cokernel of 
the morphism
$L(F_{\Stilde,\gammatilde},D_{\Sigma,\Stilde})(!\tau)
\lrarr
 L(F_{\Stilde,\gammatilde},D_{\Sigma,\Stilde})(\ast\tau)$
are naturally isomorphic to
$\iota_{1+}
 L_{!}\bigl(
 F_{S,\gamma},D_{\Sigma,S}\bigr)$
and
$\iota_{1+}
 L_{\ast}\bigl(
 F_{S,\gamma},D_{\Sigma,S}\bigr)$,
respectively.

Let $\pi_{\Sigma}:X_{\Sigma}\times S\lrarr S$
be the projection.
By Proposition \ref{prop;14.7.4.20},
we have
$\pi_{\Sigma}^i
 L_{\star}(F_{S,\gamma},D_{\Sigma,S})=0$ for $i\neq 0$.
We also have the vanishing in the first claim.
Then, we obtain the second claim of
Proposition \ref{prop;14.10.18.10}.
\hfill\qed

\vspace{.1in}

We have the canonical nilpotent map $N$ on
$\psi_{\tau}\pitilde^0_{+}L_{\ast}(\nbigatilde,\Stilde,\gammatilde)
\simeq
 \psi_{\tau}\pitilde_{+}^0L_!(\nbigatilde,\Stilde,\gammatilde)$
induced by $t\del_t$.

\begin{cor}
\label{cor;14.11.28.21}
We have the following commutative diagram:
\[
 \begin{CD}
 \Ker N @>>>
 \Cok N\\
 @V{\simeq}VV @V{\simeq}VV \\
 \pi^0_{+}L_{!}(\nbiga,S,\gamma)
@>>>
 \pi^0_{+}L_{\ast}(\nbiga,S,\gamma). 
 \end{CD}
\]
The lower horizontal arrow is the natural morphism.
\end{cor}
\pf
It follows from Proposition \ref{prop;14.10.18.20}.
\hfill\qed

\subsubsection{Graded polarized mixed twistor $\nbigd$-modules}
\label{subsection;14.11.25.100}

We continue to use the setting in 
\S\ref{subsection;14.7.4.40},
and we consider the associated mixed twistor $\nbigd$-modules.
Let us show the relations between the mixed twistor $\nbigd$-modules
associated 
$(\nbiga,S,\gamma)$ 
and $(\nbigatilde,\Stilde,\gammatilde)$,
although we have already used it in \S\ref{subsection;14.11.24.121}
implicitly.

\begin{prop}
\mbox{{}}
\label{prop;14.11.17.1}
\begin{itemize}
\item
We have 
 $\nbigt_{!}(\nbigatilde,\Stilde,\gammatilde)
\simeq
 \nbigt_{\ast}(\nbigatilde,\Stilde,\gammatilde)
\simeq
 \nbigt_{\min}(\nbigatilde,\Stilde,\gammatilde)$,
and 
$\pitilde_{\dagger}^i\nbigt_{\min}(\nbigatilde,\Stilde,\gammatilde)=0$
for $i\neq 0$.
\item
Let $\nbigk(\nbigatilde,\Stilde,\gammatilde)$ 
and $\nbigc(\nbigatilde,\Stilde,\gammatilde)$
denote the kernel and the cokernel of the following morphism
of the mixed twistor $\nbigd$-modules:
\[
 \pitilde^0_{\dagger}\nbigt_{\ast}(\nbigatilde,\Stilde,\gammatilde)[!\tau]
\lrarr
\pitilde^0_{\dagger}\nbigt_{\ast}(\nbigatilde,\Stilde,\gammatilde)[\ast \tau].
\]
Then, we have the following natural isomorphisms:
\[
 \nbigk(\nbigatilde,\Stilde,\gammatilde)
\simeq
\iota_{\dagger}\pi_{\dagger}^0
 \nbigt_{!}(\nbiga,S,\gamma),
\quad
 \nbigc(\nbigatilde,\Stilde,\gammatilde)
\simeq 
  \iota_{\dagger}\pi_{\dagger}^0
 \nbigt_{\ast}(\nbiga,S,\gamma)
 \otimes\newTate(-1).
\]
The isomorphisms are compatible
with the real structures.
\end{itemize}

\end{prop}
\pf
We obtain the first claim from Proposition \ref{prop;14.10.18.10}.
As in the proof of Proposition \ref{prop;14.10.18.10},
by using the vanishing of the cohomology 
in Proposition \ref{prop;14.11.8.110}
and the first claim of this proposition,
we obtain the isomorphisms from Proposition \ref{prop;14.10.18.11}.
We can easily compare the real structures
by using Proposition \ref{prop;14.9.12.11}.
\hfill\qed

\vspace{.1in}

On one hand,
we have the graded polarizations
$\vecnbigs_{\star}(\nbiga,S,\gamma,\veca_{m+1})$
on $\pi_{\dagger}^0\nbigt_{\star}(\nbiga,S,\gamma)$
as explained in \S\ref{subsection;14.11.24.121}.
They induce graded sesqui-linear dualities on 
$\iota_{\dagger}\pi_{\dagger}^0
 \nbigt_!(\nbiga,S,\gamma)$
and 
$\iota_{\dagger}\pi_{\dagger}^0
 \nbigt_{\ast}(\nbiga,S,\gamma)
 \otimes\newTate(-1)$.
On the other hand,
as explained in \S\ref{subsection;14.11.24.121},
the pure twistor $\nbigd$-module
$\pi_{\dagger}^0\nbigt_{\min}(\nbigatilde,\Stilde,\gammatilde)$
is equipped with the induced polarization.
It induces the graded polarizations
$\vecnbigs_{\nbigk}$ and $\vecnbigs_{\nbigc}$ of 
 $\nbigk(\nbigatilde,\Stilde,\gammatilde)$ 
and $\nbigc(\nbigatilde,\Stilde,\gammatilde)$,
as explained in \S\ref{subsection;14.11.24.1}.

\begin{prop}
\label{prop;14.11.17.3}
The isomorphisms in Proposition 
{\rm\ref{prop;14.11.17.1}}
are compatible with the graded polarizations.
\end{prop}
\pf
We use the notation in the proof of 
Proposition \ref{prop;14.10.18.10}.
Let $\nbigk_1$ and $\nbigc_1$
be the kernel and the cokernel of
the morphism of the mixed twistor $\nbigd$-modules
$\nbigt_!(F_{\gammatilde,\Sigma},D_{\Sigma,\Stilde})
\lrarr
 \nbigt_{\ast}(F_{\gammatilde,\Sigma},D_{\Sigma,\Stilde})$.
According to Proposition \ref{prop;14.9.7.42},
the isomorphisms
$\nbigk_1\simeq
 \iota_{1\dagger}
 \nbigt_{!}(F_{\gamma,\Sigma},D_{\Sigma,S})$
and 
$\nbigc_1\simeq
 \iota_{1\dagger}
 \nbigt_{\ast}(F_{\gamma,\Sigma},D_{\Sigma,S})$
are compatible with the graded sesqui-linear dualities.
Then, the claims of the propositions
follow from Corollary \ref{cor;14.11.17.10}
and the construction of the graded polarizations.
\hfill\qed

\vspace{.1in}

Let $\nbign:
 \psitilde_{\tau,-\vecdelta}
 \pitilde_{\dagger}^0\nbigt(\nbigatilde,\Stilde,\gammatilde)
 \otimes\nbigu(-1,0)
\lrarr
 \psitilde_{\tau,-\vecdelta}
 \pitilde_{\dagger}^0\nbigt(\nbigatilde,\Stilde,\gammatilde)
\otimes\nbigu(0,-1)$
be the morphism given by the pair
$(-t\del_t,-t\del_t)$.
\begin{cor}
\label{cor;14.11.28.50}
We have the following commutative diagram:
\begin{equation}
 \label{eq;14.11.28.20}
\begin{CD}
  \Ker(\nbign)
@>>>
 \Cok(\nbign)\otimes\newTate(1)\\
 @V{\simeq}VV @V{\simeq}VV \\
 \pi_{\dagger}^0\nbigt_{!}(\nbiga,S,\gamma)
@>>>
 \pi_{\dagger}^0\nbigt_{\ast}(\nbiga,S,\gamma)
\end{CD}
\end{equation}
The horizontal arrows are the natural morphisms.
\end{cor}
\pf
Because 
$\iota_{\dagger}\Cok(\nbign)$
and 
$\iota_{\dagger}\Ker(\nbign)$
are naturally isomorphic to
$\nbigc(\nbigatilde,\Stilde,\gammatilde)$
and 
$\nbigk(\nbigatilde,\Stilde,\gammatilde)$,
we obtain the vertical isomorphisms
 in (\ref{eq;14.11.28.20})
from Proposition \ref{prop;14.11.17.3}.
The commutativity of the diagram follows from 
Corollary \ref{cor;14.11.28.21}.
\hfill\qed

\subsubsection{The induced mixed TEP-structures}
\label{subsection;14.11.18.51}

We set 
$\nbigv:=
 \pitilde_{\dagger}^0\nbigl(\nbigatilde,\Stilde,\gammatilde)(\ast\tau)$.
Let $\lambda N:\psi_{\tau,-\vecdelta}(\nbigv)
\lrarr
 \psi_{\tau,-\vecdelta}(\nbigv)$
be the morphism induced by
$\lambda\tau\del_{\tau}$.
We restate the result in 
\S\ref{subsection;14.11.28.30} in this situation.

\begin{assumption}
$\nbigv$ is regular singular along $\tau=0$.
There exists an open subset $\Btilde$ in $\Stilde$
such that $\nbigv_{|\Btilde}$ is a locally free
$\nbigo_{\cnum_{\lambda}\times \Btilde}(\ast \tau)$-module.
\hfill\qed
\end{assumption}

We set $B:=\iota^{-1}(\Btilde)$
which is an open subset in $S$.
Under the assumption
$\Cok(\lambda N)_{|B}$
is a locally free $\nbigo_{\cnum_{\lambda}\times B}$-module.
We also have the pairing 
$P_{\nbigatilde}$ of weight $n$ on $\nbigv_{|\Btilde}$
induced by the real structure and the polarization
of $\pitilde_{\dagger}^0\nbigt(\nbigatilde,\Stilde,\gammatilde)$.
By the procedure in \S\ref{subsection;14.11.28.130},
we have the weight filtration $\Wtilde^{(1)}$
and the graded pairing
$(P^{(1)}_{\nbigatilde,k}\,|\,k\in\seisuu)$
on $\Cok(\lambda N)_{|B}$
and $\Ker(\lambda N)_{|B}$.

\vspace{.1in}
Note that 
$\Cok(\lambda N)$ is the underlying
$\nbigr$-module of
$\Cok(\nbign)\otimes\newTate(1)
\simeq
\pi_{\dagger}^0\nbigt_{\ast}(\nbiga,S,\gamma)$,
i.e.,
$\Cok(\lambda N)$
is isomorphic to
$\pi_{\dagger}^0\nbigl_{\ast}(\nbiga,S,\gamma)$.
As explained in \S\ref{subsection;14.11.26.131},
we obtain the weight filtration
$\Wtilde^{(2)}$ and the graded pairing 
$(P^{(2)}_{\nbigatilde,k}\,|\,k\in\seisuu)$ on
$\Cok(\lambda N)_{|B}$
from the real structure and the graded polarization 
of $\pi_{\dagger}^0\nbigt(\nbiga,S,\gamma)$.
Similarly,
$\Ker(\lambda N)$ is the underlying $\nbigr$-module of
$\Ker(\nbign)\simeq 
 \pi_{\dagger}^0\nbigt_!(\nbiga,S,\gamma)$,
i.e.,
$\Ker(\lambda N)$ is isomorphic to
$\pi_{\dagger}^0\nbigl_{!}(\nbiga,S,\gamma)$.
As explained in \S\ref{subsection;14.11.26.131},
we have the filtration $\Wtilde^{(2)}$
and the graded pairing
$(P^{(2)}_{\nbigatilde,k}\,|\,k\in\seisuu)$
on $\Ker(\lambda N)$.

\begin{prop}
\label{prop;14.11.18.52}
On $\Cok(\lambda N)$ and $\Ker(\lambda N)$,
we have $\Wtilde^{(1)}=\Wtilde^{(2)}$
and 
$P^{(1)}_{\nbigatilde,k}
=P^{(2)}_{\nbigatilde,k}$
for any $k\in\seisuu$.
\end{prop}
\pf
It follows from 
Proposition \ref{prop;14.11.28.32}
and the isomorphisms
$\Cok(\nbign)\otimes\newTate(1)\simeq
\pi_{\dagger}^0\nbigt_{\ast}(\nbiga,S,\gamma)$
and 
$\Ker(\nbign)\simeq
\pi_{\dagger}^0\nbigt_{!}(\nbiga,S,\gamma)$
in Corollary \ref{cor;14.11.28.50}.
\hfill\qed

\subsection{A regular singular case}
\label{subsection;14.12.3.120}

\subsubsection{The $\nbigd$-modules}
\label{subsection;14.7.1.13}

Let $\nbiga=\{\veca_1,\ldots,\veca_m\}$
be as in \S\ref{subsection;14.7.8.1}.
Suppose that there exists
$\alpha\in (\seisuu^n)^{\lor}$ such that 
$\alpha(\veca_i)=1$ for any $\veca_i\in\nbiga$.
In \cite{Batyrev-VMHS},
\cite{Konishi-Minabe}
and \cite{Stienstra},
the $\nbigd$-modules are described
in terms of the Gauss-Manin connection
for the relative cohomology groups
of some families.
We review it 
in the language of $\nbigd$-modules,
which fits to the theory of 
mixed twistor $\nbigd$-modules.

With an appropriate choice of a frame of $\seisuu^n$,
we may assume that
$\veca_i=(\vecb_i,1)\in\seisuu^{n-1}\times\seisuu$
for each $\veca_i\in\nbiga$.
Let $\nbigb\subset\seisuu^{n-1}$ 
denote the image of $\nbiga$ via 
the projection of
$\seisuu^{n-1}\times\seisuu\lrarr\seisuu^{n-1}$.
In other words,
we have
$\nbiga=\{(\vecb,1)\,|\,\vecb\in\nbigb\}$.
For simplicity,
we impose the following
in this subsection.

\begin{assumption}
We assume that $\nbigb$ generates $\seisuu^{n-1}$,
and $0$ is contained in the interior part of
$\Conv(\nbigb)$.
\hfill\qed
\end{assumption}

Let $X_{\nbigb}$ denote the closure of the image of
the morphism $\psi_{\nbigb}:T^{n-1}\lrarr \proj^{m}$
as in \S\ref{subsection;14.7.8.1}.
We set
$W:=\Bigl\{\sum_{i=1}^m \alpha_iz_i\,\Big|\,
 \alpha_i\in\cnum
 \Bigr\}
\subset
 H^0\bigl(\proj^m,\nbigo(1)\bigr)
=\Bigl\{
 \sum_{i=0}^m\alpha_iz_i 
 \Bigr\}$.
We use the notation in \S\ref{subsection;14.11.23.240}.

\begin{lem}
\label{lem;14.7.1.11}
Let $\gamma:S\lrarr W$ be
a morphism of complex manifolds.
The following conditions are equivalent.
\begin{itemize}
\item
$\gamma$ is non-degenerate at $\infty$ for $X_{\nbiga}$.
\item
If $\gamma=\sum_{i=1}^m \gamma_iz_i$,
the family of Laurent polynomials
$F_{\gamma,\nbigb}
 =\sum_{i=1}^m \gamma_it^{\vecb_i}$
is $\Conv(\nbigb)$-regular
in the sense that 
$F_{\gamma,\nbigb,\sigma}^{-1}(0)
\cap
 (dF_{\gamma,\nbigb,\sigma})^{-1}(0)
 =\emptyset$
for any face $\sigma$ of $\Conv(\nbigb)$.
\end{itemize}
\end{lem}
\pf
We have the family of Laurent polynomials
$F_{\gamma}=\sum \gamma_i t^{\veca_i}$
associated to $\nbiga$ and $\gamma$.
For any face $\sigma$ of $\nbigb$,
we have the face 
$\sigma(\nbiga)$ of $\nbiga$
given by
$\sigma(\nbiga):=\bigl\{
 (1,\vecc)\,\big|\,\vecc\in\sigma
 \bigr\}$.
We have
$F_{\gamma,\sigma(\nbiga)}
=t_nF_{\gamma,\nbigb,\sigma}$
and 
$dF_{\gamma,\sigma(\nbiga)}=
t_ndF_{\gamma,\nbigb,\sigma}
+dt_n F_{\gamma,\nbigb,\sigma}$.
We have
$F_{\gamma,\nbigb,\sigma}^{-1}(0)
=F_{\gamma,\sigma(\nbiga)}^{-1}(0)$
and 
$(dF_{\gamma,\sigma(\nbiga)})^{-1}(0)
=
 F_{\gamma,\nbigb,\sigma}^{-1}(0)
\cap
 (dF_{\gamma,\nbigb,\sigma})^{-1}(0)$.
Then, the claim is clear.
\hfill\qed

\vspace{.1in}

Let $\gamma:S\lrarr W$
be any holomorphic map.
Let us study 
the $\nbigd$-modules
$L_{\star}(\nbiga,S,\gamma)$
and 
$\pi^i_{+}\bigl(
 L_{\star}(\nbiga,S,\gamma)
 \bigr)$.
It is convenient to consider
a toric compactification
which is different from that in \S\ref{subsection;14.7.4.10}.
We set 
$X'_{\nbigb}:=X_{\nbigb}\times\proj^1$
which is naturally a toric variety.
We take a toric desingularization
$\varphi_{\Sigma_0}:
 X_{\Sigma_0}\lrarr X_{\nbigb}$.
We naturally have a fan $\Sigma_0'$
such that 
$X_{\Sigma_0'}=X_{\Sigma_0}\times\proj^1$.
We set
$\varphi_{\Sigma_0'}:=
\varphi_{\Sigma_0}\times\id:
 X_{\Sigma'_0}\lrarr
 X'_{\nbigb}$.
Let $\varphibar_{\Sigma_0',S}$
denote the composite of
$\varphi_{\Sigma_0'}\times\id_S$
and the inclusion
$X_{\nbigb,S}'\lrarr
 \proj^m\times\proj^1\times S$.
We have the meromorphic function
$F_{\gamma,\Sigma_0'}$ on
$(X_{\Sigma_0',S},D_{\Sigma_0',S})$
determined by 
the family of Laurent polynomials $F_{\gamma}$.
(See \S\ref{subsection;14.11.23.210}
for $F_{\gamma}$.)
For $\star=\ast,!$,
we have the $\nbigd$-modules
$L_{\star}(F_{\gamma,\Sigma'_0},D_{\Sigma'_0,S})$
on $X_{\Sigma_0',S}$.
We set
\[
 L'_{\star}(\nbigb,S,\gamma):=
 \varphibar_{\Sigma_0',S+}
  L_{\star}(F_{\gamma,\Sigma'_0},D_{\Sigma'_0,S})
\]
on $\proj^m\times\proj^1\times S$.
It is easy to see that they are independent of
the choice of a desingularization $\varphi_{\Sigma_0}$
up to natural isomorphisms.
Let $\pi_{\nbigb}:\proj^m\times\proj^1\times S \lrarr S$
and 
$\pi_{\Sigma_0'}:X_{\Sigma_0',S}\lrarr S$
denote the projections.
We have the natural isomorphisms
$\pi^i_{\nbigb+}
 L'_{\star}(\nbigb,S,\gamma)
\simeq
 \pi^i_{\Sigma'_0+}
 L_{\star}(F_{\gamma,\Sigma_0'},D_{\Sigma_0',S})$
for $\star=\ast,!$.

\begin{lem}
We have natural isomorphisms
$\pi^i_{\nbigb+}L'_{\star}(\nbigb,S,\gamma)
\simeq
 \pi^i_+L_{\star}(\nbiga,S,\gamma)$.
\end{lem}
\pf
We take a toric desingularization
$\varphi_{\Sigma_1}:X_{\Sigma_1}\lrarr X_{\nbiga}$.
We may assume that we have
a toric morphism
$\rho:X_{\Sigma_1}\lrarr X_{\Sigma_0'}$
which is birational.
We have natural isomorphisms
$\rho_+L_{\star}(F_{\gamma,\Sigma_1},D_{\Sigma_1,S})
\simeq
 L_{\star}(F_{\gamma,\Sigma_0'},D_{\Sigma_0',S})$.
They induce the desired isomorphisms.
\hfill\qed

\vspace{.1in}

The family of Laurent polynomials
$F_{\gamma}$ is described as
\[
 F_{\gamma}=t_nf_{\gamma},
\quad\quad
 f_{\gamma}:=\sum_{i=1}^m \gamma_it^{\vecb_i}.
\]
Let $Z_{f_{\gamma}}\subset X_{\Sigma_0,S}$
denote the zero set of $f_{\gamma}$.
Suppose that $\gamma$ is non-degenerate at $\infty$
for $X_{\nbiga}$.
As remarked in Lemma \ref{lem;14.7.1.11},
the family of Laurent polynomials
$f_{\gamma}$ is $\Conv(\nbigb)$-regular.
According to Proposition \ref{prop;14.7.1.1},
we have
\begin{equation}
\label{eq;14.6.28.20}
 \pi^0_{+}L_{\ast}(\nbiga,S,\gamma)
\simeq
 \pi^0_{\Sigma'_0+}
 L_{\ast}(F_{\gamma,\Sigma_0'},D_{\Sigma_0',S})
\simeq
 \pi^0_{\Sigma_0+}
 \bigl(
 \nbigo_{X_{\Sigma_0,S}}
 (!Z_{f_{\gamma}})
(\ast D_{\Sigma_0,S})
 \bigr).
\end{equation}
\begin{equation}
\label{eq;14.11.26.1}
 \pi_{+}^0L_{!}(\nbiga,S,\gamma)
\simeq
 \pi_{\Sigma_0'+}^0
 L_{!}(F_{\gamma,\Sigma_0'},D_{\Sigma_0',S})
\simeq
  \pi^0_{\Sigma_0+}
 \bigl(
 \nbigo_{X_{\Sigma_0,S}}
 (\ast Z_{f_{\gamma}})
(!D_{\Sigma_0,S})
 \bigr).
\end{equation}
We have
$\pi^i_{\Sigma_0+}
 \bigl(
 \nbigo_{X_{\Sigma_0,S}}
 (!Z_{f_{\gamma}})(\ast D_{\Sigma_0,S})
 \bigr)=0$
and 
$\pi^i_{\Sigma_0+}
 \bigl(
 \nbigo_{X_{\Sigma_0,S}}
 (\ast Z_{f_{\gamma}})(!D_{\Sigma_0,S})
 \bigr)=0$
for $i\neq 0$,
which follow from
$\pi^i_{\Sigma+}
 L_{\star}(\nbiga,S,\gamma)=0$
for $i\neq 0$.

If moreover the image of $\gamma$
is contained in $W\cap U_{\nbiga}^{\reg}$,
the $\nbigd_S$-modules
(\ref{eq;14.6.28.20})
and (\ref{eq;14.11.26.1})
are flat bundles on $S$.
The fiber of (\ref{eq;14.6.28.20})
over $s\in S$ is 
the relative cohomology group
$H^{n-1}(T^{n-1},Z^{\circ}_{f_{\gamma(s)}})$
with $\cnum$-coefficient,
where 
$Z^{\circ}_{f_{\gamma(s)}}:=Z_{f_{\gamma(s)}}\cap T^{n-1}$.
The fiber of (\ref{eq;14.11.26.1}) over $s\in S$ is
$H^{n-1}(X_{\Sigma_0}\setminus Z_{f_{\gamma(s)}},
 D_{\Sigma_0}\setminus Z_{f_{\gamma(s)}})$.

\begin{rem}
Under the identification
$X_{\Sigma_0'}=\proj^1_s\times X_{\Sigma_0}$,
as remarked in 
{\rm \S\ref{subsection;14.7.1.12}},
{\rm(\ref{eq;14.6.28.20})}
 is also naturally isomorphic to
$\pi^0_{\Sigma_0'+}
 \bigl(
 \nbigo_{X_{\Sigma_0}',S}
 (!Z_{f_{\gamma}+s})
 (\ast D_{\Sigma_0'})
 \bigr)$.
\hfill\qed
\end{rem}

\subsubsection{The mixed twistor $\nbigd$-modules}

We continue to use the notation in \S\ref{subsection;14.7.1.13}.
Suppose that $\gamma$ is non-degenerate at $\infty$
for $X_{\nbiga}$.
We have the mixed twistor $\nbigd$-modules
$\nbigt_{\star}(F_{\gamma,\Sigma_0'},D_{\Sigma_0',S})$
$(\star=\ast,!)$
on $X_{\Sigma_0',S}$.
We obtain
the mixed twistor $\nbigd$-modules
on $\proj^m\times\proj^1\times S$:
\[
 \nbigt'_{\star}(\nbigb,S,\gamma)
:=\varphibar_{\Sigma_0',S\dagger}
 \nbigt_{\star}(F_{\gamma,\Sigma_0'},D_{\Sigma_0',S})
\]
They are independent of the choice of
a toric desingularization
$\varphi_{\Sigma_0}:X_{\Sigma_0}\lrarr X_{\nbigb}$.
We have the natural isomorphisms
$\pi^i_{\nbigb\dagger}
 \nbigt'_{\star}(\nbigb,S,\gamma)
\simeq
 \pi^i_{\Sigma'_0\dagger}
 \nbigt_{\star}(F_{\gamma,\Sigma_0'},D_{\Sigma_0',S})$
$(\star=\ast,!)$.
As in the case of $\nbigd$-modules,
we have the following.

\begin{lem}
We have natural isomorphisms
$\pi^i_{\nbigb\dagger}
 \nbigt'_{\star}(\nbigb,S,\gamma)
\simeq
 \pi^i_{\dagger}
 \nbigt_{\star}(\nbiga,S,\gamma)$
for $\star=\ast,!$.
\hfill\qed
\end{lem}

Let $\nbigl_{\star}(F_{\gamma,\Sigma_0'},D_{\Sigma_0',S})$
$(\star=\ast,!)$
be the $\nbigr$-modules underlying 
$\nbigt_{\star}(F_{\gamma,\Sigma_0'},D_{\Sigma_0',S})$
on $X_{\Sigma_0',S}$.
Set
$\nbigl_{\star}'(\nbigb,S,\gamma):=
 \varphibar^0_{\Sigma_0',S\dagger}
\nbigt_{\star}(F_{\gamma,\Sigma_0'},D_{\Sigma_0',S})$
on $\proj^m\times\proj^1\times S$,
which underlie $\pi_{\dagger}^0\nbigt'_{\star}(\nbigb,S,\gamma)$.
We consider the standard $\cnum^{\ast}$-action
on $\proj^1$ given by $a\cdot [y_0:y_1]=[ay_0:y_1]$.
It induces $\cnum^{\ast}$-actions 
on $X_{\Sigma_0',S}$
and $\proj^m\times \proj^1\times S$.

\begin{lem}
\label{lem;14.11.26.100}
The $\nbigr$-modules
$\nbigl'_{\star}(\nbigb,S,\gamma)$
and 
$\nbigl_{\star}(F_{\gamma,\Sigma_0'},D_{\Sigma_0',S})$
are $\cnum^{\ast}$-homogeneous.
(See {\rm\S\ref{subsection;14.11.25.10}}
for the notion of homogeneity.)
\end{lem}
\pf
We consider the $\cnum^{\ast}$-action on $T^n$
given by $a\cdot (t_1,\ldots,t_{n-1},t_n)=(t_1,\ldots,t_{n-1},at_n)$.
It induces a $\cnum^{\ast}$-action on
$X_{\Sigma_0',S}$.
We have 
$a^{\ast}F_{\gamma,\Sigma_0'}=a\cdot F_{\gamma,\Sigma_0'}$.
By Proposition \ref{prop;14.11.25.20},
the $\nbigr$-modules
$\nbigl_{\star}(F_{\gamma,\Sigma_0'},D_{\Sigma_0',S})$
are $\cnum^{\ast}$-homogeneous.
Then, we obtain that 
$\nbigl'_{\star}(\nbigb,S,\gamma)$
are also $\cnum^{\ast}$-homogeneous.
\hfill\qed

\vspace{.1in}
Let $q_{\nbigb}:\proj^m\times\proj^1\times S\lrarr \proj^m\times S$
and $q_{\Sigma_0}:X_{\Sigma_0',S}\lrarr X_{\Sigma_0,S}$
be the projections.

\begin{cor}
\label{cor;14.11.26.111}
The $\nbigr$-modules
$q_{\nbigb\dagger}\nbigl'_{\star}(\nbigb,S,\gamma)$
and 
$q_{\Sigma_0\dagger}
 \nbigl_{\star}(F_{\gamma,\Sigma_0'},D_{\Sigma_0',S})$
are $\cnum^{\ast}$-homogeneous
for the trivial $\cnum^{\ast}$-actions
on $\proj^m\times S$ and 
$X_{\Sigma_0}\times S$.
Hence,
$q_{\nbigb\dagger}\nbigt'_{\star}(\nbigb,S,\gamma)$
and
$q_{\Sigma_0\dagger}
 \nbigt_{\star}(F_{\gamma,\Sigma_0'},D_{\Sigma_0',S})$
come from mixed Hodge modules.
In particular,
the underlying $\nbigd$-modules
are equipped with the good filtrations
such that the analytification of the Rees modules
are isomorphic to 
$q_{\nbigb\dagger}\nbigl'_{\star}(\nbigb,S,\gamma)$
and 
$q_{\Sigma_0\dagger}
 \nbigl_{\star}(F_{\gamma,\Sigma_0'},D_{\Sigma_0',S})$.
\end{cor}
\pf
Because
$q_{\nbigb}$ and $q_{\Sigma_0}$
are $\cnum^{\ast}$-equivariant,
the first claim follows from Lemma \ref{lem;14.11.26.100}.
The second follows from 
Theorem \ref{thm;15.11.17.20}.
\hfill\qed

\begin{cor}
\label{cor;14.11.26.110}
The $\nbigr$-modules
$\pi^0_{\nbigb\dagger}\nbigl'_{\star}(\nbigb,S,\gamma)
\simeq 
\pi^0_{\Sigma'_0\dagger}
 \nbigl_{\star}(F_{\gamma,\Sigma_0'},D_{\Sigma_0',S})$
are $\cnum^{\ast}$-homogeneous
for the trivial $\cnum^{\ast}$-actions on $S$.
Hence,
$\pi_{\nbigb\dagger}^0\nbigt'_{\star}(\nbigb,S,\gamma)$
come from mixed Hodge modules.
In particular,
the underlying $\nbigd$-modules
are equipped with the good filtrations
such that the analytification of the Rees modules
are isomorphic to 
$\pi^0_{\nbigb\dagger}\nbigl'_{\star}(\nbigb,S,\gamma)$.
\hfill\qed
\end{cor}

\begin{prop}
\label{prop;14.11.18.2}
Suppose that $\gamma:S\lrarr W$ is non-degenerate
at $\infty$ for $X_{\nbiga}$.
Then, we have the following natural isomorphisms
of the integrable mixed twistor $\nbigd$-modules
compatible with the real structure:
\begin{equation}
 \label{eq;14.11.18.3}
 \pi^0_{\dagger}\nbigt_{\ast}(\nbiga,S,\gamma)
\simeq
 \pi^0_{\Sigma_0\dagger}\Bigl(
 \bigl(
 \nbigu_{X_{\Sigma_0,S}}(n+d_S,-1)
 \bigr)
 [!Z_{f_{\gamma}}][\ast D_{\Sigma_0,S}]
\Bigr).
\end{equation}
\begin{equation}
 \label{eq;14.11.29.100}
 \pi^0_{\dagger}\nbigt_!(\nbiga,S,\gamma)
\simeq
 \pi^0_{\Sigma_0\dagger}\Bigl(
 \bigl(
 \nbigu_{X_{\Sigma_0,S}}(n+d_S-1,0)
 \bigr)
 [\ast Z_{f_{\gamma}}][!D_{\Sigma_0,S}]
\Bigr).
\end{equation}
\end{prop}
\pf
We give an argument for (\ref{eq;14.11.18.3}).
The other case is similar.
According to Proposition \ref{prop;14.11.7.21},
we have the isomorphism of 
the integrable mixed twistor $\nbigd$-modules
with real structure:
\begin{equation}
 \label{eq;14.11.18.1}
q^0_{\Sigma_0\dagger}\nbigt_{\ast}
 (F_{\gamma,\Sigma_0'},D_{\Sigma_0',S})
\simeq
 \bigl(
 \nbigu_{X_{\Sigma_0,S}}(n-1+d_S,0)
\otimes\newTate(-1)
 \bigr)
 [!Z_{f_{\gamma}}][\ast D_{\Sigma_0,S}].
\end{equation}
Thus, we obtain the following isomorphisms:
\begin{equation}
\label{eq;14.11.11.110}
 \pi^0_{\dagger}\nbigt_{\ast}(\nbiga,S,\gamma)
\simeq
 \pi^0_{\Sigma_0'\dagger}\nbigt_{\ast}
 (F_{\gamma,\Sigma_0'},D_{\Sigma_0',S})
\simeq
 \pi^0_{\Sigma_0\dagger}\Bigl(
 \bigl(
 \nbigu_{X_{\Sigma_0,S}}(n-1+d_S,0)
\otimes\newTate(-1)
 \bigr)
 [!Z_{f_{\gamma}}][\ast D_{\Sigma_0,S}]
\Bigr).
\end{equation}
Thus, we obtain the claim of the proposition.
\hfill\qed

\vspace{.1in}

Let $M$ be the pure Hodge module on $X_{\Sigma_0,S}$
of weight $\dim S+n-1$,
corresponding to the constant sheaf.
Let $M(-1)$ be the $(-1)$-th Tate twist of $M$.
We obtain the mixed Hodge module
$\pi_{\Sigma_0\dagger}^0M(-1)[!Z_{f_{\gamma}}][\ast D_{\Sigma_0,S}]$
which induces the mixed twistor $\nbigd$-module
in the right hand side of (\ref{eq;14.11.18.3}).
We also have the mixed Hodge module
$\pi_{\Sigma_0\dagger}^0M[\ast Z_{f_{\gamma}}][!D_{\Sigma_0,S}]$.
\begin{cor}
The Hodge filtrations of 
$\pi_{\Sigma_0\dagger}^0M(-1)[!Z_{f_{\gamma}}][\ast D_{\Sigma_0,S}]$
and 
$\pi_{\Sigma_0\dagger}^0M[\ast Z_{f_{\gamma}}][!D_{\Sigma_0,S}]$
are equal to 
the good filtrations in Corollary {\rm\ref{cor;14.11.26.110}}.
\hfill\qed
\end{cor}

Let $\iota:Z_{f_{\gamma}}\lrarr X_{\Sigma_0,S}$
be the inclusion of the complex manifolds.
We have the following natural morphisms 
of integrable mixed twistor $\nbigd$-modules on $S$:
\begin{multline}
 \pi_{\Sigma_0\dagger}^0
 \nbigu_{X_{\Sigma_0,S}}(n+d_S-1,0)[\ast Z_{f_{\gamma}}][!D_{\Sigma_0,S}]
\stackrel{a_1}{\lrarr}
 \pi_{\Sigma_0\dagger}^0
 \iota_{\dagger}\nbigu_{Z_{f_{\gamma}}}(n+d_S-1,-1)
 \\
\stackrel{a_2}{\lrarr}
 \pi_{\Sigma_0\dagger}^0
 \nbigu_{X_{\Sigma_0,S}}(n+d_S,-1)[!Z_{f_{\gamma}}][\ast D_{\Sigma_0,S}]
\end{multline}

\begin{prop}
\label{prop;14.12.3.110}
We have 
$\pi_{\Sigma_0\dagger}^0
 \iota_{\dagger}\nbigu_{Z_{f_{\gamma}}}(n+d_S-1,-1)
=\Image a_1\oplus\Ker a_2$.
We also have
\begin{multline}
 \Gr^W_{n+d_S}
 \Bigl(
  \pi_{\Sigma_0\dagger}^0
 \nbigu_{X_{\Sigma_0,S}}(n+d_S-1,0)
 [\ast Z_{f_{\gamma}}][!D_{\Sigma_0,S}]
 \Bigr)
\simeq  
 \Image a_1
 \\
\simeq
 \Gr^{W}_{n+d_S}\Bigl(
 \pi_{\Sigma_0\dagger}^0
 \nbigu_{X_{\Sigma_0,S}}(n+d_S,-1)
 [!Z_{f_{\gamma}}][\ast D_{\Sigma_0,S}]
 \Bigr)
\end{multline}
\end{prop}
\pf
According to
Corollary \ref{cor;14.11.29.400},
the morphism $a_2\circ a_1$ induces
an isomorphism 
of the weight $(n+d_S)$-parts.
Hence,
we obtain $\Image (a_2)=\Image(a_2\circ a_2)$,
which implies the claims of the proposition.
\hfill\qed

\subsubsection{Graded polarizations and induced pairings}
\label{subsection;14.12.8.20}

Let $\vecb=(0,-1)\in\seisuu^{n-1}\times\seisuu$.
Note that $0$ is contained in the interior part of
$\Conv(\nbiga\cup\{\vecb\})$.

Let $\gamma:S\lrarr W$ be a morphism
which is non-degenerate at $\infty$ for $X_{\nbiga}$.
As explained in \S\ref{subsection;14.11.24.121},
we obtain graded polarizations 
$\vecnbigs_{\star}(\nbiga,S,\gamma,\vecb)$ of 
$\pi_{\dagger}^0\nbigt_{\star}(\nbiga,S,\gamma)$.
Suppose moreover that 
the image of $\gamma$ is contained in $W\cap U_{\nbiga}^{\reg}$.
Then, 
$\nbigv'_{\star}(\nbigb,S,\gamma):=
 \pi^0_{\Sigma_0'\dagger}
 \nbigl_{\star}(F_{\gamma,\Sigma_0'},D_{\Sigma_0',S})$
are locally free $\nbigo_{\cnum_{\lambda}\times S}$-modules.
Let $W$ denote the filtration of 
$\nbigv'_{\star}(\nbigb,S,\gamma)$
underlying the weight filtration of
$\pi_{\Sigma_0'\dagger}^0
 \nbigt_{\star}(F_{\gamma,\Sigma_0'},D_{\Sigma_0',S})$.
We set 
$\Wtilde_k
 \nbigv'_{\star}(\nbigb,S,\gamma):=
  W_{d_S+k}
 \nbigv'_{\star}(\nbigb,S,\gamma)$.
Then, the polarization and the real structure
of 
$\Gr^W_{d_S+k}
 \pi_{\Sigma_0',\dagger}^0
 \nbigt_{\star}(F_{\gamma,\Sigma_0'},D_{\Sigma_0',S})$
induce a pairing $\Ptilde_k$ of weight $k$
on $\Gr^{\Wtilde}_k\nbigv'_{\star}(\nbigb,S,\gamma)$,
as explained in \S\ref{subsection;14.11.26.131}.
Let us give a description
of the pairings $\Ptilde_{k}$ in the case $\star=\ast$.

\vspace{.1in}
We have the following isomorphism:
\begin{equation}
\label{eq;14.11.26.130}
\nbigv'_{\ast}(\nbigb,S,\gamma)
\simeq
 \pi^0_{\Sigma_0\dagger}
 \nbigo_{\cnum_{\lambda}\times X_{\Sigma_0,S}}
 [!Z_{f_{\gamma}}][\ast D_{\Sigma_0,S}]
\end{equation}
Note that 
$\nbigo_{\cnum_{\lambda}\times X_{\Sigma_0,S}}
 [!Z_{f_{\gamma}}][\ast D_{\Sigma_0,S}]$
is equipped with the filtration $W$
which underlies the weight filtration of
$\nbigu_{X_{\Sigma_0,S}}(n+d_S,-1)[!Z_{f_{\gamma}}][\ast D_{\Sigma_0,S}]$.
By using the spectral sequence,
we have the following complex
\begin{multline}
\label{eq;14.12.13.40}
  \pi^{-1}_{\Sigma_0\dagger}
 \Gr^W_{d_S+n+k+1}\Bigl(
 \lambda^{-1}
 \nbigo_{\cnum\times X_{\Sigma_0,S}}
 [!Z_{f_{\gamma}}][\ast D_{\Sigma_0,S}]
 \Bigr)
\stackrel{\alpha_{k}}{\lrarr}
 \\
 \pi^0_{\Sigma_0\dagger}
 \Gr^W_{d_S+n+k}\Bigl(
 \lambda^{-1}
 \nbigo_{\cnum\times X_{\Sigma_0,S}}
 [!Z_{f_{\gamma}}][\ast D_{\Sigma_0,S}]
 \Bigr)
\stackrel{\beta_{k}}{\lrarr}
 \\
 \pi^1_{\Sigma_0\dagger}
 \Gr^W_{d_S+n+k-1}\Bigl(
 \lambda^{-1}
 \nbigo_{\cnum\times X_{\Sigma_0,S}}
 [!Z_{f_{\gamma}}][\ast D_{\Sigma_0,S}]
 \Bigr),
\end{multline}
and the cohomology is isomorphic to
$\Gr^{\Wtilde}_{n+k}\nbigv'_{\ast}(\nbigb,S,\gamma)$.
The real structure and the polarization of
\[
 \pi^0_{\Sigma_0\dagger}\Gr^{W}_{d_S+n+k}
 \nbigu_{X_{\Sigma_0,S}}(n+d_S,-1)
 [!Z_{f_{\gamma}}][\ast D_{\Sigma_0,S}] 
\]
induces 
a pairing $\nbigptilde_{n+k}$ of weight $n+k$ on 
$\pi^0_{\Sigma_0\dagger}
 \Gr^W_{d_S+n+k}\Bigl(
 \lambda^{-1}
 \nbigo_{\cnum\times X_{\Sigma_0,S}}
 [!Z_{f_{\gamma}}][\ast D_{\Sigma_0,S}]
 \Bigr)$.
The pairing $\nbigptilde_{n+k}$ induces
the following isomorphism:
\[
 \Phi_{\nbigptilde_{n+k}}:
 \pi^0_{\Sigma_0\dagger}
 \Gr^W_{d_S+n+k}\Bigl(
 \lambda^{-1}
 \nbigo_{\cnum\times X_{\Sigma_0,S}}
 [!Z_{f_{\gamma}}][\ast D_{\Sigma_0,S}]
 \Bigr)
\simeq
 \lambda^{-n-k}
 \pi^0_{\Sigma_0\dagger}
 \Gr^W_{d_S+n+k}\Bigl(
 \lambda^{-1}
 \nbigo_{\cnum\times X_{\Sigma_0,S}}
 [!Z_{f_{\gamma}}][\ast D_{\Sigma_0,S}]
 \Bigr)^{\lor}
\]
As the dual of (\ref{eq;14.12.13.40}),
we have the following:
\begin{multline}
 \pi^1_{\Sigma_0\dagger}
 \Gr^W_{d_S+n+k-1}\Bigl(
 \lambda^{-1}
 \nbigo_{\cnum\times X_{\Sigma_0,S}}
 [!Z_{f_{\gamma}}][\ast D_{\Sigma_0,S}]
 \Bigr)^{\lor}
\stackrel{\beta^{\lor}_{k}}{\lrarr} \\
 \pi^0_{\Sigma_0\dagger}
 \Gr^W_{d_S+n+k}\Bigl(
 \lambda^{-1}
 \nbigo_{\cnum\times X_{\Sigma_0,S}}
 [!Z_{f_{\gamma}}][\ast D_{\Sigma_0,S}]
 \Bigr)^{\lor}
\stackrel{\alpha^{\lor}_{k}}{\lrarr}
 \\
  \pi^{-1}_{\Sigma_0\dagger}
 \Gr^W_{d_S+n+k+1}\Bigl(
 \lambda^{-1}
 \nbigo_{\cnum\times X_{\Sigma_0,S}}
 [!Z_{f_{\gamma}}][\ast D_{\Sigma_0,S}]
 \Bigr)^{\lor}
\end{multline}
We obtain the subsheaf
$\Ker\beta_{k}\cap
 \Ker\bigl(
 \alpha_{k}^{\lor}\circ\Phi_{\nbigptilde_{n+k}}
 \bigr)$
in 
$\pi^0_{\Sigma_0\dagger}
 \Gr^W_{d_S+n+k}\Bigl(
 \lambda^{-1}
 \nbigo_{\cnum\times X_{\Sigma_0,S}}
 [!Z_{f_{\gamma}}][\ast D_{\Sigma_0,S}]
 \Bigr)$.
Condition {\bf(A)} for $\pi_{\Sigma_0}$
and $\nbigu_{X_{\Sigma_0,S}}(n+d_S,-1)
 [!Z_{f_{\gamma}}][\ast D_{\Sigma_0,S}]$
with the graded polarization
implies that the induced morphism
\begin{equation}
\label{eq;14.12.13.30}
 \Ker\beta_{k}\cap
 \Ker\bigl(
 \alpha_{k}^{\lor}\circ\Phi_{\nbigptilde_{n+k}}
 \bigr)
\lrarr
 \Gr^{\Wtilde}_{n+k}
 \nbigv'_{\ast}(\nbigb,S,\gamma)
\end{equation}
is an epimorphism.
The following proposition holds
by the construction of $\Ptilde_{n+k}$.

\begin{prop}
\label{prop;14.12.13.31}
The pull back of $\Ptilde_{n+k}$ by {\rm(\ref{eq;14.12.13.30})}
is equal to the restriction of $\nbigptilde_{n+k}$ to 
$\Ker\beta_{k}\cap
 \Ker\bigl(
 \alpha_{k}^{\lor}\circ\Phi_{\nbigptilde_{n+k}}
 \bigr)$.

\hfill\qed
\end{prop}

In the sense of Proposition \ref{prop;14.12.13.31},
a description of 
the complex (\ref{eq;14.12.13.40}) and the pairing $\nbigptilde_{n+k}$
gives a description of $\Ptilde_{n+k}$.

\paragraph{A Hodge theoretic description of the pairings}
We shall give a Hodge theoretic description of
(\ref{eq;14.12.13.40}) and $\nbigptilde_{n+k}$
in the case where $S$ is a point.
According to Proposition \ref{prop;14.11.18.41},
it is enough to consider the case where $S$ is a point.

To simplify the description,
we denote $f_{\gamma}$ by $f$,
and the zero-set is denoted by $Z_f$.
We have the irreducible decomposition
$D_{\Sigma_0}=\bigcup_{i\in \Lambda} D_{\Sigma_0,i}$.
For $I\subset \Lambda$,
we set $D_{\Sigma_0,I}:=\bigcap_{i\in I}D_{\Sigma_0,i}$.
Formally, $D_{\Sigma_0,\emptyset}:=X_{\Sigma_0}$.
Let $\iota_{I}:D_{\Sigma_0,I}\lrarr X_{\Sigma_0}$
be the inclusion.
Let $Z_{f,I}:=Z_{f}\cap D_{\Sigma_0,I}$
which are also smooth.
Let $\iota_{Z_{f},I}$ denote the inclusion
$Z_{f,I}\lrarr X_{\Sigma_0}$.

We set $\Lambda(k):=
 \bigl\{I\subset\Lambda\,\big|\,|I|=\Lambda\bigr\}$.
Let $W$ be the filtration on 
$\lambda^{-1}
 \nbigo_{\cnum\times X_{\Sigma_0}}
[!Z_f][\ast D_{\Sigma_0,I}]$
underlying the weight filtration of
the mixed Hodge module
$\nbigu_{X_{\Sigma_0}}(n,-1)[!Z_{f}][\ast D_{\Sigma_0}]$.
We have the following natural isomorphism:
\[
 \Gr^W_{n+k}\Bigl(
 \lambda^{-1}
 \nbigo_{\cnum\times X_{\Sigma_0}}
 [!Z_{f}][\ast D_{\Sigma_0}]
 \Bigr)
\simeq
 \bigoplus_{I\in\Lambda(k-1)}
 \iota_{I\dagger}
 \bigl(
 \lambda^{-k}
 \nbigo_{\cnum\times D_{\Sigma_0,I}}
 \bigr)
\oplus
  \bigoplus_{I\in\Lambda(k)}
 \iota_{Z_{f},I\dagger}
 \bigl(
 \lambda^{-k-1}
 \nbigo_{\cnum\times Z_{f,I}}
 \bigr)
\]

Let $F^{\bullet}$ be the Hodge filtration on $H^{\ast}(D_{\Sigma_0,I})$
which is a decreasing filtration.
We set 
$F_{j}H^m(D_{\Sigma_0,I}):=
 F^{-j}H^m(D_{\Sigma_0,I})$.
Let $R_FH^m(D_{\Sigma_0,I})$ be the Rees module:
\[
 R_FH^m(D_{\Sigma_0,I})
=\sum F_jH^m(D_{\Sigma_0,I})\lambda^{j}
=\sum F^jH^m(D_{\Sigma_0,I})\lambda^{-j}
\]
The associated $\nbigo_{\cnum_{\lambda}}$-module
is also denoted by the same notation.
We have the following natural isomorphism:
\begin{multline}
 \pi^i_{\Sigma_0\dagger}
 \Gr^W_{n+k}\Bigl(
 \lambda^{-1}
 \nbigo_{\cnum\times X_{\Sigma_0}}
 [!Z_{f}][\ast D_{\Sigma_0}]
 \Bigr)
\simeq
 \\
 \bigoplus_{I\in\Lambda(k-1)} 
 \lambda^{-k}
 R_F H^{n-k+i}(D_{\Sigma_0,I})
 \oplus
 \bigoplus_{I\in\Lambda(k)}
 \lambda^{-k-1}
 R_FH^{n-2-k+i}(Z_{f,I})
\end{multline}
We obtain the complex (\ref{eq;14.12.13.40})
by applying the Rees construction
to the following exact sequence of
pure Hodge structures:
\begin{multline}
\label{eq;14.12.13.50}
\bigoplus_{I\in\Lambda(k)}
 H^{n-k-2}(D_{\Sigma_0,I})\otimes\rnum(-k-1)
\oplus
\bigoplus_{I\in\Lambda(k+1)}
 H^{n-4-k}(Z_{f,I})\otimes\rnum(-k-2)
\lrarr \\
\bigoplus_{I\in\Lambda(k-1)}
 H^{n-k}(D_{\Sigma_0,I})\otimes\rnum(-k)
\oplus
\bigoplus_{I\in\Lambda(k)}
 H^{n-2-k}(Z_{f,I})\otimes\rnum(-k-1)
\lrarr \\
\bigoplus_{I\in\Lambda(k-2)}
 H^{n-k+2}(D_{\Sigma_0,I})\otimes\rnum(-k+1)
\oplus
\bigoplus_{I\in\Lambda(k-1)}
 H^{n-k}(Z_{f,I})\otimes\rnum(-k)
\end{multline}
Here, $\rnum(i)$ denote the $i$-th Tate object
in the category of Hodge structures.
We can also obtain the complex
(\ref{eq;14.12.13.50})
by using the theory of mixed Hodge modules
\cite{saito2},
or more explicitly by using the theory of mixed Hodge complexes
(see \cite{peters-steenbrink}).

\vspace{.1in}

We have the following pairings on
$\nbigp_I$ on $H^{n-1-|I|}(D_{\Sigma_0,I})$
and 
$\nbigp_{Z_f,I}$ on $H^{n-|I|-2}(Z_{f,I})$:
\[
\nbigp_I(a,b)=
 \left(\frac{1}{2\pi\sqrt{-1}}\right)^{n-1-|I|}
 \int_{D_{\Sigma_0,I}}ab,
\quad\quad
\nbigp_{Z_f,I}(a,b)=
  \left(\frac{1}{2\pi\sqrt{-1}}\right)^{n-2-|I|}
\int_{Z_{f,I}}ab.
\]
For $i\in\Lambda$,
let $p_i(f)$ be the pole order of $f$ along $D_i$.
For any integer $\ell$,
let $\epsilon(\ell):=(-1)^{\ell(\ell-1)/2}$.
We set 
\[
 \nbigptilde_I:=
 \epsilon(n-|I|-1)\cdot \prod_{i\in I}p_i(f)^{-1}
 \cdot 
 \nbigp_I,
\quad\quad
 \nbigptilde_{Z_f,I}:=
\epsilon(n-2-|I|)\cdot
 \prod_{i\in I}p_i(f)^{-1}\cdot \nbigp_{Z_f,I}.
\]
We have the pairing of weight $n-|I|-1$
on $R_FH^{n-|I|-1}(D_{\Sigma_0,I})$
induced by $\nbigptilde_{I}$,
and the pairing of weight $n-|I|-2$
on $R_FH^{n-|I|-2}(Z_{f,I})$
induced by $\nbigptilde_{Z_f,I}$.
The induced pairings are also denoted by
$\nbigptilde_I$
and 
$\nbigptilde_{Z_f,I}$
respectively.

\begin{prop}
We have
$\nbigptilde_{n+k}
=\bigoplus_{I\in\Lambda(k-1)} \nbigptilde_I
\oplus
 \bigoplus_{I\in\Lambda(k)}\nbigptilde_{Z_f,I}$.
\end{prop}
\pf
It follows from
Proposition \ref{prop;14.9.28.1}
and Proposition \ref{prop;14.11.18.41}.
\hfill\qed

\subsection{Mixed TEP-structures on
reduced GKZ-hypergeometric systems}
\label{subsection;14.11.26.200}

\subsubsection{Preliminary}

As in \S\ref{subsection;14.7.8.1},
we consider 
$\nbiga=\{\veca_1,\ldots,\veca_m\}\subset\seisuu^n$
which generates $\seisuu^n$.
We set $N_{\nbiga}:=\seisuu^n$
and $M_{\nbiga}:=\seisuu^m$.
Let $e_1,\ldots,e_m$ be the standard basis of $M_{\nbiga}$.
We have the surjective morphism
$\Xi_{\nbiga}:M_{\nbiga}\lrarr N_{\nbiga}$.
Let $L_{\nbiga}:=\Ker\Xi_{\nbiga}$.
We obtain the exact sequence
\[
\begin{CD}
 0 @>>>
 L_{\nbiga}
 @>{\Theta_{\nbiga}}>>
 M_{\nbiga}
 @>{\Xi_{\nbiga}}>>
 N_{\nbiga}
 @>>> 0.
\end{CD}
\]
By taking the dual,
we obtain the exact sequence
$0\lrarr 
 N_{\nbiga}^{\lor}
 \stackrel{\Xi_{\nbiga}^{\lor}}\lrarr
 M_{\nbiga}^{\lor}
 \stackrel{\Theta_{\nbiga}^{\lor}}{\lrarr}
 L_{\nbiga}^{\lor}\lrarr 0$.
For a finitely generated free abelian group $A$,
we set 
$A_{\cnum^{\ast}}:=
 \cnum^{\ast}\otimes_{\seisuu} A$.
We can naturally regard it as a complex algebraic variety
or a complex manifold.
Particularly,
we set 
$S_{\nbiga}:=L^{\lor}_{\nbiga,\cnum^{\ast}}$.
We have the natural surjection
$\Theta^{\lor}_{\nbiga}:
 M^{\lor}_{\nbiga,\cnum^{\ast}}
 \lrarr
 S_{\nbiga}$.

We have the action $\rho_1$ of
$N^{\lor}_{\nbiga,\cnum^{\ast}}$
on $M^{\lor}_{\nbiga,\cnum^{\ast}}$
induced by $-\varphi$.
We also have 
the natural action $\rho_0$
of $N^{\lor}_{\nbiga,\cnum^{\ast}}$
on itself by the multiplication.
We can naturally regard
$\Theta^{\lor}_{\nbiga}:
 M^{\lor}_{\nbiga,\cnum^{\ast}}
 \lrarr
 S_{\nbiga}$
as the quotient of 
the projection
$N^{\lor}_{\nbiga,\cnum^{\ast}}\times
 M^{\lor}_{\nbiga,\cnum^{\ast}}
\lrarr
 M^{\lor}_{\nbiga,\cnum^{\ast}}$
via the above actions of $N^{\lor}_{\nbiga,\cnum^{\ast}}$.

The identifications
$M_{\nbiga}=\seisuu^m$
and $N_{\nbiga}=\seisuu^n$
induce the coordinate systems
$(w_1,\ldots,w_m)$ on $M_{\nbiga,\cnum^{\ast}}^{\lor}$
and 
$(t_1,\ldots,t_n)$ on $N^{\lor}_{\nbiga,\cnum^{\ast}}$.
We set $G_{\nbiga}:=\sum_{i=1}^m w_i$
on $M^{\lor}_{\nbiga,\cnum^{\ast}}$.
The algebraic function 
$F_{\nbiga}=\sum_{i=1}^m w_i t^{\veca_i}$
on
$N^{\lor}_{\nbiga,\cnum^{\ast}}
 \times
 M^{\lor}_{\nbiga,\cnum^{\ast}}$
is 
$N^{\lor}_{\nbiga,\cnum^{\ast}}$-invariant,
and $G_{\nbiga}$ is the descent of $F_{\nbiga}$.
Any splitting 
$\gamma_{\nbiga}:
 L_{\nbiga}^{\lor}\lrarr M_{\nbiga}^{\lor}$
of $\Theta_{\nbiga}^{\lor}$
induces an algebraic morphism
$\gamma_{\nbiga}:
 L_{\nbiga,\cnum^{\ast}}^{\lor}
\lrarr 
 M^{\lor}_{\nbiga,\cnum^{\ast}}$,
and we obtain
$\id\times\gamma_{\nbiga}:
 N^{\lor}_{\nbiga,\cnum^{\ast}}
 \times
 L^{\lor}_{\nbiga,\cnum^{\ast}}
\lrarr
  N^{\lor}_{\nbiga,\cnum^{\ast}}
 \times
 M^{\lor}_{\nbiga,\cnum^{\ast}}$.
The splitting $\gamma_{\nbiga}$ also 
gives an isomorphism
$M^{\lor}_{\nbiga}
\simeq
 N^{\lor}_{\nbiga}\times
 L^{\lor}_{\nbiga}$,
and hence
$M^{\lor}_{\nbiga,\cnum^{\ast}}
\simeq
 N^{\lor}_{\nbiga,\cnum^{\ast}}
 \times
 L^{\lor}_{\nbiga,\cnum^{\ast}}$.
The pull back of $F_{\nbiga}$
by $\id\times\gamma_{\nbiga}$
is equal to $G_{\nbiga}$
under the identification.

\subsubsection{Mixed twistor $\nbigd$-modules
and the induced mixed TEP-structures}
\label{subsection;14.12.1.10}

We naturally regard
$M_{\nbiga,\cnum^{\ast}}^{\lor}
=(\cnum^{\ast})^m
=\bigl\{
 \sum_{i=1}^m\alpha_iz_i\,\big|\,
 \alpha_i\neq 0
 \bigr\}$,
and 
$N^{\lor}_{\nbiga,\cnum^{\ast}}=T^n$
in \S\ref{subsection;14.7.8.1}.
Let $\gamma_{\nbiga}:
 L_{\nbiga,\cnum^{\ast}}^{\lor}
 \lrarr
 M_{\nbiga,\cnum^{\ast}}^{\lor}$
be the splitting as above.
Note that it is non-degenerate at $\infty$
for $X_{\nbiga}$,
as remarked in Example \ref{example;14.11.24.101}.
We obtain
the following integrable mixed twistor $\nbigd$-modules
with real structure on $S_{\nbiga}$:
\[
 \nbigt_{\nbiga,\star}:=
 \pi_{\dagger}^0
 \nbigt_{\star}\bigl(\nbiga,
 S_{\nbiga},\gamma_{\nbiga}\bigr)
\quad
 (\star=\ast,!).
\]
The underlying $\nbigd$-modules are 
the reduced GKZ-hypergeometric systems.
As in Lemma \ref{lem;14.12.30.10},
$\nbigt_{\nbiga,\star}$ are independent 
of the choice of a splitting $\gamma_{\nbiga}$.

We set 
$\nbigt_{\nbiga,\min}
:=\pi_{\dagger}^0
 \nbigt_{\min}(\nbiga,S_{\nbiga},\gamma_{\nbiga})$.
By Corollary \ref{cor;14.11.29.400},
it is isomorphic to the image of the natural morphism
$\nbigt_{\nbiga,!}
\lrarr
 \nbigt_{\nbiga,\ast}$,
and 
$\Gr^W_{m}\nbigt_{\nbiga,\ast}
=\Gr^W_m\nbigt_{\nbiga,!}
=\Gr^W_m\nbigt_{\nbiga,\min}
=\nbigt_{\nbiga,\min}$.
If $0$ is contained in the interior part of 
$\Conv(\nbiga)$,
we have
$\nbigt_{\nbiga!}
=\nbigt_{\nbiga,\ast}
=\nbigt_{\nbiga,\min}$.

We set 
$\nbigl_{\nbiga,\star}:=\nbigl_{\star}\bigl(\nbiga,
 S_{\nbiga},\gamma_{\nbiga}\bigr)$
and 
$\nbigv_{\nbiga,\star}:=
 \pi_{\dagger}^0\nbigl_{\nbiga,\star}$.
The $\nbigrtilde$-modules
$\nbigv_{\nbiga,\star}$ underlie
$\nbigt_{\nbiga,\star}$.
More precisely,
$\nbigt_{\nbiga,\star}$
is expressed as 
a pair of $\nbigr$-modules
$\lambda^{m}\nbigv_{\nbiga,!}$
and 
$\nbigv_{\nbiga,\ast}$
with the induced sesqui-linear pairing,
and 
$\nbigt_{\nbiga,!}$
is expressed as 
a pair of $\nbigr$-modules
$\lambda^{m}\nbigv_{\nbiga,\ast}$
and 
$\nbigv_{\nbiga,!}$
with the induced sesqui-linear pairing.
Let $W$ denote the filtration of $\nbigv_{\nbiga,\star}$
underlying the weight filtration of $\nbigt_{\nbiga,\star}$.
We set
$\Wtilde_{k}\nbigv_{\nbiga,\star}:=
 W_{k+n-m}\nbigv_{\nbiga,\star}$.
Note $\dim S_{\nbiga}=m-n$.

Let $\nbigl_{\nbiga,\min}$ be the image of
$\nbigl_{\nbiga,!}\lrarr\nbigl_{\nbiga,\ast}$,
and set
$\nbigv_{\nbiga,\min}:=
 \pi_{\dagger}^0\nbigl_{\nbiga,\min}$.
It is naturally isomorphic to the image of
$\nbigv_{\nbiga,!}\lrarr\nbigv_{\nbiga,\ast}$ 
by Corollary \ref{cor;14.11.29.400}.
The $\nbigrtilde$-module 
$\nbigv_{\nbiga,\min}$
underlies $\nbigt_{\nbiga\min}$.
We have
$\Gr^{\Wtilde}_n\nbigv_{\nbiga,!}
\simeq
 \Gr^{\Wtilde}_n\nbigv_{\nbiga,\ast}
\simeq
 \nbigv_{\nbiga,\min}$.

If $0$ is an interior point of $\Conv(\nbiga)$,
then $\nbigt_{\nbiga\ast}=\nbigt_{\nbiga,\min}=\nbigt_{\nbiga,!}$
is pure of weight $m$,
and it is equipped with the canonical polarization.
As explained in \S\ref{subsection;14.11.24.121},
even if $0$ is not an interior point of $\Conv(\nbiga)$,
we have the graded polarizations 
$\vecnbigs_{\nbiga,\star,\vecb}$
of $\nbigt_{\nbiga,\star}$
depending on the choice of $\vecb\in\seisuu^n$
such that $0$ is an interior point of $\Conv(\nbiga\cup\{\vecb\})$.
The weight $m$-part of $\vecnbigs_{\nbiga,\star,\vecb}$
are independent of $\vecb$.

Let $S_{\nbiga}^{\reg}:=\gamma_{\nbiga}^{-1}(U_{\nbiga}^{\reg})$.
The restriction 
$\nbigv_{\nbiga,\star|S_{\nbiga}^{\reg}}$ are locally free
$\nbigo_{\cnum_{\lambda}\times S_{\nbiga}^{\reg}}$-modules.
Take a hypersurface $Y$ such that
$Y\supset S_{\nbiga}\setminus S_{\nbiga}^{\reg}$.
The real structure and the graded polarization of
$\nbigt_{\nbiga,\star}$ induce
graded pairing $P_{\nbiga,\star,\vecb}$ of
$(\nbigv_{\nbiga,\star}(\ast Y),\Wtilde)$.
In this way,
we obtain mixed TEP-structures
$(\nbigv_{\nbiga,\star}(\ast Y),\Wtilde,
 P_{\nbiga,\star,\vecb})$.

\begin{rem}
Let $\cnum[M_{\nbiga}]$ denote the group ring 
of $M_{\nbiga}$ over $\cnum$.
We may naturally regard
$M^{\lor}_{\nbiga,\cnum^{\ast}}$
as the algebraic variety
$\Spec \cnum[M_{\nbiga}]$.
Then, $\nbigt_{\nbiga,\star}$ is algebraic
in the sense of {\rm\S\ref{subsection;14.11.24.120}}.
\hfill\qed
\end{rem}

\begin{rem}
If $0$ is contained in the interior part of $\Conv(\nbiga)$,
we shall often omit the subscripts $\ast,!$
because 
$\nbigt_{\nbiga!}=\nbigt_{\nbiga\ast}=\nbigt_{\nbiga,\min}$.
\hfill\qed
\end{rem}

\subsubsection{Homogeneity}
\label{subsection;14.11.26.201}

Let $e_1^{\lor},\ldots,e_m^{\lor}$ denote the dual basis
of $M_{\nbiga}^{\lor}$.
We have the morphism
$\seisuu\lrarr M_{\nbiga}^{\lor}$
given by
$1\longmapsto
 \gminiv:=\sum_{i=1}^m e_i^{\lor}$.
It induces a $\cnum^{\ast}$-action
on $M^{\lor}_{\nbiga,\cnum^{\ast}}$.
By the composition
$\seisuu\lrarr
M^{\lor}_{\nbiga}
\lrarr
 L^{\lor}_{\nbiga}$,
we obtain a $\cnum^{\ast}$-action
on $S_{\nbiga}$.
The map
$\Theta_{\nbiga}^{\lor}:
 M^{\lor}_{\nbiga,\cnum^{\ast}}\lrarr S_{\nbiga}$ 
is $\cnum^{\ast}$-equivariant.
For the $\cnum^{\ast}$-action,
we have $a^{\ast}G_{\nbiga}=a\cdot G_{\nbiga}$
for any $a\in \cnum^{\ast}$.

We consider the action of 
$\cnum^{\ast}$ on $\cnum_{\lambda}$
given by the multiplication.
For the diagonal $\cnum^{\ast}$-action on
 $\cnum_{\lambda}\times M^{\lor}_{\nbiga,\cnum^{\ast}}$,
we have
$a^{\ast}(\lambda^{-1}G_{\nbiga})
 =\lambda^{-1}G_{\nbiga}$ for any $a\in\cnum^{\ast}$.

\begin{lem}
The $\nbigrtilde$-modules
$\nbigv_{\nbiga,\star}$ on $S_{\nbiga}$
are $\cnum^{\ast}$-homogeneous 
in the sense of {\rm\S\ref{subsection;14.11.25.10}}.
\end{lem}
\pf
Take any toric desingularization
$\varphi_{\Sigma}:X_{\Sigma}\lrarr X_{\nbiga}$.
We have the $\cnum^{\ast}$-action on
$X_{\Sigma}\times S_{\nbiga}$
which is the extension of
$T^n\times S_{\nbiga}\simeq
 M^{\lor}_{\nbiga,\cnum^{\ast}}$.
For the action,
we have 
$a^{\ast}(\lambda^{-1}F_{\gamma_{\nbiga},\Sigma})
=\lambda^{-1}F_{\gamma_{\nbiga},\Sigma}$.

We set $Y:=X_{\Sigma}\times S_{\nbiga}$
and $D_Y:=D_{\Sigma}\times S_{\nbiga}$.
By Proposition \ref{prop;14.11.25.20},
the $\nbigrtilde$-modules
$\nbigl_{\star}(F_{\gamma_{\nbiga},\Sigma},D_{Y})$
are $\cnum^{\ast}$-homogeneous.
Then, 
$\nbigv_{\nbiga,\star}
\simeq
 \pi^0_{\Sigma\dagger}
 \nbigl_{\star}(F_{\gamma_{\nbiga},\Sigma},D_Y)$
are also $\cnum^{\ast}$-homogeneous.
It is easy to check that the $\cnum^{\ast}$-actions
are independent of the choice of $\varphi_{\Sigma}$.
\hfill\qed

\subsubsection{Variation of Hodge structure}

Let us consider the special case
$\Theta_{\nbiga}^{\lor}(\gminiv)=0$,
i.e., the natural $\cnum^{\ast}$-action on $S_{\nbiga}$ is trivial.
By the general result in Theorem \ref{thm;15.11.17.20},
we have the following.
\begin{prop}
If $\Theta_{\nbiga}^{\lor}(\gminiv)=0$,
then the integrable mixed twistor $\nbigd$-modules
$\nbigt_{\nbiga,\star}$ come from mixed $\real$-Hodge modules.
\hfill\qed
\end{prop}

\begin{rem}
In particular, the underlying $\nbigd$-modules
of $\nbigt_{\nbiga,\star}$ are regular singular.
Note that the condition $\Theta^{\lor}_{\nbiga}(\gminiv)=0$
is equivalent to the standard criterion
for the GKZ-system to be regular singular.
Indeed, we have $\Theta_{\nbiga}^{\lor}(\gminiv)=0$
if and only if
there exists $\alpha\in N_{\nbiga}^{\lor}$
such that $\alpha(\Xi_{\nbiga}(e_i))=1$
$(i=1,\ldots,m)$.
If $\Theta_{\nbiga}^{\lor}(\gminiv)=0$,
we have $\alpha\in N_{\nbiga}^{\lor}$
such that
$\Xi_{\nbiga}^{\lor}(\alpha)(\gminiv)$,
and hence we have
$\langle
 \alpha,\Xi_{\nbiga}(e_i)
 \rangle
=\langle
 \Xi_{\nbiga}^{\lor}(\alpha),e_i
 \rangle
=\langle\gminiv,e_i\rangle
=1$.
Conversely, suppose that
there exists $\alpha\in N_{\nbiga}^{\lor}$
such that
$\alpha(\Xi_{\nbiga}(e_i))=1$.
Because 
$\langle
 \Xi_{\nbiga}^{\lor}(\alpha),e_i
 \rangle=1$,
we have
$\Xi_{\nbiga}^{\lor}(\alpha)=\gminiv$.
\hfill\qed
\end{rem}

\begin{cor}
$\nbigt_{\nbiga,\star|S^{\reg}_{\nbiga}}$
comes from 
a graded polarizable variation of mixed Hodge structure
on $S^{\reg}_{\nbiga}$.
For any algebraic embedding
$S^{\reg}_{\nbiga}\subset Z$,
the graded polarizable variation of mixed Hodge structure
is admissible along $Z\setminus S^{\reg}_{\nbiga}$.
\end{cor}
\pf
It follows from the algebraicity 
in Proposition \ref{prop;14.11.25.40}
and a general property of mixed twistor $\nbigd$-modules.

\hfill\qed

\subsubsection{Appendix: Comparison with the construction of Reichelt-Sevenheck}
\label{subsection;14.12.3.101}

As in \S\ref{subsection;14.7.9.20},
we set $V:=H^0(\proj^m,\nbigo(1))$.
Let $Z\subset\proj^m\times V$
be the $0$-set of the universal section of 
$\nbigo_{\proj^m}(1)\boxtimes\nbigo_V$.
We decompose $V=V_1\times V_2$,
where 
$V_1:=\bigl\{\alpha_0z_0\,\big|\,\alpha_0\in\cnum\bigr\}$
and 
$V_2:=\bigl\{\sum_{i=1}^m\alpha_iz_i\,\big|\,\alpha_i\in\cnum\bigr\}$.
We identify $V_2=M_{\nbiga}^{\lor}\otimes\cnum$.
We have the splitting 
$\gamma_{\nbiga}:
 S_{\nbiga}=L_{\nbiga}^{\lor}\otimes\cnum^{\ast}
 \lrarr
 M_{\nbiga}^{\lor}\otimes\cnum^{\ast}
\subset
 V_2$.
We obtain 
$\id\times\gamma_{\nbiga}:
 V_1\times S_{\nbiga}
\lrarr
 V_1\times V_2$.
Let $Z_{\nbiga}$ be the fiber product of
$Z$ and $V_1\times S_{\nbiga}$
over $V_1\times V_2$.
We have the naturally induced morphisms
 $q_1:Z_{\nbiga}\lrarr \proj^m$
and $q_2:Z_{\nbiga}\lrarr V_1\times S_{\nbiga}$.

Recall that we have the morphism
$\psi_{\nbiga}:T^n\lrarr \proj^m$
induced by $\nbiga$.
Let $U:=\psi_{\nbiga}(T^n)$.
We set $Z_{\nbiga,U}:=Z_{\nbiga}\times_{\proj^m}U$.
Let $\iota_{Z_{\nbiga,U}}:
 Z_{\nbiga,U}\lrarr \proj^m\times (V_1\times S_{\nbiga})$
be the inclusion.

Let us consider the pure Hodge module
$\bigl(
 \nbigo_{Z_{\nbiga,U}},F
 \bigr)$
and the mixed Hodge modules
$\iota_{Z_{\nbiga,U}\star}\bigl(
 \nbigo_{Z_{\nbiga,U}},F
 \bigr)$ $(\star=\ast,!)$
as in \S\ref{subsection;14.12.1.2}.
Let $\nbigm^{IC}(Z_{\nbiga,U})$
denote the image of the morphism
$\iota_{Z_{\nbiga,U}!}\bigl(
 \nbigo_{Z_{\nbiga,U}},F
 \bigr)
\lrarr
\iota_{Z_{\nbiga,U}\ast}\bigl(
 \nbigo_{Z_{\nbiga,U}},F
 \bigr)$.
Let $\pi_{V_1\times S_{\nbiga}}:
 \proj^m\times (V_1\times S_{\nbiga})
\lrarr V_1\times S_{\nbiga}$
be the projection.
Then, we obtain the mixed Hodge modules
$\pi^0_{V_1\times S_{\nbiga}\dagger}
 \iota_{Z_{\nbiga,U}\star}\bigl(
 \nbigo_{Z_{\nbiga,U}},F
 \bigr)$ $(\star=!,\ast)$
and 
$\pi^0_{V_1\times S_{\nbiga}\dagger}
\nbigm^{IC}(Z_{\nbiga,U})$.
By applying the procedure in \S\ref{subsection;14.7.6.1},
we obtain the following 
$\nbigrtilde_{S_{\nbiga}}$-modules:
\[
 G_0\FL^{\loc}_{S_{\nbiga}}
 \bigl(\pi^0_{V_1\times S_{\nbiga}\dagger}
 \iota_{Z_{\nbiga,U}\star}
 (\nbigo_{Z_{\nbiga,U}},F)\bigr),
\quad\quad
 G_0\FL^{\loc}_{S_{\nbiga}}
 \bigl(\pi^0_{V_1\times S_{\nbiga}\dagger}
 \nbigm^{IC}(Z_{\nbiga,U})\bigr).
\]

\begin{prop}
\label{prop;14.12.15.1}
We have isomorphisms of $\nbigrtilde_X$-modules:
\[
 \lambda \cdot \nbigv_{\nbiga\star}
\simeq
 G_0\FL^{\loc}_{S_{\nbiga}}
 \bigl(\pi^0_{V_1\times S_{\nbiga}\dagger}
 \iota_{Z_{\nbiga,U}\star}
 (\nbigo_{Z_{\nbiga,U}},F)
 \bigr)
\]
We also have the following commutative diagram:
\[
 \begin{CD}
 \lambda\cdot\nbigv_{\nbiga!}
 @>{\simeq}>>
 G_0\FL^{\loc}_{S_{\nbiga}}
 \bigl(\pi^0_{V_1\times S_{\nbiga}\dagger}
 \iota_{Z_{\nbiga,U}!}
 (\nbigo_{Z_{\nbiga,U}},F)\bigr)
 \\
 @VVV @VVV \\
 \lambda\cdot\nbigv_{\nbiga\ast}
 @>{\simeq}>>
 G_0\FL^{\loc}_{S_{\nbiga}}
 \bigl(\pi^0_{V_1\times S_{\nbiga}\dagger}
 \iota_{Z_{\nbiga,U}\ast}
 (\nbigo_{Z_{\nbiga,U}},F)\bigr)
 \end{CD}
\]
\end{prop}
\pf
We can obtain the claim 
from Proposition \ref{prop;14.11.12.3}
by using the non-characteristic pull back.
We can also prove it directly
by the argument in the proof of
Proposition \ref{prop;14.11.12.3}.
\hfill\qed

\begin{cor}
\label{cor;14.12.3.102}
We naturally have 
$\lambda\cdot\nbigv_{\nbiga,\min}
\simeq
 G_0\FL^{\loc}_{S_{\nbiga}}
 \bigl(\pi^0_{V_1\times S_{\nbiga}\dagger}
 \nbigm^{IC}(Z_{\nbiga,U})\bigr)$.
\end{cor}
\pf
Reichelt-Sevenheck proved that 
$G_0\FL^{\loc}_{S_{\nbiga}}
 \bigl(\pi^0_{V_1\times S_{\nbiga}\dagger}
 \nbigm^{IC}(Z_{\nbiga,U})\bigr)$
is naturally isomorphic to the image of
$G_0\FL^{\loc}_{S_{\nbiga}}
 \bigl(\pi^0_{V_1\times S_{\nbiga}\dagger}
 \iota_{Z_{\nbiga,U}!}
 (\nbigo_{Z_{\nbiga,U}},F)\bigr)
\lrarr
 G_0\FL^{\loc}_{S_{\nbiga}}
 \bigl(\pi^0_{V_1\times S_{\nbiga}\dagger}
 \iota_{Z_{\nbiga,U}\ast}
 (\nbigo_{Z_{\nbiga,U}},F)\bigr)$.
As mentioned in \S\ref{subsection;14.12.1.10},
$\nbigv_{\nbiga\min}$ is naturally isomorphic
to the image of 
$\nbigv_{\nbiga!}\lrarr\nbigv_{\nbiga\ast}$
by Corollary \ref{cor;14.11.29.400}.
Hence, we have the desired isomorphism.

\hfill\qed

\section{Quantum $\nbigd$-modules}
\subsection{Some mixed twistor $\nbigd$-modules
 in mirror symmetry}
\label{subsection;15.1.6.1}

\subsubsection{Toric data}
\label{subsection;14.11.27.30}

Let $X$ be an $n$-dimensional smooth weak Fano toric variety.
Let $\Sigma$ be a fan of $X$.
Let $\Sigma(1)=\{\rho_1,\ldots,\rho_m\}$
denote the set of the $1$-dimensional cones in $\Sigma$.
Let $[\rho_i]\in\seisuu^n$
be the primitive generator of
$\rho_i\cap\seisuu^n$.

Let $K_X$ denote the canonical bundle of $X$.
Let $\nbigl_j$ $(j=1,\ldots,r)$ be nef line bundles on $X$
such that
$\bigl(
 K_X\otimes \bigotimes_{j=1}^{r}\nbigl_j
 \bigr)^{\lor}$
is nef.
We may assume that
$\nbigl_j=\nbigo(\sum_{i=1}^m\beta_{ji}D_i)$
$(j=1,\ldots,r)$
for some $\beta_{ji}\in\seisuu_{\geq 0}$,
where $D_i$ are the hypersurfaces corresponding to 
the cones $\rho_i$.
Let $\seisuu^{n+r}=\seisuu^n\oplus\seisuu^r$.
Let $n_1,\ldots,n_r$ denote the standard basis of $\seisuu^r$.
We set $\veca_i\in\seisuu^n\oplus\seisuu^r$
as follows:
\[
 \veca_i:=
 \left\{
 \begin{array}{ll}
 [\rho_i]+\sum_{j=1}^r\beta_{ji}n_j
 & (i=1,\ldots,m),\\
 n_{i-m} & (i=m+1,\ldots,m+r).
 \end{array}
 \right.
\]
We put $\nbiga:=\{\veca_1,\ldots,\veca_{r+m}\}$.
We also set $\veca_{r+m+1}:=-\sum_{j=1}^r n_j$
and $\nbigatilde:=\nbiga\cup\{\veca_{r+m+1}\}$.
Note that
$0$ is contained in the interior part of
the convex hull of $\Conv(\nbigatilde)$.
(We shall use $\nbigatilde$ in \S\ref{subsection;14.12.29.1}.)

\subsubsection{$\nbigrtilde$-modules
 with filtration and graded pairing
 associated to $\nbiga$}
\label{subsection;14.12.2.30}

Applying the construction in \S\ref{subsection;14.11.26.200},
we obtain the integrable mixed twistor $\nbigd$-modules
$\nbigt_{\nbiga,\star}$ $(\star=\ast,!)$
with real structure on $S_{\nbiga}$.
They are equipped with the graded polarization
$\vecnbigs_{\nbiga,\star,\veca_{r+m+1}}$.
Let $\nbigv_{\nbiga,\star}$ denote their underlying
$\nbigrtilde$-modules.
They are $\cnum^{\ast}$-homogeneous
as in \S\ref{subsection;14.11.26.201}.
By the procedure in 
\S\ref{subsection;14.11.8.120}--\ref{subsection;14.11.28.10},
we obtain the filtration $\Wtilde$
and the graded pairing 
$(P_{\nbiga,\star,k}\,|\,k\in\seisuu)$
on $\nbigv_{\nbiga,\star|S_{\nbiga}^{\reg}}$.

\vspace{.1in}

We apply the same construction to
$\nbigt_{\nbiga,\star}\otimes\newTate(n+r)$.
Let $\gbigv_{\nbiga,\star}$ be the underlying
$\nbigrtilde$-modules.
They are $\cnum^{\ast}$-homogeneous.
The restriction $\gbigv_{\nbiga,\star|S_{\nbiga}^{\reg}}$
is equipped with the filtration $\Wtilde$
and the graded pairings 
$\gbigp_{\nbiga,\star}=(\gbigp_{\nbiga,\star,k}\,\big|\,k\in\seisuu)$.
Thus, we obtain mixed TEP-structures
$\bigl(
 \gbigv_{\nbiga,\star|S_{\nbiga}^{\reg}},
 \Wtilde,
 \gbigp_{\nbiga,\star}
 \bigr)$
on $S_{\nbiga}^{\reg}$.

By construction,
we have
$\gbigv_{\nbiga,\star}=
 \lambda^{n+r}\nbigv_{\nbiga,\star}$
and 
$\Wtilde_k\gbigv_{\nbiga,\star}=
 \lambda^{n+r}
 \Wtilde_{k+2(n+r)}\nbigv_{\nbiga,\star}$.
We have
$\Gr^{\Wtilde}_k\gbigv_{\nbiga,\star}
=\lambda^{n+r}
 \Gr^{\Wtilde}_{k+2(n+r)}\nbigv_{\nbiga,\star}$.
The graded pairings on 
$\Gr^{\Wtilde}_k\gbigv_{\nbiga,\star}(\ast\lambda)
=\Gr^{\Wtilde}_{k+2(n+r)}\nbigv_{\nbiga,\star}(\ast\lambda)$
are equal.

\subsubsection{Expression as the systems of differential equations}

Let $[z_0:\cdots:z_{m+r}]$ be 
the homogeneous coordinate system on $\proj^{m+r}$.
We set  $W:=
 \bigl\{\sum_{i=1}^{m+r}\alpha_iz_i \bigr\}
\subset
 H^0\bigl(\proj^{m+r},\nbigo(1)\bigr)
=M^{\lor}_{\nbiga}\otimes\cnum$.
The splitting 
$S_{\nbiga}=L^{\lor}_{\nbiga}\otimes\cnum^{\ast}
\lrarr
 M^{\lor}_{\nbiga}\otimes\cnum^{\ast}$
induces a morphism
$\gamma_{\nbiga}:S_{\nbiga}^{\reg}\lrarr W$.
We set $W^{\reg}:=W\cap U_{\nbiga}^{\reg}$.
We have $\gamma_{\nbiga}(S^{\reg}_{\nbiga})\subset W^{\reg}$.

We use the notation in \S\ref{section;15.5.2.2}.
By applying the construction in \S\ref{subsection;15.5.11.10}
to $\nbiga=\{\veca_{1},\ldots,\veca_{m+r}\}\subset\seisuu^{n+r}$
with $\vecbeta=0\in\cnum^{n+r}$,
we obtain the $\nbigrtilde_{\cnum^{m+r}}$-modules
$\nbigm^{\GKZ}(\nbiga,K(\nbiga),0)$
and 
$\nbigm^{\GKZ}(\nbiga,K(\nbiga)^{\circ},0)$
on $W=\cnum^{m+r}=\{(\alpha_1,\ldots,\alpha_{m+r})\}$.

We have the induced morphism
$\id\times\gamma_{\nbiga}:
 \cnum_{\lambda}\times S_{\nbiga}^{\reg}
\lrarr 
 \cnum_{\lambda}\times W^{\reg}$.
The restrictions of 
$\nbigm^{\GKZ}(\nbiga,K(\nbiga),0)$
and 
$\nbigm^{\GKZ}(\nbiga,K(\nbiga)^{\circ},0)$
to $W^{\reg}$ give locally free 
$\nbigo_{\cnum\times W^{\reg}}$-modules.
We take the pull back of them by
$\id\times\gamma_{\nbiga}$
as $\nbigo$-modules.
Then, 
$(\id\times\gamma_{\nbiga})^{\ast}
 \nbigm^{\GKZ}(\nbiga,K(\nbiga),0)$
and 
$(\id\times\gamma_{\nbiga})^{\ast}
 \nbigm^{\GKZ}(\nbiga,K(\nbiga)^{\circ},0)$
are naturally equipped with the meromorphic flat connection
with which they are 
$\nbigrtilde_{S_{\nbiga}^{\reg}}$-modules.

\begin{thm}
\label{thm;15.5.2.4}
We have the following commutative diagram.
\[
 \begin{CD}
  \lambda^n(\id\times\gamma_{\nbiga})^{\ast}
 \nbigm^{\GKZ}(\nbiga,K(\nbiga)^{\circ},0)
 @>{\simeq}>>
 \lambda^{-r}\gbigv_{\nbiga!}
 \\
 @VVV @VVV \\
 \lambda^n(\id\times\gamma_{\nbiga})^{\ast}
 \nbigm^{\GKZ}(\nbiga,K(\nbiga),0)
 @>{\simeq}>>
 \lambda^{-r}\gbigv_{\nbiga\ast}
 \end{CD}
\]
\end{thm}
\pf
By Proposition \ref{prop;14.12.3.20}
and 
Proposition \ref{prop;15.5.11.1},
we obtain the following commutative diagram:
\begin{equation}
\label{eq;15.5.2.3}
 \begin{CD}
  (\id\times\gamma_{\nbiga})^{\ast}
 \nbigm^{\GKZ}\bigl(\nbiga,K(\nbiga)^{\circ},0
 \bigr)
 @>{\simeq}>>
 \nbigv_{\nbiga!}\\
 @VVV @VVV \\
 (\id\times\gamma_{\nbiga})^{\ast}
 \nbigm^{\GKZ}(\nbiga,K(\nbiga),0)
 @>{\simeq}>>
 \nbigv_{\nbiga\ast}
 \end{CD}
\end{equation}
Then, we obtain the claim of the theorem
by construction.
\hfill\qed

\vspace{.1in}

According to \cite{Reichelt-Sevenheck2},
we have $K(\nbiga)=\seisuu_{\geq 0}\nbiga$
and $K(\nbiga)^{\circ}=K(\nbiga)+\sum_{i=1}^r\veca_{n+i}$.
We set $\vecc_1=\sum_{i=1}^r\veca_{n+i}$.
Let $\nbigi(\nbiga,0)$
and $\nbigi(\nbiga,-\vecc_1)$
be the left ideals of $\nbigrtilde_W$
given as in \S\ref{subsection;15.5.11.20}.
As remarked in \S\ref{subsection;15.5.11.20},
we have the following commutative diagram
on $W=\cnum^{m+r}$:
\[
 \begin{CD}
 \nbigm^{\GKZ}(\nbiga,K(\nbiga)^{\circ},0)
 @>{\simeq}>>
 \nbigrtilde_{W}/\nbigi(\nbiga,-\vecc_1)
 \\
 @VVV @VVV \\
 \nbigm^{\GKZ}(\nbiga,K(\nbiga),0)
 @>{\simeq}>>
 \nbigrtilde_{W}\big/\nbigi(\nbiga,0)
 \end{CD}
\]
Here, the right vertical arrow is induced by the multiplication of
$\prod_{i=1}^r(\lambda\del_{m+i})$.

\subsubsection{Relation with quantum $\nbigd$-modules
 of complete intersections (Appendix)}

Inspired by the work of 
E. Mann and T. Mignon \cite{Mignon-Mann},
Reichelt and Sevenheck 
{\rm\cite{Reichelt-Sevenheck2}} constructed 
 the $\nbigrtilde$-modules
 ${}^{\ast}_0\widehat{\nbign}^{(0,\underline{0},\underline{0})}_{\nbiga}$
and 
 ${}^{\ast}_0\widehat{\nbigm}^{(-r,\underline{0},\underline{0})}_{\nbiga}$
on $W\setminus\{0\}$
in terms of the differential systems.
(See {\rm\cite{Reichelt-Sevenheck2}} for
the precise definition of the $\nbigrtilde$-modules.
But, we remark the commutativity (\ref{eq;15.11.22.1})
below.)
See \cite{Mignon-Mann} and \cite{Reichelt-Sevenheck2}
for the precise relation of the $\nbigrtilde$-modules
and the reduced quantum $\nbigd$-module
of complete intersection of general sections of
$\nbigl_i$ $(i=1,\ldots,r)$.
By formal computations,
we can obtain the following commutative diagram 
on $(\cnum^{\ast})^{m+r}=\bigl\{
 (\alpha_1,\ldots,\alpha_{m+r})\,\big|\,
 \alpha_i\neq 0
 \bigr\}$:
\begin{equation}
\label{eq;15.11.22.1}
 \begin{CD}
 \bigl(
  \lambda^n
 \nbigrtilde_{W}\big/
 \nbigi(\nbiga,-\vecc_1)
 \bigr)_{|(\cnum^{\ast})^{m+r}}
 @>{\simeq}>>
 {}^{\ast}_0\widehat{\nbign}^{(0,\underline{0},\underline{0})}
 _{\nbiga|(\cnum^{\ast})^{m+r}}
 \\
 @VVV @VVV \\
  \bigl(
  \lambda^n
 \nbigrtilde_{W}\big/
 \nbigi(\nbiga,0)
 \bigr)_{|(\cnum^{\ast})^{m+r}}
 @>{\simeq}>>
 {}^{\ast}_0\widehat{\nbigm}^{(-r,\underline{0},\underline{0})}
 _{\nbiga|(\cnum^{\ast})^{m+r}}
 \end{CD}
\end{equation}
The vertical arrows are the natural morphisms.

\begin{cor}
We have the following commutative diagram:
\[
 \begin{CD}
 \lambda^{-r}\gbigv_{\nbiga!}
 @>{\simeq}>>
 (\id\times\gamma_{\nbiga})^{\ast}\bigl(
 {}^{\ast}_0\widehat{\nbign}^{(0,\underline{0},\underline{0})}
 _{\nbiga|(\cnum^{\ast})^{m+r}}\bigr)
 \\
 @VVV @VV{a_1}V \\
 \lambda^{-r}\gbigv_{\nbiga\ast}
 @>{\simeq}>>
 (\id\times\gamma_{\nbiga})^{\ast}\bigl(
 {}^{\ast}_0\widehat{\nbigm}^{(-r,\underline{0},\underline{0})}
 _{\nbiga|(\cnum^{\ast})^{m+r}}\bigr)
 \end{CD}
\]
In particular,
the $\nbigrtilde$-modules
$(\id\times\gamma_{\nbiga})^{\ast}\bigl(
 {}^{\ast}_0\widehat{\nbign}^{(0,\underline{0},\underline{0})}
 _{\nbiga|(\cnum^{\ast})^{m+r}}\bigr)$
and 
$(\id\times\gamma_{\nbiga})^{\ast}\bigl(
 {}^{\ast}_0\widehat{\nbigm}^{(-r,\underline{0},\underline{0})}
 _{\nbiga|(\cnum^{\ast})^{m+r}}\bigr)$
underlie mixed twistor $\nbigd$-modules,
and the $\Image a_1$ underlies
a pure twistor $\nbigd$-module.
\hfill\qed
\end{cor}

\begin{rem}
In {\rm\cite{Reichelt-Sevenheck2}},
Reichelt and Sevenheck constructed
$\nbigrtilde$-modules by using
the partial Fourier-Laplace transform 
of GKZ-hypergeometric systems
and the Brieskorn lattices
associated to the Hodge filtrations
of the GKZ-hypergeometric systems.
They conjectured in
{\rm\cite[Conjecture {\rm 6.13}]{Reichelt-Sevenheck2}}
that the $\nbigrtilde$-modules are isomorphic to
$(\id\times\gamma_{\nbiga})^{\ast}\bigl(
 {}^{\ast}_0\widehat{\nbign}^{(0,\underline{0},\underline{0})}
 _{\nbiga|(\cnum^{\ast})^{m+r}}\bigr)$
and 
$(\id\times\gamma_{\nbiga})^{\ast}\bigl(
 {}^{\ast}_0\widehat{\nbigm}^{(-r,\underline{0},\underline{0})}
 _{\nbiga|(\cnum^{\ast})^{m+r}}\bigr)$.

In {\rm \S\ref{section;14.12.27.1}},
we review their construction of $\nbigrtilde$-modules.
In Proposition {\rm\ref{prop;14.11.12.3}},
we compare their $\nbigrtilde$-modules
with the underlying $\nbigrtilde$-modules 
of mixed twistor $\nbigd$-modules
in a general setting.
See also Proposition {\rm\ref{prop;14.12.15.1}}.
So, Theorem {\rm \ref{thm;15.5.2.4}}
also verifies their conjecture.
Reichelt and Sevenheck also verified
their conjecture {\rm \cite{Reichelt-Sevenheck3}}
in a different way.
\hfill\qed
\end{rem}

\begin{rem}
\label{rem;15.11.22.2}
In the first version of this preprint,
we used some results in {\rm\cite{Reichelt-Sevenheck2}}
for the comparison with
some $\nbigrtilde$-modules
associated to differential systems,
which required us to go a roundabout way
with careful comparison of dualities.
In the second version,
as explained \S\ref{section;15.5.2.2},
we give a direct comparison of
$\nbigrtilde$-modules associated to differential systems
and $\nbigrtilde$-modules 
underlying mixed twistor $\nbigd$-modules.
Hence, the argument is made more transparent.
\hfill\qed
\end{rem}

\subsection{Description as a specialization}
\label{subsection;14.12.29.1}

This subsection is the continuation of
\S\ref{subsection;14.11.27.30}--\S\ref{subsection;14.12.2.30},
and it is the preparation for Theorem \ref{thm;14.11.21.10}.

\subsubsection{The $\nbigrtilde$-module with pairing
associated to $\nbigatilde$}
\label{subsection;14.12.30.1}

Applying the construction in \S\ref{subsection;14.11.26.200},
we obtain the integrable pure twistor $\nbigd$-module
$\nbigt_{\nbigatilde}$ on $S_{\nbigatilde}$.
Let $\nbigv_{\nbigatilde}$
denote the underlying $\nbigrtilde$-module
which is $\cnum^{\ast}$-homogeneous.
We have the filtration $\Wtilde$ and 
the graded pairing $P_{\nbigatilde}$
induced by the real structure and the graded polarization
on $\nbigv_{\nbigatilde}$.
In this case,  $\Wtilde$ is pure of weight $n+r$.

We apply the same construction
to $\nbigt_{\nbigatilde}\otimes\newTate(n+r)$.
We obtain the underlying $\nbigrtilde$-module
$\gbigv_{\nbigatilde}$ on $S_{\nbigatilde}$
which is $\cnum^{\ast}$-homogeneous.
The restriction 
$\gbigv_{\nbigatilde|S^{\reg}_{\nbigatilde}}$ is equipped 
with the filtration $\Wtilde$ which is pure of weight $-n-r$,
and the pairing $\gbigp_{\nbigatilde}$ of weight $-n-r$.
In other words,
we obtain a ${\rm TEP}(n+r)$-structure
$\bigl(\gbigv_{\nbigatilde|S^{\reg}_{\nbigatilde}},
 \gbigp_{\nbigatilde}
 \bigr)$.
By construction,
we have
$\gbigv_{\nbigatilde}=
 \lambda^{n+r}\nbigv_{\nbigatilde}$,
and 
$P_{\nbigatilde}=\gbigp_{\nbigatilde}$
on 
$\nbigv_{\nbigatilde}(\ast\lambda)
 =\gbigv_{\nbigatilde}(\ast\lambda)$.

\subsubsection{Relation of $\gbigv_{\nbiga}$ and $\gbigv_{\nbigatilde}$}

We set $\Stilde_{\nbiga}:=S_{\nbiga}\times\cnum_{\tau}$.
Applying the construction in \S\ref{subsection;14.11.21.21},
we obtain the integrable pure twistor $\nbigd$-module
$\nbigt(\nbigatilde,\Stilde_{\nbiga},\gammatilde_{\nbiga})$
and the underlying $\nbigrtilde$-module
$\nbigl(\nbigatilde,\Stilde_{\nbiga},\gammatilde_{\nbiga})$
on $\proj^{m+r+1}\times \Stilde_{\nbiga}$.
Let $\pitilde:
 \proj^{m+r+1}\times\Stilde_{\nbiga}
\lrarr
 \Stilde_{\nbiga}$
be the projection.
We set
$\gbigvtilde_{\nbiga}:=
 \lambda^{n+r}\pitilde_{\dagger}^0
 \nbigl(\nbigatilde,\Stilde_{\nbiga},\gammatilde_{\nbiga})(\ast\tau)$.
It is equipped with the pairing 
$\gbigptilde_{\nbiga}$ of weight $-n-r$
as in the case of \S\ref{subsection;14.12.30.1}.

\begin{lem}
\label{lem;14.11.27.100}
We can naturally identify
$S_{\nbigatilde}$ with
$S_{\nbiga}\times\cnum_{\tau}^{\ast}
\subset \Stilde_{\nbiga}$
under which
$\gbigvtilde_{\nbiga|S_{\nbigatilde}}
=\gbigv_{\nbigatilde}$.
In other words,
$\gbigvtilde_{\nbiga}$
is a meromorphic extension of 
$\gbigv_{\nbigatilde}$
on $\Stilde_{\nbiga}$.
We also have the $\cnum^{\ast}$-action 
on $\Stilde_{\nbiga}$
for which the inclusions
$S_{\nbiga}\simeq S_{\nbiga}\times\{0\}\subset\Stilde_{\nbiga}$
and 
$S_{\nbigatilde}\lrarr \Stilde_{\nbiga}$
are equivariant.
\end{lem}
\pf
We use the notation in \S\ref{subsection;14.11.26.200}.
We naturally have 
$M_{\nbigatilde}=M_{\nbiga}\oplus\seisuu e_{m+r+1}$
and $N_{\nbigatilde}=N_{\nbiga}$.
The natural inclusion $M_{\nbiga}\lrarr M_{\nbigatilde}$
induces $L_{\nbiga}\lrarr L_{\nbigatilde}$,
and hence
$L_{\nbigatilde}^{\lor}\lrarr L_{\nbiga}^{\lor}$.

We also consider $M_X:=\seisuu^m$
and $N_X:=\seisuu^n$
and a morphism
$\Theta_X:M_X\lrarr N_X$ determined by
$\Theta_X(e_i)=[\rho_i]$,
where the tuple $e_1,\ldots,e_m$  is the standard basis of $M_X$,
and $[\rho_i]$ are as in \S\ref{subsection;14.11.27.30}.
Let $L_X$ be the kernel of $\Theta_X$.
The projections
$M_{\nbigatilde}\lrarr M_X$
and $N_{\nbigatilde}\lrarr N_X$
induces
$L_{\nbigatilde}\lrarr L_X$.
It is well known and easy that
the composite of
$L_{\nbiga}\lrarr L_{\nbigatilde}\lrarr L_X$
is an isomorphism.

Let $\Upsilon$ denote the image of
$e_{m+r+1}^{\lor}$ by
$M_{\nbigatilde}^{\lor}\lrarr L_{\nbigatilde}^{\lor}$.
Then, $\Upsilon$ is the frame of
the kernel of $L_{\nbigatilde}^{\lor}\lrarr L_{\nbiga}^{\lor}$.
Hence, $L_{\nbigatilde}^{\lor}$
is the direct sum of the image of $L_X^{\lor}$
and $\seisuu\Upsilon$.
The direct sum induces
\[
 S_{\nbigatilde}
=L^{\lor}_{\nbigatilde}\otimes\cnum^{\ast}
=\bigl(
 L^{\lor}_{\nbiga}\otimes\cnum^{\ast}
 \bigr)
\times
 \cnum^{\ast}
=S_{\nbiga}\times\cnum^{\ast}.
\]
We may assume that the splitting
$\kappa_{\nbigatilde}:L_{\nbigatilde}^{\lor}\lrarr
M_{\nbigatilde}^{\lor}$
is given by the direct sum of
$\kappa_{\nbiga}:L_{\nbiga}^{\lor}\lrarr M_{\nbiga}^{\lor}$
and $\seisuu\Upsilon\lrarr \seisuu e_{m+r+1}^{\lor}$.
Then, 
we have
$F_{\gamma_{\nbigatilde}}(x,\tau,t)
=F_{\gamma_{\nbiga}}(x,t)+
 \tau t^{\veca_{m+r+1}}$.
Then, by the construction of $\gammatilde_{\nbiga}$,
we have
$\gammatilde_{\nbiga|S_{\nbigatilde}}
=\gamma_{\nbigatilde}$.
We obtain
$\nbigt(\nbigatilde,\Stilde_{\nbiga},\gammatilde_{\nbiga})
_{|\proj^m\times S_{\nbigatilde}}
=\nbigt(\nbigatilde,S_{\nbigatilde},\gamma_{\nbigatilde})$,
and hence
$\gbigvtilde_{\nbiga|S_{\nbigatilde}}
=\gbigv_{\nbigatilde}$.
The claim on the $\cnum^{\ast}$-action is also easy to see.
\hfill\qed

\vspace{.1in}
Recall the following due to Reichelt-Sevenheck
who proved it in a more general situation.

\begin{prop}[\cite{Reichelt-Sevenheck1}]
$\gbigvtilde_{\nbiga}$ is regular along $\tau$.
\hfill\qed
\end{prop}

\begin{lem}
We consider the morphism
$\lambda N:\psi_{\tau,-\vecdelta}\gbigvtilde_{\nbiga}
\lrarr
 \psi_{\tau,-\vecdelta}\gbigvtilde_{\nbiga}$,
where $N$ is induced by $\tau\del_{\tau}$.
\begin{itemize}
\item
The kernel and the cokernel of $\lambda N$
are isomorphic to
$\gbigv_{\nbiga,!}$
and 
$\gbigv_{\nbiga,\ast}$
respectively.
\item
Let $W(N)$ denote the weight filtration of $\lambda N$
on $\psi_{\tau,-\vecdelta}\gbigvtilde_{\nbiga}$.
The induced filtrations on
$\Ker(\lambda N)$ and $\Cok(\lambda N)$ are also denoted by
$W(N)$.
Then,
\[
 W(N)_{k}\bigl( \Cok(\lambda N)\bigr)
\simeq
 \Wtilde_{k-n-r}\bigl(
 \gbigv_{\nbiga,\ast}
 \bigr)
\]
\[
  W(N)_{k}\bigl(\Ker(\lambda N)\bigr)
\simeq
 \Wtilde_{k-n-r}\bigl(
 \gbigv_{\nbiga,!}
 \bigr)
\]
\end{itemize}
\end{lem}
\pf
The first follows from the isomorphisms in
Corollary \ref{cor;14.11.28.50}
and the constructions of
$\gbigv_{\nbiga}$
and $\gbigvtilde_{\nbiga}$.
The second also follows from 
the isomorphisms in Corollary \ref{cor;14.11.28.50},
the comparison of the filtrations
$W(\nbign)$ and $W$ on $\psi_{\tau,-\vecdelta}(\nbigt)$,
and the constructions of the filtrations $\Wtilde$.
\hfill\qed

\vspace{.1in}

Let $\Btilde$ be any open set in $\Stilde_{\nbiga}$
such that 
$\gbigvtilde_{\nbiga}$
is locally free 
$\nbigo_{\cnum_{\lambda}\times \Btilde}(\ast\tau)$-module,
and regular singular along $\tau=0$.
The polarization and the real structure of
$\pi_{\dagger}^0
 \nbigt(\nbigatilde,\Stilde_{\nbiga},\gammatilde_{\nbiga})$
induce a pairing $\gbigptilde_{\nbiga}$ of weight $-n-r$ 
on $\gbigvtilde_{\nbiga|\Btilde}$.
We clearly have
$\gbigptilde_{\nbiga|\Btilde\setminus\{\tau=0\}}
=\gbigp_{\nbigatilde|\Btilde\setminus\{\tau=0\}}$.
We obtain the induced graded pairing
$\sp_{\tau}(\gbigptilde_{\nbiga})$
on $(\Cok(\lambda N),\Wtilde)$
by the procedure in \S\ref{subsection;14.11.28.130}.
We also have 
$\sp_{\tau}(\gbigptilde_{\nbiga})$
on $(\Ker(\lambda N),\Wtilde)$.

\begin{lem}
We have 
$\sp_{\tau}(\gbigptilde_{\nbiga})
=\gbigp_{\nbiga,\ast}$
under the isomorphism
$(\Cok(\lambda N),\Wtilde)
 \simeq
 (\gbigv_{\nbiga\ast},\Wtilde)$.
We also have
$\sp_{\tau}(\gbigptilde_{\nbiga})
=\gbigp_{\nbiga,!}$
under the isomorphism
$(\Ker(\lambda N),\Wtilde)
 \simeq
 (\gbigv_{\nbiga!},\Wtilde)$.
\end{lem}
\pf
It follows from 
Proposition \ref{prop;14.11.18.52}.
\hfill\qed

\subsubsection{Cylindrical ends}
\label{subsection;14.11.22.10}

We set 
$Y:=\proj\bigl(\nbigo_X\oplus\bigoplus_{j=1}^{r}\nbigl_j\bigr)$
with the projection $\pi:Y\lrarr X$.
It is the projective completion of 
the vector bundle $\bigoplus_{j=1}^{r}\nbigl_j^{\lor}$.
Recall that 
$Y$ is also weak Fano.
Indeed, for the canonical bundle $K_Y$ of $Y$,
we have
$K_{Y}^{\lor}
=\pi^{\ast}\bigl(
 K_X\otimes\bigotimes_{i=1}^r\nbigl_i\bigr)^{\lor}
\otimes
 \nbigo_Y(1)^{\otimes (r+1)}$.
By the assumption,
$\pi^{\ast}(K_X\otimes 
 \bigotimes_{i=1}^r\nbigl_i )^{\lor}$
is nef.
Because $\nbigl_i$ are nef line bundles
on a toric variety,
they are globally generated.
Then, it is easy to see that
$\nbigo_Y(1)$ is also globally generated.
In particular,
it is nef.
As a result, $K_Y^{\lor}$ is nef.

We have the natural morphisms:
\[
 \begin{CD}
 L_{X}^{\lor} @>>>
 L_{\nbigatilde}^{\lor} @>>>
 L_{\nbiga}^{\lor}\\
 @V{\simeq}VV @V{\simeq}VV @V{\simeq}VV \\
 H^2(X,\seisuu) @>{\pi^{\ast}}>>
 H^2(Y,\seisuu) @>{i_0^{\ast}}>>
 H^2(X,\seisuu)
 \end{CD}
\]
We take a frame 
$\eta_1,\ldots,\eta_{\ell}$ of
$H^2(X,\seisuu)$
such that each $\eta_i$ is the first Chern class
of a nef line bundle.
Let $\eta_{\ell+1}\in H^2(Y,\seisuu)$ be the first Chern class 
of the tautological line bundle $\nbigo(1)$ of $Y$ over $X$.
By the assumption on the line bundles $\nbigl_j$,
$\eta_{\ell+1}$ is also a nef class.
The tuple $\eta_1,\ldots,\eta_{\ell+1}$
gives a frame of $H^2(Y,\seisuu)$.
Let $\xi_i\in L_{\nbiga}^{\lor}$
be the elements corresponding to $\eta_i$.
We have $\xi_{\ell+1}=\Upsilon$.

The tuple $\xi_1,\ldots,\xi_{\ell}$ gives
a frame of $L_{\nbiga}^{\lor}$.
It gives a coordinate system
$(x_1,\ldots,x_{\ell})$ on $S_{\nbiga}$
with which 
$S_{\nbiga}\simeq (\cnum^{\ast})^{\ell}$.
Similarly,
the tuple $\xi_1,\ldots,\xi_{\ell+1}$ gives
a frame of $L_{\nbigatilde}^{\lor}$.
It gives a coordinate system
$(x_1,\ldots,x_{\ell+1})$ on $S_{\nbigatilde}$
with which 
$S_{\nbigatilde}\simeq (\cnum^{\ast})^{\ell+1}$.
We have 
$\tau=x_{\ell+1}$
under the identification in Lemma \ref{lem;14.11.27.100}.
Recall the following due to Iritani.

\begin{lem}[\cite{Iritani}]
\label{lem;14.12.2.1}
There exists $\epsilon>0$
such that the following holds:
\begin{itemize}
\item
$B_{\nbigatilde}=\bigl\{
 (x_1,\ldots,x_{\ell+1})\,\big|\,
 0<|x_i|<\epsilon\,\,(i=1,\ldots,\ell+1)
 \bigr\}$
is contained in
$S_{\nbigatilde}^{\reg}$.
\item
$B_{\nbiga}=\bigl\{
 (x_1,\ldots,x_{\ell})\,\big|\,
 0<|x_i|<\epsilon\,\,(i=1,\ldots,\ell)
 \bigr\}$
is contained in
$S_{\nbiga}^{\reg}$.
\end{itemize}
\end{lem}
\pf
The first claim is proved  in \cite[Appendix A1]{Iritani}.
Let $\sigma$ be a face of $\Conv(\nbiga\cup\{0\})$
such that $0$ is not contained in $\sigma$.
Then, $\sigma$ is a face of $\Conv(\nbigatilde)$.
We have the equality
$F_{\gamma_{\nbiga},\sigma}
=F_{\gamma_{\nbigatilde},\sigma}$,
and 
$F_{\gamma_{\nbiga},\sigma}$
does not contain $x_{\ell+1}$.
Hence, the second claim follows.
\hfill\qed

\vspace{.1in}

We consider the following subset $\Btilde_{\nbiga}$
in $\Stilde_{\nbiga}$:
\[
 \Btilde_{\nbiga}=
 \bigl\{
 (x_1,\ldots,x_{\ell+1})\,\big|\,
 0<|x_i|<\epsilon\,\,(i=1,\ldots,\ell),\,\,
 |x_{\ell+1}|<\epsilon
 \bigr\}
\]
The restriction of $\gbigvtilde_{\nbiga|B_{\nbigatilde}}$
is a locally free 
$\nbigo_{\cnum_{\lambda}\times B_{\nbigatilde}}
 (\ast x_{\ell+1})$-module,
and it is regular along $x_{\ell+1}$.

\subsubsection{Coverings}

We set $\gbigs_{\nbiga}:=
\cnum\otimes L_{\nbiga}^{\lor}$.
The exponential map
$\cnum\lrarr \cnum^{\ast}$
induces 
the covering map
$\chi_{\nbiga}:
 \gbigs_{\nbiga} \lrarr
 S_{\nbiga}$.
We obtain an $\nbigrtilde$-module
$\chi^{\ast}_{\nbiga}\gbigv_{\nbiga\star}$
equipped with the filtration $\chi^{\ast}_{\nbiga}\Wtilde$
on $\gbigs_{\nbiga}$.
The restriction to
$\gbigs^{\reg}_{\nbiga}:=
 \chi^{-1}_{\nbiga}(S^{\reg}_{\nbiga})$
are equipped with graded pairings 
$\chi_{\nbiga}^{\ast}\gbigp_{\nbiga\star}$.

We set 
$\widetilde{\gbigs}_{\nbiga}:=
 \cnum\otimes L^{\lor}_{\nbigatilde}$.
By the decomposition
$L_{\nbigatilde}^{\lor}=L_{\nbiga}^{\lor}
\oplus\seisuu\Upsilon$
in the proof of Lemma \ref{lem;14.11.27.100},
we naturally have
$\widetilde{\gbigs}_{\nbiga}
=\gbigs_{\nbiga}\times\cnum$.
The map $\chi_{\nbiga}$
and the identity on $\cnum$
induces 
$\chitilde_{\nbiga}:\widetilde{\gbigs}_{\nbiga}\lrarr \Stilde_{\nbiga}$.

Take a frame $\xi_1,\ldots,\xi_{\ell}$ of $L_{\nbiga}^{\lor}$
and an open subset $\Btilde_{\nbiga}$
as in \S\ref{subsection;14.11.22.10}.
The restriction 
of $\chitilde_{\nbiga}^{\ast}\gbigvtilde_{\nbiga}$
to $\chitilde_{\nbiga}^{-1}(\Btilde_{\nbiga})$
is equipped with 
the pairing
$\chitilde_{\nbiga}^{\ast}\gbigptilde_{\nbiga}$
of weight $-n-r$.

\begin{lem}
\label{lem;14.11.22.1}
The restriction of
the mixed TEP-structure
$\chi^{\ast}_{\nbiga}\bigl(
 \gbigv_{\nbiga\star|S_{\nbiga}^{\reg}},
 \Wtilde,\gbigp_{\nbiga\ast}
 \bigr)$
to $\chitilde_{\nbiga}^{-1}(B_{\nbiga})$
$(\star=\ast,!)$
are obtained from 
$\chitilde^{\ast}_{\nbiga}(\gbigvtilde_{\nbiga},\gbigptilde_{\nbiga})
 _{|\Btilde_{\nbiga}}$
by the procedure in {\rm\S\ref{subsection;14.11.28.130}}.
\hfill\qed
\end{lem}

\subsubsection{Logarithmic extension and 
endomorphisms of $\gbigv_{\nbigatilde}$}

This is the continuation of \S\ref{subsection;14.11.22.10}.
Let us recall 
the explicit description of 
the logarithmic extension of $\gbigv_{\nbigatilde}$
due to Reichelt and Sevenheck in \cite{Reichelt-Sevenheck1}.
See \cite{Reichelt-Sevenheck1} for more details.
We also give an easy remark on uniqueness of automorphisms,
although it is also essentially implied in \cite{Reichelt-Sevenheck1}.

Let $B_{\nbigatilde}$ be as in Lemma \ref{lem;14.12.2.1}.
We set 
\[
 \Bbar_{\nbigatilde}:=
 \bigl\{
 (x_1,\ldots,x_{\ell+1})\,\big|\,
 |x_i|<\epsilon\,\,(i=1,\ldots,\ell+1)
 \bigr\}.
\]
First, the following holds.
\begin{prop}[\S3.2 of \cite{Reichelt-Sevenheck1}]
There exists a locally free
$\nbigo_{\cnum_{\lambda}\times \Bbar_{\nbigatilde}}$-module
$\gbigq_{\nbigatilde}$
with an isomorphism
\[
\gbigq
 _{\nbigatilde|\cnum_{\lambda}\times B_{\nbigatilde}}
\simeq
\gbigv_{\nbigatilde|B_{\nbigatilde}}
\]
such that the following holds:
\begin{itemize}
\item
The meromorphic flat connection $\nabla$ 
of $\gbigv_{\nbigatilde|B_{\nbigatilde}}$
gives a meromorphic flat connection 
$\nabla^{\gbigq_{\nbigatilde}}$ of $\gbigq_{\nbigatilde}$
which is logarithmic along $\{x_i=0\}$ $(i=1,\ldots,\ell+1)$.
\item
The residues 
$\Res_{x_i}(\nabla^{\gbigq_{\nbigatilde}})
 _{|\cnum_{\lambda}\times \{0\}}$ 
are nilpotent.
\end{itemize}
Such $\gbigq_{\nbigatilde}$
is unique up to canonical isomorphisms.
\hfill\qed
\end{prop}

They considered the residual connection $\nabla^E$
on $E:=\gbigq_{\nbigatilde|\cnum_{\lambda}\times\{0\}}$.
They give a quite explicit description
of the meromorphic flat bundle
$\bigl(
 E,\nabla^E
 \bigr)$.
Let $N_{i}:E\lrarr \lambda^{-1}E$ $(i=1,\ldots,\ell+1)$
denote 
$\Res_{x_i}(\nabla^{\gbigq_{\nbigatilde}})
 _{|\cnum_{\lambda}\times\{0\}}$.
We have the fundamental vector fields
$\underline{\gminiv}=\sum_{i=1}^{\ell+1} k_ix_i\del_{x_i}$
of the $\cnum^{\ast}$-action on $S_{\nbigatilde}$
(see \S\ref{subsection;14.11.26.201}).
Note that for the $\cnum^{\ast}$-homogeneity of 
$\gbigv_{\nbigatilde}$,
the fundamental vector field of the $\cnum^{\ast}$-action
on $\cnum_{\lambda}\times S_{\nbigatilde}$
is $\lambda\del_{\lambda}+\underline{\gminiv}$.
Let us consider the $\cnum$-homomorphism
$L:E\lrarr E$
given by
$\nabla^E(\lambda\del_{\lambda})
+\sum_{i=1}^{\ell+1} k_iN_i$.
Because $[L,N_i]=0$,
we have $[L,\lambda N_i]=\lambda N_i$.
For the following proposition,
and more detailed description 
in terms of the cohomology group of $Y$,
see Lemma 3.8 and its proof in \cite{Reichelt-Sevenheck1}.

\begin{prop}[\S3.2 of \cite{Reichelt-Sevenheck1}]
\label{prop;14.12.2.10}
There exists a frame $u_1,\ldots,u_{\rank E}$
for which we have the following:
\begin{itemize}
\item
$L(u_j)=c_ju_j$ for $c_j\in\seisuu$.
\item
We have $c_1<c_i$ $(i\neq 1)$.
\item
We have $N_iu_j=\sum \alpha^i_{kj}u_k$
for some $\alpha^i_{kj}\in\seisuu$.
In particular,
$\Sym\Bigl(
 \bigoplus_{i=1}^{\ell+1} \cnum N_i
 \Bigr)$
acts on $F:=\bigoplus_{i=1}^{\rank E} \cnum u_i$.
\item
The map
$\Sym\Bigl(
 \bigoplus_{i=1}^{\ell+1} \cnum N_i
 \Bigr)
\lrarr F$
given by 
$P(N_1,\ldots,N_{\ell+1})
\longmapsto
 P(N_1,\ldots,N_{\ell+1})u_{1}$
is surjective.
\hfill\qed
\end{itemize}
\end{prop}

By the first property in Proposition \ref{prop;14.12.2.10},
we have the $\cnum^{\ast}$-action on $E$
which induces $L$.
As a result, we have a vector space 
$H$ with an increasing filtration $F$ indexed by integers
such that the $\cnum^{\ast}$-equivariant bundle $E$ 
is isomorphic to the analytification
of the Rees module of $(H,F)$.
Let $d_0:=\min\{d\,|\,\Gr^F_d(H)\neq 0\}$.
The second property implies that
$\dim F_{d_0}(H)=1$.
The third and the fourth properties imply that 
$F_{d_0}(H)$ generates $H$
over the induced actions of the residues.
The following is also essentially implied in 
the proof of \cite[Lemma 3.8]{Reichelt-Sevenheck1}.

\begin{cor}
\label{cor;14.12.2.21}
For any endomorphism $\varphi$
of the $\nbigrtilde$-module
$\gbigv_{\nbigatilde|B_{\nbigatilde}}$,
there exists a complex number $\alpha$
such that $\varphi=\alpha\id$.
\end{cor}
\pf
Because 
$\gbigq_{|\cnum_{\lambda}^{\ast}\times \Bbar_{\nbigatilde}}$
is the Deligne extension of
$\gbigv_{\nbigatilde|\cnum_{\lambda}^{\ast}\times B_{\nbigatilde}}$
to $\cnum_{\lambda}^{\ast}\times \Bbar_{\nbigatilde}$,
the morphism $\varphi$ is extended to an endomorphism of
the $\nbigo$-module $\gbigq_{\nbigatilde}$
compatible with $\nabla^{\gbigq_{\nbigatilde}}$.
Hence, it induces an endomorphism $\varphi^E$ of
$(E,\nabla^E)$
compatible with the actions of the residues.
Because the $\cnum^{\ast}$-action
is determined by $\nabla^E$ and the residues,
$\varphi^E$ is $\cnum^{\ast}$-equivariant.
Hence, it induces an endomorphism $\varphi^H$ of
$(H,F)$ compatible with the induced actions of residues.
Because $F_{d_0}$ generates $H$ over the actions of
the residues,
$\varphi^H$ is uniquely determined by 
the restriction $\varphi^H_{|F_{d_0}}$.
Because $\dim F_{d_0}=1$,
we have $\alpha\in\cnum$ such that 
$\varphi^H_{|F_{d_0}}$ is the multiplication of $\alpha$.
If $\alpha=0$, then $\varphi=0$.
Hence, we have
$\varphi=\alpha\id$.
\hfill\qed

\begin{cor}
\label{cor;14.12.2.22}
Let $U\subset S_{\nbigatilde}$ be any connected open subset
such that $B_{\nbigatilde}\subset U$.
For any endomorphism $\varphi$
of the $\nbigrtilde$-module
$\gbigv_{\nbigatilde|\cnum_{\lambda}\times U}$,
there exists a complex number $\alpha$
such that $\varphi=\alpha\id$.
\hfill\qed
\end{cor}

\begin{cor}
Let 
$P':\gbigv_{\nbigatilde|S^{\reg}_{\nbigatilde}}
 \otimes
 j^{\ast}\gbigv_{\nbigatilde|S^{\reg}_{\nbigatilde}}
\lrarr
 \lambda^{n+r}\nbigo_{\cnum_{\lambda}\times S^{\reg}_{\nbiga}}$
be a morphism of $\nbigrtilde$-modules.
Then, there exists a complex number $\alpha$
such that 
$P'=\alpha \gbigp_{\nbigatilde}$.
\hfill\qed
\end{cor}

\subsubsection{Logarithmic extension
and endomorphisms of $\gbigv_{\nbiga}$}

We can easily deduce the same property for $\gbigv_{\nbiga\star}$
($\star=\ast,!$).
Let $B_{\nbiga}$ be as in Lemma \ref{lem;14.12.2.1}.
We set
\[
 \Bbar_{\nbiga}:=
 \bigl\{
 (x_1,\ldots,x_{\ell})\,\big|\,
 |x_i|<\epsilon\,\,(i=1,\ldots,\ell)
 \bigr\}.
\]
We regard $\Bbar_{\nbiga}$
as $\bigl\{
 (x_1,\ldots,x_{\ell+1})\in\Bbar_{\nbigatilde}\,\big|\,
 x_{\ell+1}=0
 \bigr\}$.
We obtain the vector bundle
$\gbigq_{\nbigatilde|\Bbar_{\nbiga}}$.
It is equipped with the endomorphism
$\lambda\Res_{x_{\ell+1}}(\nabla^{\gbigq_{\nbigatilde}})$,
and the induced connection.

\begin{lem}
The conjugacy classes of
$\lambda\Res_{x_{\ell+1}}(\nabla^{\gbigq_{\nbigatilde}})_{|(\lambda,P)}$
are independent of the choice of 
$(\lambda,P)\in \cnum_{\lambda}\times \Bbar_{\nbiga}$.
\end{lem}
\pf
By the construction,
$\gbigv_{\nbigatilde}$ comes from
a harmonic bundle.
It is easy to see that 
$\gbigq_{\nbigatilde}$ is equal to 
the prolongment of the family of $\lambda$-flat bundles
$\nbigp_0\bigl(
 \gbigv_{\nbigatilde|B_{\nbigatilde}}
 \bigr)$
associated to tame harmonic bundles,
as remarked in \cite{Reichelt-Sevenheck1}.
Then, the claim follows from a general theory of
harmonic bundles \cite{Mochizuki-tame}.
\hfill\qed

\vspace{.1in}
Let $\gbigq_{\nbiga}$ 
be the cokernel
of $\lambda\Res_{x_{\ell+1}}(\nabla^{\gbigq})$.
It is a locally free
$\nbigo_{\cnum_{\lambda}\times\Bbar_{\nbiga}}$-module.
The following is clear 
by the relation of $\gbigv_{\nbiga\ast}$
and $\gbigv_{\nbigatilde}$.

\begin{lem}
We have natural isomorphisms
of $\nbigrtilde$-modules
$\gbigq_{\nbiga|\cnum_{\lambda}\times B_{\nbiga}}
\simeq
 \gbigv_{\nbiga\ast|\cnum_{\lambda}\times B_{\nbiga}}$.
\hfill\qed
\end{lem}

Set 
$E_{\nbiga}:=
 \gbigq_{\nbiga|\cnum_{\lambda}\times\{0\}}$.
It is equipped with the induced connection
$\nabla^{E_{\nbiga}}$.
We obtain an explicit description of 
$(E_{\nbiga},\nabla^{E_{\nbiga}})$
from the description of
$(E,\nabla^E)$ in \cite{Reichelt-Sevenheck1}.
Let $N_{i}:E_{\nbiga}\lrarr \lambda^{-1}E_{\nbiga}$
$(i=1,\ldots,\ell)$
denote 
$\Res_{x_i}(\nabla^{\gbigq_{\nbiga}})
 _{|\cnum_{\lambda}\times\{0\}}$.
We have the fundamental vector fields
$\underline{\gminiv}=\sum_{i=1}^{\ell} k_ix_i\del_{x_i}$
of the $\cnum^{\ast}$-action on $S_{\nbiga}$
(see \S\ref{subsection;14.11.26.201}).
Note that for the $\cnum^{\ast}$-homogeneity of
$\gbigv_{\nbiga\ast}$,
the fundamental vector field of the $\cnum^{\ast}$-action
on $\cnum_{\lambda}\times S_{\nbiga}$
is $\lambda\del_{\lambda}+\underline{\gminiv}$.
Let us consider the $\cnum$-homomorphism
$L:E_{\nbiga}\lrarr E_{\nbiga}$
given by
$\nabla^E(\lambda\del_{\lambda})
+\sum_{i=1}^{\ell} k_iN_i$.
Because $[L,N_i]=0$,
we have $[L,\lambda N_i]=\lambda N_i$.
The following is a direct consequence of 
Proposition \ref{prop;14.12.2.10}.

\begin{prop}
\label{prop;14.12.2.20}
There exists a frame $u_1,\ldots,u_{\rank E_{\nbiga}}$
for which we have the following:
\begin{itemize}
\item
$L(u_j)=c_ju_j$ for $c_j\in\seisuu$.
\item
We have $c_1<c_i$ $(i\neq 1)$.
\item
We have $N_iu_j=\sum \beta^{i}_{kj}u_k$
for $\beta^{i}_{kj}\in\seisuu$.
In particular,
$\Sym\Bigl(
 \bigoplus_{i=1}^{\ell} \cnum N_i
 \Bigr)$
acts on $F_{\nbiga}:=\bigoplus_{i=1}^{\rank E_{\nbiga}} \cnum u_i$.
\item
The map
$\Sym\Bigl(
 \bigoplus_{i=1}^{\ell} \cnum N_i
 \Bigr)
\lrarr F_{\nbiga}$
given by 
$P(N_1,\ldots,N_{\ell})
\longmapsto
 P(N_1,\ldots,N_{\ell})u_{1}$
is surjective.
\hfill\qed
\end{itemize}
\end{prop}
We omit a more detailed description 
in terms of the cohomology group of $X$.

By the first property in Proposition \ref{prop;14.12.2.20},
we have the $\cnum^{\ast}$-action on $E_{\nbiga}$
which induces $L$.
As a result, we have a vector space $H_{\nbiga}$
with an increasing filtration $F$ indexed by integers
such that the $\cnum^{\ast}$-equivariant bundle $E_{\nbiga}$ 
is isomorphic to the analytification
of the Rees module of $(H_{\nbiga},F)$.
Let $d_0:=\min\{d\,|\,\Gr^F_d(H)\neq 0\}$.
The second property implies that
$\dim F_{d_0}(H_{\nbiga})=1$.
The third and the fourth properties imply that 
$F_{d_0}(H_{\nbiga})$ generates $H_{\nbiga}$
over the induced actions of the residues.

\begin{cor}
Let $U\subset S_{\nbiga}$ be any connected open subset
such that $B_{\nbiga}\subset U$.
For any endomorphism $\varphi$
of the $\nbigrtilde$-module
$\gbigv_{\nbiga\star|\cnum_{\lambda}\times U}$ $(\star=\ast,!)$,
there exists a complex number $\alpha$
such that $\varphi=\alpha\id$.
\end{cor}
\pf
We obtain the claim for
$\gbigv_{\nbiga,\ast}$
by the argument in the proof of 
Corollary \ref{cor;14.12.2.21}
and Corollary \ref{cor;14.12.2.22}.
Because
$\gbigv_{\nbiga,!}\simeq
 \lambda^{2n-m+r}j^{\ast}\DDD\gbigv_{\nbiga\ast}$,
we obtain the claim for 
$\gbigv_{\nbiga,!}$.
\hfill\qed

\vspace{.1in}
We have an isomorphism
$\Psi:
 j^{\ast}\bigl(
 \gbigv_{\nbiga,\ast|S_{\nbiga}^{\reg}}
 \bigr)^{\lor}
\simeq
 \lambda^{-n-r}
 \gbigv_{\nbiga,!|S_{\nbiga}^{\reg}}$
induced by the duality isomorphism
$\gbigv_{\nbiga,!}
\simeq
\lambda^{2n-m+r}j^{\ast}\DDD\gbigv_{\nbiga,\ast}$.

\begin{cor}
Let $U\subset S_{\nbiga}^{\reg}$ be any open subset
such that $B_{\nbiga}\subset U$.
For any morphism of $\nbigrtilde$-modules
$\Psi':
 j^{\ast}\bigl(
 \gbigv_{\nbiga,\ast|\cnum_{\lambda}\times U}
 \bigr)^{\lor}
\lrarr
 \lambda^{-n-r}
 \gbigv_{\nbiga,!|\cnum_{\lambda}\times U}$,
there exists a complex number $\alpha$
such that $\Psi'=\alpha\Psi$.
\hfill\qed
\end{cor}

\subsection{Preliminary on the A-side}

\label{subsection;14.11.19.50}

We recall some basic matters on the quantum products of toric varieties.
For example,
see \cite{Iritani-convergence}, \cite{Iritani},
\cite{Konishi-Minabe}, \cite{Konishi-Minabe2},
\cite{Reichelt-Sevenheck1}, \cite{Reichelt-Sevenheck2}
for more details.

\subsubsection{The cohomology ring of projective bundles}
\label{subsection;14.9.16.1}

Let $X$ be a complex manifold.
For simplicity, we assume 
$H^{2i+1}(X,\seisuu)=0$ for any $i$.
We shall omit the coefficient of the cohomology group
if the coefficient is $\cnum$.
Let $\nbige$ be a locally free $\nbigo_X$-module
of rank $r$.
Let $\nbige^{\lor}$ denote the dual.
Let $Y:=\proj(\nbige\oplus\nbigo_X)=
 \Proj\bigl(
 \Sym^{\bullet}(\nbige\oplus\nbigo_X)\bigr)$
denote the projective completion of $\nbige^{\lor}$.
We have the natural inclusions
$i_0:X\simeq \proj(\nbigo)\lrarr Y$
and 
$i_{\infty}:
 H_{\infty}:=\proj(\nbige)\lrarr Y$.
We have
$i_0(X)\subset \nbige^{\lor}$
and 
$Y=\nbige^{\lor}\cup H_{\infty}$.

Let $\pi:Y\lrarr X$ denote the projection.
Let $\pi^{\ast}$ and $i_0^{\ast}$
denote the pull back of the cohomology.
We have
$H^{\ast}(Y)
=\Ker i_0^{\ast}\oplus
 \Image \pi^{\ast}$.
Let $i_{0\ast}:H^i(X)\lrarr H^{i+2r}(Y)$ 
and $\pi_{\ast}:H^i(Y)\lrarr H^{i-2r}(X)$
denote the Gysin map.
We have
$H^{\ast}(Y)
=\Image i_{0\ast}
\oplus
 \Ker \pi_{\ast}$.

Let $\nbigo_Y(1)$ denote 
the tautological line bundle
of $Y$ over $X$.
We set
$\gamma:=c_1\bigl(\nbigo_Y(1)\bigr)
 \in H^2(Y)$.
We have
$H^{\ast}(Y)
=\bigoplus_{j=0}^r
 \gamma^j\cdot
 \pi^{\ast}\bigl(
 H^{\ast}(X)\bigr)$.
Note that $\gamma$ is the cohomology class
representing $H_{\infty}$.
In particular, we have $i_0^{\ast}\gamma=0$.
Hence, $\Ker i_0^{\ast}=
 \bigoplus_{j=1}^{r}
 \pi^{\ast}H^{\ast}(X)\gamma^j$.
Because 
$c_{r+1}(\nbige^{\lor}\oplus \nbigo)=0$,
we have
\begin{equation}
\label{eq;14.9.17.1}
 \gamma\cdot
 \sum_{j=0}^{r}
 \gamma^{r-j}
 \pi^{\ast}c_j(\nbige^{\lor})
=0. 
\end{equation}
Let $N:H^{\ast}(Y)\lrarr H^{\ast}(Y)$
be determined by
$N(\sigma)=\gamma\sigma$.
By (\ref{eq;14.9.17.1}),
we have 
$\Image N=
 \bigoplus_{i=1}^r 
 \gamma^i\cdot
 \pi^{\ast}\bigl(H^{\ast}(X)\bigr)$,
and $\Cok N\simeq \Image\pi^{\ast}$.
\begin{lem}
We have $\Image i_{0\ast}=\Ker N$.
\end{lem}
\pf
Because $H_{\infty}\cap i_0(X)=0$,
we have $\Image (i_{0\ast})\subset\Ker N$.
By comparison of the dimension,
we obtain $\Image i_{0\ast}=\Ker N$.
\hfill\qed

\vspace{.1in}

In particular,
$\Ker(N)\cap H^{2r}(Y)
 =i_{0\ast}H^0(X,\cnum)$
is one dimensional,
and that
$\sum_{j=0}^r\gamma^{r-j}
 \pi^{\ast}c_j(\nbige^{\lor})$
is a base.
By considering the restriction to
a ball in $X$,
we obtain that 
$i_{0\ast}(1)
=\sum_{j=0}^r\gamma^{r-j}
 \pi^{\ast}c_j(\nbige^{\lor})$.
In particular, we have
$i_{0\ast}\sigma
=
 \pi^{\ast}\sigma
 \sum_{j=0}^r\gamma^{r-j}
 \pi^{\ast}c_j(\nbige^{\lor})$.

\subsubsection{Quantum products of $Y$ in the weak Fano toric case}
\label{subsection;14.11.20.1}

Suppose that 
(i) $X$ is a smooth weak Fano toric variety,
(ii) $\nbige$ is the direct sum of 
 holomorphic line bundles
 $\nbigl_i$ $(i=1,\ldots,r)$,
(iii) $\nbigl_i$ $(i=1,\ldots,r)$
 and  $(\bigotimes_{i=1}^r\nbigl_i
 \otimes K_X)^{\lor}$
 are nef.
Then, $Y$ is also a weak Fano toric manifold.

Let $\Eff^{\ast}(Y)\subset H_2(Y,\seisuu)$
denote the subset of 
the homology classes of 
non-empty algebraic curves.
We set $\Eff(Y):=\Eff^{\ast}(Y)\cup\{0\}$.
We take a homogeneous base 
$\phi_1,\ldots,\phi_{\chi(Y)}$ of $H^{\ast}(Y,\cnum)$.
Let $\phi^1,\ldots,\phi^{\chi(Y)}$
be the dual base 
of $H^{\ast}(Y)$
with respect to the Poincar\'e pairing.

For $d\in\Eff(Y)$,
let $Y_{0,3,d}$ denote the moduli stack of stable maps
$f:C\lrarr Y$
where $C$ denotes a $3$-pointed pre-stable curve
with genus $0$, and $f$ denotes a morphism such that 
the homology class of $f(C)$ is $d$.
Let $p:\nbigc_{0,3,d}\lrarr Y_{0,3,d}$ be the universal curve,
and let $\ev_i:\nbigc_{0,3,d}\lrarr X$ be the evaluation map
at the $i$-th marked point.
Let $[Y_{0,3,d}]^{\vir}$ denote
the virtual fundamental class of $Y_{0,3,d}$.
For $\alpha_i\in H^{\ast}(Y)$ $(i=1,2,3)$,
we obtain the number
$\langle
 \alpha_1,\alpha_2,\alpha_3
 \rangle^Y_{0,3,d}
:=
\int_{[Y_{0,3,d}]^{\vir}}
 \prod_{i=1,2,3}\ev_i^{\ast}\alpha_i$
called the Gromov-Witten invariants of $Y$.

According to a result of Iritani \cite{Iritani-convergence}
(see also \cite{Iritani} and \cite{Reichelt-Sevenheck1}),
there exists an open subset
 $\nbigu\subset H^{2}(Y,\cnum)$ 
such that the following holds:
\begin{itemize}
\item $\nbigu$ is of the form
 $\bigl\{
 u\in H^2(Y,\cnum)\,\big|\,
 \Re\langle u,d\rangle<-M,\,\,\,
 \forall d\in \Eff^{\ast}(Y)
 \bigr\}$ for some $M>0$.
\item
 For any $u\in\nbigu$,
 the quantum product $\bullet_u$
 on $H^{\ast}(Y,\cnum)$ is convergent:
\[
 \alpha\bullet_u\beta
:=\sum_{d\in\Eff(Y)}
\sum_{i=1}^{\chi(Y)}
 \langle
 \alpha,\beta,\phi_i
 \rangle^Y_{0,3,d}
 \cdot
 \phi^i
 \cdot
 e^{\langle u,d\rangle}
\]
\end{itemize}

We have the decomposition
$H^2(Y,\seisuu)
=\pi^{\ast}\bigl(
 H^2(X,\seisuu)
 \bigr)\oplus \seisuu\gamma$,
which induces
$H^2(Y)
=\pi^{\ast}\bigl( H^2(X)\bigr)
\oplus \cnum\gamma$.

\begin{lem}
For a large $M'>0$,
we have 
\[
  \pi^{\ast}
 \bigl\{
 \sigma\in H^2(X,\cnum)\,\big|\,
 \Re(\langle\sigma,d\rangle)
<-M',\,\,\forall d\in\Eff^{\ast}(X)
 \bigr\}
\times
 \{c\cdot \gamma\,|\,
 c \in\cnum,\,\,\Re(c)<-M'\}
\subset\nbigu.
\]
\end{lem}
\pf
Any algebraic curve $C$
is homologous to 
a curve $C_0\cup 
 \bigcup_{j=1}^r C_j
 \cup\bigcup_{i=1}^{\ell} F_i$
where $C_0$ is contained in $i_0(X)$,
$C_{j}$ are contained in $\proj(\nbigl_j)\subset Y$,
and $F_i$ $(i=1,\ldots,\ell)$
are contained in fibers of 
$\pi:Y\lrarr X$.
We have
\[
\sum_{j=0}^r
 \Re
 \langle \sigma+c\gamma,[C_j]
 \rangle
+\sum_{i=1}^{\ell}
 \Re\langle\sigma+c\gamma,[F_i]
 \rangle
\leq
 \sum_{j=0}^r
\Re \langle
 \sigma,[C_j]
 \rangle
+
 \sum_{i=1}^{\ell}
 \Re(c)\langle \gamma,[F_i]\rangle
\]
Then, the claim of the lemma is clear.
\hfill\qed

\vspace{.1in}

We set
$\nbigu_X:=
 \bigl\{
 \sigma\in H^2(X,\cnum)\,\big|\,
 \Re(\langle\sigma,d\rangle)
<-M'
 \bigr\}$,
$\nbigu_{\gamma}:=
 \bigl\{
 c\cdot\gamma\,\big|\,
 c\in\cnum,\,\,\, \Re(c)<-M'
 \bigr\}$
and $\nbigu_Y:=\nbigu_X\times\nbigu_{\gamma}$.
A point $\sigma+c\gamma\in \nbigu_Y$
is denoted by $(\sigma,c)$.
We have a natural action of
$2\pi\sqrt{-1}H^2(Y,\seisuu)$
on $\nbigu_{Y}$
by the addition.
The quantum product
$\bullet_{\sigma,c}$
depends only on 
the equivalence class of $(\sigma,c)$
in $\nbigu_Y\big/2\pi\sqrt{-1}H^2(Y,\seisuu)$.

\subsubsection{Degenerated quantum products on 
 $H^{\ast}(Y)$ and $H^{\ast}(X)$}
\label{subsection;15.5.16.30}

Let $\Eff(Y,\gamma)$ denote the set of
$[C]\in \Eff(Y)$ such that 
$\langle \gamma,[C]\rangle=0$.
We set 
$\Eff^{\ast}(Y,\gamma):=\Eff(Y,\gamma)\cap\Eff^{\ast}(Y)$.
Let $\sigma\in\nbigu_X$.
As in \cite{Konishi-Minabe2},
by taking the limit of
$\bullet_{(\sigma,c)}$ for $\Re(c)\lrarr -\infty$,
we obtain the following product
$\bullet_{\sigma}$ on $H^{\ast}(Y)$:
\[
\alpha\bullet_{\sigma}\beta
:=\sum_{d\in\Eff(Y,\gamma)}
 \sum_{i=1}^{\chi(Y)}
 \langle
 \alpha,\beta,\phi_i
 \rangle^{Y}_{0,3,d}
 \cdot
 \phi^i
 \cdot
 e^{\langle \sigma,\pi_{\ast}d\rangle}
\]
As in \S\ref{subsection;14.9.16.1},
let $N:H^{\ast}(Y)\lrarr H^{\ast}(Y)$
be the endomorphism given by
the cup product of $\gamma$,
i.e., $N(\beta):=\gamma\cup\beta$ 
for $\beta\in H^{\ast}(Y)$.
Note that
$\gamma\bullet_{\sigma}\beta
=\beta\bullet_{\sigma}\gamma=
\beta\cup\gamma$
for any $\beta\in H^{\ast}(Y)$,
which follows from the divisor axiom.
Hence, 
$\Image(N)$ and $\Ker(N)$
are the ideal of 
the algebra $(H^{\ast}(Y),\bullet_{\sigma})$.
By the natural isomorphism of 
vector spaces $H^{\ast}(X)\simeq \Cok(N)$,
we obtain the induced product 
$\bullet^{\nbige^{\lor}}_{\sigma}$
on $H^{\ast}(X)$.
The algebra $(H^{\ast}(X),\bullet^{\nbige^{\lor}}_{\sigma})$
is denoted by
$Q_{\sigma}H^{\ast}(X,\nbige^{\lor})$ in this paper.
The multiplication
$\Cok(N)\times \Ker(N)\lrarr \Ker(N)$
and the natural isomorphism
$\Ker(N)\simeq H^{\ast}(X)$
induces a structure of
$Q_{\sigma}H^{\ast}(X,\nbige^{\lor})$-module
on $H^{\ast}(X)$.
The $Q_{\sigma}H^{\ast}(X,\nbige^{\lor})$-module
is denoted by 
$K_{\sigma}H^{\ast}(X,\nbige^{\lor})$.
The multiplication of
$b_1\in Q_{\sigma}H^{\ast}(X,\nbige^{\lor})$
and $b_2\in K_{\sigma}H^{\ast}(X,\nbige^{\lor})$
is denoted by
$b_1\bullet^{\nbige^{\lor}}_{K,\sigma}b_2$.

\subsubsection{Filtrations and the graded pairings}
\label{subsection;15.5.16.41}

Let $W(N)$ be the weight filtration of 
the nilpotent map $N$ on $H^{\ast}(Y)$.
For any $\sigma\in\nbigu_X$,
by shifting the filtrations 
of Konishi-Minabe in \cite{Konishi-Minabe2},
we set
\[
 \Wtilde_{k}Q_{\sigma}H^{\ast}(X,\nbige^{\lor}):=
 \Image\Bigl(
 W(N)_{k+n+r}H^{\ast}(Y)\lrarr
 Q_{\sigma}H^{\ast}(X,\nbige^{\lor})
 \Bigr).
\]
Here, $Q_{\sigma}H^{\ast}(X,\nbige^{\lor})$
is identified with $\Cok(N)$.
Note that
\[
\Gr^{\Wtilde}_{k-n-r}Q_{\sigma}H^{\ast}(X,\nbige^{\lor})
\simeq
 \left\{
 \begin{array}{ll}
 P\Gr^{W(N)}_kH^{\ast}(Y) & (k\geq 0) \\
 0 & (k<0) 
 \end{array}
 \right.
\]
Here, 
$P\Gr^{W(N)}_kH^{\ast}(Y)$ is the primitive part,
i.e., the kernel of
$N^{k+1}:\Gr^{W(N)}_kH^{\ast}(Y)
\lrarr \Gr^{W(N)}_{-k-2}H^{\ast}(Y)$.
The Poincar\'e pairing $(\cdot,\cdot)_Y$ on $H^{\ast}(Y)$
and the nilpotent map $N$ induces the following symmetric pairing
on $\Gr^{\Wtilde}_{k-n-r}Q_{\sigma}H^{\ast}(X,\nbige^{\lor})$:
\[
 P^{\nbige^{\lor}}_{k-n-r}
 \bigl(a,b \bigr):=
 \bigl(
 a,N^{k}b
 \bigr)_Y
\]
It is easy and standard to see that the pairings are non-degenerate.

We also have the filtration $W$
on $K_{\sigma}H^{\ast}(X,\nbige^{\lor})$
given as follows:
\[
 \Wtilde_kK_{\sigma}H^{\ast}(X,\nbige^{\lor}):=
 K_{\sigma}H^{\ast}(X,\nbige^{\lor})
\cap
 W(N)_{k+n+r}H^{\ast}(Y).
\]
Here 
$K_{\sigma}H^{\ast}(X,\nbige^{\lor})$
is identified with $\Ker(N)$.
We have
\[
 \Gr^{\Wtilde}_{k-n-r}K_{\sigma}H^{\ast}(X,\nbige^{\lor})
\simeq
 \left\{
 \begin{array}{ll}
 P'\Gr_{k}^{W(N)}H^{\ast}(Y) & (k\leq 0),\\
 0 & (k>0).
 \end{array}
 \right.
\]
Here,
$P'\Gr_{k}^{W(N)}H^{\ast}(Y)$
is the image of 
$P\Gr_{-k}^{W(N)}H^{\ast}(Y)$
by $N^{-k}$.
Hence, we have the induced pairings
$\bigl(\cdot,\cdot\bigr)^{\nbige^{\lor}}_{k-n-r}$
on 
$\Gr^{\Wtilde}_{k-n-r}K_{\sigma}H^{\ast}(X,\nbige^{\lor})$.

\vspace{.1in}

Note that the filtrations $\Wtilde$
and the symmetric pairings 
$P^{\nbige^{\lor}}_k$ are independent of
the choice of $\sigma$.
Namely, we have the filtrations $\Wtilde$ 
on $\Cok(N)$ and $\Ker(N)$
such that
$\Wtilde_k\Cok(N)=
 \Wtilde_kQ_{\sigma}H^{\ast}(X,\nbige^{\lor})$
and 
$\Wtilde_k\Ker(N)=
 \Wtilde_kK_{\sigma}H^{\ast}(X,\nbige^{\lor})$.
We also have symmetric pairings
$P^{\nbige^{\lor}}_k$
on $\Gr^{\Wtilde}_k\Cok(N)$
and $\Gr^{\Wtilde}_k\Ker(N)$
which are equal to the pairings on
$\Gr^{\Wtilde}_kQ_{\sigma}H^{\ast}(X,\nbige^{\lor})$
and 
$\Gr^{\Wtilde}_kK_{\sigma}H^{\ast}(X,\nbige^{\lor})$
for any $\sigma$.

\subsubsection{
Degenerated quantum products
and local Gromov-Witten invariants
(Appendix)}
\label{subsection;15.5.16.50}

Let us recall the relation between
the degenerated quantum products
and local Gromov-Witten invariants
in some special cases.
We essentially follow \cite{Konishi-Minabe2}.

Let $d\in \Eff^{\ast}(X)$.
Let $X_{0,3,d}$ be the moduli stack of 
stable maps $f:C\lrarr X$
where $C$ denotes a $3$-pointed pre-stable curve
with genus $0$, and $f$ denotes a morphism such that
the homology class of $f(C)$ is $d$.
Let $p:\nbigc_{0,3,d}\lrarr X_{0,3,d}$ be the universal curve,
and let $\mu:\nbigc_{0,3,d}\lrarr X$ be the universal map.
Let us consider the following concavity condition for $d$:
\begin{description}
\item[(B)]
 $H^0(C,f^{\ast}\nbige^{\lor})=0$ for any $(C,f)\in X_{0,3,d}$.
\end{description}

Recall that $\nbige$ is the direct sum of the nef line bundles
$\nbigl_i$ $(i=1,\ldots,r)$.
We recall the following standard lemma.
\begin{lem}
Suppose that $\langle c_1(\nbigl_i^{\lor}),d\rangle<0$
for any $i$.
Then, Condition {\rm(B)} holds for $d$.
\end{lem}
\pf
Take $d\in \Eff^{\ast}(X)$
and  $(C,f)\in X_{0,3,d}$.
Let $C=\bigcup C_a$ be the irreducible decomposition.
Because $\nbigl_i$ is nef,
the degree of $f^{\ast}(\nbigl_i^{\lor})_{|C_a}$
is non-positive for any $a$.
By the assumption
$\langle c_1(\nbigl_i^{\lor},d\rangle)<0$,
we have an irreducible component $C_{a_0}$ of $C$
such that 
the degree of $f^{\ast}(\nbigl_i^{\lor})_{|C_a}$
is strictly negative.
We also remark that $C$ is connected.
Then, we have $H^0(C,f^{\ast}\nbigl_i)=0$,
i.e., Condition (B) holds for $d$.
\hfill\qed

\vspace{.1in}

If Condition (B) holds for $d$,
we obtain the locally free sheaf
$R^1p_{\ast}\mu^{\ast}\nbige^{\lor}$
on $X_{0,3,d}$.
Let $[X_{0,3,d}]^{\vir}$ denote the virtual fundamental class
of $X_{0,3,d}$.
Let $\ev_i:X_{0,3,d}\lrarr X$ be the evaluation map
at the $i$-th marked point.
For any vector bundle $E$,
let $e(E)$ denote the Euler class.
Then, for any $\alpha_i\in H^{\ast}(X)$ $(i=1,2,3)$,
we set
\[
 \bigl\langle
 \alpha_1,\alpha_2,\alpha_3
 \bigr\rangle^{X,\nbige^{\lor},e^{-1}}_{0,3,d}
 :=
 \int_{[X_{0,3,d}]^{\vir}}
 \prod_{i=1,2,3}\ev_i^{\ast}(\alpha_i)
\cdot
 e(Rp_{\ast}\mu^{\ast}\nbige^{\lor}).
\]
They are called 
$e^{-1}$-twisted Gromov-Witten invariants
or local Gromov-Witten invariants.

For any $\alpha,\beta\in H^{\ast}(X)$,
let $\alpha\cup\beta$ denote the ordinary cup product.
We take a homogeneous base $\rho_1,\ldots,\rho_{\ell}$
of $H^{\ast}(X)$.
Let $\rho^{1},\ldots,\rho^{\ell}$ denote the dual base of $H^{\ast}(X)$
with respect to the Poincar\'e pairing of $X$.

\begin{prop}
\label{prop;14.11.19.11}
Suppose one of the following holds:
\begin{description}
\item[Case 1.]
$X$ is Fano, and $\nbigl_i$ $(i=1,\ldots,r)$ are ample.
Note that Condition {\rm(B)} holds
for any $d\in\Eff^{\ast}(X)$.
\item[Case 2.]
$X$ is a weak Fano surface, $r=1$,
and $\nbigl_1=K_X^{-1}$.
Note that Condition {\rm(B)} holds
for any $d\in\Eff^{\ast}(X)$
with $\langle c_1(K_X^{\lor}),d\rangle\neq 0$.
\end{description}
Then, we have the following 
for $\sigma\in \nbigu_X$
and for $\alpha,\beta\in Q_{\sigma}H^{\ast}(X,\nbige^{\lor})$:
\begin{equation}
 \label{eq;14.11.19.20}
 \alpha\bullet^{\nbige^{\lor}}_{\sigma}\beta
=
\left\{
\begin{array}{ll}
{\displaystyle
 \alpha\cup\beta
+\sum_{d\in\Eff^{\ast}(X)}
 \sum_{i=1}^{\ell}
 \bigl\langle
 \alpha,\beta,e(\nbige^{\lor})\rho_i
 \bigr\rangle^{X,\nbige^{\lor},e^{-1}}_{0,3,d}
 e^{\langle \sigma,d \rangle}\rho^i
 }
 &
 (\mbox{\rm Case 1})
 \\
{\displaystyle
 \alpha\cup\beta
+\sum_{
 \substack{d\in\Eff^{\ast}(X),\\
 \langle c_1(K_X),d\rangle\neq 0}}
 \sum_{i=1}^{\ell}
 \bigl\langle
 \alpha,\beta,e(\nbige^{\lor})\rho_i
 \bigr\rangle^{X,\nbige^{\lor},e^{-1}}_{0,3,d}
 e^{\langle \sigma,d \rangle}\rho^i
 }
 &
 (\mbox{\rm Case 2})
\end{array}
\right.
\end{equation}
We also have the following
for $\sigma\in\nbigu_X$
and 
for $\alpha\in Q_{\sigma}H^{\ast}(X,\nbige^{\lor})$,
and $\beta\in K_{\sigma}H^{\ast}(X,\nbige^{\lor})$:
\begin{equation}
 \label{eq;14.11.19.21}
 \alpha\bullet^{\nbige^{\lor}}_{K,\sigma}\beta
=
\left\{
 \begin{array}{ll}
{\displaystyle
 \alpha\cup\beta
+\sum_{d\in\Eff^{\ast}(X)}
 \sum_{i=1}^{\ell}
 \bigl\langle
 \alpha,e(\nbige^{\lor})\beta,\rho_i
 \bigr\rangle^{X,\nbige^{\lor},e^{-1}}_{0,3,d}
 e^{\langle \sigma,d \rangle}\rho^i
 }
 & (\mbox{\rm Case 1})\\
 {\displaystyle
  \alpha\cup\beta
+\sum_{\substack{d\in\Eff^{\ast}(X)\\
 \langle c_1(K_X),d\rangle\neq 0
 \\}}
 \sum_{i=1}^{\ell}
 \bigl\langle
 \alpha,e(\nbige^{\lor})\beta,\rho_i
 \bigr\rangle^{X,\nbige^{\lor},e^{-1}}_{0,3,d}
 e^{\langle \sigma,d \rangle}\rho^i
 }
 & (\mbox{\rm Case 2})
 \end{array}
\right.
\end{equation}
\end{prop}
\pf
For the proof of the proposition,
we study the degenerated quantum products
$\bullet_{\sigma}$ on $H^{\ast}(Y)$.
Any element $\alpha\in H^{\ast}(Y)$
has the expression
$\alpha=\sum_{j=0}^r\pi^{\ast}\alpha_j\cdot\gamma^j$
where $\alpha_j\in H^{\ast}(X)$.

\begin{lem}
\label{lem;14.11.19.10}
Let $\sigma\in\nbigu_X$.
Suppose Case {\rm 1} or Case {\rm 2}.
For $\alpha=\sum_{j=0}^r \pi^{\ast}(\alpha_j)\gamma^j$
and $\beta=\sum_{j=0}^r \pi^{\ast}(\beta_j)\gamma^j$,
we have the following:
\begin{equation}
 \label{eq;14.11.19.12}
 \alpha\bullet_{\sigma}\beta
=
\left\{
 \begin{array}{ll}
 {\displaystyle
\alpha\cup\beta
+\sum_{d\in\Eff^{\ast}(X)}
 \sum_{i=1}^{\ell}
 \langle
 \alpha_0,
 \beta_0,\rho_i
 \rangle^{X,\nbige^{\lor},e^{-1}}_{0,3,d}
 e^{\langle \sigma,d\rangle}
 i_{0\ast}(\rho^i)}
 & \mbox{\rm(Case 1)}
 \\
 {\displaystyle
\alpha\cup\beta
+\sum_{\substack{d\in\Eff^{\ast}(X)\\
 \langle c_1(K_X),d\rangle\neq 0
 }}
 \sum_{i=1}^{\ell}
 \langle
 \alpha_0,
 \beta_0,\rho_i
 \rangle^{X,\nbige^{\lor},e^{-1}}_{0,3,d}
 e^{\langle \sigma,d\rangle}
 i_{0\ast}(\rho^i)
 }
 & \mbox{\rm(Case 2)}
\end{array}
\right.
\end{equation}
\end{lem}

Before giving a proof of Lemma \ref{lem;14.11.19.10},
we deduce Proposition \ref{prop;14.11.19.11}
from Lemma \ref{lem;14.11.19.10}.
We give an argument in Case 1.
The other case can be given similarly.
Because $i_0^{\ast}i_{0\ast}(\rho^j)=e(\nbige^{\lor})\rho^j$,
we obtain the following for $\alpha,\beta\in H^{\ast}(X)$
from (\ref{eq;14.11.19.12}):
\begin{multline}
 \alpha\bullet^{\nbige^{\lor}}_{\sigma}\beta
=\alpha\cup\beta
+\sum_{d\in\Eff^{\ast}(X)}\sum_{i=1}^{\ell}
 \bigl\langle
 \alpha,\beta,\rho_i
 \bigr\rangle^{X,\nbige^{\lor},e^{-1}}_{0,3,d}
 e^{\langle\sigma,d\rangle}
 e(\nbige^{\lor})\rho^i
 \\
=\alpha\cup\beta
+\sum_{d\in\Eff^{\ast}(X)}\sum_{i=1}^{\ell}
 \bigl\langle
 \alpha,\beta, e(\nbige^{\lor})\rho_i
 \bigr\rangle^{X,\nbige^{\lor},e^{-1}}_{0,3,d}
 e^{\langle\sigma,d\rangle}
 \rho^i
\end{multline}
Thus, we obtain (\ref{eq;14.11.19.20}).
We identify
$\Ker (N)=i_{0\ast}H^{\ast}(X)$.
We have
$i_{0\ast}(\beta)=
 e(\nbige^{\lor})\beta
+\sum_{j=1}^{r}c_j\gamma^j$
for some $c_j\in H^{\ast}(X)$.
Hence
we have the following
for $\alpha\in H^{\ast}(X)$
and $\beta\in H^{\ast}(X)$:
\begin{multline}
 i_{0\ast}\bigl(
 \alpha\bullet^{\nbige^{\lor}}_{K,\sigma}
 \beta\bigr)
=\pi^{\ast}(\alpha)\bullet_{\sigma}
 i_{0\ast}(\beta)
=\pi^{\ast}(\alpha)\cup i_{0\ast}\beta
+\sum_{d\in \Eff^{\ast}(X)}
 \sum_{i=1}^{\ell}
 \bigl\langle
 \alpha,e(\nbige^{\lor})\beta,\rho_i
 \bigr\rangle^{X,\nbige^{\lor},e^{-1}}_{0,3,d}
 e^{\langle\sigma,d\rangle}
 i_{0\ast}\rho^i
 \\
=i_{0\ast}\Bigl(
 \alpha\cup\beta
+\sum_{d\in \Eff^{\ast}(X)}
 \sum_{i=1}^{\ell}
 \bigl\langle
 \alpha,e(\nbige^{\lor})\beta,\rho_i
 \bigr\rangle^{X,\nbige^{\lor},e^{-1}}_{0,3,d}
 e^{\langle\sigma,d\rangle}
 \rho^i
 \Bigr)
\end{multline}
Thus, we obtain (\ref{eq;14.11.19.21}).
It remains to prove
Lemma \ref{lem;14.11.19.10}.

\vspace{.1in}

Let us prove Lemma \ref{lem;14.11.19.10} in Case 1.
For any algebraic curve $C$ in $Y$,
let $[C]\in H_2(Y)$ denote the homology class of $C$.

\begin{lem}
\label{lem;14.11.19.30}
Let $C$ be a non-empty algebraic curve in $Y$.
Then, 
we have $\langle \gamma,[C]\rangle=0$
if and only if  $C\subset i_0(X)$.
\end{lem}
\pf
Because $\nbigo_Y(1)=\nbigo(H_{\infty})$,
the second condition implies the first.
Let us prove that the first condition implies the second.
We may assume that $C$ is irreducible.
If $C\subset H_{\infty}$,
we clearly have $\langle \gamma,[C]\rangle>0$.
Hence, we have $C\not\subset H_{\infty}$.
If $C\cap H_{\infty}\neq\emptyset$,
we have $\langle \gamma,[C]\rangle>0$.
Hence, we obtain that $C\cap H_{\infty}=\emptyset$.
Let $\varphi:\Ctilde\lrarr C\subset\nbige^{\lor}$
be the normalization.
Let $q:\nbige^{\lor}\lrarr X$ be the projection.
Note that $q(C)$ is not a point.
We obtain the morphism
$q\circ\varphi:\Ctilde\lrarr X$
and a section $s$ of
$(q\circ\varphi)^{\ast}\nbige^{\lor}$.
Because $(q\circ\varphi)^{\ast}\nbigl_i$
are ample,
we have $s=0$.
It means $C\subset i_0(X)$.
Thus, we obtain Lemma \ref{lem;14.11.19.30}.
\hfill\qed

\begin{lem}
\label{lem;14.11.19.31}
For $d\in\Eff^{\ast}(Y,\gamma)$
and for $i>0$,
we have
$\langle
 \alpha_1,\alpha_2,\gamma^i\alpha_3
 \rangle^Y_{0,3,d}=0$.
\end{lem}
\pf
The support of the cohomology class
$\gamma^i\alpha_3$
is contained in $H_{\infty}$.
For $(C,f)\in Y_{0,3,d}$,
we have $f(C)\subset i_{0}(X)$.
Then, the claim of Lemma \ref{lem;14.11.19.31}
is clear.
\hfill\qed

\vspace{.1in}

By Lemma \ref{lem;14.11.19.30},
we have
$\Eff^{\ast}(Y,\gamma)=\Eff^{\ast}(X)$.
By Lemma \ref{lem;14.11.19.31},
for $d\in \Eff^{\ast}(Y,\gamma)$
and for $\alpha_i=\sum_{j=0}^{r} 
 \pi^{\ast}\alpha_{i,j}\gamma^j$,
we have
\begin{multline}
 \langle
 \alpha_1,\alpha_2,\alpha_3
 \rangle^{Y}_{0,3,d}
=\langle
 \pi^{\ast}\alpha_{1,0},
 \pi^{\ast}\alpha_{2,0},
 \pi^{\ast}\alpha_{3,0}
 \rangle^{Y}_{0,3,d}
=
 \int_{[X_{0,3,d}]^{\vir}}
 \prod \ev_i^{\ast}(\alpha_{i,0})
 \cdot e\bigl(R^1p_{\ast}\mu^{\ast}\nbige^{\lor}\bigr)
 \\
=
 \langle
 \alpha_{1,0},\alpha_{2,0},
 \alpha_{3,0}
 \rangle^{X,\nbige^{\lor},e^{-1}}_{0,3,d}
\end{multline}

We set 
$\phi_{ij}:=\rho_i\gamma^j$
$(1\leq i\leq \ell,\,\,
 0\leq j\leq r)$
which give a frame of $H^{\ast}(Y)$.
Let $\phi^{ij}$ 
$(1\leq i\leq \ell,\,\,
 0\leq j\leq r)$
be the dual base
with respect to the Poincar\'e pairing of $Y$.
Note that 
$\phi^{i0}=i_{0\ast}(\rho^i)$.
Indeed, it is enough to check
$\langle
 \phi_{kj},
 i_{0\ast}(\rho^i)
 \rangle=0$ $(j>0)$,
$\langle
 \phi_{k0},
 i_{0\ast}(\rho^i)
 \rangle=0$ $(k\neq i)$
and 
$\langle
 \phi_{k0},
 i_{0\ast}(\rho^i)
 \rangle=1$,
which can be checked by direct computations.
Then, we obtain the equality (\ref{eq;14.11.19.12})
in Case 1.

\vspace{.1in}

Let us prove (\ref{eq;14.11.19.12}) in Case 2.
It is exactly the case given in \cite{Konishi-Minabe2}.
We just revisit it
in a slightly different way.
Let $d\in\Eff(Y,\gamma)$.
Note that 
 the virtual dimension of
 $Y_{0,3,d}$ is $3$.
Indeed, because we have
$K_Y=\nbigo_Y(2)$,
the virtual dimension of
$Y_{0,n,d}$ is
$n-3+\dim Y-\langle c_1(K_Y),d\rangle=n$.

\begin{lem}
\label{lem;14.11.19.40}
For 
$\alpha_i=\alpha_{i0}+\alpha_{i1}\gamma$
$(\alpha_{ij}\in H^{\ast}(X),\,i=1,2,3)$,
we have
\[
 \langle
 \alpha_1,\alpha_2,\alpha_3
 \rangle^{Y}_{0,3,d}
=\langle
 \alpha_{10},\alpha_{20},\alpha_{30}
 \rangle^Y_{0,3,d}.
\]
\end{lem}
\pf
By the dimension reason,
we have 
$\langle\alpha_1,\alpha_2,\alpha_3
 \rangle^{Y}_{0,3,d}=0$
unless the cohomological degree of $\alpha_i$
are $2$.
Moreover,  by the divisor axiom,
we have
$\langle\alpha_1,\alpha_2,\alpha_3
 \rangle^{Y}_{0,3,d}=0$
if one of $\alpha_i$ is $\gamma$.
Then, the claim of Lemma \ref{lem;14.11.19.40} follows.
\hfill\qed

\begin{lem}
\label{lem;14.11.19.41}
Suppose that 
$d\in\Eff^{\ast}(Y,\gamma)$
satisfies $\langle \pi^{\ast}c_1(K_X),d\rangle=0$.
Then, we have 
$\langle \alpha_1,\alpha_2,\alpha_3\rangle_{0,3,d}^Y=0$
for any $\alpha_i\in H^{\ast}(Y)$.
\end{lem}
\pf
It is enough to consider the case where
$\alpha_i\in \pi^{\ast}H^2(X)$.
Let us observe that
the induced morphism
$Y_{0,3,d}\lrarr X_{0,3,d}$ is smooth
and the fibers are $\proj^1$.
Take $(C,g)\in X_{0,3,d}$.
Because $K_X^{\lor}$ is nef,
and because we have $\langle c_1(K_X),[g(C)]\rangle=0$,
the restriction of $g^{\ast}K_X$ 
to each irreducible component of $C$
is of degree $0$.
In particular, we have
$H^1(C,g^{\ast}K_X)=0$,
and 
$\dim H^0(C,g^{\ast}K_X)=1$.
Let $\nbigc_{0,3,d}$ denote the universal curve
over $X_{0,3,d}$.
Let $\mu:\nbigc_{0,3,d}\lrarr X$ be the universal morphism.
Let $p:\nbigc_{0,3,d}\lrarr X_{0,3,d}$ 
denote the projection.
Then, $p_{\ast}\mu^{\ast}(\nbigo\oplus K_X^{\lor})$
is a locally free sheaf,
and the projectivization is
naturally isomorphic to $Y_{0,3,d}$.
Then, it is easy to see that
$[Y_{0,3,d}]^{\vir}$
is the pull back of 
$[X_{0,3,d}]^{\vir}$
via the natural morphism
$q:Y_{0,3,d}\lrarr X_{0,3,d}$.
Because 
$q_{\ast}([Y_{0,3,d}]^{\vir})=0$,
we obtain 
$\langle \alpha_1,\alpha_2,\alpha_3\rangle_{0,3,d}^Y=0$.
Thus, we obtain Lemma \ref{lem;14.11.19.41}.
\hfill\qed

\begin{lem}
\label{lem;14.11.19.42}
Let $d\in \Eff(Y,\gamma)$.
 If $\langle \pi^{\ast}c_1(K_X),d\rangle< 0$,
 then  we have $f(C)\subset i_0(X)$
for any $(C,f)\in Y_{0,3,d}$.
 In particular,
 $d$ comes from $\Eff^{\ast}(X)$.
\end{lem}
\pf
Let $(C,f)\in Y_{0,3,d}$.
We have $f(C)=\sum m_i[C_i]$ $(m_i>0)$
for algebraic curves in $Y$.
Note that $C$ is connected.
Suppose 
$f(C)\cap H_{\infty}\neq \emptyset$.
If $f(C)\not\subset H_{\infty}$,
we have
$\langle \gamma,[C_i]\rangle\geq 0$
for any $i$,
and 
$\langle\gamma,[C_{i_0}]\rangle>0$
for some $i_0$,
and hence
$\langle \gamma,[f(C)]\rangle>0$
which contradicts with the assumption.
If $f(C)\subset H_{\infty}$,
we have
$\langle
 \gamma,f(C)
 \rangle
=\langle
 c_1(K_X^{\lor}),
 d\rangle>0$
which also contradicts with the assumption.
Hence, we have $f(C)\cap H_{\infty}=\emptyset$.
Then, $f$ is equivalent to a morphism
$g:C\lrarr X$
and a section
$s$ of $g^{\ast}(K_X)$.
But, because $\langle c_1(K_X),d\rangle<0$,
we have $s=0$.
Thus, we obtain
Lemma \ref{lem;14.11.19.42}.
\hfill\qed

\vspace{.1in}
We can deduce the equality (\ref{eq;14.11.19.12}) in Case 2
from Lemma \ref{lem;14.11.19.40},
Lemma \ref{lem;14.11.19.41}
and Lemma \ref{lem;14.11.19.42}
as in Case 1.
Thus, Lemma \ref{lem;14.11.19.10}
and Proposition \ref{prop;14.11.19.11} are proved.
\hfill\qed

\subsection{Associated quantum $\nbigd$-modules}
\label{subsection;14.12.8.10}

We continue to use the notation in \S\ref{subsection;14.11.19.50}.
We recall some basic matters on the quantum $\nbigd$-modules
of toric varieties.
Again, see \cite{Iritani-convergence}, \cite{Iritani},
\cite{Konishi-Minabe}, \cite{Konishi-Minabe2},
\cite{Reichelt-Sevenheck1}, \cite{Reichelt-Sevenheck2}
for more details.

\subsubsection{Quantum $\nbigd$-modules
associated to degenerated quantum products on $H^{\ast}(X)$}

We set 
$\QDM(X,\nbige^{\lor}):=
 \nbigo_{\cnum_{\lambda}\times\nbigu_X}
 \otimes H^{\ast}(X)$.
As in \cite{Konishi-Minabe2},
we have the meromorphic connection $\nabla^{\nbige^{\lor}}$
on $\QDM(X,\nbige^{\lor})$
which is a variant of the connections of Dubrovin and Givental.
It is given as follows.
Let $\mu_X$ be the grading operator on $H^{\ast}(X)$
defined by $\mu(a):=ka$ for $a\in H^{2k}(X)$.
We can naturally identify
$H^{2}(X)\otimes\nbigo_{\nbigu_X}$
with the tangent sheaf of $\nbigu_X$,
and so we can regard
$\xi\in H^{2}(X)$ as a vector field on $\nbigu_X$.
We can naturally regard
$b\in H^{\ast}(X)$ as a section of $\QDM(X,\nbige^{\lor})$.
Then, 
$\nabla^{\nbige^{\lor}}$ is determined by the following 
on $\cnum^{\ast}\times \nbigu_X$:
\[
 \bigl(
\nabla^{\nbige^{\lor}}_{\xi}b
 \bigr)_{(\lambda,\sigma)}
 =-\lambda^{-1}\xi\bullet^{\nbige^{\lor}}_{\sigma}b,
\quad
 \bigl(
 \nabla^{\nbige^{\lor}}_{\lambda\del_{\lambda}}b
 \bigr)_{(\lambda,\sigma)}
=\lambda^{-1}E\bullet^{\nbige^{\lor}}_{\sigma}b
+\mu_X(b),
\quad\quad
((\lambda,\sigma)\in \cnum^{\ast}\times \nbigu_X).
\]
The filtrations $\Wtilde$ of 
$Q^{\nbige^{\lor}}_{\sigma}H^{\ast}(X)$ $(\sigma\in\nbigu_X)$
give a filtration $\Wtilde$ on $\QDM(X,\nbige^{\lor})$.
It is preserved by the connection $\nabla^{\nbige^{\lor}}$.
We have the induced meromorphic connection
on $\Gr^{\Wtilde}_k\QDM(X,\nbige^{\lor})$,
denoted by the same notation.

We regard $b_i\in \Gr^{\Wtilde}_{k}\Cok(N)$ $(i=1,2)$
as sections of 
$\Gr^{\Wtilde}_k\QDM(X,\nbige^{\lor})$.
We set
\[
 \nbigp^{\nbige^{\lor}}_k(b_1,j^{\ast}b_2)
:=\lambda^{-k}P^{\nbige^{\lor}}_k(b_1,b_2).
\]
Thus, we obtain a morphism
$\nbigp^{\nbige^{\lor}}_k:
\Gr^{\Wtilde}_k\QDM(X,\nbige^{\lor})
\otimes
j^{\ast}
\Gr^{\Wtilde}_k\QDM(X,\nbige^{\lor})
\lrarr
\lambda^{-k}
\nbigo_{\cnum\times\nbigu_X}$.

We will give a proof of the following lemma
by using the description
as a specialization of the quantum $\nbigd$-module of $Y$
in \S\ref{subsection;14.11.28.100}.
\begin{lem}
\label{lem;14.11.28.200}
$\nabla^{\nbige^{\lor}}$ is flat,
and 
$\nbigp^{\nbige^{\lor}}_k$ are  pairings of weight $k$
on 
$\bigl(
 \Gr^{\Wtilde}_k\QDM(X,\nbige^{\lor}),\nabla^{\nbige^{\lor}}
 \bigr)$.
In other words,
$\bigl(
 \QDM(X,\nbige^{\lor}),\Wtilde,
\bigl\{
\nbigp^{\nbige^{\lor}}_k
 \bigr\} 
 \bigr)$
is a mixed TEP-structure.
\end{lem}

We also set 
$\nbigk(X,\nbige^{\lor}):=
 \nbigo_{\cnum_{\lambda}\times\nbigu_X}
 \otimes
 H^{\ast}(X)$.
As in the case of $\QDM(X,\nbige^{\lor})$,
it is equipped with the meromorphic flat connection
$\nabla^{\nbige^{\lor}}$
satisfying the following:
\[
 \bigl(
 \nabla^{\nbige^{\lor}}_{\xi}b
 \bigr)_{(\lambda,\sigma)}
=-\lambda^{-1}\xi\bullet^{\nbige^{\lor}}_{K,\sigma}b,
\quad
 \bigl(
 \nabla^{\nbige^{\lor}}_{\lambda\del_{\lambda}}b
 \bigr)_{(\lambda,\sigma)}
=\lambda^{-1}E\bullet^{\nbige^{\lor}}_{K,\sigma}b
+\mu_X(b)+rb,
\quad
\quad
((\lambda,\sigma)\in \cnum^{\ast}\times\nbigu_X)
\]
Here, $b\in H^{\ast}(X)$ and $\xi\in H^{2}(X)$.
The filtration $\Wtilde$ is preserved by 
$\nabla^{\nbige^{\lor}}$.
We obtain a morphism
$\nbigp_k^{\nbige^{\lor}}:
 \Gr^{\Wtilde}_k\nbigk(X,\nbige^{\lor})
\times
 j^{\ast}
 \Gr^{\Wtilde}_k\nbigk(X,\nbige^{\lor})
\lrarr
 \lambda^{-k}\nbigo_{\cnum\times\nbigu_X}$
by the formula:
\[
 \nbigp^{\nbige^{\lor}}_k(b_1,j^{\ast}b_2)
=(-1)^{k-n-r}\lambda^{-k}
 P^{\nbige^{\lor}}_k(b_1,b_2)
\]

We will give a proof of the following lemma
by using the description
as a specialization of the quantum $\nbigd$-module of $Y$
in \S\ref{subsection;14.11.28.100}.
\begin{lem}
\label{lem;14.11.29.30}
$\nabla^{\nbige^{\lor}}$ is flat,
and $\nbigp_k^{\nbige^{\lor}}$
are pairings of weight $k$
on 
$\bigl(
 \Gr^{\Wtilde}_k\nbigk(X,\nbige^{\lor}),
 \nabla^{\nbige^{\lor}}
 \bigr)$.
In other words,
$\bigl(
 \nbigk(X,\nbige^{\lor}),\Wtilde,
 \{\nbigp_k^{\nbige^{\lor}}\}
 \bigr)$
is a mixed TEP-structure.
\end{lem}

\subsubsection{Quantum $\nbigd$-module of $Y$}

We recall the quantum $\nbigd$-module 
associated to the quantum products of $Y$.
We set 
 $\nbigu_{\gamma 1}:=
 \nbigu_{\gamma}/2\pi\sqrt{-1}\seisuu\gamma$.
We set 
 $\nbigu_{Y1}:=\nbigu_X\times\nbigu_{\gamma 1}$
which is the quotient of $\nbigu_{Y}$
by the action of $2\pi\sqrt{-1}\seisuu\gamma$.
We embed $\nbigu_{\gamma 1}\lrarr \cnum^{\ast}$
induced by $c\longmapsto \exp(c)$.
The quantum products $\bullet_{\sigma,c}$
depends only on $(\sigma,e^c)$.
So, we denote them by $\bullet_{\sigma,e^c}$.

We set 
 $\QDM(Y)
 :=\nbigo_{\cnum_{\lambda}\times\nbigu_{Y1}}\otimes H^{\ast}(Y)$.
Recall that we have the meromorphic flat connection $\nabla^G$
on $\QDM(Y)$ due to Dubrovin and Givental.
Let $\mu_Y$ be the grading operator on $H^{\ast}(Y)$,
given by $\mu_Y(b):=kb$ for $b\in H^{2k}(Y)$.
We naturally regard
$\xi\in H^{2}(Y)$ as a vector field on $\nbigu_Y$,
which induces a vector field on $\nbigu_{Y1}$.
We also naturally regard $b\in H^{\ast}(Y)$ as a section of 
$\QDM(Y)$.
Then, $\nabla^G$ is determined by the following on 
$\cnum_{\lambda}^{\ast}\times\nbigu_{Y1}$:
\[
 \bigl(
 \nabla^G_{\xi}b
 \bigr)_{(\lambda,\sigma,t)}
=-\lambda^{-1}\xi\bullet_{(\sigma,t)}b,
\quad\quad
 \bigl(
 \nabla^G_{\lambda\del_{\lambda}}b
 \bigr)_{(\lambda,\sigma,t)}
=\lambda^{-1}E\bullet_{(\sigma,t)}b
+\mu_Y(b),
\quad\quad
 ((\lambda,\sigma,t)\in\cnum^{\ast}\times\nbigu_{Y1}).
\]
It is equipped with the pairing $\nbigp_Y$ of weight $-n-r$
induced by the Poincar\'e pairing:
\begin{equation}
 \label{eq;14.11.20.2}
 \nbigp_Y(b_1,j^{\ast}b_2):=\lambda^{n+r}(b_1,b_2)
\end{equation}
It is well known that $\nbigp_Y$ is $\nabla^G$-flat.
It is clearly equivariant with respect to the action of
$2\pi\sqrt{-1}H^2(X,\seisuu)$.

\subsubsection{Logarithmic extension and the specialization}
\label{subsection;14.11.28.100}

Put $\nbigubar_{\gamma 1}:=\nbigu_{\gamma 1}\cup\{0\}$
which is a neighbourhood of $0$ in $\cnum$.
We set $\nbigubar_{Y1}:=\nbigu_X\times\nbigubar_{\gamma 1}$.
We set 
$\QDMbar(Y):=
 \nbigo_{\cnum_{\lambda}\times\nbigubar_{Y1}}\otimes H^{\ast}(Y)$.
We naturally regard
$\QDMbar(Y)_{|\nbigu_{Y1}}=\QDM(Y)$.
The meromorphic flat connection $\nabla^G$ of $\QDM(Y)$
naturally gives a meromorphic flat connection
on $\QDMbar(Y)$
which is logarithmic along $t=0$.
The pairing $\nbigp_Y$ is naturally extended to
the following morphism
which is also denoted by $\nbigp_Y$:
\[
 \QDMbar(Y)
\otimes
 j^{\ast}\QDMbar(Y)
\lrarr
 \lambda^{n+r}\nbigo_{\cnum_{\lambda}\times\nbigubar_{Y1}}
\]

We set 
$\QDM(Y)^{\sp}:=
 \nbigo_{\cnum_{\lambda}\times\nbigu_X}
 \otimes H^{\ast}(Y)$.
We naturally have 
\[
\QDM(Y)^{\sp}=\QDMbar(Y)_{|\cnum_{\lambda}\times\nbigu_X\times{0}}.
\]
The residue
$\Res_{t}(\nabla^{G}):
 \lambda\QDM(Y)^{\sp}
\lrarr 
 \QDM(Y)^{\sp}$
is given by the multiplication of $-\lambda^{-1}\gamma$.
The specialization  of the connection 
$\nabla^{G\sp}$
is induced by the degenerated quantum products:
\[
 \bigl(
 \nabla^{G\sp}_{\xi}b
 \bigr)_{(\lambda,\sigma)}
=-\lambda^{-1}\xi\bullet_{\sigma}b,
\quad\quad
  \bigl(
 \nabla^{G\sp}_{\lambda\del_{\lambda}}b
 \bigr)_{(\lambda,\sigma)}
=\lambda^{-1}E\bullet_{\sigma}b
+\mu_Y(b),
\quad\quad
 \bigl(
(\lambda,\sigma)\in \cnum_{\lambda}^{\ast}\times\nbigu_X
 \bigr)
\]
It is equipped with the pairing $\nbigp_Y$ of weight $-n-r$
given by the formula (\ref{eq;14.11.20.2}).

\vspace{.1in}
By the construction
$\psi_{t,-\vecdelta}\bigl(\QDMbar(Y)(\ast t)\bigr)$
is naturally isomorphic to
$\QDM(Y)^{\sp}$ with $\nabla^{G\sp}$.
By the construction,
we have the following isomorphisms
which are compatible with the meromorphic connections:
\[
 \Cok\bigl(\lambda\Res_t(\nabla^G)\bigr)
 \simeq
 \QDM(X,\nbige^{\lor}),
\quad\quad
 \Ker\bigl(\lambda\Res_t(\nabla^G)\bigr)
 \simeq
 \nbigk(X,\nbige^{\lor})
\]

By the construction,
we have the following
which implies Lemma \ref{lem;14.11.28.200}
and Lemma \ref{lem;14.11.29.30}.
\begin{lem}
\label{lem;14.11.21.20}
The weight filtration $\Wtilde$ on 
$\QDM(X,\nbige^{\lor})$
and the pairings $\nbigp^{\nbige^{\lor}}_k$ 
on $\Gr^{W(N)}_k\QDM(X,\nbige^{\lor})$
are constructed 
from $\bigl(\QDM(Y),\nabla^G,\nbigp_Y\bigr)$
by the procedure in 
{\rm \S\ref{subsection;14.11.28.130}}.
\hfill\qed
\end{lem}

\subsection{Local mirror isomorphism}

Recall that 
the isomorphism
$L_X^{\lor}\simeq H^2(X,\seisuu)$
is induced by
$M_X^{\lor}\lrarr H^2(X,\seisuu)$
given by 
$\sum_{i=1}^m a_ie_i^{\lor}\longmapsto \sum a_i[D_i]$,
where $e_i^{\lor}$ $(i=1,\ldots,m)$
is the dual frame of the standard basis $e_1,\ldots,e_m$.

Let $\eta_1,\ldots,\eta_{\ell}$ be a frame of 
$H^2(X,\seisuu)$ such that
the first Chern class of
$K_X^{\lor}\otimes
 \bigotimes_{i=1}^r\nbigl_i^{\lor}$
is contained in the cone generated by
$\eta_i$ $(i=1,\ldots,\ell)$.
We can check the existence of such a frame
in an elementary way.
We have a frame $\xi_1,\ldots,\xi_{\ell}$ of $L_X^{\lor}$
corresponding to $\eta_1,\ldots,\eta_{\ell}$ 
by the above isomorphism.
For $M\in\real$,
we set 
\[
 \nbigb_M(\veceta):=
 \Bigl\{
 \sum_{i=1}^{\ell}
 \alpha_i\eta_i\in H^2(X,\cnum)
 \,\Big|\,
 \Re(\alpha_i)<-M\,\,(i=1,\ldots,\ell)
 \Bigr\},
\]
\[
 \nbigb'_M(\vecxi):=\Bigl\{
 \sum_{i=1}^{\ell}
 \beta_i\xi_i\in \gbigs_A=L_X^{\lor}\otimes\cnum
 \,\Big|\,
 \Re(\beta_i)<-M\,\,(i=1,\ldots,\ell)
 \Bigr\}.
\]

We obtain the following theorem
as a corollary of the theorem of Givental,
according to which
we have isomorphisms between
the mixed TEP-structures
associated to local A-models and local B-models
in \cite{Konishi-Minabe, Konishi-Minabe2}.
It is also related with 
\cite[Conjecture 6.14]{Reichelt-Sevenheck2}.
We use the notation
in \S\ref{subsection;14.12.29.1}
and \S\ref{subsection;14.12.8.10}.

\begin{thm}
\label{thm;14.11.21.10}
We have the following:
\begin{itemize}
\item
 An open subset $U_1\subset \nbigu_X\subset H^2(X,\cnum)$;
 It contains $\nbigb_{M_1}(\veceta)$ 
 for a positive number $M_1$;
 It is preserved by the natural action of
 $2\pi\sqrt{-1}H^2(X,\seisuu)$.
\item
 An open subset $U_2\subset 
 \gbigs_{\nbiga}^{\reg}$;
 It contains $\nbigb_{M_2}(\vecxi)$ 
 for a positive number $M_2$;
 It is preserved by the natural action of
 $(2\pi\sqrt{-1}\seisuu) L_X^{\lor}$.
\item
 A holomorphic isomorphism $\varphi:U_1\simeq U_2$;
 It preserves the Euler vector fields.
\item
Isomorphism of mixed TEP-structures
\[
 \Phi_{\ast}:
 \bigl(
 \QDM(X,\nbige^{\lor}),
 \Wtilde,
 \{\alpha \nbigp^{\nbige^{\lor}}_k\,|\,k\in\seisuu\}
 \bigr)_{|U_1}
\simeq
 \varphi^{\ast}
 \chi_{\nbiga}^{\ast}
 \bigl(
 \gbigv_{\nbiga,\ast},
 \Wtilde,
 \{\gbigp_{\nbiga,\ast,k}\,|\,k\in\seisuu\}
 \bigr)_{|U_2}
\]
\[
 \Phi_!:\bigl(
 \nbigk(X,\nbige^{\lor}),
 \Wtilde,
 \{\alpha \nbigp^{\nbige^{\lor}}_k\,|\,k\in\seisuu\}
 \bigr)_{|U_1}
\simeq
 \varphi^{\ast}
 \chi_{\nbiga}^{\ast}
 \bigl(
 \gbigv_{\nbiga,!},
 \Wtilde,
 \{\gbigp_{\nbiga,!,k}\,|\,k\in\seisuu\}
 \bigr)_{|U_2}
\]
Here, $\alpha$ is a non-zero complex number.
We also have the following commutative diagram
of mixed TE-structures:
\[
 \begin{CD}
\bigl(
\nbigk(X,\nbige^{\lor}),\Wtilde
 \bigr)_{|U_1}
 @>{\Phi_{!}}>{\simeq}>
 \varphi^{\ast}\chi_{\nbiga}^{\ast}
 (\gbigv_{\nbiga!},\Wtilde)_{|U_2}
 \\
 @VVV @VVV \\
 \bigl(
\QDM(X,\nbige^{\lor}),\Wtilde
 \bigr)_{|U_1}
 @>{\Phi_{\ast}}>{\simeq}>
 \varphi^{\ast}\chi_{\nbiga}^{\ast}
 (\gbigv_{\nbiga\ast},\Wtilde)_{|U_2}
 \end{CD}
\]
\end{itemize}
The isomorphisms
$\varphi$ and $\Phi_{\star}$ $(\star=\ast,!)$
are equivariant with respect to
 the actions of
 $2\pi\sqrt{-1}H^2(X,\seisuu)
 \simeq
 (2\pi\sqrt{-1}\seisuu)L_X^{\lor}$.
\end{thm}
\pf
Let $B_{\epsilon}(t):=\bigl\{|t|<\epsilon\bigr\}$.
According to \cite{Givental},
as reformulated
in \cite{Iritani} and \cite{Reichelt-Sevenheck1},
we have the following.
\begin{itemize}
\item
 An open subset $U_Y\subset \nbigu_Y$;
 It contains $\nbigb_{M_1}(\veceta)\times B_{\epsilon_1}(t_1)$
 for a large $M_1$ and a small $\epsilon_1>0$;
 It is preserved by the action of $2\pi\sqrt{-1}H^2(X,\seisuu)$.
\item
 An open subset 
 $U_{\nbigatilde}\subset
 \gbigstilde_{\nbiga}
=L_{\nbigatilde}\otimes\cnum$;
 It contains $\nbigb_{M_2}(\vecxi)\times B_{\epsilon_2}(t_2)$;
 It is preserved by the action of $2\pi\sqrt{-1}L_X^{\lor}$;
 The restriction of $\nbigv_{\nbiga}$ to $U_{\nbigatilde}$
 is locally free 
 $\nbigo_{\cnum_{\lambda}\times U_{\nbigatilde}}(\ast t_2)$-module.
\item
 A holomorphic isomorphism
 $\widehat{\varphi}:U_{Y}\simeq U_{\nbigatilde}$;
 It preserves the Euler vector fields;
 We have $\widehat{\varphi}^{-1}\{t_2=0\}=\{t_1=0\}$.
\item
An isomorphism of TEP-structures
$\widehat{\Phi}:
 \bigl(
 \QDMbar(Y)(\ast t),
 \alpha\nbigp_Y
 \bigr)_{|U_Y}
\simeq
 \widehat{\varphi}^{\ast}
 \chitilde_{\nbiga}^{\ast}
\bigl(
 \gbigvtilde_{\nbiga},\gbigptilde_{\nbiga}
\bigr)_{|U_{\nbigatilde}}$.
\item
The isomorphisms $\widehat{\varphi}$
and $\widehat{\Phi}$ are equivariant
with respect to the action of
$(2\pi\sqrt{-1})L_X^{\lor}
\simeq
 H^2(X,\seisuu)$.
\end{itemize}

We set
$U_1:=U_Y\cap\{t_1=0\}$
and 
$U_2:=U_{\nbigatilde}\cap\{t_2=0\}$.
As the restriction of $\widehat{\varphi}$,
we obtain a holomorphic isomorphism
$\varphi:U_1\simeq U_2$
which is equivariant with respect to 
the action of the lattices.
It also preserves the Euler vector fields.
Then, the claim follows from
Lemma \ref{lem;14.11.21.20}
and Lemma \ref{lem;14.11.22.1}.
\hfill\qed

\begin{rem}
If the Euler vector fields are $0$,
then the $\nbigrtilde$-modules
are $\cnum^{\ast}$-homogeneous 
with respect to the trivial action on the base space.
So, they are equipped with the Hodge filtration.
Because the Hodge filtration can be recovered 
from the connection,
the isomorphisms in Theorem {\rm\ref{thm;14.11.21.10}}
also preserve the Hodge filtrations.
\hfill\qed
\end{rem}

\begin{rem}
If $r=1$ and $\nbigl_1=K_X^{\lor}$,
the mixed TEP-structure is expressed 
as in {\rm\S\ref{subsection;14.12.3.120}}
in terms of the variation of Hodge structure
associated to $\nbigb$.
Here, $\nbigb$ is related to $\nbiga$
as in {\rm\S\ref{subsection;14.7.1.13}.}
\hfill\qed
\end{rem}

\begin{rem}
We set 
$U_1':=U_1\big/\bigl(2\pi\sqrt{-1}H^2(X,\seisuu)\bigr)$
and $U_2':=U_2\big/\bigl(2\pi\sqrt{-1}L_X^{\lor}\bigr)$.
We have the induced isomorphism
$\varphi':U_1'\simeq U_2'$.
We may regard 
$U_1'\subset H^2(X,\seisuu)\otimes\cnum^{\ast}$
and 
$U_2'\subset L_X^{\lor}\otimes\cnum^{\ast}$.
They are naturally extended to 
neighbourhoods
$\overline{U}_1'\subset H^2(X,\cnum)$ 
and 
$\overline{U}_2'\subset L_X^{\lor}\otimes\cnum$
of $0$.
According to
{\rm \cite{Reichelt-Sevenheck1}},
the isomorphism $\varphi'$
is extended to an isomorphism
$\overline{U}_1'\simeq
 \overline{U}_2'$.

We can consider the descent of the objects 
and the isomorphisms in Theorem {\rm\ref{thm;14.11.21.10}}
to $U_i'$.
The objects are extended on $\overline{U}_i'$
in the logarithmic way
such that the residues endomorphisms are nilpotent.
Then, 
by the uniqueness of such extensions,
the isomorphisms given on $U_i'$
are extended to the isomorphisms
on $\overline{U}_i'$.
\hfill\qed
\end{rem}

\begin{rem}
In {\rm\cite{Reichelt-Sevenheck2}},
it is announced that
a result in {\rm \cite{Iritani-Mann-Mignon}}
implies the comparison of the TE-structures
in Theorem {\rm\ref{thm;14.11.21.10}}.
At this moment, it is not clear to the author 
if we could also deduce 
the comparison of the weight filtrations
from {\rm \cite{Iritani-Mann-Mignon}}.
\hfill\qed
\end{rem}

\appendix
\section{Duality of meromorphic flat bundles}
\label{section;14.12.26.10}
\subsection{Isomorphism}
\label{subsection;14.11.14.11}

Let $X$ be a complex manifold
with a hypersurface $H$.
Let $d_X:=\dim X$.
Let $M$ be a reflexive meromorphic flat connection
on $(X,H)$,
i.e., $M$ is a $\nbigd_X$-module
and that $M$ is 
a coherent and reflexive $\nbigo_X(\ast H)$-module.
We set
$M^{\lor}:=
 \nhom_{\nbigo_X(\ast H)}
 \bigl(M,\nbigo_X(\ast H)\bigr)$
which is naturally a reflexive meromorphic flat connection.
In this section,
let $\DDD$ denote the duality functor
on the category of holonomic $\nbigd_X$-modules.
As well known,
$\DDD_{X(\ast H)}(M):=
 \DDD(M)(\ast H)$ is isomorphic to
$M^{\lor}$
as a meromorphic flat connection.
Let us recall the construction of 
the isomorphism
$\nu_{M}:
 \DDD_{X(\ast H)}(M)\stackrel{\simeq}{\lrarr}
 M^{\lor}$
as in \cite{Mochizuki-MTM}.

Let $\Theta_X$ denote the tangent sheaf of $X$.
We set $\Theta^{-p}_X:=\bigwedge^p\Theta_X$
for $p\geq 0$.
We have the Spencer resolution
$\nbigd_X\otimes\Theta^{\bullet}_X$
of $\nbigo_X$
by a locally free left $\nbigd$-modules.
We set $\Omega_X:=\Omega_X^{d_X}$.
The dual line bundle is denoted by 
$\Omega_X^{-1}$.
We have the following natural identification:
\[
 \DDD_{X(\ast H)}(M)
\simeq
 \nhom_{\nbigd_X(\ast H)}
 \bigl(
 \nbigd_X\otimes\Theta_X^{\bullet}\otimes M,
 \nbigd_X\otimes\Omega^{-1}_X(\ast H)
 \bigr)[d_X].
\]
The degree $0$ term in the right hand side is
$\nhom_{\nbigd_X(\ast H)}
 (\nbigd_X\otimes\Omega_X^{-1}\otimes M,
 \nbigd_X\otimes\Omega_X^{-1}(\ast H))
\simeq
 \nbigd_X\otimes M^{\lor}$.
The isomorphism $\nu_{M}$
is induced by the canonical isomorphism
$\DDD_{X(\ast H)}(M)
\simeq
 \nbigd_X\otimes\Theta^{\bullet}\otimes M^{\lor}$
whose degree $0$-term is the identity.

We have another description of the isomorphism $\nu_{M}$
in the case $H=\emptyset$.
Let $L_{M^{\lor}}$ be the sheaf of flat sections of $M^{\lor}$
which is a $\cnum$-local system.
We have a natural isomorphism
$\DR\DDD_{X} M=
 \nrhom_{\nbigd_X}(M,\nbigo_X[d_X])
\simeq
 L_{M^{\lor}}[d_X]$.
We have $\DR(\nu_{M})=(-1)^{d_X}\id$
which characterizes $\nu_{M}$.

Under the natural identification
$\DDD_{X(\ast H)}\circ
 \DDD_{X(\ast H)}(M)=M$,
we have the induced morphism
$\DDD_{X(\ast H)}\nu_{M}:
 \DDD_{X(\ast H)}(M^{\lor})\simeq M$.
We shall use the following lemma implicitly.
(See \cite{Mochizuki-MTM}.)
\begin{lem}
We have 
 $\DDD_{X(\ast H)}\nu_{M}
 =(-1)^{d_X}\nu_{M^{\lor}}$.
\hfill\qed
\end{lem}

\subsubsection{A commutative diagram}

For any morphism of meromorphic flat bundles
$f:M_1\lrarr M_2$ on $(X,H)$,
the following diagram is commutative:
\[
 \begin{CD}
 \DDD_{X(\ast H)} M_2@>{\DDD_{X(\ast H)} f}>>
 \DDD_{X(\ast H)} M_1\\
 @V{\nu_{M_2}}VV @V{\nu_{M_1}}VV \\
 M_2^{\lor}
 @>{f^{\lor}}>>
 M_1^{\lor}
 \end{CD}
\]

Let $M_i$ be reflexive meromorphic flat connections on $(X,H)$.
Let $\rho:M_1\lrarr \DDD_{X(\ast H)} M_2$
be a morphism of $\nbigd$-modules.
We have the induced morphisms
$\nu_{M_2}\circ\rho:
 M_1\lrarr M_2^{\lor}$
and 
$\nu_{M_1}\circ\DDD_{X(\ast H)}\rho:
 M_2\lrarr M_1^{\lor}$.
We shall use the following lemma 
in Proposition \ref{prop;14.11.8.40}.
\begin{lem}
\label{lem;14.11.8.41}
We have
$\nu_{M_1}\circ\DDD_{X(\ast H)}\rho
=(-1)^{d_X}(\nu_{M_2}\circ\rho)^{\lor}$.
\end{lem}
\pf
We have the following commutative diagram:
\[
 \begin{CD}
 \DDD_{X(\ast H)} (M_1)
 @<{\DDD_{X(\ast H)}(\nu_{M_2}\circ\rho)}<<
 \DDD_{X(\ast H)} (M_2^{\lor})
 \\
 @V{\nu_{M_1}}VV @VVV \\
 M_1^{\lor}
 @<{(\nu_{M_2}\circ\rho)^{\lor}}<<
 M_2
 \end{CD}
\]
Here, the right vertical arrow is
the composite of
$\DDD_{X(\ast H)} (M_2^{\lor})
\stackrel{\nu_{M_2^{\lor}}}\simeq
(M_2^{\lor})^{\lor}
\simeq M_2$.
Note that
$\nu_{M_1}\circ\DDD_{X(\ast H)} \rho$
precisely means the composite of
$M_2
 \simeq
 \DDD_{X(\ast H)} (\DDD_{X(\ast H)} M_2)
\stackrel{\DDD_{X(\ast H)} \rho}{\lrarr}
 \DDD_{X(\ast H)} M_1
\stackrel{\nu_{M_1}}{\simeq} M_1^{\lor}$.
We remark that the composite of
the following natural isomorphisms
is the multiplication of $(-1)^{d_X}$.
\[
\begin{CD}
 M_2
 @>>>
 (M_2^{\lor})^{\lor}
 @>>>
 \DDD_{X(\ast H)}(M_2^{\lor})
 @>{\DDD_{X(\ast H)}\nu_{M_2}}>>
 \DDD_{X(\ast H)} (\DDD_{X(\ast H)} M_2)
 @>>>M_2
\end{CD}
\]
Then, the claim of Lemma \ref{lem;14.11.8.41} follows.
\hfill\qed

\subsection{Duality and localizations of flat bundles}

This subsection is a preliminary for 
\S\ref{subsection;14.11.14.2}
and \S\ref{subsection;14.12.10.1}.

\subsubsection{Statement}
\label{subsection;14.11.14.1}

Let $X$ be a complex manifold
with a smooth hypersurface $Y$.
Set $d_X:=\dim X$.
Let $\iota:Y\lrarr X$ be the inclusion.
Let $M$ be a flat bundle on $X$,
i.e., a locally free $\nbigo_X$-module
with a flat connection.
The pull back of $M$ to $Y$ is denoted by $M_Y$.
For any $\nbigd_X$-module $N$,
we set $N(\ast Y):=N\otimes\nbigo_X(\ast Y)$
and $N(!Y):=\DDD_X\bigl(\DDD_X(N)(\ast Y)\bigr)$.

We have the natural morphisms
$\rho_{M,1}:M(\ast Y)
\lrarr
 \iota_{+}M_Y$
and 
$\rho_{M,2}:\iota_{+}M_Y
\lrarr M(!Y)$.
Locally, they are given as follows.
Let $(x_1,\ldots,x_n)$ be 
a holomorphic local coordinate system
of $X$ such that $Y=\{x_1=0\}$.
We have
$\iota_{+}M_Y
 \simeq
 \bigoplus_{n=0}^{\infty}
 \iota_{\ast}M_Y(dx_1)^{-1}\,\del_{x_1}^n$.
Then, we have
$\rho_{M,1}(f\cdot x_1^{-1})
=\iota_{\ast}(f_{|Y}(dx_1)^{-1})$
for any section $f$ of $M$.
And, we have
 $\rho_{M,2}(g\iota_{\ast}(dx_1)^{-1})
=-\del_{x_1}(\gtilde)$ in $M(!Y)$
for any section $g$ of $M_Y$,
where $\gtilde$ is a section of $M$
satisfying $\gtilde_{|Y}=g$ and $\del_{x_1}(\gtilde)=0$ in $M$.
They are independent of the choice of coordinate systems,
and we have the global morphisms.

The isomorphism
$\nu_M:\DDD M\simeq M^{\lor}$
induces
$\nu_{M\ast}:
\DDD\bigl(
 M(\ast Y)\bigr)
\simeq
 M^{\lor}(!Y)$
and 
$\nu_{M!}:
\DDD\bigl(
 M(!Y)\bigr)
\simeq
 M^{\lor}(\ast Y)$.
We have
$\nu_{M_Y}:\DDD M_Y\lrarr M^{\lor}_Y$
determined by a similar condition.
The composite 
$\DDD\iota_{+}M_Y
\simeq
 \iota_{+}\DDD M_Y
\simeq
 \iota_{+}M_Y^{\lor}$
is denoted by
$\iota_{+}\nu_{M_Y}$.
We shall prove the following proposition
in \S\ref{subsection;14.11.8.1}--\ref{subsection;14.11.8.2}.

\begin{prop}
\label{prop;14.9.12.10}
We have 
$\nu_{M\ast}\circ\DDD \rho_{M,1}
=-\rho_{M^{\lor},2}\circ\iota_{+}\nu_{M_Y}$
and 
$\iota_{+}\nu_{M_Y}\circ \DDD\rho_{M,2}
=\rho_{M^{\lor},1}\circ \nu_{M!}$.
\end{prop}

Clearly, it is enough to prove the case 
$M=\nbigo_X$
for which we naturally identify $M^{\lor}=\nbigo_X$.
We shall $\rho_{M,i}$ by $\rho_i$,
and $\nu_{M_X}$ and $\nu_{M_Y}$ by 
$\nu_X$ and $\nu_Y$, respectively.

We set $d_Y:=\dim Y=d_X-1$.
Because
$\DDD\nu_X=(-1)^{d_X}\nu_X$
and 
$\DDD\nu_Y=(-1)^{d_Y}\nu_Y$,
the two equalities are equivalent.
Hence, it is enough to prove
$\iota_{+}\nu_{Y}\circ \DDD\rho_2
=\rho_1\circ \nu_{X!}$.

We have only to prove the claim
locally around any point of $Y$.
We may assume that
$X=\cnum_t\times Y$.
We have natural isomorphisms
$\nbigo_X
\simeq
 \nbigo_{\cnum_t}\boxtimes\nbigo_Y$
and 
$\DDD\nbigo_X
\simeq \DDD\nbigo_{\cnum_t}
 \boxtimes
 \DDD\nbigo_Y$.
The equality
$\iota_{+}\nu_{Y}\circ \DDD\rho_2
=\rho_1\circ \nu_{X!}$
means the commutativity of the following diagram:
\[
 \begin{CD}
 \DDD\bigl(
 \nbigo_{\cnum_t}(!t)
 \bigr)
 \boxtimes
 \DDD\nbigo_Y
 @>{\nu_{\cnum_t!}\boxtimes\nu_Y}>>
 \nbigo_{\cnum_t}(\ast t)
 \boxtimes
 \nbigo_Y\\
 @VV{\DDD\rho_2\boxtimes\id}V
 @VV{\rho_1\boxtimes\id}V \\
 \DDD(\iota_{+}\nbigo_{\{0\}})
 \boxtimes
 \DDD\nbigo_Y
 @>{\iota_{+}\nu_{\{0\}}\boxtimes\nu_Y}>>
 \iota_{+}\nbigo_{\{0\}}
 \boxtimes
 \nbigo_Y
 \end{CD}
\]
Hence, it is enough to consider the case
$\dim X=1$,
i.e.,
$X=\cnum_t$
and $Y=\{0\}$,
which we assume in the following.

\subsubsection{Commutativity
of duality and push-forward}
\label{subsection;14.11.8.1}

First, let us describe the isomorphism
$\DDD\iota_{+}\nbigo_Y
\simeq
 \iota_{+}\DDD\nbigo_Y$.
We naturally have 
$\iota_{+}\nbigo_Y
=\bigoplus_{j=0}^{\infty}
\iota_{\ast} \cnum \cdot\del_t^{j}(dt)^{-1}$.
We have the exact sequence
$0\lrarr\nbigd_X\stackrel{\kappa_0}{\lrarr}
 \nbigd_X
 \stackrel{\kappa_1}{\lrarr}
 \iota_{+}\nbigo_Y
 \lrarr 0$,
where 
$\kappa_0(P)=P\cdot t$
and $\kappa_1(1)=\iota_{\ast}(dt)^{-1}$.
It induces the exact sequence:
\[
 0\lrarr
 \nhom_{\nbigd_X}(\nbigd_X,\nbigd_X\otimes\Omega_X^{-1})
\stackrel{\kappa_0'}\lrarr
 \nhom_{\nbigd_X}(\nbigd_X,\nbigd_X\otimes\Omega_X^{-1})
\stackrel{\kappa_1'}\lrarr
 \DDD(\iota_{+}\nbigo_Y)
\lrarr 0
\]
Here, $\kappa_0'(g)(P)=g(Pt)$.
Under the natural identifications
$\nhom_{\nbigd_X}(\nbigd_X,\nbigd_X\otimes\Omega_X^{-1})
\simeq
 \nbigd_X\otimes\Omega_X^{-1}$,
we have $\kappa_0'(Q)=tQ$.

We naturally have
$\iota_{+}\DDD\nbigo_Y
=\bigoplus_{j=0}^{\infty}
 \iota_{\ast}\Hom(\cnum,\cnum)\cdot
 \del_t^j\iota_{\ast}(dt)^{-1}$.
\begin{lem}
\label{lem;14.9.12.2}
The isomorphism
$\DDD(\iota_{+}\nbigo_Y)
\simeq
 \iota_{+}\DDD(\nbigo_Y)$
is given by
$\kappa_1'((dt)^{-1})
\longmapsto
 -\iota_{\ast}\bigl(\id\cdot(dt)^{-1}\bigr)$.
\end{lem}
\pf
Note that the following diagram is commutative,
due to M. Saito \cite{saito4}
(see also \cite[Proposition 12.4.25]{Mochizuki-MTM}):
\[
\begin{CD}
 \DR(\DDD\iota_{+}\nbigo_Y)
@>{\simeq}>>
 \DR(\iota_{+}\DDD\nbigo_Y) \\
 @VVV @VVV\\
 \DDD \iota_{\ast}\DR\nbigo_Y
@>{\simeq}>>
 \iota_{\ast}\DDD\DR\nbigo_Y
\end{CD}
\]
Here, the vertical arrows are 
the exchange of the de Rham functors
with the others,
and the horizontal arrows are the exchange of
the duality functor and the push-forward.
We describe the induced morphism
$\DR\DDD\iota_{+}\nbigo_Y
\simeq
 \iota_{\ast}\DDD\DR\nbigo_Y$.
We set
$\omega_X^{top}=\distribution_X^{\bullet}[2]$
and 
$\omega_Y^{top}=\cnum_Y$.
Then, it is described as
the composite of the following:
\begin{multline}
\nrhom_{\nbigd_X}\bigl(
 \iota_{+}\nbigo_Y,\nbigo_X[1]
 \bigr)
\stackrel{a}{\simeq}
 \nrhom_{\cnum_X}\bigl(
 \DR\iota_{+}\nbigo_Y,
 \omega_X^{top}
 \bigr)
\stackrel{b}{\simeq}
 \nrhom_{\cnum_X}
 \bigl(
 \iota_{\ast}\DR\nbigo_Y,
 \omega_X^{top}
 \bigr)
 \\
\stackrel{c}{\simeq}
 \nhom_{\cnum_X}
 (\iota_{\ast}\DR \nbigo_Y,\iota_{\ast}\omega_Y^{top})
=\iota_{\ast}\Hom(\cnum,\cnum).
\end{multline}
Let us describe the morphisms more explicitly.
We have
$\nbigo_X[1]=\distribution_X^{0,\bullet}[1]
=\bigl(
 \distribution_X^{0,0}
 \stackrel{-\delbar}{\lrarr}
 \distribution_X^{0,1}
 \bigr)$,
where $\distribution_X^{0,0}$ sits in the degree $-1$.
We have
\[
 \nrhom_{\nbigd_X}(\iota_{+}\nbigo_Y,\nbigo_X[1])
\simeq
 \nhom_{\nbigd_X}(\iota_{+}\nbigo_Y,
 \distribution_X^{0,\bullet}[1])
\simeq \distribution_{X,Y}^{0,\bullet}[1],
\]
where
$\distribution_{X,Y}^{0,i}$ denotes
the sheaf of $(0,i)$-currents on $X$
whose supports are contained in $Y$.
The last isomorphism is given by
$\Psi\longmapsto
 \Psi(\iota_{\ast}(dt)^{-1})$.
We have the section 
$-\delbar(t^{-1})$
of  $\distribution^{0,1}_{X,Y}$.
Recall that,
for a complex manifold $Z$
with $\dim Z=d_Z$,
we identify the total complex of
$\Omega_Z^{\bullet}[d_Z]
 \otimes\distribution_Z^{0,\bullet}[d_Z]$
with
$\distribution_Z^{\bullet}[2d_Z]$
by the correspondence
$\eta^{d_Z+i}\otimes
 \etabar^{d_Z+j}
\longmapsto
 (-1)^{d_Z(d_Z-1)/2+d_Zi}
 \eta^{d_Z+i}\wedge\etabar^{d_Z+j}$.
We have
\[
 a\bigl(
 -\delbar(1/t)
 \bigr)
 \bigl(
 \del_t^j\iota_{\ast}(dt)^{-1}dt
 \bigr)
=-\del_t^jdt\delbar(t^{-1})
=\del_t^j\delbar(dt/t),
\quad
  a\bigl(
 -\delbar(1/t)
 \bigr)
 \bigl(
 \del_t^j\iota_{\ast}(dt)^{-1}
 \bigr)
=\del_t^j\delbar(t^{-1}).
\]
Hence, we have
$(b\circ a)(-\delbar(t^{-1}))\bigl(
 \iota_{\ast}(dt)^{-1}dt
 \bigr)
=\delbar(dt/t)$.
Note that $c$ is induced by the trace morphism
$\iota_{\ast}\omega_Y^{top}\lrarr \omega_X^{top}$,
for which we have
$\iota_{\ast}(1)
\longmapsto \delbar(t^{-1}dt)$.
Hence, we have
$(c\circ b\circ a)(-\delbar(t^{-1}))=\iota_{\ast}1$.

\vspace{.1in}

We have the quasi-isomorphism:
\begin{equation}
\label{eq;14.9.12.1}
 \nhom(\iota_{+}\nbigo_Y,\distribution_X^{0,\bullet}[1])
\simeq
 \DR\DDD(\iota_{+}\nbigo_Y)
\simeq
 \Bigl(
  \nhom_{\nbigd_X}(\nbigd_X,\nbigo)
\stackrel{t}{\lrarr}
 \nhom_{\nbigd_X}(\nbigd_X,\nbigo)
 \Bigr)
\end{equation}
The complexes (\ref{eq;14.9.12.1})
are naturally quasi-isomorphic
to the total complex associated to the following double complex:
\[
 \begin{CD}
 \nhom_{\nbigd_X}(\nbigd_X,\distribution_X^{0,0})
 @>{t}>>
 \nhom_{\nbigd_X}(\nbigd_X,\distribution_X^{0,0})
 \\
 @V{-\delbar}VV @V{-\delbar}VV \\
 \nhom_{\nbigd_X}(\nbigd_X,\distribution_X^{0,1})
 @>{t}>>
  \nhom_{\nbigd_X}(\nbigd_X,\distribution_X^{0,1})
 \end{CD}
\quad\quad
 \left(
 \begin{CD}
 (0,-1) @>>> (1,-1)\\
 @VVV @VVV \\
 (0,0) @>>> (1,0)
 \end{CD}
\right)
\]
Here, the right square in the parenthesis indicates
the degrees of the double complex.
Then, we can easily deduce
$\delbar(t^{-1})$ in 
$\nhom(\iota_+\nbigo_Y,
 \distribution_X^{0,1})$
is equal to 
$1$ in
$\nhom_{\nbigd_X}(\nbigd_X,\nbigo)$
in the cohomology of (\ref{eq;14.9.12.1}).

Hence, the isomorphism
$\DR\iota_{+}\DDD\nbigo_Y
\simeq
 \DR\DDD(\iota_{+}\nbigo_Y)$
is identified with
$\iota_{\ast}\Hom(\cnum,\cnum)
\simeq
 \Cok(\nbigo\stackrel{t}{\lrarr}\nbigo)$
for which
$\iota_{\ast}\id$
is mapped to 
the image of $-1$
in 
$\Cok(\nbigo\stackrel{t}{\lrarr}\nbigo)$.
Then, the claim of Lemma \ref{lem;14.9.12.2} follows.
\hfill\qed

\vspace{.1in}

The isomorphism
$\nu_Y:\DDD\nbigo_Y\lrarr\nbigo_Y$
is given by
$\id\longmapsto 1$.
We have
the induced isomorphism
$\iota_{+}\nu_Y:
 \DDD(\iota_{+}\nbigo_Y)
\simeq
\iota_{+}\nbigo_Y$.
The following lemma is an immediate corollary
of Lemma \ref{lem;14.9.12.2}.
\begin{lem}
\label{lem;14.10.4.20}
We have
$\iota_{+}\nu_Y\bigl(
 \kappa'_1((dt)^{-1})\bigr)
=-\kappa_1(1)$.
\hfill\qed
\end{lem}

\subsubsection{A description of $\nu_X$}

Let us describe 
the isomorphism
$\nu_X:\DDD\nbigo_X\simeq\nbigo_X$.
We use the following natural free resolution of $\nbigo_X$:
\[
 \nbigo_X\simeq
 \bigl(
 \nbigd_X\stackrel{\cdot\del_t}{\lrarr}
 \nbigd_X
 \bigr)
\]
It induces the following free resolution
of $\DDD\nbigo_X$:
\[
 \DDD\nbigo_X\simeq
 \bigl(
 \nbigd_X\otimes\Omega_X^{-1}
\stackrel{\del_t\cdot}{\lrarr}
 \nbigd_X\otimes\Omega_X^{-1}
 \bigr) 
\]
We have the morphism of
$\nbigd_X$-modules
$\mu:
 \nbigd_X
\lrarr
 \nbigd_X\otimes\Omega_X^{-1}$
given by
$a_j\del_t^j
\longmapsto
 (-1)^j\del_t^ja_j\otimes(dt)^{-1}$.
We have the following commutative diagram:
\begin{equation}
\label{eq;14.10.4.11}
 \begin{CD}
 \nbigd_X\otimes\Omega_X^{-1}
 @>{\del_t\cdot}>>
 \nbigd_X\otimes\Omega_X^{-1} \\
 @V{-\mu^{-1}}VV @V{\mu^{-1}}VV \\
 \nbigd_X
 @>{\cdot\del_t}>>
 \nbigd_X
 \end{CD}
\end{equation}
It induces
$\DDD\nbigo_X
\lrarr
 \nbigo_X$.
We can check that
the induced morphism
$\DR\nbigo_X\simeq
 \DR\DDD\nbigo_X$
is given by
$\cnum\simeq \Hom(\cnum,\cnum)$,
$1\longmapsto -\id$.
Hence, it is $\nu_X$.

\subsubsection{End of the proof of 
 Proposition \ref{prop;14.9.12.10}}
\label{subsection;14.11.8.2}

We have the following commutative diagram:
\begin{equation}
\label{eq;14.9.12.3}
 \begin{CD}
 \nbigd_X
 @>{\cdot\del_t}>>
 \nbigd_X
 @>>>
 \nbigo_X
 \\ 
 @V{\id}VV @V{\cdot t}VV @VVV 
 \\
  \nbigd_X
 @>{\cdot \del_tt}>>
 \nbigd_X
 @>{\varphi_1}>>
 \nbigo_X(\ast t)
 \\
 @V{\cdot \del_t}VV @V{\id}VV @V{\rho_1}VV
 \\
 \nbigd_X
 @>{\cdot t}>>
 \nbigd_X
 @>{\kappa_1}>>
 \iota_{+}\nbigo_Y
 \end{CD}
\end{equation}
Here $\varphi_1(1)=t^{-1}$.
We also have the following commutative diagram:
\begin{equation}
\label{eq;14.9.12.5}
 \begin{CD}
 \nbigd_X
 @>{\cdot t}>>
 \nbigd_X
 @>{\kappa_1}>>
 \iota_{+}\nbigo_Y
 \\
 @V{-\id}VV @V{\cdot(-\del_t)}VV @V{\rho_2}VV\\ 
 \nbigd_X
 @>{\cdot t\del_t}>>
 \nbigd_X
 @>{\varphi_2}>>
 \nbigo_X(!t)
 \\
 @V{\cdot t}VV @V{\id}VV @VVV \\
 \nbigd_X
 @>{\cdot \del_t}>>
 \nbigd_X
 @>>>
 \nbigo_X
 \end{CD}
\end{equation}
Here $\varphi_2(1)=1$.
As the dual of (\ref{eq;14.9.12.5}),
we obtain the following:
\begin{equation}
 \label{eq;14.10.4.10}
 \begin{CD}
 \DDD(\iota_{+}\nbigo_Y)
@<{\kappa_1'}<<
 \nbigd_X\otimes\Omega_X^{-1}
 @<{t\cdot }<<
 \nbigd_X\otimes\Omega_X^{-1}
 \\
@A{\DDD\rho_2}AA @A{-\id}AA @A{(-\del_t)\cdot}AA 
 \\ 
 \DDD\bigl(\nbigo_X(!t)\bigr)
@<<< 
 \nbigd_X\otimes\Omega_X^{-1}
 @<{t\del_t\cdot }<<
 \nbigd_X\otimes\Omega_X^{-1}
 \\
@AAA @A{t}AA @A{\id}AA 
\\
 \DDD\nbigo_X
 @<<<
 \nbigd_X\otimes\Omega_X^{-1}
 @<{\del_t\cdot}<<
 \nbigd_X\otimes\Omega_X^{-1}
 \end{CD}
\end{equation}
We have the following commutative diagram:
\[
 \begin{CD}
 \nbigd_X\otimes\Omega_X^{-1}
 @>{(t\del_t)\cdot}>>
 \nbigd_X
 @>>>
 \DDD(\nbigo_X(!t))
  \\
 @V{-\mu^{-1}}VV @V{\mu^{-1}}VV
 @V{\nu_{X!}}VV \\
 \nbigd_X
 @>{\cdot(\del_tt)}>>
 \nbigd_X
 @>{\varphi_1}>>
 \nbigo_X(\ast t)
 \end{CD}
\]
Because $\varphi_1(1)=t^{-1}$,
we have
$\DDD\rho_2\circ
 \nu_{X!}^{-1}(t^{-1})
=-\kappa'_1\bigl((dt)^{-1}\bigr)$.
We also have $\rho_1(t^{-1})=\kappa_1(1)$.
Then, we obtain 
$\DDD\rho_2\circ\nu_{X!}^{-1}
=(\iota_{+}\nu_Y)^{-1}\circ\rho_1$
by Lemma \ref{lem;14.10.4.20}.
Thus, the proof of Proposition \ref{prop;14.9.12.10}
is finished.
\hfill\qed

\subsection{Duality and nearby cycle functors}
\label{subsection;14.11.14.2}

We shall generalize Proposition \ref{prop;14.9.12.10}.
This subsection is a preliminary for
\S\ref{subsection;14.11.28.30}.

\subsubsection{Some induced isomorphisms}
\label{subsection;14.11.15.20}

Let $t$ be a function on $X$
such that the zero divisor $(t)_0$ is smooth and reduced,
i.e., the zero set $H$ of $t$ is smooth
and reduced as a subscheme.
Let $M$ be a regular singular meromorphic
flat bundle on $(X,H)$.

We use the natural identification
$\bigl(
 \Pi^{a,b}_tM
 \bigr)^{\lor}
\simeq
 \Pi_t^{-b+1,-a+1}M^{\lor}$.
(See \cite[\S2.3]{Mochizuki-MTM}
where the identification is explained for $\nbigr$-modules.
The identification for flat bundles is obtained
as the specialization at $\lambda=1$.)
We have the isomorphisms
$\nu_{\Pi^{a,b}_tM}:(\DDD_{X(\ast H)} \Pi_t^{a,b}M)
 \simeq \Pi^{-b+1,-a+1}_tM^{\lor}$,
which induce
\[
 \nu_{\Pi^{a,b}_tM,\ast}:
 \DDD_X \bigl(
\Pi_{t\ast}^{a,b}M
 \bigr)
 \simeq
 \Pi^{-b+1,-a+1}_{t!}M^{\lor},
\quad\quad
 \nu_{\Pi^{a,b}_tM,!}:
 \DDD_X\bigl(\Pi_{t!}^{a,b}M\bigr)
 \simeq
 \Pi^{-b+1,-a+1}_{t\ast}M^{\lor}.
\]
We have the following isomorphism
for a sufficiently large $N$:
\[
 \psi^{(a)}_t(M)\simeq
 \Ker\Bigl(
  \Pi^{a-N,a}_{t!}(M)
\lrarr
 \Pi^{a-N,a}_{t\ast}(M)
 \Bigr)
\simeq
 \Cok\Bigl(
 \Pi^{a,a+N}_{t!}(M)
\lrarr
 \Pi^{a,a+N}_{t\ast}(M)
 \Bigr),
\]
Hence, we have the following canonical isomorphisms:
\[
 \nu_{M,\psi_t^{(a)}}:
 \DDD_H\psi^{(a)}_t(M)
\simeq
 \psi^{(-a+1)}_t(M^{\lor})
\]

We have the $V$-filtration $U$ of $M$ along $t$.
Here, we adopt the condition that
$\del_tt+\alpha$ is nilpotent on
$\psi_{t,\alpha}(M):=U_{\alpha}(M)\big/U_{<\alpha}(M)$.
The natural perfect pairing
$M\otimes M^{\lor}\lrarr \nbigo_X(\ast H)$
naturally induces a perfect pairing of
$\nbigo_X$-modules
$V_{-1}(M)
\otimes
 V_{-1}(M^{\lor})
\lrarr 
\nbigo_X$.
It induces a perfect pairing
$\psi_{t,-1}(M)\times
 \psi_{t,-1}(M^{\lor})
\lrarr
 \nbigo_H$
of flat bundles.
It gives the following natural identification:
\begin{equation}
 \psi_{t,-1}(M)^{\lor}
\simeq
 \psi_{t,-1}(M^{\lor}).
\end{equation}

We have the canonical isomorphisms
$\xi^{(0)}_M:\psi_t^{(0)}(M)
\simeq
 \iota_{+}\psi_{t,-1}(M)$
and 
$\xi^{(1)}_M:
\psi_t^{(1)}(M)\simeq
\iota_{+}\psi_{t,-1}(M)$.
The construction of the isomorphisms
is explained in the case of
$\nbigr$-modules in \cite[S4.3]{Mochizuki-MTM}.
The construction for $\nbigd$-modules
is obtained as the specialization to $\lambda=1$.

\subsubsection{Some commutative diagrams}

\begin{prop}
\label{prop;14.10.4.40}
The following diagram is commutative:
\begin{equation}
 \label{eq;14.9.29.12}
 \begin{CD}
 \DDD_X\psi_t^{(1)}(M)
 @>{\nu_{M,\psi_t^{(1)}}}>>
 \psi_t^{(0)}(M^{\lor})
 \\
 @A{\simeq}A{\DDD_X \xi_M^{(1)}}A
 @V{\simeq}V{\xi_{M^{\lor}}^{(0)}}V \\
 \DDD_X\iota_{+} \psi_{t,-1}(M)
 @>{\simeq}>>
 \iota_{+}
 \psi_{t,-1}(M^{\lor})
 \end{CD}
\end{equation}
Here, the lower horizontal arrow is
induced as the composite of 
the natural isomorphisms
$\DDD_X\iota_{+} \psi_{t,-1}(M)
\simeq
 \iota_{+}
 \DDD_H\psi_{t,-1}(M)
\simeq
 \iota_{+}\psi_{t,-1}(M)^{\lor}
\simeq
 \iota_{+}
 \psi_{t,-1}(M^{\lor})$.

The following diagram is also commutative.
\begin{equation}
\label{eq;14.9.29.13}
 \begin{CD}
 \DDD_X\psi_t^{(0)}(M)
 @>{-\nu_{M,\psi^{(0)}}}>>
 \psi_t^{(1)}(M^{\lor})
 \\
 @A{\simeq}A{\DDD_X\xi_M^{(0)}}A 
 @V{\simeq}V{\xi_{M^{\lor}}^{(1)}}V \\
 \DDD_X\iota_{+}
 \psi_{t,-1}(M)
@>{\simeq}>>
  \iota_{+}\psi_{t,-1}(M^{\lor})
 \end{CD}
\end{equation}
Here, the lower horizontal arrow is induced 
as the composite of the natural isomorphisms
$ \DDD_X\iota_{+}
 \psi_{t,-1}(M)
\simeq
 \iota_{+}\DDD_H\psi_{t,-1}(M)
\simeq
 \iota_{+}\psi_{t,-1}(M)^{\lor}
\simeq
  \iota_{+}\psi_{t,-1}(M^{\lor})$.

\end{prop}
\pf
Let $V$ be any regular singular meromorphic flat bundle
on $(X,H)$.
We set $V_{\star}:=V(\star H)$.
Let $K_V$ and $C_V$
denote the kernel and the cokernel of
$V_!\lrarr V_{\ast}$.
We have the following natural commutative diagram:
\[
 \begin{CD}
 0@>>> \DDD_X C_V
 @>>> \DDD_X V_{\ast}
 @>>> \DDD_XV_!
 @>>> \DDD_X K_V
 @>>> 0\\
 @. 
 @V{\nu_{V,C}}VV 
 @V{\nu_{V\ast}}VV 
 @V{\nu_{V!}}VV 
 @V{\nu_{V,K}}VV 
 @.\\
 0@>>> K_{V^{\lor}}
 @>>>
 (V^{\lor})_!
 @>>>
 (V^{\lor})_{\ast}
 @>>>
 C_{V^{\lor}}
 @>>> 0
 \end{CD}
\]

Let $K_{0V}$ and $C_{0V}$
denote the kernel and the cokernel of
the morphism
$\psi_{t,-1}(V)
 \stackrel{N}{\lrarr}
 \psi_{t,-1}(V)$.
We have natural isomorphisms
$\xi_{C,V}:C_V\simeq
 \iota_{+}C_{0V}$
and 
$\xi_{K,V}:K_V\simeq
 \iota_{+}K_{0V}$.
The construction of the isomorphisms
is explained in \cite[\S4.3.2]{Mochizuki-MTM}
for $\nbigr$-modules.
The construction for $\nbigd$-modules
is obtained as the specialization at $\lambda=1$.
Under the natural identification
$\psi_{t,-1}(V)^{\lor}
=\psi_{t,-1}(V^{\lor})$,
we have
$C_{0V}^{\lor}
=K_{0V^{\lor}}$
and 
$K_{0V}^{\lor}
=C_{0V^{\lor}}$.

\begin{lem}
\label{lem;14.10.4.50}

The following diagram is commutative:
\begin{equation}
\label{eq;14.9.29.10}
 \begin{CD}
 \DDD_X K_V
 @>{\nu_{V,K}}>>
 C_{V^{\lor}} \\
 @A{\simeq}A{\DDD_X\xi_{K,V}}A
 @V{\simeq}V{\xi_{C,V^{\lor}}}V \\
 \DDD_X(\iota_{+}K_{0V})
 @>>>
 \iota_{+}C_{0V^{\lor}}
 \end{CD}
\end{equation}
The lower horizontal arrow is
the composite of
$\DDD_X(\iota_{+}K_{0V})
\simeq
 \iota_{+}\DDD_H K_{0V}
\simeq
 \iota_{+}K_{0V}^{\lor}
\simeq
 \iota_{+}C_{0V^{\lor}}$.

The following diagram is commutative:
\begin{equation}
 \label{eq;14.9.29.11}
 \begin{CD}
 \DDD_X C_V
 @>{-\nu_{V,C}}>>
 K_{V^{\lor}} \\
 @A{\simeq}A{\DDD_X\xi_{C,V}}A 
 @V{\simeq}V{\xi_{K,V^{\lor}}}V \\
 \DDD_X(\iota_{+}C_{0V})
 @>>>
 \iota_{+}K_{0V^{\lor}}
 \end{CD}
\end{equation}
The lower horizontal arrow is
the composite of
$\DDD_X(\iota_{+}C_{0V})
\simeq
 \iota_{+}\DDD_H C_{0V}
\simeq
 \iota_{+}C_{0V}^{\lor}
\simeq
 \iota_{+}K_{0V^{\lor}}$.

\end{lem}
\pf
It is enough to consider the case that
the monodromy of $V$ along the loop around
$t=0$ is unipotent.
Moreover, we may assume that
the logarithm of the monodromy is
a Jordan block.
Then, we have a flat subbundle
$L\subset V$ of rank one
such that  $K_L\simeq K_V$
and $K_{0L}\simeq K_{0V}$.
We also have
$V^{\lor}\lrarr L^{\lor}$
which induces
$C_{0V^{\lor}}\simeq C_{0L^{\lor}}$
and $C_{V^{\lor}}\simeq C_{L^{\lor}}$.
Then, we obtain (\ref{eq;14.9.29.10})
from the claim for $L$,
which follows from Proposition \ref{prop;14.9.12.10}.
We obtain (\ref{eq;14.9.29.11})
similarly.
\hfill\qed

\vspace{.1in}

We have 
$K_{0,\Pi^{-N,1}M}\simeq
  \psi_{t,-1}(M)$
and 
$C_{0,\Pi^{0,N+1},M^{\lor}}
\simeq
 \psi_{t,-1}(M^{\lor})$.
Hence, we obtain the commutativity of
(\ref{eq;14.9.29.12}) from 
Lemma \ref{lem;14.10.4.50}.
We obtain the commutativity of
(\ref{eq;14.9.29.13}) similarly.
\hfill\qed

\subsubsection{Specialization of pairings}

Let $M_i$ be regular singular meromorphic flat bundles
on $(X,H)$.
Let $P:M_1\otimes M_2\lrarr \nbigo_X(\ast H)$
be a morphism of $\nbigd$-modules.
We have the naturally induced pairing
$\psi_{t,-1}(P):
 \psi_{t,-1}(M_1)\otimes\psi_{t,-1}(M_2)\lrarr\nbigo_H$.

The pairing $P$ and the isomorphism $\nu_{M_2}$ induce
a morphism $\Psi_P:M_1\simeq \DDD_{X(\ast H)} M_2$.
For $a=0,1$, we have the induced morphisms
\[
 \psi^{(a)}(M_1)
\lrarr
 \psi^{(a)}(\DDD_X M_2)
\simeq
 \DDD_X\psi^{(-a+1)}(M_2).
\]
By the natural isomorphisms
$\iota_{+}\psi_{t,-1}(M_i)
 \simeq
 \psi^{(a)}(M_i)$
and the isomorphism $\nu_{\psi_{t,-1}(M_2)}$,
we obtain
\[
 \iota_{+}\psi_{t,-1}(M_1)
\simeq
 \DDD_X \iota_{+}\psi_{t,-1}(M_2)
\simeq
 \iota_{+}\DDD_H\psi_{t,-1}(M_2)
\simeq
 \iota_{+}\psi_{t,-1}(M_2)^{\lor}.
\]
It comes from a morphism
$\rho_a:\psi_{t,-1}(M_1)\simeq\psi_{t,-1}(M_2)^{\lor}$,
which is equivalent to a pairing
\[
P_a:\psi_{t,-1}(M_1)\times\psi_{t,-1}(M_2)
\lrarr\nbigo_H.
\]

\begin{cor}
\label{cor;14.11.15.32}
We have $P_0=\psi_{t,-1}(P)$
and $P_1=-\psi_{t,-1}(P)$.
\end{cor}
\pf
Let $\mu:\psi_{t,-1}(M_1)\lrarr \psi_{t,-1}(M_2)^{\lor}$
be the morphism induced by $\psi_{t,-1}(P)$.
The push-forward
$\iota_{+}\mu:
 \iota_{+}\psi_{t,-1}(M_1)
\lrarr
 \iota_{+}\psi_{t,-1}(M_2)^{\lor}$
is the composite of the following morphisms:
\[
 \iota_{+}\psi_{t,-1}(M_1)
\simeq
 \psi^{(a)}(M_1)
\lrarr
 \psi^{(a)}(M_2^{\lor})
\simeq
 \iota_{+}\psi_{t,-1}(M_2^{\lor})
\]
By the construction,
$\iota_{+}\rho_a$ is composite of the following morphisms:
\begin{multline}
 \iota_{+}\psi_{t,-1}(M_1)
\simeq
 \psi^{(a)}(M_1)
\lrarr
  \psi^{(a)}(M_2^{\lor})
 \simeq
 \psi^{(a)}(\DDD_X M_2)
\simeq
 \DDD_X\psi^{(-a+1)}(M_2)
\simeq
 \DDD_X\iota_{+}\psi_{t,-1}(M_2)
 \\
\simeq
 \iota_{+}\DDD_X\psi_{t,-1}(M_2)
\simeq
 \iota_{+}\psi_{t,-1}(M_2)^{\lor}
\end{multline}
By Proposition \ref{prop;14.10.4.40},
we have 
$\iota_{+}\mu
=\iota_{+}\rho_0$
and 
$\iota_{+}\mu
=-\iota_{+}\rho_1$.
\hfill\qed

\subsection{Push-forward}

This subsection is a preliminary for
\S\ref{subsection;14.12.10.2}.

\subsubsection{Statement}

Let $X$ be a complex manifold.
Let $M_i$ be flat bundles on $X$.
Let $P:M_1\otimes M_2\lrarr\nbigo_X$
be a morphism of $\nbigd$-modules.
By the isomorphism
$\nu_{M_2}:\DDD_X M_2\simeq M_2^{\lor}$,
$P$ is equivalent to a morphism
$\varphi:M_1\lrarr \DDD_X M_2$.

Let $F:X\lrarr S$ be a smooth proper morphism
of complex manifolds.
We have the flat bundles $F^j_{+}M_i$ on $S$.
We obtain the following morphism:
\[
\begin{CD}
 F_{+}^jM_1
 @>>>
 F_{+}^j\DDD_X M_2
 @>>>
 \DDD_S F_{+}^{-j}M_2
\end{CD}
\]
The composite is denoted by
$F^j_{+}\varphi$.
By the isomorphism
$\nu_{F_{+}^{-j}M_2}:
  \DDD_S F_{+}^{-j}M_2
\simeq 
\bigl(
 F_{+}^{-j}M_2
\bigr)^{\lor}$,
we obtain a flat morphism
\[
 P(F_{+}^{j}\varphi):
 F_{+}^jM_1
\otimes
 F_{+}^{-j}M_2
\lrarr
 \nbigo_S.
\]

We give a more direct expression of the pairing
$P(F_{+}^j\varphi)$.
For each $Q\in S$,
we set $X_Q:=F^{-1}(Q)$,
and let $M_{i,Q}$ denote the restriction 
of the flat bundles $M_i$ to $X_Q$.
Let $d:=\dim X-\dim S$.
Let $H_{\DR}^{\bullet}(X_Q,M_{i,Q})$
denote the de Rham cohomology of
the flat bundle $M_{i,Q}$.
The fiber of $F_{+}^j(M_i)_{|Q}$ is 
naturally identified with 
$H^{d+j}_{\DR}\bigl(X_Q,M_{i,Q}\bigr)$.
The pairing $P$ naturally induces
the following pairing 
$P_{F,j,Q}$:
\[
 H^{d+j}_{\DR}\bigl(X_Q,M_{1,Q}\bigr)
\times
 H^{d-j}_{\DR}\bigl(X_Q,M_{2,Q}\bigr)
\lrarr
 H_{\DR}^{2d}(X_Q,\cnum)
\lrarr \cnum
\]
Here, $H_{\DR}^{2d}(X_Q,\cnum)\lrarr\cnum$
is induced by $(2\pi\sqrt{-1})^{-d}\int_{X_Q}$.
The family 
$\{P_{F,j,Q}\,|\,Q\in S\}$
gives a flat morphism
$P_{F,j}:
 F_{+}^jM_1\otimes
 F_{+}^{-j}M_2\lrarr\nbigo_S$.
We shall prove the following proposition
in \S\ref{subsection;14.11.15.12}--\ref{subsection;14.11.15.13}.

\begin{prop}
\label{prop;14.11.15.1}
We have 
$P(F_{+}^j\varphi)
=\epsilon(d)(-1)^{dj}P_{F,j}$,
where $\epsilon(d)=(-1)^{d(d-1)/2}$.
\end{prop}

\subsubsection{Preliminary}
\label{subsection;14.11.15.12}

We have only to prove the claim
locally around any point of $S$.
So, we may assume that $S$ is a multi-disc.
Let $t$ be a holomorphic function on $S$
such that the zero divisor $(t)_0$ is smooth and reduced.
Let $S_0:=|(t)_0|$.
The pull back $t\circ F$ is also denoted by $t$.
Set $X_0:=\{Q\in X\,|\,t(Q)=0\}$.
The induced morphism
$X_0\lrarr S_0$ is denoted by $F_0$.
We naturally identify
the restriction of the flat bundles $M_i$
to $X_0$
with $\psi_{t,-1}(M_i)$.
We also identify 
the restriction of $F_{+}^jM_i$
to $S_0$
with
$\psi_{t,-1}F_{+}^j(M_i)$,
which is also naturally isomorphic to
$F_{0+}^j\psi_{t,-1}(M_i)$.

We obtain the following pairings 
as the restriction of $P$,
$P(F^j_{+}\varphi)$
and $P_{F,j}$:
\[
 \psi_{t,-1}P:\psi_{t,-1}(M_1)\otimes \psi_{t,,-1}(M_2)
\lrarr \nbigo_{X_0}
\]
\[
 \psi_{t,-1}P(F^j_{+}\varphi):
 \psi_{t,-1}F_{+}^j(M_1)
\otimes
 \psi_{t,-1}F_{+}^{-j}(M_2)
\lrarr
 \nbigo_{S_0}
\]
\[
 \psi_{t,-1}P_{F,j}:
 \psi_{t,-1}F_{+}^j(M_1)
\otimes
 \psi_{t,-1}F_{+}^{-j}(M_2)
\lrarr
 \nbigo_{S_0}
\]
By the construction,
we have
$(\psi_{t,-1}P)_{F_0,j}
=\psi_{t,-1}P_{F,j}$.

Let $\varphi_0:\psi_{t,-1}M_1\lrarr
 \DDD_{X_0}\psi_{t,-1}M_2$
be the morphism corresponding to $\psi_{t,-1}P$.
\begin{lem}
\label{lem;14.11.15.20}
We have
$\psi_{t,-1}P(F^j_{+}\varphi)
=P(F^j_{+}\varphi_0)$.
As a consequence,
to prove Proposition {\rm\ref{prop;14.11.15.1}},
it is enough to consider the case
where $S$ is a point.
\end{lem}
\pf
Let $\iota_{X_0}:X_0\lrarr X$ be the inclusion.
By Proposition \ref{prop;14.10.4.40},
the morphism
$\iota_{X_0+}\varphi_0:
 \iota_{X_0+}\psi_{t,-1}M_1
\lrarr
 \iota_{X_0+}\DDD_{X_0}\psi_{t,-1}M_2$
is identified with the composite of the following morphism:
\begin{equation}
\label{eq;14.11.15.10}
\begin{CD}
 \psi^{(0)}_{t}(M_1)
 @>{\psi^{(0)}_t\varphi}>>
 \psi^{(0)}_t(\DDD_X M_2)
 \simeq
 \DDD_X \psi^{(1)}_t(M_2)
\end{CD}
\end{equation}
We have the morphism
$\rho_j:
 \psi_{t,-1}F_{+}^jM_1
\lrarr
 \DDD_{S_0}\psi_{t,-1}F_{+}^jM_2$
corresponding to
$\psi_{t,-1}P(F_{+}^j\varphi)$.
Let $\iota_{S_0}:S_0\lrarr S$ be the inclusion.
By Proposition \ref{prop;14.10.4.40},
$\iota_{S_0+}\rho_j$
is identified with 
\begin{equation}
 \label{eq;14.11.15.11}
\begin{CD}
 \psi_{t}^{(0)}F_{+}^jM_1
@>{\psi_t^{(0)}F_{+}^j\varphi}>>
 \psi_{t}^{(0)}\DDD_S F_{+}^{-j}M_2
\simeq
 \DDD_S\psi_t^{(1)}F_{+}^{-j}M_2.
\end{CD}
\end{equation}
Note that the following diagram 
of the natural isomorphisms is commutative,
which can be checked easily:
\[
 \begin{CD}
 \DDD_S \psi_t^{(1)} F_{+}^{-j}M_2
 @>{\simeq}>>
 \psi_t^{(0)}\DDD_S F_{+}^{-j}M_2
 @>{\simeq}>>
  \psi_t^{(0)}F_{+}^{-j}\DDD_X M_2
 \\
 @V{\simeq}VV @. @V{\simeq}VV \\
  \DDD_S F_{+}^{-j}\psi_{t}^{(1)}M_2
 @>{\simeq}>>
F_{+}^{-j}  \DDD_X \psi_t^{(1)}M_2
 @>{\simeq}>>
  F_{+}^{-j}\psi_t^{(0)}
 \DDD_X M_2
 \end{CD}
\]
Then, the morphism (\ref{eq;14.11.15.11})
is identified with the push-forward of
(\ref{eq;14.11.15.10}) by $F$.
Hence, we have 
$\rho_j=F_{+}^j\varphi_0$,
which implies the first claim of the lemma.
The second claim follows from
the first claim
and the equality
$(\psi_{t,-1}P)_{F_0,j}
=\psi_{t,-1}P_{F,j}$.
\hfill\qed

\subsubsection{The case where $S$ is a point}
\label{subsection;14.11.15.13}

Let $d:=\dim X$.
Let $L_i$ be the sheaf of flat sections of $M_i$.
Let $L_i^{\lor}$ denote the local system
obtained as the dual of $L_i$.
The pairing $P$ induces
$\langle\cdot,\cdot\rangle:
 L_1\times L_2\lrarr \cnum_X$,
which induces 
$\upsilon:L_1\lrarr L_2^{\lor}$.
We have the natural isomorphism
$\DR M_1\simeq L_1[d]$
and $\DR \DDD M_2\simeq 
 \nrhom_{\nbigd_X}(M_2,\nbigo_X[d])
\simeq L_2^{\lor}[d]$.
The morphism $\varphi$
induces 
$(-1)^d\upsilon:
 L_1[d]\lrarr L^{\lor}_2[d]$.

Let us describe the induced morphism
$L_1[d]\lrarr
 \nrhom_{\cnum_X}\bigl(
 \DR M_2,
 \DR \nbigo_X[d]
 \bigr)$.
A section $e$ of $L_1$ naturally gives
a section of $(L_1[d])^{-d}$,
which is denoted by $e'$.
We have the section $(-1)^d\upsilon(e')$
of $\nhom_{\nbigd_X}(M_2,\nbigo[d])^{-d}$.
The image of $(-1)^d\upsilon(e')$ in
$\nhom_{\cnum_X}(\DR M_2,\DR\nbigo[d])^{-d}$
is denoted by $F_e$. 
Then, we have
\[
 F_e(\eta^{d+j}m)
=(-1)^{dj}\eta^{d+j}\cdot
 (-1)^d\upsilon(e')(m)
=(-1)^{dj+d}\eta^{d+j}\langle e,m\rangle'
\]
Here, $\eta^{d+j}$ is a section of
$(\Omega^{\bullet}[d])^j=\Omega^{d+j}$,
$m$ is a section of $M_2$,
and $\langle e,m\rangle'$ denotes
the section of $(\nbigo[d])^{-d}$
corresponding to $\langle e,m\rangle$.

We have the natural quasi isomorphisms
$\DR\nbigo[d]
\simeq
 \Tot\bigl(
 \Omega^{\bullet}[d]\otimes
 \distribution_X^{0,\bullet}[d]
 \bigr)
\simeq
 \distribution_X^{\bullet}[2d]$,
where the latter is given by 
$\eta^{d+j}\otimes
 \etabar^{d+i}
\longmapsto
 \epsilon(d)(-1)^{dj}
 \eta^{d+j}\wedge\etabar^{d+i}$.
Let $\Ftilde_e$ denote the image of $F_e$
in $\nhom_{\cnum_X}(\DR M_2,\distribution_X^{\bullet}[2d])^{-d}$.
It is given by
\[
 \Ftilde_e(\eta^{d+j}m)
=\epsilon(d)(-1)^d
 \eta^{d+j}
 \langle e,m\rangle
\]
which is naturally regarded as a section of
$(\distribution^{\bullet}_X[2d])^{j-d}$.
The morphism
$L_1[d]\lrarr
\nhom_{\cnum_X}(\DR M_2,\distribution_X^{\bullet}[2d])^{-d}$
is extended to the morphism
$\DR M_1
\lrarr
\nhom_{\cnum_X}(\DR M_2,\distribution_X^{\bullet}[2d])^{-d}$
induced by the multiplication
$\DR M_1\times \DR M_2
\lrarr
 \distribution_X^{\bullet}[2d]$
given by
$(\xi^{d+k}e,\eta^{d+j}m)
\longmapsto
 \epsilon(d)(-1)^{kd}
 \xi^{d+k}\eta^{d+j}
 \langle e,m\rangle$.
Then,
we obtain  the claim of Proposition \ref{prop;14.11.15.1}
in the case where $S$ is a point.
By Lemma \ref{lem;14.11.15.20},
the proof of Proposition \ref{prop;14.11.15.1}
is finished.
\hfill\qed

\section{Pairings and their functoriality}
\label{section;14.12.26.11}
\subsection{Mixed TEP-structures and their functoriality}
\label{subsection;14.12.9.2}

\subsubsection{Pairing of weight $w$ and  graded pairing}
\label{subsection;14.11.28.1}

Let $X$ be a complex manifold with a hypersurface $H$.
Set $\nbigx:=\cnum_{\lambda}\times X$.
For any $\nbigo_{\nbigx}$-module $N$,
let $N(\ast H):=
 N\otimes
 \nbigo_{\nbigx}\bigl(
 \ast(\cnum_{\lambda}\times H)
 \bigr)$.
An $\nbigr_X(\ast H)$-module $\nbigv$
is called smooth
if (i) $\nbigv_{|\cnum_{\lambda}\times (X\setminus H)}$
is locally free 
as an $\nbigo_{\cnum_{\lambda}\times (X\setminus H)}$-module,
(ii) $\nbigv$ is $\nbigo_{\nbigx}(\ast H)$-coherent and reflexive.
For a smooth $\nbigr_X(\ast H)$-module $\nbigv$,
a pairing of weight $w$ on  $\nbigv$ 
is defined to be an $\nbigr_X$-homomorphism 
 $P:\nbigv\otimes j^{\ast}\nbigv\lrarr 
 \lambda^{-w}\nbigo_{\nbigx}(\ast H)$
such that 
(i) $j^{\ast}P=(-1)^{w}P\circ\exchange$,
(ii) $P$ is non-degenerate,
i.e.,
the induced morphism
$\nbigv\lrarr
 \nhom_{\nbigo_{\nbigx}(\ast H)}
 \bigl(j^{\ast}\nbigv,\lambda^{-w}\nbigo_{\nbigx}(\ast H)\bigr)$
is an isomorphism.
Such a pair $(\nbigv,P)$ is a ${\rm TP}(-w)$-structure
in the sense of \cite{Hertling}.
Note that $P$ naturally induces
$P(\ast\lambda):\nbigv(\ast\lambda)\otimes
 j^{\ast}\nbigv(\ast\lambda)
\lrarr \nbigo_{\nbigx}(\ast H)(\ast\lambda)$,
and $P$ is the restriction of $P(\ast\lambda)$
to $\nbigv\otimes j^{\ast}\nbigv$.
We also obtain pairings
$P^{(a)}$ of weight $w+2a$ on $\lambda^{-a}\nbigv$
as the restriction of $P(\ast\lambda)$.
We have
$P^{(a)}(\lambda^{-a}v_1,j^{\ast}(\lambda^{-a}v_2))
=(-1)^a\lambda^{-2a}P(v_1,j^{\ast}v_2)$.

Suppose that 
a smooth $\nbigr_X(\ast H)$-module
$\nbigv$ is equipped with 
an exhaustive increasing filtration $\Wtilde$ 
of $\nbigr_X(\ast H)$-submodules
indexed by integers
such that $\Gr^{\Wtilde}_w(\nbigv)$ are 
also smooth $\nbigr_X(\ast H)$-modules.
A graded pairing of $(\nbigv,\Wtilde)$ is defined to be
a tuple of pairings $P_w$ $(w\in\seisuu)$
of weight $w$ on $\Gr^{\Wtilde}_w(\nbigv)$.
Such a filtered smooth $\nbigr_X(\ast H)$-module
$(\nbigv,W)$ with a graded pairing $\{P_w\,|\,w\in\seisuu\}$
is called a mixed TP-structure.
For any integer $a$,
we set
$\Wtilde_{k+2a}(\lambda^{-a}\nbigv)
:=\lambda^{-a}\Wtilde_k\nbigv$.
We naturally have
$\Gr^{\Wtilde}_{k+2a}(\lambda^{-a}\nbigv)
=\lambda^{-a}\Gr^{\Wtilde}_k(\nbigv)$.
The pairing $P_k$ on $\Gr^{\Wtilde}_k(\nbigv)$
induces 
$P_{k+2a}^{(a)}$
on $\Gr^{\Wtilde}_{k+2a}(\lambda^{-a}\nbigv)$.
Thus, we obtain another mixed TP-structure
which consists of 
a filtered  $\nbigr_X(\ast H)$-module
$(\lambda^{-a}\nbigv,\Wtilde)$
and a graded pairing  $(P^{(a)}_w\,|\,w\in\seisuu)$.

A mixed TP-structure 
$(\nbigv,\Wtilde,P)$ is called integrable if 
(i) $\nbigv$ is an $\nbigrtilde_X(\ast H)$-module,
(ii) $\Wtilde$ is a filtration by 
$\nbigrtilde_X(\ast H)$-submodules,
(iii) $P_w$ are also 
$\nbigrtilde_X(\ast H)$-homomorphisms.
In that case,
$(\nbigv,\Wtilde,P)$ is also called a mixed TEP-structure.

\subsubsection{Specialization of pairings}
\label{subsection;14.11.28.110}

Suppose that $Y$ is a hypersurface of $X$
given as $\{t=0\}$ for a holomorphic function $t$
such that $dt$ is nowhere vanishing.
We introduce a procedure for pairings
in the context of $\nbigrtilde_X(\ast Y)$-modules,
which we call specialization.
The procedure can also work in the non-integrable case.

We set $\nbigx:=\cnum_{\lambda}\times X$
and $\nbigy:=\cnum_{\lambda}\times Y$.
Let $p_{\lambda}:\nbigx\lrarr X$ be the projection.
Let $V\nbigrtilde_X\subset\nbigrtilde_X$
denote the sheaf of subalgebras generated by
$\lambda p_{\lambda}^{\ast}\Theta_X(\log Y)$
and $\lambda^2\del_{\lambda}$
over $\nbigo_{\nbigx}$.

Let $\nbigv$ be a smooth $\nbigrtilde_X(\ast Y)$-module.
We assume the following:
\begin{description}
\item[(B0)]
 $\nbigv$  is regular along $t$,
 i.e., there exists a filtration 
 $V_{\bullet}(\nbigv)=(V_a(\nbigv)\,|\,a\in\real)$ 
 by $V\nbigrtilde_X$-submodules
 such that 
 (i) $V_a(\nbigv)$ are locally free $\nbigo_{\nbigx}$-submodules
 and $\Gr^V_a(\nbigv)$ are locally free
 $\nbigo_{\nbigy}$-modules,
 (ii) $\lambda\del_tV_{a}\subset V_{a+1}$
 and $tV_a(\nbigv)= V_{a-1}(\nbigv)$,
 (iii) $\lambda \del_tt+\lambda a$ are nilpotent
 on $\Gr^V_a(\nbigv)$.
\end{description}
We obtain $\nbigrtilde_Y$-modules $\Gr^V_a(\nbigv)$
which are smooth.
We are particularly interested in 
$\Gr^V_{-1}(\nbigv)$.
It is also denoted as $\psi_{t,-\vecdelta}(\nbigv)$.
We also assume the following:
\begin{description}
\item[(B1)]
 The conjugacy classes of 
$(\lambda N)_{|(\lambda,W)}:
 \psi_{t,-\vecdelta}(\nbigv)_{(\lambda,Q)}
\lrarr
 \psi_{t,-\vecdelta}(\nbigv)_{(\lambda,Q)}$
are independent of the choice of 
$(\lambda,Q)\in \nbigy$.
Here, 
$N:\psi_{t,-\vecdelta}(\nbigv)
 \lrarr \psi_{t,-\vecdelta}(\nbigv)\lambda^{-1}$
denote the morphism induced by 
the multiplication of $t\del_t$.
\end{description}

Let $P$ be a pairing of weight $w$ on $\nbigv$.
We have the induced morphism
$P:V_{-1}\nbigv\otimes j^{\ast}V_{-1}\nbigv
\lrarr \lambda^{-w}\nbigo_{\nbigx}$.
Then, we have the naturally induced pairing of weight $w$:
\[
 \psi_{t,-\vecdelta}(P):
 \psi_{t,-\vecdelta}(\nbigv)
\otimes
 j^{\ast}
 \psi_{t,-\vecdelta}(\nbigv)
\lrarr
 \lambda^{-w}\nbigo_{\nbigy}
\]
Let $W(N)$ denote the monodromy weight filtration of
$\lambda N$ on $\psi_{t,-\vecdelta}(\nbigv)$.
By the assumption,
$W(N)$ is a filtration by subbundles.
We have the induced non-degenerate pairings:
\[
 \psi_{t,-\vecdelta}(P)_k:
 \Gr^{W(N)}_{-k}\psi_{t,-\vecdelta}(\nbigv)
\times
 j^{\ast}\Gr^{W(N)}_{k}\psi_{t,-\vecdelta}(\nbigv)
\lrarr
 \lambda^{-w}\nbigo_{\nbigy}
\]
For $k\geq 0$,
let $P\Gr^W_k\psi_{t,-\vecdelta}(\nbigv)$
denote the primitive part,
i.e.,
the kernel of 
$N^{k+1}:
 \Gr^W_k\psi_{t,-\vecdelta}(\nbigv)
\lrarr
 \lambda^{-k-1}\Gr^W_{-k-2}\psi_{t,-\vecdelta}(\nbigv)$.
We have the induced pairings:
\[
 \sp_t(P)_{w+k}:=\psi_{t,-\vecdelta}(P)_k\circ
 (N^k\times \id):
 P\Gr^{W(N)}_k\psi_{t,-\vecdelta}(\nbigv)
\times
 j^{\ast}P\Gr^{W(N)}_{k}\psi_{t,-\vecdelta}(\nbigv)
\lrarr
 \lambda^{-w-k}\nbigo_{\nbigy}
\]
By using
$\psi_{t,-\vecdelta}(P)\circ(\id\times j^{\ast}N)
=-\psi_{t,-\vecdelta}(P)\circ(N\times\id)$,
we can easily check that
$\sp_t(P)_{w+k}$ is a pairing of weight $w+k$
on 
$P\Gr^{W(N)}_k\psi_{t,-\vecdelta}(\nbigv)$.

\vspace{.1in}
Let $k\leq 0$.
Let $P'\Gr^W_{k}\psi_{t,-\vecdelta}(\nbigv)\subset
 \Gr^W_{k}\psi_{t,-\vecdelta}(\nbigv)$
denote the image of 
$\lambda^{-k}P\Gr^W_{-k}\psi_{t,-\vecdelta}(\nbigv)$
by $N^{-k}$.
We have the isomorphism
\[
 \id\times j^{\ast}N^{-k}:
 P'\Gr^W_{k}\psi_{t,-\vecdelta}(\nbigv)
 \times
 j^{\ast} \bigl(
 \lambda^{-k}P\Gr^W_{-k}\psi_{t,-\vecdelta}(\nbigv) 
 \bigr)
\simeq
 P'\Gr^W_{k}\psi_{t,-\vecdelta}(\nbigv)
 \times
 j^{\ast}P'\Gr^W_{k}\psi_{t,-\vecdelta}(\nbigv).
\]
We obtain the induced pairing of weight $w+k$.
\[
 \sp_t(P)_{w+k}:=
 \psi_{t,-\vecdelta}(P)_{k}
 \circ(\id\times j^{\ast}N^{-k})^{-1}:
  P'\Gr^W_{k}\psi_{t,-\vecdelta}(\nbigv)
 \times
 j^{\ast}P'\Gr^W_{k}\psi_{t,-\vecdelta}(\nbigv)
\lrarr
 \lambda^{-k-w}\nbigo_{\nbigy}
\]
It is also induced by
$(-1)^k\psi_{t,-\vecdelta}(P)_{-k}\circ
 (N^k\times\id)^{-1}$.

\subsubsection{Filtration and graded pairings}
\label{subsection;14.11.28.130}

We continue to use the notation in \S\ref{subsection;14.11.28.110}.
Let $W(N)$ denote the filtration on 
$\Cok(\lambda N)$
induced by the filtration $W(N)$
on $\psi_{t,-\vecdelta}(\nbigv)$.
It is equipped with the flat connection.
We set
\[
 \Wtilde_{k}\Cok(\lambda N):=
 W(N)_{k-w}\Cok(\lambda N).
\]
In other words,
$\Wtilde_k\Cok(\lambda N)$
is the image of 
$W(N)_{k-w}\psi_{t,-\vecdelta}(\nbigv)
\lrarr
 \Cok(\lambda N)$.
We have 
\[
\Gr^{\Wtilde}_k\Cok(\lambda N)\simeq
\left\{
\begin{array}{ll}
 P\Gr^{W(N)}_{k-w}\psi_{t,-\vecdelta}(\nbigv) 
 &  (k\geq w)\\
 0 & (k<w)
\end{array}
\right.
\]
Hence, we have the pairing of weight $k$ 
on $\Gr^{\Wtilde}_k\Cok(\lambda N)$
induced by $\sp_t(P)_{k}$.
The induced pairing is also denoted by
$\sp_t(P)_{k}$.
The tuple $\bigl(\sp_t(P)_k\,\big|\,k\in\seisuu \bigr)$
is denoted by $\sp_t(P)$.
Thus, we obtain a mixed TEP-structure
$(\Cok(\lambda N),\Wtilde,\sp_t(P))$.

Similarly,
we set 
$\Wtilde_k\Ker(\lambda N):=
 W(N)_{k-w}\cap
 \Ker(\lambda N)$.
We have
\[
 \Gr^{\Wtilde}_k\Ker(\lambda N)
\simeq
 \left\{
 \begin{array}{ll}
 P'\Gr^{W(N)}_{k-w}\psi_{t,-\vecdelta}(\nbigv)
 & (k\leq w)\\
 0 & (k>w)
 \end{array}
 \right.
\]
Hence, we have the pairing $\sp_t(P)_k$ of weight $k$
on $\Gr^{\Wtilde}_k\Ker(\lambda N)$.
The induced pairing is also denoted by 
$\sp_t(P)_k$.
Thus, we obtain a mixed TEP-structure
$(\Ker(\lambda N),\Wtilde,\sp_t(P))$,
where 
$\sp_t(P):=\bigl(\sp_t(P)_k\,\big|\,k\in\seisuu\bigr)$.

\subsubsection{Dependence on defining functions}

Let $\tau$ be a holomorphic function on $X$
such that
(i) $Y=\{\tau=0\}$,
(ii) $d\tau$ is nowhere vanishing.
The $V$-filtrations for $t$ and $\tau$ are the same,
and we naturally have
$\psi_{\tau,-\vecdelta}(\nbigv)
=\psi_{t,-\vecdelta}(\nbigv)$
as $\nbigo_{\nbigy}$-modules.
Let $N_t$ and $N_{\tau}$ denote the morphisms 
$\psi_{\tau,-\vecdelta}(\nbigv)
\lrarr
 \lambda^{-1}\psi_{\tau,-\vecdelta}(\nbigv)$
induced by $t\del_t$ and $\tau\del_{\tau}$.
We remark the following standard lemma.
\begin{lem}
We have $N_t=N_{\tau}$
and $\psi_{t,-\vecdelta}(P)=\psi_{\tau,-\vecdelta}(P)$.
As a result,
we have $\Cok(N_t)=\Cok(N_{\tau})$
and $\Ker(N_t)=\Ker(N_{\tau})$.
On $\Cok(N_t)=\Cok(N_{\tau})$
and $\Ker(N_t)=\Ker(N_{\tau})$,
the induced connections are also equal.
Moreover,
the induced filtrations and the graded pairings
are also equal.
\hfill\qed
\end{lem}

In other words,
the induced mixed TEP-structures
are well defined in the sense that
they are independent of the choice of
a defining function of $Y$.
Note that the induced connections are not the same
on $\psi_{t,-\vecdelta}(\nbigv)=\psi_{\tau,-\vecdelta}(\nbigv)$,
in general.

\subsection{Pairings associated with real structure and graded
 sesqui-linear duality} 
\label{subsection;14.12.9.1}

\subsubsection{Compatibility of sesqui-linear duality and real structure}

Let $X$ be a complex manifold.
Let $\nbigt=(\nbigm',\nbigm'',C)$
be a pure twistor $\nbigd$-module of weight $w$
on $X$.
Let $\nbigs:\nbigt\lrarr \nbigt^{\ast}(-w)$ be a sesqui-linear duality.
Let $\kappa:\gammatilde^{\ast}\nbigt\lrarr \nbigt$
be a real structure.
We say that $\nbigs$ and $\kappa$ are compatible
if the following diagram is commutative.
\[
 \begin{CD}
 \gammatilde^{\ast}\nbigt
 @>>>
 \gammatilde^{\ast}(\nbigt^{\ast})(-w)
 \\
 @V{\kappa}VV @A{\kappa^{\ast}}AA
 \\
 \nbigt
 @>{\nbigs}>>
 \nbigt^{\ast}(-w)
 \end{CD}
\]
Recall that 
$\nbigs$ is expressed as a pair of morphisms
$\nbigm'\stackrel{\nbigs'}{\llarr}\nbigm''\lambda^w$
and 
$\nbigm''\stackrel{\nbigs''}{\lrarr}\nbigm'\lambda^{-w}$
satisfying $\nbigs'=\nbigs''$,
and that $\kappa$ is expressed as a pair of morphisms
$\DDD_X j^{\ast}\nbigm''\stackrel{\kappa'}{\llarr}\nbigm'$
and 
$\DDD_X j^{\ast}\nbigm'\stackrel{\kappa''}{\lrarr}\nbigm''$
satisfying
$j^{\ast}\DDD_X\kappa''\circ\kappa'=\id$
and 
$\kappa''\circ j^{\ast}\DDD_X\kappa'=\id$.
Here, $\DDD_X$ denotes the duality functor
for $\nbigr_X$-modules.
Note that 
$\gammatilde^{\ast}\bigl(\nbigt(-w)\bigr)
\simeq
 \gammatilde^{\ast}(\nbigt)(-w)$
is given by
$\bigl((-1)^w,(-1)^w\bigr)$.
Then, the commutativity of the above diagram
means the commutativity of the following:
\[
 \begin{CD}
 j^{\ast}\DDD_X(\nbigm')
 @>{(-1)^wj^{\ast}\DDD_X\nbigs''}>>
 j^{\ast}\DDD_X(\nbigm'')\lambda^{-w}
 \\
 @V{\kappa''}VV @A{\kappa'}AA \\
 \nbigm''
 @>{\nbigs''}>>
 \nbigm'\lambda^{-w}
 \end{CD}
\]
Because
$\kappa''=(j^{\ast}\DDD_X\kappa')^{-1}$,
the commutativity means
$\kappa'\circ\nbigs''
 =(-1)^wj^{\ast}\DDD_X(\kappa'\circ\nbigs'')$.

\vspace{.1in}
Let $(\nbigt,W)$ be a mixed twistor $\nbigd$-module
with a real structure $\kappa$
and a graded sesqui-linear duality
$(\nbigs_w\,|\,w\in\seisuu)$.
We say that the real structure and the graded sesqui-linear duality
are compatible
if the induced real structure and the sesqui-linear duality
on $\Gr^W_w(\nbigt)$ are compatible for any $w$.

\subsubsection{The associated TEP-structure in the pure case}
\label{subsection;14.11.8.120}

Let $\nbigt=(\nbigm',\nbigm'',C)$ be 
an integrable pure twistor $\nbigd$-module
of weight $w$
on a complex manifold $X$
with a real structure $\kappa$
and a polarization $\nbigs$ which are compatible
and integrable.
We have the following isomorphisms
of $\nbigrtilde_X$-modules:
\[
\begin{CD}
 \nbigm''
 @>{\nbigs''}>>
 \nbigm'\lambda^{-w}
 @>{\kappa'}>>
 \lambda^{-w}j^{\ast}\DDD_X\nbigm''
\end{CD}
\]

Let $H$ be a hypersurface of $X$.
Suppose that
$\nbigm'(\ast H)$ and $\nbigm''(\ast H)$
are smooth $\nbigr_X(\ast H)$-modules.
We have a natural isomorphism
\[
  \nu_{\nbigm''}:
 (\DDD_X\nbigm'')(\ast H)
 \simeq 
\lambda^{d_X}\cdot \nbigm''(\ast H)^{\lor}
\]
whose restriction to $\{\lambda_0\}\times X$ $(\lambda_0\neq 0)$
is given by the morphism in \S\ref{subsection;14.11.14.11}.
We consider the pairing $P$
induced by $(-1)^{d_X}\nu_{\nbigm''}$
and $\kappa'\circ\nbigs''$:
\[
\begin{CD}
P:
 \nbigm''(\ast H)\times j^{\ast}\nbigm''(\ast H)
@>>>
 \lambda^{d_X-w}\nbigo_{\nbigx}(\ast H)
\end{CD}
\]

\begin{prop}
\label{prop;14.11.8.40}
We have
$j^{\ast}P\circ\exchange
=(-1)^{d_X-w}P$.
In other words,
$P$ is a pairing of weight $w-d_X$
on $\nbigm''(\ast H)$.
In other words,
$(\nbigm''(\ast H),P)$
is a $TEP(d_X-w)$-structure.
\end{prop}
\pf
The claim of this proposition
follows from the relation
$\kappa'\circ\nbigs''
=(-1)^{w}j^{\ast}\DDD(\kappa'\circ\nbigs'')$
and Lemma \ref{lem;14.11.8.41}.
\hfill\qed

\vspace{.1in}

For any $n\in\seisuu$,
the pure twistor $\nbigd$-module
$\nbigt\otimes\newTate(-n)$
of weight $w+2n$ naturally equipped with 
the polarization and the real structure
which are compatible.
The polarization is given by
$((-1)^n\nbigs',(-1)^n\nbigs'')$.
The real structure is given by
$((-1)^n\kappa',(-1)^n\kappa'')$.
Hence, the induced pairing is $P^{(n)}$
which was introduced in \S\ref{subsection;14.11.28.1}.

\begin{example}
Let $\nbigt(F)=\bigl(
 \lambda^{d_X}\nbigl(F),\nbigl(F),C
 \bigr)$ 
be the integrable pure twistor $\nbigd$-module
of weight $d_X$ associated
to a holomorphic function $F$
on $X$,
where $d_X:=\dim X$.
The polarization $\nbigs$ is given by
$\bigl((-1)^{d_X},(-1)^{d_X}\bigr)$.
The real structure 
$\gammatilde^{\ast}\nbigt(F)
 \simeq\nbigt(F)$
is given by 
$\bigl(\nu^{-1},(-1)^{d_X}\nu\bigr)$.
The induced pairing
$\nbigl(F)\times
 j^{\ast}\nbigl(F)
\lrarr \nbigo_{\nbigx}$
is just the multiplication.
\hfill\qed
\end{example}

\subsubsection{The associated mixed TEP-structure}
\label{subsection;14.11.28.10}

Let $(\nbigt,W)$ be a mixed twistor $\nbigd$-module on $X$.
Let $(\nbigm',\nbigm'',C)$ be the underlying $\nbigr_X$-triple
of $\nbigt$.
We often pick $\nbigm''$,
and we say that $\nbigm''$ is the underlying
$\nbigr_X$-module of $\nbigt$.
We have the filtration $W$ on $\nbigm''$
such that $W_j\nbigm''$ are the underlying 
$\nbigr$-modules of $W_j\nbigt$,
which we call the filtration of $\nbigm''$
underlying the weight filtration of $(\nbigt,W)$.
In the following, we shall often use
the shifted filtration $\Wtilde$ on $\nbigm''$
given by
$\Wtilde_{k}(\nbigm''):=W_{k+\dim X}\nbigm''$.

\vspace{.1in}

Let $(\nbigt,W)$ be an integrable mixed twistor $\nbigd$-module
with real structure and a graded sesqui-linear duality
which are compatible and integrable.
Let $H$ be a hypersurface of $X$.
Let $\nbigm''$ be the underlying $\nbigrtilde_X$-module.
Suppose that $\nbigm''(\ast H)$ is 
a smooth $\nbigr_X(\ast H)$-module.
Let $W$ be the filtration of $\nbigm''$
underlying the weight filtration of $\nbigt$.
We set 
$\Wtilde_k\bigl(\nbigm''(\ast H)\bigr):=
 W_{k+d_S}\nbigm''(\ast H)$.
Each $\Gr^{\Wtilde}_k(\nbigm''(\ast H))$
is equipped with the pairing $P_k$ of weight $k$
induced by the real structure and the sesqui-linear duality.
In this way,
we obtain a mixed TEP-structure
$(\nbigm''(\ast H),\Wtilde,\{P_k\})$.

\begin{rem}
The underlying $\nbigr$-module of
$\nbigt\otimes\newTate(-n)$
is $\lambda^{-n}\nbigm''(\ast H)$.
The filtration 
$\Wtilde(\lambda^{-n}\nbigm''(\ast H))$
induced by the mixed twistor structure
$(\nbigt,W)\otimes\newTate(-n)$
is given by
$\Wtilde_k(\lambda^{-n}\nbigm''(\ast H))
=\lambda^{-n}\Wtilde_{k-2n}\nbigm''(\ast H)$,
which is equal to the shift of the filtration
in {\rm\S\ref{subsection;14.11.28.1}}.
The graded pairing on 
$(\lambda^{-n}\nbigm''(\ast H),\Wtilde)$
induced by the real structure and the graded sesqui-linear duality
on $(\nbigt,W)\otimes\newTate(-n)$
is equal to that induced by the graded pairing
on $(\nbigm''(\ast H),\Wtilde)$ by the procedure 
in {\rm\S\ref{subsection;14.11.28.1}}.
\hfill\qed
\end{rem}

\subsection{Comparison of the specializations}
\label{subsection;14.11.28.30}

\subsubsection{Statement}
\label{subsection;14.11.28.40}

Let $X$ be a complex manifold
with a smooth hypersurface $Y$.
Suppose that $Y$ is given as  $\{t=0\}$
for a holomorphic function $t$
such that $dt$ is nowhere vanishing.
Let $\nbigt=(\nbigm_1,\nbigm_2,C)$ 
be an integrable pure twistor $\nbigd$-module of weight $w$
on $X$ with a sesqui-linear duality $\nbigs$
and a real structure $\kappa$
which are compatible and integrable.

\begin{assumption}
We assume that 
 $\nbigm_i(\ast Y)$
are smooth $\nbigr_X(\ast Y)$-modules,
and regular along $Y$.
(See the condition {\bf (B0)} in
{\rm\S\ref{subsection;14.11.28.110}}
 for the regularity.)
\hfill\qed
\end{assumption}
Set $d:=\dim X$.
We have the pairing $P$ of weight $w-d$
on $\nbigm_2(\ast Y)$.
We have the endomorphism
$\lambda N:
 \psi_{t,-\vecdelta}(\nbigm_2)
 \lrarr
 \psi_{t,-\vecdelta}(\nbigm_2)$,
where $N$ is the induced by $t\del_t$.
Note that the condition {\bf (B1)} is satisfied
because $\nbigm_2$ comes from 
a pure twistor $\nbigd$-module.
By the procedure in \S\ref{subsection;14.11.28.130},
we obtain filtrations $\Wtilde^{[1]}$
and graded pairings
$(P^{[1]}_k\,|\,k\in\seisuu):=\sp_t(P)$ 
on $\Cok(\lambda N)$
and $\Ker(\lambda N)$.

\vspace{.1in}
Let $\nbign:
 \psitilde_{t,-\vecdelta}(\nbigt)\otimes\nbigu(-1,0)
\lrarr
 \psitilde_{t,-\vecdelta}(\nbigt)\otimes\nbigu(0,-1)$
be given by $(-N,-N)$,
where $N:
 \psi_{t,-\vecdelta}(\nbigm_i)
\lrarr
 \lambda^{-1}\psi_{t,-\vecdelta}(\nbigm_i)$
are induced by $t\del_t$.
According to the isomorphism in 
Proposition 4.3.1 of \cite{Mochizuki-MTM},
we have the following commutative diagram:
\[
 \begin{CD}
 \psi^{(1)}_t(\nbigt)
@>>>
  \psi^{(0)}_t(\nbigt)
 \\
 @V{\simeq}VV @V{\simeq}VV \\
 \iota_{\dagger}\psitilde_{t,-\vecdelta}(\nbigt)
 \otimes\nbigu(-1,0)
 @>{\nbign}>>
  \iota_{\dagger}\psitilde_{t,-\vecdelta}(\nbigt)
 \otimes\nbigu(0,-1)
 \end{CD}
\]
Here, the upper horizontal arrow is the canonical morphism.
Hence, we have the following natural isomorphisms
of mixed twistor $\nbigd$-modules:
\[
 \iota_{\dagger}\Cok(\nbign)
 \simeq
 \Cok\bigl(
 \psi^{(1)}_t(\nbigt)
\lrarr
 \psi^{(0)}_t(\nbigt)
 \bigr)
\simeq
 \Cok\Bigl(
 \nbigt[!t]
\lrarr
  \nbigt[\ast t]
 \Bigr)
\]
\[
 \iota_{\dagger}\Ker(\nbign)
 \simeq
 \Ker\bigl(
 \psi^{(1)}_t(\nbigt)
\lrarr
 \psi^{(0)}_t(\nbigt)
 \bigr)
\simeq
 \Ker\Bigl(
 \nbigt[!t]
\lrarr
  \nbigt[\ast t]
 \Bigr)
\]
So, $\Cok(\nbign)\otimes\newTate(1)$
and $\Ker(\nbign)$
are equipped with the real structure
and the graded sesqui-linear duality
which are compatible and integrable.
The underlying $\nbigrtilde_Y$-modules of
$\Cok(\nbign)\otimes\newTate(1)$
and $\Ker(\nbign)$
are $\Cok(\lambda N)$ and $\Ker(\lambda N)$,
respectively.
Applying the procedure in \S\ref{subsection;14.11.28.10},
we obtain filtrations $\Wtilde^{[2]}$
and graded pairings $(P^{[2]}_k\,|\,k\in\seisuu)$
on $\Cok(\lambda N)$
and $\Ker(\lambda N)$.

\begin{prop}
\label{prop;14.11.28.32}
We have 
$\Wtilde^{[1]}=\Wtilde^{[2]}$
and $P^{[1]}_k=P^{[2]}_k$
for any $k\in\seisuu$.
In other words,
the induced mixed TEP-structures are the same.
\end{prop}

For the monodromy weight filtration $W(\nbign)$ of $\nbign$
on $\psitilde_{t,-\vecdelta}(\nbigt)$,
we have
$W(\nbign)_{j}\psitilde_{t,-\vecdelta}(\nbigt)
 =W_{j+w}\psitilde_{t,-\vecdelta}(\nbigt)$.
Then, we can check the equality for the filtrations easily.
We shall compare the pairings in 
\S\ref{subsection;14.11.8.10}--\ref{subsection;14.11.8.11}.

\subsubsection{Preliminary}
\label{subsection;14.11.8.10}

Let $\nbigv$ be a  smooth $\nbigrtilde_X(\ast Y)$-module
which is regular along $t=0$.
We have a natural isomorphism
$\nu_{\nbigv}:\DDD(\nbigv)(\ast Y)
\simeq
 \lambda^d\cdot\nbigv^{\lor}$
whose specialization at
$\{\lambda_0\}\times X$ $(\lambda_0\neq 0)$
is equal to the morphism 
in \S\ref{subsection;14.11.14.11}.
By a procedure similar to that
in \S\ref{subsection;14.11.15.20},
it induces the following isomorphisms
of $\nbigrtilde_X$-modules for $a=0,1$:
\[
 \nu_{\nbigv,\psi^{(a)}_t}:
 \DDD_X \psi_t^{(a)}(\nbigv)
\simeq
 \psi_t^{(-a+1)}(\nbigv^{\lor})
 \lambda^d
\]

\begin{lem}
\label{lem;14.11.15.31}
The following diagram is commutative:
\begin{equation}
\label{eq;14.11.15.30}
 \begin{CD}
 \DDD_X\psi_t^{(1)}(\nbigv)
 @>{\nu_{\nbigv,\psi_t^{(1)}}}>{\simeq}>
 \psi_t^{(0)}(\nbigv^{\lor})\lambda^d
 \\
 @V{\simeq}VV @V{\simeq}VV \\
 \DDD_X\iota_{\dagger}\psi_{t,-\vecdelta}(\nbigv)
@>{\simeq}>>
 \iota_{\dagger}\psi_{t,-\vecdelta}(\nbigv^{\lor})
 \lambda^{d-1}
 \end{CD}
\end{equation}
Here, the vertical arrows
are the isomorphisms in 
{\rm\S4.3} of {\rm\cite{Mochizuki-MTM}},
and the lower horizontal arrow is induced as
the composite of the isomorphisms
$\DDD_X(\iota_{\dagger}\psi_{t,-\vecdelta}(\nbigv))
\simeq
 \iota_{\dagger}\DDD_Y \psi_{t,-\vecdelta}(\nbigv)
\simeq
  \iota_{\dagger}\psi_{t,-\vecdelta}(\nbigv)^{\lor}\lambda^{d-1}
\simeq
 \iota_{\dagger}\psi_{t,-\vecdelta}(\nbigv^{\lor})
 \lambda^{d-1}$.

The following diagram is commutative:
\begin{equation}
\label{eq;14.11.15.31}
 \begin{CD}
 \DDD_X\bigl(
 \psi^{(0)}_t(\nbigv)
 \bigr)
 @>{-\nu_{\nbigv,\psi^{(0)}_t}}>{\simeq}>
 \psi_t^{(1)}(\nbigv^{\lor})\lambda^d
 \\
 @V{\simeq}VV @V{\simeq}VV \\
 \DDD_X(\iota_{\dagger}
 \psi_{t,-\vecdelta}(\nbigv) \lambda^{-1})
 @>{\simeq}>>
 \iota_{\dagger}
 \psi_{t,-\vecdelta}(\nbigv^{\lor})
 \lambda^d
 \end{CD}
\end{equation}
The vertical arrows and the lower horizontal arrows
are given as in the diagram 
{\rm (\ref{eq;14.11.15.30})}.
\end{lem}
\pf
We have only to prove the commutativity
of the diagrams
after taking the specializations
to $\{\lambda_0 \}\times X$
for any generic $\lambda_0\neq 0$,
which follows from Proposition \ref{prop;14.10.4.40}.
\hfill\qed

\vspace{.1in}

Let $\nbigv_i$ $(i=1,2)$
be smooth $\nbigrtilde_X(\ast Y)$-modules
which are regular along $t=0$.
Let $P':\nbigv_1\otimes j^{\ast}\nbigv_2\lrarr 
 \lambda^{-m}\nbigo_{\nbigx}(\ast Y)$
be a morphism of $\nbigrtilde_X(\ast Y)$-modules.
We have the induced morphism
of $\nbigrtilde_Y$-modules
\[
\psi_{t,-\vecdelta}(P'):
 \psi_{t,-\vecdelta}(\nbigv_1)
\otimes
 j^{\ast}
 \psi_{t,-\vecdelta}(\nbigv_2)
\lrarr
 \lambda^{-m}
 \nbigo_{\nbigy}.
\]

The pairing $P'$ and 
the isomorphism $(-1)^{d}\nu_{\nbigv_2}$
induce
an $\nbigrtilde_X$-homomorphism 
$\nbigv_1\lrarr 
 \lambda^{-m-d}(j^{\ast}\DDD_X\nbigv_2)(\ast Y)$.
It induces the following morphisms
for $a=0,1$:
\[
 \psi^{(a)}_t(\nbigv_1)
\lrarr
 \lambda^{-m-d}
 j^{\ast}\DDD_X\psi^{(1-a)}_t(\nbigv_2)
\]
By the isomorphisms
$\lambda^{-1+a}\cdot
 \iota_{\dagger}\psi_{t,-\vecdelta}(\nbigv_1)
\simeq
 \psi^{(a)}_t(\nbigv_1)$
and 
$\lambda^{-a}\iota_{\dagger}\psi_{t,-\vecdelta}(\nbigv_2)
\simeq
 \psi^{(1-a)}_t(\nbigv_2)$,
we obtain the following:
\begin{equation}
 \label{eq;14.11.15.40}
 \iota_{\dagger}\psi_{t,-\vecdelta}(\nbigv_1)
\lrarr 
 \lambda^{-m-d+1}
 j^{\ast}\DDD_X\bigl(
 \iota_{\dagger}\psi_{t,-\vecdelta}(\nbigv_2)
 \bigr)
\simeq
 \lambda^{-m-d+1}
 \iota_{\dagger}j^{\ast}\DDD_Y\psi_{t,-\vecdelta}(\nbigv_2)
\end{equation}
The morphism (\ref{eq;14.11.15.40})
and $(-1)^{d-1}\nu_{\psi_{t,-\vecdelta}(\nbigv_2)}$
induce a pairing
$P'_a:\psi_{t,-\vecdelta}(\nbigv_1)
 \times
 j^{\ast}\psi_{t,-\vecdelta}(\nbigv_2)
\lrarr
 \lambda^{-m}\nbigo_{\nbigy}$.
We obtain the following lemma from 
Lemma \ref{lem;14.11.15.31},
as in the case of Corollary \ref{cor;14.11.15.32}.
We remark the twist of the signatures,
i.e.,
we use $(-1)^{d}\nu_{\nbigv}$
and $(-1)^{d-1}\nu_{\psi_{t,-\vecdelta}(\nbigv_2)}$.
\begin{lem}
\label{lem;14.11.15.100}
We have
$\psi_{t,-\vecdelta}(P')=P'_1$
and 
$\psi_{t,-\vecdelta}(P')=-P'_0$.
\hfill\qed
\end{lem}

\subsubsection{Proof of Proposition  \ref{prop;14.11.28.32}}
\label{subsection;14.11.8.11}

For $a=0,1$,
from the morphisms
 $\psi^{(a)}_t(\nbigs):
 \psi_t^{(a)}(\nbigt)
\lrarr
 \psi_t^{(a)}(\nbigt)^{\ast}\otimes\newTate(-w-1+2a)$
and 
$\psi_t^{(a)}(\gammatilde^{\ast}):
 \gammatilde^{\ast}\psi_t^{(a)}(\nbigt)
 \simeq
\psi_t^{(a)}(\nbigt)$,
we have the following $\nbigrtilde_X$-homomorphisms:
\begin{equation}
 \label{eq;14.11.15.110}
 \psi^{(a)}_t(\nbigm_2)
\stackrel{b^{(a)}_1}{\lrarr}
 \lambda^{-w-1+2a}\psi^{(1-a)}_t(\nbigm_1)
\stackrel{b^{(a)}_2}{\lrarr}
 \lambda^{-w-1+2a}
 j^{\ast}\DDD_X\psi^{(a)}_t(\nbigm_2)
\end{equation}
The composite 
$b_2^{(a)}\circ b_1^{(a)}$
comes from an isomorphism
$\psi_{t,-\vecdelta}(\nbigm_2)\lambda^{-1}
\simeq
  \lambda^{-w-1}
 j^{\ast}\DDD_X
 \bigl(\psi_{t,-\vecdelta}(\nbigm_2)
 \lambda^{-1}\bigr)$,
in the case $a=0$,
or 
$\psi_{t,-\vecdelta}(\nbigm_2)
\simeq
  \lambda^{-w+1}
 j^{\ast}\DDD_X
 \bigl(\psi_{t,-\vecdelta}(\nbigm_2)\bigr)$
in the case $a=1$.
Together with
$(-1)^{d-1}\nu_{\psi_{t,-\vecdelta}(\nbigm_2)}$,
we obtain a pairing
$\Ptilde_a:
\psi_{t,-\vecdelta}(\nbigm_2)
\times
j^{\ast}\psi_{t,-\vecdelta}(\nbigm_2)
\lrarr
 \lambda^{d-w}\nbigo_{\nbigy}$.

\begin{lem}
\label{lem;14.11.15.120}
We have $\psi_{t,-\vecdelta}(P)=\Ptilde_0=\Ptilde_1$.
\end{lem}
\pf
According to \S4.3 of \cite{Mochizuki-MTM},
we have the following commutative diagram:
\begin{equation}\label{eq;14.11.29.10}
\begin{CD}
 \psi^{(0)}_t(\nbigt)
@>{\psi_t^{(0)}(\nbigs)}>>
 \psi^{(0)}_t(\nbigt)^{\ast}
 \otimes\newTate(-w-1)
 \\
 @V{\simeq}VV @V{\simeq}VV \\
 \iota_{\dagger}\psi_{t,-\vecdelta}(\nbigt)
 \otimes\nbigu(0,-1)
 @>>>
 \Bigl(
 \iota_{\dagger}\psi_{t,-\vecdelta}(\nbigt)
 \otimes\nbigu(0,-1)
 \Bigr)^{\ast}
\otimes\newTate(-w-1)
\end{CD}
\end{equation}
The lower horizontal arrow is induced by
$\psi_{t,-\vecdelta}(\nbigs)$.
It means that the following diagram is commutative:
\begin{equation}
 \label{eq;14.11.15.121}
 \begin{CD}
 \psi^{(0)}_t(\nbigm_2)
 @>{b_1^{(0)}}>>
 \psi^{(1)}_t(\nbigm_1) \lambda^{-w-1}
 \\
 @V{\simeq}VV @V{\simeq}VV \\
 \iota_{\dagger}
 \lambda^{-1}\psi_{t,-\vecdelta}(\nbigm_2)
 @>{\psi_{t,-\vecdelta}(\nbigs'')}>>
 \iota_{\dagger}\psi_{t,-\vecdelta}(\nbigm_1)
 \lambda^{-w-1}
 \end{CD}
\end{equation}

We have the pairing
$P':\nbigm_1(\ast Y)
 \times
 j^{\ast}\nbigm_2(\ast Y)
\lrarr
 \lambda^{d}
 \nbigo_{\nbigx}(\ast Y)$
induced by
$\kappa':\nbigm_1\lrarr j^{\ast}\DDD\nbigm_2$
and 
$(-1)^{d}\nu_{\nbigm_2(\ast Y)}$.
The morphism $b^{(0)}_2$
comes from
$\psi_{t,-\vecdelta}(\nbigm_1)
\lrarr
 j^{\ast}\DDD(\psi_{t,-\vecdelta}(\nbigm_2)\lambda^{-1})$.
Together with the morphism
$(-1)^{d-1}\nu_{\psi_{t,-\vecdelta}(\nbigm_2)}$,
we obtain a pairing
$P'_1:
 \psi_{t,-\vecdelta}(\nbigm_1)
\times
 j^{\ast}\psi_{t,-\vecdelta}(\nbigm_2)
\lrarr \lambda^{d}\nbigo_{\nbigy}$.
We obtain $\psi_{t,-\vecdelta}(P)=\Ptilde_0$
from the equality 
$\psi_{t,-\vecdelta}(P')=P'_1$
in Lemma \ref{lem;14.11.15.100},
and the commutativity 
(\ref{eq;14.11.15.121}).

We also have the following commutative diagram
from (\ref{eq;14.11.29.10}):
\begin{equation}
 \label{eq;14.11.29.11}
 \begin{CD}
 \psi^{(1)}_t(\nbigm_2)
 @>{b_1^{(1)}}>>
 \psi^{(0)}_t(\nbigm_1) \lambda^{-w+1}
 \\
 @V{\simeq}VV @V{\simeq}VV \\
 \iota_{\dagger}
  \psi_{t,-\vecdelta}(\nbigm_2)
 @>{-\psi_{t,-\vecdelta}(\nbigs'')}>>
 \iota_{\dagger}\psi_{t,-\vecdelta}(\nbigm_1)
 \lambda^{-w}
 \end{CD}
\end{equation}
The morphism  $b_2^{(1)}$ comes from
a morphism 
$\lambda^{-1}\psi_{t,-\vecdelta}(\nbigm_1)
\lrarr
 j^{\ast}\DDD\bigl(
 \psi_{t,-\vecdelta}(\nbigm_2)
 \bigr)$.
It induces $P_0':
\psi_{t,-\vecdelta}(\nbigm_1)
\times
 j^{\ast}\psi_{t,-\vecdelta}(\nbigm_2)
\lrarr \lambda^d\nbigo_{\nbigy}$.
From the equality
$\psi_{t,-\vecdelta}(P')=-P_0'$
and the commutativity of the diagram (\ref{eq;14.11.29.11}),
we obtain 
$\psi_{t,-\vecdelta}(P)=\Ptilde_1$.
Thus, the proof of 
Lemma \ref{lem;14.11.15.120}
is finished.
\hfill\qed

\vspace{.1in}

The graded sesqui-linear duality
on $\Gr^W_{w+1+k}\iota_{\dagger}\Cok(\nbign)$
is induced by
$\psi_t^{(0)}(\nbigs)\circ(-\nbign)^k$.
Hence, by Lemma \ref{lem;14.11.15.120},
we obtain $P^{[1]}=P^{[2]}$ on 
$\Cok(\nbign)$.
The graded sesqui-linear duality 
on $\Gr^{W}_{w-1-k}\iota_{\dagger}\Ker(\nbign)$
is induced by
$\psi_t^{(1)}(\nbigs)\circ \nbign^k$.
Hence, by Lemma \ref{lem;14.11.15.120},
we obtain 
$P^{[1]}=P^{[2]}$ on 
$\Ker(\nbign)$.
Thus the proof of Proposition \ref{prop;14.11.28.32}
is finished.
\hfill\qed

\subsection{Push-forward}
\label{subsection;14.12.10.2}

Let $X$ be a complex manifold.
We set $d_X:=\dim X$.
Let $\nbigt=(\nbigm_1,\nbigm_2,C)$
be an integrable pure twistor $\nbigd$-module 
of weight $w$ on $X$
with a sesqui-linear duality $\nbigs$ and a real structure $\kappa$
which are compatible and integrable.
Suppose that $\nbigt$ is smooth,
i.e., $\nbigm_i$ are locally free $\nbigo_{\nbigx}$-modules.
We have the associated pairing
$P(\nbigs,\kappa)$ on $\nbigm_2$
of weight $w-d_X$.

Let $F:X\lrarr Y$ be a smooth projective morphism.
We set $d_Y:=\dim Y$.
We have the pure twistor $\nbigd$-modules
$F_{\dagger}^i(\nbigt)
=(F_{\dagger}^{-i}\nbigm_1,F_{\dagger}^i\nbigm_2,F_{\dagger}^iC)$.
They are equipped with the induced real structure $\kappa_i$.
We also have the induced morphisms
$\nbigs_i:
 F_{\dagger}^i\nbigt
\lrarr
\bigl(F_{\dagger}^{-i}\nbigt\bigr)^{\ast}
 \otimes\newTate(-w)$.
Note that $F_{\dagger}^j\nbigm_i$ are 
locally free $\nbigo_{\nbigy}$-modules.
As in \S\ref{subsection;14.11.8.120},
we have the associated morphism
\[
 P(F_{\dagger}\nbigs,F_{\dagger}\kappa)_i:
 F_{\dagger}^i\nbigm_2
\otimes
 j^{\ast}F_{\dagger}^{-i}\nbigm_2
\lrarr
 \lambda^{-w+d_Y}
 \nbigo_{\nbigy}
\]
induced by the composite of the following morphisms:
\begin{equation}
 \label{eq;14.11.16.1}
\begin{CD}
 F_{\dagger}^i\nbigm_2
 @>{F_{\dagger}^i\nbigs}>>
 \lambda^{-w}F_{\dagger}^i\nbigm_1
 @>{F_{\dagger}^i\kappa}>>
 \lambda^{-w}F_{\dagger}^ij^{\ast}\DDD\nbigm_2
\simeq
 \lambda^{-w}
 j^{\ast}\DDD F_{\dagger}^{-i}\nbigm_2
 \simeq
 \lambda^{-w+d_Y}
 j^{\ast}\bigl(
 F_{\dagger}^{-i}\nbigm_2
 \bigr)^{\lor}
\end{CD}
\end{equation}
The last isomorphism is given by
$(-1)^{d_Y}\nu_{F_{\dagger}^{-i}\nbigm_2}$.
We give an expression of
$P(F_{\dagger}\nbigs,F_{\dagger}\kappa)_i$
in terms of $P$.

\vspace{.1in}

Let $Q\in Y$ be any point of $Y$.
We set $X_Q:=F^{-1}(Q)$.
Let $F_Q:X_Q\lrarr \{Q\}$ be the restriction of $F$.
Let $\nbigt_Q:=(\nbigm_{1Q},\nbigm_{2Q},C_Q)$
be the restriction of $\nbigt$ to $X_Q$
as smooth $\nbigr$-triple.
Let $F_{\dagger}^i(\nbigt)_Q$ be the restriction of
$F_{\dagger}^i(\nbigt)$ to $Q$
as smooth $\nbigr$-triple.
We naturally have
$F_{Q\dagger}^i(\nbigt_Q)
=F_{\dagger}^i(\nbigt)_Q$.

We set 
$\Omegabar_{X_Q}^1:=
 \lambda^{-1}p_{\lambda}^{\ast}\Omega_{X_Q}^1$
and 
$\Omegabar_{X_Q}^p:=\bigwedge^p\Omegabar_{X_Q}^1$.
We have the de Rham complex 
$\Omegabar_{X_Q}^{\bullet}\otimes
 (\nbigm_{2Q})$ for $\nbigm_{2Q}$.
We have
\[
 F_{Q\dagger}^i(\nbigm_{2Q})=
 R^{i+d}F_{\ast}
 \Bigl(
 \Omegabar_{X_Q}^{\bullet}\otimes(\nbigm_{2Q}).
 \Bigr)
\]
Set $d:=d_X-d_Y=\dim X_Q$.
Let $\Omega^{0,q}_{X_Q}$ denote the sheaf of
$C^{\infty}$ $(0,q)$-forms on $X_Q$.
Let $\Omega_{\cnum_{\lambda}\times X_Q/\cnum_{\lambda}}^{p,q}:=
 p_{\lambda}^{\ast}(\Omega_{X_Q}^p)
 \otimes_{p_{\lambda}^{-1}\nbigo_{X_Q}}
 p_{\lambda}^{-1}\bigl(\Omega_{X_Q}^{0,q}\bigr)$.
The pairing $P$ induces the following:
\[
\Bigl(
 \lambda^{-k}\Omega^{k,p}_{\cnum\times X_Q/\cnum}
 \otimes\nbigm_2
\Bigr)
\times
j^{\ast}
\Bigl(
 \lambda^{-(d-k)}\Omega^{d-k,d-p}_{\cnum\times X_Q/\cnum}
 \otimes\nbigm_2
\Bigr)
\lrarr
 \lambda^{-d-w+d_X}
 \Omega^{d,d}_{\cnum\times X_Q/\cnum}
\]
given by
$(\xi^{k,p}m_1,\xi^{d-k,d-p}m_2)
\longmapsto
 \xi^{k,p}\xi^{d-k,d-p}
 P(m_1,m_2)$.
Note that $-d-w+d_X=-w+d_Y$.
It induces
\[
 P_{i|Q}:
 R^{d+i}F_{\ast}\bigl(
 \Omegabar^{\bullet}_{X_Q}(\nbigm_{2Q})
 \bigr)
\times
 j^{\ast}R^{d-i}F_{\ast}\bigl(
 \Omega^{\bullet}_{X_Q}(\nbigm_2)
 \bigr)
\lrarr
 \lambda^{-w+d_Y}
 \nbigo_{\cnum_{\lambda}}
\]
By varying $Q$,
we obtain 
$P_{F,i}:
 F_{\dagger}^i\nbigm_2
 \otimes
 F_{\dagger}^{-i}\nbigm_2
\lrarr
 \lambda^{-w+d_Y}
 \nbigo_{\nbigy}$.

\begin{prop}
\label{prop;14.11.18.41}
We have
$P(F_{\dagger}\nbigs,F_{\dagger}\kappa)_i
=\epsilon(d)(-1)^{di}P_{F,i}$.
\end{prop}
\pf
It is enough to compare the specializations
along $\{\lambda\}\times X$ for $\lambda\neq 0$,
which follows from Proposition \ref{prop;14.11.15.1}.
\hfill\qed

\section{Some functoriality of $\nbigrtilde$-modules}
\label{section;14.12.26.12}
In an earlier version,
in the proof of Theorem \ref{thm;15.5.2.4},
we used the comparison 
of the duality functors for
$\nbigr$-modules and $\nbigrtilde$-modules.
After the proof of the theorem has been simplified,
we do not need it in the current version.
(See Remark \ref{rem;15.11.22.2}.)
Hence, the most results in this section
are not used in the other part of this paper.
However,
because the author hopes that 
the results in this section would be useful
for our understanding of the duality of
integrable mixed twistor $\nbigd$-modules,
we keep this section.

\subsection{Preliminary}

\subsubsection{Some sheaves}

For any complex manifold $X$,
let $\nbigx:=\cnum_{\lambda}\times X$.
Let $p_{\lambda}:\nbigx\lrarr X$ be the projection.
Set $d_X:=\dim X$.
Let $\nbigo_{\cnum_{\lambda}}$
denote the sheaf of functions on $\nbigx$
which are locally constant in the $X$-direction,
and holomorphic in the $\cnum_{\lambda}$-direction.

Let $\Omega_X^1$ denote the sheaf of holomorphic
$1$-forms on $X$.
Let $\Omegabar^1_X:=
 \lambda^{-1}p_{\lambda}^{\ast}\Omega^{1}_X$
and $\Omegabar^p_X:=
 \bigwedge^p\Omegabar_X^1$.
In particular, we set
$\Omegabar_X:=\Omegabar_X^{d_X}$.
As in the case of $\nbigd$-modules,
for any left $\nbigr_X$-module $\nbign^{\ell}$,
we naturally obtain a right $\nbigr_X$-module
$\Omegabar_X\otimes\nbign^{\ell}$.
Conversely,
for any right $\nbigr_X$-module $\nbign^r$,
we naturally obtain a left $\nbigr_X$-module
$\nbign^r\otimes\Omegabar_X^{-1}$.

Let $\Theta_X$ be the sheaf of holomorphic vector fields
on $X$ .
We set 
$\Thetabar_X:=\lambda p_{\lambda}^{\ast}\Theta_X$.
We set 
$\Thetabar^{p}_X:=\bigwedge^{-p}\Thetabar_X$.
We obtain the Spencer resolution
$\nbigr_X\otimes\Thetabar^{\bullet}_X$
of $\nbigo_{\nbigx}$
by locally free left $\nbigr_X$-modules.
(See \cite{sabbah2}.)

Let $\nbigx^0:=\{0\}\times X$.
Let $\Omega^1_{\nbigx}(\log \nbigx^0)$
be the sheaf of logarithmic holomorphic
$1$-forms on $(\nbigx,\nbigx^0)$.
We set 
$\Omegatilde^1_{\nbigx}:=
 \Omega^1_{\nbigx}(\log\nbigx^0)
 \otimes\nbigo(\nbigx^0)$,
and 
$\Omegatilde^p_{\nbigx}:=
 \bigwedge^p\Omegatilde^1_{\nbigx}$.
We set 
$\Omegatilde_{\nbigx}:=
 \Omegatilde_{\nbigx}^{\dim X+1}$.
For any left $\nbigrtilde_X$-module $\nbign^{\ell}$,
we naturally obtain a right $\nbigrtilde_X$-module
$\Omegatilde_{\nbigx}\otimes\nbign^{\ell}$.
Conversely, 
for any right $\nbigrtilde_X$-module $\nbign^{r}$,
we naturally obtain a left $\nbigrtilde_X$-module
$\nbign^{r}\otimes\Omegatilde_{\nbigx}^{-1}$.
We have the natural isomorphism
$\lambda^{-2}\Omegabar_X
\simeq \Omegatilde_{\nbigx}$
given by $\tau\longmapsto \tau d\lambda$.

Let $\Theta_{\nbigx}(\log \nbigx^0)$ be 
the sheaf of tangent vectors 
on $\nbigx$ which are logarithmic along $\nbigx^0$.
Set $\Thetatilde_{\nbigx}:=\lambda \Theta_{\nbigx}(\log\nbigx^0)$
and $\Thetatilde^p_{\nbigx}:=\bigwedge^{-p}\Thetatilde_{\nbigx}$.
We obtain the Spencer resolution
$\nbigrtilde_X\otimes\Thetatilde_{\nbigx}^{\bullet}$
of $\nbigo_{\nbigx}$
by locally free left $\nbigr_{\nbigx}$-modules
in the standard way.

\subsubsection{Partial Spencer resolution}

Let $\nbign$ be any $\nbigr_X$-module.
We set
$\nbign\langle
 \lambda^2\del_{\lambda}\rangle:=
 \bigoplus_{j=0}^{\infty}
 \del_{\lambda}^j\otimes
 \lambda^{2j}\nbign$.
We have the isomorphism of sheaves
$\Phi:
 \nbigrtilde_X\otimes_{\nbigr_X}\nbign
\lrarr
\bigl(
\nbign\langle\lambda^2\del_{\lambda}\rangle,
\triangledown
\bigr)$
given by
$\sum (\del_{\lambda}^j\lambda^{2j})\otimes m_j
\longmapsto
 \sum \del_{\lambda}^j\otimes (\lambda^{2j}m_j)$.
The natural $\nbigrtilde_X$-action on
$\nbigrtilde_X\otimes_{\nbigr_X}\nbign$
is transferred onto
$\nbign\langle\lambda^2\del_{\lambda}\rangle$.
We denote it by $\triangledown$
which is explicitly described as follows.
\begin{itemize}
\item
Let $P$ and $m$ denote local sections of
$\nbigr_X$ and $\nbign$, respectively.
For any $j\geq 0$,
we have local sections $P_{ij}$ $(i=0,\ldots,j)$
of $\nbigr_X$
such that 
$P\del_{\lambda}^j\lambda^{2j}
=\sum_{i=0}^j \del_{\lambda}^iP_{ij}\lambda^{2i}$
in $\nbigrtilde_X$.
Then, 
$P\triangledown(\del_{\lambda}^j\otimes\lambda^{2j}m)=
 \sum_{i=0}^j\del_{\lambda}^i\otimes(P_{ij}\lambda^{2i}m)$.
\item
We also have
$(\del_{\lambda}\lambda^2)\triangledown
 (\del_{\lambda}^j\otimes \lambda^{2j}m)
=\del_{\lambda}^{j+1}\otimes(\lambda^{2j+2}m)
-2j\del_{\lambda}^j\otimes \lambda^{2j+1}m
+j(j-1)\del_{\lambda}^{j-1}\otimes \lambda^{2j}m$.
\end{itemize}
We have the inclusion
$\iota:
 \nbign\lrarr 
 \nbign\langle\lambda^2\del_{\lambda}\rangle$
given by 
$\iota(m)=1\otimes m$.
It is an $\nbigr_X$-homomorphism.
We have the decomposition:
\begin{equation}
\label{eq;14.12.20.10}
 \nbign\langle\lambda^2\del_{\lambda}\rangle
=\bigoplus_{j=0}^{\infty}
 (\del_{\lambda}\lambda^2)^j\triangledown\iota(\nbign)
\end{equation}

We also have the $\nbigr_X$-action 
$\blacktriangledown$ on 
$\nbign\langle\lambda^2\del_{\lambda}\rangle$
given as follows.
\begin{itemize}
\item
For any $P$ and $m$ as above,
we set
$P\blacktriangledown(\del_{\lambda}^j\otimes\lambda^{2j}m)
=\del_{\lambda}^j\otimes P(\lambda^{2j}m)$.
\end{itemize}
If moreover $\nbign$ is an $\nbigrtilde_X$-module,
the $\nbigr_X$-action $\blacktriangledown$
is extended to the $\nbigrtilde_X$-action
$\blacktriangledown$
on $\nbign\langle\lambda^2\del_{\lambda}\rangle$
given as follows:
\[
 \del_{\lambda}\lambda^2\blacktriangledown
 (\del_{\lambda}^j\otimes \lambda^{2j}m)
=\del_{\lambda}^j\otimes
 \del_{\lambda}\lambda^{2+2j}m
-\del_{\lambda}^{j+1}
 \otimes\lambda^{2j+2}m.
\]
It is easy to check that
$(\nbign\langle\lambda^2\del_{\lambda}\rangle,\blacktriangledown)$
is an $\nbigrtilde_X$-module.
(See also Lemma \ref{lem;14.12.21.130} below.)
We have the decomposition:
\begin{equation}
\label{eq;14.12.20.21}
 \nbign\langle\lambda^2\del_{\lambda}\rangle
=\bigoplus_{j=0}^{\infty}
 (\del_{\lambda}\lambda^2)^j\blacktriangledown
 \iota(\nbign).
\end{equation}

For any $\ell\in\seisuu$,
we set 
$\lambda^{\ell}\nbign
 \langle\lambda^2\del_{\lambda}\rangle
:=\bigoplus_{j=0}^{\infty}
\del_{\lambda}^{j}\otimes \lambda^{2j+\ell}\nbign$.
It is equipped with
the natural $\nbigrtilde$-actions $\triangledown$
and $\blacktriangledown$,
obtained as above.
We naturally have
$\lambda^2\nbign\langle\lambda^2\del_{\lambda}\rangle
=\lambda^2\triangledown
 \nbign\langle\lambda^2\del_{\lambda}\rangle
=\lambda^2\blacktriangledown
 \nbign\langle\lambda^2\del_{\lambda}\rangle$.
We have the well defined morphisms of sheaves
$\del_{\lambda}\triangledown:
 \lambda^2\nbign\langle\lambda^2\del_{\lambda}\rangle
\lrarr
 \nbign\langle\lambda^2\del_{\lambda}\rangle$
and 
$\del_{\lambda}\blacktriangledown:
 \lambda^2\nbign\langle\lambda^2\del_{\lambda}\rangle
\lrarr
 \nbign\langle\lambda^2\del_{\lambda}\rangle$.

It is easy to check the following lemma.
\begin{lem}
\mbox{{}}\label{lem;14.12.21.130}
Let $\nbign$ be an $\nbigrtilde_X$-module.
\begin{itemize}
\item
For any local section $m$ of $\nbign$,
we have the following relation
\begin{equation}
\label{eq;14.12.20.11}
 (\del_{\lambda}\lambda^2)\triangledown
 \iota(m)
+
 (\del_{\lambda}\lambda^2)\blacktriangledown
 \iota(m)
=\iota(\del_{\lambda}\lambda^2m).
\end{equation}
\item
For the above isomorphism
$\Phi:\nbigrtilde_X\otimes_{\nbigr_X}\nbign
\simeq
 \nbign\langle\lambda^2\del_{\lambda}\rangle$,
we have
\[
 \del_{\lambda}\lambda^2\blacktriangledown
 \Phi\bigl(
 \del_{\lambda}^j\lambda^{2j}
\otimes m
 \bigr)
=\Phi\bigl(
\del_{\lambda}^j\lambda^{2j+2}(-\del_{\lambda})
\otimes m
 \bigr)
+\Phi\bigl(
 \del_{\lambda}^j\lambda^{2j+2}
\otimes
 \del_{\lambda}m
 \bigr).
\]
For any $P\in\nbigr_X$,
we have the following:
\[
 P\blacktriangledown
 \Phi(\del_{\lambda}^j\lambda^{2j}\otimes m)
=\Phi(\del_{\lambda}^j\lambda^{2j}\otimes Pm).
\]
\item
$\del_{\lambda}\triangledown$
gives a morphism
$\bigl(
 \lambda^2\nbign\langle\lambda^2\del_{\lambda}\rangle,
 \blacktriangledown
 \bigr)
\lrarr
 \bigl(
 \nbign\langle\lambda^2\del_{\lambda}\rangle,
 \blacktriangledown
 \bigr)$,
and 
$\del_{\lambda}\blacktriangledown$
gives a morphism
$\bigl(
 \lambda^2\nbign\langle\lambda^2\del_{\lambda}\rangle,
 \triangledown
 \bigr)
\lrarr
 \bigl(
 \nbign\langle\lambda^2\del_{\lambda}\rangle,
 \triangledown
 \bigr)$.
\hfill\qed
\end{itemize}
\end{lem}

By the lemma, for any $\nbigrtilde_X$-module $\nbign$,
we obtain the following exact sequence:
\[
\begin{CD}
 0 @>>>
 (\lambda^2\nbign\langle\lambda^2\del_{\lambda}\rangle,
 \blacktriangledown)
 @>{\del_{\lambda}\triangledown}>>
 (\nbign\langle\lambda^2\del_{\lambda}\rangle,
 \blacktriangledown)
 @>{a_1}>>
 \nbign
 @>>> 0
\end{CD}
\]
Here, 
$a_1$ is the projection 
onto the $0$-th component with respect to
the decomposition (\ref{eq;14.12.20.10}).
We also have the following exact sequence:
\[
\begin{CD}
 0 @>>>
 (\lambda^2\nbign\langle\lambda^2\del_{\lambda}\rangle,
 \triangledown)
 @>{\del_{\lambda}\blacktriangledown}>>
 (\nbign\langle\lambda^2\del_{\lambda}\rangle,
 \triangledown)
 @>{a_2}>>
 \nbign
 @>>> 0
\end{CD}
\]
Here, $a_2$ is the projection
onto the $0$-th component with respect to
the decomposition (\ref{eq;14.12.20.21}).
We obtain 
the following $\nbigrtilde_X$-resolutions of $\nbign$:
\[
S_{\lambda}^{\blacktriangledown}(\nbign):=
\Bigl(
 \bigl(
 \lambda^2\nbign\langle\lambda^2\del_{\lambda}\rangle,
 \blacktriangledown
 \bigr)
\stackrel{\del_{\lambda}\triangledown}{\lrarr}
  \bigl(
 \nbign\langle\lambda^2\del_{\lambda}\rangle,
 \blacktriangledown
 \bigr)
\Bigr)
\]
\[
S_{\lambda}^{\triangledown}(\nbign):=
\Bigl(
 \bigl(
 \lambda^2\nbign\langle\lambda^2\del_{\lambda}\rangle,
 \triangledown
 \bigr)
\stackrel{\del_{\lambda}\blacktriangledown}{\lrarr}
  \bigl(
 \nbign\langle\lambda^2\del_{\lambda}\rangle,
 \triangledown
 \bigr)
\Bigr)
\]

\begin{example}
\label{example;14.12.21.3}
Let us consider the case $\nbign=\nbigr_X$.
We consider the $\nbigr_X$-action given by
the left multiplication.
We have the isomorphism
$\nbigr_X\langle\lambda^2\del_{\lambda}\rangle
\simeq
 \nbigrtilde_X$
given by
$\sum \del_{\lambda}^j\otimes \lambda^{2j}m_j
\longmapsto
 \sum \del_{\lambda}^j\lambda^{2j}m_j$.
Under the isomorphism,
$(\del_{\lambda}\lambda^2)\triangledown$
is the left multiplication of
$\del_{\lambda}\lambda^2$,
and 
$(\del_{\lambda}\lambda^2)\blacktriangledown$
is the right multiplication of
$-\lambda^2\del_{\lambda}$.
Hence, we naturally have 
$\nbigrtilde_X\otimes\Thetatilde_{\nbigx}^{\bullet}
\simeq
 S_{\lambda}^{\triangledown}
 \bigl(
 \nbigr_X\otimes\Thetabar^{\bullet}_X\bigr)$.
\hfill\qed
\end{example}

\subsubsection{Exchange}
\label{subsection;14.12.22.1}

We define
$\Psi:
 \nbign\langle\lambda^2\del_{\lambda}\rangle
\lrarr
 \nbign\langle\lambda^2\del_{\lambda}\rangle$
by 
$\Psi\bigl(
 (\del_{\lambda}\lambda^2)^j\triangledown\iota(m)
 \bigr)
:=
  (\del_{\lambda}\lambda^2)^j\blacktriangledown\iota(m)$.
By the construction,
$\Psi$ gives an isomorphism
$\bigl(
 \nbign\langle\lambda^2\del_{\lambda}\rangle,
 \triangledown
 \bigr)
\lrarr
 \bigl(
 \nbign\langle\lambda^2\del_{\lambda}\rangle,
 \blacktriangledown
 \bigr)$.

\begin{lem}
$\Psi$ also gives an isomorphism
$\Psi:
 \bigl(
 \nbign\langle\lambda^2\del_{\lambda}\rangle,
 \blacktriangledown\bigr)
\lrarr
 \bigl(
 \nbign\langle\lambda^2\del_{\lambda}\rangle,
 \triangledown\bigr)$.
\end{lem}
\pf
It is enough to check that
$\Psi
 \bigl(
 (\del_{\lambda}\lambda^2)\blacktriangledown
 (\del_{\lambda}\lambda^2)^j\triangledown
 \iota(m)
 \bigr)
= (\del_{\lambda}\lambda^2)\triangledown
 \Psi\bigl(
 (\del_{\lambda}\lambda^2)^j\triangledown
 \iota(m)
 \bigr)$
due to the decomposition (\ref{eq;14.12.20.10}).
By the commutativity of
$(\del_{\lambda}\lambda^2)\blacktriangledown$
and 
$(\del_{\lambda}\lambda^2)\triangledown$,
we have only to check
$\Psi\bigl(
 (\del_{\lambda}\lambda^2)\blacktriangledown
\iota(m)
 \bigr)
=(\del_{\lambda}\lambda^2)\triangledown
 \iota(m)$.
It follows from (\ref{eq;14.12.20.11}).
\hfill\qed

\vspace{.1in}
We have the induced isomorphisms
$\Psi:
 S_{\lambda}^{\triangledown}(\nbign)
\lrarr
 S_{\lambda}^{\blacktriangledown}(\nbign)$
and 
$\Psi:
 S_{\lambda}^{\blacktriangledown}(\nbign)
\lrarr
 S_{\lambda}^{\triangledown}(\nbign)$.

\begin{example}
If $\nbign=\nbigr_X$,
under the isomorphism
$\nbigrtilde_X\simeq\nbigr_X\langle\lambda^2\del_{\lambda}\rangle$,
the isomorphism
$\Psi:\nbigrtilde_X\lrarr\nbigrtilde_X$
is given by
$\Psi\bigl(\sum \del_{\lambda}^jm_j\bigr)
=\sum m_j(-\del_{\lambda})^j$.
\hfill\qed
\end{example}

\subsubsection{Bi-modules}
\label{subsection;14.12.24.2}

Let $\nbigo_{\cnum_{\lambda}}
 \langle\lambda^2\del_{\lambda}\rangle$
denote the sheaf of subalgebras in $\nbigrtilde_X$
generated by $\lambda^2\del_{\lambda}$
over $\nbigo_{\cnum_{\lambda}}$.
Let $\nbigb$ be an $\nbigo_{\cnum_{\lambda}}$-module.
An action of $\lambda^2\del_{\lambda}$
on $\nbigb$ means
a morphism of sheaves 
$\rho(\lambda^2\del_{\lambda}):\nbigb\lrarr \nbigb$
such that 
$\rho(\lambda^2\del_{\lambda})(fP)
=f\rho(\lambda^2\del_{\lambda})(P)
+(\lambda^2\del_{\lambda}f)\,P$
for any local sections $f\in\nbigo_{\cnum_{\lambda}}$
and $P\in\nbigb$.
If $\nbigb$ is a sheaf of algebras over 
$\nbigo_{\cnum_{\lambda}}$,
we moreover impose the Leibniz rule
$\rho(\lambda^2\del_{\lambda})(PQ)
=\rho(\lambda^2\del_{\lambda})(P)Q
+P\rho(\lambda^2\del_{\lambda})(Q)$.
If we are given a sheaf of algebras $\nbigb$
over $\nbigo_{\cnum}$
equipped with an action of $\lambda^2\del_{\lambda}$,
the sheaf
$\nbigb\otimes_{\nbigo_{\cnum_{\lambda}}}
 \nbigo\langle\lambda^2\del_{\lambda}\rangle$
is naturally a sheaf of algebras.
The multiplication is given by
$(P\otimes (\lambda^2\del_{\lambda})^k)
\cdot
 (Q\otimes (\lambda^2\del_{\lambda})^{\ell})
=\sum_{j=0}^k (k:j)
 P\cdot
 \rho(\lambda^2\del_{\lambda})^j(Q)\otimes
 (\lambda^2\del_{\lambda})^{k+\ell-j}$,
where $(k:j)$ denote the binomial coefficients.

We can naturally regard
$\nbigr_X$ as a left 
$\nbigo_{\cnum_{\lambda}}
 \langle\lambda^2\del_{\lambda}\rangle$.
The action of $\lambda^2\del_{\lambda}$
on $\nbigr_X$ is given by
$\rho(\lambda^2\del_{\lambda})(P)
=[\lambda^2\del_{\lambda},P]$
in $\nbigrtilde_X$
which satisfies the Leibniz rule
$\rho(\lambda^2\del_{\lambda})(PQ)
=\rho(\lambda^2\del_{\lambda})(P)\cdot Q
+P\cdot \rho(\lambda^2\del_{\lambda})(Q)$.
The sheaf
$\nbigr_X\otimes_{\nbigo_{\cnum_{\lambda}}}
 \nbigo_{\cnum_{\lambda}}\langle\lambda^2\del_{\lambda}\rangle$
is naturally a sheaf of algebras.
We have the isomorphism
$\nbigr_X\otimes_{\nbigo_{\cnum_{\lambda}}}
 \nbigo_{\cnum_{\lambda}}\langle\lambda^2\del_{\lambda}\rangle
\simeq
 \nbigrtilde_X$
given by
$\sum P_j\otimes (\lambda^2\del_{\lambda})^j
\longmapsto
 \sum P_j(\lambda^2\del_{\lambda})^j$.

Because $\nbigo_{\cnum_{\lambda}}$ is the center of
$\nbigr_X$,
we have the naturally defined sheaf of algebras
$\nbigr_X\otimes_{\nbigo_{\cnum_{\lambda}}}
 \nbigr_X$ on $\nbigx$.
It is naturally an 
$\nbigo_{\cnum_{\lambda}}
\langle\lambda^2\del_{\lambda}\rangle$-module.
The action of $\lambda^2\del_{\lambda}$
is given by
$\rho_1(\lambda^2\del_{\lambda})
(P\otimes Q)
=\rho(\lambda^2\del_{\lambda})(P)\otimes Q
+P\otimes\rho(\lambda^2\del_{\lambda})(Q)$.
It satisfies the Leibniz rule.
We set $\nbiga_X:=
 \nbigr_X\otimes_{\nbigo_{\cnum_{\lambda}}}
 \nbigr_X\otimes_{\nbigo_{\cnum_{\lambda}}}
 \nbigo_{\cnum_{\lambda}}\langle\lambda^2\del_{\lambda}\rangle$.
It is naturally a sheaf of algebras.
Let $k_i:\nbigr_X\lrarr \nbiga_X$ $(i=1,2)$
be the morphism of algebras 
given by the inclusion to the $i$-th component.
They are extended to the morphism of algebras
$k_i:\nbigrtilde_X\lrarr \nbiga_X$.

Let $\nbign$ be an $\nbiga_X$-module.
Let $\ell$ (resp. $r$)
denote the $\nbigr_X$-action or $\nbigrtilde_X$-action
on $\nbign$ induced by $k_1$ (resp. $k_2$).
We set
$\nbign\langle\lambda^2\del_{\lambda}\rangle:=
 \bigoplus_{j=0}^{\infty}
 \del_{\lambda}^j\otimes\lambda^{2j}\nbign$.
The $\nbigrtilde_X$-action $\ell$ on $\nbign$
and the construction $\blacktriangledown$
give an $\nbigrtilde_X$-action $\ell\blacktriangledown$
on $\nbign\langle\lambda^2\del_{\lambda}\rangle$.
The induced multiplication of
$g\in\nbigrtilde_X$ and 
$c\in\nbign\langle\lambda^2\del_{\lambda}\rangle$
is denoted by
$g(\ell\blacktriangledown)c$.
We obtain an $\nbigr_X$-action
by taking the restriction,
which is also denoted by $\ell\blacktriangledown$.
For $a\in\nbigrtilde_X$,
let $\ell\blacktriangledown(a)$ denote the endomorphism of
the sheaves $\nbign\langle\lambda^2\del_{\lambda}\rangle$
given by $m\longmapsto a(\ell\blacktriangledown)m$.
We use the notation
$\ell\triangledown$,
$r\triangledown$
and $r\blacktriangledown$
similarly.

\begin{lem}
\label{lem;14.12.21.201}
We have
$r\triangledown(a_1)
 \circ
 \ell\blacktriangledown(a_2)
=\ell\blacktriangledown (a_2)
 \circ
 r\triangledown (a_1)$
for any $a_1,a_2\in\nbigrtilde_X$.
\end{lem}
\pf
By using the explicit construction of
$\triangledown$ and $\blacktriangledown$,
we can observe that
$r\triangledown (a_1)\circ 
 \ell\blacktriangledown(a_2)
=\ell\blacktriangledown (a_2)
 \circ
 r\triangledown (a_1)$
if $a_1,a_2\in\nbigr_X$.
By the assumption,
we have $r(\lambda^2\del_{\lambda})=\ell(\lambda^2\del_{\lambda})$.
Hence,
we have 
$r\triangledown(\del_{\lambda}\lambda^2)
=\ell\triangledown(\del_{\lambda}\lambda^2)$
and 
$r\blacktriangledown(\del_{\lambda}\lambda^2)
=\ell\blacktriangledown(\del_{\lambda}\lambda^2)$.
By using Lemma \ref{lem;14.12.21.130},
we obtain 
$r\triangledown(a_1)
 \circ 
 \ell\blacktriangledown(a_2)
=\ell\blacktriangledown(a_2)
 \circ
 r\triangledown(a_1)$
if either one of $a_1$ or $a_2$ is 
$\del_{\lambda}\lambda^2$.
\hfill\qed

\vspace{.1in}
The following lemma is clear by construction.

\begin{lem}
\label{lem;14.12.21.200}
We have
$r\triangledown(a_1)
\circ
 \ell\triangledown(a_2)
=\ell\triangledown(a_2)
 \circ
 r\triangledown(a_1)$
and 
$r\blacktriangledown(a_1)
 \circ
 \ell\blacktriangledown(a_2)
=\ell\blacktriangledown(a_2)
 \circ
 r\blacktriangledown(a_1)$
for any $a_1,a_2\in\nbigr_X$.
\hfill\qed
\end{lem}

As mentioned in the proof of Lemma \ref{lem;14.12.21.201},
we have
$r\triangledown(\del_{\lambda}\lambda^2)
=\ell\triangledown(\del_{\lambda}\lambda^2)$
and 
$r\blacktriangledown(\del_{\lambda}\lambda^2)
=\ell\blacktriangledown(\del_{\lambda}\lambda^2)$,
the morphism $\Psi$ in \S\ref{subsection;14.12.22.1}
is independent of the choices of $\ell$ and $r$.
It exchanges
$(r\triangledown,\ell\blacktriangledown)$
and 
$(r\blacktriangledown,\ell\triangledown)$.

\begin{example}
\label{example;14.12.22.20}
We consider 
$\nbign:=\nbigr_X\otimes\Omegabar_X^{-1}$
as an $\nbiga_X$-module.
Here, the $\nbigr_X$-actions $\ell$ and $r$
are induced by the left and the right multiplications.
The action of $\lambda^2\del_{\lambda}$
is given by
$\lambda^2\del_{\lambda} (P)
=[\lambda^2\del_{\lambda},P]$
in $\nbigrtilde_X$.

We have the isomorphism
$\bigl(
 \nbign\langle\lambda^2\del_{\lambda}\rangle,
 \ell\triangledown
\bigr)
\simeq
\nbigrtilde_X\otimes_{\nbigr_X}
 \bigl(
 \nbign,\ell
 \bigr)
=\nbigrtilde_X\otimes\Omegabar_X^{-1}
\simeq
 \nbigrtilde_X\otimes
 \Omegabar_X^{-1}(d\lambda)^{-1}$.
Under the isomorphism,
$\ell\triangledown$ is equal to the left multiplication,
and $r\blacktriangledown$
is equal to the $\nbigrtilde_X$-action
induced by the right multiplication.

We have the isomorphism
$\bigl(
 \nbign\langle\lambda^2\del_{\lambda}\rangle,
 r\triangledown
\bigr)
\simeq
\nbigrtilde_X\otimes_{\nbigr_X}
 \bigl(
 \nbign,r
 \bigr)
=\nbigrtilde_X\otimes\Omegabar_X^{-1}
\simeq
 \nbigrtilde_X\otimes
 \Omegabar_X^{-1}(d\lambda)^{-1}$.
Under the isomorphism,
$r\triangledown$ is equal to
the $\nbigrtilde_X$-action induced by the right multiplication,
and $\ell\blacktriangledown$ is equal to
the left multiplication.
\hfill\qed
\end{example}

\subsubsection{Functoriality for inner homomorphisms}

Let $\nbign_1$ be an $\nbigrtilde_X$-module. 
Let $\nbign_2$ be an $\nbiga_X$-module.
Let $\nhom_{\nbigr_X}(\nbign_1,\nbign_2^{\ell})$
be the sheaf of $\nbigr_X$-homomorphisms from
$\nbign_1$ to $(\nbign_2,\ell)$.
It is equipped with the $\nbigr_X$-action induced by $r$.
It is naturally extended to an $\nbigrtilde_X$-action,
where the action of $\lambda^2\del_{\lambda}$
is given by 
$\lambda^2\del_{\lambda}(f)(m)=
 \lambda^2\del_{\lambda}(f(m))
-f\bigl(\lambda^2\del_{\lambda}m\bigr)$.
We have the $\nbigrtilde_X$-actions
$\triangledown$ and $\blacktriangledown$
on $\nhom_{\nbigr_X}(\nbign_1,\nbign_2^{\ell})
 \langle\lambda^2\del_{\lambda}\rangle$.

\vspace{.1in}

Let $\nhom_{\nbigr_X}\bigl(\nbign_1,
 \nbign_2\langle\lambda^2\del_{\lambda}\rangle
 ^{\ell\blacktriangledown}
 \bigr)$
denote the sheaf of $\nbigr_X$-homomorphisms
$\nbign_1\lrarr
 \bigl(
 \nbign_2\langle\lambda^2\del_{\lambda}\rangle,
 \ell\blacktriangledown
\bigr)$.
By Lemma \ref{lem;14.12.21.201},
we have the $\nbigrtilde_X$-action
on 
$\nhom_{\nbigr_X}\bigl(\nbign_1,
 \nbign_2\langle\lambda^2\del_{\lambda}\rangle
 ^{\ell\blacktriangledown}
 \bigr)$
induced by $r\triangledown$ 
on 
$\nbign_2\langle\lambda^2\del_{\lambda}\rangle$.
The induced $\nbigrtilde_X$-action
is also denoted by $r\triangledown$.
By Lemma \ref{lem;14.12.21.200},
we have the $\nbigr_X$-action
on 
$\nhom_{\nbigr_X}\bigl(\nbign_1,
 \nbign_2\langle\lambda^2\del_{\lambda}\rangle
 ^{\ell\blacktriangledown}
 \bigr)$
induced by $r\blacktriangledown$
on 
$\nbign_2\langle\lambda^2\del_{\lambda}\rangle$.
It is extended to an $\nbigrtilde_X$-action.
The action of $\lambda^2\del_{\lambda}$
is given by
$(\lambda^2\del_{\lambda}g)(m)
=\lambda^2\del_{\lambda}(\ell\blacktriangledown)(g(m))
-g(\lambda^2\del_{\lambda}m)$.
The induced $\nbigrtilde_X$-action is denoted by
$r\blacktriangledown$.

\vspace{.1in}

We define a morphism of sheaves
$F:
 \nhom_{\nbigr_X}\bigl(
 \nbign_1,\nbign_2^{\ell}
 \bigr)\langle
 \lambda^2\del_{\lambda}\rangle
\lrarr
 \nhom_{\nbigr_X}\bigl(
 \nbign_1,\nbign_2^{\ell\blacktriangledown}
 \langle\lambda^2\del_{\lambda}\rangle
 \bigr)$
given as follows.
For any local section
$\sum \del_{\lambda}^j\otimes\lambda^{2j}g_j$
of  
$\nhom_{\nbigr_X}\bigl(
 \nbign_1,\nbign_2^{\ell}
 \bigr)\langle
 \lambda^2\del_{\lambda}\rangle$
and any local section
$m$ of $\nbign_1$,
we set
\[
 F\Bigl(
 \sum\del_{\lambda}^j\otimes \lambda^{2j}g_j
 \Bigr)(m)
:=\sum \del_{\lambda}^j\otimes \lambda^{2j}g_j(m).
\]
If $\nbign_1$ is $\nbigr_X$-coherent,
then $F$ is an isomorphism.
We can check the following lemma
by a direct computation.
\begin{lem}
\label{lem;14.12.21.51}
For local sections $a\in \nbigrtilde_X$
and 
$b\in
 \nhom_{\nbigr_X}(\nbign_1,\nbign_2^{\ell})
 \langle\lambda^2\del_{\lambda}\rangle$,
we have 
$F(a \triangledown b)
=a(r\triangledown) F(b)$
and 
$F(a \blacktriangledown b)
=a(r\blacktriangledown) F(b)$.
\hfill\qed
\end{lem}

We have the following natural isomorphism
by the decomposition (\ref{eq;14.12.20.10}):
\begin{equation}
\label{eq;14.12.21.40}
 \nhom_{\nbigr_X}\bigl(\nbign_1,
 \nbign_2
  \langle\lambda^2\del_{\lambda}
 \rangle^{\ell\blacktriangledown}
 \bigr)
\simeq
 \nhom_{\nbigrtilde_X}\bigl(
 \nbign_1\langle\lambda^2\del_{\lambda}
 \rangle^{\blacktriangledown},
 \nbign_2
  \langle\lambda^2\del_{\lambda}
 \rangle^{\ell\blacktriangledown}
 \bigr)
\end{equation}

\begin{lem}
\label{lem;14.12.21.50}
We have the following commutative diagram:
\[
 \begin{CD}
 \nhom_{\nbigr_X}\bigl(\nbign_1,
 \nbign_2
  \langle\lambda^2\del_{\lambda}\rangle^{\ell\blacktriangledown}
 \bigr)
 @>{b}>>
 \nhom_{\nbigrtilde_X}\bigl(
 \nbign_1\langle\lambda^2\del_{\lambda}\rangle^{\blacktriangledown},
 \nbign_2
  \langle\lambda^2\del_{\lambda}\rangle^{\ell\blacktriangledown}
 \bigr)
 \\
 @V{a_1}VV @V{a_2}VV \\
\nhom_{\nbigr_X}\bigl(
 \lambda^2\nbign_1,
\nbign_2
  \langle\lambda^2\del_{\lambda}\rangle^{\ell\blacktriangledown}
 \bigr)
 @>{b}>>
 \nhom_{\nbigrtilde_X}\bigl(
 \lambda^2\nbign_1
 \langle\lambda^2\del_{\lambda}\rangle^{\blacktriangledown},
 \nbign_2
  \langle\lambda^2\del_{\lambda}\rangle^{\ell\blacktriangledown}
 \bigr)
 \end{CD}
\]
The horizontal arrows $b$ are 
{\rm(\ref{eq;14.12.21.40})}.
The morphism $a_2$ is given by
$\del_{\lambda}\triangledown$ 
on $\nbign_1\langle\lambda^2\del_{\lambda}\rangle$.
The morphism $a_1$ is determined by
$a_1(g)(m)=
-\del_{\lambda}(\ell\blacktriangledown)\bigl(g(m)\bigr)
+g(\del_{\lambda}m)$
for any 
$g\in \nhom_{\nbigr_X}\bigl(\nbign_1,
 \nbign_2\langle\lambda^2\del_{\lambda}
 \rangle^{\ell\blacktriangledown}\bigr)$
and 
$m\in \lambda^2\nbign_1$.
\end{lem}
\pf
For any 
$g\in \nhom_{\nbigr_X}\bigl(\nbign_1,
 \nbign_2\langle\lambda^2\del_{\lambda}
 \rangle^{\ell\blacktriangledown}\bigr)$
and $m=\lambda^2m_1\in\lambda^2\nbign_1$,
we have
\begin{multline}
 (a_2\circ b)(g)\bigl(\iota(\lambda^2 m_1)\bigr)
=b(g)\Bigl(
 (\del_{\lambda}\lambda^2)
 \triangledown 
 \iota(m_1)
\bigr)
 \Bigr)
=b(g)\Bigl(
 -(\del_{\lambda}\lambda^2)
 \blacktriangledown
 \iota(m_1)
+\iota(\del_{\lambda}\lambda^2m_1)
 \Bigr)
 \\
=(-\del_{\lambda}\lambda^2)
 (\ell\blacktriangledown)
\bigl(g(m_1)\bigr)
+g\bigl(\del_{\lambda}\lambda^2m_1\bigr)
 =
-\del_{\lambda} (\ell\blacktriangledown)
 \bigl(g(m)\bigr)
+g(\del_{\lambda}m)
\end{multline}
Thus, we obtain the claim of the lemma.
\hfill\qed

\begin{cor}
\label{cor;14.12.21.101}
We have the following natural isomorphism
of $\nbigrtilde_X$-complexes:
\begin{equation}
\label{eq;14.12.22.10}
 S_{\lambda}^{\triangledown}
 \nhom_{\nbigr_X}\bigl(
 \nbign_1,\nbign_2^{\ell}
 \bigr)
\simeq
 \nhom_{\nbigrtilde_X}
\bigl(S_{\lambda}^{\blacktriangledown}\nbign_1,
 \lambda^2\nbign_2\langle\lambda^2\del_{\lambda}\rangle
 ^{\ell\blacktriangledown}
 \bigr)[1]
\end{equation}
\begin{equation}
\label{eq;14.12.22.11}
 S_{\lambda}^{\blacktriangledown}
 \nhom_{\nbigr_X}\bigl(
 \nbign_1,\nbign_2^{\ell}
 \bigr)
\simeq
 \nhom_{\nbigrtilde_X}
\bigl(S_{\lambda}^{\triangledown}\nbign_1,
 \lambda^2\nbign_2\langle\lambda^2\del_{\lambda}\rangle
 ^{\ell\triangledown}
 \bigr)[1]
\end{equation}
\end{cor}
\pf
We obtain (\ref{eq;14.12.22.10})
from Lemma \ref{lem;14.12.21.51}
and Lemma \ref{lem;14.12.21.50}.
We obtain (\ref{eq;14.12.22.11})
by exchanging $\triangledown$ and $\blacktriangledown$.
\hfill\qed

\subsection{Push-forward}
\label{subsection;15.1.4.21}

\subsubsection{Push-forward of $\nbigr$-modules
and $\nbigrtilde$-modules}

\label{subsection;14.12.21.1}

We recall the push-forward of $\nbigr$-modules
and $\nbigrtilde$-modules.
See \cite{sabbah2} for more details.
Let $X$ be any complex manifold.
Let $\nbigc^{\infty}_{\nbigx/\cnum_{\lambda}}$
denote the sheaf of $C^{\infty}$-functions $F$ on $\nbigx$
such that $\delbar_{\lambda}F=0$.
Let $\Omega^{0,q}_X$ denote the sheaf of
$C^{\infty}$ $(0,q)$-forms on $X$.
Let $\nbigc^{\infty}_X$ denote the sheaf of
$C^{\infty}$-functions on $X$.
We define
\[
 \nbigcbar_X^{\infty,j}:=
 \bigoplus_{p+q=j}
 p_{\lambda}^{-1}\Omega_{X}^{0,q}
 \otimes_{p_{\lambda}^{-1}\nbigc^{\infty}_X}
 \bigl(
 \Omegabar^{p}_X
\otimes_{\nbigo_{\nbigx}}
 \nbigc^{\infty}_{\nbigx/\cnum_{\lambda}}
 \bigr).
\]
With the exterior derivative in the $X$-direction,
we obtain the complexes
$\bigl(
 \Omegabar^{\bullet}_X,d
 \bigr)$
and 
$\bigl(
 \nbigcbar_X^{\infty,\bullet},d
 \bigr)$.
The natural inclusion gives a quasi-isomorphism
$\bigl(
 \Omegabar^{\bullet}_X,d
 \bigr)
\lrarr
\bigl(
 \nbigcbar_X^{\infty,\bullet},d
 \bigr)$.

Let $f:X\lrarr Y$ be any morphism of complex manifolds.
The induced morphism
$\nbigx\lrarr\nbigy$ is also denoted by $f$.
We set 
$\nbigr_{X,Y,f}:=
 \nbigr_X\otimes_{\nbigo_{\cnum_{\lambda}}}
 f^{-1}\nbigr_Y$.
It is naturally a sheaf of algebras.
It is also naturally 
an $\nbigo_{\cnum_{\lambda}}
 \langle\lambda^2\del_{\lambda}\rangle$-module,
and the action of $\lambda^{2}\del_{\lambda}$
satisfies the Leibniz rule.
Hence, we obtain the sheaf of algebras
$\nbiga_{X,Y,f}:=
 \nbigr_{X,Y,f}\otimes_{\nbigo_{\cnum_{\lambda}}}
 \nbigo_{\cnum_{\lambda}}
 \langle\lambda^2\del_{\lambda}\rangle$
as in the case of $\nbiga_X$.
We have a natural inclusion
$\nbigrtilde_X\lrarr \nbiga_{X,Y,f}$
and 
$f^{-1}\nbigrtilde_Y\lrarr \nbiga_{X,Y,f}$.
We set
$\nbigr_{Y\larr X}:=
 \Omegabar_{X}
 \otimes_{\nbigo_{\nbigx}}
 f^{-1}(\nbigr_Y\otimes
 \Omegabar_Y^{-1})$.
It is naturally 
a left $\nbigr_{X,Y,f}$-module.
We have the natural action of $\lambda^2\del_{\lambda}$
on $f^{-1}(\nbigr_Y)$ given by 
$\lambda^2\del_{\lambda}(Q):=[\lambda^2\del_{\lambda},Q]$
in $f^{-1}(\nbigrtilde_Y)$.
It induces an action of 
$\lambda^2\del_{\lambda}$
on $\nbigr_{Y\larr X}$,
with which
$\nbigr_{Y\larr X}$
is an $\nbiga_{X,Y,f}$-module.
We have a canonical 
locally free 
$\nbigr_{X,Y,f}$-resolution
$\bigl(\Omegabar_X^{\bullet}[d_X]
 \otimes_{\nbigo_{\nbigx}}\nbigr_X\bigr)
 \otimes_{f^{-1}\nbigo_{\nbigy}}
 f^{-1}(\nbigr_Y\otimes\Omegabar_Y^{-1})$
of $\nbigr_{Y\larr X}$.
It is also an 
$\nbiga_{X,Y,f}$-resolution.
We have the following natural quasi-isomorphism
of $\nbiga_{X,Y,f}$-complexes:
\[
 \bigl(\Omegabar_X^{\bullet}[d_X]
 \otimes_{\nbigo_{\nbigx}}\nbigr_X\bigr)
 \otimes_{f^{-1}\nbigo_{\nbigy}}
 f^{-1}(\nbigr_Y\otimes\Omegabar_Y^{-1})
\lrarr
 \bigl(
 \nbigcbar^{\infty\bullet}_X[d_X]
 \otimes_{\nbigo_{\nbigx}}\nbigr_X
\bigr)
\otimes_{f^{-1}\nbigo_{\nbigy}}
 f^{-1}\bigl(\nbigr_Y\otimes\Omegabar_Y^{-1}\bigr) 
\]

Let $\nbigm^{\bullet}$ be any $\nbigr_X$-complex.
We obtain the following $\nbigr_Y$-complex:
\begin{multline}
\label{eq;14.12.20.1}
 f_{\dagger}\nbigm^{\bullet}:=
 f_{!}\left(
 \Bigl(\bigl(
\nbigcbar^{\infty,\bullet}_X[d_X]
 \otimes_{\nbigo_{\nbigx}}\nbigr_X
\bigr)
\otimes_{f^{-1}\nbigo_{\nbigy}}
 f^{-1}\bigl(\nbigr_Y\otimes\Omegabar_Y^{-1}\bigr)
\Bigr)
 \otimes_{\nbigr_X}
 \nbigm^{\bullet}
\right)
 \\
=f_{!}\left(
\Bigl(
 \nbigcbar^{\infty,\bullet}_X[d_X]
\otimes_{f^{-1}\nbigo_{\nbigy}}
 f^{-1}\bigl(\nbigr_Y\otimes\Omegabar_Y^{-1}\bigr)
\Bigr)
 \otimes_{\nbigo_{\nbigx}}
 \nbigm^{\bullet}
\right)
\end{multline}
It induces a functor 
$D^b(\nbigr_X)\lrarr D^b(\nbigr_Y)$.
We denote the right hand side of (\ref{eq;14.12.20.1}) by 
$Rf_!\bigl(\nbigr_{Y\larr X}
 \otimes_{\nbigr_X}^L\nbigm^{\bullet}\bigr)$.
If $\nbigm^{\bullet}$ is cohomologically $\nbigr_X$-coherent
and good relative to $f$,
then $f_{\dagger}\nbigm^{\bullet}$ is also cohomologically 
$\nbigr_X$-coherent.
If $\nbigm^{\bullet}$ is an $\nbigrtilde_X$-complex,
$f_{\dagger}\nbigm^{\bullet}$ is naturally
an $\nbigrtilde_Y$-complex.
Hence, (\ref{eq;14.12.20.1})
gives a functor 
$D^b(\nbigrtilde_X)\lrarr D^b(\nbigrtilde_Y)$.

\subsubsection{Another expression of the push-forward of
   $\nbigrtilde$-modules}

We put
\[
 \nbigctilde^{\infty,j}_{\nbigx}
 :=\bigoplus_{p+q=j}
 p_{\lambda}^{-1}\Omega_X^{0,q}
 \otimes_{p_{\lambda}^{-1}\nbigc^{\infty}_X}
 \bigl(
 \Omegatilde_{\nbigx}^{p}
 \otimes_{\nbigo_{\nbigx}}
 \nbigc^{\infty}_{\nbigx/\cnum_{\lambda}}
 \bigr).
\]
With the exterior derivative $d$,
we obtain the complexes
$\bigl(
 \Omegatilde_{\nbigx}^{\bullet},d
 \bigr)$
and 
$\bigl(
 \nbigctilde_{\nbigx}^{\infty,\bullet},d
 \bigr)$,
and we have the natural quasi-isomorphism
$\bigl(\Omegatilde_{\nbigx}^{\bullet},d\bigr)
\lrarr
\bigl(\nbigctilde_{\nbigx}^{\infty,\bullet},
 d\bigr)$.

Let $f:X\lrarr Y$ be a morphism of complex manifolds.
We set
$\nbigrtilde_{X,Y,f}:=
 \nbigrtilde_X\otimes_{\cnum}
 f^{-1}\nbigrtilde_Y$.
We set 
$\nbigrtilde_{Y\larr X}:=
 \Omegatilde_{\nbigx}\otimes_{f^{-1}\nbigo_{\nbigy}}
 f^{-1}\bigl(
 \nbigrtilde_Y\otimes\Omegatilde_{\nbigy}^{-1}
 \bigr)$.
It is naturally an $\nbigrtilde_{X,Y,f}$-module.
We have a canonical locally free
$\nbigrtilde_{X,Y,f}$-resolution
$\bigl(
 \Omegatilde_{\nbigx}^{\bullet}[d_X+1]
 \otimes_{\nbigo_{\nbigx}}
 \nbigrtilde_X
 \bigr)
 \otimes_{f^{-1}\nbigo_{\nbigy}}
 f^{-1}(\nbigrtilde_Y\otimes\Omegatilde_{\nbigy}^{-1})$
of $\nbigrtilde_{Y\larr X}$.
We have the following natural quasi-isomorphism
of $\nbigrtilde_X\otimes f^{-1}\nbigrtilde_Y$-complexes:
\[
 \bigl(\Omegatilde_{\nbigx}^{\bullet}[d_X+1]
 \otimes_{\nbigo_{\nbigx}}\nbigrtilde_X\bigr)
 \otimes_{f^{-1}\nbigo_{\nbigy}}
 f^{-1}(\nbigrtilde_Y\otimes\Omegatilde_{\nbigy}^{-1})
\lrarr
 \bigl(
 \nbigctilde^{\infty,\bullet}_{\nbigx}[d_X+1]
 \otimes_{\nbigo_{\nbigx}}\nbigrtilde_X
\bigr)
\otimes_{f^{-1}\nbigo_{\nbigx}}
 f^{-1}\bigl(\nbigrtilde_Y\otimes
 \Omegatilde_{\nbigy}^{-1}\bigr) 
\]

Let $\nbigm^{\bullet}$ be any $\nbigrtilde_X$-complex.
We obtain the following $\nbigrtilde_Y$-complex:
\begin{multline}
 \label{eq;14.12.21.2}
 \ftilde_{\dagger}\nbigm^{\bullet}:=
 f_!\left(
 \Bigl(
 \bigl(
 \nbigctilde^{\infty,\bullet}_{\nbigx}[d_X+1]
 \otimes_{\nbigo_{\nbigx}}\nbigrtilde_X
\bigr)
\otimes_{f^{-1}\nbigo_{\nbigy}}
 f^{-1}\bigl(\nbigrtilde_Y\otimes
 \Omegatilde_{\nbigy}^{-1}\bigr) 
\Bigr)
\otimes_{\nbigrtilde_X}
 \nbigm^{\bullet}
 \right)
\\
=f_!\left(
 \Bigl(
 \nbigctilde^{\infty,\bullet}_{\nbigx}[d_X+1]
\otimes_{f^{-1}\nbigo_{\nbigy}}
 f^{-1}\bigl(\nbigrtilde_Y\otimes
 \Omegatilde_{\nbigy}^{-1}\bigr) 
\Bigr)
\otimes_{\nbigo_{\nbigx}}
 \nbigm^{\bullet}
 \right)
\end{multline}
It induces a functor
$D^b(\nbigrtilde_X)\lrarr D^b(\nbigrtilde_Y)$.
We denote the right hand side of (\ref{eq;14.12.21.2})
by 
$Rf_!\bigl(
 \nbigrtilde_{Y\larr X}\otimes^L_{\nbigrtilde_X}\nbigm
 \bigr)$.

\begin{lem}
\label{lem;14.12.23.1}
We have the following natural isomorphism:
\begin{multline}
\label{eq;14.12.21.10}
 \Bigl(
 \nbigctilde^{\infty,\bullet}_{\nbigx}[d_X+1]
\otimes_{f^{-1}\nbigo_{\nbigy}}
 f^{-1}\bigl(\nbigrtilde_Y\otimes
 \Omegatilde_{\nbigy}^{-1}\bigr) 
\Bigr)
\otimes_{\nbigo_{\nbigx}}
 \nbigm^{\bullet}
\simeq \\
 S_{\lambda}^{\triangledown}
\left(
 \Bigl(
 \nbigcbar^{\infty\bullet}_X[d_X]
\otimes_{f^{-1}\nbigo_{\nbigy}}
 f^{-1}\bigl(\nbigr_Y\otimes\Omegabar_Y^{-1}\bigr)
\Bigr)
 \otimes_{\nbigo_{\nbigx}}
 \nbigm^{\bullet} 
\right)
\end{multline}
\end{lem}
\pf
Let $\ell$ and $r$ denote 
the $\nbigr_Y$-action 
(resp. $\nbigrtilde_Y$-action)
on
$\nbigr_{Y}\otimes\Omegabar_Y^{-1}$
(resp. $\nbigrtilde_Y\otimes\Omegatilde_{\nbigy}^{-1}$)
induced by the left and right multiplications.
Let $\ell\blacktriangledown$
(resp. $r\triangledown$)
denote the $\nbigrtilde_Y$-action on 
$\bigl(
 \nbigr_Y\otimes\Omegabar_Y^{-1}\bigr)
 \langle\lambda^2\del_{\lambda}\rangle$
induced by
$\ell$ and $\blacktriangledown$
(resp. $r$ and $\triangledown$).
We have the isomorphism
$\Omegabar_Y^{-1}\simeq
 \lambda^{-2}\Omegatilde_{\nbigy}^{-1}$
given by 
$\tau\longmapsto
 \tau(d\lambda)^{-1}$.
As mentioned in Example \ref{example;14.12.22.20},
they induce the following isomorphism:
\begin{equation}
\label{eq;14.12.21.4}
 \Bigl(
 \nbigrtilde_Y\otimes
 (\lambda^{-2}\Omegatilde_{\nbigy}^{-1}),
 \ell,r
\Bigr)
 \simeq
 \Bigl(
 \bigl(
 \nbigr_Y\otimes\Omegabar_Y^{-1}
 \bigr)
 \langle\lambda^2\del_{\lambda}\rangle,
 \ell\blacktriangledown,
 r\triangledown
 \Bigr).
\end{equation}
We obtain (\ref{eq;14.12.21.10})
from (\ref{eq;14.12.21.4}).
\hfill\qed

\begin{cor}
\label{cor;14.12.22.141}
We have the natural isomorphism
$\ftilde_{\dagger}\nbigm^{\bullet}
\simeq
 S_{\lambda}^{\triangledown}f_{\dagger}\nbigm^{\bullet}$
of $\nbigrtilde_Y$-complexes.
In particular,
we have a natural isomorphism
$\ftilde_{\dagger}\nbigm^{\bullet}
\lrarr
 f_{\dagger}\nbigm^{\bullet}$
in the derived category of $\nbigrtilde_Y$-modules.
\hfill\qed
\end{cor}

\subsubsection{Trace morphisms}
\label{subsection;14.12.22.130}

Let $\distribution_{\nbigx/\cnum_{\lambda}}$
be the sheaf of distributions $F$ on $\nbigx$
such that $\delbar_{\lambda} F=0$.
We set 
\[
 \distributionbar^j_X:=
 \bigoplus_{p+q=j}
 p_{\lambda}^{-1}\Omega^{0,q}_X
 \otimes_{p_{\lambda}^{-1}\nbigc^{\infty}_X}
 \bigl(
 \Omegabar_X^p\otimes_{\nbigo_{\nbigx}}
 \distribution_{\nbigx/\cnum_{\lambda}}
 \bigr).
\]
With the exterior derivative $d$ in the $X$-direction,
we obtain the complex
$\bigl(
 \distributionbar^{\bullet}_X,d
 \bigr)$.
We have the natural quasi-isomorphism
$(\Omegabar^{\bullet}_X,d)
\lrarr
 (\distributionbar^{\bullet}_X,d)$.
We also set
\[
 \distributiontilde^j_{\nbigx}:=
 \bigoplus_{p+q=j}
 p_{\lambda}^{-1}\Omega^{0,q}_X
 \otimes_{p_{\lambda}^{-1}\nbigc^{\infty}_X}
 \bigl(
 \Omegatilde_{\nbigx}^p\otimes_{\nbigo_{\nbigx}}
 \distributiontilde_{\nbigx/\cnum_{\lambda}}
 \bigr).
\]
With the exterior derivative,
we obtain the complex
$\bigl(
 \distributiontilde^{\bullet}_{\nbigx},d
 \bigr)$.
The natural inclusion 
$(\Omegatilde_{\nbigx}^{\bullet},d)
\lrarr
 (\distributiontilde_{\nbigx}^{\bullet},d)$
is a quasi-isomorphism.

\vspace{.1in}

Let $f:X\lrarr Y$ be a proper morphism.
Recall that we have the trace morphism
$\tr_f:\lambda^{d_X}f_{\dagger}\nbigo_{\nbigx}[d_X]
\lrarr
 \lambda^{d_Y}\nbigo_{\nbigy}[d_Y]$
in the derived category of
$\nbigrtilde_Y$-modules:
\begin{multline}
\label{eq;14.12.23.11}
 \lambda^{d_X}f_{\dagger}\nbigo_{\nbigx}[d_X]
\simeq
 f_!\Bigl(
 \lambda^{d_X}
 \distributionbar_X^{\bullet}[2d_X]
 \otimes_{f^{-1}\nbigo_{\nbigy}}
 f^{-1}\bigl(
 \nbigr_Y\otimes\Omegabar_Y^{-1}\bigr)
 \Bigr)
\simeq
 f_!\bigl(
 \lambda^{d_X}
  \distributionbar_X^{\bullet}[2d_X]
 \bigr)
\otimes_{\nbigo_{\nbigy}}
 \bigl(
 \nbigr_Y\otimes\Omegabar_Y^{-1}
 \bigr)
 \\
\stackrel{a}{\lrarr}
 \lambda^{d_Y}
 \distributionbar_Y^{\bullet}[2d_Y]
\otimes_{\nbigo_{\nbigy}}
 \bigl(
 \nbigr_Y\otimes\Omegabar_Y^{-1}
 \bigr)
\simeq
 \lambda^{d_Y}\nbigo_{\nbigy}[d_Y].
\end{multline}
Here, $a$ is given by the integration of currents along $f$
multiplied by $(2\pi\sqrt{-1})^{-d_X+d_Y}$.
We also have the trace morphism
$\trtilde_f:
 \lambda^{d_X}\ftilde_{\dagger}\nbigo_{\nbigx}[d_X]
\lrarr
 \lambda^{d_Y}\nbigo_{\nbigy}[d_Y]$:
\begin{multline}
\label{eq;14.12.23.10}
 \lambda^{d_X}\ftilde_{\dagger}\nbigo_{\nbigx}[d_X]
\simeq
 f_!\Bigl(
 \lambda^{d_X}\distributiontilde_{\nbigx}^{\bullet}[2d_X+1]
\otimes_{f^{-1}\nbigo_{\nbigy}}
 f^{-1}\bigl(
 \nbigrtilde_Y\otimes\Omegatilde_{\nbigy}^{-1}
 \bigr)
 \Bigr)
\simeq
  f_!\bigl(
 \lambda^{d_X}\distributiontilde_{\nbigx}^{\bullet}[2d_X+1]
 \bigr)
\otimes_{\nbigo_{\nbigy}}
 \bigl(
 \nbigrtilde_Y\otimes\Omegatilde_{\nbigy}^{-1}
 \bigr)
\\
\stackrel{b}{\lrarr}
 \lambda^{d_Y}\distributiontilde_{\nbigy}^{\bullet}[2d_Y+1]
\otimes_{\nbigo_{\nbigy}}
 \bigl(
 \nbigrtilde_Y\otimes\Omegatilde_{\nbigy}^{-1}
 \bigr)
\simeq
 \lambda^{d_Y}
 \nbigo_{\nbigy}[d_Y]
\end{multline}
Here, $b$ is given by the integration of currents along $f$
multiplied by $(2\pi\sqrt{-1})^{-d_X+d_Y}$.
The following is clear by the construction of the morphisms.

\begin{lem}
\label{lem;14.12.23.32}
$\trtilde_f$ is equal to the composite of
$\lambda^{d_X}\ftilde_{\dagger}\nbigo_{\nbigx}[d_X]
 \lrarr
 \lambda^{d_X}f_{\dagger}\nbigo_{\nbigx}[d_X]$
and $\tr_f$.
\end{lem}
\pf
As in Lemma \ref{lem;14.12.23.1},
we can identify the morphism (\ref{eq;14.12.23.10})
with that obtained from (\ref{eq;14.12.23.11})
by the functor $S_{\lambda}^{\triangledown}$.
Hence, the claim is clear.
\hfill\qed

\subsubsection{Complement}
\label{subsection;14.12.23.210}

Let $\nbign$ be any $\nbigrtilde_X$-module.
We have the morphism of sheaves
\begin{multline}
\label{eq;14.12.23.100}
 \Bigl(
 \nbigcbar^{\infty\,\bullet}_X[d_X]
 \otimes_{f^{-1}\nbigo_{\nbigy}}
 f^{-1}(\nbigr_Y\otimes\Omegabar^{-1}_Y)
 \Bigr)
 \otimes_{\nbigo_{\nbigx}}
 \bigl(
 \nbign\langle\lambda^2\del_{\lambda}\rangle,
 \blacktriangledown
 \bigr)
\lrarr
 \\
 \Bigl(
 \bigl(
 \nbigcbar^{\infty\,\bullet}_X[d_X]
 \otimes_{f^{-1}\nbigo_{\nbigy}}
 f^{-1}(\nbigr_Y\otimes\Omegabar^{-1}_Y)
 \bigr)
 \otimes_{\nbigo_{\nbigx}}
 \nbign
 \Bigr)
\langle\lambda^2\del_{\lambda}\rangle
\end{multline}
by 
$g\otimes(\del_{\lambda}^j\otimes\lambda^{2j}m)
\longmapsto
 \del_{\lambda}^j\otimes (g\otimes\lambda^{2j}m)$.
By construction, the following holds.
\begin{lem}
The morphism {\rm(\ref{eq;14.12.23.100})}
is an $\nbigrtilde_X$-homomorphism
with respect to the natural $f^{-1}\nbigrtilde_Y$-action
on the left hand side,
and the $f^{-1}\nbigrtilde_Y$-action $\blacktriangledown$ 
on the right hand side.
The morphism {\rm(\ref{eq;14.12.23.100})}
is compatible with 
the actions $(\del_{\lambda}\lambda^2)\triangledown$
on the both sides.
\hfill\qed
\end{lem}

\begin{cor}
The morphism {\rm(\ref{eq;14.12.23.100})}
induces an isomorphism
$f_{\dagger}\bigl(
 \nbigm^{\bullet}\langle\lambda^2\del_{\lambda}\rangle,
 \blacktriangledown
 \bigr)
\lrarr
\bigl(
 f_{\dagger}(\nbigm^{\bullet})\langle\lambda^2\del_{\lambda}\rangle,
 \blacktriangledown
\bigr)$.
It also induces an isomorphism
$f_{\dagger}S_{\lambda}^{\blacktriangledown}(\nbigm)
\simeq
 S_{\lambda}^{\blacktriangledown}
 f_{\dagger}(\nbigm)$.
\hfill\qed
\end{cor}

\subsection{Duality}
\label{subsection;15.1.4.20}

\subsubsection{Duality of $\nbigr$-modules and $\nbigrtilde$-modules}

Let $X$ be a complex manifold
with a hypersurface $H$. Set $d_X:=\dim X$.
We naturally regard 
$\lambda^{d_X}\nbigr_X(\ast H)\otimes\Omegabar_X^{-1}$
as an $\nbiga_{X}(\ast H)$-module.
We take an $\nbiga_{X}(\ast H)$-injective resolution
$\nbigg_0^{\bullet}$ of 
$\lambda^{d_X}\nbigr_X(\ast H)
 \otimes\Omegabar_X^{-1}[d_X]$.
For any $\nbigr_X(\ast H)$-module $\nbign$,
let $\nhom_{\nbigr_X(\ast H)}(\nbign,\nbigg_0^{p,\ell})$
denote the sheaf of $\nbigr_X(\ast H)$-homomorphisms
$\nbign\lrarr (\nbigg_0^{p},\ell)$.
It is naturally an $\nbigr_X(\ast H)$-module
by $r$.
If $\nbign$ is an $\nbigrtilde_X(\ast H)$-module,
it is naturally an $\nbigrtilde_X(\ast H)$-module.
The action of $\lambda^2\del_{\lambda}$
is given by
$(\lambda^2\del_{\lambda}f)(m)
=\lambda^2\del_{\lambda}(f(m))
-f(\lambda^2\del_{\lambda}m)$.

For any $\nbigr_X(\ast H)$-complex $\nbigm^{\bullet}$,
we obtain the following $\nbigr_X(\ast H)$-complex:
\begin{equation}
\label{eq;14.12.21.11}
 \DDD_{X(\ast H)}\nbigm^{\bullet}:=
 \nhom_{\nbigr_X(\ast H)}\bigl(
 \nbigm^{\bullet},
 \nbigg_0^{\bullet\,\ell}
 \bigr)
\end{equation}
If $\nbigm^{\bullet}$ is an $\nbigrtilde_X(\ast H)$-complex,
then $\DDD_{X(\ast H)}\nbigm^{\bullet}$ is naturally
an $\nbigrtilde_X(\ast H)$-complex.
We shall often denote the right hand side of
(\ref{eq;14.12.21.11})
by 
$\nrhom_{\nbigr_X(\ast H)}\bigl(\nbigm^{\bullet},
 \lambda^{d_X}\nbigr_X(\ast H)\otimes\Omegabar_X^{-1}
 \bigr)[d_X]$
if there is no risk of confusion.

\subsubsection{Another expression of the duality of
 $\nbigrtilde$-modules}
\label{subsection;14.12.23.200}

We give another expression of the duality functor
for $\nbigrtilde$-modules.
We naturally regard 
$\nbigrtilde_X(\ast H)\otimes\Omegatilde_{\nbigx}^{-1}$
as $\nbigrtilde_X(\ast H)\otimes_{\cnum}
 \nbigrtilde_X(\ast H)$-module.
Let $\ell$ (resp. $r$) denote
the $\nbigrtilde_X(\ast H)$-actions
induced by the left multiplication
(resp. the right multiplication).
We take an injective
$\nbigrtilde_X(\ast H)
 \otimes_{\cnum}
 \nbigrtilde_X(\ast H) $-resolution 
$\nbigg_1^{\bullet}$ of
$\lambda^{d_X}\nbigrtilde_X(\ast H)
 \otimes\Omegatilde_{\nbigx}^{-1}[d_X+1]$.
For any $\nbigrtilde_X(\ast H)$-module $\nbign$,
let 
$\nhom_{\nbigrtilde_X(\ast H)}(\nbign,\nbigg_1^{p,\ell})$
denote the sheaf of 
$\nbigrtilde_X(\ast H)$-homomorphisms
$\nbign\lrarr (\nbigg_1^{p},\ell)$.
It is naturally an $\nbigrtilde_X(\ast H)$-module
induced by the $\nbigrtilde_X(\ast H)$-action $r$.

For any $\nbigrtilde_X(\ast H)$-complex 
$\nbigm^{\bullet}$,
we obtain the following $\nbigrtilde_X(\ast H)$-complex:
\begin{equation}
\label{eq;14.12.21.20}
 \DDDtilde_{X(\ast H)}(\nbigm^{\bullet}):=
 \nhom_{\nbigrtilde_X(\ast H)}
 \bigl(
 \nbigm^{\bullet},
 \nbigg_1^{\bullet\,\ell}
 \bigr)
\end{equation}
We shall denote the right hand side of
(\ref{eq;14.12.21.20})
by 
$\nrhom_{\nbigrtilde_X(\ast H)}\bigl(
 \nbigm^{\bullet},
 \nbigrtilde_X\otimes\Omegatilde_{\nbigx}^{-1}
 \bigr)[d_X+1]$
if there is no risk of confusion.

\vspace{.1in}
We have the $\nbigrtilde_X(\ast H)$-action
$\ell\blacktriangledown$
on $(\nbigr_X(\ast H)\otimes\Omegabar_X^{-1})
 \langle\lambda^2\del_{\lambda}\rangle$
induced by
the $\nbigrtilde_X(\ast H)$-action $\ell$
and the construction $\blacktriangledown$.
We also have the $\nbigrtilde_X(\ast H)$-action
$r\triangledown$
on $(\nbigr_X(\ast H)\otimes\Omegabar_X^{-1})\langle
 \lambda^2\del_{\lambda}\rangle$
induced by the $\nbigr(\ast H)$-action $r$
and the construction $\triangledown$.
Similarly, we have the $\nbigrtilde_X(\ast H)$-actions
$\ell\blacktriangledown$
and $r\triangledown$ 
on $\nbigg^{\bullet}_0\langle\lambda^2\del_{\lambda}\rangle$.
As in (\ref{eq;14.12.21.4}),
we have the isomorphism
\begin{equation}
\label{eq;14.12.21.31}
\Bigl(
\bigl(
 \nbigr_X(\ast H)\otimes\Omegabar_X^{-1}
 \bigr)
 \langle\lambda^2\del_{\lambda}
 \rangle,
 \ell\blacktriangledown,
 r\triangledown
\Bigr)
\simeq
 \Bigl(
 \nbigrtilde_X(\ast H)\otimes
 (\lambda^{-2}\Omegatilde_{\nbigx}^{-1}),
 \ell,r
\Bigr).
\end{equation}
Hence, we may assume 
to have a quasi-isomorphism
of $\nbigrtilde_X(\ast H)
 \otimes_{\cnum}\nbigrtilde_X(\ast H)$-complexes
\begin{equation}
\label{eq;14.12.21.100}
 \lambda^{2}\nbigg_0^{\bullet}
 \langle\lambda^2\del_{\lambda}\rangle[1]
\lrarr
 \nbigg_1^{\bullet}.
\end{equation}
From Corollary \ref{cor;14.12.21.101}
and (\ref{eq;14.12.21.100}),
we have the following morphisms
of $\nbigrtilde_X(\ast H)$-complexes:
\begin{equation}
\label{eq;14.12.21.112}
S_{\lambda}^{\triangledown}
\nhom_{\nbigr_X(\ast H)}\bigl(
 \nbigm^{\bullet},
 \nbigg_0^{\bullet\,\ell}
\bigr)
\simeq
\nhom_{\nbigrtilde_X(\ast H)}\bigl(
 S_{\lambda}^{\blacktriangledown}
 \nbigm^{\bullet},
 \lambda^2\nbigg_0^{\bullet}
 \langle\lambda^2\del_{\lambda}\rangle
 ^{\ell\blacktriangledown}
 [1]
 \bigr)
\lrarr
\nhom_{\nbigrtilde_X(\ast H)}\bigl(
 S_{\lambda}^{\blacktriangledown}
 \nbigm^{\bullet},
 \nbigg_1^{\bullet\,\ell}
 \bigr)
\end{equation}
We also have the following natural quasi-isomorphisms:
\begin{equation}
\label{eq;14.12.21.110}
 S_{\lambda}^{\triangledown}
\nhom_{\nbigr_X(\ast H)}\bigl(
 \nbigm^{\bullet},
 \nbigg_0^{\bullet\,\ell}
\bigr)
\lrarr
 \nhom_{\nbigr_X(\ast H)}\bigl(
 \nbigm^{\bullet},
 \nbigg_0^{\bullet,\ell}
 \bigr)
\end{equation}
\begin{equation}
\label{eq;14.12.21.111}
\nhom_{\nbigrtilde_X(\ast H)}\bigl(
 \nbigm^{\bullet},
 \nbigg_1^{\bullet\,\ell}
 \bigr)
\lrarr
\nhom_{\nbigrtilde_X(\ast H)}\bigl(
 S_{\lambda}^{\blacktriangledown}
 \nbigm^{\bullet},
 \nbigg_1^{\bullet\,\ell}
 \bigr)
\end{equation}

\begin{prop}
\label{prop;14.12.17.1}
If $\nbigm^{\bullet}$ is 
cohomologically $\nbigr_X(\ast H)$-coherent,
then {\rm(\ref{eq;14.12.21.112})}
is a quasi-isomorphism.
Together with {\rm(\ref{eq;14.12.21.110})}
and {\rm(\ref{eq;14.12.21.111})},
we obtain the following isomorphisms
in the derived category of 
cohomologically
$\nbigrtilde_X(\ast H)$-coherent complexes:
\[
\begin{CD}
 \DDDtilde_{X}\nbigm^{\bullet}
@>{\simeq}>>
 \DDDtilde_X S_{\lambda}^{\blacktriangledown}\nbigm^{\bullet}
@>{\simeq}>>
S_{\lambda}^{\triangledown}
 \DDD_X\nbigm^{\bullet}
@>{\simeq}>>
 \DDD_X\nbigm^{\bullet}
\end{CD}
\]
\end{prop}
\pf
For simplicity of the description,
we omit to denote $(\ast H)$
in the proof.
It is enough to prove the claim for 
an $\nbigrtilde_X$-module $\nbigm$
which is coherent over $\nbigr_X$.
We begin with the following lemma.

\begin{lem}
\label{lem;14.12.16.40}
Let $0\lrarr\nbign_1\lrarr\nbign_2\lrarr\nbign_3\lrarr 0$
be an exact sequence of
coherent $\nbigr_X$-modules.
Let $\nbigi$ be 
any injective $\nbigr_X$-module.
Then, the following is exact.
\begin{multline}
\label{eq;14.12.16.20}
 0\lrarr
 \nhom_{\nbigrtilde_X}
 \bigl(
 \nbign_3\langle\lambda^2\del_{\lambda}\rangle,
 \nbigi\langle\lambda^2\del_{\lambda}\rangle
 \bigr)
\lrarr
 \nhom_{\nbigrtilde_X}
 \bigl(
 \nbign_2\langle\lambda^2\del_{\lambda}\rangle,
 \nbigi\langle\lambda^2\del_{\lambda}\rangle
 \bigr)
 \\
\lrarr
 \nhom_{\nbigrtilde_X}
 \bigl(
 \nbign_1\langle\lambda^2\del_{\lambda}\rangle,
 \nbigi\langle\lambda^2\del_{\lambda}\rangle
 \bigr)
\lrarr 0
\end{multline}
\end{lem}
\pf
The sequence (\ref{eq;14.12.16.20})
is rewritten as follows:
\begin{equation}
\label{eq;14.12.16.21}
 0\lrarr
 \nhom_{\nbigr_X}
 \bigl(
 \nbign_3,
 \nbigi\langle\lambda^2\del_{\lambda}\rangle
 \bigr)
\lrarr
 \nhom_{\nbigr_X}
 \bigl(
 \nbign_2,
 \nbigi\langle\lambda^2\del_{\lambda}\rangle
 \bigr)
\stackrel{\beta}{\lrarr}
 \nhom_{\nbigr_X}
 \bigl(
 \nbign_1,
 \nbigi\langle\lambda^2\del_{\lambda}\rangle
 \bigr)
\lrarr 0
\end{equation}
It is enough to prove that $\beta$ is an epimorphism.
For any non-negative integer $L$,
we set $F_L\nbigi\langle\lambda^2\del_{\lambda}\rangle:=
 \bigoplus_{j=0}^L\del_{\lambda}^j\otimes\lambda^{2j}\nbigi$.
They are naturally an $\nbigr_X$-submodules
of $\nbigi\langle\lambda^2\del_{\lambda}\rangle$.
Let $F_L\nhom_{\nbigr_X}
 \bigl(\nbign_i,
 \nbigi\langle\lambda^2\del_{\lambda}\rangle
 \bigr)$
be the image of 
$\nhom_{\nbigr_X}
 \bigl(\nbign_i,
 F_L\nbigi\langle\lambda^2\del_{\lambda}\rangle
 \bigr)$.
Because $\nbign_i$ is $\nbigr_X$-coherent,
the filtration $F$ is exhaustive on
$\nhom_{\nbigr_X}
 \bigl(\nbign_i,
 \nbigi\langle\lambda^2\del_{\lambda}\rangle
 \bigr)$.
By construction,
we have 
$\Gr^L_j\nhom_{\nbigr_X}\bigl(
 \nbign_i,
 \nbigi\langle\lambda^2\del_{\lambda}\rangle
\bigr)
\simeq
 \nhom_{\nbigr_X}\bigl(
 \nbign_i,
 \lambda^{2j}\nbigi
\bigr)$.
Because $\nbigi$ is $\nbigr$-injective,
we have the surjectivity of
$\Gr^L_j\nhom_{\nbigr_X}\bigl(
 \nbign_2,
 \nbigi\langle\lambda^2\del_{\lambda}\rangle
\bigr)
\lrarr
\Gr^L_j\nhom_{\nbigr_X}\bigl(
 \nbign_1,
 \nbigi\langle\lambda^2\del_{\lambda}\rangle
\bigr)$.
Then, we obtain the desired surjectivity.
\hfill\qed

\vspace{.1in}
The claim of Proposition \ref{prop;14.12.17.1}
is reduced to the following lemma.

\begin{lem}
\label{lem;14.12.16.33}
For any coherent $\nbigr_X$-module
$\nbign$,
the natural morphism
\begin{equation}
\label{eq;14.12.16.32}
 \nhom_{\nbigrtilde_X}
 \bigl(
 \nbign\langle\lambda^2\del_{\lambda}\rangle,
 \lambda^2\nbigg_0^{\bullet}
 \langle\lambda^2\del_{\lambda}\rangle
 ^{\ell\blacktriangledown}
 [1]
 \bigr)
\lrarr 
 \nhom_{\nbigrtilde_X}
 \bigl(
 \nbign\langle\lambda^2\del_{\lambda}\rangle,
 \nbigg_1^{\bullet\,\ell}
 \bigr)
\end{equation}
is a quasi-isomorphism.
\end{lem}
\pf
The morphism (\ref{eq;14.12.16.32}) is a quasi-isomorphism
if $\nbign=\nbigr_X$.
It is enough to prove that 
(\ref{eq;14.12.16.32}) is a quasi-isomorphism
locally around any point of $X$.
We may assume to have an $\nbigr_X$-free resolution
$\nbigp^{\bullet}$ of $\nbign$.
We have the following commutative diagram:
\[
 \begin{CD}
 \nhom_{\nbigrtilde_X}
 \bigl(
 \nbign\langle\lambda^2\del_{\lambda}\rangle,
 \lambda^2\nbigg_0^{\bullet}
 \langle\lambda^2\del_{\lambda}\rangle
 ^{\ell\blacktriangledown}
 [1]
 \bigr)
 @>{\alpha_1}>>
 \nhom_{\nbigrtilde_X}
 \bigl(
 \nbign\langle\lambda^2\del_{\lambda}\rangle,
 \nbigg_1^{\bullet\,\ell}
 \bigr)
 \\
 @V{\alpha_2}VV @V{\alpha_3}VV \\
\Tot
 \nhom_{\nbigrtilde_X}
 \bigl(
 \nbigp^{\bullet}\langle\lambda^2\del_{\lambda}\rangle,
 \lambda^2\nbigg_0^{\bullet}
 \langle\lambda^2\del_{\lambda}\rangle
 ^{\ell\blacktriangledown}
 [1]
 \bigr)
 @>{\alpha_4}>>
 \Tot
 \nhom_{\nbigrtilde_X}
 \bigl(
 \nbigp^{\bullet}\langle\lambda^2\del_{\lambda}\rangle,
 \nbigg_1^{\bullet\,\ell}
 \bigr)
 \end{CD}
\]
Here, $\Tot\nbigc^{\bullet\bullet}$ denote the total complex
of a double complex $\nbigc^{\bullet\bullet}$.
As remarked above,
$\alpha_4$ is a quasi-isomorphism.
Because $(\nbigg_1^{\bullet},\ell)$ is $\nbigrtilde_X$-injective,
$\alpha_3$ is a quasi-isomorphism.
As for $\alpha_2$,
note that $(\lambda^2\nbigg_0^{\bullet},\ell)$ is 
$\nbigr_X$-injective.
We have the isomorphism
exchanging $\triangledown$ and $\blacktriangledown$
as in \S\ref{subsection;14.12.22.1}.
Hence, $\alpha_2$ is a quasi-isomorphism
by Lemma \ref{lem;14.12.16.40}.
Then, we obtain that $\alpha_1$ is a quasi-isomorphism.
Thus we obtain Lemma \ref{lem;14.12.16.33},
and hence Proposition \ref{prop;14.12.17.1}.
\hfill\qed

\subsubsection{Smooth case}

Let $\nbigm$ be a smooth $\nbigrtilde_X(\ast H)$-module.
Let $\nbigm^{\lor}$ be the smooth $\nbigrtilde_X(\ast H)$-module
$\nhom_{\nbigo_{\nbigx}(\ast H)}
 \bigl(\nbigm,\nbigo_{\nbigx}(\ast H)\bigr)$
with the induced meromorphic connection.
As remarked in \cite{Mochizuki-MTM},
we have a natural isomorphism
$\DDD_{X(\ast H)}\nbigm
 \simeq \lambda^{d_X}\nbigm^{\lor}$
in the derived category.
Together with Proposition \ref{prop;14.12.17.1},
we obtain
$\DDDtilde_{X(\ast H)
 }\nbigm\simeq \lambda^{d_X}\nbigm^{\lor}$.

We also have the isomorphism
$\DDDtilde_{X(\ast H)}\nbigm
\simeq
 \lambda^{d_X}\nbigm^{\lor}$
given in a more direct and standard way,
by using the Spencer resolution
$\nbigrtilde_X(\ast H)
 \otimes\Thetatilde^{\bullet}_{\nbigx}$
of $\nbigo_{\nbigx}(\ast H)$:
\begin{multline}
\label{eq;14.12.17.102}
 \nrhom_{\nbigrtilde_X(\ast H)}
 \bigl(
 \nbigm,
 \lambda^{d_X}\nbigrtilde_X(\ast H)
 \otimes \Omegatilde_{\nbigx}^{-1}
 \bigr)[d_X+1]
\simeq
 \\
 \nhom_{\nbigrtilde_X(\ast H)}
 \bigl(
 \nbigrtilde_X(\ast H)
 \otimes\Thetatilde^{\bullet}_{\nbigx}
 \otimes\nbigm,
\lambda^{d_X}\nbigrtilde_X(\ast H)
 \otimes \Omegatilde_{\nbigx}^{-1}
 \bigr)[d_X+1]
\simeq 
\\
 \nbigrtilde_X(\ast H)
 \otimes\Thetatilde^{\bullet}_{\nbigx}
 \otimes\lambda^{d_X}\nbigm^{\lor}
\simeq
 \lambda^{d_X}\nbigm^{\lor}
\end{multline}

\begin{lem}
\label{lem;14.12.23.220}
The isomorphisms
$\DDDtilde_{X(\ast H)}\nbigm
 \simeq\lambda^{d_X}\nbigm^{\lor}$
are equal.
\end{lem}
\pf
We have the natural quasi-isomorphisms
$S^{\blacktriangledown}_{\lambda}
 \bigl((\nbigr_X(\ast H)
 \otimes_{\nbigo_{\nbigx}}\Thetabar_{X}^{\bullet})
 \otimes_{\nbigo_{\nbigx}}\nbigm\bigr)
\lrarr
S^{\blacktriangledown}_{\lambda}(\nbigm)
\lrarr \nbigm$.
We have the following natural isomorphisms:
\[
S^{\blacktriangledown}_{\lambda}
 \bigl((\nbigr_X(\ast H)
 \otimes_{\nbigo_{\nbigx}}\Thetabar_{\nbigx}^{\bullet})
 \otimes_{\nbigo_{\nbigx}}\nbigm\bigr)
\simeq
S^{\triangledown}_{\lambda}
 \bigl((\nbigr_X(\ast H)
 \otimes_{\nbigo_{\nbigx}}\Thetabar_{\nbigx}^{\bullet})
 \otimes_{\nbigo_{\nbigx}}\nbigm\bigr)
\simeq
 \bigl(
 \nbigrtilde_X(\ast H)
 \otimes_{\nbigo_{\nbigx}}
 \Thetatilde^{\bullet}_{\nbigx}
 \bigr)
 \otimes_{\nbigo_{\nbigx}}\nbigm
\]
The isomorphism in Proposition \ref{prop;14.12.17.1}
is expressed as follows:
\begin{multline}
\nhom_{\nbigrtilde_X(\ast H)}
 \bigl(
 \nbigrtilde_{X}(\ast H)
 \otimes\Thetatilde^{\bullet}_{\nbigx}
 \otimes\nbigm,
  \nbigrtilde_{X}(\ast H)
 \otimes \Omegatilde_{\nbigx}^{-1}
 \bigr)[d_X+1]
\simeq
 \\
\nhom_{\nbigrtilde_X(\ast H)}
 \bigl(
 S^{\blacktriangledown}_{\lambda}\bigl(
 \nbigr_X(\ast H)\otimes\Thetabar^{\bullet}_X\otimes
 \nbigm
 \bigr),
  \nbigrtilde_{X}(\ast H)
 \otimes \Omegatilde_{\nbigx}^{-1}
 \bigr)[d_X+1]
\simeq \\
 S^{\triangledown}_{\lambda}
 \nhom_{\nbigr_X(\ast H)}
\bigl(
 \nbigr_X(\ast H)\otimes\Thetabar^{\bullet}_X\otimes
 \nbigm,
  \nbigr_{X}(\ast H)
 \otimes \Omegabar_{X}^{-1}
 \bigr)[d_X]
\simeq
 \\
 \nhom_{\nbigr_X(\ast H)}\bigl(
 \nbigr_X(\ast H)\otimes\Thetabar^{\bullet}_X\otimes
 \nbigm,
  \nbigr_{X}(\ast H)
 \otimes \Omegabar_{X}^{-1}
 \bigr)[d_X]
\end{multline}
Then, the claim is clear by the constructions
of the isomorphisms
$\DDDtilde_{X(\ast H)}\nbigm
 \simeq\lambda^{d_X}\nbigm^{\lor}$.
\hfill\qed

\subsection{Compatibility of push-forward and duality}

\subsubsection{Compatibility for $\nbigr$-modules}
\label{subsection;14.12.22.110}

Let us recall the compatibility of the push-forward
and the duality for $\nbigr$-modules studied in
\cite{Mochizuki-MTM}.
Let $Y$ be a complex manifold.

We set 
$\nbigm_{Y0}:=(\nbigr_Y\otimes\Omegabar_Y^{-1})
 \lefttop{\ell}\otimes^{\ell}_{\nbigo_{\nbigy}}
 (\nbigr_Y\otimes\Omegabar_Y^{-1})$,
where 
$\lefttop{\ell}\otimes^{\ell}$
means that we consider the actions $\ell$
for the tensor product over $\nbigo_{\nbigy}$.
It is equipped with the $\nbigr_Y$-actions $r_i$ $(i=1,2)$
induced by the action $r$ on the $i$-th factor.
We also have the $\nbigr_Y$-action
$\ell\otimes\ell$ induced by the actions $\ell$.
We have the action of $\lambda^2\del_{\lambda}$
on $\nbigm_{Y0}$
induced by the actions on $\nbigr_Y\otimes\Omegabar_Y^{-1}$
with the Leibniz rule. 
Each $\nbigr_Y$-action with the action of $\lambda^2\del_{\lambda}$
gives an $\nbigrtilde_Y$-action.

We set 
$\nbigm_{Y1}:=
 (\nbigr_Y\otimes\Omegabar_Y^{-1})
 \lefttop{r}\otimes^{\ell}_{\nbigo_{\nbigy}}
 (\nbigr_Y\otimes\Omegabar_Y^{-1})$,
where 
$\lefttop{r}\otimes^{\ell}$ means that
we consider the action $r$ on the first factor
and the action $\ell$ on the second factor
for the tensor product over $\nbigo_{\nbigy}$.
It is equipped with 
the $\nbigr_Y$-action $\ell_1$
induced by the action $\ell$ on the first factor,
and the $\nbigr_Y$-action $r_2$
induced by the action $r$ on the second factor.
We also have the $\nbigr_Y$-action $r\otimes\ell$
induced by the action $r$ on the first factor 
and the action $\ell$ on the second factor.
As in the case of $\nbigm_{Y0}$,
we have the action of $\lambda^2\del_{\lambda}$
on $\nbigm_{Y1}$
with which each $\nbigr_Y$-action is extended 
to an $\nbigrtilde_Y$-action.

Recall that we have a unique isomorphism of sheaves
$\Psi_1:\nbigm_{Y1}\lrarr\nbigm_{Y0}$
such that 
(i) $\Psi_1\circ r_2=r_2\circ \Psi_1$,
$\Psi_1\circ\ell_1=(\ell\otimes\ell)\circ\Psi_1$
and 
$\Psi_1\circ(r\otimes\ell)=r_1\circ\Psi_1$,
(ii) $\Psi_1$ induces the identity on
$\Omegabar_Y^{-1}\otimes\Omegabar_Y^{-1}$.
Here, 
$\Psi_1\circ r_2=r_2\circ\Psi_1$ means
$\Psi_1(r_2(P)m)=r_2(P)\Psi_1(m)$
for any local sections $P\in\nbigr_X$ and $m\in \nbigm_{Y1}$,
and 
the meaning of $\Psi_1\circ\ell_1=(\ell\otimes\ell)\circ\Psi_1$
and $\Psi_1\circ(r\otimes\ell)=r_1\circ\Psi_1$ are similar.
The uniqueness is clear by the conditions
(i) and (ii).
Locally, $\Psi$ is given as follows.
For any local frame $\tau$ of $\Omegabar_X^{-1}$,
we have 
$\Psi_1
\bigl((Q\otimes\tau)\otimes m \bigr)
=(\ell\otimes\ell)(Q)\bigl((1\otimes\tau)\otimes m\bigr)$
for $Q\in\nbigr_X$ and $m\in \nbigr_X\otimes\Omegabar_X^{-1}$.
It is compatible with the actions of $\lambda^{2}\del_{\lambda}$.

\vspace{.1in}

Let $f:X\lrarr Y$ be a proper morphism.
We obtain the morphism 
\begin{equation}
\label{eq;14.12.22.100}
 \lambda^{d_X}
 f_{\dagger}\Bigl(
 \nbigo_{\nbigx}
 \otimes_{f^{-1}\nbigo_{\nbigy}}
 f^{-1}\bigl(\nbigr_Y\otimes\Omegabar_Y^{-1}\bigr)
 \Bigr)[d_X]
\lrarr
 \lambda^{d_Y}
 \nbigr_Y\otimes\Omegabar_Y^{-1}[d_Y]
\end{equation}
in the derived category of $\nbiga_Y$-modules
as the composite of the following morphisms
(see \S\ref{subsection;14.12.24.2}
for $\nbiga_Y$):
\begin{multline}
\lambda^{d_X}
 f_{\dagger}\Bigl(
 \nbigo_{\nbigx}
 \otimes_{f^{-1}\nbigo_{\nbigy}}
 f^{-1}\bigl(\nbigr_Y\otimes\Omegabar_Y^{-1}\bigr)
 \Bigr)[d_X]
 \\
\simeq
 \lambda^{d_X}
 f_{!}\Bigl(
\Bigl(
 \bigl(
 \distributionbar^{\bullet}_X[2d_X]
 \otimes \nbigr_X
 \bigr)
\otimes_{f^{-1}\nbigo_{\nbigy}}
 f^{-1}(\nbigr_Y\otimes\Omegabar_Y^{-1})
\Bigr)\otimes_{\nbigr_X}
 \bigl(
 \nbigo_{\nbigx}
 \otimes_{f^{-1}\nbigo_{\nbigy}}
 f^{-1}\bigl(\nbigr_Y\otimes\Omegabar_Y^{-1}\bigr)
 \bigr)
 \Bigr)
\\
\simeq
 \lambda^{d_X}
 f_!\Bigl(
 \distributionbar^{\bullet}_X[2d_X]
\otimes_{f^{-1}\nbigo_{\nbigy}}
f^{-1}\nbigm_{Y0}
 \Bigr)
\stackrel{a_1}{\simeq}
 \lambda^{d_X}
 f_!\Bigl(
 \distributionbar^{\bullet}_X[2d_X]
\otimes_{f^{-1}\nbigo_{\nbigy}}
f^{-1}\nbigm_{Y1}
 \Bigr)
 \\
\simeq
 \lambda^{d_X}
 f_!\Bigl(
  \distributionbar^{\bullet}_X[2d_X]
 \otimes_{f^{-1}\nbigo_{\nbigy}}
 f^{-1}(\nbigr_Y\otimes\Omegabar_Y^{-1})
 \Bigr)
\otimes_{\nbigo_{\nbigy}}
 \bigl(
 \nbigr_Y\otimes\Omegabar_{Y}^{-1}
 \bigr)
\\
\simeq
 \lambda^{d_X}f_{\dagger}(\nbigo_{\nbigx}[d_X])
 \otimes_{\nbigo_{\nbigy}}
 \bigl(\nbigr_Y\otimes\Omegabar_Y^{-1}\bigr)
\stackrel{a_2}{\lrarr}
 \lambda^{d_Y}
 \bigl(\nbigr_Y\otimes\Omegabar_Y^{-1}\bigr)[d_Y]
\end{multline}
Here, $a_1$ is induced by $\Psi_1^{-1}$,
and $a_2$ is induced by the trace morphism
(see \S\ref{subsection;14.12.22.130}).

Let $\nbigm$ be a coherent $\nbigr_X$-module
which is good relative to $f$.
Then, we have the natural morphism
\begin{equation}
\label{eq;14.12.22.101}
 f_{\dagger}\DDD_X\nbigm\lrarr
 \DDD_Y f_{\dagger}\nbigm
\end{equation}
obtained as the composite of the following morphisms:
\begin{multline}
f_{\dagger}\DDD_X\nbigm
\simeq
 Rf_!\nrhom_{\nbigr_X}\Bigl(
 \nbigm,
 \lambda^{d_X}\nbigo_{\nbigx}
 \otimes_{f^{-1}\nbigo_{\nbigy}}
 f^{-1}(\nbigr_Y\otimes\Omegabar_Y^{-1})[d_X]
 \Bigr)
 \\
\lrarr
 Rf_!\nrhom_{f^{-1}\nbigr_Y}\Bigl(
 \nbigr_{Y\larr X}\otimes^L_{\nbigr_X}\nbigm,
 \nbigr_{Y\larr X}\otimes^L_{\nbigr_X}
 \bigl(
 \lambda^{d_X}\nbigo_{\nbigx}
 \otimes_{f^{-1}\nbigo_{\nbigy}}
 f^{-1}(\nbigr_Y\otimes\Omegabar_Y^{-1})
 \bigr)[d_X]
 \Bigr)
 \\
\lrarr
 \nrhom_{\nbigr_Y}\Bigl(
 f_{\dagger}\nbigm,
 \lambda^{d_X}f_{\dagger}\Bigl(
 \nbigo_{\nbigx}\otimes_{f^{-1}\nbigo_{\nbigy}}
 f^{-1}(\nbigr_Y\otimes\Omegabar_Y^{-1})
 \Bigr)[d_X]
 \Bigr)
 \\
\lrarr
 \nrhom_{\nbigr_Y}\Bigl(
 f_{\dagger}\nbigm,
 \lambda^{d_Y}\nbigr_Y\otimes\Omegabar_Y^{-1}[d_Y]
 \Bigr)
\simeq
 \DDD_Y f_{\dagger}\nbigm
\end{multline}
If $\nbigm$ is also an $\nbigrtilde_X$-module,
then (\ref{eq;14.12.22.101})
is a morphism in the derived category of
$\nbigrtilde_Y$-modules.

\subsubsection{Compatibility for $\nbigrtilde$-modules}
\label{subsection;14.12.22.140}

We set 
$\nbigmtilde_{Y0}:=
 (\nbigrtilde_Y\otimes\Omegatilde_{\nbigy}^{-1})
 \lefttop{\ell}\otimes^{\ell}_{\nbigo_{\nbigy}}
 (\nbigrtilde_Y\otimes\Omegatilde_{\nbigy}^{-1})$.
The tensor product is taken as in the case of
$\nbigm_{Y0}$.
We have the $\nbigrtilde_X$-actions
$r_i$ $(i=1,2)$ and $\ell\otimes\ell$
on $\nbigmtilde_{Y0}$
as in the case of $\nbigm_{Y0}$.
We set 
$\nbigmtilde_{Y1}:=
 (\nbigrtilde_Y\otimes\Omegatilde_{\nbigy}^{-1})
 \lefttop{r}\otimes^{\ell}_{\nbigo_{\nbigy}}
 (\nbigrtilde_Y\otimes\Omegatilde_{\nbigy}^{-1})$,
where the tensor product is taken 
as in the case of $\nbigm_{Y1}$.
We have the $\nbigrtilde_X$-actions
$r_2$, $\ell_1$ and $r\otimes\ell$
on $\nbigmtilde_{Y1}$
as in the case of $\nbigm_{Y1}$.

We have a unique isomorphism of sheaves
$\Psitilde_1:\nbigmtilde_{Y1}\lrarr\nbigmtilde_{Y0}$
such that 
(i) $\Psitilde_1\circ r_2=r_2\circ \Psitilde_1$,
$\Psitilde_1\circ\ell_1=(\ell\otimes\ell)\circ\Psitilde_1$,
(ii) $\Psitilde_1$ induces the identity on
$\Omegatilde_{\nbigx}^{-1}
 \otimes
 \Omegatilde_{\nbigx}^{-1}$.
For any local frame $\tau$ of
$\Omegatilde_{\nbigy}^{-1}$,
we have
$\Psitilde_1\bigl(
 (Q\otimes\tau)\otimes m
 \bigr)=
 (\ell\otimes\ell)(Q)
 \bigl((1\otimes\tau)\otimes m\bigr)$.

\begin{lem}
\label{lem;14.12.24.1}
We have 
$\Psitilde_1\circ(r\otimes\ell)
=r_1\circ\Psitilde_1$.
\end{lem}
\pf
It is enough to prove the equality
locally around any point of 
$\nbigx\setminus\nbigx^0$.
We may assume to have a frame $\tau$
of $\Omegatilde_{\nbigy}^{-1}$
such that
$r(v)(\tau^{-1})=0$
for any $v\in\Thetatilde_{\nbigx}$.
Let $Q\in\nbigrtilde_Y$ and 
$m\in \nbigrtilde_Y\otimes\Omegatilde_{\nbigy}^{-1}$.
For any $g\in\nbigo_{\nbigy}$,
we have
\begin{multline}
 r_1(g)\Psitilde_1\bigl((Q\otimes\tau)\otimes m\bigr)
=r_1(g)(\ell\otimes\ell)(Q)\bigl((1\otimes\tau)\otimes m\bigr)
=(\ell\otimes\ell)(Q)\bigl((g\otimes\tau)\otimes m\bigr)
 \\
=(\ell\otimes\ell)(Q)\bigl((1\otimes\tau)\otimes gm\bigr)
=\Psitilde_1\bigl((Q\otimes\tau)\otimes gm\bigr)
=\Psitilde_1\bigl(
 (r\otimes\ell)(g)\bigl((Q\otimes\tau)\otimes m\bigr)
 \bigr).
\end{multline}
For any $v\in \Thetatilde_{\nbigx}$,
we have
\begin{multline}
 \Psitilde_1\Bigl(
 (r\otimes\ell)(v)\bigl((Q\otimes\tau)\otimes m\bigr)
 \Bigr)
=\Psitilde_1\bigl((-Qv\otimes\tau)\otimes m\bigr)
+\Psitilde(Q\otimes vm)
 \\
=-(\ell\otimes\ell)(Qv)\bigl((1\otimes\tau)\otimes m\bigr)
+(\ell\otimes\ell)(Q)\bigl((1\otimes\tau)\otimes vm\bigr)
=(\ell\otimes\ell)(Q)\bigl((-v\otimes\tau)\otimes m\bigr)
\\
=(\ell\otimes\ell)(Q)r_1(v)
 \bigl((1\otimes\tau)\otimes m\bigr)
=r_1(v)\Psitilde_1\bigl((Q\otimes\tau)\otimes m\bigr)
\end{multline}
Thus, we are done.
\hfill\qed

\vspace{.1in}

Let $f:X\lrarr Y$ be a proper morphism
as in \S\ref{subsection;14.12.22.110}.
We obtain the morphism
\begin{equation}
\label{eq;14.12.22.120}
 \lambda^{d_X}
 \ftilde_{\dagger}\Bigl(
 \nbigo_{\nbigx}
 \otimes_{f^{-1}\nbigo_{\nbigy}}
 f^{-1}\bigl(\nbigrtilde_Y\otimes\Omegatilde_{\nbigy}^{-1}\bigr)
 \Bigr)[d_X+1]
\lrarr
 \lambda^{d_Y}
 \nbigrtilde_Y\otimes\Omegatilde_{\nbigy}^{-1}[d_Y+1]
\end{equation}
in the derived category of 
$\nbigrtilde_Y\otimes_{\cnum}
 \nbigrtilde_Y$-modules
as the composite of the following morphisms:
\begin{multline}
\lambda^{d_X}
 \ftilde_{\dagger}\Bigl(
 \nbigo_{\nbigx}
 \otimes_{f^{-1}\nbigo_{\nbigy}}
 f^{-1}\bigl(\nbigrtilde_Y\otimes
 \Omegatilde_{\nbigy}^{-1}\bigr)
 \Bigr)[d_X+1]
 \\
\simeq
 \lambda^{d_X}
 f_{!}\Bigl[
\Bigl(
 \bigl(
 \distributiontilde^{\bullet}_{\nbigx/\cnum_{\lambda}}[2d_X+2]
 \otimes \nbigrtilde_X
 \bigr)
\otimes_{f^{-1}\nbigo_{\nbigy}}
 f^{-1}(\nbigrtilde_Y\otimes\Omegatilde_{\nbigy}^{-1})
\Bigr)\otimes_{\nbigrtilde_X}
 \bigl(
 \nbigo_{\nbigx}
 \otimes_{f^{-1}\nbigo_{\nbigy}}
 f^{-1}\bigl(\nbigrtilde_Y\otimes\Omegatilde_{\nbigy}^{-1}\bigr)
 \bigr)
\Bigr]
\\
\simeq
 \lambda^{d_X}
 f_!\Bigl(
 \distributiontilde^{\bullet}_{\nbigx/\cnum_{\lambda}}[2d_X+2]
\otimes_{f^{-1}\nbigo_{\nbigy}}
f^{-1}\nbigmtilde_{Y0}
 \Bigr)
\stackrel{\atilde_1}{\simeq}
 \lambda^{d_X}
 f_!\Bigl(
 \distributiontilde^{\bullet}_{\nbigx/\cnum_{\lambda}}[2d_X+2]
\otimes_{f^{-1}\nbigo_{\nbigy}}
f^{-1}\nbigmtilde_{Y1}
 \Bigr)
 \\
\simeq
 \lambda^{d_X}
 f_!\Bigl(
  \distributiontilde^{\bullet}_{\nbigx/\cnum_{\lambda}}[2d_X+2]
 \otimes_{f^{-1}\nbigo_{\nbigy}}
 f^{-1}(\nbigrtilde_Y\otimes\Omegatilde_{\nbigy}^{-1})
 \Bigr)
\otimes_{\nbigo_{\nbigy}}
 \bigl(
 \nbigrtilde_Y\otimes\Omegatilde_{\nbigy}^{-1}
 \bigr)
\\
\simeq
 \lambda^{d_X}\ftilde_{\dagger}(\nbigo_{\nbigx}[d_X+1])
 \otimes_{\nbigo_{\nbigy}}
 \bigl(\nbigrtilde_X\otimes\Omegatilde_{\nbigx}^{-1}\bigr)
\stackrel{\atilde_2}{\lrarr}
 \lambda^{d_Y}
 \bigl(\nbigrtilde_Y\otimes\Omegatilde_{\nbigy}^{-1}\bigr)[d_Y+1]
\end{multline}
Here, $\atilde_1$ is induced by $\Psitilde_1^{-1}$,
and $\atilde_2$ is induced by the trace morphism
(see \S\ref{subsection;14.12.22.130}).

Let $\nbigm$ be a coherent $\nbigrtilde_X$-module
which is good relative to $f$.
Then, we have the natural isomorphism
\begin{equation}
\label{eq;14.12.22.131}
 \ftilde_{\dagger}\DDDtilde_X\nbigm\lrarr
 \DDDtilde_Y \ftilde_{\dagger}\nbigm
\end{equation}
in the derived category of $\nbigrtilde_Y$-modules,
obtained as the composite of the following morphisms:
\begin{multline}
\label{eq;14.12.23.40}
\ftilde_{\dagger}\DDDtilde_X\nbigm
\simeq
 Rf_!\nrhom_{\nbigrtilde_X}\Bigl(
 \nbigm,
 \lambda^{d_X}\nbigo_{\nbigx}
 \otimes_{f^{-1}\nbigo_{\nbigy}}
 f^{-1}(\nbigrtilde_Y\otimes\Omegatilde_{\nbigy}^{-1})[d_X+1]
 \Bigr)
 \\
\lrarr
 Rf_!\nrhom_{f^{-1}\nbigrtilde_Y}\Bigl(
 \nbigrtilde_{Y\larr X}\otimes^L_{\nbigrtilde_X}\nbigm,
 \nbigrtilde_{Y\larr X}\otimes^L_{\nbigrtilde_X}
 \bigl(
 \lambda^{d_X}\nbigo_{\nbigx}
 \otimes_{f^{-1}\nbigo_{\nbigy}}
 f^{-1}(\nbigrtilde_Y\otimes\Omegatilde_{\nbigy}^{-1})
 \bigr)[d_X+1]
 \Bigr)
 \\
\lrarr
 \nrhom_{\nbigrtilde_Y}\Bigl(
 \ftilde_{\dagger}\nbigm,
 \lambda^{d_X}\ftilde_{\dagger}\Bigl(
 \nbigo_{\nbigx}\otimes_{f^{-1}\nbigo_{\nbigy}}
 f^{-1}(\nbigrtilde_{Y}\otimes\Omegatilde_{\nbigy}^{-1})
 \Bigr)[d_X+1]
 \Bigr)
 \\
\lrarr
 \nrhom_{\nbigrtilde_Y}\Bigl(
 \ftilde_{\dagger}\nbigm,
 \lambda^{d_Y}\nbigrtilde_Y\otimes\Omegatilde_{\nbigy}^{-1}[d_Y+1]
 \Bigr)
\simeq
 \DDDtilde_Y \ftilde_{\dagger}\nbigm.
\end{multline}

\subsubsection{Comparison}

Let $f:X\lrarr Y$ be a proper morphism of
complex manifolds as above.
We shall prove the following proposition
in \S\ref{subsection;14.12.23.120}--\ref{subsection;14.12.24.20}.

\begin{prop}
\label{prop;14.12.23.210}
Let $\nbigm$ be a coherent $\nbigrtilde_X$-module
which is good relative to $f$.
The following diagram is commutative:
\begin{equation}
\begin{CD}
 \ftilde_{\dagger}\DDDtilde_X \nbigm
 @>{\simeq}>>
 \DDDtilde_Y \ftilde_{\dagger}\nbigm
 \\ @V{\simeq}VV @V{\simeq}VV \\
 f_{\dagger}\DDD_X\nbigm
 @>{\simeq}>>
 \DDD_Y f_{\dagger}\nbigm
\end{CD}
\end{equation}
Here, 
the horizontal arrows are given in 
{\rm\S\ref{subsection;14.12.22.110}}
and {\rm\S\ref{subsection;14.12.22.140}},
and the vertical arrows are given 
by the isomorphisms
in Corollary {\rm\ref{cor;14.12.22.141}}
and Proposition {\rm\ref{prop;14.12.17.1}}.
\end{prop}

\subsubsection{Step 1}
\label{subsection;14.12.23.120}

We construct a natural transform
$\ftilde_{\dagger}\DDDtilde_X
\lrarr
 \DDDtilde_Y f_{\dagger}$
modifying the constructions in 
\S\ref{subsection;14.12.22.110}
and \S\ref{subsection;14.12.22.140}.

We set 
$\nbigm'_{Y0}:=
 (\nbigr_Y\otimes\Omegabar_Y^{-1})
 \lefttop{\ell}\otimes^{\ell}_{\nbigo_{\nbigy}}
 \bigl(
 \nbigrtilde_Y\otimes\Omegatilde_{\nbigy}^{-1}
 \bigr)$.
The tensor product is taken as in the case of $\nbigm_{Y0}$.
It is equipped with the $\nbigrtilde_X$-action $r_2$
induced by the $\nbigrtilde_X$-action $r$
on the second factor.
It is equipped with the $\nbigrtilde_X$-action
$\ell\otimes\ell$ induced by the $\nbigrtilde_X$-actions
$\ell$ on the first and second factors.
It is equipped with the $\nbigr_X$-action $r_1$
induced by the $\nbigr_X$-action $r$
on the first factor.
Together with the action of $\lambda^2\del_{\lambda}$
by $\ell\otimes\ell$,
it is extended to an $\nbigrtilde_X$-action
which is also denoted by $r_1$.

We set 
$\nbigm'_{Y1}:=
 (\nbigr_Y\otimes\Omegabar_Y^{-1})
 \lefttop{r}\otimes^{\ell}_{\nbigo_{\nbigy}}
 \bigl(
 \nbigrtilde_Y\otimes\Omegatilde_{\nbigy}^{-1}
 \bigr)$.
The tensor product is taken as in the case of $\nbigm_{Y1}$.
It is equipped with the $\nbigrtilde_X$-action $r_2$
induced by the $\nbigrtilde_X$-action $r$
on the second factor.
It is equipped with the $\nbigrtilde_X$-action
$r\otimes\ell$ induced by the $\nbigrtilde_X$-actions
$r$ and $\ell$ on the first and second factors,
respectively.
It is equipped with the $\nbigr_X$-action $\ell_1$
induced by the $\nbigr_X$-action $\ell$
on the first factor.
Together with the action of $\lambda^2\del_{\lambda}$
by $r\otimes\ell$,
it is extended to an $\nbigrtilde_X$-action
which is also denoted by $\ell_1$.

We have the isomorphism
$\Psi'_1:\nbigm'_{1Y}\simeq\nbigm'_{0Y}$
determined by the conditions
(i) $\Psi'_1\circ\ell_1=(\ell\otimes\ell)\circ\Psi_1$
 and $\Psi'_1\circ r_2=r_2\circ\Psi_1$,
(ii) $\Psi'_1$ induces the identity on 
 $\Omegabar_{Y}^{-1}\otimes\Omegatilde_{\nbigy}^{-1}$.
It also satisfies the condition
$\Psi'_1\circ (r\otimes\ell)=r_1\circ\Psi'_1$,
which can be checked by the argument in 
Lemma \ref{lem;14.12.24.1}.

Then, we obtain the morphism
\begin{equation}
\label{eq;14.12.23.20}
 \lambda^{d_X}f_{\dagger}\Bigl(
 \nbigo_{\nbigx}\otimes_{f^{-1}\nbigo_{\nbigy}}
 f^{-1}\bigl(
 \nbigrtilde_Y\otimes\Omegatilde_{\nbigy}^{-1}
 \bigr)
 \Bigr)[d_X+1]
\lrarr
 \lambda^{d_Y}\nbigrtilde_Y\otimes
 \Omegatilde_{\nbigy}^{-1}[d_Y+1]
\end{equation}
in the derived category of 
$\nbigrtilde_Y\otimes\nbigrtilde_Y$-complexes,
obtained as the composite of the following:
\begin{multline}
 \lambda^{d_X}f_{\dagger}\Bigl(
 \nbigo_{\nbigx}\otimes_{f^{-1}\nbigo_{\nbigy}}
 f^{-1}(\nbigrtilde_Y\otimes\Omegatilde_{\nbigy}^{-1})
 \Bigr)[d_X+1]
 \\
\simeq
\lambda^{d_X}
 f_{!}\Bigl(
\Bigl(
 \bigl(
 \distributionbar^{\bullet}_X[2d_X+1]
 \otimes \nbigr_X
 \bigr)
\otimes_{f^{-1}\nbigo_{\nbigy}}
 f^{-1}(\nbigr_Y\otimes\Omegabar_Y^{-1})
\Bigr)\otimes_{\nbigr_X}
 \bigl(
 \nbigo_{\nbigx}
 \otimes_{f^{-1}\nbigo_{\nbigy}}
 f^{-1}\bigl(\nbigrtilde_Y\otimes
 \Omegatilde_{\nbigy}^{-1}\bigr)
 \bigr)
 \Bigr)
\\
\simeq
 \lambda^{d_X}
 f_!\Bigl(
 \distributionbar^{\bullet}_X[2d_X+1]
\otimes_{f^{-1}\nbigo_{\nbigy}}
f^{-1}\nbigm'_{Y0}
 \Bigr)
\stackrel{a'_1}{\simeq}
 \lambda^{d_X}
 f_!\Bigl(
 \distributionbar^{\bullet}_X[2d_X+1]
\otimes_{f^{-1}\nbigo_{\nbigy}}
f^{-1}\nbigm'_{Y1}
 \Bigr)
 \\
\simeq
 \lambda^{d_X}
 f_!\Bigl(
  \distributionbar^{\bullet}_X[2d_X+1]
 \otimes_{f^{-1}\nbigo_{\nbigy}}
 f^{-1}(\nbigr_Y\otimes\Omegabar_{Y}^{-1})
 \Bigr)
\otimes_{\nbigo_{\nbigy}}
 \bigl(
 \nbigrtilde_Y\otimes\Omegatilde_{\nbigy}^{-1}
 \bigr)
\\
\simeq
 \lambda^{d_X}f_{\dagger}(\nbigo_{\nbigx}[d_X+1])
 \otimes_{\nbigo_{\nbigy}}
 \bigl(\nbigrtilde_Y\otimes\Omegatilde_{\nbigy}^{-1}\bigr)
\stackrel{a'_2}{\lrarr}
 \lambda^{d_Y}
 \bigl(\nbigrtilde_Y\otimes\Omegatilde_{\nbigy}^{-1}\bigr)[d_Y+1]
\end{multline}
Here $a_1'$ is induced by $\Psi_1'$,
and $a_2'$ is induced by the trace morphism
in \S\ref{subsection;14.12.22.130}.
Then, for any $\nbigrtilde_X$-module $\nbigm$
which is good relative to $f$ and coherent over $\nbigr_X$,
we obtain the morphism
\begin{equation}
 \label{eq;14.12.23.2}
 \ftilde_{\dagger}\DDDtilde_X\nbigm
\lrarr
 \DDDtilde_Y f_{\dagger}\nbigm
\end{equation}
as the composite of the following morphisms:
\begin{multline}
\label{eq;14.12.23.41}
\ftilde_{\dagger}\DDDtilde_X\nbigm
\simeq
 Rf_!\nrhom_{\nbigrtilde_X}\Bigl(
 \nbigm,
 \lambda^{d_X}\nbigo_{\nbigx}
 \otimes_{f^{-1}\nbigo_{\nbigy}}
 f^{-1}(\nbigrtilde_Y\otimes\Omegatilde_{\nbigy}^{-1})[d_X+1]
 \Bigr)
 \\
\lrarr
 Rf_!\nrhom_{f^{-1}\nbigrtilde_Y}\Bigl(
 \nbigr_{Y\larr X}\otimes^L_{\nbigr_X}\nbigm,
 \nbigr_{Y\larr X}\otimes^L_{\nbigr_X}
 \bigl(
 \lambda^{d_X}\nbigo_{\nbigx}
 \otimes_{f^{-1}\nbigo_{\nbigy}}
 f^{-1}(\nbigrtilde_Y\otimes\Omegatilde_{\nbigy}^{-1})
 \bigr)[d_X+1]
 \Bigr)
 \\
\lrarr
 \nrhom_{\nbigrtilde_Y}\Bigl(
 f_{\dagger}\nbigm,
 \lambda^{d_X}f_{\dagger}\Bigl(
 \nbigo_{\nbigx}\otimes_{f^{-1}\nbigo_{\nbigy}}
 f^{-1}(\nbigrtilde_{Y}\otimes\Omegatilde_{\nbigy}^{-1})
 \Bigr)[d_X+1]
 \Bigr)
 \\
\lrarr
 \nrhom_{\nbigrtilde_Y}\Bigl(
 f_{\dagger}\nbigm,
 \lambda^{d_Y}\nbigrtilde_Y\otimes\Omegatilde_{\nbigy}^{-1}[d_Y+1]
 \Bigr)
\simeq
 \DDDtilde_Y f_{\dagger}\nbigm
\end{multline}

\begin{lem}
\label{lem;14.12.23.34}
The following diagram is commutative:
\begin{equation}
\label{eq;14.12.23.33}
 \begin{CD}
\lambda^{d_X}\ftilde_{\dagger}\Bigl(
 \nbigo_{\nbigx}\otimes_{f^{-1}\nbigo_{\nbigy}}
 f^{-1}\bigl(
 \nbigrtilde_Y\otimes\Omegatilde_{\nbigy}^{-1}
 \bigr)
 \Bigr)[d_X+1]
 @>{c_1}>>
 \lambda^{d_Y}\nbigrtilde_Y\otimes
 \Omegatilde_{\nbigy}^{-1}[d_Y+1]
 \\
 @VVV @V{=}VV \\
 \lambda^{d_X}f_{\dagger}\Bigl(
 \nbigo_{\nbigx}\otimes_{f^{-1}\nbigo_{\nbigy}}
 f^{-1}\bigl(
 \nbigrtilde_Y\otimes\Omegatilde_{\nbigy}^{-1}
 \bigr)
 \Bigr)[d_X+1]
 @>{c_2}>>
 \lambda^{d_Y}\nbigrtilde_Y\otimes
 \Omegatilde_{\nbigy}^{-1}[d_Y+1]
 \end{CD}
\end{equation}
Here, $c_1$ and $c_2$
are given by {\rm(\ref{eq;14.12.22.120})}
and {\rm(\ref{eq;14.12.23.20})},
and the left vertical arrow is 
given in Proposition {\rm\ref{prop;14.12.17.1}}.
\end{lem}
\pf
We can regard $\nbigm'_{Y0}$
as an $\nbiga_{Y}$-module 
by $\ell\otimes\ell$ and $r_1$ with
the action of $\lambda^2\del_{\lambda}$.
We can also regard $\nbigm'_{Y1}$
as an $\nbiga_{Y}$-module 
by $r\otimes\ell$ and $\ell_1$ with
the action of $\lambda^2\del_{\lambda}$.
As mentioned in Example \ref{example;14.12.22.20},
we have the following isomorphisms:
\[
 \bigl(
 \nbigr_Y\otimes\Omegabar_Y^{-1}
 \langle\lambda^2\del_{\lambda}\rangle,
 \ell\blacktriangledown,
 r\triangledown
 \bigr)
\simeq
 \bigl(
\nbigrtilde_Y\otimes
  \lambda^{-2}\Omegatilde_{\nbigy}^{-1},
 \ell, r\bigr)
\]
It induces the following isomorphisms:
\begin{equation}
 \label{eq;14.12.23.30}
\Bigl(
 \nbigm'_{Y0}\langle\lambda^2\del_{\lambda}\rangle,
 (\ell\otimes\ell)\blacktriangledown,
 r_1\triangledown
\Bigr)
\simeq
 \Bigl(
 \lambda^{-2}\nbigmtilde_{Y0},
 \ell\otimes\ell,
 r_1
\Bigr)
\end{equation}
\begin{equation}
\label{eq;14.12.23.31}
 \Bigl(
 \nbigm'_{Y1}\langle\lambda^2\del_{\lambda}\rangle,
 (r\otimes\ell)\triangledown,
 \ell_1\blacktriangledown
\Bigr)
\simeq
 \Bigl(
 \lambda^{-2}\nbigmtilde_{Y1},
 r\otimes\ell,
 \ell_1
\Bigr)
\end{equation}
The isomorphisms are compatible with the $\nbigrtilde_X$-actions
$r_2$.
We have the following commutative diagram:
\[
 \begin{CD}
  \nbigm'_{Y1}\langle\lambda^2\del_{\lambda}\rangle
 @>>>
 \lambda^{-2}\nbigmtilde_{Y1}
 \\ 
 @VVV @VVV\\
 \nbigm'_{Y0}\langle\lambda^2\del_{\lambda}\rangle
 @>>>
 \lambda^{-2}\nbigmtilde_{Y0}
 \end{CD}
\]
Here, the vertical arrows are induced by
$\Psi'_1$ and $\Psitilde_1$,
and the horizontal arrows are as in 
(\ref{eq;14.12.23.30})
and (\ref{eq;14.12.23.31}).
Hence, the following diagram is commutative:
\[
 \begin{CD}
 \lambda^{d_X}
 f_!\Bigl(
 \distributiontilde^{\bullet}_{\nbigx/\cnum_{\lambda}}[2d_X+2]
\otimes_{f^{-1}\nbigo_{\nbigy}}
f^{-1}\nbigmtilde_{Y0}
 \Bigr)
 @>{\atilde_1}>{\simeq}>
 \lambda^{d_X}
 f_!\Bigl(
 \distributiontilde^{\bullet}_{\nbigx/\cnum_{\lambda}}[2d_X+2]
\otimes_{f^{-1}\nbigo_{\nbigy}}
f^{-1}\nbigmtilde_{Y1}
 \Bigr)
 \\
 @VVV @VVV\\
 \lambda^{d_X}
 f_!\Bigl(
 \distributionbar^{\bullet}_X[2d_X+1]
\otimes_{f^{-1}\nbigo_{\nbigy}}
f^{-1}\nbigm'_{Y0}
 \Bigr)
 @>{a'_1}>{\simeq}>
 \lambda^{d_X}
 f_!\Bigl(
 \distributionbar^{\bullet}_X[2d_X+1]
\otimes_{f^{-1}\nbigo_{\nbigy}}
f^{-1}\nbigm'_{Y1}
 \Bigr)
 \end{CD}
\]
We also have the compatibility of the trace morphism
in Lemma \ref{lem;14.12.23.32}.
Then, we obtain the claim of Lemma \ref{lem;14.12.23.34}
by construction of the morphisms.
\hfill\qed

\vspace{.1in}

The natural transform $\ftilde_{\dagger}\lrarr f_{\dagger}$
in Corollary \ref{cor;14.12.22.141}
induces
$\DDDtilde_Y f_{\dagger}\lrarr
\DDDtilde_Y \ftilde_{\dagger}$.

\begin{lem}
The composite of 
{\rm(\ref{eq;14.12.23.2})}
and $\DDDtilde_Y f_{\dagger}\lrarr \DDDtilde_Y \ftilde_{\dagger}$
is equal to 
{\rm(\ref{eq;14.12.22.131})}.
\end{lem}
\pf
We obtain the claim
by comparing the constructions 
(\ref{eq;14.12.23.40})
and (\ref{eq;14.12.23.41})
and by using Lemma \ref{lem;14.12.23.34}.
\hfill\qed

\subsubsection{Step 2}
\label{subsection;14.12.24.20}

We consider the following diagram:
\begin{equation}
\label{eq;14.12.24.3}
\begin{CD}
 \ftilde_{\dagger}
 \DDDtilde_X
@>{\alpha_1}>>
 \ftilde_{\dagger}\DDDtilde_X S_{\lambda}^{\blacktriangledown}
@>{\alpha_2}>>
 f_{\dagger}\DDD_X
 \\
 @V{\beta_1}VV @V{\beta_2}VV @V{\beta_3}VV 
\\
\DDDtilde_Y f_{\dagger}
 @>{\gamma_1}>>
 \DDDtilde_Y f_{\dagger}S_{\lambda}^{\blacktriangledown}
 @>{\gamma_2}>>
 \DDD_Y f_{\dagger}
\end{CD}
\end{equation}
The morphisms $\alpha_1$ and $\gamma_1$ 
are the natural ones.
The morphisms $\beta_i$ $(i=1,2)$ are given in 
\S\ref{subsection;14.12.23.120}.
The morphism $\beta_3$ is given in 
\S\ref{subsection;14.12.23.210}.
The morphism
$\alpha_2$ is 
induced by
$\ftilde_{\dagger}\lrarr f_{\dagger}$
in Corollary \ref{cor;14.12.22.141},
and
$\DDDtilde_X S_{\lambda}^{\blacktriangledown}
\simeq
 S_{\lambda}^{\triangledown}\DDD_X
\lrarr \DDD_X$
in \S\ref{subsection;14.12.23.200}.
The morphism $\gamma_2$
is induced by
$f_{\dagger}S_{\lambda}^{\blacktriangledown}
\simeq
 S_{\lambda}^{\blacktriangledown}f_{\dagger}$
in \S\ref{subsection;14.12.23.210},
and $\DDDtilde_Y S_{\lambda}^{\blacktriangledown}
\simeq S_{\lambda}^{\triangledown}\DDD_Y
\lrarr \DDD_Y$
in \S\ref{subsection;14.12.23.200}.
The morphism $\beta_3$
is given in \S\ref{subsection;14.12.22.110}.
Clearly, the left square is commutative.
For the proof of Proposition \ref{prop;14.12.23.210},
it is enough to prove the commutativity of
the right square.

We have the following expression:
\begin{equation}
\label{eq;14.12.24.10}
 \ftilde_{\dagger}\DDDtilde
 S_{\lambda}^{\blacktriangledown}\nbigm
\simeq 
\lambda^{d_X}Rf_!
\Bigl[
\Bigl(
 \Omegatilde_{\nbigx}^{\bullet}[d_X+1]
 \otimes_{f^{-1}\nbigo_{\nbigy}}
 f^{-1}(\nbigrtilde_{Y}\otimes\Omegatilde_{\nbigy}^{-1})
\Bigr)
\otimes_{\nbigo_{\nbigx}}
\nrhom_{\nbigrtilde_X}\bigl(
S_{\lambda}^{\blacktriangledown}\nbigm,
\nbigrtilde_X\otimes\Omegatilde_{\nbigx}^{-1}
\bigr)[d_X+1]
\Bigr]
\end{equation}
It is rewritten as follows:
\begin{equation}
\label{eq;14.12.23.200}
 \lambda^{d_X}
 Rf_!\Bigl[
 \Bigl(
 \Omegatilde_{\nbigx}^{\bullet}[d_X+1]
 \otimes_{f^{-1}\nbigo_{\nbigy}}
 f^{-1}(\nbigrtilde_{Y}\otimes\Omegatilde_{\nbigy}^{-1})
\Bigr)
\otimes_{\nbigo_{\nbigx}}
 S_{\lambda}^{\triangledown}
\nrhom_{\nbigr_X}\Bigl(
 \nbigm,
\nbigr_X\otimes\Omegabar_{X}^{-1}
\bigr)[d_X]
\Bigr]
\end{equation}
We have the natural morphism from
(\ref{eq;14.12.23.200})
to the following:
\begin{multline}
\label{eq;14.12.23.201}
\lambda^{d_X}Rf_!
\Bigl[
 \bigl(
 \Omegabar^{\bullet}_{X}[d_X]
 \otimes_{f^{-1}\nbigo_{\nbigy}}
  f^{-1}(\nbigr_Y\otimes\Omegabar_Y^{-1})
\bigr)
\otimes_{\nbigo_{\nbigx}}
\nrhom_{\nbigr_X}\Bigl(
 \nbigm,
\nbigr_X\otimes\Omegabar_X^{-1} 
\Bigr)[d_X]
\Bigr]
\simeq
 \\
\lambda^{d_X}Rf_!\nrhom_{\nbigr_X}\Bigl(
 \nbigm,
 \nbigo_{\nbigx}
 \otimes_{f^{-1}\nbigo_{\nbigy}}
  f^{-1}(\nbigr_Y\otimes\Omegabar_Y^{-1})
\Bigr)[d_X]
\end{multline}
Indeed, we have 
\begin{multline}
\label{eq;14.12.24.100}
\lambda^{d_X}
Rf_!\Bigl[
 \Bigl(
 \Omegatilde_{\nbigx}^{\bullet}[d_X+1]
 \otimes_{f^{-1}\nbigo_{\nbigy}}
 f^{-1}(\nbigrtilde_{Y}\otimes\Omegatilde_{\nbigy}^{-1})
\Bigr)
\otimes_{\nbigo_{\nbigx}}
 S_{\lambda}^{\triangledown}
\nrhom_{\nbigr_X}\Bigl(
 \nbigm,
\nbigr_X\otimes\Omegabar_{X}^{-1}
\bigr)[d_X]
\Bigr]
\lrarr
 \\
\lambda^{d_X}
 Rf_!\Bigl[
 \Bigl(
 \Omegatilde_{\nbigx}^{\bullet}[d_X+1]
 \otimes_{f^{-1}\nbigo_{\nbigy}}
 f^{-1}(\nbigrtilde_{Y}\otimes\Omegatilde_{\nbigy}^{-1})
\Bigr)
\otimes_{\nbigo_{\nbigx}}
\nrhom_{\nbigr_X}\Bigl(
 \nbigm,
\nbigr_X\otimes\Omegabar_{X}^{-1}
\bigr)[d_X]
\Bigr]
\\
\lambda^{d_X}Rf_!
\Bigl[
 \bigl(
 \Omegabar^{\bullet}_{X}[d_X]
 \otimes_{f^{-1}\nbigo_{\nbigy}}
  f^{-1}(\nbigr_Y\otimes\Omegabar_Y^{-1})
\bigr)
\otimes_{\nbigo_{\nbigx}}
\nrhom_{\nbigr_X}\Bigl(
 \nbigm,
\nbigr_X\otimes\Omegabar_X^{-1} 
\Bigr)[d_X]
\Bigr]
\end{multline}
We also have
\begin{multline}
\label{eq;14.12.24.101}
\lambda^{d_X}
Rf_!\Bigl[
 \Bigl(
 \Omegatilde_{\nbigx}^{\bullet}[d_X+1]
 \otimes_{f^{-1}\nbigo_{\nbigy}}
 f^{-1}(\nbigrtilde_{Y}\otimes\Omegatilde_{\nbigy}^{-1})
\Bigr)
\otimes_{\nbigo_{\nbigx}}
 S_{\lambda}^{\triangledown}
\nrhom_{\nbigr_X}\Bigl(
 \nbigm,
\nbigr_X\otimes\Omegabar_{X}^{-1}
\bigr)[d_X]
\Bigr]
\lrarr
 \\
\lambda^{d_X}Rf_!
\Bigl[
 \bigl(
 \Omegabar^{\bullet}_{X}[d_X]
 \otimes_{f^{-1}\nbigo_{\nbigy}}
  f^{-1}(\nbigr_Y\otimes\Omegabar_Y^{-1})
\bigr)
\otimes_{\nbigo_{\nbigx}}
 S_{\lambda}^{\triangledown}
\nrhom_{\nbigr_X}\Bigl(
 \nbigm,
\nbigr_X\otimes\Omegabar_X^{-1} 
\Bigr)[d_X]
\Bigr]
\\
\lambda^{d_X}Rf_!
\Bigl[
 \bigl(
 \Omegabar^{\bullet}_{X}[d_X]
 \otimes_{f^{-1}\nbigo_{\nbigy}}
  f^{-1}(\nbigr_Y\otimes\Omegabar_Y^{-1})
\bigr)
\otimes_{\nbigo_{\nbigx}}
\nrhom_{\nbigr_X}\Bigl(
 \nbigm,
\nbigr_X\otimes\Omegabar_X^{-1} 
\Bigr)[d_X]
\Bigr]
\end{multline}
The composite of the morphisms (\ref{eq;14.12.24.100})
and (\ref{eq;14.12.24.101}) are equal.
We have the following morphisms:
\begin{multline}
 Rf_!\nrhom_{\nbigr_X}\Bigl(
 \nbigm,
 \lambda^{d_X}\nbigo_{\nbigx}[d_X]
 \otimes_{f^{-1}\nbigo_{\nbigy}}
  f^{-1}(\nbigr_Y\otimes\Omegabar_Y^{-1})
\Bigr)
\lrarr
 \\
 \nrhom_{\nbigr_Y}\Bigl(
 f_{\dagger}\nbigm,
 f_{\dagger}\Bigl(
 \lambda^{d_X}\nbigo_{\nbigx}[d_X]
 \otimes_{f^{-1}\nbigo_{\nbigy}}
  f^{-1}(\nbigr_Y\otimes\Omegabar_Y^{-1})
 \Bigr)
\Bigr)
\\
\simeq
 \nrhom_{\nbigr_Y}\Bigl(
 f_{\dagger}\nbigm,
 \lambda^{d_X}
 f_{\dagger}(\nbigo_{\nbigx})[d_X]
 \otimes_{\nbigo_{\nbigy}}
\bigl(
  \nbigr_Y\otimes\Omegabar_Y^{-1}
\bigr)
\Bigr)
 \\
\lrarr 
 \nrhom_{\nbigr_Y}\Bigl(
 f_{\dagger}\nbigm,
 \lambda^{d_Y}
  \nbigr_Y\otimes\Omegabar_Y^{-1}[d_Y]
\Bigr)
\lrarr
 \DDD f_{\dagger}\nbigm
\end{multline}
As the composite of the above morphisms,
we obtain a morphism
$G:\ftilde_{\dagger}\DDDtilde S_{\lambda}^{\blacktriangledown}
 \nbigm
\lrarr
 \DDD f_{\dagger}\nbigm$.
By construction,
we can check that
both $\beta_3\circ\alpha_2$
and $\gamma_2\circ\beta_2$
are equal to $G$.
Indeed, $\beta_3\circ\alpha_2$ is induced by
(\ref{eq;14.12.24.100}),
and $\gamma_2\circ\beta_2$ is induced by
(\ref{eq;14.12.24.101}).
Thus, we obtain the commutativity of the right square of
(\ref{eq;14.12.24.3}),
and the proof of Proposition \ref{prop;14.12.23.210}
is finished.
\hfill\qed

\section{Comparison of some $\nbigrtilde$-modules}
\label{section;14.12.27.1}
\subsection{Preliminary}

\subsubsection{A construction}
\label{subsection;14.7.6.1}

Let $W_t=\cnum[t]\langle\del_t\rangle$ 
be the Weyl algebra in the variable $t$.
We set $W_t^{\loc}:=W_t\otimes_{\cnum[t]}\cnum[t,t^{-1}]$.
Let $X$ be any smooth algebraic variety.
Let $p:\cnum_t\times X\lrarr X$ be the projection.
In the algebraic setting,
we naturally identify
$p_{\ast}\nbigd_{\cnum_t\times X}$
with $\nbigd_X\otimes W_t$.
We naturally identify
algebraic $\nbigd_{\cnum_t\times X}$-modules
and
algebraic $\nbigd_X\otimes W_t$-modules.
We also naturally identify
algebraic 
$\nbigd_{\cnum_t^{\ast}\times X}$-modules
and
algebraic $\nbigd_X\otimes W^{\loc}_t$-modules.

We regard $\cnum_{\tau}$ as the dual space of
$\cnum_t$ by the pairing 
$(\tau,t)\longmapsto\tau t$.
Let $M$ be an algebraic 
$\nbigd_{\cnum_t\times X}$-module.
We have its partial Fourier-Laplace transform
$\FL_X(M)$ on $\cnum_{\tau}\times X$.
As a transform from
$\nbigd_{X}\otimes W_t$-modules to 
$\nbigd_{X}\otimes W_{\tau}$-modules,
it is described as the quotient of
$\del_t-\tau:
 M[\tau]\lrarr M[\tau]$.
Here, the action of 
$\nbigd_X\otimes W_{\tau}$
on $M[\tau]$ is given as follows.
We have $v(m\tau^j)=(vm)\tau^j$
for a vector field $v$ on $X$,
and 
$\del_{\tau}(m\tau^j)=
 jm\tau^{j-1}
-tm\tau^{j}$.

By taking the restriction,
we have the localized partial Fourier transform
$\FL^{\loc}_X(M)$
on $\cnum_{\tau}^{\ast}\times X$.
We set $\lambda:=\tau^{-1}$.
We may naturally regard
$\FL_X^{\loc}(M)$
as a $\nbigd_{X}\otimes W_{\lambda}^{\loc}$-module.
Let us explicit it more explicitly.
We have the action of
$\nbigd_{\cnum_{\lambda}^{\ast}\times X}$
on $M[\lambda,\lambda^{-1}]$
given as follows.
We have $v(m\lambda^j)=(vm)\lambda^j$
for a vector field $v$ on $X$,
and 
$\del_{\lambda}(m\lambda^j)=
 tm\lambda^{j-2}+jm\lambda^{j-1}$.
Then, the $\nbigd_{\cnum_{\lambda}^{\ast}\times X}$-module
$\FL_X^{\loc}(M)$ is described
as the quotient of
\[
 M[\lambda,\lambda^{-1}]
 \stackrel{\del_t-\lambda^{-1}}{\lrarr}
 M[\lambda,\lambda^{-1}].
\]
The natural inclusion $M\lrarr M[\lambda,\lambda^{-1}]$
induces a $\nbigd_X$-homomorphism
$\loc:M\lrarr \FL^{\loc}_X(M)$.

Let $(M,F)$ be an algebraic filtered 
$\nbigd_{\cnum_t\times X}$-module,
i.e., $M$ is an algebraic $\nbigd_X$-module,
and $F$ is a good filtration of $M$.
Following \cite{sabbah-Fourier-Laplace},
we set
\[
G_0\FL_X^{\loc}(M,F):=
 \sum_{j\in\seisuu} \lambda^j\loc(F_jM)
\subset
 \FL_X^{\loc}(M)
\]
on $X$.
It is naturally an algebraic $\nbigr_X$-module.
It is also equipped with
the action of $\lambda^{2}\del_{\lambda}$.
Thus, we obtain an algebraic $\nbigrtilde_X$-module.
Because 
$\del_{t}F_j(M)\subset F_{j+1}(M)$,
we have
$\loc(F_jM)\subset
 \lambda\loc(F_{j+1}M)$,
and hence
$G_0\FL_X^{\loc}(M,F)
=\sum_{j\geq N}
 \lambda^j\loc(F_jM)$
for any $N$.

We have the Rees module
$R(M,F)=\sum F_jM\lambda^j$.
We have the $\nbigr_X$-homomorphism
$R(M,F)\rarr
 \FL_X^{\loc}(M,F)$,
and the image is
$G_0\FL_X^{\loc}(M,F)$.

The quotient of the following is denoted by
$G_0'\FL_X^{\loc}(M,F)$:
\[
 \lambda\,R(M,F)
\stackrel{\del_t-\lambda^{-1}}\lrarr
 R(M,F)
\]
It is naturally an algebraic $\nbigrtilde_X$-module.

\begin{lem}
\label{lem;14.12.11.2}
The image of the induced morphism
$G'_0\FL_X^{\loc}(M,F)
\lrarr
 \FL_X^{\loc}(M)$
is $G_0\FL_X^{\loc}(M,F)$.
If the multiplication of $\lambda$
on $G_0'\FL_X^{\loc}(M,F)$
is injective,
the natural morphism
$G_0'\FL_X^{\loc}(M,F)
\lrarr
 G_0\FL_X^{\loc}(M,F)$
is an isomorphism
of $\nbigrtilde_X$-modules.
\end{lem}
\pf
The first claim is clear by the construction.
The morphism
$G_0'\FL_X^{\loc}(M,F)
\lrarr
 G_0\FL_X^{\loc}(M,F)$
is an isomorphism after the localization
with respect to $\lambda$.
Then, the second claim follows.
\hfill\qed

\vspace{.1in}
We may interpret the construction $G_0'\FL^{\loc}_X$
in terms of $\nbigr$-modules.
Let $p$ denote the projection of
$\cnum_t\times X$ onto $X$, as above.
We have the algebraic $\nbigr_{\cnum_t\times X}$-module
$\nbigl(-t)=\bigl(
 \nbigo_{\cnum_{\lambda}\times \cnum_t\times X},\,
 d-d(\lambda^{-1}t)
 \bigr)$.
Then, by definition,
$G_0'\FL^{\loc}_X(M,F)$
is isomorphic to
\[
 \lambda\cdot p^0_{\dagger}
 \bigl(
 R(M,F)\otimes\nbigl(-t)
 \bigr)
\simeq
 \lambda\cdot
 \hyperr^1p_{\ast}\Bigl(
 R(M,F)\otimes\nbigl(-t)
 \stackrel{b}{\lrarr}
 R(M,F)\otimes\nbigl(-t)
 \otimes
 \bigl(
 \lambda^{-1}
 \Omega^{1}_{\cnum_{\lambda}\times X\times\cnum_t
 /\cnum_{\lambda}\times X}
 \bigr)
 \Bigr)
\]
Here, $b$ is induced by the meromorphic flat connection 
of $R(M,F)\otimes\nbigl(-t)$.
(See \S\ref{subsection;14.11.6.11} for example.)

\begin{prop}
\label{prop;14.10.16.1}
Suppose that $(M,F)$ is 
a filtered $\nbigd_{\cnum_t\times X}$-module
underlying an algebraic mixed Hodge module.
Then, 
the natural morphism
$G_0'\FL^{\loc}_X(M,F)
\lrarr
 G_0\FL^{\loc}_X(M,F)$
is an isomorphism.
\end{prop}
\pf
By the assumption,
we have a mixed Hodge module on $\proj^1_t\times X$
whose underlying filtered $\nbigd$-module
$(M',F')$ on $\proj^1_t\times X$
satisfies
$(M',F')_{|\cnum_t\times X}=(M,F)$.
As explained in \cite[\S13.5]{Mochizuki-MTM},
the analytification of $R(M',F')$ 
underlies an integrable mixed twistor $\nbigd$-module
$(\nbigt,W)$ on $\proj^1_t\times X$.
Let $\nbigv(-t)=\bigl(\nbigl(-t),\nbigl(-t),C(-t)\bigr)$
denote the $\nbigr_{\proj^1_t\times X}$-triple
as in \S\ref{subsection;14.12.11.1}.
As explained in \cite[\S11.3.2]{Mochizuki-MTM},
we have a mixed twistor $\nbigd$-module
$\bigl((\nbigt,W)\otimes\nbigv(-t)\bigr)[\ast (t)_{\infty}]$
on $\proj^1_t\times X$.

\begin{lem}
\label{lem;14.12.11.10}
The underlying $\nbigr_{\proj^1_t\times X}$-module
of $\bigl((\nbigt,W)\otimes\nbigv(-t)\bigr)[\ast (t)_{\infty}]$
is the analytification of $R(M,F)\otimes \nbigl(-t)$.
\end{lem}
\pf
The claim is clear on $\cnum_t\times X$.
Let $Q$ be any point of $X$.
It is enough to check the claim locally around
$(\infty,Q)$.
Let $\nbigm$ and $\nbigmtilde$
be the  $\nbigr_{\proj^1_t\times X}$-modules
underlying $\nbigt$
and 
$\bigl((\nbigt,W)\otimes\nbigv(-t)\bigr)[\ast (t)_{\infty}]$,
respectively.
We have 
$\nbigmtilde(\ast (t)_{\infty})
=\nbigm\otimes\nbigl(-t)$.

Let $\kappa:=t^{-1}$.
We consider 
$V\nbigr_{\cnum_{\kappa}\times X}
\subset 
 \nbigr_{\cnum_{\kappa}\times X}$
which is generated by
$\lambda
 \Theta_{\cnum_{\kappa}\times X}(\log \kappa)$
over $\nbigo_{\cnum_{\lambda}\times\cnum_{\kappa}\times X}$.
We have the $V$-filtration $V(\nbigm)$
of $\nbigm$ along $\kappa$.
Then,
$V_{c}(\nbigm)$ is
$V\nbigr_{\cnum_{\kappa}\times X}$-coherent
for any $c\in\real$.
Let $q:\cnum_{\kappa}\times X\lrarr X$ be
the projection.
Note that
$V_{c}(\nbigm)$ are $q^{\ast}\nbigr_X$-coherent
because $\nbigm$ is regular.
We take sections $s_1,\ldots,s_m$
of $V_c(\nbigm)$
which generates $V_c(\nbigm)$
over $q^{\ast}\nbigr_X$.
Let $\upsilon$ be the generator of $\nbigl(-t)$
such that $\nabla\upsilon=\upsilon d(-\lambda^{-1}t)$.

Let $\nbigm_1$ be the $\nbigr_{\cnum_{\kappa}\times X}$-submodule
of $\nbigm\otimes\nbigl(-t)$
generated by 
$s_i\otimes\upsilon$ $(i=1,\ldots,m)$.
Let us observe
$\nbigm_1=\nbigm\otimes\nbigl(-t)$.
It is enough to prove that 
$\kappa^{-n}(V_c\nbigm\otimes \upsilon)
\subset
 \nbigm_1$ for any $n$.
We have 
$V_c\nbigm\otimes \upsilon
\subset
 \nbigm_1$
because $V_c(\nbigm)$
is generated by $s_1,\ldots,s_m$
over $q^{\ast}\nbigr_X$.
Suppose 
$\kappa^{-n}(V_c\nbigm\otimes \upsilon)
\subset
 \nbigm_1$.
Let $f$ be any section of 
$V_c\nbigm$.
We have 
\[
\kappa\del_{\kappa}
\bigl(
 \kappa^{-n}f\otimes\upsilon
\bigr)
=-n\kappa^{-n}f\otimes\upsilon
+\kappa^{-n}(\kappa\del_{\kappa}f)\otimes\upsilon
-\kappa^{-n-1}f\otimes\upsilon. 
\]
Hence, we obtain
$\kappa^{-n-1}f\otimes\upsilon
\in\nbigm_1$.

Suppose $c<0$.
For a sufficiently large $N$,
$t^Ns_i\otimes \upsilon$ $(i=1,\ldots,m)$
are sections of $\nbigmtilde$.
We also have that 
$t^Ns_i$ generates
$V_{c-N}\nbigm$.
Hence, we obtain that 
$t^Ns_i\otimes \upsilon$ $(i=1,\ldots,m)$
generates $\nbigm\otimes\nbigl(-t)$
over $\nbigr_X$.
It implies that
$\nbigmtilde=\nbigm\otimes\nbigl(-t)$.
\hfill\qed

\vspace{.1in}

Hence,
the analytification of
$p^0_{\dagger}
 \bigl(
 R(M,F)\otimes\nbigl(-t)
 \bigr)$
underlies a mixed twistor $\nbigd$-module.
It implies that
the multiplication of $\lambda$
on 
$p^0_{\dagger}
 \bigl(
 R(M,F)\otimes\nbigl(-t)
 \bigr)$
is injective.
Then, the claim follows from the Lemma \ref{lem;14.12.11.2}.
\hfill\qed

\subsubsection{Comparison}
\label{subsection;14.7.16.4}

Let $Y$ be a smooth projective variety.
Let $X$ be any smooth algebraic variety.
Let $\pi_1:Y\times \cnum_{t}\times X
\lrarr \cnum_{t}\times X$
and 
$\pi_2:Y\times\cnum_{t}\times X
\lrarr X$ be the projections.

Let $(M,F)$ be a filtered $\nbigd$-module
on $Y\times\cnum_{t}\times X$
which underlies a mixed Hodge module.
We obtain filtered $\nbigd$-modules
$\pi^j_{1\dagger}(M,F)$
on $\cnum_{t}\times X$,
which underlie mixed Hodge modules.
Applying the construction in \S\ref{subsection;14.7.6.1},
we obtain $\nbigr_X$-modules
$G_0\FL^{\loc}_X\pi^j_{1\dagger}(M,F)$.

\vspace{.1in}

Let $R(M,F)$ be the Rees module of $(M,F)$.
It is naturally an algebraic
$\nbigrtilde_{Y\times\cnum_{t}\times X}$-module.
We have the $\nbigr_{Y\times\cnum_t\times X}$-module
$\nbigl(-t)$.
We obtain an algebraic
$\nbigrtilde_{Y\times\cnum_t\times X}$-module
$R(M,F)\otimes \nbigl(-t)$.
It underlies a mixed twistor $\nbigd$-module
(Lemma \ref{lem;14.12.11.10}).
We obtain $\nbigr_{X}$-modules
$\pi^j_{2\dagger}\bigl(
 R(M,F)\otimes\nbigl(-t)
\bigr)$.

\begin{prop}
\label{prop;14.7.9.30}
We have natural isomorphisms
$\lambda\cdot\pi^j_{2\dagger}\bigl(
 R(M,F)\otimes\nbigl(-t)
 \bigr)
\simeq
 G_0\FL^{\loc}_X\pi^j_{1\dagger}(M,F)$.
\end{prop}
\pf
Note that we have
$\pi_{1\dagger}^j\bigl(R(M,F)\bigr)
\simeq
 R\bigl(\pi_{1\dagger}^j(M,F)\bigr)$
by the theory of mixed Hodge modules.
Let $p_X:\cnum_t\times X\lrarr X$ denote the projection.
We have $p_X\circ\pi_1=\pi_2$.
Let $L(-t)$ denote the $\nbigd$-module associated to $-t$
on $\proj^1_t\times X$
as in \S\ref{subsection;14.10.14.10}.
Let us recall the following lemma.
\begin{lem}
\label{lem;14.11.5.1}
For regular holonomic $\nbigd$-modules $M$
on $\cnum_t\times X$,
we have
$p_{X\dagger}^jM\otimes L(-t)=0$ for $j\neq 0$.
\end{lem}
\pf
It is enough to consider the case $j=-1$.
Let $s$ be any section of $M$ 
such that $\del_ts-s=0$,
and let us check that $s=0$.
We have the $V$-filtration of $V_{\bullet}M$ 
along $u=t^{-1}$.
Suppose $s\neq 0$ around 
$(\infty,P)\in\proj^1\times X$.
Because $M$ is regular along $u$,
we have 
$a_0:=\min\{a\,|\,s\in V_{a}M\}$.
Because $\del_ts\in V_{<a_0}$,
we also have $s\in V_{<a_0}$
which contradicts with the choice of $a_0$.
Hence, we have $s=0$ around $(\infty,P)$.
Let $M_1\subset M$ be the $\nbigd$-submodule 
generated by $s$.
We have $\Supp(M_1)\cap \{\infty\}\times X=\emptyset$.
Let $P\in X$.
By shrinking $X$ appropriately,
we take an algebraic function $f$ on $X$
such that $f(P)=0$ and $df$ is nowhere vanishing.
We have the $V$-filtrations $V_f$ of $M_1$ and $M_1\otimes L(-t)$
along $f$.
It is easy to see that
$\Gr^{V_f}(M\otimes L(-t))
 \simeq
 \Gr^{V_f}(M)\otimes L(-t)$.
If $s$ is non-zero around $(t,P)$,
there exists $b$ such that 
$s$ gives a non-zero section $[s]$ of $\Gr_b^{V_f}(M)$
satisfying $\del_t[s]-[s]=0$.
Hence, we can reduce the claim to the case $\dim X=0$,
where it can be checked easily.
Thus, we obtain Lemma \ref{lem;14.11.5.1}.
\hfill\qed

\vspace{.1in}

By the lemma,
we have
$p^0_{X\dagger}
 \pi_{1\dagger}^j\bigl(
 R(M,F)\otimes \nbigl(-t)
 \bigr)
\simeq
 \pi_{2\dagger}^j\bigl(
 R(M,F)\otimes\nbigl(-t)
 \bigr)$.
Then, the claim of Proposition \ref{prop;14.7.9.30}
follows from Proposition \ref{prop;14.10.16.1}.
\hfill\qed

\subsection{$\nbigrtilde$-modules
associated to subvarieties of $\proj^m$}
\label{subsection;14.7.9.20}

\subsubsection{Setting}

We fix a homogeneous coordinate system
$[z_0:\cdots:z_m]$ of $\proj^m$.
We set $V:=H^0\bigl(\proj^m,\nbigo(1)\bigr)$.
We put 
$V_1:=\{\alpha_0 z_0\,|\,\alpha_0\in\cnum\}\subset V$
and 
$V_2:=\bigl\{\sum_{i=1}^m\alpha_iz_i\,\big|\,
 \alpha_i\in\cnum
 \bigr\}
 \subset V$.
We may regard $\alpha_0$
and $(\alpha_1,\ldots,\alpha_m)$
as coordinate systems of $V_1$ and $V_2$,
respectively.

We set $H_0:=\{z_0=0\}$.
Let $U$ be a smooth quasi-projective variety
with an immersion
$\iota_U:U\lrarr\proj^m\setminus H_0$.
For simplicity,
we assume that there exists 
a hypersurface 
$H\subset\proj^m$
for which $\iota_U(U)$
is a closed subset of 
$\proj^m\setminus(H_0\cup H)$.
We do not assume that $H$ is smooth.

We set $V_2^{\ast}:=V_2\setminus\{0\}$.
We shall construct some $\nbigrtilde$-modules on $V_2^{\ast}$
associated to $U$.

\subsubsection{A construction}
\label{subsection;14.7.16.10}

We take a smooth projective variety $Y$
with an open immersion
$\iota_1:U\lrarr Y$
and a morphism
$\iota_2:Y\lrarr \proj^m$
such that 
(i) $D_Y:=Y\setminus U$ is normal crossing,
(ii) $\iota_2\circ\iota_1=\iota_U$.
We have
$D_Y=\iota_{2}^{\ast}(H\cup H_0)$.
We set $D_{Y,V_2^{\ast}}:=D_Y\times V_2^{\ast}$.
We set $\iota_{2,V_2^{\ast}}:=\iota_2\times \id_{V_2^{\ast}}:
 Y\times V_2^{\ast}\lrarr \proj^m\times V_2^{\ast}$.
We have the meromorphic function
\[
 F^u_{Y,V_2^{\ast}}:=
 \iota_{2,V_2^{\ast}}^{\ast}
 \Bigl(
  \sum_{i=1}^m \alpha_iz_i/z_0\Bigr)
\]
on $(Y,D_Y)\times V_2^{\ast}$.
We obtain the integrable mixed twistor $\nbigd$-modules
$\nbigt_{\star}(F^u_{Y,V_2^{\ast}},D_{Y,V_2^{\ast}})$
$(\star=\ast,!)$
and the underlying $\nbigrtilde_{Y\times V_2^{\ast}}$-modules
$\nbigl_{\star}(F^u_{Y,V_2^{\ast}},D_{Y,V_2^{\ast}})$.
Let $\pi_{Y,V_2^{\ast}}:Y\times V_2^{\ast}\lrarr V_2^{\ast}$ 
be the projection.
We obtain the integrable mixed twistor $\nbigd$-modules
$\pi_{Y,V_2^{\ast}\dagger}^0
 \nbigt_{\star}(F^u_{Y,V_2^{\ast}},D_{Y,V_2^{\ast}})$
and the underlying $\nbigrtilde_{V_2^{\ast}}$-modules
$\pi_{Y,V_2^{\ast}\dagger}^0
 \nbigl_{\star}(F^u_{Y,V_2^{\ast}},D_{Y,V_2^{\ast}})$.

\vspace{.1in}
If $Y'$ and $\iota_{i}'$ be another choice,
then we take $Y''$ with morphisms
$\iota_1'':U\lrarr Y''$,
$a:Y''\lrarr Y'$ and $b:Y''\lrarr Y$
such that
$a\circ\iota_1''=\iota'_1$,
$b\circ\iota_1''=\iota_1$,
$\iota_2'\circ a=\iota_2\circ b=:\iota_{2}''$.
Then, we have natural commutative diagram of
isomorphisms:
\[
\begin{array}{ccccc}
 \pi^0_{Y,V_2^{\ast}\dagger}
 \nbigt_{!}(F^u_{Y,V_2^{\ast}},D_{Y,V_2^{\ast}})
 & \simeq &
 \pi^0_{Y'',V_2^{\ast}\dagger}
 \nbigt_{!}(F^u_{Y'',V_2^{\ast}},D_{Y'',V_2^{\ast}})
 & \simeq  &
 \pi^0_{Y',V_2^{\ast}\dagger}
 \nbigt_{!}(F^u_{Y',V_2^{\ast}},D_{Y',V_2^{\ast}})\\
 \darr & & \darr & & \darr \\
 \pi^0_{Y,V_2^{\ast}\dagger}
 \nbigt_{\ast}(F^u_{Y,V_2^{\ast}},D_{Y,V_2^{\ast}})
 & \simeq &
 \pi^0_{Y'',V_2^{\ast}\dagger}
 \nbigt_{\ast}(F^u_{Y'',V_2^{\ast}},D_{Y'',V_2^{\ast}})
 & \simeq  &
 \pi^0_{Y',V_2^{\ast}\dagger}
 \nbigt_{\ast}(F^u_{Y',V_2^{\ast}},D_{Y',V_2^{\ast}})\\
\end{array}
\]
In this sense, 
$\pi^0_{Y,V_2^{\ast}\dagger}
 \nbigt_{\star}(F^u_{Y,V_2^{\ast}},D_{Y,V_2^{\ast}})$
and the underlying $\nbigrtilde_{V_2^{\ast}}$-modules
$\pi^0_{Y,V_2^{\ast}\dagger}
 \nbigl_{\star}(F^u_{Y,V_2^{\ast}},D_{Y,V_2^{\ast}})$
are independent of the choice of $Y$.

We shall observe that 
$\pi^0_{Y,V_2^{\ast}\dagger}
 \nbigl_{\star}(F^u_{Y,V_2^{\ast}},D_{Y,V_2^{\ast}})$
are isomorphic to 
the $\nbigrtilde_{V_2^{\ast}}$-modules
given in \cite{Reichelt-Sevenheck2}.

\subsubsection{Comparison with a construction of Reichelt-Sevenheck}
\label{subsection;14.12.1.2}

Let $Z\subset\proj^m\times V$
be the $0$-set of 
the universal section $\sum \alpha_iz_i$ of 
$\nbigo_{\proj^m}(1)\boxtimes\nbigo_V$.
The projections $q_1:Z\lrarr \proj^m$
and $q_2:Z\lrarr V$ are smooth.
We set $Z_U:=Z\times_{\proj^m}U$.
Let $\iota_{Z_U}:Z_U\lrarr \proj^m\times V$
denote the natural inclusion.
Let $\pi_V:\proj^m\times V\lrarr V$
be the projection.

We have the variation of pure Hodge structure
$(\nbigo_{Z_U,}F)$ with the canonical real structure,
where the Hodge filtration $F$
is given by $F_0=\nbigo_{Z_U}$ and $F_{-1}=0$.
We have the associated pure Hodge module
which is also denoted by $(\nbigo_{Z_U},F)$. 
We obtain the mixed Hodge modules
$\iota_{Z_U\star}(\nbigo_{Z_U},F)$
$(\star=\ast,!)$
on $\proj^m\times V$,
and then
$\pi^0_{V\dagger}\iota_{Z_U\star}(\nbigo_{Z_U},F)$
on $V=V_1\times V_2$.
By applying the procedure in \S\ref{subsection;14.7.6.1},
and by taking the restriction to $V_2^{\ast}\subset V_2$,
we obtain the following 
$\nbigrtilde_{V_2^{\ast}}$-modules:
\[
 G_0\FL^{\loc}_{V_2}
 \pi^0_{V\dagger}\iota_{Z_U\star}(\nbigo_{Z_U},F)_{|V_2^{\ast}}
\]

\begin{prop}
\label{prop;14.11.12.3}
We have isomorphisms
of $\nbigrtilde_X$-modules:
\[
\lambda\cdot \pi^0_{Y,V_2^{\ast}\dagger}
 \nbigl_{\star}(F_Y^u,D_{Y,V_2^{\ast}})
\simeq
 G_0\FL^{\loc}_{V_2}
 \pi^0_{V\dagger}\iota_{Z_U\star}(\nbigo_{Z_U},F)_{|V_2^{\ast}}
\]
We also have the following commutative diagram:
\[
 \begin{CD}
 \lambda\cdot \pi^0_{Y,V_2^{\ast}\dagger}
 \nbigl_{!}(F^u_{Y,V_2^{\ast}},D_{Y,V_2^{\ast}})
 @>{\simeq}>>
 G_0\FL^{\loc}_{V_2}
 \pi^0_{V\dagger}\iota_{Z_U!}(\nbigo_{Z_U},F)_{|V_2^{\ast}}
 \\
 @VVV @VVV \\
 \lambda\cdot \pi^0_{Y,V_2^{\ast}\dagger}
 \nbigl_{\ast}(F^u_{Y,V_2^{\ast}},D_{Y,V_2^{\ast}})
 @>{\simeq}>>
 G_0\FL^{\loc}_{V_2}
 \pi^0_{V\dagger}\iota_{Z_U\ast}(\nbigo_{Z_U},F)_{|V_2^{\ast}}
 \end{CD}
\]
\end{prop}
\pf
We set $Z_0:=
 Z\cap\bigl(\proj^m\times V_1\times V_2^{\ast}
 \bigr)$.
Let $\proj_{V_1}$ denote the projective completion of $V_1$,
i.e.,
$\proj_{V_1}=\proj(V_1^{\lor}\oplus\cnum)$.
Let $\Zbar_0$ be the closure of
$Z_0$ in
$\proj^m\times\proj_{V_1}\times V_2^{\ast}$.
Let $\iota_{\Zbar_0}:
 \Zbar_0\lrarr
 \proj^m\times\proj_{V_1}\times V_2^{\ast}$
denote the natural inclusion.
By the construction of $\Zbar_0$,
we have the following equality
of meromorphic functions on $\Zbar_0$:
\[
 \iota_{\Zbar_0}^{\ast}\alpha_0
=-\iota_{\Zbar_0}^{\ast}\Bigl(
 \sum_{i=1}^m\alpha_iz_i/z_0
 \Bigr)
\]

Let $Z_{0U}:=Z_U\cap Z_0$,
and let $\Zbar_{0U}$ denote the closure of
$Z_{0U}$ in $\Zbar_0$.
Let $q:\Zbar_{0U}\lrarr \proj^m\times V_2^{\ast}$
denote the naturally induced morphism.
Note that
$Z_{0U}$ is naturally isomorphic to $U\times V_2^{\ast}$,
and that the restriction
$q_{|Z_{0U}}$ is an immersion.
We can take a smooth complex algebraic variety $B$
with projective morphisms
$\varphi_1:B\lrarr \Zbar_{0U}$ and
$\varphi_2:B\lrarr Y\times V_2^{\ast}$,
and an open immersion
$j:U\times V_2^{\ast}\subset B$
such that 
(i)
$q\circ\varphi_1
=\iota_{Y,V_2^{\ast}}\circ\varphi_2$,
(ii) $\varphi_a\circ j$ $(a=1,2)$
are the identity of $U\times V_2^{\ast}$,
(iii) $D_B:=B\setminus j(U\times V_2^{\ast})$
is a normal crossing hypersurface.
We set 
$G:=-\varphi_1^{\ast}\alpha_0
 =\varphi_2^{\ast}F_Y^u$.
We have the integrable mixed twistor $\nbigd$-modules
$\nbigt_{\star}(G,D_B)$ on $B$.
Let $\nbigl_{\star}(G,D_B)$
denote the underlying $\nbigrtilde_B$-modules.
Then, we have
the following natural isomorphisms
and the commutative diagrams:
\[
 \begin{CD}
 \varphi_{2\dagger}\nbigl_{!}(G,D_B)
 @>{\simeq}>>
 \nbigl_{!}(F^u_{Y,V_2^{\ast}},D_{Y,V_2^{\ast}})\\
 @VVV @VVV \\ 
\varphi_{2\dagger}\nbigl_{\ast}(G,D_B)
 @>{\simeq}>>
 \nbigl_{\ast}(F^u_{Y,V_2^{\ast}},D_{Y,V_2^{\ast}})
 \end{CD}
\]
We also have the following:
\[
 \begin{CD}
 \varphi_{1\dagger}\nbigl_{!}(G,D_B)
 @>{\simeq}>>
 R\bigl(
 \iota_{Z_U!}(\nbigo_{Z_U},F)
 \bigr)
 \otimes\nbigl(-\alpha_0)
 \\
 @VVV @VVV \\
 \varphi_{1\dagger}\nbigl_{\ast}(G,D_B)
 @>{\simeq}>>
 R\bigl(
 \iota_{Z_U\ast}(\nbigo_{Z_U},F)
 \bigr)
 \otimes\nbigl(-\alpha_0)
 \end{CD}
\]
Then, the claim of the theorem follows from
Proposition \ref{prop;14.7.9.30}.
\hfill\qed

\vspace{.1in}

Reichelt and Sevenheck considered 
the image $(M,F)$ of
$\iota_{Z_U!}(\nbigo_{Z_U},F)
\lrarr
 \iota_{Z_U\ast}(\nbigo_{Z_U},F)$,
that is the minimal extension of
$Z_U$ in $\proj^m\times V$.
They proved that
$G_0\FL^{\loc}_{V_2}\bigl(
 \pi^0_{V\dagger}(M,F)
 \bigr)_{|V_2^{\ast}}$
is isomorphic to the image of the natural morphism
$G_0\FL^{\loc}_{V_2}
 \pi^0_{V\dagger}\iota_{Z_U!}
 \bigl(\nbigo_{Z_U},F
 \bigr)_{|V_2^{\ast}}
\lrarr
G_0\FL^{\loc}_{V_2}
 \pi^0_{V\dagger}\iota_{Z_U\ast}
 \bigl(
 \nbigo_{Z_U},F
 \bigr)_{|V_2^{\ast}}$.

\begin{cor}
\label{cor;14.12.3.100}
$G_0\FL^{\loc}_{V_2}\bigl(
 \pi_{V\dagger}(M,F)
 \bigr)_{|V_2^{\ast}}$
is isomorphic to the image of 
the natural morphism
\[
  \lambda\cdot 
 \pi_{Y,V_2^{\ast}\dagger}^0
 \nbigl_{!}(F^u_{Y,V_2^{\ast}},D_{Y,V_2^{\ast}})
\lrarr
\lambda\cdot 
 \pi_{Y,V_2^{\ast}\dagger}^0
 \nbigl_{\ast}(F^u_{Y,V_2^{\ast}},D_{Y,V_2^{\ast}}).
\]
In particular,
$G_0\FL^{\loc}_{V_2}\bigl(
 \pi_{V\dagger}(M,F)
\bigr)_{|V_2^{\ast}}$
underlies a polarizable integrable pure twistor $\nbigd$-module.
\end{cor}
\pf
The first claim follows from Proposition \ref{prop;14.11.12.3}
and the result of Reichelt and Sevenheck mentioned above.
Because the morphism of the mixed twistor $\nbigd$-modules
$\pi_{Y,V_2^{\ast}\dagger}^0
 \nbigt_{!}(F^u_{Y,V_2^{\ast}},D_{Y,V_2^{\ast}})
\lrarr
 \pi_{Y,V_2^{\ast}\dagger}^0
 \nbigt_{\ast}(F^u_{Y,V_2^{\ast}},D_{Y,V_2^{\ast}})$
factors through the morphism of the pure twistor $\nbigd$-modules
\[
 \Gr^W_{\dim U+m}\pi_{Y,V_2^{\ast}\dagger}^0
 \nbigt_{!}(F^u_{Y,V_2^{\ast}},D_{Y,V_2^{\ast}})
\lrarr
 \Gr^W_{\dim U+m}\pi_{Y,V_2^{\ast}\dagger}^0
 \nbigt_{\ast}(F^u_{Y,V_2^{\ast}},D_{Y,V_2^{\ast}}),
\]
the second claim follows.
\hfill\qed

\section{Better behaved GKZ-systems and de Rham complexes}
\label{section;15.5.2.2}
\subsection{$\nbigd$-modules}

\label{subsection;15.4.26.1}

\subsubsection{The $\nbigd$-modules
associated to some better behaved GKZ-systems}
\label{subsection;15.4.26.2}

We recall a special version of better behaved GKZ-systems
introduced by L. Borisov and P. Horja \cite{Borisov-Horja}.
Let $\nbiga=\{\veca_1,\ldots,\veca_m\}\subset \seisuu^n$
be a finite subset generating $\seisuu^n$.
Let $K_{\real}(\nbiga)\subset\real^n$ denote the cone
$\bigl\{
 \sum_{j=1}^m r_j\veca_j\,\big|\,
 r_j\geq 0
 \bigr\}$
generated by $\nbiga$.
We set $K(\nbiga):=K_{\real}(\nbiga)\cap\seisuu^n$.
The semigroup $K(\nbiga)$ is the saturation of
$\seisuu_{\geq 0}\nbiga
=\bigl\{
\sum_{j=1}^m n_{j}\veca_j\,\big|\, 
 n_{j}\in\seisuu_{\geq 0}
\bigr\}$.
We denote 
$\veca_j=(a_{j1},\ldots,a_{jn})$.

Let $\Gamma\subset K(\nbiga)$ be any subset
such that
$\Gamma+\veca\subset \Gamma$ for any $\veca\in\nbiga$.
Let $\vecbeta\in\cnum^n$.
The following system of differential equations
$\GKZ(\nbiga,\Gamma,\vecbeta)$
for tuples of holomorphic functions 
$(\Phi_{\vecc}\,|\,\vecc\in \Gamma)$
on any open subset $\cnum^m$ 
is called the better behaved GKZ-hypergeometric system
associated to $(\nbiga,\Gamma,\vecbeta)$:
\[
 \del_{x_j}\Phi_{\vecc}=\Phi_{\vecc+\veca_j}
\quad(\forall \vecc\in \Gamma,\,\forall j=1,\ldots,m)
\]
\[
 \Bigl(
 \sum_{j}a_{ji}x_j\del_{x_j}+c_i-\beta_i
 \Bigr)
 \Phi_{\vecc}=0
\quad
 (\forall \vecc\in \Gamma,\,\forall i=1,\ldots,n)
\]
Here, $(x_1,\ldots,x_m)$ denotes the standard coordinate system
of $\cnum^m$.

Let us describe the $\nbigd_{\cnum^m}$-module
$M^{\GKZ}(\nbiga,\Gamma,\vecbeta)$
corresponding to the system $\GKZ(\nbiga,\Gamma,\vecbeta)$.
We introduce a free $\nbigo_{\cnum^m}$-module
$Q(\nbiga,\Gamma)$ generated by $\Gamma$.
Let $e(\vecc)$ denote the element corresponding to
$\vecc\in\Gamma$.
So, we have
$Q(\nbiga,\Gamma)=
 \bigoplus_{\vecc\in\Gamma}
 \nbigo_{\cnum^m}\cdot e(\vecc)$.
We introduce the action of 
$\nbigd_{\cnum^m}$ on
$Q(\nbiga,\Gamma)$
by $\del_{x_j}e(\vecc)=e(\vecc+\veca_j)$.
Let $J(\nbiga,\Gamma,\vecbeta)$
denote the $\nbigd_{\cnum^m}$-submodule
of $Q(\nbiga,\Gamma)$
generated by
$\bigl(
 \sum_{j}a_{ji}x_j\del_{x_j}+c_i-\beta_i
\bigr)e(\vecc)$
for $\vecc\in \Gamma$
and $i=1,\ldots,n$.
We set 
$M^{\GKZ}(\nbiga,\Gamma,\vecbeta):=
 Q(\nbiga,\Gamma)/J(\nbiga,\Gamma,\vecbeta)$.

Let $U\subset \cnum^m$ be any open subset.
An $\nbigo$-homomorphism
$\Phi:M^{\GKZ}(\nbiga,\Gamma,\vecbeta)_{|U}
\lrarr
 \nbigo_U$
is uniquely determined by
$\Phi_{\vecc}:=\Phi(e(\vecc))$ $(\vecc\in \Gamma)$.
Conversely, any tuple of holomorphic functions
$(\Phi_{\vecc}\,|\,\vecc\in\Gamma)$
determines an $\nbigo$-morphism
$M^{\GKZ}(\nbiga,\Gamma,\vecbeta)_{|U}
\lrarr
 \nbigo_U$.
Then, the following lemma is clear by the construction.

\begin{lem}
\label{lem;15.5.15.1}
Let $U\subset \cnum^m$ be any open subset.
The above correspondence induces
a bijective correspondence
between
morphisms of $\nbigd$-modules
$M^{\GKZ}(\nbiga,\Gamma,\vecbeta)_{|U}
\lrarr
 \nbigo_U$
and 
solutions of $\GKZ(\nbiga,\Gamma,\vecbeta)$.
\hfill\qed
\end{lem}

\begin{rem}
In {\rm\cite{Borisov-Horja}},
Borisov and Horja considered 
an $m$-tuple of a finitely generated abelian group,
instead of a finite subset of $\seisuu^n$.
In that sense, we consider the only special case.
However,
we do not impose the existence of an element
$\vecalpha\in(\seisuu^n)^{\lor}$
such that
$\vecalpha(\veca_j)=1$ $(j=1,\ldots,m)$.
\hfill\qed
\end{rem}

\subsubsection{Twisted de Rham complexes}
\label{subsection;15.4.26.3}

We recall some notation used in \S\ref{section;14.12.8.11}.
Let $T^n:=(\cnum^{\ast})^n$.
We consider the morphism
$\psi^{\aff}_{\nbiga}:T^n\lrarr \cnum^m$
given by 
$\psi^{\aff}_{\nbiga}(t_1,\ldots,t_n)
=(t^{\veca_1},\ldots,t^{\veca_m})$,
where $t^{\veca_j}=\prod_{i=1}^nt_i^{a_{ji}}$.
Let $X_{\nbiga}^{\aff}$ denote the closure
of the image of $\psi^{\aff}_{\nbiga}$.
Let $\check{X}^{\aff}_{\nbiga}\lrarr X^{\aff}_{\nbiga}$
be the normalization.
Let $\check{D}^{\aff}_{\nbiga}$
denote the complement of $T^n$ in
$\check{X}^{\aff}_{\nbiga}$.
Note that 
$\check{X}^{\aff}_{\nbiga}
=\Spec \cnum[K(\nbiga)]$
and 
$X^{\aff}_{\nbiga}
=\Spec\cnum[\seisuu_{\geq 0}\nbiga]$.

The family of Laurent polynomials 
$\sum_{j=1}^m x_jt^{\veca_j}$
induces a meromorphic function
$F_{\nbiga}$ on 
$(\check{X}^{\aff}_{\nbiga},\check{D}^{\aff}_{\nbiga})\times\cnum^m$.
We also have the logarithmic closed one form
$\kappa(\vecbeta)=\sum_{i=1}^n \beta_i dt_i/t_i$.
We obtain the relative algebraic de Rham complexes:
\[
 \nbigc^{\bullet}(\nbiga,\vecbeta)_{\bullet}:=
 \Bigl(
 \Omega^{\bullet}_{(\check{X}^{\aff}_{\nbiga}\times \cnum^m)/\cnum^m}
 (\log\check{D}^{\aff}_{\nbiga}\times\cnum^m),
 d+dF_{\nbiga}-\kappa(\vecbeta)
 \Bigr)
\]
\[
 \nbigc^{\bullet}(\nbiga,\vecbeta)_{\circ}:=
 \Bigl(
 \Omega^{\bullet}_{(\check{X}^{\aff}_{\nbiga},\check{D}^{\aff}_{\nbiga})
 \times \cnum^m/\cnum^m},
 d+dF_{\nbiga}-\kappa(\vecbeta)
 \Bigr)
\]

Let $\check{\pi}_{\nbiga}:
 \check{X}^{\aff}_{\nbiga}\times\cnum^m
\lrarr \cnum^m$ denote the projection.
Each $\check{\pi}_{\nbiga\ast}\nbigc^{k}(\nbiga,\vecbeta)_{\star}$
is naturally a $\nbigd_{\cnum^m}$-module
by $\del_{x_j}\bullet(g)=\del_{x_j}g+\del_{x_j}F_{\nbiga}\cdot g$.
The differential $d+dF_{\nbiga}-\kappa(\vecbeta)$ of
the complexes are compatible
with the actions of $\nbigd_{\cnum^m}$.
So, we obtain the following $\nbigd_{\cnum^m}$-modules
\[
 M_{\nbiga,\vecbeta,\star}:=
 \hyperr^n\check{\pi}_{\nbiga\ast}
 \nbigc^{\bullet}(\nbiga,\vecbeta)_{\star}
\quad
 (\star=\bullet,\circ).
\]

\begin{rem}
In the notation of {\rm \S\ref{subsection;14.11.23.221}}
(see also Example {\rm\ref{example;14.11.24.100}}),
we have the following commutative diagram:
\[
 \begin{CD}
 M_{\nbiga,0,\bullet}
 @>>> 
 M_{\nbiga,0,\circ}
 \\
 @V{\simeq}VV @V{\simeq}VV \\
 \pi_{\dagger}L_{\ast}(\nbiga,\cnum^m,\id)
 @>>>
 \pi_{\dagger}L_{!}(\nbiga,\cnum^m,\id)
 \end{CD}
\]
\hfill\qed
\end{rem}

\subsubsection{Comparison}

Let $K(\nbiga)^{\circ}$ denote the intersection of
$\seisuu^n$ and the interior part of $K_{\real}(\nbiga)$.

\begin{prop}
\label{prop;15.4.25.2}
We have the following natural isomorphisms
of $\nbigd_{\cnum^m}$-modules:
\[
 M^{\GKZ}(\nbiga,K(\nbiga),\vecbeta)
\simeq
 M_{\nbiga,\vecbeta,\bullet},
\quad\quad
  M^{\GKZ}(\nbiga,K(\nbiga)^{\circ},\vecbeta)
\simeq
 M_{\nbiga,\vecbeta,\circ}.
\]
\end{prop}
\pf
Although the claim is rather obvious by definition,
we give some more details.
Let $\psi_{\Sigma_1}:X_{\Sigma_1}\lrarr \check{X}^{\aff}_{\nbiga}$
be any toric projective resolution.
Set $D_{\Sigma_1}:=X_{\Sigma_1}\setminus T$.
According to \cite{Batyrev-VMHS},
we have 
$R\psi_{\Sigma_1\ast}\Omega^p_{X_{\Sigma_1}}(\log D_{\Sigma_1})
\simeq
\Omega^p_{\check{X}^{\aff}_{\nbiga}}(\log \check{D}^{\aff}_{\nbiga})$.
Let $\nbigi(p)$ denote the set of tuples
of integers $I=(i_1,\ldots,i_p)$ with
$1\leq i_1<\cdots <i_p\leq n$.
For any $I\in\nbigi(p)$,
we set 
$\tau_I:=
\prod_{j=1}^pt_{i_j}^{-1}
dt_{i_1}\wedge\cdots\wedge dt_{i_p}$.
Hence, we have
$\Omega_{\check{X}^{\aff}_{\nbiga}}^p
 (\log\check{D}^{\aff}_{\nbiga})$
is a free sheaf over $\nbigo_{\check{X}^{\aff}_{\nbiga}}$
with a basis
$\tau_I$ $(I\in S(p))$.
Hence, we have
\[
\check{\pi}_{\nbiga\ast}
\nbigc^p(\nbiga)_{\ast}
\simeq
 \bigoplus_{I\in S(p)}
 \bigoplus_{\vecc\in K(\nbiga)}
 \nbigo_{\cnum^m}
 t^{\vecc}\tau_I.
\]
We set 
$\omega:=\tau_{\{1,\ldots,n\}}$.
We have
$\del_{x_j}(t^{\vecc}\omega)
=t^{\vecc+\veca_j}\omega$.
Hence, we have the isomorphism of
$\nbigd_{\cnum^m}$-modules
$Q(\nbiga,K(\nbiga))\simeq
 \check{\pi}_{\nbiga\ast}
 \nbigc^n(\nbiga)_{\bullet}$
given by 
$e(\vecc)\longmapsto
 t^{\vecc}\omega$.

Let $\omega_i$ denote the inner product of
$\omega$ and $t_i\del/\del t_i$.
We have 
$\omega=(dt_i/t_i)\omega_i$.
We have
\[
 \bigl(d+dF_{\nbiga}-\kappa(\vecbeta)\bigr)
 (t^{\vecc}\omega_i)
=\Bigl(
 c_i-\beta_i+\sum_{j=1}^na_{ji}x_jt^{\veca_j}
 \Bigr)(t^{\vecc}\omega)
=\Bigl(
 c_i-\beta_i+\sum a_{ji}x_j\del_{x_j}
 \Bigr)(t^{\vecc}\omega).
\]
Hence, the image of 
$\check{\pi}_{\nbiga\ast}
 \nbigc^{n-1}(\nbiga)_{\bullet}
\lrarr
 \check{\pi}_{\nbiga\ast}
 \nbigc^{n}(\nbiga)_{\bullet}$
is identified with
$J(\nbiga,K(\nbiga),\vecbeta)$.
Hence, we obtain the desired isomorphism
$M^{\GKZ}(\nbiga,K(\nbiga),\vecbeta)\simeq
 M_{\nbiga,\vecbeta,\bullet}$.

\vspace{.1in}
According to \cite{Danilov-de-Rham},
the space of global sections of 
$\Omega^p_{(\check{X}^{\aff}_{\nbiga},\check{D}^{\aff}_{\nbiga})}$
is
$\bigoplus_{I\in S(p)}
 \bigoplus_{\vecc\in K(\nbiga)^{\circ}}
 \cnum\cdot t^{\vecc}\tau_I$.
Hence, we have
\[
 \check{\pi}_{\nbiga\ast}
 \nbigc^p(\nbiga)_{\circ}
=\bigoplus_{I\in S(p)}
 \bigoplus_{\vecc\in K(\nbiga)^{\circ}}
 \nbigo_{\cnum^m}
 t^{\vecc}\tau_I.
\]
Then, as in the case of $\bullet$,
we obtain the desired isomorphism 
$M^{\GKZ}(\nbiga,K(\nbiga)^{\circ},\vecbeta)
\simeq
 M_{\nbiga,\vecbeta,\circ}$.
\hfill\qed

\subsubsection{Special cases}
\label{subsection;15.5.11.2}

For any $\vecp=(p_1,\ldots,p_m)\in \seisuu^m$,
we put $\supp_+(\vecp):=\{j\,|\,p_j\geq 0\}$
and $\supp_-(\vecp):=\{j\,|\,p_j\leq 0\}$.
We set
\[
 \square_{\vecp}
=\prod_{j\in\supp_+(\vecp)}\del_{x_j}^{p_j}
-\prod_{j\in\supp_-(\vecp)}\del_{x_j}^{-p_j}.
\]
We have the morphism
$\seisuu^m\lrarr \seisuu^n$
given by
$\vecp=(p_1,\ldots,p_m)\longmapsto
\sum_{j=1}^m p_j\veca_j$.
Let $L_{\nbiga}$ denote the kernel.

Suppose that $\Gamma=\seisuu_{\geq 0}\nbiga+\vecc_0$
for an element $\vecc_0\in K(\nbiga)$.
Then, as remarked in \cite{Borisov-Horja},
$\GKZ(\nbiga,\Gamma,\vecbeta)$
is equivalent to 
the following ordinary GKZ-hypergeometric system
$\GKZ^{\ord}(\nbiga,\vecbeta-\vecc_0)$
for holomorphic functions $\Phi_{\vecc_0}$
on any open subset of $\cnum^m$:
\[
 \square_{\vecp}\Phi_{\vecc_0}=0
\quad
 (\forall\vecp\in L_{\nbiga})
\]
\[
 \Bigl(c_{0i}-\beta_i+\sum_{j=1}^ma_{ji}x_j\del_{x_j}\Bigr)
 \Phi_{\vecc_0}=0
\quad
 (i=1,\ldots,n)
\]
For any $\vecgamma\in\cnum^n$,
let $I(\nbiga,\vecgamma)$ denote
the left ideal of $\nbigd_{\cnum^m}$
generated by 
$\square_{\vecp}$ $(\vecp\in L_{\nbiga})$
and $-\gamma_i+\sum_{j=1}^ma_{ji}x_j\del_{x_j}$
$(i=1,\ldots,n)$.
Then, we have a natural isomorphism
$M^{\GKZ}(\nbiga,\Gamma,\vecbeta)
\simeq
 \nbigd_{\cnum^m}/I(\nbiga,\vecbeta-\vecc_0)$.

Suppose that we are given 
$\Gamma_1\subset\Gamma_2\subset K(\nbiga)$
such that
$\Gamma_i=\vecc_i+\seisuu_{\geq 0}\nbiga$.
We have the natural morphism
$M^{\GKZ}(\nbiga,\Gamma_1,\vecbeta)
\lrarr
 M^{\GKZ}(\nbiga,\Gamma_2,\vecbeta)$
induced by the inclusion
$Q(\nbiga,\Gamma_1)\subset Q(\nbiga,\Gamma_2)$.
The above isomorphisms induce
$g_{\vecc_2,\vecc_1}:
 \nbigd_{\cnum^m}/I(\nbiga,\vecc_1,\vecbeta)
\lrarr
 \nbigd_{\cnum^m}/I(\nbiga,\vecc_2,\vecbeta)$.

We have an expression 
$\vecc_1=\vecc_2+\sum_{j=1}^m b_j\veca_j$
$(b_j\in\seisuu_{\geq 0})$.
We have the morphism of $\nbigd_{\cnum^m}$-modules
$\htilde_{\vecc_2,\vecc_1}:
 \nbigd_{\cnum^m}\lrarr \nbigd_{\cnum^m}$
induced by the right multiplication of
$B_{\vecc_2,\vecc_1}=
 \prod_{j=1}^m\del_j^{b_j}$.
We have
$\square_{\vecp}\cdot B_{\vecc_2,\vecc_1}
=B_{\vecc_2,\vecc_1}\cdot \square_{\vecp}$.
We also have
\[
 \Bigl(
 c_{1i}-\beta_i+\sum_{j=1}^ma_{ji}x_j\del_{x_j}
 \Bigr)
 \cdot\prod_{j=1}^m\del_j^{b_j}
=
 \prod_{j=1}^m\del_j^{b_j}\cdot
  \Bigl(
 c_{2i}-\beta_i+\sum_{j=1}^ma_{ji}x_j\del_{x_j}
 \Bigr).
\]
Hence, we have the induced morphism
$h_{\vecc_2,\vecc_1}:
  \nbigd_{\cnum^m}/I(\nbiga,\vecc_1,\vecbeta)
 \lrarr
 \nbigd_{\cnum^m}/I(\nbiga,\vecc_2,\vecbeta)$.
It is easy to check
$g_{\vecc_2,\vecc_1}
=h_{\vecc_2,\vecc_1}$.

\vspace{.1in}

Let us consider some more special cases.
Suppose that
$K(\nbiga)=\seisuu_{\geq 0}\nbiga$,
i.e.,
$X_{\nbiga}^{\aff}
=\check{X}^{\aff}_{\nbiga}$.
Then, 
the system $\GKZ(\nbiga,K(\nbiga),0)$
is equivalent to
$\GKZ^{\ord}(\nbiga,0)$,
and we have
$M^{\GKZ}(\nbiga,K(\nbiga),0)
\simeq
 \nbigd_{\cnum^m}/I(\nbiga,0)$.

Suppose moreover that
$K(\nbiga)^{\circ}=K(\nbiga)+\vecc_1$
for some $\vecc_1\in K(\nbiga)$.
Then, 
$\GKZ(\nbiga,K(\nbiga)^{\circ},0)$
is equivalent to
$\GKZ^{\ord}(\nbiga,-\vecc_1)$,
and we have
$M^{\GKZ}(\nbiga,K(\nbiga)^{\circ},0)
\simeq
 \nbigd_{\cnum^m}/I(\nbiga,-\vecc_1)$.

The morphism
$M^{\GKZ}(\nbiga,K(\nbiga)^{\circ},0)
\lrarr
 M^{\GKZ}(\nbiga,K(\nbiga),0)$
induces
$\nbigd_{\cnum^m}/I(\nbiga,-\vecc_1)
\lrarr
 \nbigd_{\cnum^m}/I(\nbiga,0)$.
If $\vecc_1=\sum_{i=1}^{\ell} \veca_i$ for some $\ell$,
it is equal to the morphism induced by
the multiplication of
$\prod_{j=1}^{\ell}\del_{x_j}$.

\subsection{$\nbigrtilde$-modules}
\label{subsection;15.5.2.1}

\subsubsection{Systems of differential equations
 and the associated $\nbigrtilde$-modules}
\label{subsection;15.5.11.10}

We continue to use the notation in \S\ref{subsection;15.4.26.1}.
Let $\Gamma\subset K(\nbiga)$ be any subset 
satisfying $\Gamma+\veca\subset \Gamma$
for any $\veca\in\nbiga$. 
We consider the following system of differential equations
$\GKZ_{\nbigrtilde}(\nbiga,\Gamma,\vecbeta)$
for a tuple $\Phi_{\Gamma}=(\Phi_{\vecc}\,|\,\vecc\in \Gamma)$
of holomorphic functions on any open subset of 
$\cnum\times\cnum^m
=\{(\lambda,x_1,\ldots,x_m)\}$:
\[
 \lambda\del_{x_j}\Phi_{\vecc}
=\Phi_{\vecc+\veca_j},
\quad
(\forall \vecc\in \Gamma,\,\,j=1,\ldots,m)
\]
\[
 \Bigl(
 \lambda^2\del_{\lambda}
+n\lambda
+\sum_{j=1}^m\lambda x_j\del_{x_j}
 \Bigr)
\Phi_{\vecc}
=0,
\quad
(\forall \vecc\in \Gamma)
\]
\[
 \Bigl(
 \lambda (c_i-\beta_i)+\sum_{j=1}^m a_{ji}\lambda x_j\del_{x_j}
 \Bigr)
 \Phi_{\vecc}=0,
\quad
(\forall\vecc\in \Gamma,\,\,\,i=1,\ldots,n)
\]

We describe the corresponding
$\nbigrtilde_{\cnum^m}$-module.
We introduce a free
$\nbigo_{\cnum\times\cnum^m}$-module
 $\nbigq(\nbiga,\Gamma)$
generated by $\Gamma$.
Let $\overline{e}(\vecc)$ denote the section
corresponding to $\vecc\in\Gamma$.
So, we have
$\nbigq(\nbiga,\Gamma)
=\bigoplus_{\vecc\in\Gamma}
 \nbigo_{\cnum\times\cnum^m}\overline{e}(\vecc)$.
It is an $\nbigrtilde_{\cnum^m}$-module
by the actions
$\lambda\del_{x_j}\overline{e}(\vecc)
=\overline{e}(\vecc+\veca_j)$
and
\[
 \lambda^2\del_{\lambda}
 \overline{e}(\vecc)
=-n\lambda \overline{e}(\vecc)
-\sum_{j=1}^mx_j\lambda\del_{x_j}\overline{e}(\vecc).
\]
Let $\nbigj(\nbiga,\Gamma,\vecbeta)$
denote the $\nbigrtilde$-submodule of
$\nbigq(\nbiga,\Gamma)$
generated by
$\bigl(
 \lambda (c_i-\beta_i)+\sum_{j=1}^m a_{ji}\lambda x_j\del_{x_j}
 \bigr)
 \overline{e}(\vecc)$
for $\vecc\in\Gamma$
and $i=1,\ldots,n$.
We set
$\nbigm^{\GKZ}(\nbiga,\Gamma,\vecbeta):=
 \nbigq(\nbiga,\Gamma)/\nbigj(\nbiga,\Gamma,\vecbeta)$.

As in \S\ref{subsection;15.4.26.2},
we have the following.
\begin{lem}
\label{lem;15.5.15.11}
Let $U$ be any open subset in $\cnum\times\cnum^m$.
 We have the natural bijective correspondence
between the $\nbigrtilde$-homomorphisms
$\nbigm^{\GKZ}(\nbiga,\Gamma,\vecbeta)_{|U}
\lrarr
 \nbigo_{U}$
and the solutions of $\GKZ_{\nbigrtilde}(\nbiga,\Gamma,\vecbeta)$
on $U$.
\end{lem}

\subsubsection{Comparison with the twisted de Rham complexes}

We use the notation in \S\ref{subsection;15.4.26.3}.
Let $q:\cnum\times \check{X}_{\nbiga}^{\aff}\times\cnum^m
\lrarr \check{X}_{\nbiga}^{\aff}$
denote the projection.
We set
\[
 \overline{\nbigc}^{k}(\nbiga,\vecbeta)_{\bullet}:=
 \lambda^{-k}
 q^{\ast}\Omega^k_{\check{X}^{\aff}_{\nbiga}}
 (\log \check{D}^{\aff}_{\nbiga}),
\quad\quad
 \overline{\nbigc}^k(\nbiga,\vecbeta)_{\circ}:=
 \lambda^{-k}
 q^{\ast}\Omega^k_{(\check{X}^{\aff}_{\nbiga},\check{D}^{\aff}_{\nbiga})}.
\]
With the differential induced by
$d+d(\lambda^{-1}F_{\nbiga})-\kappa(\vecbeta)$,
we obtain the complexes of sheaves
$\overline{\nbigc}^{\bullet}(\nbiga,\vecbeta)_{\star}$
$(\star=\bullet,\circ)$.

Let $\check{\pi}_{\nbiga}:
 \cnum\times\check{X}^{\aff}_{\nbiga}\times\cnum^m
\lrarr
 \cnum\times \cnum^m$
denote the projection.
Each 
$\check{\pi}_{\nbiga\ast}
 \overline{\nbigc}^{k}(\nbiga,\vecbeta)_{\star}$
is naturally an $\nbigrtilde_{\cnum^m}$-module
by $\lambda\del_{x_j}\bullet(g)=
 \lambda\del_{x_j}g+\del_{x_j}F_{\nbiga}\cdot g$
 $(j=1,\ldots,m)$
and 
 $\lambda^2\del_{\lambda}\bullet(g)
 =\lambda^2\del_{\lambda}g
-F_{\nbiga}\cdot g$.
The differentials 
$d+d(\lambda^{-1}F_{\nbiga})-\kappa(\vecbeta)$
of the complexes are compatible
with the actions of $\nbigrtilde_{\cnum^m}$.
So, we obtain the following $\nbigrtilde$-modules
\[
 \nbigm_{\nbiga,\vecbeta,\star}:=
 \hyperr^n\check{\pi}_{\nbiga\ast}
 \overline{\nbigc}^{\bullet}(\nbiga,\vecbeta)_{\star}
\quad
 (\star=\bullet,\circ).
\]

\begin{prop}
\label{prop;15.5.11.1}
We have the following natural isomorphisms
of $\nbigrtilde_{\cnum^m}$-modules:
\[
 \nbigm^{\GKZ}(\nbiga,K(\nbiga),\vecbeta)\simeq
 \nbigm_{\nbiga,\vecbeta,\bullet},
\quad
 \nbigm^{\GKZ}(\nbiga,K(\nbiga)^{\circ},\vecbeta)
\simeq
 \nbigm_{\nbiga,\vecbeta,\circ}.
\]
\end{prop}
\pf
The argument is similar to that in
the proof of Proposition \ref{prop;15.4.25.2}.
We give just an indication.
We have
\[
 \check{\pi}_{\nbiga\ast}
 \overline{\nbigc}^p(\nbiga,\vecbeta)_{\bullet}
=\bigoplus_{I\in S(p)}
 \bigoplus_{\vecc\in K(\nbiga)}
 \nbigo_{\cnum\times\cnum^m}
 t^{\vecc}\lambda^{-p}\tau_I.
\]
We have the following equalities.
\[
 \lambda\del_{x_j}(t^{\vecc}\lambda^{-n}\omega)
=t^{\vecc+\veca_j}\lambda^{-n}\omega
\]
\[
 \Bigl(
 \lambda^2\del_{\lambda}
+n\lambda+\sum_{j=1}^m \lambda x_j\del_{x_j}
 \Bigr)
 t^{\vecc}\lambda^{-n}\omega
=0
\]
We also have the following:
\[
 \Bigl(\lambda (c_i-\beta_i)+\sum_{j=1}^m  a_{ji}\lambda x_j\del_{x_j}\Bigr)
 (t^{\vecc}\lambda^{-n}\omega)
=
 \Bigl(
 d+d(\lambda^{-1}F_{\nbiga})-\kappa(\vecbeta)
 \Bigr)
 (t^{\vecc}\lambda^{-n+1}\omega_i)
\]
Then, we obtain the claim in the case for $\bullet$
as in the case of Proposition \ref{prop;15.4.25.2}.
The claim for $\circ$ is obtained in a similar way.
\hfill\qed

\begin{rem}
\label{rem;15.5.15.20}
As studied in {\rm\S\ref{subsection;14.12.15.100}},
 we have the following commutative diagram:
\[
 \begin{CD}
 \pi_{\dagger}^0\nbigl_{!}(\nbiga,\cnum^m,\id)
@>>>
 \pi_{\dagger}^0\nbigl_{\ast}(\nbiga,\cnum^m,\id)\\
@V{\simeq}VV @V{\simeq}VV  \\
 \check{\pi}_{\ast}
 \overline{\nbigc}^{\bullet}(\nbiga,0)_{\circ}
@>>>
  \check{\pi}_{\ast}
 \overline{\nbigc}^{\bullet}(\nbiga,0)_{\bullet} 
 \end{CD}
\]
Hence, Proposition {\rm\ref{prop;15.5.11.1}}
gives expressions of 
$\pi_{\dagger}^0\nbigl_{\star}(\nbiga,\cnum^m,\id)$
$(\star=\ast,!)$
as systems of differential equations.
\hfill\qed
\end{rem}

\subsubsection{Special case}
\label{subsection;15.5.11.20}

We use the notation in \S\ref{subsection;15.5.11.2}.
For any $\vecp\in \seisuu^m$,
we set
\[
 \overline{\square}_{\vecp}:=
 \prod_{j\in \supp_+(\vecp)}
 (\lambda\del_j)^{p_j}
-\prod_{i\in\supp_-(\vecp)}
 (\lambda\del_j)^{-p_j}.
\]

Suppose that $\Gamma=\seisuu_{\geq 0}\nbiga+\vecc_0$
for an element $\vecc_0\in K(\nbiga)$.
Then,
$\GKZ_{\nbigrtilde}(\nbiga,\Gamma,\vecbeta)$
is equivalent to the following system
$\GKZ_{\nbigrtilde}^{\ord}(\nbiga,\vecbeta-\vecc_0)$
for $\Phi_{\vecc_0}$:
\[
 \overline{\square}_{\vecp}\Phi_{\vecc_0}=0
\quad
 (\forall\vecp\in L_{\nbiga})
\]
\[
 \Bigl(
 \lambda^2\del_{\lambda}
 +n\lambda+\sum_{j=1}^mx_j\lambda\del_{x_j}
 \Bigr)\Phi_{\vecc_0}=0
\]
\[
 \Bigl(
 \lambda(c_{0i}-\beta_i)
+\sum_{j=1}^ma_{ji}x_j\lambda\del_{x_j}\Bigr)
 \Phi_{\vecc_0}=0
\quad
 (i=1,\ldots,n)
\]
For any $\vecgamma\in\seisuu^m$,
let $\nbigi(\nbiga,\vecgamma)$ denote
the left ideal of $\nbigrtilde_{\cnum^m}$
generated by 
$\overline{\square}_{\vecp}$ $(\vecp\in L_{\nbiga})$,
$\lambda^2\del_{\lambda}+n\lambda
+\sum_{j=1}^m\lambda x_j\del_{x_j}$,
and $-\lambda\gamma_i+\sum_{j=1}^ma_{ji}\lambda x_j\del_{x_j}$
$(i=1,\ldots,n)$.
Then, we have a natural isomorphism
\[
 \nbigm^{\GKZ}(\nbiga,\Gamma,\vecbeta)
\simeq
 \nbigrtilde_{\cnum^m}/\nbigi(\nbiga,\vecbeta-\vecc_0).
\]

Suppose that we are given 
$\Gamma_1\subset\Gamma_2\subset K(\nbiga)$
such that
$\Gamma_i=\vecc_i+\seisuu_{\geq 0}\nbiga$.
We have the natural morphism
$\nbigm^{\GKZ}(\nbiga,\Gamma_1,\beta)
\lrarr
 \nbigm^{\GKZ}(\nbiga,\Gamma_2,\beta)$
induced by the inclusion
$\nbigq(\nbiga,\Gamma_1)\subset 
 \nbigq(\nbiga,\Gamma_2)$.
The above isomorphisms induce
$g_{\vecc_2,\vecc_1}:
 \nbigrtilde_{\cnum^m}/\nbigi(\nbiga,\vecc_1,\beta)
\lrarr
 \nbigrtilde_{\cnum^m}/\nbigi(\nbiga,\vecc_2,\beta)$.

We have an expression $\vecc_1=\vecc_2+\sum_{i=1}^m b_i\veca_i$
$(b_i\in\seisuu_{\geq 0})$.
We have the morphism of
$\nbigrtilde_{\cnum^m}$-modules
$\htilde_{\vecc_2,\vecc_1}:
 \nbigrtilde_{\cnum^m}
\lrarr
 \nbigrtilde_{\cnum^m}$
induced by the right multiplication of
$\overline{B}_{\vecc_2,\vecc_1}:=
 \prod_{j=1}^m
 (\lambda\del_{x_j})^{b_j}$.
We have
$\square_{\vecp}\cdot
 \overline{B}_{\vecc_2,\vecc_1}
=\overline{B}_{\vecc_2,\vecc_1}\cdot
 \square_{\vecp}$.
We have
\[
 \Bigl(
 \lambda(c_{2i}-\beta_i)
+\sum_{j=1}^ma_{ji}x_j\lambda\del_{x_j}
 \Bigr)
 \prod_{k=1}^m(\lambda\del_{x_k})^{b_k}
= \prod_{k=1}^m(\lambda\del_{x_k})^{b_k}
  \Bigl(
 \lambda(c_{1i}-\beta_i)
+\sum_{j=1}^ma_{ji}x_j\lambda\del_{x_j}
 \Bigr)
\]
We have
$\bigl(\lambda^2\del_{\lambda}+x_j\lambda\del_{x_j}\bigr)
 \lambda \del_{x_j}
=\lambda\del_{x_j}
\bigl(\lambda^2\del_{\lambda}+x_j\lambda\del_{x_j}\bigr)$.
Hence, 
\[
 \Bigl(
 \lambda^2\del_{\lambda}
+n\lambda
+\sum_{j=1}^mx_j\lambda\del_{x_j}
 \Bigr)\cdot
 \prod_{k=1}^m(\lambda\del_k)^{b_k}
=
 \prod_{k=1}^m(\lambda\del_k)^{b_k}
\cdot
 \Bigl(
 \lambda^2\del_{\lambda}
+n\lambda
+\sum_{j=1}^mx_j\lambda\del_{x_j}
 \Bigr)
\]
Hence, we have the induced morphism
$h_{\vecc_2,\vecc_1}:
 \nbigrtilde_{\cnum^m}/\nbigi(\nbiga,\vecc_1,\vecbeta)
\lrarr
 \nbigrtilde_{\cnum^m}/\nbigi(\nbiga,\vecc_2,\vecbeta)$.
It is easy to see that
$g_{\vecc_2,\vecc_1}=h_{\vecc_2,\vecc_1}$.

\vspace{.1in}

Let us consider some more special cases.
Suppose that
$K(\nbiga)=\seisuu_{\geq 0}\nbiga$,
i.e.,
$X_{\nbiga}^{\aff}
=\check{X}^{\aff}_{\nbiga}$.
Then, 
the system
$\GKZ_{\nbigrtilde}\bigl(\nbiga,K(\nbiga),0\bigr)$
is equivalent to
$\GKZ_{\nbigrtilde}^{\ord}(\nbiga,0)$,
and we have
$\nbigm^{\GKZ}(\nbiga,K(\nbiga),0)
\simeq
 \nbigrtilde_{\cnum^m}/\nbigi(\nbiga,0)$.

Suppose moreover that
$K(\nbiga)^{\circ}=K(\nbiga)+\vecc_1$
for some $\vecc_1\in K(\nbiga)$.
Then, 
$\GKZ_{\nbigrtilde}(\nbiga,K(\nbiga)^{\circ},0)$
is equivalent to
$\GKZ_{\nbigrtilde}^{\ord}(\nbiga,-\vecc_1)$,
and we have
$\nbigm^{\GKZ}(\nbiga,K(\nbiga)^{\circ},0)
\simeq
 \nbigrtilde_{\cnum^m}/\nbigi(\nbiga,-\vecc_1)$.

The morphism
$\nbigm^{\GKZ}(\nbiga,K(\nbiga)^{\circ},0)
\lrarr
 \nbigm^{\GKZ}(\nbiga,K(\nbiga),0)$
induces
$\nbigrtilde_{\cnum^m}/\nbigi(\nbiga,-\vecc_1)
\lrarr
 \nbigrtilde_{\cnum^m}/\nbigi(\nbiga,0)$.
If $\vecc_1=\sum_{j=1}^{\ell} \veca_j$ for some $\ell$,
it is equal to the morphism induced by
the multiplication of
$\prod_{j=1}^{\ell}(\lambda\del_{x_j})$.

\vspace{.1in}

\noindent
{\em Address\\
Research Institute for Mathematical Sciences,
Kyoto University,
Kyoto 606-8502, Japan\\
takuro@kurims.kyoto-u.ac.jp
}


\begin{thebibliography}{99}

\bibitem{Adolphson}
A. Adolphson,
{\em Hypergeometric functions and rings generated by monomials},
Duke Math. J. {\bf 73}, (1994), 269--290

\bibitem{Adolphson-Sperber}
A. Adolphson, 
S. Sperber,
{\em A-hypergeometric systems that come from geometry},
Proc. Amer. Math. Soc. {\bf 140}, (2012), 2033--2042

\bibitem{Batyrev-VMHS}
V. V. Batyrev, 
{\em Variations of the mixed Hodge structure 
of affine hypersurfaces in algebraic tori},
Duke Math. J. {\bf 69}, (1993), 349--409.

\bibitem{beilinson2}
A. Beilinson,
{\em How to glue perverse sheaves},
in; {\em $K$-theory, arithmetic and geometry
 (Moscow, 1984--1986)}, 
Lecture Notes in Math., {\bf 1289},
Springer, Berlin, (1987),
42--51.

\bibitem{Borisov-Horja}
L. Borisov, P. Horja,
{\em On the better behaved version of 
the GKZ hypergeometric system},
Math. Ann. {\bf 357} (2013), 585--603.

\bibitem{Chiang-Klemm-Yau-Zaslow}
T.-M. Chiang, 
A. Klemm,
S.-T. Yau,
E. Zaslow, 
{\em Local mirror symmetry: calculations 
and interpretations},
Adv. Theor. Math. Phys. {\bf 3}, (1999), 495--565.

\bibitem{Cox-et-al-toric}
D. A. Cox,
J. B. Little,
H. K. Schenck,
{\em Toric varieties},
Graduate Studies in Mathematics, {\bf 124},
American Mathematical Society, Providence, RI, 2011.

\bibitem{Danilov-de-Rham}
V. I. Danilov,
{\em de Rham complex on toroidal variety},
Algebraic geometry (Chicago, IL, 1989), 
Lecture Notes in Math., {\bf 1479}, Springer, Berlin, (1991),
26--38.

\bibitem{Douai-Sabbah1}
A. Douai, C. Sabbah, 
{\em Gauss-Manin systems, 
Brieskorn lattices and Frobenius structures. I.},
Ann. Inst. Fourier (Grenoble), {\bf 53}, (2003), 1055--1116.

\bibitem{Douai-Sabbah2}
A. Douai, C. Sabbah,
{\em Gauss-Manin systems, Brieskorn lattices and Frobenius structures. II.},
Frobenius manifolds,  
Aspects Math., E36, Vieweg, Wiesbaden, (2004),
1--18.

\bibitem{Esnault-Sabbah-Yu}
H. Esnault, 
C. Sabbah, 
J.-D. Yu,
(with an appendix by M. Saito),
{\em $E_1$-degeneration of the irregular Hodge filtration},
arXiv:1302.4537

\bibitem{Fan}
H. Fan,
{\em Schr\"{o}dinger equations, deformation theory and
$tt^*$-geometry},
arXiv:1107.1290 

\bibitem{Fulton-toric}
W. Fulton,
{\em Introduction to toric varieties},
Annals of Mathematics Studies, {\bf 131},
Princeton University Press, Princeton, NJ, 1993.

\bibitem{Hertling}
C. Hertling,
{\em $tt^{\ast}$geometry, 
Frobenius manifolds, their connections, 
and the construction for singularities},
J. Reine Angew. Math. {\bf 555}, (2003), 77--161.

\bibitem{GKZ}
I. M. Gelʹfand, 
M. M. Kapranov, 
A. V. Zelevinsky, 
{\em Generalized Euler integrals and 
A-hypergeometric functions.}
Adv. Math. {\bf 84}, (1990), 255--271.

\bibitem{GKZ-book}
I. M. Gelʹfand, 
M. M. Kapranov,
A. V. Zelevinsky,
{\em Discriminants, resultants, and multidimensional determinants},
Birkh\"{a}user Boston, Inc., Boston, MA, 1994.

\bibitem{Givental}
A. Givental,
{\em A mirror theorem for toric complete intersections},
In 
{\em Topological field theory, primitive forms and related topics
(Kyoto, $1996$)}, 
Progr. Math., {\bf 160},
Birkh\"auser Boston, Boston, MA, (1998)
141--175.

\bibitem{Hertling-Manin}
C. Hertling, Y. Manin,
{\em Unfoldings of meromorphic connections and 
a construction of Frobenius manifolds},
in {\em Frobenius manifolds}, 
Aspects Math., {\bf E36}, Vieweg, Wiesbaden, (2004),
113--144.

\bibitem{Hotta-Takeuchi-Tanisaki}
R. Hotta,
K. Takeuchi,
T. Tanisaki,
{\em D-modules, perverse sheaves, and representation theory},
Progress in Mathematics, {\bf 236},
Birkh\"{a}user Boston, Inc., Boston, MA, 2008. 

\bibitem{Iritani-convergence}
H. Iritani,
{\em Convergence of quantum cohomology 
by quantum Lefschetz},
J. Reine Angew. Math. {\bf 610}, (2007), 29--69.

\bibitem{Iritani}
H. Iritani,
{\em An integral structure in quantum cohomology 
and mirror symmetry for toric orbifolds},
Adv. Math. {\bf 222}, (2009), 1016--1079.

\bibitem{Iritani-quantum-cohomology-period}
H. Iritani,
{\em Quantum cohomology and periods},
Ann. Inst. Fourier (Grenoble) {\bf 61} (2011),
2909--2958.

\bibitem{Iritani-shift-operators}
H. Iritani,
{\em Shift operators and toric mirror theorem},
arXiv:1411.6840

\bibitem{Iritani-Mann-Mignon}
H. Iritani,
E. Mann, 
T. Mignon,
{\em
 Quantum Serre in terms of quantum D-modules},
arXiv:1412.4523




\bibitem{kashiwara-mixed-Hodge}
M. Kashiwara,
{\em A study of variation of mixed Hodge structure},
Publ. Res. Inst. Math. Sci. {\bf 22}, (1986)
991--1024. 

\bibitem{Kashiwara-Schapira}
M. Kashiwara and P. Schapira,
{\em Sheaves on manifolds},
Springer-Verlag, Berlin, 1990

\bibitem{kashiwara_text}
M. Kashiwara,
{\em $D$-modules and microlocal calculus,}
Translations of Mathematical Monographs, {\bf 217},
American Mathematical Society, Providence, (2003).

\bibitem{Katzarkov-Kontsevich-Pantev}
L. Katzarkov, 
M. Kontsevich,
T. Pantev, 
{\em Hodge theoretic aspects of mirror symmetry},
In 
{\em From Hodge theory to integrability and TQFT $tt^{\ast}$-geometry} 
Proc. Sympos. Pure Math., {\bf 78}, Amer. Math. Soc., Providence, RI, (2008)
87--174.

\bibitem{Konishi-Minabe-cubic}
Y. Konishi and S. Minabe, 
{\em Local Gromov-Witten invariants of cubic surfaces 
via nef toric degeneration},
Ark. Mat. {\bf 47}, (2009),  345--360.

\bibitem{Konishi-Minabe}
Y. Konishi,
S. Minabe,
{\em Local B-model and mixed Hodge structure},
Adv. Theor. Math. Phys. {\bf 14}, (2010), 1089--1145.

\bibitem{Konishi-Minabe2}
Y. Konishi,
S. Minabe,
{\em Mixed Frobenius Structure and Local A-model},
arXiv:1209.5550

\bibitem{Konishi-Minabe3}
Y. Konishi,
S. Minabe,
{\em Local Quantum Cohomology and Mixed Frobenius Structure},
arXiv:1405.7476 

\bibitem{Kouchnirenko}
A. G. Kouchnirenko,
{Poly\`{e}dres de Newton et nombres de Milnor},
Invent. Math. {\bf 32}, (1976), 1--31. 

\bibitem{Mignon-Mann}
E. Mann,
T. Mignon,
{\em Quantum D-modules for toric nef complete intersections},
arXiv:1112.1552 

\bibitem{Mochizuki-tame}
T. Mochizuki,
{\em Asymptotic behaviour of tame harmonic bundles
and an application to pure twistor $D$-modules I, II},
Mem. AMS. {\bf 185}, (2007)

\bibitem{mochi8}
T. Mochizuki,
{\em Asymptotic behaviour of variation of
pure polarized TERP structures},
Publ. RIMS. {\bf 47} (2011), 419--534.

\bibitem{Mochizuki-wild}
T. Mochizuki,
{\em Wild harmonic bundles and wild pure twistor $D$-modules},
Ast\'{e}risque {\bf 340},
Soci\'{e}t\'{e} Math\'{e}matique de France, Paris,
2011.

\bibitem{Mochizuki-Toda-lattice}
T. Mochizuki,
{\em Harmonic bundles and Toda lattices with opposite sign II}, 
Comm. Math. Phys. {\bf 328}, 1159--1198,
DOI:10.1007/s00220-014-1994-0,
the second part of arXiv:1301.1718

\bibitem{Mochizuki-Betti}
T. Mochizuki,
{\em Holonomic D-modules with Betti structure},
arXiv:1001.2336,
to apper in
{\em M\'{e}m. Soc. Math. France}


\bibitem{Mochizuki-MTM}
T. Mochizuki,
{\em Mixed twistor $D$-modules},
Springer, 2015.

\bibitem{Mochizuki-K-complex}
T. Mochizuki,
{\em A twistor approach to
the Kontsevich complexes},
arXiv:1501.0415

\bibitem{peters-steenbrink}
C. Peters and
J. Steenbrink, 
{\em Mixed Hodge structure},
Springer-Verlag, Berlin, 2008. 

\bibitem{Reichelt}
T. Reichelt,
{\em Laurent polynomials, GKZ-hypergeometric systems 
and mixed Hodge modules},
 Compos. Math. {\bf 150}, (2014), 911--941

\bibitem{Reichelt-Sevenheck1}
T. Reichelt, 
C. Sevenheck,
{\em Logarithmic Frobenius manifolds, 
hypergeometric systems and quantum D-modules},
arXiv:1010.2118,
to appear in J. Alg. Geom.


\bibitem{Reichelt-Sevenheck2}
T. Reichelt, 
C. Sevenheck,
{\em Non-affine Landau-Ginzburg models and intersection cohomology},
 \\
arXiv:1210.6527

\bibitem{Reichelt-Sevenheck3}
T. Reichelt, C. Sevenheck,
{\em Hypergeometric Hodge modules},
arXiv:1503.01004

\bibitem{sabbah2}
C. Sabbah,
{\em Polarizable twistor $D$-modules}
Ast\'{e}risque, {\bf 300}, (2005)

\bibitem{Sabbah-tame-polynomial}
C. Sabbah, 
{\em Hypergeometric periods for a tame polynomial},
Port. Math. (N.S.) {\bf 63}, (2006), 173--226.

\bibitem{sabbah-Fourier-Laplace}
C. Sabbah,
{\em Fourier-Laplace transform of 
a variation of polarized complex Hodge structure},
J. Reine Angew. Math. {\bf 621}, (2008), 123--158

\bibitem{sabbah5}
C. Sabbah,
{\em Wild twistor $D$-modules},
in
{\em Algebraic analysis and around},
Adv. Stud. Pure Math., {\bf 54},
Math. Soc. Japan, Tokyo, (2009),
293--353.

\bibitem{Sabbah-Yu}
C. Sabbah,
J.-D. Yu.
{\em On the irregular Hodge filtration
of exponentially twisted mixed Hodge modules},
arXiv:1406.1339

\bibitem{Saito-Takahashi-survey}
K. Saito,
A. Takahashi, 
{\em From primitive forms to Frobenius manifolds},
in
{\em From Hodge theory to integrability and TQFT tt*-geometry,}
31--48, Proc. Sympos. Pure Math., {\bf 78},
Amer. Math. Soc., Providence, RI, 2008. 


\bibitem{saito1}
M. Saito,
{\em Modules de Hodge polarisables},
Publ. RIMS., {\bf 24},
(1988), 849--995.

\bibitem{saito-Brieskorn}
M. Saito,
{\em On the structure of Brieskorn lattice},
Ann. Inst. Fourier (Grenoble), {\bf 39}, (1989), 27--72.

\bibitem{saito3}
M. Saito,
{\em Duality for vanishing cycle functors},
Publ. Res. Inst. Math. Sci.  {\bf 25},  (1989),
889--921.

\bibitem{saito4}
M. Saito,
{\em Induced $D$-modules and 
differential complexes},
Bull. Soc. Math. France {\bf 117}, (1989),
361--387. 

\bibitem{saito2}
M. Saito,
{\em Mixed Hodge modules},
Publ. RIMS., {\bf 26}, (1990),
221--333.

\bibitem{Schulze-Walther}
M. Schulze and U. Walther, 
{\em Hypergeometric D-modules and twisted 
Gauss–Manin systems}, J. Algebra {\bf 322}, (2009), 3392--3409.




\bibitem{s3}
C. Simpson,
{\it Mixed twistor structures},
math.AG/9705006.

\bibitem{steenbrink-zucker}
J. Steenbrink and S. Zucker,
{\em Variation of mixed Hodge structure. I},
Invent. Math.  {\bf 80},  (1985),  489--542.

\bibitem{Stienstra}
J. Stienstra,
{\em Resonant hypergeometric systems 
and mirror symmetry},
Integrable systems and algebraic geometry (Kobe/Kyoto, 1997), 
World Sci. Publ., River Edge, NJ, (1998),
412--452.



\bibitem{Yu}
J.-D. Yu,
{\em Irregular Hodge filtration on twisted de Rham cohomology},
Manuscripta Math. {\bf 144}, (2014),  99--133



\end{thebibliography}
\end{document}